\documentclass[reqno,12pt]{amsart}
\usepackage[T1]{fontenc}
\usepackage{amsfonts}
\usepackage{amsmath}
\usepackage{amsthm}
\usepackage{amssymb}
\usepackage{graphicx}
\usepackage{subcaption}
\usepackage{url}
\usepackage{hyperref}
\usepackage{geometry}


\newtheorem{theorem}[subsection]{Theorem}
\newtheorem{definition}{Definition}[section]
\newtheorem{lemma}{Lemma}[section]
\newtheorem{proposition}{Proposition}[section]

\newtheorem{remark}{Remark}[section]

\newcommand{\weak}{\rightharpoonup}
\newcommand{\weakstar}{\overset{*}{\rightharpoonup}}

\title[A doubly nonlinear parabolic problem]{A doubly nonlinear parabolic problem with variable exponents, homogeneous Neumann boundary conditions and generalized logistic source}
\author{Bogdan Maxim}
\email{maxim.bogdan.n6h@student.ucv.ro}
\date{\today}

\begin{document}

\maketitle

	\begin{center} \footnotesize{Department of Mathematics, University of Craiova, Al. I. Cuza Street, no. 13, 200585, Craiova, Romania}
\end{center}

\begin{abstract}
	The aim of this work is to develop a self-contained existence and uniqueness theory for a doubly nonlinear parabolic problem with variable exponents and homogeneous Neumann boundary conditions. Our approach is based on Rothe's time-discretization method and builds on the results established by the author in \cite{max4}. We also investigate the asymptotic behavior of the solution. All arguments are presented in full detail, so that the paper is self-contained.
\end{abstract}

\bigskip


\smallskip


\tableofcontents

\section{Introduction}

\noindent We consider, for an arbitrary $T\in (0,\infty)$, the following doubly nonlinear parabolic problem:

\begin{equation}\tag{$DNP$}\label{eqdpg}
	\begin{cases}\dfrac{\partial b(x,u(t,x))}{\partial t}-\operatorname{div}\mathbf{a}(x,\nabla u)=f\big (x,u(t,x)\big ), & (t,x)\in (0,T)\times\Omega\\[3mm] \mathbf{a}(x,\nabla u)\cdot\nu=0, & (t,x)\in (0,T)\times\partial\Omega\\[3mm] u(0,x)=u_0(x)\in [\varepsilon,\delta], & x\in\Omega\end{cases}
\end{equation}

\noindent together with the associated stationary problem:

\begin{equation}\tag{$E$}\label{eqedg}
	\begin{cases}-\operatorname{div}\mathbf{a}\big (x,\nabla U(x)\big )=f\big (x,U(x)\big ), & x\in\Omega\\[3mm] \mathbf{a}(x,\nabla U)\cdot\nu =0, & x\in \partial\Omega\\[3mm] \varepsilon\leq U(x)\leq \delta, & x\in\Omega\end{cases}
\end{equation}

\noindent The present work extends the results obtained in \cite{giaco2} and \cite{Giaco1} to a doubly nonlinear parabolic setting. A central feature of the problem is the spatial heterogeneity of the nonlinearity appearing in the time derivative: we allow $b=b(x,s)$, rather than restricting the analysis to the spatially homogeneous case $b=b(s)$. This additional dependence on the spatial variable introduces substantial difficulties and requires a number of arguments that are absent in the homogeneous setting.

\medskip

\noindent In Section \ref{s2}, we introduce the notation and state the hypotheses used throughout the paper. We work on a bounded Lipschitz domain $\Omega\subset\mathbb{R}^N$. Of particular interest is hypothesis \textbf{(H7)}, which generalizes the usual \textbf{Leray-Lions structural condition}:

\begin{equation}\label{eqstruct}
	A(x,\xi)\geq \delta |\xi|^{p(x)}-\tilde{\delta},\ \text{for a.a.}\ x\in\Omega,\ \forall\ \xi\in\mathbb{R}^N,
\end{equation}

\noindent by replacing it with the weaker requirement that there exists a \textbf{coercivity profile}, namely, a nondecreasing function $\gamma:[0,\infty)\to [0,\infty)$, with $\lim\limits_{s\to\infty} \gamma(s)=+\infty$, such that:

\begin{equation}
	\mathcal{A}(V)\geq\gamma\left(\rho_{p(x)}\bigl(|\nabla V|\bigr)\right),\ \forall\ V\in W^{1,p(x)}(\Omega).
\end{equation}

\noindent A further feature of our framework is that we no longer impose assumptions \textbf{(A1)}, \textbf{(A2)}, \textbf{(A3)}, \textbf{(A4)}, \textbf{(A5)}, \textbf{(A6)} and \textbf{(A7)} on $A$ or $\Phi$, which are assumed in the related paper \cite{Giaco1}. That work establishes qualitative results for \eqref{eqdpg} in the homogeneous linear case $b(x,s)=s$. Concerning the flux operator $\mathbf{a}$, we assume only strict monotonicity together with the requirement that its Nemytskii operator maps $L^{p(x)}(\Omega)^N$ into $L^{p'(x)}(\Omega)^N$. This mapping property is essential, since without it the integral appearing in the definition of a weak solution would not, in general, be well defined.

\medskip

\noindent Section \ref{s3} develops a collection of auxiliary results that are used throughout the paper. The definition of a \textbf{weak solution} to \eqref{eqdpg}, introduced in Section \ref{s4}, is inspired by the PDE course \cite{Leoni2014PDEII} and differs from the formulations commonly used in this area, for instance in \cite{antontsev2012existence} and \cite[Definition 2.1]{giaco2}. Whereas those formulations employ test functions belonging to Bochner spaces, our definition uses test functions only from the Sobolev space $W^{1,p(x)}(\Omega)$. A similar formulation is also used in \cite[Theorem 5.1]{Goro2014}. This choice is particularly useful in treating the parabolic problem \eqref{eqdpg} with irregular initial data $u_0\in L^{\infty}(\Omega)$, a situation not covered in \cite{Goro2014}.

\medskip

\noindent In Section \ref{s5}, we introduce one of the main tools of the paper: a \textit{weak parabolic comparison principle} for \eqref{eqdpg}. Its proof is inspired by \cite[Proposition 6]{chaouai2018qualitative}.

\medskip

\noindent Section \ref{s6}, which constitutes a substantial part of the paper, is devoted to the auxiliary problem

\begin{equation}
	\begin{cases}\dfrac{\partial b(x,u(t,x))}{\partial t}-\operatorname{div}\mathbf{a}(x,\nabla u)=g(t,x), & (t,x)\in (0,T)\times\Omega\\[3mm] \mathbf{a}(x,\nabla u)\cdot\nu=0, & (t,x)\in (0,T)\times\partial\Omega\\[3mm] u(0,x)=u_0(x)\in [\varepsilon,\delta], & x\in\Omega\end{cases}
\end{equation}

\noindent where $g$ belongs to a suitable class of functions in $L^{\infty}\bigl((0,T)\times\Omega\bigr)$. Two of the main results of the paper are established in this section. Theorem \ref{thmauxiliar} concerns regular initial data $u_0\in W^{1,p(x)}(\Omega)$. In this case, the solution has the regularity 

$$u\in C\bigl([0,T];L^2(\Omega)\bigr) \cap H^1\bigl((0,T); L^2(\Omega)\bigr)\cap L^{\infty}\bigl(0,T;W^{1,p(x)}(\Omega)\bigr).$$ 

\noindent To prove this theorem, we use a \textbf{Minty-type argument}, originating in the classical work \cite{minty1962monotone}, which allows us to avoid unnecessary additional hypotheses on the flux operator $\textbf{a}$. The second major result of this section is Theorem \ref{theoremverygeneral}, which deals with irregular initial data $u_0\in L^{\infty}(\Omega)$. In this case, we prove that 

$$v\in C\big ( [0,T]; L^2(\Omega)\big )\cap H^1_{\textnormal{loc}}((0,T);L^2(\Omega)\big )\cap L^{1}\big (0,T;W^{1,p(x)}(\Omega)\big ).$$

\medskip

\noindent At this point, it is important to explain why we do not formulate the notion of weak solution exclusively through integrals over the space-time cylinder $(0,T)\times\Omega$. For a general $v\in L^{1}\big (0,T;W^{1,p(x)}(\Omega)\big )$, the integral $\displaystyle\int_{t_1}^{t_2}\int_{\Omega} \mathbf{a}(x,\nabla v(t,x))\cdot\nabla \psi(t,x)\ dx\ dt$ need not be finite, and hence need not be well defined as a Lebesgue integral, even if $\psi\in C^{\infty}_c\bigl((0,T)\times\Omega\bigr)$.

\medskip

\noindent In Section \ref{s7}, we establish in Theorem \ref{theexistenceresult} the existence result for the original problem \eqref{eqdpg} by means of a convergent monotone scheme.

\medskip

\noindent Finally, Section \ref{s8} is devoted to the long-time behavior of solutions. We prove that, if the stationary problem \eqref{eqedg} admits a unique weak solution $U$, then the solution of \eqref{eqdpg} satisfies $u(t,\cdot)\to U$ in $L^r(\Omega)$ as $t\to\infty$ for every $r\in [1,\infty)$; see Theorem \ref{eqasym}.

\medskip

\noindent The results developed in this paper apply, in particular, to the porous medium equation with variable exponents and heterogeneous nonlinear diffusion

\begin{equation}
	\begin{cases}
		\dfrac{\partial}{\partial t}\left( u^{\theta(x)}(t,x)\right)
		-\operatorname{div}\mathbf{a}(x,\nabla u)
		=
		f\bigl(x,u(t,x)\bigr),
		& (t,x)\in(0,T)\times\Omega,
		\\[3mm]
		\mathbf{a}(x,\nabla u)\cdot\nu=0,
		& (t,x)\in(0,T)\times\partial\Omega,
		\\[3mm]
		u(0,x)=u_0(x)\in[\varepsilon,\delta],
		& x\in\Omega.
	\end{cases}
\end{equation}

\noindent The Appendix collects a number of general results used in the paper. Whenever we were unable to find a clear and complete proof in the literature, we provide a full proof for completeness.

\medskip

\noindent The bibliography has been selected with particular emphasis on clear, modern references in which the relevant arguments are presented in an accessible and detailed form.

\medskip

\noindent The level of detail and, consequently, the length of the paper are deliberate. Our aim is to provide a transparent and self-contained treatment, based on modern tools from functional analysis, of the existence theory for weak solutions to a broad class of nonlinear parabolic problems under general structural assumptions. In this way, the paper is also intended to serve as a useful reference for researchers and PhD students working on nonlinear evolutionary PDEs.

\section{Hypotheses and notations}\label{s2}

\noindent We consider the following hypotheses and notations:

\begin{enumerate}
	\item[\textbf{(H1)}] $T>0$ and $\Omega\subset\mathbb{R}^N,\ N\geq 2$ is an open, bounded and connected Lipschitz domain.

	\item[\textbf{(H2)}] $p:\overline{\Omega}\to (1,\infty)$ is a continuous variable exponent with $p^->\dfrac{2N}{N+2}$. In this situation we have the following compact embedding $W^{1,p(x)}(\Omega)\stackrel{c}{\hookrightarrow} L^2(\Omega)$.\footnote{See \cite[Proposition 2.2]{fan2009remarks}, \cite[Theorem 2.2 (c)]{Dinca2} and \cite[Theorem 1.3]{fan2001sobolev}.}
	
	\bigskip
	
	\item[$\bullet$] Denote $p^{-}=\displaystyle\min_{x\in\overline{\Omega}}\ p(x) >1$ and $p^+=\displaystyle\max_{x\in\overline{\Omega}}\  p(x)<\infty$. Let also $p'(x)=\dfrac{p(x)}{p(x)-1}$ be the conjugate variable exponent of $p(x)$.

	\bigskip
	
	\item[\textbf{(H3)}] $\Psi:\overline{\Omega}\times (0,\infty)\to (0,\infty)$ with $\Psi(\cdot, s)$ measurable for each $s\in (0,\infty)$ and $\Psi(x,\cdot)\in\operatorname{AC}_{\text{loc}}\big ((0,\infty)\big )$ for a.e. $x\in\Omega$.

	\item[\textbf{(H4)}] $\lim\limits_{s\to 0^+} \Psi(x,s)s=0$ for a.e. $x\in\Omega$.
	
	\item[\textbf{(H5)}] For a.e. $x\in\Omega$ we have that $(0,\infty)\ni s\longmapsto \Psi(x,s)s\in (0,\infty)$ is a strictly increasing function.
	
	\bigskip
	
	\begin{remark}\label{remmax4}
		From Proposition 3.1 \textbf{(1)} given in \cite{max4} we have that $\mathbf{a}$ is a \textbf{strictly monotone flux operator}, i.e.:
		
		\begin{equation}
			\bigl(\mathbf{a}(x,\xi_1)-\mathbf{a}(x,\xi_2)\bigr)\cdot (\xi_1-\xi_2)\geq 0,\ \forall\ x\in\overline{\Omega}\ \text{and}\ \forall\ \xi_1,\xi_2\in\mathbb{R}^N,
		\end{equation}
		
		\noindent with equality iff $\xi_1=\xi_2$.
	\end{remark}

	\item[$\bullet$] We take $\mathbf{a}:\overline{\Omega}\times\mathbb{R}^N\to\mathbb{R}^N,\ \mathbf{a}(x,\xi)=\Psi(x,|\xi|)\xi$ if $\xi\neq \mathbf{0}$ and $\mathbf{a}(x,\textbf{0})=\textbf{0}$. 
	\bigskip
	
	\item[$\bullet$] Define $\Phi:\overline{\Omega}\times\mathbb{R}\to [0,\infty),\ \Phi(x,s)=\begin{cases} \Psi(x,|s|)|s|, & s\neq 0\\ 0, & s=0\end{cases}$.
	
	\bigskip
	
	\item[\textbf{(H6)}] $L^{p(x)}(\Omega)\ni v\longmapsto \Phi\big(\cdot,v(\cdot)\big )\in L^{p'(x)}(\Omega)$. This is the same as saying that the Nemytskii operator of $\Phi$ can be defined as follows $\mathcal{N}_{\Phi}:L^{p(x)}(\Omega)\to L^{p'(x)}(\Omega)$. An other equivalent version would be: $L^{p(x)}(\Omega)^N\ni\mathbf{v}\longmapsto\mathbf{a}\big (\cdot,\mathbf{v}(\cdot)\big )\in L^{p'(x)}(\Omega)^N$.\footnote{On the space $L^{p(x)}(\Omega)^N$ we consider the norm $\Vert \mathbf{v}\Vert_{L^{p(x)}(\Omega)^N}:=\Vert |\mathbf{v}|\Vert_{L^{p(x)}(\Omega)}$, for any $\mathbf{v}\in L^{p(x)}(\Omega)$.}
	
	\bigskip
	
	\item[$\bullet$] $A:\overline{\Omega}\times\mathbb{R}^N\to [0,\infty),\ A(x,\xi):=\displaystyle\int_{0}^{|\xi|}\Phi(x,s)\ ds$. Note that $A$ is well-defined since $\Phi(x,\cdot)$ is continuous on $[0,\infty)$ for a.e. $x\in\Omega$.

	\item[$\bullet$] $\mathcal{A}:W^{1,p(x)}(\Omega)\to [0,\infty)$ given by $\mathcal{A}(V)=\displaystyle\int_{\Omega} A(x,\nabla V(x))\ dx$.\footnote{Look at Proposition 3 \textbf{(6)} from \cite{max2} to see why $\mathcal{A}$ is well-defined.}
	
	\begin{remark}\label{remmathcalA}
		For $V\in W^{1,p(x)}(\Omega)$ we have that
		
		\begin{align*}
			&\mathcal{A}(V)=0\ \Longleftrightarrow\ \int_{\Omega} \underbrace{A(x,\nabla V(x))}_{\geq 0}\ dx=0\ \Longleftrightarrow\ A(x,\nabla V(x))=0,\ \text{for a.e.}\ x\in\Omega\\
			&\Longleftrightarrow\ \int_{0}^{|\nabla V(x)|}\underbrace{\Phi(x,s)}_{>0,\ \forall\ s>0}\ ds=0,\ \text{for a.e.}\ x\in\Omega\ \Longleftrightarrow\ |\nabla V(x)|=0\ \text{for a.e.}\ x\in\Omega\\
			&\Longleftrightarrow\ \nabla V(x)=0\ \text{for a.e.}\ x\in\Omega.
		\end{align*}
		
		\noindent Since $\Omega$ is an open and connected domain we deduce, using Corollary 2.1.9 from \cite[page 47]{ziemer1989weakly} and the fact that $W^{1,p(x)}(\Omega)\subset W^{1,1}(\Omega)$,\footnote{See \cite[relation (3.2), page 604]{kovacik}.} that there is some constant $c\in\mathbb{R}$ with $V(x)=c$ for a.e. $x\in\Omega$.
	\end{remark}
	
	\bigskip
	
	\item[$\bullet$] $\overline{\mathcal{A}}:L^2(\Omega)\to [0,\infty]$ given by $\overline{\mathcal{A}}(V):=\begin{cases} \mathcal{A}(V)=\displaystyle\int_{\Omega} A(x,\nabla V(x))\ dx,& V\in W^{1,p(x)}(\Omega)\\[3mm] +\infty,& V\in L^2(\Omega)\setminus W^{1,p(x)}(\Omega)\end{cases}$.

	\bigskip
	
	\item[\textbf{(H7)}] For any sequence $(V_n)_{n\geq 1}\subset W^{1,p(x)}(\Omega)$ with $\displaystyle\int_{\Omega}|\nabla V_n|^{p(x)}\ dx\to\infty$ we have that: $\lim\limits_{n\to\infty}\mathcal{A}(V_n)=\lim\limits_{n\to\infty} \displaystyle\int_{\Omega} A(x,\nabla V_n) =\infty$.\footnote{See Proposition \ref{propomathcalA1} \textbf{(5)} for a characterization of \textbf{(H7)}.}
	
	\bigskip
	
	\begin{remark}\label{remboundednormmodular} Note that if $(V_n)_{n\geq 1}\subset W^{1,p(x)}(\Omega)$ has the property that $\bigl(\mathcal{A}(V_n)\bigr)_{n\geq 1}$ is bounded, then the sequence $\bigl(\rho_{p(x)}(|\nabla V_n|)\bigr)_{n\geq 1}$ is also bounded, i.e. the sequence $\left (\Vert \nabla V_n\Vert_{L^{p(x)}(\Omega)^N}\right )_{n\geq 1}$ is bounded.\footnote{Indeed, if we suppose that $\bigl(\rho_{p(x)}(|\nabla V_n|)\bigr)_{n\geq 1}$ is unbounded, we will get a subsequence with $\rho_{p(x)}(|\nabla V_{n_k}|)\to \infty$ as $k\to\infty$. But from \textbf{(H7)}, this implies that $\mathcal{A}(V_{n_k})\to\infty$, which is a contradiction. See also Lemma 3.2.5 from \cite[page 73]{Hasto}.} 
	\end{remark}
	
	\bigskip
	
	\item[\textbf{(H8)}] $f:\overline{\Omega}\times [\varepsilon,\delta]\to\mathbb{R}$ is a 
	\textbf{generalized logistic source}, meaning that:
	
	\medskip
	
	\begin{itemize}
		
		\item $f$ is a Carath\'{e}odory function, i.e. $[\varepsilon,\delta]\ni s\mapsto f(x,s)$ is continuous for a.e. $x\in\Omega$ and $\Omega\ni x\mapsto f(x,s)$ is measurable for any $s\in [\varepsilon,\delta]$. 
		
		\item  $0\leq\varepsilon<\delta$.\footnote{It is very important that $\varepsilon$ is allowed to be $0$.}
		
		\item $f(x,\varepsilon)\geq 0$ and $f(x,\delta)\leq 0$ a.e. on $\Omega$.
	\end{itemize}
	
	\medskip
	
	
	\item[\textbf{(H9)}] $b\in C^1\big (\Omega\times (0,\delta)\big )\cap C\big (\overline{\Omega}\times [0,\delta] \big )$.
	
	\item[\textbf{(H10)}] $b(x,0)=0,\ \forall\ x\in\overline{\Omega}$.\footnote{As we are interested in only in the derivative with respect to the second argument of $b$, we can assume this condition without losing the generality of the problem. Otherwise we may take $b(x,s)-b(x,0)$ instead of $b(x,s)$.}
	
	\begin{remark}\label{remdbdxi} It is important to point out that, when the second argument is $0$, all partial derivatives of $b$ with respect to the spatial variable exist and are zero, i.e. for each $i\in\{1,2,...,N\}$: $\dfrac{\partial b}{\partial x_i}(x,0)=\lim\limits_{h\to 0}\dfrac{b(x+he_i,0)-b(x,0)}{h}=0$, for any $x\in\Omega$.
	\end{remark}
	
	\item[\textbf{(H11)}] For any $(x_0,s_0)\in\partial\big (\Omega\times (0,\delta)\big )$ the following limit exists\footnote{We have that $\partial\big (\Omega\times (0,\delta)\big )=\big (\overline{\Omega}\times\{0,\delta\}\big )\cup \big (\partial\Omega\times [0,\delta]\big )$.} and it is finite $\lim\limits_{\underset{(x,s)\in\Omega\times (0,\delta)}{(x,s)\to (x_0,s_0)}}\dfrac{\partial b}{\partial s}(x,s)$.
	
	\item[\textbf{(H12)}] For each $i\in\overline{1,N}$ and for any $(x_0,s_0)\in\partial\big (\Omega\times (0,\delta)\big )$ the following limit exists and it is finite $\lim\limits_{\underset{(x,s)\in\Omega\times (0,\delta)}{(x,s)\to (x_0,s_0)}}\dfrac{\partial b}{\partial x_i}(x,s)$.

	\item[\textbf{(H13)}] There is some $\ell_0>0$ such that: $\dfrac{\partial b}{\partial s}(x,s)\geq \ell_0$ for any $(x,s)\in\overline{\Omega}\times [\varepsilon,\delta]$.

\begin{remark}\label{rem23}
	In particular we have for any $x\in\overline{\Omega}$ that the function $[\varepsilon,\delta]\ni s\mapsto b(x,s)$ is strictly increasing. We also get from the Mean Value Theorem, that:
	
	\begin{equation}
		|b(x,s_1)-b(x,s_2)|\geq \ell_0 |s_1-s_2|,\ \forall\ x\in\overline{\Omega},\ \forall\ s_1,s_2\in [\varepsilon,\delta].
	\end{equation}
\end{remark}

\begin{remark} Even though from \textnormal{\textbf{(H10)}} we know that $b(x,0)=0$ on $\overline{\Omega}$ and from Remark \ref{rem23} we know that $[\varepsilon,\delta]\ni s\mapsto b(x,s)$ is strictly increasing, we may have $b(x,\varepsilon)\leq 0$ or $b(x,\delta)\leq 0$ in some subsets of $\overline{\Omega}$ or even on the entire set $\overline{\Omega}$. Therefore we have that:
	
	\begin{equation}
		b(x,\varepsilon)\leq b(x,s)\leq b(x,\delta),\ \forall\ s\in [\varepsilon,\delta], 
	\end{equation}
	
	\noindent from where, taking into account that $b(\cdot,\varepsilon),b(\cdot,\delta)\in C(\overline{\Omega})\subset L^{\infty}(\Omega)$, we can write:
	
	\begin{align}
		|b(x,s)|&\leq\max\{|b(x,\varepsilon)|,|b(x,\delta)|\}\nonumber\\
		&\leq \max\bigl\{\Vert b(\cdot,\varepsilon)\Vert_{L^{\infty}(\Omega)},\Vert b(\cdot,\delta)\Vert_{L^{\infty}(\Omega)}\bigr\},\ \forall\ s\in [\varepsilon,\delta].\label{bequ1}
	\end{align}
	
\end{remark}
	
	\item[\textbf{(H14)}] There is a constant $\lambda_0>0$ such that the function $[\varepsilon,\delta]\ni s\mapsto f(x,s)+\lambda_0 b(x,s)$ is strictly increasing for a.e. $x\in\Omega$.
	
	\item[\textbf{(H15)}] $u_0\in\mathcal{U}_{[\varepsilon,\delta]}:=\{U\in L^\infty(\Omega)\ |\ \varepsilon\leq U\leq \delta\ \text{a.e. on}\ \Omega\}$.
%
%
	
		\bigskip
	
		\noindent In what follows we will introduce some extra hypothesis that shall be useful only when we are dealing with the uniqueness problem for \eqref{eqedg} and \eqref{eqdpg}:
	
	\bigskip


	\item[\textbf{(EH$_\Phi$)}] There is some $\alpha\in (1,2)$ with $p^-\geq\alpha$ such that for a.e. $x\in\Omega$ the function $(0,\infty)\ni s\mapsto\dfrac{\Phi(x,s)}{s^{\alpha-1}}$ is \textbf{strictly increasing} and the function $(0,\delta^{\alpha}]\ni s\longmapsto \dfrac{f(x,\sqrt[\alpha]{s})}{\sqrt[\alpha]{s^{\alpha-1}}}$ is decreasing.
	
	\item[\textbf{(EH$_f$)}] There is some $\alpha\in (1,2)$ with $p^-\geq\alpha$ such that for a.e. $x\in\Omega$ the function $(0,\delta^{\alpha}]\ni s\longmapsto \dfrac{f(x,\sqrt[\alpha]{s})}{\sqrt[\alpha]{s^{\alpha-1}}}$ is \textbf{strictly decreasing} and the function $(0,\infty)\ni s\mapsto\dfrac{\Phi(x,s)}{s^{\alpha-1}}$ is increasing.
	
	\item[\textbf{(EHU)}] There is some $L_f>0$ such that for a.e. $x\in\Omega$ we have that:
	
	\begin{equation}
		|f(x,s_2)-f(x,s_1)|\leq L_f |s_2-s_1|,\ \forall\ s_1,s_2\in [\varepsilon,\delta].
	\end{equation}

\end{enumerate}
	
	\bigskip
	\bigskip
	
	\noindent In this article we will use the following notations:
	
	\bigskip
	
	\begin{enumerate}
	
	\item[$\bullet$] $\mathcal{U}_{[\varepsilon,\delta]}=\big\{U\in L^{\infty}(\Omega)\ |\ \varepsilon\leq U\leq \delta\ \text{a.e. on}\ \Omega \big \}$.
	
	\medskip
	
	\item[$\bullet$] $\mathcal{V}_{[\varepsilon,\delta]}=\big\{v\in L^{\infty}\big ((0,T)\times\Omega\big )\ \big |\ \varepsilon\leq v\leq\delta\ \text{a.e. on}\ (0,T)\times\Omega\big \}$.
	
	\medskip
	
	\noindent For any $\lambda>0$ we denote:
	
	\medskip
	
	\item[$\bullet$] $\mathcal{M}^{[\varepsilon,\delta]}_{\lambda}:=\big \{g:\Omega\to\mathbb{R}\ \text{measurable}\ |\ \lambda b(x,\varepsilon)\leq g(x)\leq \lambda b(x,\delta)\ \text{for a.e.}\ x\in\Omega\big \}$.
	
	\medskip
	
	\item[$\bullet$] $\mathfrak{M}^{[\varepsilon,\delta]}_{\lambda}:=\big \{g:(0,T)\times\Omega\to\mathbb{R}\ \text{measurable}\ |\ \lambda b(x,\varepsilon)\leq g(t,x)\leq \lambda b(x,\delta)\ \text{a.e. on}\ (0,T)\times\Omega\big \}$.
	
	\begin{remark}\label{rem221}
		From hypothesis \textnormal{\textbf{(H9)}} we get that $\mathcal{M}_{\lambda}^{[\varepsilon,\delta]}\subseteq L^{\infty}(\Omega)$ and $\mathfrak{M}_{\lambda}^{[\varepsilon,\delta]}\subseteq L^{\infty}\big ((0,T)\times\Omega\big )$, because $b:\overline{\Omega}\times [0,\delta]\to\mathbb{R}$ is a continuous function on a compact set, and hence it is bounded.\footnote{See the Weierstrass theorem in \cite[Theorem 5.1.1, page 162]{kurdila2005convex}.}
	\end{remark}
	
	\begin{remark}\label{remark22} Note that if $\lambda\geq \lambda_0$ then for each $U\in\mathcal{U}_{[\varepsilon,\delta]}$ we get that $\Omega\ni x\mapsto g(x):=f(x,U(x))+\lambda b(x,U(x))\in\mathcal{M}^{[\varepsilon,\delta]}_{\lambda}$. Indeed, for a.e. $x\in\Omega$ we have that:
		
		\begin{equation}
			\lambda b(x,\varepsilon)\leq f(x,\varepsilon)+\lambda b(x,\varepsilon)\leq g(x)\leq f(x,\delta)+\lambda b(x,\delta)\leq \lambda b(x,\delta).
		\end{equation}
	\end{remark}
	
		\begin{remark}\label{remark222} Similarly, if $\lambda\geq \lambda_0$ then for each $u\in\mathcal{V}_{[\varepsilon,\delta]}$ we get that $(0,T)\times\Omega\ni (t,x)\mapsto g(t,x):=f(x,u(t,x))+\lambda b(x,u(t,x))\in\mathfrak{M}^{[\varepsilon,\delta]}_{\lambda}$. Indeed, for a.e. $x\in\Omega$ we have that:
		
		\begin{equation}
			\lambda b(x,\varepsilon)\leq f(x,\varepsilon)+\lambda b(x,\varepsilon)\leq g(t,x)\leq f(x,\delta)+\lambda b(x,\delta)\leq \lambda b(x,\delta).
		\end{equation}
	\end{remark}
	
%
	
	\item[$\bullet$] $\overline{f}:\overline{\Omega}\times\mathbb{R}\to\mathbb{R}$, $\overline{f}(x,s)=\begin{cases} f(x,\varepsilon)-\dfrac{\lambda_0}{2}(s-\varepsilon), & s\in (-\infty,\varepsilon)\\[3mm] f(x,s), & s\in [\varepsilon,\delta] \\[3mm] f(x,\delta)-\tilde{\lambda}_0(s-\delta), & s\in (\delta,\infty) \end{cases}$, where $0\leq\tilde{\lambda}_0<\lambda_0$ is any fixed constant.
	
	\item[$\bullet$] $\overline{b}:\overline{\Omega}\times\mathbb{R}\to\mathbb{R},\  \overline{b}(x,s):=\begin{cases}b(x,\varepsilon)+s-\varepsilon, & s<\varepsilon\\[3mm] b(x,s), & s\in [\varepsilon,\delta]\\[3mm] b(x,\delta)+s-\delta, & s>\delta \end{cases}$.
	
	\item[$\bullet$] $g_0:\overline{\Omega}\times\mathbb{R}\to\mathbb{R},\ g_0(x,s)=\overline{f}(x,s)+\lambda_0 \overline{b}(x,s)$.
	
	\medskip
	
	\item[$\bullet$] $\frak{B}:\overline{\Omega}\times [\varepsilon,\delta]\to\mathbb{R},\ \frak{B}(x,s)=\displaystyle\int_{\varepsilon}^sr\cdot\dfrac{\partial b}{\partial s}(x,r)\ dr$.
	
\end{enumerate}
\section{Preliminary results}\label{s3}

\subsection{Properties of $L^{p(x)}(\Omega)$, $L^{p'(x)}(\Omega)$ and $W^{1,p(x)}(\Omega)$}

\begin{proposition}\label{propospatii} The follwoing facts are valid:
	
\begin{enumerate}
	\item[\textnormal{\textbf{(1)}}] $L^{p(x)}(\Omega)$, $L^{p'(x)}(\Omega)$ and $W^{1,p(x)}(\Omega)$ are reflexive and separable Banach spaces.
	
	\item[\textnormal{\textbf{(2)}}] $L^{p(x)}(\Omega)^N$ and $L^{p'(x)}(\Omega)^N$ are also reflexive and separable Banach spaces.
	
	\item[\textnormal{\textbf{(3)}}] There is some constant $\frak{C}>0$ such that:
	
	\begin{equation}\label{eqcontrolWnorm}
		\Vert V\Vert_{W^{1,p(x)}(\Omega)}\leq\frak{C}\bigl(\Vert \nabla V\Vert_{L^{p(x)}(\Omega)^N}+\Vert V\Vert_{L^2(\Omega)}\bigr),\ \forall\ V\in W^{1,p(x)}(\Omega).
	\end{equation}
\end{enumerate}
	
\end{proposition}

\begin{proof} \noindent\textbf{(1)} From $1<p^-\leq p^+<\infty$ we deduce from \cite[Theorem 2.7]{kovacik} that $L^{p(x)}(\Omega)$ is a reflexive Banach space. The exponent $p$ being bounded (because it is a continuous function defined on the compact set $\overline{\Omega}$) we also obtain from \cite[Corollary 2.12]{kovacik} that $L^{p(x)}(\Omega)$ is separable.

\noindent In a similar fashion, from $+\infty>p^+\geq p(x)\geq p^{-}>1$ for any $x\in\overline{\Omega}$, we infer from \textbf{(H2)} that $p'(x):\overline{\Omega}\to (1,\infty),\ p'(x):=\dfrac{p(x)}{p(x)-1}$ is a continuous (and hence bounded) exponent, with $(p')^{-}\geq \dfrac{p^+}{p^+-1}>1$ and $(p')^+\leq\dfrac{p^-}{p^--1}<\infty$. Therefore, from \cite[Corollary 2.12]{kovacik}, we have that $L^{p'(x)}(\Omega)$ is separable, and from \cite[Corollary 2.7]{kovacik} we get that $L^{p'(x)}(\Omega)$ is also reflexive.  

\noindent From \cite[Theorem 3.1]{kovacik} we get that $W^{1,p(x)}(\Omega)$ is a reflexive and separable Banach space.
	
\bigskip

\noindent\textbf{(2)} From \cite[Corollary 1.11.20, page 105]{megginson2012introduction} it follows that $L^{p(x)}(\Omega)^N$ and $L^{p'(x)}(\Omega)^N$ are reflexive. In the same way, from \cite[Proposition 1.12.9 (f)]{megginson2012introduction}, we obtain that $L^{p(x)}(\Omega)^N$ and $L^{p'(x)}(\Omega)^N$ are separable Banach spaces.

\bigskip

\noindent\textbf{(3)} Suppose that we cannot find such $\frak{C}$. Therefore, for each $n\geq 1$ there is some element $V_n\in W^{1,p(x)}(\Omega)$ such that:

\begin{equation}
	\Vert V_n\Vert_{W^{1,p(x)}(\Omega)}>n\bigl(\Vert \nabla V_n\Vert_{L^{p(x)}(\Omega)^N}+\Vert V_n\Vert_{L^2(\Omega)} \bigr).
\end{equation}	

\noindent We denote $U_n=\dfrac{V_n}{\Vert V_n\Vert_{W^{1,p(x)}(\Omega)}},\ n\geq 1$. Hence the above inequality can be rewritten for each $n\geq 1$ as:

\begin{equation}\label{taticu1}
	\Vert \nabla U_n\Vert_{L^{p(x)}(\Omega)^N}+\Vert U_n\Vert_{L^2(\Omega)}<\dfrac{1}{n},\ \text{where}\ \Vert U_n\Vert_{L^{p(x)}(\Omega)}=1.
\end{equation}

\noindent This proves that $\lim\limits_{n\to\infty} \Vert U_n\Vert_{L^2(\Omega)}=0$, i.e. $U_n\to 0$ in $L^2(\Omega)$. Also from \eqref{taticu1} we get that $\lim\limits_{n\to\infty}	\Vert \nabla U_n\Vert_{L^{p(x)}(\Omega)^N}=0$, which in particular means that the sequence $\left( \bigl\Vert \nabla U_n\bigr\Vert_{L^{p(x)}(\Omega)^N}\right)_{n\geq 1}$ is bounded. Therefore, for each $n\geq 1$ we have that:

%

\begin{equation}
	\Vert U_n\Vert_{W^{1,p(x)}(\Omega)}=\Vert U_n\Vert_{L^{p(x)}(\Omega)}+\Vert \nabla U_n\Vert_{L^{p(x)}(\Omega)^N}=1+\Vert \nabla U_n\Vert_{L^{p(x)}(\Omega)^N}\leq 1+\sup_{n\geq 1}\Vert \nabla U_n\Vert_{L^{p(x)}(\Omega)^N}<\infty,
\end{equation}

\noindent which means that $\bigl(U_n\bigr)_{n\geq 1}$ is bounded in $W^{1,p(x)}(\Omega)$.

\medskip

\noindent Note that from \cite[Proposition 2.2]{fan2009remarks} we have that $W^{1,p(x)}(\Omega)\stackrel{c}{\hookrightarrow} L^{p(x)}(\Omega)$. As $\bigl(U_n\bigr)_{n\geq 1}$ is bounded in $W^{1,p(x)}(\Omega)$, we get from the compactness of this embedding that there is some $U\in L^{p(x)}(\Omega)$ and a subsequence $\bigl(U_{n_k}\bigr)_{k\geq 1}$ such that:\footnote{For the definition and properties of compact operators see \cite[Section 3.1, page 266]{gasinski2005nonlinear}.}

\begin{equation}
	U_{n_k}\stackrel{k\to\infty}{\longrightarrow} U,\ \text{in}\ L^{p(x)}(\Omega).
\end{equation}

\noindent But now we have that $\begin{cases} U_{n_k}\to 0, & \text{in}\ L^2(\Omega)\hookrightarrow L^1(\Omega)\\  U_{n_k}\to U, & \text{in}\ L^{p(x)}(\Omega)\hookrightarrow L^1(\Omega)\end{cases}$, which gives us that $\begin{cases} U_{n_k}\to 0, & \text{in}\ L^1(\Omega)\\  U_{n_k}\to U, & \text{in}\ L^{1}(\Omega)\end{cases}$. Therefore $U=0$ a.e. on $\Omega$.

\medskip

\noindent We now make the final step in order to obtain the desired contradiction. Since $\Vert U_{n_k}\Vert_{L^{p(x)}(\Omega)}=1,\ \forall\ k\geq 1$ and $U_{n_k}\to U$ in $L^{p(x)}(\Omega)$, we get, using the continuity of the norm\footnote{In any normed vector space $\bigl(X,\Vert\cdot\Vert_X\bigr)$ the norm is 1-Lipschitz continuous, because from the triangle inequality one has $\big|\Vert x_1\Vert_X-\Vert x_2\Vert_X \big|\leq \Vert x_1-x_2\Vert_X,\ \forall\ x_1,x_2\in X$.}, that 

\begin{equation}
0=\Vert U\Vert_{L^{p(x)}(\Omega)}=\lim\limits_{k\to\infty}  \Vert U_{n_k}\Vert_{L^{p(x)}(\Omega)}=1,\ \text{false}.
\end{equation}

\end{proof}

\subsection{Properties of $\mathbf{a}$}

\begin{proposition}\label{propoa0} For each $\mathbf{v}\in L^{p(x)}(\Omega)^N$ we define the function $\mathcal{N}_{\mathbf{a}}(\mathbf{v}):\Omega\to\mathbb{R}^N$ by:
	
	\begin{equation}
		\mathcal{N}_{\mathbf{a}}(\mathbf{v})(x):=\mathbf{a}(x,\mathbf{v}(x)),\ x\in\Omega.
	\end{equation}

\noindent Then the following assertions hold:

\begin{enumerate}
	\item[\textnormal{\textbf{(1)}}] $\mathbf{a}:\overline{\Omega}\times\mathbb{R}^N\to\mathbb{R}^N$ is a Carath\'{e}odory function.
	
	\item[\textnormal{\textbf{(2)}}] $\mathcal{N}_{\mathbf{a}}:L^{p(x)}(\Omega)^N\to L^{p'(x)}(\Omega)^N$ is well-defined.
	
	\item[\textnormal{\textbf{(3)}}] There is some non-negative function $a\in L^{p'(x)}(\Omega)^+$ and some constant $b\geq 0$ such that
	
	\begin{equation}\label{neimportant7}
		|\mathbf{a}(x,\xi)|\leq a(x)+b|\xi|^{p(x)-1},\ \forall\ \xi\in\mathbb{R}^N\ \text{and for a.e.}\ x\in\Omega.
	\end{equation}
	
	\item[\textnormal{\textbf{(4)}}] $\mathcal{N}_{\mathbf{a}}$ is a bounded nonlinear operator, i.e. it maps bounded subsets of $L^{p(x)}(\Omega)^N$ into bounded subsets of $L^{p'(x)}(\Omega)^N$. Moreover, the following inequality takes place:
	
	\begin{equation}\label{impo1}
		\Vert \mathbf{a}(\cdot,\mathbf{v}(\cdot))\Vert_{L^{p'(x)}(\Omega)^N}\leq \Vert a\Vert_{L^{p'(x)}(\Omega)}+b\max\bigl\{\Vert \mathbf{v}\Vert_{L^{p(x)}(\Omega)^N}^{p^+-1},\Vert \mathbf{v}\Vert_{L^{p(x)}(\Omega)^N}^{p^--1}\bigr\},\ \forall\ \mathbf{v}\in L^{p(x)}(\Omega)^N.
	\end{equation}
	
	\item[\textnormal{\textbf{(5)}}] $\mathcal{N}_{\mathbf{a}}$ is norm-to-norm continuous. Namely,
	
	\begin{equation}
		\mathbf{v}_n\to \mathbf{v}\ \text{in}\ L^{p(x)}(\Omega)^N\ \Longrightarrow\ \mathbf{a}(\cdot,\mathbf{v}_n)\to\mathbf{a}(\cdot,\mathbf{v})\ \text{in}\ L^{p'(x)}(\Omega)^N.
	\end{equation}

\end{enumerate}
	
\end{proposition}

\begin{proof}\noindent\textbf{(1)} For a.e. $x\in\Omega$, from \textbf{(H3)}, we have that $\Psi(x,\cdot)$ is continuous on $(0,\infty)$. Therefore $\xi\mapsto\mathbf{a}(x,\xi)=\Psi(x,|\xi|)\xi$ is continuous at every $\xi\neq \mathbf{0}$. At $\mathbf{0}$, from \textbf{(H4)} we deduce that:
	
	\begin{equation}
		|\mathbf{a}(x,\xi)|=\Psi(x,|\xi|)|\xi|\stackrel{\xi\to\mathbf{0}}{\longrightarrow}\mathbf{0}=\mathbf{a}(x,\mathbf{0}).
	\end{equation}
	
	\noindent Hence for a.e. $x\in\Omega$ the mapping $\mathbb{R}^N\ni\xi\mapsto\mathbf{a}(x,\xi)$ is continuous on all $\mathbb{R}^N$.
	
	\noindent Fix now any $\xi\in\mathbb{R}^N$. If $\xi=\mathbf{0}$ then $\Omega\ni x\mapsto\mathbf{a}(x,\mathbf{0})\equiv\mathbf{0}$ is measurable, being identically zero. If $\xi\neq\mathbf{0}$ then, from \textbf{(H3)}, the mapping $\Omega\ni x\mapsto\Psi(x,|\xi|)$ is measurable. Thus the mapping $\Omega\ni x\mapsto\mathbf{a}(x,\xi)=\Psi(x,|\xi|)\xi$ is also measurable. In conclusion $\mathbf{a}:\overline{\Omega}\times\mathbb{R}^N\to\mathbb{R}^N$ is a Carath\'{e}odory function.

	\bigskip
	
	\noindent\textbf{(2)} Since $\mathbf{a}:\Omega\times\mathbb{R}^N\to\mathbb{R}^N$ is a Carath\'{e}odory function, using Corollary 2.5.24 from the book \cite[page 190]{denkowski2013introduction}, we get that $\mathcal{N}_{\mathbf{a}}(\mathbf{v})$ is a measurable function.
	
	\noindent Now, since for any $\mathbf{v}\in L^{p(x)}(\Omega)^N$ we have that $|\mathbf{v}|\in L^{p(x)}(\Omega)$, we get from \textbf{(H6)} that
	
	\[
	|\mathcal{N}_{\mathbf{a}}(\mathbf{v})|=|\mathbf{a}(\cdot,\mathbf{v}(\cdot))|=\Phi(\cdot,|\mathbf{v}(\cdot)|)=\mathcal{N}_{\Phi}(|\mathbf{v}|)\in L^{p'(x)}(\Omega).
	\]
	
	\noindent This means that $\mathcal{N}_{\mathbf{a}}(\mathbf{v})\in L^{p'(x)}(\Omega)^N$ as needed.

	\bigskip
	
	\noindent\textbf{(3)} Note that for each $\xi\in\mathbb{R}^N$ we have that $|\mathbf{a}(x,\xi)|=\Phi(x,|\xi|),\ \forall\ x\in\Omega$. From \textbf{(H6)} we know that $\mathcal{N}_{\Phi}:L^{p(x)}(\Omega)\to L^{p'(x)}(\Omega)$ and therefore, using Theorem 1.16 from \cite{fan2001spaces} we get that there is some function $a\in L^{p'(x)}(\Omega)$ and some constant $b\geq 0$ such that:
	
	\begin{equation}
		\Phi(x,s)\leq a(x)+b|s|^{p(x)-1},\ \forall\ s\in\mathbb{R}\ \text{and for a.e.}\ x\in\Omega.
	\end{equation}
	
	\noindent So, for each $\xi\in\mathbb{R}^N$ we have that
	
	\begin{equation}
	|\mathbf{a}(x,\xi)|=\Phi(x,|\xi|)\leq a(x)+b|\xi|^{p(x)-1},\ \text{and for a.e.}\ x\in\Omega.
	\end{equation}
	
	\bigskip
	
	\noindent\textbf{(4)} Take any $\mathbf{v}\in L^{p(x)}(\Omega)^N$. This means that $|\mathbf{v}|\in L^{p(x)}(\Omega)$. Just remark that:
	
	\begin{align*}
	\Vert\mathcal{N}_{\mathbf{a}}(\mathbf{v})\Vert_{L^{p'(x)}(\Omega)^N}&=\Vert \mathbf{a}(\cdot,\mathbf{v}(\cdot))\Vert_{L^{p'(x)}(\Omega)^N}=\Vert |\mathbf{a}(\cdot,\mathbf{v}(\cdot))|\Vert_{L^{p'(x)}(\Omega)}\\
	\text{Proposition \ref{propu2u1} and \eqref{neimportant7}}\ \ \ &\leq \left\Vert a(x)+b|\mathbf{v}|^{p(x)-1}\right\Vert_{L^{p'(x)}(\Omega)}\\
\text{(Norm inequality)}\ \ \ 	&\leq \Vert a\Vert_{L^{p'(x)}(\Omega)}+b\left\Vert |\mathbf{v}|^{p(x)-1} \right\Vert_{L^{p'(x)}(\Omega)}\\
\text{Proposition \ref{lpxprimineq}}\ \ \ &\leq \Vert a\Vert_{L^{p'(x)}(\Omega)}+b\max\left\{\Vert |\mathbf{v}|\Vert_{L^{p(x)}(\Omega)}^{p^+-1},\Vert |\mathbf{v}|\Vert_{L^{p(x)}(\Omega)}^{p^--1}\right\}\\
&=\Vert a\Vert_{L^{p'(x)}(\Omega)}+b\max\left\{\Vert \mathbf{v}\Vert_{L^{p(x)}(\Omega)^N}^{p^+-1},\Vert \mathbf{v}\Vert_{L^{p(x)}(\Omega)^N}^{p^--1}\right\}.
	\end{align*}
	
	\noindent This inequality shows in particular that bounded sets from $L^{p(x)}(\Omega)^N$ are taken into bounded set of $L^{p'(x)}(\Omega)^N$ by the Nemytskii operator $\mathcal{N}_{\mathbf{a}}$.
	
		\bigskip
	
	\noindent\textbf{(5)} \noindent  Define $\Phi_\mathbf{v}:\overline{\Omega}\times\mathbb{R}^N\to\mathbb{R},\ \Phi_{\mathbf{v}}(x,\mathbf{\xi})=|\mathbf{a}(x,\mathbf{\xi})-\mathbf{a}(x,\mathbf{v}(x))|,\ \forall\ x\in\overline{\Omega},\ \forall\ \mathbf{\xi}\in\mathbb{R}^N$ which is clearly a Carath\'{e}odory function. The Nemytskii operator associated to $\Phi_{\mathbf{v}}$ is $\mathcal{N}_{\Phi_\mathbf{v}}:L^{p(x)}(\Omega)^N\to L^{p'(x)}(\Omega)$. Indeed, since for all $\mathbf{w}\in L^{p(x)}(\Omega)^N$ we have from \textbf{(H6)} that:
	
	\begin{equation*}
		\mathcal{N}_{\Phi_{\mathbf{v}}}(\mathbf{w})=|\mathbf{a}(\cdot,\mathbf{w})-\mathbf{a}(\cdot,\mathbf{v})|\leq |\mathbf{a}(\cdot,\mathbf{w})|+|\mathbf{a}(\cdot,\mathbf{v})|=\mathcal{N}_{\Phi}(|\mathbf{w}|)+\mathcal{N}_{\Phi}(|\mathbf{v}|)\in L^{p'(x)}(\Omega),
	\end{equation*}
	\noindent we deduce from Proposition \ref{propu2u1} in the Appendix that $\mathcal{N}_{\Phi_{\mathbf{v}}}(\mathbf{w})\in L^{p'(x)}(\Omega)$. So, using Theorem 4.2 from \cite{kovacik}, we get that $\mathcal{N}_{\Phi_\mathbf{v}}$ is norm-to-norm continuous. Now, for $(\mathbf{v}_n)_{n\geq 1}\subset L^{p(x)}(\Omega)^N$ we know that $\Vert \mathbf{v}_n-\mathbf{v}\Vert_{L^{p(x)}(\Omega)^N}=\Vert |\mathbf{v}_n-\mathbf{v}|\Vert_{L^{p(x)}(\Omega)}\stackrel{n\to\infty}{\longrightarrow 0}$. Thus,

	\begin{equation*}\big \Vert \mathbf{a}(\cdot,\mathbf{v}_n)-\mathbf{a}(\cdot,\mathbf{v})\big \Vert_{L^{p'(x)}(\Omega)^N}=\big \Vert |\mathbf{a}(\cdot,\mathbf{v}_n)-\mathbf{a}(\cdot,\mathbf{v})|\big \Vert_{L^{p'(x)}(\Omega)}=\Vert\mathcal{N}_{\Phi_\mathbf{v}}(\mathbf{v}_n)-\mathcal{N}_{\Phi_\mathbf{v}}(\mathbf{v})\Vert_{L^{p'(x)}(\Omega)}\longrightarrow 0.
	\end{equation*}
	
	\noindent For another proof see \cite[Theorem 3.4.4, page 407]{gasinski2005nonlinear}.
	
\end{proof}

\begin{proposition}\label{propotheta} For any $r\in [1,\infty]$ and $t_1<t_2$ two real numbers, if $\Theta\in L^r\bigl(t_1,t_2;W^{1,p(x)}(\Omega)\bigr)$ then $\nabla\Theta\in L^r\bigl(t_1,t_2;L^{p(x)}(\Omega)^N\bigr)$ and moreover:
	
	\begin{equation}\label{import2}
		\Vert \nabla\Theta\Vert_{L^r(t_1,t_2;L^{p(x)}(\Omega)^N)}\leq \Vert\Theta\Vert_{L^r(t_1,t_2;W^{1,p(x)}(\Omega))}.
	\end{equation}
\end{proposition}

\begin{proof} 
	
	\noindent Consider the spatial-gradient operator $G:W^{1,p(x)}(\Omega)\to L^{p(x)}(\Omega)^N$ given by $G(w):=\nabla w$. It is easy to check that it is a linear operator and moreover:
	
	\begin{equation}
		\Vert \nabla w\Vert_{L^{p(x)}(\Omega)^N}=\Vert |\nabla w|\Vert_{L^{p(x)}(\Omega)}\leq |w|_{L^{p(x)}(\Omega)}+\Vert |\nabla w|\Vert_{L^{p(x)}(\Omega)}=\Vert w\Vert_{W^{1,p(x)}(\Omega)}.
	\end{equation} 
	
	\noindent This shows that $G$ is a bounded linear operator and $\Vert G\Vert_{\mathcal{L}(W^{1,p(x)}(\Omega),L^{p(x)}(\Omega)^N)}\leq 1$. Now using Proposition \ref{inducedbochner} we get that $G$ induces a bounded linear operator $\overline{G}:L^r\bigl(t_1,t_2;W^{1,p(x)}(\Omega)\bigr)\to L^r\bigl(t_1,t_2;L^{p(x)}(\Omega)^N\bigr)$ via the relation $\overline{G}(w)(t):=G(w(t))$ and moreover for each $w\in L^r\bigl(t_1,t_2;W^{1,p(x)}(\Omega)\bigr)$ one has:
	
	\begin{equation}
		\Vert \overline{G}(w)\Vert_{L^r(t_1,t_2;L^{p(x)}(\Omega)^N)}\leq \Vert G\Vert_{\mathcal{L}(W^{1,p(x)}(\Omega),L^{p(x)}(\Omega)^N)}\Vert w\Vert_{L^r(t_1,t_2;W^{1,p(x)}(\Omega))}\leq \Vert w\Vert_{L^r(t_1,t_2;W^{1,p(x)}(\Omega))}.
	\end{equation}
	
	\noindent By setting $w:=\Theta\in L^r\bigl(t_1,t_2;W^{1,p(x)}(\Omega)\bigr)$ we conclude that $\nabla\Theta=\overline{G}(\Theta)\in L^r\bigl(t_1,t_2;L^{p(x)}(\Omega)^N\bigr)$ and $\Vert \nabla\Theta\Vert_{L^r(t_1,t_2;L^{p(x)}(\Omega)^N)}\leq \Vert \Theta\Vert_{L^r(t_1,t_2;W^{1,p(x)}(\Omega))}$.
	
\end{proof}

\begin{proposition}\label{propolplprim} Let $r\in [1,\infty]$, $t_1<t_2$ two real numbers and two functions $\mathbf{h}_1\in L^{\infty}\bigl(t_1,t_2;L^{p'(x)}(\Omega)^N\bigr)$ and $\mathbf{h}_2\in L^r\bigl(t_1,t_2;L^{p(x)}(\Omega)^N\bigr)$. Define the function $h:(t_1,t_2)\times\Omega\to\mathbb{R},\ h(t,x)= \mathbf{h}_1(t,x)\cdot\mathbf{h}_2(t,x)$. Then $h\in L^r\bigl(t_1,t_2;L^1(\Omega)\bigr)$, and moreover:

\begin{equation}\label{import1}
	\Vert h\Vert_{L^r(t_1,t_2;L^1(\Omega))}\leq 2\Vert \mathbf{h}_1\Vert_{L^{\infty}(t_1,t_2;L^{p'(x)}(\Omega)^N)}\Vert\mathbf{h}_2\Vert_{L^r(t_1,t_2;L^{p(x)}(\Omega)^N)}.
\end{equation}
	
\end{proposition}

\begin{proof} We divide the proof into several steps:
	
\bigskip
	
\noindent\textbf{Step I:} Define $\mathcal{G}:L^{p'(x)}(\Omega)^N\times L^{p(x)}(\Omega)^N\to L^1(\Omega)$ by $\mathcal{G}(\mathbf{u},\mathbf{v})=\mathbf{u}\cdot\mathbf{v}$ (the scalar product of the two vectorial functions). Then $\mathcal{G}$ is a continuous bilinear form and:

\begin{equation}\label{neimportant9}
	\Vert \mathcal{G}(\mathbf{u},\mathbf{v})\Vert_{L^1(\Omega)}\leq 2\Vert\mathbf{u}\Vert_{L^{p'(x)}(\Omega)^N}\Vert\mathbf{v}\Vert_{L^{p(x)}(\Omega)^N},\ \forall\ \mathbf{u}\in L^{p'(x)}(\Omega)^N\ \text{and}\ \mathbf{v}\in L^{p(x)}(\Omega)^N.
\end{equation}

\bigskip

\noindent$\bullet$ We first show that $\mathcal{G}$ is well-defined. Fix any $\mathbf{u}\in L^{p'(x)}(\Omega)^N$ and $\mathbf{v}\in L^{p(x)}(\Omega)^N$. Notice that $\mathbf{u}\cdot\mathbf{v}$ is a measurable function on $\Omega$ being a finite sum of products of measurable functions. Moreover:

\begin{align*}
\Vert\mathcal{G}(\mathbf{u},\mathbf{v})\Vert_{L^1(\Omega)}	=\Vert\mathbf{u}\cdot\mathbf{v}\Vert_{L^1(\Omega)}&=\int_{\Omega} |\mathbf{u}(x)\cdot\mathbf{v}(x)|\ dx\\
\text{(Cauchy ineq. for scalar product)}\ \ \ 	&\leq \int_{\Omega} |\mathbf{u}(x)|\cdot|\mathbf{v}(x)|\ dx\\
\text{(H\"{o}lder ineq.)}\ \ \ &\leq 2\Vert|\mathbf{u}|\Vert_{L^{p'(x)}(\Omega)}\cdot \Vert|\mathbf{v}|\Vert_{L^{p(x)}(\Omega)}\\
&=2\Vert\mathbf{u}\Vert_{L^{p'(x)}(\Omega)^N}\cdot \Vert\mathbf{v}\Vert_{L^{p(x)}(\Omega)^N}<\infty.
\end{align*}

\noindent This shows that $\mathcal{G}(\mathbf{u},\mathbf{v})=\mathbf{u}\cdot\mathbf{v}\in L^1(\Omega)$. It is straightforward to check that $\mathcal{G}$ is a biliniar form. Also from \eqref{neimportant9} it follows that $\mathcal{G}$ is a continuous bilinear form.\footnote{For this classical result about bilinear form one can follow \cite[Proposition 1, Section 10, page 56]{zorich2016mathematical}.}

\bigskip

\noindent\textbf{Step II:} The mapping $(t_1,t_2)\ni t\mapsto \mathcal{G}(\mathbf{h}_1(t,\cdot),\mathbf{h}_2(t,\cdot))=\mathbf{h}_1(t,\cdot)\cdot\mathbf{h}_2(t,\cdot)=h(t,\cdot)\in L^1(\Omega)$ is strongly measurable.

\bigskip 

\noindent$\bullet$ Indeed, from  $\mathbf{h}_1\in L^{\infty}\bigl(t_1,t_2;L^{p'(x)}(\Omega)^N\bigr)$ and $\mathbf{h}_2\in L^{r}\bigl(t_1,t_2;L^{p(x)}(\Omega)^N\bigr)$ we get that $\mathbf{h}_1:(t_1,t_2)\to L^{p'(x)}(\Omega)^N$ and $\mathbf{h}_2:(t_1,t_2)\to L^{p(x)}(\Omega)^N$ are both strongly measurable. Therefore, using Lemma 4.49 from \cite{Alip} we get that the mapping $H:(t_1,t_2)\to L^{p'(x)}(\Omega)^N\times L^{p(x)}(\Omega)^N$ given by $H(t)=\bigl(\mathbf{h}_1(t,\cdot),\mathbf{h}_2(t,\cdot)\bigr)$ is strongly measurable. 

\noindent Now since $\mathcal{G}:L^{p'(x)}(\Omega)^N\times L^{p(x)}(\Omega)^N\to L^1(\Omega)$ is continuous and $H:(t_1,t_2)\to L^{p'(x)}(\Omega)^N\times L^{p(x)}(\Omega)^N$ is strongly measurable, we get that their composition $\mathcal{G}\circ H: (t_1,t_2)\to L^1(\Omega),\ \bigl(\mathcal{G}\circ H\bigr)(t)=\mathbf{h}_1(t,\cdot)\cdot\mathbf{h}_2(t,\cdot)=h(t,\cdot)$ is also strongly measurable.\footnote{See \cite[Corollary 1.1.11, page 7]{hytonen2016analysis}.}

\bigskip

\noindent\textbf{Step III:} The mapping $(t_1,t_2)\ni t\mapsto h(t,\cdot)=\mathbf{h}_1(t,\cdot)\cdot\mathbf{h}_2(t,\cdot)\in L^1(\Omega)$ is from $L^r(t_1,t_2;L^1(\Omega))$.

\bigskip

\noindent$\bullet$ From the fact that $\mathbf{h}_1\in L^{\infty}\bigl(t_1,t_2;L^{p'(x)}(\Omega)^N\bigr)$, we get that $\Vert \mathbf{h}_1(t,\cdot)\Vert_{L^{p'(x)}(\Omega)^N}\leq \Vert \mathbf{h}_1\Vert_{L^{\infty}(t_1,t_2;L^{p'(x)}(\Omega)^N)}$, for a.e. $t\in (t_1,t_2)$. For $r\in [1,\infty)$ we may write that:

\begin{align*}
	\Vert h\Vert_{L^r(t_1,t_2;L^1(\Omega))} &=\left (\int_{t_1}^{t_2} \Vert h(t,\cdot)\Vert_{L^1(\Omega)}^r\ dt\right )^{\frac{1}{r}}=\left (\int_{t_1}^{t_2} \Vert \mathbf{h}_1(t,\cdot)\cdot\mathbf{h}_2(t,\cdot)\Vert_{L^1(\Omega)}^r\ dt\right )^{\frac{1}{r}}\\
\text{(H\"{o}lder ineq.)}\ \ \ &\leq \left (\int_{t_1}^{t_2} 2^r\Vert \mathbf{h}_1(t,\cdot)\Vert_{L^{p'(x)}(\Omega)^N}^r\cdot\Vert\mathbf{h}_2(t,\cdot)\Vert_{L^{p(x)}(\Omega)^N}^r\ dt\right )^{\frac{1}{r}}\\
&\leq \left (\int_{t_1}^{t_2} 2^r\Vert \mathbf{h}_1\Vert_{L^{\infty}(t_1,t_2;L^{p'(x)}(\Omega)^N)}^r\cdot\Vert\mathbf{h}_2(t,\cdot)\Vert_{L^{p(x)}(\Omega)^N}^r\ dt\right )^{\frac{1}{r}}\\
&=2\Vert \mathbf{h}_1\Vert_{L^{\infty}(t_1,t_2;L^{p'(x)}(\Omega)^N)}\left (\int_{t_1}^{t_2}\Vert\mathbf{h}_2(t,\cdot)\Vert_{L^{p(x)}(\Omega)^N}^r\ dt\right )^{\frac{1}{r}}\\
&=2\Vert \mathbf{h}_1\Vert_{L^{\infty}(t_1,t_2;L^{p'(x)}(\Omega)^N)}\cdot \Vert \mathbf{h}_2\Vert_{L^{r}(t_1,t_2;L^{p(x)}(\Omega)^N)}.
\end{align*}

\noindent This proves that $h\in L^r\bigl(t_1,t_2;L^1(\Omega)\bigr)$ and also the required estimate.

\noindent If $r=\infty$ the proof goes in a similar fashion:

\begin{align*}
	\Vert h\Vert_{L^\infty(t_1,t_2;L^1(\Omega))} &= \underset{t\in (t_1,t_2)}{\operatorname{ess\ sup}}\ \Vert h(t,\cdot)\Vert_{L^1(\Omega)}=\underset{t\in (t_1,t_2)}{\operatorname{ess\ sup}}\  \Vert \mathbf{h}_1(t,\cdot)\cdot\mathbf{h}_2(t,\cdot)\Vert_{L^1(\Omega)}\\ 
	\text{(H\"{o}lder ineq.)}\ \ \ &\leq 2\cdot\underset{t\in (t_1,t_2)}{\operatorname{ess\ sup}}\ \Vert \mathbf{h}_1(t,\cdot)\Vert_{L^{p'(x)}(\Omega)^N}\cdot\Vert\mathbf{h}_2(t,\cdot)\Vert_{L^{p(x)}(\Omega)^N}\\
	&\leq 2\cdot\underset{t\in (t_1,t_2)}{\operatorname{ess\ sup}}\ \Vert \mathbf{h}_1(t,\cdot)\Vert_{L^{p'(x)}(\Omega)^N}\cdot\underset{t\in (t_1,t_2)}{\operatorname{ess\ sup}}\ \Vert \mathbf{h}_2(t,\cdot)\Vert_{L^{p(x)}(\Omega)^N}\\
	&=2\Vert \mathbf{h}_1\Vert_{L^{\infty}(t_1,t_2;L^{p'(x)}(\Omega)^N)}\cdot \Vert \mathbf{h}_2\Vert_{L^{\infty}(t_1,t_2;L^{p(x)}(\Omega)^N)}.
\end{align*}

\noindent This proves that $h\in L^\infty\bigl(t_1,t_2;L^1(\Omega)\bigr)$ and also the required estimate. The proof is now complete.

\end{proof}

\begin{proposition}\label{propoa1} The following assertions hold:
	
\begin{enumerate}
	\item[\textnormal{\textbf{(1)}}] If $v\in L^\infty\bigl(0,T;W^{1,p(x)}(\Omega)\bigr)$ then $\mathbf{a}(\cdot,\nabla v(\cdot,\cdot))\in L^\infty\bigl(0,T;L^{p'(x)}(\Omega)^N\bigr)$ and the following inequality takes place:
	
	\begin{equation}\label{impo2}
		\Vert \mathbf{a}(\cdot,\nabla v(\cdot,\cdot))\Vert_{L^{\infty}(0,T;L^{p^\prime(x)}(\Omega)^N)}\leq\Vert a\Vert_{L^{p'(x)}(\Omega)}+b\max\bigl\{\Vert v\Vert_{L^{\infty}(0,T;W^{1,p(x)}(\Omega))}^{p^+-1},\Vert v\Vert_{L^{\infty}(0,T;W^{1,p(x)}(\Omega))}^{p^--1}\bigr\}.
	\end{equation}
	
	\noindent In particular, for any $r\in [1,\infty)$, we will get that $\mathbf{a}(\cdot,\nabla v(\cdot,\cdot))\in L^r\bigl(0,T;L^{p'(x)}(\Omega)^N\bigr)$ and:
	
	\begin{equation}\label{impo3}
		\Vert \mathbf{a}(\cdot,\nabla v(\cdot,\cdot))\Vert_{L^{r}(0,T;L^{p^\prime(x)}(\Omega)^N)}\leq T^{\frac{1}{r}}\Vert a\Vert_{L^{p'(x)}(\Omega)}+T^{\frac{1}{r}}b\max\bigl\{\Vert v\Vert_{L^{\infty}(0,T;W^{1,p(x)}(\Omega))}^{p^+-1},\Vert v\Vert_{L^{\infty}(0,T;W^{1,p(x)}(\Omega))}^{p^--1}\bigr\}.
	\end{equation}
	
	\item[\textnormal{\textbf{(2)}}] If $(v_n)_{n\geq 1}\subset L^\infty\bigl(0,T;W^{1,p(x)}(\Omega)\bigr)$ is a bounded sequence then the sequence $\bigl(\mathbf{a}(\cdot,\nabla v_n(\cdot,\cdot)) \bigr)_{n\geq 1}$ is bounded in $L^\infty\bigl(0,T;L^{p'(x)}(\Omega)^N\bigr)$. Moreover, if $\displaystyle\sup_{n\geq 1}\Vert v_n\Vert_{L^{\infty}(0,T;W^{1,p(x)}(\Omega))}=M$, then:
	
	\begin{equation} 
		\Vert \mathbf{a}(\cdot,\nabla v_n(\cdot,\cdot))\Vert_{L^{\infty}(0,T;L^{p^\prime(x)}(\Omega)^N)}\leq \Vert a\Vert_{L^{p'(x)}(\Omega)}+b\max\{M^{p^+-1},M^{p^--1}\},\ \forall\ n\geq 1.
	\end{equation}
\end{enumerate}

\end{proposition}

\begin{proof} \noindent\textbf{(1)} Since $v\in L^{\infty}\bigl(0,T;W^{1,p(x)}(\Omega)\bigr)$, we get from Proposition \ref{propotheta} -- applied for $r=\infty$ -- that $\nabla v\in L^{\infty}\bigl(0,T;L^{p(x)}(\Omega)^N\bigr)$ and:
	
\begin{equation}\label{neimportant8}
	\Vert \nabla v\Vert_{L^{\infty}(0,T;L^{p(x)}(\Omega)^N)}\leq \Vert v\Vert_{L^{\infty}(0,T;W^{1,p(x)}(\Omega))}.
\end{equation}

\noindent Also, from Proposition \ref{propoa0} \textbf{(4),(5)} we know that the (nonlinear) operator $\mathcal{N}_{\mathbf{a}}:L^{p(x)}(\Omega)^N\to L^{p'(x)}(\Omega)^N$ is continuous and bounded (in the nonlinear sense).

\noindent Therefore, we can apply Proposition \ref{nonlinearlinfty} \textbf{(1)} for $\nabla v\in L^{\infty}\bigl(0,T;L^{p(x)}(\Omega)^N\bigr)$, and deduce that $\mathbf{a}(\cdot,\nabla v(\cdot,\cdot))=\mathcal{N}_{\mathbf{a}}\circ\nabla v\in L^{\infty}\bigl(0,T;L^{p'(x)}(\Omega)^N\bigr)$, as needed.

\noindent Furthermore, since from \eqref{impo1} we know $\forall\ \mathbf{v}\in L^{p(x)}(\Omega)^N$ that

\begin{equation}
	\Vert \mathcal{N}_{\mathbf{a}}(\mathbf{v}(\cdot))\Vert_{L^{p'(x)}(\Omega)^N}\leq \Vert a\Vert_{L^{p'(x)}(\Omega)}+b\max\bigl\{\Vert \mathbf{v}\Vert_{L^{p(x)}(\Omega)^N}^{p^+-1},\Vert \mathbf{v}\Vert_{L^{p(x)}(\Omega)^N}^{p^--1}\bigr\}=h\left(\Vert\mathbf{v}\Vert_{L^{p(x)}(\Omega)^N}\right),
\end{equation}

\noindent where the scalar function $h:[0,\infty)\to [0,\infty),\  h(s)=\Vert a\Vert_{L^{p'(x)}(\Omega)}+b\max\bigl\{s^{p^+-1},s^{p^--1}\bigr\}$ is strictly increasing, we may apply Proposition \ref{nonlinearlinfty} \textbf{(2)} for $\nabla v\in L^{\infty}\bigl(0,T;L^{p(x)}(\Omega)^N\bigr)$ to deduce that:

\begin{align*}
	\Vert \mathbf{a}(\cdot,\nabla v(\cdot,\cdot))\Vert_{L^{\infty}(0,T;L^{p^\prime(x)}(\Omega)^N)}&=\Vert \mathcal{N}_{\mathbf{a}}\circ\nabla v\Vert_{L^{\infty}(0,T;L^{p'(x)}(\Omega)^N)}\leq h\left(\Vert\nabla v\Vert_{L^{\infty}(0,T;L^{p(x)}(\Omega)^N)}\right)\\
(h\ \text{increasing and \eqref{neimportant8}})\ \ \	&\leq h\left (\Vert v\Vert_{L^{\infty}(0,T;W^{1,p(x)}(\Omega))} \right )\\
&=\Vert a\Vert_{L^{p'(x)}(\Omega)}+b\max\bigl\{\Vert v\Vert_{L^{\infty}(0,T;W^{1,p(x)}(\Omega))}^{p^+-1},\Vert v\Vert_{L^{\infty}(0,T;W^{1,p(x)}(\Omega))}^{p^--1}\bigr\}.
\end{align*}

\noindent Taking into account that $L^{\infty}(0,T;L^{p^\prime(x)}(\Omega)^N)\hookrightarrow L^{r}(0,T;L^{p^\prime(x)}(\Omega)^N)$ we get that $\mathbf{a}(\cdot,\nabla v(\cdot,\cdot))\in L^{r}(0,T;L^{p^\prime(x)}(\Omega)^N)$ and:

\begin{align*}
	\Vert \mathbf{a}(\cdot,\nabla v(\cdot,\cdot))\Vert_{L^{r}(0,T;L^{p^\prime(x)}(\Omega)^N)}&=\left (\int_{0}^T \Vert \mathbf{a}(\cdot,\nabla v(t,\cdot))\Vert_{L^{p'(x)}(\Omega)^N}^{r} \right )^{\frac{1}{r}}\\
	&\leq T^{\frac{1}{r}}\underset{t\in (0,T)}{\operatorname{ess\ sup}}\ \Vert \mathbf{a}(\cdot,\nabla v(t,\cdot))\Vert_{L^{p'(x)}(\Omega)^N}\\
	&=T^{\frac{1}{r}}\Vert \mathbf{a}(\cdot,\nabla v(\cdot,\cdot))\Vert_{L^{\infty}(0,T;L^{p^\prime(x)}(\Omega)^N)}.
\end{align*}

\noindent The conclusion follows.
\bigskip

\noindent\textbf{(2)} Applying \eqref{impo2} for each $v_n\in L^\infty\bigl(0,T;W^{1,p(x)}(\Omega)\bigr),\ n\geq 1$ will give us that:

\begin{align*}
		\Vert \mathbf{a}(\cdot,\nabla v_n(\cdot,\cdot))\Vert_{L^{\infty}(0,T;L^{p^\prime(x)}(\Omega)^N)}&\leq\Vert a\Vert_{L^{p'(x)}(\Omega)}+b\max\bigl\{\Vert v_n\Vert_{L^{\infty}(0,T;W^{1,p(x)}(\Omega))}^{p^+-1},\Vert v_n\Vert_{L^{\infty}(0,T;W^{1,p(x)}(\Omega))}^{p^--1}\bigr\}\\
		&=h\left (\Vert v_n\Vert_{L^{\infty}(0,T;W^{1,p(x)}(\Omega))} \right )\\
(h\ \text{increasing})\ \ \  &\leq h(M)=\Vert a\Vert_{L^{p'(x)}(\Omega)}+b\max\{M^{p^+-1},M^{p^--1}\},\ \forall\ n\geq 1.
\end{align*}

\noindent The proof is now complete.
	
\end{proof}

\begin{proposition}\label{propoa3} Assume we have a sequence $(\mathbf{w}_n)_{n\geq 1}\subset L^{\infty}\bigl(0,T;L^{p(x)}(\Omega)^N\bigr)$ and a function $\mathbf{w}\in L^{\infty}\bigl(0,T;L^{p(x)}(\Omega)^N\bigr)$ such that $\mathbf{w}_n\to \mathbf{w}$ in $L^{\infty}\bigl(0,T;L^{p(x)}(\Omega)^N\bigr)$. Then for every $r\in [1,\infty)$ one has that
	
	\begin{equation}
		\mathbf{a}(\cdot,\mathbf{w}_n)\to\mathbf{a}(\cdot,\mathbf{w})\ \text{in}\ L^r\bigl(0,T;L^{p'(x)}(\Omega)^N\bigr).
	\end{equation}
	
\end{proposition}

\begin{proof} From Proposition \ref{propoa0} \textbf{(4),(5)} we have that $\mathcal{N}_{\mathbf{a}}:L^{p(x)}(\Omega)^N\to L^{p'(x)}(\Omega)^N$ is a continuous and bounded (nonlinear) operator. Since $\mathbf{w}_n\to \mathbf{w}$ in $L^{\infty}\bigl(0,T;L^{p(x)}(\Omega)^N\bigr)$ we can apply Proposition \ref{propocompositionprinc} to conclude that $\mathbf{a}(\cdot,\mathbf{w}_n)=\mathcal{N}_{\mathbf{a}}\circ\mathbf{w}_n\to \mathcal{N}_{\mathbf{a}}\circ \mathbf{w}=\mathbf{a}(\cdot,\mathbf{w})$ in $L^{r}\bigl(0,T;L^{p'(x)}(\Omega)^N\bigr)$, for any $r\in [1,\infty)$, as required.

\end{proof}

\begin{proposition}\label{propoa2} Let $v\in L^\infty\bigl(0,T;W^{1,p(x)}(\Omega)\bigr)$ and fix $0\leq t_1<t_2\leq T$. Define now $\ell_v:L^2\bigl(0,T;W^{1,p(x)}(\Omega)\bigr)\to\mathbb{R}$ by:
	
	\begin{equation}
		\ell_v(w):=\int_{t_1}^{t_2}\int_{\Omega} \mathbf{a}(x,\nabla v(t,x))\cdot\nabla w(t,x)\ dx\ dt.
	\end{equation}
	
\noindent Therefore $\ell_v$ is a well-defined bounded linear functional, i.e. $\ell_v\in L^2\bigl(0,T;W^{1,p(x)}(\Omega)\bigr)^*$. Moreover the following inequality holds for every $w\in L^2\bigl(0,T;W^{1,p(x)}(\Omega)\bigr)$:

\begin{equation}\label{import3}
	|\ell_v(w)|\leq \sqrt{t_2-t_1}\left (\Vert a\Vert_{L^{p'(x)}(\Omega)}+b\max\bigl\{\Vert v\Vert_{L^{\infty}(0,T;W^{1,p(x)}(\Omega))}^{p^+-1},\Vert v\Vert_{L^{\infty}(0,T;W^{1,p(x)}(\Omega))}^{p^--1}\bigr\}\right )\Vert w\Vert_{L^2(0,T;W^{1,p(x)}(\Omega))}.
\end{equation}

\end{proposition}

\begin{proof} Since $v\in L^\infty\bigl(0,T;W^{1,p(x)}(\Omega)\bigr)$ we deduce from Proposition \ref{propoa1} \textbf{(1)} that $\mathbf{a}(\cdot,\nabla v(\cdot,\cdot))\in L^{\infty}\bigl(0,T;L^{p'(x)}(\Omega)^N\bigr)$ and
	
\begin{equation}
	\Vert \mathbf{a}(\cdot,\nabla v(\cdot,\cdot))\Vert_{L^{\infty}(0,T;L^{p^\prime(x)}(\Omega)^N)}\leq\Vert a\Vert_{L^{p'(x)}(\Omega)}+b\max\bigl\{\Vert v\Vert_{L^{\infty}(0,T;W^{1,p(x)}(\Omega))}^{p^+-1},\Vert v\Vert_{L^{\infty}(0,T;W^{1,p(x)}(\Omega))}^{p^--1}\bigr\}:=C.
\end{equation}

\noindent It is obvious to see that $\Vert \mathbf{a}(\cdot,\nabla v(\cdot,\cdot))\Vert_{L^{\infty}(t_1,t_2;L^{p^\prime(x)}(\Omega)^N)}\leq \Vert \mathbf{a}(\cdot,\nabla v(\cdot,\cdot))\Vert_{L^{\infty}(0,T;L^{p^\prime(x)}(\Omega)^N)}\leq C$.

\noindent Also, for any $w\in L^2\bigl(0,T;W^{1,p(x)}(\Omega)\bigr)$ we have from Proposition \ref{propotheta} that $\nabla w\in L^2\bigl(0,T;L^{p(x)}(\Omega)^N\bigr)$. Therefore, from Proposition \ref{propolplprim} we deduce that the function $(0,T)\times\Omega\ni (t,x)\mapsto h(t,x):=\mathbf{a}(x,\nabla v(t,x))\cdot\nabla w(t,x)\in\mathbb{R}$ is in $L^2\bigl(0,T;L^1(\Omega)\bigr)\subset L^2\bigl(t_1,t_2;L^1(\Omega)\bigr)\subset L^1\bigl(t_1,t_2;L^1(\Omega)\bigr)$. This shows that $\ell_v$ is well-defined and moreover

\begin{align*}
	|\ell_v(w)|&=\Vert h\Vert_{L^1(t_1,t_2;L^1(\Omega))}\stackrel{\text{Cauchy}}{\leq}\sqrt{t_2-t_1}\Vert h\Vert_{L^2(t_1,t_2;L^1(\Omega))} \\
	&\stackrel{\eqref{import1}}{\leq} \sqrt{t_2-t_1}\Vert \mathbf{a}(\cdot,\nabla v(\cdot,\cdot))\Vert_{L^{\infty}(t_1,t_2;L^{p^\prime(x)}(\Omega)^N)}\cdot \Vert\nabla w\Vert_{L^2(t_1,t_2;L^{p(x)}(\Omega)^N)}\\
	&\stackrel{\eqref{import2}}{\leq}\sqrt{t_2-t_1}\Vert \mathbf{a}(\cdot,\nabla v(\cdot,\cdot))\Vert_{L^{\infty}(t_1,t_2;L^{p^\prime(x)}(\Omega)^N)}\cdot \Vert w\Vert_{L^2(t_1,t_2;W^{1,p(x)}(\Omega))}\\
	&\leq C\sqrt{t_2-t_1}\Vert w\Vert_{L^2(t_1,t_2;W^{1,p(x)}(\Omega))}\\
	&\leq C\sqrt{t_2-t_1}\Vert w\Vert_{L^2(0,T;W^{1,p(x)}(\Omega))}.
\end{align*}

\noindent It is straightforward to see that $\ell_v$ is a linear operator, so it is a bounded linear operator from the previous relation. Therefore $\ell_v\in L^2\bigl(0,T;W^{1,p(x)}(\Omega)\bigr)^*$ and the required estimate holds.

\end{proof}

\subsection{Properties of $\mathcal{A}$}

\begin{proposition}\label{propomathcalA1} The following assertions about $\mathcal{A}:W^{1,p(x)}(\Omega)\to [0,\infty)$ are valid:
	
\begin{enumerate}
	\item[\textnormal{\textbf{(1)}}] $\mathcal{A}$ is a convex functional. In particular, since $\mathcal{A}(0)=0$, we have that
	
	\begin{equation}\label{Asubhomogeneity}
		\mathcal{A}(\theta V)\leq \theta\mathcal{A}(V),\ \forall\ V\in W^{1,p(x)}(\Omega),\ \forall\ \theta\in [0,1].
	\end{equation}
	
	\item[\textnormal{\textbf{(2)}}] $\mathcal{A}\in C^1\bigl(W^{1,p(x)}(\Omega)\bigr)$ and:
	
	\begin{equation}\label{laKAUF16}
		\langle\mathcal{A}'(u),\phi\rangle=\int_{\Omega}\mathbf{a}(x,\nabla u(x))\cdot\nabla\phi(x)\ dx,\ \forall\ u,\phi\in W^{1,p(x)}(\Omega).
	\end{equation}
	
	\noindent In particular $\mathcal{A}':W^{1,p(x)}(\Omega)\to W^{1,p(x)}(\Omega)^*$ is a continuous (nonlinear) operator.
	
	\item[\textnormal{\textbf{(3)}}] $\mathcal{A}$ is weakly lower semicontinuous, meaning that:
	
	\begin{equation}
		\text{If}\ u_n\weak u\ \text{in}\ W^{1,p(x)}(\Omega)\ \text{then}\ \liminf\limits_{n\to\infty}\mathcal{A}(u_n)\geq \mathcal{A}(u).
	\end{equation}
	
	\item[\textnormal{\textbf{(4)}}] For any $u_1,u_2\in W^{1,p(x)}(\Omega)$ the following inequality holds:
	
	\begin{equation}
		\int_{\Omega}\mathbf{a}(x,\nabla u_2(x))\cdot \bigl(\nabla u_2(x)-\nabla u_1(x)\bigr)\ dx.\geq\mathcal{A}(u_2)-\mathcal{A}(u_1)\geq \int_{\Omega}\mathbf{a}(x,\nabla u_1(x))\cdot \bigl(\nabla u_2(x)-\nabla u_1(x)\bigr)\ dx.
	\end{equation}
	
	\item[\textnormal{\textbf{(5)}}] If we define the \textbf{coercivity profile} $\gamma:[0,\infty)\to [0,\infty)$ by
	
	\begin{equation}
		\gamma(s):=\inf\bigl\{\mathcal{A}(V)\ |\ V\in W^{1,p(x)}(\Omega)\ \text{and}\ \rho_{p(x)}\bigl(|\nabla V|\bigr)\geq s\bigr\},
	\end{equation}
	
	\noindent then $\gamma$ is a well-defined, nondecreasing function with $\lim\limits_{s\to\infty} \gamma(s)=+\infty$. Moreover the following inequality holds:\footnote{In fact the converse is also true: if there is a function $\gamma:[0,\infty)\to [0,\infty)$ that is nondecreasing, has $\lim\limits_{s\to \infty} \gamma(s)=+\infty$ and $\mathcal{A}(V)\geq\gamma\left(\rho_{p(x)}\bigl(|\nabla V|\bigr)\right),\ \forall\ V\in W^{1,p(x)}(\Omega)$, then \textbf{(H7)} holds. So this is a characterization of the hypothesis \textbf{(H7)}.}
	
	\begin{equation}
		\mathcal{A}(V)\geq\gamma\left(\rho_{p(x)}\bigl(|\nabla V|\bigr)\right),\ \forall\ V\in W^{1,p(x)}(\Omega).
	\end{equation}
	
	\item[\textnormal{\textbf{(6)}}] The following inequality holds:
	
	\begin{equation}
		\Vert\nabla V\Vert_{L^{p(x)}(\Omega)^N}:=\bigl\Vert|\nabla V|\bigr\Vert_{L^{p(x)}(\Omega)}\leq C_*\bigl(1+\mathcal{A}(V)\bigr), \ \forall\ V\in W^{1,p(x)}(\Omega),
	\end{equation}
	
	\noindent where the constant $C_*\in (0,\infty)$ is given by 
	
	\begin{equation}
		C_*=\sup\bigl\{\Vert\nabla V\Vert_{L^{p(x)}(\Omega)^N}\ |\ V\in W^{1,p(x)}(\Omega)\ \text{and}\ \mathcal{A}(V)\leq 1 \bigr\}.
	\end{equation}
	
\end{enumerate}
	
\end{proposition}

\begin{proof} \noindent\textbf{(1)} From the proof of Proposition 3.1 \textbf{(6)} given in \cite{max4} we have that $\mathcal{A}$ is well-defined. 
	
\noindent Set any $V_1,V_2\in W^{1,p(x)}(\Omega)$ and any $\theta\in (0,1)$ we have from the (strict) convexity of $A(x,\cdot)$ for a.e. $x\in\Omega$ proved in \cite[Proposition 3.1 \textbf{(4)}]{max4}, we get that 

\begin{equation}
	A(x,\theta \nabla V_1(x)+(1-\theta)\nabla V_2(x)\bigr)\leq \theta A(x,\nabla V_1(x))+(1-\theta)A(x,\nabla V_2(x))\ \text{for a.e.}\ x\in\Omega. 
\end{equation}

\noindent Integrating this inequality on $\Omega$ we get that:

\begin{align}\label{convexitynotstrict}
	\mathcal{A}(\theta V_1+(1-\theta)V_2)&=\int_{\Omega} A(x,\theta \nabla V_1(x)+(1-\theta)\nabla V_2(x)\bigr)\ dx\\
	&\leq \int_{\Omega} \theta A(x,\nabla V_1(x))+(1-\theta)A(x,\nabla V_2(x))\ dx\\
	&=\theta\int_{\Omega}  A(x,\nabla V_1(x))\ dx+(1-\theta)\int_{\Omega} A(x,\nabla V_2(x))\ dx\\
	&=\theta \mathcal{A}(V_1)+(1-\theta)\mathcal{A}(V_2),
\end{align} 

\noindent i.e. $\mathcal{A}$ is a convex functional.\footnote{However, it is easy to remark that $\mathcal{A}$ is \textbf{not} a strictly convex functional. More precisely equality in \eqref{convexitynotstrict} occurs iff $\nabla V_1=\nabla V_2$ a.e. on $\Omega$. Taking into account that $\Omega$ is a connected domain this means that $V_1-V_2\equiv$ constant.} For the last inequality just remark that for any $\theta\in [0,1]$ we may write

\begin{equation}
	\mathcal{A}(\theta V)=\mathcal{A}(\theta V+(1-\theta)0)\leq \theta\mathcal{A}(V)+(1-\theta)\mathcal{A}(0)=\theta\mathcal{A}(V),\ \forall\ V\in W^{1,p(x)}(\Omega).
\end{equation}

	\bigskip	
	\noindent\textbf{(2)} The complete proof can be found in \cite[Proposition 3.1 \textbf{(6)}]{max4}.

	\bigskip	
	\noindent\textbf{(3)} Since $\mathcal{A}$ is Fr\'{e}chet differentiable we have from \cite[pages 59-60]{coleman} that $\mathcal{A}$ is also G\^{a}teaux differentiable. Because $W^{1,p(x)}(\Omega)$ is a reflexive Banach space and $\mathcal{A}:W^{1,p(x)}(\Omega)\to [0,\infty)$ is a convex functional (from \textbf{(1)}), we deduce from \cite[Theorem 6.2.1]{kurdila2005convex} that:
	
	\begin{equation}
		\mathcal{A}(V_1)-\mathcal{A}(V_2)\geq \langle\mathcal{A}'(V_2),V_2-V_1\rangle,\ \forall\ V_1, V_2\in W^{1,p(x)}(\Omega).
	\end{equation}
	
	\noindent Therefore we may write for each $n\geq 1$ that:
	
		\begin{equation}
		\mathcal{A}(u_n)-\mathcal{A}(u)\geq \langle\mathcal{A}'(u),u_n-u\rangle\ \Longrightarrow\ \mathcal{A}(u_n)\geq \mathcal{A}(u)+\langle\mathcal{A}'(u),u_n-u\rangle.
	\end{equation}
	
	\noindent From the fact that $\mathcal{A}'(u)\in W^{1,p(x)}(\Omega)^*$ and $u_n\weak u$ in $W^{1,p(x)}(\Omega)$ we obtain that 
	
	$$\lim\limits_{n\to\infty} \mathcal{A}'(u)(u_n-u)=\lim\limits_{n\to\infty}\langle\mathcal{A}'(u),u_n-u\rangle=0.$$
	
	\noindent Finally:
	
	\begin{align*}
	\liminf\limits_{n\to\infty} \mathcal{A}(u_n)&\geq \liminf\limits_{n\to\infty}\mathcal{A}(u)+\langle\mathcal{A}'(u),u_n-u\rangle=\mathcal{A}(u)+\liminf\limits_{n\to\infty}\langle\mathcal{A}'(u),u_n-u\rangle\\
	&=\mathcal{A}(u)+\lim\limits_{n\to\infty}\langle\mathcal{A}'(u),u_n-u\rangle=\mathcal{A}(u),
	\end{align*}
	
	\noindent as wanted.
	
	\bigskip
	
	\noindent\textbf{(4)} Using Proposition 3.1 \textbf{(3)} from \cite{max4} we get for any $x\in\overline{\Omega}$ that:
	
	\begin{equation}
			\mathbf{a}(x,\nabla u_2(x))\cdot \bigl(\nabla u_2(x)-\nabla u_1(x)\bigr)\geq A(x,\nabla u_2(x))-A(x,\nabla u_1(x))\geq\mathbf{a}(x,\nabla u_1(x))\cdot \bigl(\nabla u_2(x)-\nabla u_1(x)\bigr).
	\end{equation}
	
	\noindent Integrating this relation on $\Omega$ gives the desired inequality.
	
	\bigskip	
	\noindent\textbf{(5)} $\bullet$ First we show that $\gamma$ is well-defined. In fact, we will show that for any $s\geq 0$ the set $E_s:=\bigl\{V\in W^{1,p(x)}(\Omega)\ |\ \rho_{p(x)}\bigl(|\nabla V|\bigr)\geq s\bigr\}$ is not empty.
	
	\medskip
	
	\noindent Indeed, if we fix any $\varphi_0\in C^{\infty}_{c}(\Omega)$ that is not a constant function, we get that $|\nabla\varphi_0(x)|>0$ on a subset of $\Omega$ with strictly positive measure. Therefore $\rho_{p(x)}(|\nabla\varphi_0|)=\displaystyle\int_{\Omega} |\nabla\varphi_0(x)|^{p(x)}\ dx>0$. Note that for any $\theta\geq 1$ we have for any $x\in\Omega$ that $|\nabla (\theta\varphi_0)(x)|^{p(x)}=\theta^{p(x)}|\nabla\varphi_0(x)|^{p(x)}\geq \theta^{p^-}|\nabla\varphi_0(x)|^{p(x)}$. From here we deduce that $\rho_{p(x)}\bigl(|\nabla (\theta\varphi_0)|\bigr)\geq \theta^{p^-}\rho_{p(x)}\bigl(|\nabla\varphi_0|\bigr)$ for any $\theta\geq 1$. But for any $\theta\geq\max\left\{1,\left [\dfrac{s}{\rho_{p(x)}\bigl(|\nabla\varphi_0|\bigr)}\right ]^{\frac{1}{p^-}}\right\}$ one has that $\theta\varphi_0\in C^{\infty}_c(\Omega)\subset W^{1,p(x)}(\Omega)$ and $\rho_{p(x)}\bigl(|\nabla (\theta\varphi_0)|\bigr)\geq s$, from which we deduce that $E_s\neq\emptyset$. But then $\bigl\{\mathcal{A}(V)\ |\ V\in W^{1,p(x)}(\Omega)\ \text{and}\ \rho_{p(x)}\bigl(|\nabla V|\bigr)\geq s\bigr\}\neq \emptyset$, and henceforth $0\leq \gamma(s)\leq \mathcal{A}(\theta\varphi_0)<\infty$. This shows that $\gamma:[0,\infty)\to [0,\infty)$ is well-defined.
		
	\medskip
	
	\noindent $\bullet$ Let any $0\leq s_1\leq s_2<\infty$. Remark that $E_{s_2}\subset E_{s_1}$. This happens because if $V\in E_{s_2}$ then $\rho_{p(x)}(|\nabla V|)\geq s_2\geq s_1$, which means that $V\in E_{s_1}$. Taking the infimum of the same functional over a smaller set can only increase the
	infimum. Hence $\gamma(s_1)=\displaystyle\inf_{V\in E_{s_1}}\mathcal{A}(V)\leq \displaystyle\inf_{V\in E_{s_2}}\mathcal{A}(V)=\gamma(s_2)$. 
	
	\medskip
	
	\noindent $\bullet$ Since $\gamma$ is a nondecreasing function the following limit exists: $\lim\limits_{s\to\infty} \gamma(s)=\ell\in [0,\infty]$. Suppose that $0\leq \ell<\infty$. We have that $\gamma(s)\leq \ell,\ \forall\ s\in [0,\infty)$. From the definition of the infimum, for every $n\geq 1$ we get that there is some $V_n\in E_n\subset W^{1,p(x)}(\Omega)$ such that: 
	
	\begin{equation}\label{capasuna1}
		0\leq \mathcal{A}(V_n)\leq \dfrac{1}{n}+\inf_{V\in E_n}\mathcal{A}(V)=\dfrac{1}{n}+\gamma(n)\leq\ell+\dfrac{1}{n}\leq\ell+1,\ \forall\ n\geq 1.
	\end{equation}
	
	\noindent But, from $V_n\in E_n$ for every $n\geq 1$, we get that $\rho_{p(x)}(|\nabla V_n|)=\displaystyle\int_{\Omega} |\nabla V_n(x)|^{p(x)}\ dx\geq n,\ \forall\ n\geq 1$. Therefore we have that $\lim\limits_{n\to\infty} \displaystyle\int_{\Omega} |\nabla V_n(x)|^{p(x)}\ dx=+\infty$. So, from \textbf{(H7)}, we get that $\lim\limits_{n\to\infty} \mathcal{A}(V_n)=+\infty$, which contradicts \eqref{capasuna1}. Thereby, we have shown that $\ell=\infty$, i.e. $\lim\limits_{s\to\infty} \gamma(s)=\infty$.
	
	\medskip
	
	\noindent $\bullet$ Let any $V\in W^{1,p(x)}(\Omega)$ and denote $s_V:=\rho_{p(x)}(|\nabla V|)$. Since $\rho_{p(x)}(|\nabla V|)=s_V\geq s_V$ it follows that $V\in E_{s_V}$. Thence, from the definition of the infimum we get that:
	
	\begin{equation}
	\gamma\bigl(\rho_{p(x)}(|\nabla V|) \bigr)=\gamma(s_V)=\inf_{\tilde{V}\in E_{s_V}}\mathcal{A}(\tilde{V})\leq \mathcal{A}(V).
	\end{equation}
	
	\bigskip	
	\noindent\textbf{(6)} $\bullet$ First we show that $C_*\in (0,\infty)$. Indeed, if we fix any $\varphi_0\in C^{\infty}_{c}(\Omega)$ that is not a constant function, we get using Remark \ref{remmathcalA} that $\mathcal{A}(\theta\varphi_0)>0$, for any $\theta\in (0,1]$. From \eqref{Asubhomogeneity} we get that for any $\theta\in (0,1]$ one has that:
	
	\begin{equation}
		\mathcal{A}(\theta\varphi_0)\leq \theta\mathcal{A}(\varphi_0).
	\end{equation}
	
	\noindent Therefore for any $0<\theta\leq \min\left\{1,\dfrac{1}{\mathcal{A}(\varphi_0)}\right\}$ we get that $\theta\varphi_0\in\{V\in W^{1,p(x)}(\Omega)\ |\ \mathcal{A}(V)\leq 1\}$ and $\mathcal{A}(\theta\varphi_0)>0$. So $C_*>0$.
	
	\medskip
	
	\noindent Suppose that $C_*=+\infty$. Therefore we can find a sequence $(V_n)_{n\geq 1}\subset W^{1,p(x)}(\Omega)$ with 
	
	\begin{equation}\label{septembrie1}
	\lim\limits_{n\to\infty}\Vert \nabla V_n\Vert_{L^{p(x)}(\Omega)^N}=+\infty\ \text{and}\ \mathcal{A}(V_n)\leq 1, \ \forall\ n\geq 1.
	\end{equation}

	\noindent Without loss of generality we may assume that $\Vert \nabla V_n\Vert >1$ for any $n\geq 1$. Therefore, from Lemma 3.2.4 given in \cite[page 73]{Hasto} we get that $\rho_{p(x)}(|\nabla V_n|)\geq \Vert \nabla V_n\Vert_{L^{p(x)}(\Omega)^N}$ for any $n\geq 1$. As a consequence $\rho_{p(x)}(|\nabla V_n|)\to \infty$ and therefore, using \textbf{(H7)}, we get that $\mathcal{A}(V_n)\to\infty$, contradicting \eqref{septembrie1}. Thus $C_*\in (0,\infty)$.
	
	\medskip
	
	\noindent $\bullet$ Fix any arbitrary $V\in W^{1,p(x)}(\Omega)$. We distinguish two separate cases:
	
	\begin{enumerate}
		\item[$\blacksquare$] If $\Vert \nabla V\Vert_{L^{p(x)}(\Omega)^N}\leq C_*$ then there is nothing to prove since:
		
		\begin{equation}
			\Vert \nabla V\Vert_{L^{p(x)}(\Omega)^N}\leq C_*\leq C_*\bigl(1+\mathcal{A}(V)\bigr).
		\end{equation}

		\medskip
		
		\item[$\blacksquare$] If $C_V:=\Vert \nabla V\Vert_{L^{p(x)}(\Omega)^N}>C_*$ we denote $\theta_*:=\dfrac{C_*}{C_V}\in (0,1)$. Then, for any $\theta\in (\theta_*,1)$ we have that

		\begin{equation}
			\Vert\nabla(\theta V)\Vert_{L^{p(x)}(\Omega)^N}=\theta\Vert\nabla V\Vert_{L^{p(x)}(\Omega)^N}=\theta C_V>\dfrac{C_*}{C_V}\cdot C_V=C_*.
		\end{equation}
		
		\noindent This shows, using the definition of $C_*$, that $\mathcal{A}(\theta V)>1$ for any $\theta\in (\theta_*,1)$. Using \eqref{Asubhomogeneity}, we get that:
		
		\begin{equation}
		1<\mathcal{A}(\theta V)\leq \theta\mathcal{A}(V)\ \Longrightarrow\ \theta\mathcal{A}(V)>1,\ \forall\ \theta\in (\theta_*,1).
		\end{equation}
		
		\noindent Making $\theta\to \theta_*$ we deduce that $\theta_*\mathcal{A}(V)\geq 1$, from where:
		
		\begin{equation}
		\Vert \nabla V\Vert_{L^{p(x)}(\Omega)^N}=C_V\leq C_*\mathcal{A}(V)\leq  C_*\bigl(1+\mathcal{A}(V)\bigr).
		\end{equation}
		
	\end{enumerate}
\end{proof}

\subsection{Properties of $\overline{\mathcal{A}}$}

\begin{proposition}\label{propoverA} The following assertions about $\overline{\mathcal{A}}:L^{2}(\Omega)\to [0,\infty]$ are valid:
	
	\begin{enumerate}
		
		\item[\textnormal{\textbf{(1)}}] $\overline{\mathcal{A}}$ is a convex functional. In particular, since $\overline{\mathcal{A}}(0)=0$, we have that
		
		\begin{equation}
			\overline{\mathcal{A}}(\theta V)\leq \theta\overline{\mathcal{A}}(V),\ \forall\ V\in L^{2}(\Omega),\ \forall\ \theta\in [0,1].
		\end{equation}
		
		\item[\textnormal{\textbf{(2)}}]  With the extra assumption that $p\in\mathcal{P}^{\text{log}}(\Omega)$ -- i.e. $p$ is a log-H\"{o}lder continuous exponent -- we have that $\overline{\mathcal{A}}$ is a lower semicontinuous functional, i.e. 
		
		\begin{equation}
			\text{if}\ u_n\to u\ \text{in}\ L^2(\Omega)\ \text{then}\ \liminf\limits_{n\to\infty} \overline{\mathcal{A}}(u_n)\geq \overline{\mathcal{A}}(u).
		\end{equation}
		
		\item[\textnormal{\textbf{(3)}}]  If $p\in\mathcal{P}^{\text{log}}(\Omega)$ then the subdifferential of $\overline{\mathcal{A}}$ at any $u\in L^2(\Omega)$ is given by:
	\end{enumerate}
	
	\begin{equation}\label{subdifoA}
		\partial\overline{\mathcal{A}}(u)=\begin{cases} \left\{w\in L^2(\Omega)\ \bigg|\ \mathcal{A}(v)\geq \mathcal{A}(u)+\displaystyle\int_{\Omega} w\cdot\bigl (v-u\bigr)\ dx,\ \forall\ v\in W^{1,p(x)}(\Omega)\right\},\  \text{if }u\in W^{1,p(x)}(\Omega)\\[3mm] \emptyset, \ \ \text{if}\ u\in L^2(\Omega)\setminus W^{1,p(x)}(\Omega). \end{cases}
	\end{equation}
	
	\item[\textnormal{\textbf{(4)}}] If $p\in\mathcal{P}^{\text{log}}(\Omega)$, then $\overline{\mathcal{A}}$ is weakly lower semicontinuous, i.e.:
	
	\begin{equation}
		\text{if}\ u_n\weak u\ \text{in}\ L^2(\Omega)\ \text{then}\ \liminf\limits_{n\to\infty} \overline{\mathcal{A}}(u_n)\geq \overline{\mathcal{A}}(u).
	\end{equation}
	
\end{proposition}

\begin{proof}\noindent\textbf{(1)} Let any $V_1,V_2\in L^2(\Omega)$ and any $\theta\in (0,1)$. We need to show that:
	
	\begin{equation}\label{convex1}
		\overline{\mathcal{A}}(\theta V_1+(1-\theta)V_2)\leq \theta\overline{\mathcal{A}}(V_1)+(1-\theta)\overline{\mathcal{A}}(V_2).
	\end{equation}
	
\noindent Note that if $V_1\notin W^{1,p(x)}(\Omega)$ or $V_2\notin W^{1,p(x)}(\Omega)$ then the right-hand side of \eqref{convex1} is $+\infty$ and therefore there is nothing to prove.

\noindent Consider now the case in which $V_1,V_2\in W^{1,p(x)}(\Omega)$. Therefore $\theta V_1+(1-\theta)V_2\in W^{1,p(x)}(\Omega)$ and \eqref{convex1} becomes:

\begin{equation}\label{convex2}
	\mathcal{A}(\theta V_1+(1-\theta)V_2)\leq \theta\mathcal{A}(V_1)+(1-\theta)\mathcal{A}(V_2),
\end{equation}

\noindent which is true from the convexity of $\mathcal{A}$, proved in Proposition \ref{propomathcalA1} \textbf{(1)}. Thus $\overline{\mathcal{A}}$ is a convex functional. For the last inequality just remark that for any $\theta\in [0,1]$ we may write

\begin{equation}
	\overline{\mathcal{A}}(\theta V)=\overline{\mathcal{A}}(\theta V+(1-\theta)0)\leq \theta\overline{\mathcal{A}}(V)+(1-\theta)\overline{\mathcal{A}}(0)=\theta\overline{\mathcal{A}}(V),\ \forall\ V\in L^{2}(\Omega).
\end{equation}

\bigskip	
\noindent\textbf{(2)} If $\liminf\limits_{n\to\infty} \overline{\mathcal{A}}(u_n)=+\infty$ then there is nothing to prove.

\noindent Suppose now that $d:=\liminf\limits_{n\to\infty} \overline{\mathcal{A}}(u_n)\in [0,\infty)$. Therefore we can find a subsequence $(u_{n_k})_{k\geq 1}\subset L^2(\Omega)$ such that $d=\lim\limits_{k\to\infty}\overline{\mathcal{A}}(u_{n_k})$.

\noindent Without losing the generality we can assume that $\overline{\mathcal{A}}(u_{n_k})<\infty$ for each $k\geq 1$. This means that $(u_{n_k})_{k\geq 1}\subset W^{1,p(x)}(\Omega)$ and, as a consequence, we can write that:

\begin{equation}\label{convex3}
	d=\lim_{k\to\infty}\overline{\mathcal{A}}(u_{n_k})=\lim_{k\to\infty}\mathcal{A}(u_{n_k}).
\end{equation}

\noindent If $(\nabla u_{n_k})_{k\geq 1}$ is unbounded in $L^{p(x)}(\Omega)^N$, then we can find a further subsequence $\bigl(u_{n_{k_\ell}}\bigr)_{\ell\geq 1}$ such that $\rho_{p(x)}\bigl(|\nabla u_{n_{k_{\ell}}}|\bigr)=\displaystyle\int_{\Omega} |\nabla u_{n_{k_{\ell}}}(x)|^{p(x)}\ dx\stackrel{\ell\to\infty}{\longrightarrow}+\infty$. Using hypothesis \textbf{(H7)} we obtain that $\mathcal{A}(u_{n_{k_{\ell}}})\stackrel{\ell\to\infty}{\longrightarrow}+\infty$ which is in contradiction with \eqref{convex3}.

\noindent Thus $(\nabla u_{n_k})_{k\geq 1}$ is a bounded sequence in $L^{p(x)}(\Omega)^N$. Now since $L^{p(x)}(\Omega)$ is a reflexive Banach space, so is $L^{p(x)}(\Omega)^N$. Applying \textit{Eberlein-\v{S}mulian theorem}\footnote{See \cite[Theorem 3.4.16, page 221]{papageorgiou2018applied}.} we deduce that there is some function $\mathbf{U}\in L^{p(x)}(\Omega)^N$ and a further subsequence such that $\nabla u_{n_{k_{\ell}}}\weak \mathbf{U}$ in $L^{p(x)}(\Omega)^N$. But we also know that $u_{n_{k_{\ell}}}\to u$ in $L^2(\Omega)$. 

\noindent Using now Lemma \ref{bijuteria} (here we need the extra assumption that $p$ is log-H\"{o}lder continuous) we get that $u\in W^{1,p(x)}(\Omega)$ and moreover $u_{n_{k_{\ell}}}\weak u$ in $W^{1,p(x)}(\Omega)$.

\noindent In conclusion from the weak lower semicontinuity of $\mathcal{A}:W^{1,p(x)}(\Omega)\to [0,\infty)$ proved at Proposition \ref{propomathcalA1} \textbf{(3)} we conclude that:

\begin{equation}
	\liminf\limits_{n\to\infty} \overline{\mathcal{A}}(u_n)=d=\lim\limits_{\ell\to\infty} \mathcal{A}(u_{n_{k_{\ell}}})=\liminf\limits_{\ell\to\infty} \mathcal{A}(u_{n_{k_{\ell}}})\geq \mathcal{A}(u)=\overline{\mathcal{A}}(u).
\end{equation}

\bigskip

\noindent\textbf{(3)} From the definition of the subdifferential of a convex and proper functional defined on a Hilbert space\footnote{Consult Chapter 16 from \cite{BauschkeCombettes2017}.} we have that:

\begin{equation}\label{subdif1a}
	\partial\overline{\mathcal{A}}(u):=\left\{w\in L^2(\Omega)\ \bigg| \ \overline{\mathcal{A}}(v)\geq \overline{\mathcal{A}}(u)+ \int_{\Omega} w(v-u)\ dx,\ \forall\ v\in L^2(\Omega) \right\}.
\end{equation}

\noindent If $u\in L^2(\Omega)\setminus W^{1,p(x)}(\Omega)$, then $\partial \overline{\mathcal{A}}(u)=\emptyset$, because otherwise we will have some $w\in \partial \overline{\mathcal{A}}(u)$ and if we choose any $v\in W^{1,p(x)}(\Omega)$ we will get the following contradiction:

 \begin{equation}
 	+\infty>\mathcal{A}(v)=\overline{\mathcal{A}}(v)\geq \overline{\mathcal{A}}(u)+\int_{\Omega} w(v-u)\ dx=+\infty.
 \end{equation}
 
 \noindent Consider now that $u\in W^{1,p(x)}(\Omega)$. For any $v\in L^2(\Omega)\setminus W^{1,p(x)}(\Omega)$ we have that:

 \begin{equation}
 	\overline{\mathcal{A}}(v)=+\infty>\mathcal{A}(u)+\int_{\Omega} w(v-u)\ dx=\overline{\mathcal{A}}(u)+\int_{\Omega} w(v-u)\ dx,\ \forall\ w\in L^2(\Omega).
 \end{equation}
 
 \noindent So, in order to have $w\in \partial\overline{\mathcal{A}}(u)$ we only need to check that:
 
 \begin{equation}
 	\overline{\mathcal{A}}(v)=\mathcal{A}(v)\geq \mathcal{A}(u)+\displaystyle\int_{\Omega} w\cdot\bigl (v-u\bigr)\ dx,\ \forall\ v\in W^{1,p(x)}(\Omega).
 \end{equation}
 
\bigskip

\noindent\textbf{(4)} This follows directly from Corollary 3.9 in \cite[page 61]{brezis2011functional}, taking into account that $\overline{\mathcal{A}}$ is a convex and lower semicontinuous functional.

\noindent The proof is now complete.

\end{proof}

\subsection{Properties of $b$}

\begin{proposition}\label{propoext} For $b:\overline{\Omega}\times [0,\delta]\to\mathbb{R}$ we can define, from \textnormal{\textbf{(H11)}} and \textnormal{\textbf{(H12)}} the following functions:
	
	\begin{equation}
		\begin{cases} L_i:\overline{\Omega}\times [0,\delta]\to\mathbb{R},\ L_i(x,s)=\begin{cases} \dfrac{\partial b}{\partial x_i}(x,s), & (x,s)\in\Omega\times (0,\delta)\\[3mm] \lim\limits_{\underset{(x,s)\in\Omega\times (0,\delta)}{(x,s)\to (x_0,s_0)}}\dfrac{\partial b}{\partial x_i}(x,s), & (x,s)\in\partial\big (\Omega\times (0,\delta)\big )\end{cases},\ \text{where}\ i\in\overline{1,N}\\[3mm] L_{N+1}:\overline{\Omega}\times [0,\delta]\to\mathbb{R},\ L_{N+1}(x,s)=\begin{cases} \dfrac{\partial b}{\partial s}(x,s), & (x,s)\in\Omega\times (0,\delta)\\[3mm] \lim\limits_{\underset{(x,s)\in\Omega\times (0,\delta)}{(x,s)\to (x_0,s_0)}}\dfrac{\partial b}{\partial s}(x,s), & (x,s)\in\partial\big (\Omega\times (0,\delta)\big )\end{cases}\\[3mm] L:\overline{\Omega}\times [0,\delta]\to\mathbb{R},\ L(x,s)=\big (L_1(x,s),L_2(x,s),\hdots, L_N(x,s),L_{N+1}(x,s)\big ).\end{cases}
	\end{equation}
	
	\noindent The following properties hold:
	
	\begin{enumerate}
		\item[\textbf{(1)}] $L\in C\big (\overline{\Omega}\times [0,\delta];\mathbb{R}^{N+1}\big )$, i.e. $L_1,L_2,\hdots,L_N, L_{N+1}\in C(\overline{\Omega}\times [0,\delta])$.
		
		\item[\textbf{(2)}] There is a function $b_{\textnormal{ext}}:\mathbb{R}^{N+1}\to\mathbb{R}$ with $b_{\textnormal{ext}}\in C^1(\mathbb{R}^{N+1})$ such that 
		
		\begin{itemize}
			\item  $b_{\textnormal{ext}}(x,s)=b(x,s),\ \forall\ (x,s)\in\overline{\Omega}\times [0,\delta]$.
			
			\item $\nabla b_{\textnormal{ext}}(x,s)=\big (L_1(x,s),L_2(x,s),\hdots, L_{N}(x,s)\big )$ and $\dfrac{\partial b_{\textnormal{ext}}}{\partial s}(x,s)=L_{N+1}(x,s)$ for any $(x,s)\in \overline{\Omega}\times [0,\delta]$.
			
			\item In particular $\nabla b_{\textnormal{ext}}(x,s)=\nabla b(x,s)$ and $\dfrac{\partial b_{\textnormal{ext}}}{\partial s}(x,s)=\dfrac{\partial b}{\partial s}(x,s)$ for any $(x,s)\in\Omega\times (0,\delta)$.
		\end{itemize}
		
		\begin{remark}\label{remb} \textbf{From now on we shall write for simplicity $b_{\textnormal{ext}}$ as $b$. So we have $b:\mathbb{R}^{N+1}\to\mathbb{R}$ with $b\in C^1(\mathbb{R}^{N+1})$, $b(x,0)=0$ for any $x\in\overline{\Omega}$ and $\dfrac{\partial b}{\partial s}(x,s)\geq \ell_0$ for any $(x,s)\in\overline{\Omega}\times [\varepsilon,\delta]$.}	
		\end{remark}
		
		\item[\textbf{(3)}] There is a constant $L_0>0$ such that: 
		
		\begin{equation}
			\begin{cases} \left |\dfrac{\partial b}{\partial s}(x,s)\right |,\ \left |\dfrac{\partial b}{\partial x_i}(x,s)\right |, |\nabla b (x,s)|\leq L_0\ \text{for any}\ (x,s)\in\overline{\Omega}\times [0,\delta], \ \text{and for each} \ i\in\overline{1,N}\\[5mm]
			|b(x_1,s_1)-b(x_2,s_2)|\leq L_0|x_1-x_2|+L_0|s_1-s_2|,\ \forall\ x_1,x_2\in\overline{\Omega},\ \forall\ s_1,s_2\in [0,\delta].\end{cases}
		\end{equation}
		
		\noindent In particular we have that:
		
		\begin{equation}\label{blipschitz}
			|b(x,s_1)-b(x,s_2)|\leq L_0|s_1-s_2|,\ \forall\ x\in\overline{\Omega},\ \forall\ s_1,s_2\in [0,\delta].
		\end{equation}
	\end{enumerate}
	
	\item[\textbf{(4)}] For each $i\in\overline{1,N}$ and for any $x\in\overline{\Omega}$ we have that $\dfrac{\partial b}{\partial x_i}(x,0)=0$.
	
	\item[\textbf{(5)}] There are some real constants $\tilde{\varepsilon}<\varepsilon<\delta<\tilde{\delta}$ such that 
	
	\begin{itemize}
	\item $\dfrac{\partial b}{\partial s}(x,s)\geq \dfrac{\ell_0}{2}$ for any $(x,s)\in\overline{\Omega}\times [\tilde{\varepsilon},\tilde{\delta}]$. In particular for each $x\in\overline{\Omega}$ the function $[\tilde{\varepsilon},\tilde{\delta}]\ni s\mapsto b(x,s)$ is strictly increasing and we can write that:
	
	\begin{equation}\label{equationell0pe2}
		|b(x,s_1)-b(x,s_2)|\geq \dfrac{\ell_0}{2}|s_1-s_2|,\ \forall\ x\in\overline{\Omega},\ \forall\ s_1,s_2\in [\tilde{\varepsilon},\tilde{\delta}].
	\end{equation}
	
	\item $\left |\dfrac{\partial b}{\partial s}(x,s)\right |,\ \left |\dfrac{\partial b}{\partial x_i}(x,s)\right |, |\nabla b (x,s)|\leq 2L_0\ \text{for any}\ (x,s)\in\overline{\Omega}\times [\tilde{\varepsilon},\tilde{\delta}], \ \text{and for each} \ i\in\overline{1,N}$.
	
	\item 	$|b(x_1,s_1)-b(x_2,s_2)|\leq 2L_0|x_1-x_2|+2L_0|s_1-s_2|,\ \forall\ x_1,x_2\in\overline{\Omega},\ \forall\ s_1,s_2\in [\tilde{\varepsilon},\tilde{\delta}]$.
	
	\end{itemize} 
	
\end{proposition}

\begin{proof} \noindent\textbf{(1), (2)} Note that from Proposition \ref{propolip1} we deduce that $\Omega\times (0,\delta)$ is an open, bounded and connected Lipschitz domain. We also know from \textbf{(H9)} that $b\in C^1(\Omega\times (0,\delta))$ and $b\in C(\overline{\Omega}\times [0,\delta])$. Therefore we are in position to apply Proposition \ref{propolip4} \textbf{(1)} and \textbf{(3)} to deduce exactly the required conclusion.
	
\medskip

\noindent\textbf{(3)} The first conclusion follows from \textit{Weierstrass theorem}, taking into account that the real functions $\left |\dfrac{\partial b}{\partial s}\right |,\ \left |\dfrac{\partial b}{\partial x_i}\right |, |\nabla b|$ are continuous on $\overline{\Omega}\times [0,\delta]$, which is a compact set. The second conclusion follows by a direct application of Proposition \ref{propolip3}.

\medskip

\noindent\textbf{(4)} For $x\in\Omega$ the statement is true from Remark \ref{remdbdxi}. Fix now any $x\in\partial\Omega$. Since $b\in C^1(\mathbb{R}^{N+1})$ we have that $\dfrac{\partial b}{\partial x_i}(x,0)=\lim\limits_{y\to x,\ y\in\Omega} \dfrac{\partial b}{\partial x_i}(y,0)=0$.

\medskip

\noindent\textbf{(5)} Consider the function $h:\mathbb{R}\to \mathbb{R},\ h(s)=\displaystyle\inf_{x\in\overline{\Omega}} \dfrac{\partial b}{\partial s}(x,s)$. Since $ \dfrac{\partial b}{\partial s}\in C\big (\overline{\Omega}\times \mathbb{R}\big )$ and $\overline{\Omega}$ is a compact set, using Theorem \ref{topothm} we deduce that $h$ is continuous. But we already know that $L_0\geq h(s)\geq \ell_0$ for any $s\in [\varepsilon,\delta]$. From the continuity of $h$ there are some $\tilde{\varepsilon}<\varepsilon<\delta<\tilde{\delta}$ such that $2L_0\geq h(s)\geq \dfrac{\ell_0}{2},\ \forall\ s\in [\tilde{\varepsilon},\tilde{\delta}]$.

\noindent Similarly, for each $i\in\overline{1,N}$, the functions $h_i:\mathbb{R}\to \mathbb{R},\ h(s)=\displaystyle\inf_{x\in\overline{\Omega}} \left |\dfrac{\partial b}{\partial x_i}(x,s)\right |$ are also continuous and $h(s)\leq L_0$ for any $s\in [\varepsilon,\delta]$. To obtain the Lipschitz continuity for $b:\overline{\Omega}\times [\tilde{\varepsilon},\tilde{\delta}]\to\mathbb{R}$ we use Proposition \ref{propolip3}. The conclusion follows with ease.
	
\end{proof}

\begin{proposition}\label{propob} The following assertions hold:
	
	\begin{enumerate}
		
		\item[\textbf{(1)}] If $v\in W^{1,p(x)}(\Omega)\cap \mathcal{U}_{[\varepsilon,\delta]}$, then $b\big (\cdot,v(\cdot)\big )\in W^{1,p(x)}(\Omega)$ and:
		
		\begin{equation}
			\dfrac{\partial b\big (x,v(x)\big )}{\partial x_i}=\dfrac{\partial b}{\partial x_i}\big (x,v(x)\big )+\dfrac{\partial b}{\partial s}\big (x,v(x)\big )\cdot\dfrac{\partial v}{\partial x_i}(x),\ \text{a.e. on}\ \Omega.
		\end{equation}
		
		\noindent Moreover, $\Vert b(\cdot,v(\cdot))\Vert_{W^{1,p(x)}(\Omega)}\leq L_0\Vert v\Vert_{W^{1,p(x)}(\Omega)}+L_0\Vert 1\Vert_{L^{p(x)}(\Omega)}$.
		
		\item[\textbf{(2)}] Let $(v_n)_{n\geq 1}\subset W^{1,p(x)}(\Omega)\cap\mathcal{U}_{[\varepsilon,\delta]}$ and $v\in W^{1,p(x)}(\Omega)$ such that $v_n\to v$ in $W^{1,p(x)}(\Omega)$. Then $v\in\mathcal{U}_{[\varepsilon,\delta]}$ and $b(\cdot,v_n)\to b(\cdot,v)$ in $W^{1,p(x)}(\Omega)$.
	\end{enumerate}
\end{proposition}

\begin{proof} \noindent\textbf{(1)} We know from Proposition \ref{propoext} \textbf{(2)} that $b\in C^1\big (\mathbb{R}^{N+1})\ \Rightarrow\ b\in C^1\big (\Omega\times (\tilde{\varepsilon},\tilde{\delta})\big )$. We also deduce from Proposition \ref{propoext} \textbf{(5)} that $b:\overline{\Omega}\times [\tilde{\varepsilon},\tilde{\delta}]\to\mathbb{R}$ is a $2L_0$ -- Lipschitz function.

\noindent Since the mapping $\Omega\ni x\mapsto (x,v(x))\in\Omega\times (\tilde{\varepsilon},\tilde{\delta})$\footnote{Here it is essential that $\Omega$ is bounded and $\delta<\tilde{\delta}$.} is in $W^{1,p(x)}(\Omega)$, we are in position to apply the \textit{Chain Rule} (given as Theorem \ref{thmchainrule} in Appendix) and deduce that $b(\cdot,v(\cdot))\in W^{1,p(x)}(\Omega)$. Moreover for every $i\in\overline{1,N}$ and a.e. $x\in\Omega$ we have that:

\begin{equation*}
\dfrac{\partial}{\partial x_i} \big [b(x,v(x))\big ]=\dfrac{\partial b}{\partial s}(x,v(x))\cdot\dfrac{\partial v}{\partial x_i}(x)+\sum_{j=1}^N \dfrac{\partial b}{\partial x_j}(x,v(x))\cdot\dfrac{\partial x_j}{\partial x_i}=\dfrac{\partial b}{\partial x_i}(x,v(x))+\dfrac{\partial b}{\partial s}(x,v(x))\cdot\dfrac{\partial v}{\partial x_i}(x).
\end{equation*}

\noindent\textbf{(2)} From $v_n\to v$ in $W^{1,p(x)}(\Omega)$ we have that $v_n\to v$ in $L^{p(x)}(\Omega)\hookrightarrow L^1(\Omega)$. Therefore there is a subsequence $(v_{n_k})_{k\geq 1}$ with $v_{n_k}\to v$ pointwise a.e. on $\Omega$\footnote{See \cite{Jones}, page 234.}. Since for each $k\geq 1$ we have that $\varepsilon\leq v_{n_k}(x)\leq \delta$ a.e. on $\Omega$, we deduce by taking $k\to\infty$ that $\varepsilon\leq v(x)\leq \delta$ a.e. on $\Omega$.

\noindent If, for $x\in\Omega$, we set $\mathbf{u}_n(x)=(x,v_n(x))$ and $\mathbf{u}(x)=(x,v(x))$, we have that $\mathbf{u}_n\to \mathbf{u}$ in $W^{1,p(x)}(\Omega)^{N+1}$. Using Theorem \ref{thmchainrule} we get that $b(\cdot,v_n)=b(\mathbf{u}_n)\to b(\mathbf{u})=b(\cdot,v)$ in $W^{1,p(x)}(\Omega)$, as needed.

\noindent The inequality follows as such: first note that from \textbf{(H10)} and Proposition \ref{propoext} \textbf{(3)} we obtain that $|b(x,v(x))|=|b(x,v(x))-b(x,0)|\leq L_0|v(x)|$ for a.e. $x\in\Omega$. Therefore, from Proposition \ref{propu2u1} we get that $\Vert b(\cdot,v)\Vert_{L^{p(x)}(\Omega)}\leq L_0\Vert v\Vert_{L^{p(x)}(\Omega)}$.

\noindent Next, we have that $\nabla [b(x,v(x))]=(\nabla b)(x,v(x))+\dfrac{\partial b}{\partial s}(x,v(x))\nabla v(x)$ for a.e. $x\in\Omega$. Therefore $|\nabla [b(x,v(x))]|\leq |(\nabla b)(x,v(x))|+\left | \dfrac{\partial b}{\partial s}(x,v(x))\right |\cdot |\nabla v(x) |\leq L_0+L_0|\nabla v(x) |$. Using one more time Proposition \ref{propu2u1} we get that $\Vert |\nabla [b(x,v(x))]|\Vert_{L^{p(x)}(\Omega)}\leq L_0\Vert 1\Vert_{L^{p(x)}(\Omega)}+L_0\Vert |\nabla v(x)|\Vert_{L^{p(x)}(\Omega)}$. Summing up leads to $\Vert b(\cdot,v(\cdot))\Vert_{W^{1,p(x)}(\Omega)}\leq L_0\Vert v\Vert_{W^{1,p(x)}(\Omega)}+L_0\Vert 1\Vert_{L^{p(x)}(\Omega)}$.
\end{proof}

\subsection{Properties of $\mathfrak{b}$}

\noindent Let's consider the restriction $b:\overline{\Omega}\times [\tilde{\varepsilon},\tilde{\delta}]\to\mathbb{R}$. For any $x\in\overline{\Omega}$ and any $\sigma\in [b(x,\tilde{\varepsilon}), b(x,\tilde{\delta})]$ there is a unique $s\in [\tilde{\varepsilon},\tilde{\delta}]$ such that $b(x,s)=\sigma$. This is because of the strict monotony and continuity of $b(x,\cdot):[\tilde{\varepsilon},\tilde{\delta}]\to [b(x,\tilde{\varepsilon}),b(x,\tilde{\delta})]$ proved in from Proposition \ref{propoext} \textbf{(5)}. Thus we may define:

\begin{equation}
	\tilde{\mathfrak{b}}:\bigcup_{x\in\overline{\Omega}} \{x\}\times [b(x,\tilde{\varepsilon}), b(x,\tilde{\delta})]\to [\tilde{\varepsilon},\tilde{\delta}],\ \tilde{\mathfrak{b}}(x,\sigma)=s,\ \text{with}\ b(x,s)=\sigma.
\end{equation}

\noindent Consider also the restriction of $\tilde{\mathfrak{b}}$:

\begin{equation}
	\mathfrak{b}:\bigcup_{x\in\overline{\Omega}} \{x\}\times [b(x,\varepsilon), b(x,\delta)]\to [\varepsilon,\delta],\ \mathfrak{b}(x,\sigma)=\tilde{\mathfrak{b}}(x,\sigma).
\end{equation}

\begin{proposition}\label{propobrussian} The following results regarding $\tilde{\mathfrak{b}}$ and $\mathfrak{b}$ hold:
	
	\begin{enumerate}
			\item[\textbf{(1)}] The sets $\widetilde{D}:=\displaystyle\bigcup_{x\in\Omega} \{x\}\times \big (b(x,\tilde{\varepsilon}), b(x,\tilde{\delta})\big )$ and $D:=\displaystyle\bigcup_{x\in\Omega} \{x\}\times \big (b(x,\varepsilon), b(x,\delta)\big )$ are open, and their closures are:
			
			\begin{equation}
				\operatorname{cl}(\widetilde{D})=\bigcup_{x\in\overline{\Omega}} \{x\}\times [b(x,\tilde{\varepsilon}), b(x,\tilde{\delta})],\ \text{and}\ \operatorname{cl}(D)=\bigcup_{x\in\overline{\Omega}} \{x\}\times [b(x,\varepsilon), b(x,\delta)].
			\end{equation}
			
			\begin{remark}\label{rembfraktilde} We may write $\begin{cases} \tilde{\mathfrak{b}}:\operatorname{cl}(\widetilde{D})\to [\tilde{\varepsilon},\tilde{\delta}]\\ \tilde{\mathfrak{b}}:\widetilde{D}\to (\tilde{\varepsilon},\tilde{\delta})\end{cases} $ and $\begin{cases} \mathfrak{b}:\operatorname{cl}(D)\to [\varepsilon,\delta]\\ \mathfrak{b}:D\to (\varepsilon,\delta)\end{cases}$. Of course we have $D\subsetneq \widetilde{D}$, and $\tilde{\mathfrak{b}}(x,\sigma)=\mathfrak{b}(x,\sigma)$ for any $(x,\sigma)\in\operatorname{cl}(D)$.
			\end{remark}

			\item[\textbf{(2)}] For any $x\in\overline{\Omega}$, the function $[b(x,\tilde{\varepsilon}), b(x,\tilde{\delta})]\ni \sigma\mapsto \tilde{\mathfrak{b}}(x,\sigma)$ is strictly increasing. In particular, for any $x\in\overline{\Omega}$ the function $[b(x,\varepsilon),b(x,\delta)]\ni\sigma\mapsto \mathfrak{b}(x,\sigma)$ is strictly increasing.

			\item[\textbf{(3)}] $\tilde{\mathfrak{b}}$ is continuous on $\operatorname{cl}(\widetilde{D})$ and $\mathfrak{b}$ is continuous on $\operatorname{cl}(D)$.

			\item[\textbf{(4)}] $\tilde{\mathfrak{b}}\in C^1(\widetilde{D})$,  $\dfrac{\partial\tilde{\mathfrak{b}}}{\partial \sigma}\in C\big (\operatorname{cl}(\widetilde{D})\big )$ and the following formulas hold:
			
			\begin{equation}
			\dfrac{\partial \tilde{\mathfrak{b}}}{\partial \sigma}(x,\sigma)=\dfrac{1}{\dfrac{\partial b}{\partial s}(x,\tilde{\mathfrak{b}}(x,\sigma))} \ \forall\ (x,\sigma)\in \operatorname{cl}(\widetilde{D}),\quad 	\dfrac{\partial \tilde{\mathfrak{b}}}{\partial x_i}(x,\sigma)=-\dfrac{\dfrac{\partial b}{\partial x_i}(x,\tilde{\mathfrak{b}}(x,\sigma))}{\dfrac{\partial b}{\partial s}(x,\tilde{\mathfrak{b}}(x,\sigma))} \ \forall\ (x,\sigma)\in \widetilde{D}.
		\end{equation}
		
		\noindent Moreover:
		
		\begin{equation}
			\dfrac{1}{2L_0}\leq \dfrac{\partial \tilde{\mathfrak{b}}}{\partial \sigma}(x,\sigma)\leq\dfrac{2}{\ell_0}.\ \forall\  (x,\sigma)\in\operatorname{cl}(\widetilde{D})\ \text{and}\ \left |\dfrac{\partial \tilde{\mathfrak{b}}}{\partial x_i}(x,\sigma) \right |\leq\dfrac{4L_0}{\ell_0},\ \forall\  (x,\sigma)\in\widetilde{D}.
		\end{equation}
		
		\noindent We may also write that:
		
		\begin{equation}
		\dfrac{1}{2L_0}|\sigma_1-\sigma_2|\leq 	|\tilde{\mathfrak{b}}(x,\sigma_1)-\tilde{\mathfrak{b}}(x,\sigma_2)|\leq\dfrac{2}{\ell_0}|\sigma_1-\sigma_2|,\ \forall\ x\in\overline{\Omega},\ \forall\ \sigma_1,\sigma_2\in [b(x,\tilde{\varepsilon}),b(x,\tilde{\delta})].
		\end{equation}

		\item[\textbf{(5)}] $\mathfrak{b}\in C^1(D)$ and if we denote:
			
			\begin{equation}
				D\subsetneq \widehat{D}:=\bigcup_{x\in\Omega} \{x\}\times [b(x,\varepsilon),b(x,\delta)]\subsetneq \operatorname{cl}(D)\subsetneq \widetilde{D},
			\end{equation}
			
			\noindent then for each $i\in\overline{1,N}$: $\dfrac{\partial \mathfrak{b}}{\partial x_i},\dfrac{\partial \mathfrak{b}}{\partial \sigma}\in C(\widehat{D})$ and $\forall\ (x,\sigma)\in \widehat{D}$ we have that:
			
				\begin{equation}
				\dfrac{\partial\tilde{\mathfrak{b}}}{\partial \sigma}(x,\sigma)=\dfrac{\partial \mathfrak{b}}{\partial \sigma}(x,\sigma)=\dfrac{1}{\dfrac{\partial b}{\partial s}(x,\mathfrak{b}(x,\sigma))},\	\dfrac{\partial\tilde{\mathfrak{b}}}{\partial x_i}(x,\sigma)=\dfrac{\partial \mathfrak{b}}{\partial x_i}(x,\sigma)=-\dfrac{\dfrac{\partial b}{\partial x_i}(x,\mathfrak{b}(x,\sigma))}{\dfrac{\partial b}{\partial s}(x,\mathfrak{b}(x,\sigma))}, 
			\end{equation}
			
			\noindent Moreover, for any $(x,\sigma)\in \widehat{D}$:
			
			\begin{equation}
				\dfrac{1}{L_0}\leq \dfrac{\partial \mathfrak{b}}{\partial \sigma}(x,\sigma)\leq\dfrac{1}{\ell_0}\ \text{and}\ \left |\dfrac{\partial \mathfrak{b}}{\partial x_i}(x,\sigma) \right |,|\nabla\mathfrak{b}(x,\sigma)|\leq\dfrac{L_0}{\ell_0}.
			\end{equation}
			
			\noindent We may also write that:
			
			\begin{equation}
				\dfrac{1}{L_0}|\sigma_1-\sigma_2|\leq 	|\mathfrak{b}(x,\sigma_1)-\mathfrak{b}(x,\sigma_2)|\leq\dfrac{1}{\ell_0}|\sigma_1-\sigma_2|,\ \forall\ x\in\overline{\Omega},\ \forall\ \sigma_1,\sigma_2\in [b(x,\varepsilon),b(x,\delta)].
			\end{equation}

			\item[\textbf{(6)}] Let any $w\in W^{1,p(x)}(\Omega)$ with $b(x,\varepsilon)\leq w(x)\leq b(x,\delta)$ a.e. on $\Omega$ and denote $v(x)=\mathfrak{b}\big (x,w(x)\big )$ for a.e. $x\in \Omega$. Then $v\in W^{1,p(x)}(\Omega)\cap\mathcal{U}_{[\varepsilon,\delta]}$ and for each $i\in\overline{1,N}$:
			
			\begin{equation}
				\dfrac{\partial v}{\partial x_i}(x)=\dfrac{\partial\mathfrak{b}}{\partial x_i}\big (x, w(x)\big )+\dfrac{\partial\mathfrak{b}}{\partial\sigma}\big (x,w(x)\big)\dfrac{\partial w}{\partial x_i}.
			\end{equation}
			
			\noindent Moreover $\Vert v\Vert_{W^{1,p(x)}(\Omega)}\leq \dfrac{1}{\ell_0}\Vert w\Vert_{W^{1,p(x)}(\Omega)}+\dfrac{L_0}{\ell_0}\Vert 1\Vert_{L^{p(x)}(\Omega)}+\dfrac{1}{\ell_0}\Vert b(\cdot,\varepsilon)-\ell_0\varepsilon\Vert_{L^{p(x)}(\Omega)}$.
						
			\item[\textbf{(7)}] If the sequence $(w_n)_{n\geq 1}\subset W^{1,p(x)}(\Omega)$ satisfies for each $n\geq 1$ the inequality $b(x,\varepsilon)\leq w_n(x)\leq b(x,\delta)$ a.e. on $\Omega$ and $w\in W^{1,p(x)}(\Omega)$ such that $w_n\to w$ in $W^{1,p(x)}(\Omega)$ then $b(x,\varepsilon)\leq w(x)\leq b(x,\delta)$ a.e. on $\Omega$. Moreover, if we denote $v=\mathfrak{b}(\cdot,w(\cdot))$ and for each $n\geq 1$ we set $v_n=\mathfrak{b}\big (\cdot,w_n(\cdot)\big )$, then $v_n\to v$ in $W^{1,p(x)}(\Omega)$.
			
			\item[\textbf{(8)}] $\tilde{\mathfrak{b}}:\operatorname{cl}\big (\widetilde{D}\big )\to [\tilde{\varepsilon},\tilde{\delta}]$ is a Lipschitz function on $\operatorname{cl}\big (\widetilde{D}\big )$, i.e. $\tilde{\mathfrak{b}}\in \operatorname{Lip}(\operatorname{cl}\big (\widetilde{D}\big ))$, and the following inequality holds
			
			\begin{equation}\label{lipschitzbrustilde}
				\big |\tilde{\mathfrak{b}}(x_1,\sigma_1)-\tilde{\mathfrak{b}}(x_2,\sigma_2)\big |\leq \dfrac{2}{\ell_0}|\sigma_1-\sigma_2|+\dfrac{4L_0}{\ell_0}|x_1-x_2|,\ \forall\ (x_1,\sigma_1),(x_2,\sigma_2)\in \operatorname{cl}\big (\widetilde{D}\big ).
			\end{equation}

			\item[\textbf{(9)}] $\mathfrak{b}:\operatorname{cl}(D)\to [\varepsilon,\delta]$ is a Lipschitz function on $\operatorname{cl}(D)$, i.e. $\mathfrak{b}\in \operatorname{Lip}(\operatorname{cl}(D))$, and the following inequality holds
			
			\begin{equation}\label{lipschitzbrus}
				\big |\mathfrak{b}(x_1,\sigma_1)-\mathfrak{b}(x_2,\sigma_2)\big |\leq \dfrac{1}{\ell_0}|\sigma_1-\sigma_2|+\dfrac{L_0}{\ell_0}|x_1-x_2|,\ \forall\ (x_1,\sigma_1),(x_2,\sigma_2)\in \operatorname{cl}(D).
			\end{equation}
			
			\item[\textbf{(10)}] There is an extension of $\tilde{\mathfrak{b}}$ denoted by $\mathfrak{b}_{\textnormal{ext}}:\mathbb{R}^{N+1}\to\mathbb{R}$ with the following properties
			
			\begin{itemize}
			\item $\mathfrak{b}_{\textnormal{ext}}\in C^1(\mathbb{R}^{N+1})$;
			
			\item $\mathfrak{b}_{\textnormal{ext}}(x,\sigma)=\tilde{\mathfrak{b}}(x,\sigma)$ for any $(x,\sigma)\in\operatorname{cl}\big (\widetilde{D}\big )$;
			
			\item $\mathfrak{b}_{\textnormal{ext}}(x,\sigma)=\mathfrak{b}(x,\sigma)$ for any $(x,\sigma)\in\operatorname{cl}(D)$;
			
			\item $\dfrac{\partial\mathfrak{b}_{\textnormal{ext}}}{\partial\sigma}(x,\sigma)=\dfrac{\partial\tilde{\mathfrak{b}}}{\partial\sigma}(x,\sigma)$ and $\dfrac{\partial\mathfrak{b}_{\textnormal{ext}}}{\partial x_i}(x,\sigma)=\dfrac{\partial\tilde{\mathfrak{b}}}{\partial x_i}(x,\sigma)$ for each $i\in\overline{1,N}$ and for every $(x,\sigma)\in\widetilde{D}$;
			
			\item $\dfrac{\partial\mathfrak{b}_{\textnormal{ext}}}{\partial\sigma}(x,\sigma)=\dfrac{\partial\mathfrak{b}}{\partial\sigma}(x,\sigma)$ and $\dfrac{\partial\mathfrak{b}_{\textnormal{ext}}}{\partial x_i}(x,\sigma)=\dfrac{\partial\mathfrak{b}}{\partial x_i}(x,\sigma)$ for each $i\in\overline{1,N}$ and for every $(x,\sigma)\in D$;
			\end{itemize}
			
			\begin{remark}\label{rembrussian}
				\textbf{From now on we shall write for simplicity $\mathfrak{b}$ instead of the extension $\mathfrak{b}_{\textnormal{ext}}$. So $\mathfrak{b}\in C^1(\mathbb{R}^{N+1})$.}
			\end{remark}
		
	\end{enumerate}
\end{proposition}

\begin{proof} \noindent\textbf{(1)} First we show that $\widetilde{D}$ is open. Let $(x_0,\sigma_0)\in\widetilde{D}$. We know that $b(x_0,\tilde{\varepsilon})<\sigma_0<b(x_0,\widetilde{\delta})$. Consider $\tau=\dfrac{1}{2}\max\{\sigma_0-b(x_0,\tilde{\varepsilon}),b(x_0,\tilde{\delta})-\sigma_0\}$. Thus $b(x_0,\tilde{\varepsilon})<\sigma_0-\tau<\sigma_0<\sigma_0+\tau< b(x_0,\widetilde{\delta})$. We know that $b(\cdot,\tilde{\varepsilon}),b(\cdot,\tilde{\delta}):\Omega\to\mathbb{R}$ are both continuous functions. From the fact that $b(x_0,\tilde{\varepsilon})<\sigma_0-\tau$ we deduce that there is some $r_1>0$ such that $b(x,\tilde{\varepsilon})<\sigma_0-\tau$ for any $x\in B(x_0,r_1)\subset\Omega$. Similarly, from the fact that $b(x_0,\tilde{\delta})>\sigma_0+\tau$ one can find some $r_2>0$ such that $b(x,\tilde{\delta})>\sigma_0+\tau$ for any $x\in B(x_0,r_2)\subset\Omega$. We used the fact that $\Omega$ is an open set and $x_0\in\Omega$. Taking $r:=\min\{r_1,r_2\}>0$ we have that for any $(x,\sigma)\in B(x_0,r)\times (\sigma_0-\tau,\sigma_0+\tau)$ the following relation holds: $b(x,\tilde{\varepsilon})<\sigma_0-\tau<\sigma<\sigma_0+\tau<b(x,\tilde{\delta})$. Hence $B(x_0,r)\times (\sigma_0-\tau,\sigma_0+\tau)\subset \widetilde{D}$. This proves that $\widetilde{D}$ is open. In the same manner we obtain that $D$ is also open.
	
\bigskip
	
\noindent Next, we show that $\operatorname{cl}(\widetilde{D})=\displaystyle\bigcup_{x\in\overline{\Omega}} \{x\}\times [b(x,\tilde{\varepsilon}), b(x,\tilde{\delta})]$ by double inclusion.

\noindent Let's say $(x_0,\sigma_0)\in \operatorname{cl}(\widetilde{D})$. Therefore we find a sequence $(x_n,\sigma_n)\in\widetilde{D}$ with $(x_n,\sigma_n)\to (x_0,\sigma_0)$. Since $x_n\in\Omega, \ \forall\ n\geq 1$ and $x_n\to x_0$ we easily deduce that $x_0\in\overline{\Omega}$. We also know that $b(x_n,\tilde{\varepsilon})<\sigma_n<b(x_n,\tilde{\delta}),\ \forall\ n\geq 1$. Making $n\to\infty$ in this relation and using the continuity of $b$, gives us that: $b(x_0,\tilde{\varepsilon})\leq\sigma_0\leq b(x_0,\tilde{\delta})$. Thus $(x_0,\sigma_0)\in \displaystyle\bigcup_{x\in\overline{\Omega}} \{x\}\times [b(x,\tilde{\varepsilon}), b(x,\tilde{\delta})]$. 

\noindent Now consider any $(x_0,\sigma_0)\in \displaystyle\bigcup_{x\in\overline{\Omega}} \{x\}\times [b(x,\tilde{\varepsilon}), b(x,\tilde{\delta})]$, i.e. $x_0\in\overline{\Omega}$ and $b(x_0,\tilde{\varepsilon})\leq \sigma_0\leq b(x_0,\tilde{\delta})$. So, there is a unique $\lambda\in [0,1]$ such that: $\sigma_0=(1-\lambda)b(x_0,\tilde{\varepsilon})+\lambda b(x_0,\tilde{\delta})$. Since $x_0\in\overline{\Omega}$ there is some sequence $(x_n)_{n\geq 1}\subset\Omega$ with $x_n\to x_0$. Now define for each $n\geq 1$:

\[
\sigma_n=(1-\lambda) \underbrace{\left [b(x_n,\tilde{\varepsilon})+\dfrac{b(x_n,\tilde{\delta})-b(x_n,\tilde{\varepsilon})}{n+1} \right ]}_{\in (b(x_n,\tilde{\varepsilon}),b(x_n,\tilde{\delta}))}+\lambda\underbrace{\left [b(x_n,\tilde{\delta})- \dfrac{b(x_n,\tilde{\delta})-b(x_n,\tilde{\varepsilon})}{n+1}\right ]}_{\in (b(x_n,\tilde{\varepsilon}),b(x_n,\tilde{\delta}))}\in (b(x_n,\tilde{\varepsilon}),b(x_n,\tilde{\delta})).
\]

\noindent Note that $\sigma_n\to(1-\lambda)b(x_0,\tilde{\varepsilon})+\lambda b(x_0,\tilde{\delta})=\sigma_0$. In conclusion, we have that $(x_n,\sigma_n)\in\widetilde{D}$ for each $n\geq 1$ and $(x_n,\sigma_n)\to (x_0,\sigma_0)$ as $n\to\infty$. This proves that $(x_0,\sigma_0)\in\operatorname{cl}(\widetilde{D})$, as needed. We have proved that $\operatorname{cl}(\widetilde{D})=\displaystyle\bigcup_{x\in\overline{\Omega}} \{x\}\times [b(x,\tilde{\varepsilon}), b(x,\tilde{\delta})]$. Exactly in the same way one can show that $\operatorname{cl}(D)=\displaystyle\bigcup_{x\in\overline{\Omega}} \{x\}\times [b(x,\varepsilon), b(x,\delta)]$. 

\medskip

\noindent\textbf{(2)} Let any $b(x,\tilde{\varepsilon})\leq \sigma_1<\sigma_2\leq b(x,\tilde{\delta})$, and denote $\begin{cases} \tilde{\mathfrak{b}}(x,\sigma_1)=s_1\in [\tilde{\varepsilon},\tilde{\delta}] \ \Longrightarrow\ b(x,s_1)=\sigma_1\\[2mm] \tilde{\mathfrak{b}}(x,\sigma_2)=s_2\in [\tilde{\varepsilon},\tilde{\delta}] \ \Longrightarrow\ b(x,s_2)=\sigma_2\end{cases}$. Since, from Proposition \ref{propoext} \textbf{(5)}, we know that $b(x,\cdot)$ is strictly increasing on $[\tilde{\varepsilon},\tilde{\delta}]$ and $b(x,s_1)<b(x,s_2)$, we conclude that $s_1<s_2$, i.e. $\tilde{\mathfrak{b}}(x,\sigma_1)<\tilde{\mathfrak{b}}(x,\sigma_2)$.

\medskip
	
	\noindent\textbf{(3)} Fix any $(x_0,\sigma_0)\in\operatorname{cl}(\widetilde{D})$ and any sequence $(x_n,\sigma_n)\in \operatorname{cl}(\widetilde{D})$ with $(x_n,\sigma_n)\to (x_0,\sigma_0)$. We want to show that $s_n\to s_0$, where $s_n=\tilde{\mathfrak{b}}(x_n,\sigma_n)\in [\tilde{\varepsilon},\tilde{\delta}]$ for each $n\geq 1$ and $s_0=\tilde{\mathfrak{b}}(x_0,\sigma_0)\in [\tilde{\varepsilon},\tilde{\delta}]$. We have that:
	
	\begin{align*}
		|\sigma_n-\sigma_0|&=|b(x_n,s_n)-b(x_0,s_0)|=|b(x_n,s_n)-b(x_n,s_0)+b(x_n,s_0)-b(x_0,s_0)|\\
		&\geq |b(x_n,s_n)-b(x_n,s_0)|-|b(x_n,s_0)-b(x_0,s_0)|\\
	\eqref{equationell0pe2}\	&\geq \dfrac{\ell_0}{2}|s_n-s_0|-|b(x_n,s_0)-b(x_0,s_0)|.
	\end{align*}
	
	\noindent Hence, we can write that $|s_n-s_0|\leq \dfrac{2}{\ell_0}|\sigma_n-\sigma_0|+\dfrac{2}{\ell_0}|b(x_n,s_0)-b(x_0,s_0)|\stackrel{n\to\infty}{\longrightarrow}0$, from the continuity of $b$. So $\tilde{\mathfrak{b}}\in C(\operatorname{cl}(\widetilde{D}))$. Using Remark \ref{rembfraktilde} we also obtain that $\mathfrak{b}\in C(\operatorname{cl}(D))$.
	
	\medskip
	
	\noindent\textbf{(4)} We start by showing that $\dfrac{\partial\tilde{\mathfrak{b}}}{\partial \sigma}(x,\sigma)$ exists for any $(x,\sigma)\in\operatorname{cl}(\widetilde{D})$. Consider any sequence $(\sigma_n)_{n\geq 1}$ with $(x,\sigma_n)\in \operatorname{cl}(\widetilde{D})$ and $\sigma_n\to\sigma$ as $n\to\infty$.
	
	\noindent We denote $\begin{cases} \tilde{\mathfrak{b}}(x,\sigma_n)=s_n \ \Longrightarrow\ b(x,s_n)=\sigma_n \\[2mm] \tilde{\mathfrak{b}}(x,\sigma)=s \ \Longrightarrow\ b(x,s)=\sigma\end{cases}$. Using the continuity of $\tilde{\mathfrak{b}}$ on $\operatorname{cl}(\widetilde{D})$ we get that $s_n\to s$ as $n\to\infty$. Therefore:
	
	\begin{align*}
		\lim\limits_{n\to\infty} \dfrac{\tilde{\mathfrak{b}}(x,\sigma_n)-\tilde{\mathfrak{b}}(x,\sigma)}{\sigma_n-\sigma}&=	\lim\limits_{n\to\infty}\dfrac{s_n-s}{b(x,s_n)-b(x,s)}=\lim\limits_{n\to\infty}\dfrac{1}{\dfrac{b(x,s_n)-b(x,s)}{s_n-s}}\\
		&=\dfrac{1}{\dfrac{\partial b}{\partial s}(x,s)}=\dfrac{1}{\dfrac{\partial b}{\partial s}(x,\tilde{\mathfrak{b}}(x,\sigma))} :=\dfrac{\partial\tilde{\mathfrak{b}}}{\partial \sigma}(x,\sigma).
	\end{align*}
	
	\noindent Since $\dfrac{\ell_0}{2}\leq \dfrac{\partial b}{\partial s}(x,s)\leq 2L_0$ we get that $\dfrac{1}{2L_0}\leq\dfrac{\partial\tilde{\mathfrak{b}}}{\partial \sigma}(x,\sigma)\leq\dfrac{2}{\ell_0},\ \forall\ (x,\sigma)\in \operatorname{cl}(\widetilde{D})$. Also, using the continuity of $\tilde{\mathfrak{b}}$ on $\operatorname{cl}(\widetilde{D})$ -- proved at \textbf{(3)} -- and the continuity of $\dfrac{\partial b}{\partial s}$ we deduce that $\dfrac{\partial\tilde{\mathfrak{b}}}{\partial \sigma}\in C\big (\operatorname{cl}(\widetilde{D})\big )$.
	
	\noindent Fix any $(x_0,\sigma_0)\in \widetilde{D}$ and any $i\in\overline{1,N}$. We denote $\tilde{\mathfrak{b}}(x_0,\sigma_0)=s_0\in (\tilde{\varepsilon},\tilde{\delta})$. Since $\Omega$ is an open set we get that there is some $\tau>0$ such that $x_0+he_i\in\Omega$ for any $h\in (-\tau,\tau)$. We define the function $G:(-\tau,\tau)\times (\tilde{\varepsilon},\tilde{\delta})\to\mathbb{R},\ G(h,s)=b(x_0+he_i,s)-\sigma_0$. It exhibits the following properties:
	
	\begin{itemize}
		\item $G(0,s_0)=b(x_0,s_0)-\sigma_0=0$.
		
		\item  $G\in C^1\big ((-\tau,\tau)\times (\tilde{\varepsilon},\tilde{\delta})\big )$.
		
		\item $\dfrac{\partial G}{\partial s}(h,s)=\dfrac{\partial b}{\partial s}(x_0+he_i,s)\geq\dfrac{\ell_0}{2}>0$, for any $(h,s)\in (-\tau,\tau)\times (\tilde{\varepsilon},\tilde{\delta})$, from Proposition \ref{propoext} \textbf{(5)}.
		
		\item $G(h,s)=0\ \Longleftrightarrow\ b(x_0+he_i,s)=\sigma_0\ \Longleftrightarrow\ s=\tilde{\mathfrak{b}}(x_0+he_i,\sigma_0):=\alpha(h)$.
	\end{itemize}
	
	\noindent Now we are in position to apply the \textit{Implicit function theorem}\footnote{See Extension 3 of the \textit{General Implicit Function Theorem}, given in \cite[page 101]{zorich2016mathematical}.} and obtain that:
	
	\begin{itemize}
		\item $\alpha\in C^1\big ((-\tau,\tau)\big )$
		
		\item $\alpha(0)=s_0$
		
		\item For any $h\in (-\tau,\tau)$:
		
		\begin{equation}
			\alpha'(h)=-\dfrac{\dfrac{\partial{G}}{\partial h}(h,\alpha(h))}{\dfrac{\partial G}{\partial s}(h,\alpha(h)) }=-\dfrac{\dfrac{\partial b}{\partial x_i}(x_0+he_i,\alpha(h))}{\dfrac{\partial b}{\partial s}(x_0+he_i,\alpha(h))}=-\dfrac{\dfrac{\partial b}{\partial x_i}(x_0+he_i,\tilde{\mathfrak{b}}(x_0+he_i,\sigma_0))}{\dfrac{\partial b}{\partial s}(x_0+he_i,\tilde{\mathfrak{b}}(x_0+he_i,\sigma_0))}.
		\end{equation}
		
		\noindent In particular $\alpha'(0)=-\dfrac{\dfrac{\partial b}{\partial x_i}(x_0,\tilde{\mathfrak{b}}(x_0,\sigma_0))}{\dfrac{\partial b}{\partial s}(x_0,\tilde{\mathfrak{b}}(x_0,\sigma_0))}$.
	\end{itemize}
	
	\noindent Next, just note that:
	
	\begin{align*}
		\dfrac{\partial\tilde{\mathfrak{b}}}{\partial x_i}(x_0,\sigma_0) &:=\lim\limits_{h\to 0}\dfrac{\tilde{\mathfrak{\mathfrak{b}}}(x_0+he_i,\sigma_0)-\tilde{\mathfrak{b}}(x_0,\sigma_0)}{h}=\lim\limits_{h\to 0}\dfrac{\alpha(h)-\alpha(0)}{h}=\alpha'(0)=-\dfrac{\dfrac{\partial b}{\partial x_i}(x_0,\tilde{\mathfrak{b}}(x_0,\sigma_0))}{\dfrac{\partial b}{\partial s}(x_0,\tilde{\mathfrak{b}}(x_0,\sigma_0))}.
	\end{align*}
	
	\noindent Moreover, since from Proposition \ref{propoext} \textbf{(5)} we have that $\left |\dfrac{\partial b}{\partial x_i}(x,s) \right |\leq 2L_0$ and $\dfrac{\partial b}{\partial s}(x,s)\geq\dfrac{\ell_0}{2}$ for any $(x,s)\in\Omega\times (\tilde{\varepsilon},\tilde{\delta})$, we deduce that: $\left |\dfrac{\partial\tilde{\mathfrak{b}}}{\partial x_i}(x,\sigma) \right |\leq \dfrac{4L_0}{\ell_0},\ \forall\ (x,\sigma)\in \widetilde{D}$.
	
	\noindent Using the continuity of $\dfrac{\partial b}{\partial x_i},\ \dfrac{\partial b}{\partial s}$ and the continuity of $\tilde{\mathfrak{b}}$ on $\widetilde{D}$, we obtain that $\dfrac{\partial\tilde{\mathfrak{b}}}{\partial x_i}\in C(\widetilde{D})$. Thus, all partial derivatives of $\tilde{\mathfrak{b}}$ exists and are continuous on $\widetilde{D}$. This means that $\tilde{\mathfrak{b}}\in C^1(\widetilde{D})$.
	
	\noindent\textbf{(5)} Since $\mathfrak{b}=\tilde{\mathfrak{b}}$ on the open set $D\subset \widetilde{D}$, from \textbf{(4)} we get that $\mathfrak{b}=\tilde{\mathfrak{b}}\in C^1(D)$. As a consequence on $D$ we have that $\dfrac{\partial\mathfrak{b}}{\partial\sigma}=\dfrac{\partial\tilde{\mathfrak{b}}}{\partial\sigma}$ and $\dfrac{\partial\mathfrak{b}}{\partial x_i}=\dfrac{\partial\tilde{\mathfrak{b}}}{\partial x_i}$ for each $i\in\overline{1,N}$. Recall that: $\tilde{\mathfrak{b}}=\mathfrak{b}\ \text{on}\ \widehat{D}$.
		 
	\noindent Let now any $i\in\overline{1,N}$ and any $(x_0,\sigma_0)\in \widehat{D}\ \Longrightarrow\ x_0\in\Omega$ and $\tilde{\mathfrak{b}}(x_0,\sigma_0)=\mathfrak{b}(x_0,\sigma_0)\in [\varepsilon,\delta]$. 
	
	\noindent So $\begin{cases} \dfrac{\partial b}{\partial s}\big (x_0,\tilde{\mathfrak{b}}(x_0,\sigma_0)\big )=\dfrac{\partial b}{\partial s}\big (x_0,\mathfrak{b}(x_0,\sigma_0)\big )\\[3mm] \dfrac{\partial b}{\partial x_i}\big (x_0,\tilde{\mathfrak{b}}(x_0,\sigma_0)\big )=\dfrac{\partial b}{\partial x_i}\big (x_0,\mathfrak{b}(x_0,\sigma_0)\big ) \end{cases}$. Therefore:

 \begin{align*}
 &\dfrac{\partial\tilde{\mathfrak{b}}}{\partial \sigma}(x_0,\sigma_0)=\dfrac{1}{\dfrac{\partial b}{\partial s}\big (x_0,\tilde{\mathfrak{b}}(x_0,\sigma_0)\big )}=\dfrac{1}{\dfrac{\partial b}{\partial s}\big (x_0,\mathfrak{b}(x_0,\sigma_0)\big )}\ \text{and} \\ &\dfrac{\partial\tilde{\mathfrak{b}}}{\partial x_i}(x_0,\sigma_0)=-\dfrac{\dfrac{\partial b}{\partial x_i}(x_0,\tilde{\mathfrak{b}}(x_0,\sigma_0))}{\dfrac{\partial b}{\partial s}(x_0,\tilde{\mathfrak{b}}(x_0,\sigma_0))}=-\dfrac{\dfrac{\partial{b}}{\partial x_i}(x_0,\mathfrak{b}(x_0,\sigma_0))}{\dfrac{\partial b}{\partial s}(x_0,\mathfrak{b}(x_0,\sigma_0))}.
 \end{align*}
 
 \noindent When $(x_0,\sigma_0)\in D$ we know that $\dfrac{\partial \mathfrak{b}}{\partial \sigma}(x_0,\sigma_0)=\dfrac{\partial \tilde{\mathfrak{b}}}{\partial \sigma}(x_0,\sigma_0)$ and  $\dfrac{\partial \mathfrak{b}}{\partial x_i}(x_0,\sigma_0)=\dfrac{\partial \tilde{\mathfrak{b}}}{\partial x_i}(x_0,\sigma_0)$, so we are done.
 
 \noindent If $\sigma_0=b(x_0,\varepsilon)$ or $\sigma_0=b(x_0,\delta)$ just remark that:
 
 \[
 \dfrac{\partial\mathfrak{b}}{\partial x_i}(x_0,\sigma_0):=\lim\limits_{h\to 0}\dfrac{\mathfrak{b}(\overbrace{x_0+he_i}^{\in\Omega},\sigma_0)-\mathfrak{b}(x_0,\sigma_0)}{h}=\lim\limits_{h\to 0}\dfrac{\tilde{\mathfrak{b}}(x_0+he_i,\sigma_0)-\tilde{\mathfrak{b}}(x_0,\sigma_0)}{h}=\dfrac{\partial\tilde{\mathfrak{b}}}{\partial x_i}(x_0,\sigma_0),
 \]
 
 \noindent and
 
 \begin{align*}
\dfrac{\partial\mathfrak{b}}{\partial\sigma}(x_0,b(x,\varepsilon))&:=\lim\limits_{\sigma\to b(x,\varepsilon)^+}\dfrac{\mathfrak{b}(\overbrace{x_0,\sigma}^{\in\widehat{D}})-\mathfrak{b}(\overbrace{x_0,b(x,\varepsilon)}^{\in\widehat{D}})}{\sigma-b(x,\varepsilon)}=\lim\limits_{\sigma\to b(x,\varepsilon)^+}\dfrac{\tilde{\mathfrak{b}}(x_0,\sigma)-\tilde{\mathfrak{b}}(x_0,b(x,\varepsilon))}{\sigma-b(x,\varepsilon)}\\
&=\dfrac{\partial\tilde{\mathfrak{b}}}{\partial\sigma}(x_0,b(x,\varepsilon))\\
\dfrac{\partial\mathfrak{b}}{\partial\sigma}(x_0,b(x,\delta))&:=\lim\limits_{\sigma\to b(x,\delta)^-}\dfrac{\mathfrak{b}(\overbrace{x_0,\sigma}^{\in\widehat{D}})-\mathfrak{b}(\overbrace{x_0,b(x,\delta)}^{\in\widehat{D}})}{\sigma-b(x,\delta)}=\lim\limits_{\sigma\to b(x,\delta)^-}\dfrac{\tilde{\mathfrak{b}}(x_0,\sigma)-\tilde{\mathfrak{b}}(x_0,b(x,\delta))}{\sigma-b(x,\delta)}\\
&=\dfrac{\partial\tilde{\mathfrak{b}}}{\partial\sigma}(x_0,b(x,\delta)).
 \end{align*}

\noindent We have proved so far that on $\widehat{D}$: $\dfrac{\partial \mathfrak{b}}{\partial \sigma}=\dfrac{\partial \tilde{\mathfrak{b}}}{\partial \sigma}\stackrel{\widehat{D}\subset \widetilde{D}}{\in} C(\widehat{D})$ and  $\dfrac{\partial \mathfrak{b}}{\partial x_i}=\dfrac{\partial \tilde{\mathfrak{b}}}{\partial x_i}\stackrel{\widehat{D}\subset \widetilde{D}}{\in} C(\widehat{D})$.

\noindent The rest of the proof is just an immediate consequence of \textbf{(H13)}, Proposition \ref{propoext} \textbf{(3)} and Remark \ref{rem23}.

	\medskip

	\noindent\textbf{(6)} Recall that $\tilde{\mathfrak{b}}:\widetilde{D}\to (\tilde{\varepsilon},\tilde{\delta}),\ \tilde{\mathfrak{b}}\in C^1(\widetilde{D})$ and $\widetilde{D}$ is an open set from $\mathbb{R}^{N+1}$. Consider $\mathbf{u}:\Omega\to\widetilde{D},\ \mathbf{u}(x)=(x,w(x))\in \widetilde{D}$. This is because for each $x\in\Omega$:
	
	\[
	b(x,\tilde{\varepsilon})<b(x,\varepsilon)\leq w(x)\leq b(x,\delta)<b(x,\tilde{\delta}).
	\]
	
	\noindent Taking into account that the canonical projections are smooth functions from $C^{\infty}(\Omega)\subset W^{1,p(x)}(\Omega)$ and $w\in W^{1,p(x)}(\Omega)$ we deduce that $\mathbf{u}\in W^{1,p(x)}(\Omega)^{N+1}$. Using now the \textit{General Chain Rule} -- i.e. Theorem \ref{thmchainrule} from the Appendix --  we get that $v=\tilde{\mathfrak{b}}\circ \mathbf{u}\in W^{1,p(x)}(\Omega)$ and for each $i\in\overline{1,N}$:
	
	\begin{equation}
		\dfrac{\partial v}{\partial x_i}(x)=\dfrac{\partial\tilde{\mathfrak{b}}}{\partial \sigma}\big (\mathbf{u}(x)\big )\dfrac{\partial w}{\partial x_i}(x)+\displaystyle\sum_{j=1}^{N} \dfrac{\partial\tilde{\mathfrak{b}}}{\partial x_i}\big (\mathbf{u}(x)\big )\dfrac{\partial x_j}{\partial x_i}=\dfrac{\partial\tilde{\mathfrak{b}}}{\partial x_i}\big (x,w(x)\big )+\dfrac{\partial\tilde{\mathfrak{b}}}{\partial \sigma}\big (x,w(x)\big )\dfrac{\partial w}{\partial x_i}(x),
	\end{equation}
	
	\noindent for almost all $x\in\Omega$. Note that $(x,w(x))\in\widehat{D}$ for a.e. $x\in\Omega$ and therefore, from \textbf{(5)} it follows that $\dfrac{\partial\tilde{\mathfrak{b}}}{\partial x_i}\big (x,w(x)\big )=\dfrac{\partial \mathfrak{b}}{\partial x_i}\big (x,w(x)\big )$ and $\dfrac{\partial\tilde{\mathfrak{b}}}{\partial \sigma}\big (x,w(x)\big )=\dfrac{\partial \mathfrak{b}}{\partial \sigma}\big (x,w(x)\big )$. In conclusion:
	
	\begin{equation}
		\dfrac{\partial v}{\partial x_i}(x)=\dfrac{\partial \mathfrak{b}}{\partial x_i}\big (x,w(x)\big )+\dfrac{\partial \mathfrak{b}}{\partial\sigma}\big (x,w(x)\big )\dfrac{\partial w}{\partial x_i}(x), \text{a.e. on}\ \Omega.
	\end{equation}
	
	\noindent We also have to show here that $v\in\mathcal{U}_{[\varepsilon,\delta]}$. This is obvious, since from \textbf{(2)} we get that $\varepsilon=\mathfrak{b}(x,b(x,\varepsilon))\leq v(x)=\mathfrak{b}(x,w(x))\leq \mathfrak{b}(x,b(x,\delta))=\delta$, a.e. on $\Omega$.
	
	\noindent The inequality can be easily derived. First, note that $\varepsilon\leq v(x)=\mathfrak{b}(x,w(x))=\mathfrak{b}(x,w(x))-\mathfrak{b}(x,b(x,\varepsilon))+\varepsilon\leq \dfrac{1}{\ell_0}\big (w(x)-b(x,\varepsilon)\big )+\varepsilon=\dfrac{1}{\ell_0}w(x)+\dfrac{1}{\ell_0}(\varepsilon\ell_0-b(x,\varepsilon))$. So, from Proposition \ref{propu2u1}: $\Vert v\Vert_{L^{p(x)}(\Omega)}\leq \dfrac{1}{\ell_0}\Vert w\Vert_{L^{p(x)}(\Omega)}+\dfrac{1}{\ell_0}\Vert b(\cdot,\varepsilon)-\varepsilon\ell_0\Vert_{L^{p(x)}(\Omega)}$. Also $\left |\nabla v(x) \right |\leq |(\nabla\mathfrak{b})(x,w(x))|+\left | \dfrac{\partial\mathfrak{b}}{\partial\sigma}(x,w(x))\right |\cdot |\nabla w(x)|\leq\dfrac{L_0}{\ell_0}+\dfrac{1}{\ell_0}|\nabla w(x)|$. Again, using Proposition \ref{propu2u1} gives us that $\Vert |\nabla v|\Vert_{L^{p(x)}(\Omega)}\leq \dfrac{L_0}{\ell_0}\Vert 1\Vert_{L^{p(x)}(\Omega)}+\dfrac{1}{\ell_0}\Vert |\nabla w|\Vert_{L^{p(x)}(\Omega)}$. Adding the two obtained inequalities gives us the desired result.
	
	\medskip
	
	\noindent\textbf{(7)} From $w_n\to w$ in $W^{1,p(x)}(\Omega)$ we get that $w_n\to w$ in $L^{p(x)}(\Omega)\hookrightarrow L^1(\Omega)$. So there is a subsequence $(w_{n_k})_{k\geq 1}$ such that $w_{n_k}\to w$ pointwise a.e. on $\Omega$. Since for each $k\geq 1$ we have that $b(x,\varepsilon)\leq w_{n_k}(x)\leq b(x,\delta)$ for a.e. $x\in\Omega$, by passing $k\to\infty$ we get that $b(x,\varepsilon)\leq w(x)\leq b(x,\delta)$ for a.e. $x\in\Omega$, as needed.

	\noindent Now defining for any $n\geq 1$: $\mathbf{u}_n:\Omega\to \widetilde{D},\ \mathbf{u}_n(x)=(x,w_n(x))$ we have that $\mathbf{u}_n\in W^{1,p(x)}(\Omega)^{N+1}$ and $\mathbf{u}_n\to \mathbf{u}$ in $W^{1,p(x)}(\Omega)^{N+1}$, where $\mathbf{u}:\Omega\to\widetilde{D},\ \mathbf{u}(x)=(x,w(x))$. Now using Theorem \ref{thmchainrule} we get that $\tilde{\mathfrak{b}}\circ\mathbf{u}_n\to \tilde{\mathfrak{b}}\circ\mathbf{u}$ in $W^{1,p(x)}(\Omega)$. Just note that $\tilde{\mathfrak{b}}(\mathbf{u}_n(x))=\tilde{\mathfrak{b}}(\underbrace{x,w_n(x)}_{\in\widehat{D}})=\mathfrak{b}(x,w_n(x))=v_n(x)$ a.e. on $\Omega$. Similarly: $\tilde{\mathfrak{b}}(\mathbf{u}(x))=\tilde{\mathfrak{b}}(\underbrace{x,w(x)}_{\in\widehat{D}})=\mathfrak{b}(x,w(x))=v(x)$. We conclude that $v_n\to v$ in $W^{1,p(x)}(\Omega)$ as stated.
	
	\medskip

	\noindent\textbf{(8)} Fix any $(x_1,\sigma_1),(x_2,\sigma_2)\in\operatorname{cl}\big (\widetilde{D}\big )$ and denote $\begin{cases}\tilde{\mathfrak{b}}(x_1,\sigma_1):=s_1\in [\tilde{\varepsilon},\tilde{\delta}]\ \Longrightarrow\ b(x_1,s_1)=\sigma_1\\ \tilde{\mathfrak{b}}(x_2,\sigma_2):=s_2\in [\tilde{\varepsilon},\tilde{\delta}]\ \Longrightarrow\ b(x_2,s_2)=\sigma_2\end{cases}$.
	
	\noindent Observe that $(x_1,s_1),(x_2,s_2)\in \overline{\Omega}\times [\tilde{\varepsilon},\tilde{\delta}]$ and
	
	\begin{align*}
		|\sigma_1-\sigma_2|&=|b(x_1,s_1)-b(x_2,s_2)|\geq |b(x_1,s_1)-b(x_1,s_2)|-|b(x_1,s_2)-b(x_2,s_2)|\\
\eqref{equationell0pe2}\ \ \ 	 &\geq \dfrac{\ell_0}{2}|s_1-s_2|-|b(x_1,s_2)-b(x_2,s_2)|\\
\text{Prop. \ref{propoext} \textbf{(5)}}\ \ \ &\geq  \dfrac{\ell_0}{2}|s_1-s_2|-2L_0|x_1-x_2|\\
&= \dfrac{\ell_0}{2}|\tilde{\mathfrak{b}}(x_1,\sigma_1)-\tilde{\mathfrak{b}}(x_2,\sigma_2)|-2L_0|x_1-x_2|.
	\end{align*}
	
	\noindent Therefore $|\tilde{\mathfrak{b}}(x_1,\sigma_1)-\tilde{\mathfrak{b}}(x_2,\sigma_2)|\leq \dfrac{2}{\ell_0}|\sigma_1-\sigma_2|+\dfrac{4L_0}{\ell_0}|x_1-x_2|$, as needed.

		\medskip

	\noindent\textbf{(9)} Fix any $(x_1,\sigma_1),(x_2,\sigma_2)\in\operatorname{cl}(D)$ and denote $\begin{cases}\mathfrak{b}(x_1,\sigma_1):=s_1\in [\varepsilon,\delta]\ \Longrightarrow\ b(x_1,s_1)=\sigma_1\\ \mathfrak{b}(x_2,\sigma_2):=s_2\in [\varepsilon,\delta]\ \Longrightarrow\ b(x_2,s_2)=\sigma_2\end{cases}$.
	
	\noindent Observe that $(x_1,s_1),(x_2,s_2)\in \overline{\Omega}\times [\varepsilon,\delta]$ and
	
	\begin{align*}
		|\sigma_1-\sigma_2|&=|b(x_1,s_1)-b(x_2,s_2)|\geq |b(x_1,s_1)-b(x_1,s_2)|-|b(x_1,s_2)-b(x_2,s_2)|\\
		\text{Remark \ref{rem23}}\ \ \ 	 &\geq \ell_0|s_1-s_2|-|b(x_1,s_2)-b(x_2,s_2)|\\
		\text{Prop. \ref{propoext} \textbf{(3)}}\ \ \ &\geq  \ell_0|s_1-s_2|-L_0|x_1-x_2|\\
		&= \ell_0|\tilde{\mathfrak{b}}(x_1,\sigma_1)-\tilde{\mathfrak{b}}(x_2,\sigma_2)|-L_0|x_1-x_2|.
	\end{align*}
	
	\noindent Therefore $|\mathfrak{b}(x_1,\sigma_1)-\mathfrak{b}(x_2,\sigma_2)|\leq \dfrac{1}{\ell_0}|\sigma_1-\sigma_2|+\dfrac{L_0}{\ell_0}|x_1-x_2|$, as needed.

		\medskip

	\noindent\textbf{(10)} For each $i\in\overline{1,N}$ we define the functions $\widetilde{L}_i:\operatorname{cl}\big (\widetilde{D}\big )\to\mathbb{R},\ \widetilde{L}_i(x,\sigma)=-\dfrac{\dfrac{\partial b}{\partial x_i}(x,\tilde{\mathfrak{b}}(x,\sigma))}{\dfrac{\partial b}{\partial s}(x,\tilde{\mathfrak{b}}(x,\sigma))}$ and then we define $\widetilde{L}_{N+1}:\operatorname{cl}\big (\widetilde{D}\big )\to\mathbb{R},\ \widetilde{L}_i(x,\sigma)=\dfrac{1}{\dfrac{\partial b}{\partial s}(x,\tilde{\mathfrak{b}}(x,\sigma))}$. Set $\widetilde{L}:\operatorname{cl}\big (\widetilde{D}\big )\to\mathbb{R}^{N+1},\ \widetilde{L}(x,\sigma)=\big (\widetilde{L}_1(x,\sigma),\widetilde{L}_2(x,\sigma),\hdots, \widetilde{L}_{N}(x,\sigma),\widetilde{L}_{N+1}(x,\sigma)\big )$ for every $(x,\sigma)\in \operatorname{cl}\big (\widetilde{D}\big )$.
	
	\noindent Notice that all these functions are well-defined, because $\dfrac{\partial b}{\partial s}(x,\tilde{\mathfrak{b}}(x,\sigma))\geq\dfrac{\ell_0}{2}>0$ for any $(x,\sigma)\in \operatorname{cl}\big (\widetilde{D}\big )$. They are also continous on $\operatorname{cl}\big (\widetilde{D}\big )$, because $\tilde{\mathfrak{b}}\in C\big (\operatorname{cl}\big (\widetilde{D}\big ) \big )$ -- from \textbf{(4)} --, $\dfrac{\partial b}{\partial x_i}\in C\big (\mathbb{R}^{N+1} \big )$ for each $i\in\overline{1,N}$ and $\dfrac{\partial b}{\partial s}\in C\big (\mathbb{R}^{N+1} \big )$, from Remark \ref{remb}. Untill now we have that $\tilde{\mathfrak{b}}\in C\big (\operatorname{cl}\big (\widetilde{D}\big ) \big )$ and $\widetilde{L}\in C\big (\operatorname{cl}\big (\widetilde{D}\big );\mathbb{R}^{N+1}\big )$. Consider now any two distinct elements $(x_1,\sigma_1)=(x_{11},x_{12},\hdots,x_{1N},\sigma_1)\neq (x_2,\sigma_2)=(x_{21},x_{22},\hdots,x_{2N},\sigma_2)$ from $\operatorname{cl}\big (\widetilde{D}\big )$ -- which is a compact set. We set $\tilde{\mathfrak{b}}(x_1,\sigma_1)=s_1\ \Rightarrow\ b(x_1,s_1)=\sigma_1$ and $\tilde{\mathfrak{b}}(x_2,\sigma_2)=s_2\ \Rightarrow\ b(x_2,s_2)=\sigma_2$. Now we may write that:
	
	\begin{align*}
		&\dfrac{\big |\tilde{\mathfrak{b}}(x_2,\sigma_2)-\tilde{\mathfrak{b}}(x_1,\sigma_1)-\widetilde{L}(x_1,\sigma_1)\cdot ((x_2,\sigma_2)-(x_1,\sigma_1))\big |}{|x_2-x_1|+|\sigma_2-\sigma_1|}=\\
		=& \dfrac{1}{|x_2-x_1|+|b(x_2,s_2)-b(x_1,s_1)|}\cdot\left |s_2-s_1-\sum_{i=1}^N\widetilde{L}_i(x_1,\sigma_1)\cdot(x_{2i}-x_{1i})-L_{N+1}(x_1,\sigma_1)\cdot (\sigma_2-\sigma_1)\right |\\
		=&\dfrac{\left |s_2-s_1+\displaystyle\sum_{i=1}^N\dfrac{\frac{\partial b}{\partial x_i}(x_1,s_1)}{\frac{\partial b}{\partial s}(x_1,s_1)}\cdot (x_{2i}-x_{1i})-\dfrac{1}{\frac{\partial b}{\partial s}(x_1,s_1)}\cdot \big (b(x_2,s_2)-b(x_1,s_1)\big ) \right |}{|x_2-x_1|+|b(x_2,s_2)-b(x_1,s_1)|} \\
		=&\dfrac{\left |b(x_2,s_2)-b(x_1,s_1)-\displaystyle\sum_{i=1}^N \frac{\partial b}{\partial x_i}(x_1,s_1)\cdot (x_{2i}-x_{1i})-\frac{\partial b}{\partial s}(x_1,s_1)\cdot (s_2-s_1) \right |}{\frac{\partial b}{\partial s}(x_1,s_1)[|x_2-x_1|+|b(x_2,s_2)-b(x_1,s_1)|]}\\
		=&\dfrac{1}{\frac{\partial b}{\partial s}(x_1,s_1)}\cdot \dfrac{|x_2-x_1|+|s_2-s_1|}{|x_2-x_1|+|b(x_2,s_2)-b(x_1,s_1)|}\cdot\dfrac{\big | b(x_2,s_2)-b(x_1,s_1)-L(x_1,s_1)\cdot ((x_2,s_2)-(x_1,s_1))\big |}{|x_2-x_1|+|s_2-s_1|}\\
		\leq &\dfrac{2}{\ell_0}\cdot \dfrac{|x_2-x_1|+|s_2-s_1|}{|x_2-x_1|+\frac{\ell_0}{2}|s_2-s_1|}\cdot\dfrac{\big | b(x_2,s_2)-b(x_1,s_1)-L(x_1,s_1)\cdot ((x_2,s_2)-(x_1,s_1))\big |}{|x_2-x_1|+|s_2-s_1|}\\
		\leq &\dfrac{2\max\{1,\frac{2}{\ell_0}\}}{\ell_0}\cdot\dfrac{\big | b(x_2,s_2)-b(x_1,s_1)-L(x_1,s_1)\cdot ((x_2,s_2)-(x_1,s_1))\big |}{|x_2-x_1|+|s_2-s_1|}.
	\end{align*} 
	
	\noindent It is important to see that from \textbf{(8)} we have that:
	
	\[
	|x_2-x_1|+|s_2-s_1|=|x_2-x_1|+|\tilde{\mathfrak{b}}(x_2,\sigma_2)-\tilde{\mathfrak{b}}(x_1,\sigma_1)|\leq \dfrac{2}{\ell_0}|\sigma_2-\sigma_1|+\dfrac{4L_0+\ell_0}{\ell_0}|x_2-x_1|.
	\]
	
	\noindent This shows that if $(x_2,\sigma_2)-(x_1,\sigma_1)\to 0$ then $(x_2,s_2)-(x_2,s_1)\to 0$.

	\noindent Since we already know that $b:\overline{\Omega}\times [\tilde{\varepsilon},\tilde{\delta}]\to\mathbb{R}$ admits a $C^1$ extension on $\mathbb{R}^{N+1}$ we deduce from the \textit{Converse of Whitney's Extension Theorem} -- i.e. Theorem \ref{thmwhit} --  that 
	
	\[
\lim\limits_{\underset{(x_2,s_2),(x_1,s_1)\in\overline{\Omega}\times [\tilde{\varepsilon},\tilde{\delta}]}{|(x_2,s_2)-(x_1,s_1)|\to 0}}	\dfrac{\big | b(x_2,s_2)-b(x_1,s_1)-L(x_1,s_1)\cdot ((x_2,s_2)-(x_1,s_1))\big |}{|x_2-x_1|+|s_2-s_1|}=0.
	\]
	
	\noindent Putting all together we get from the \textit{squeezing principle} that
	
	\[
	\lim\limits_{\underset{(x_2,\sigma_2),(x_1,\sigma_1)\in \operatorname{cl}(\widetilde{D})}{|(x_2,\sigma_2)-(x_1,\sigma_1)|\to 0}}\dfrac{\big |\tilde{\mathfrak{b}}(x_2,\sigma_2)-\tilde{\mathfrak{b}}(x_1,\sigma_1)-\widetilde{L}(x_1,\sigma_1)\cdot ((x_2,\sigma_2)-(x_1,\sigma_1))\big |}{|x_2-x_1|+|\sigma_2-\sigma_1|}=0.
	\]
	
	\noindent We are now in position to apply \textit{Whitney's Extension Theorem} and deduce that there is a function $\mathfrak{b}_{\text{ext}}:\mathbb{R}^{N+1}\to\mathbb{R}$ with $\mathfrak{b}_{\text{ext}}\in C^1\big (\mathbb{R}^{N+1}\big )$ that extends $\tilde{\mathfrak{b}}$. The listed properties are straightforward to deduce.
\end{proof}

\subsection{Properties of $\overline{b}$}

\begin{proposition}\label{propoverb} The following properties of $\overline{b}:\overline{\Omega}\times\mathbb{R}\to \mathbb{R}$ hold:
	
	\begin{enumerate}
		\item[\textbf{(1)}] $\overline{b}\in C\big (\overline{\Omega}\times\mathbb{R} \big )$ and for any $x\in\overline{\Omega}$, the function $\mathbb{R}\ni s\mapsto \overline{b}(x,s)$ strictly increasing.
		
		\item[\textbf{(2)}] For any $x\in\overline{\Omega}$ and for any $s_1,s_2\in\mathbb{R}$: $|\overline{b}(x,s_1)-\overline{b}(x,s_2)|\leq \overline{L}_0|s_1-s_2|$, where $\overline{L}_0=\max\{1,L_0\}$.
		
		\item[\textbf{(3)}] We have that the Nemytskii operator $\mathcal{N}_{\overline{b}}:L^2(\Omega)\to L^2(\Omega)$ is an $\overline{L}_0$-Lipschitz operator, i.e.:
		
		\begin{equation}
			\Vert\mathcal{N}_{\overline{b}}(V_1)-\mathcal{N}_{\overline{b}}(V_2) \Vert_{L^2(\Omega)}\leq \overline{L}_0\Vert V_1-V_2\Vert_{L^2(\Omega)},\ \forall\ V_1,V_2\in L^2(\Omega).
		\end{equation}
		
		\noindent Furthermore, for any $r\in [1,\infty)$, $\mathcal{N}_{\overline{b}}:L^r(\Omega)\to L^r(\Omega)$ is an $\overline{L}_0$ -- Lipschitz operator.
		
		\item[\textbf{(4)}]  We have that the Nemytskii operator $\mathcal{N}_{\overline{b}}:L^2\big ((0,T)\times \Omega\big )\to L^2\big ((0,T)\times \Omega\big )$ is an $\overline{L}_0$-Lipschitz operator, i.e.:
		
		\begin{equation}
			\Vert\mathcal{N}_{\overline{b}}(v_1)-\mathcal{N}_{\overline{b}}(v_2) \Vert_{L^2((0,T)\times\Omega)}\leq \overline{L}_0\Vert v_1-v_2\Vert_{L^2((0,T)\times \Omega)},\ \forall\ v_1,v_2\in L^2\big ((0,T)\times\Omega\big ).
		\end{equation}
		
		\noindent Similarly, for any $r\in [1,\infty)$, $\mathcal{N}_{\overline{b}}:L^r((0,T)\times \Omega)\to L^r((0,T)\times\Omega)$ is an $\overline{L}_0$ -- Lipschitz operator.
	\end{enumerate}
	
\end{proposition}

\begin{proof}\noindent \textbf{(1)} It follows from \textbf{(H9)} that $\overline{b}$ is continuous on $\overline{\Omega}\times \bigl [(-\infty,\varepsilon)\cup (\varepsilon,\delta)\cup (\delta,\infty)\bigr ]$. We check the continuity of $\overline{b}$ at $(x_0,\varepsilon)$ and $(x_0,\delta)$ for any $x_0\in\overline{\Omega}$. Indeed:
	
	\begin{equation*}
		\begin{cases} \lim\limits_{(x,s)\to (x_0,\varepsilon^-)} \overline{b}(x,s)=\lim\limits_{(x,s)\to (x_0,\varepsilon^-)} b(x,\varepsilon)+s-\varepsilon\stackrel{\textbf{(H9)}}{=}b(x_0,\varepsilon)=\overline{b}(x_0,\varepsilon)\\ \lim\limits_{(x,s)\to (x_0,\varepsilon^+)} \overline{b}(x,s)=\lim\limits_{(x,s)\to (x_0,\varepsilon^+)} b(x,s)\stackrel{\textbf{(H9)}}{=}b(x_0,\varepsilon)=\overline{b}(x_0,\varepsilon)\end{cases},
	\end{equation*}
	
	\noindent and:
	
	\begin{equation*}
		\begin{cases} \lim\limits_{(x,s)\to (x_0,\delta^+)} \overline{b}(x,s)=\lim\limits_{(x,s)\to (x_0,\delta^+)} b(x,\delta)+s-\delta\stackrel{\textbf{(H9)}}{=}b(x_0,\delta)=\overline{b}(x_0,\delta)\\ \lim\limits_{(x,s)\to (x_0,\delta^-)} \overline{b}(x,s)=\lim\limits_{(x,s)\to (x_0,\delta^-)} b(x,s)\stackrel{\textbf{(H9)}}{=}b(x_0,\delta)=\overline{b}(x_0,\delta)\end{cases}.
	\end{equation*}
	
	\noindent Finally, since $\mathbb{R}\ni s\mapsto \overline{b}(x,s)$ is continuous and strictly increasing on $(-\infty,\varepsilon)$, $(\varepsilon,\delta)$, $(\delta,\infty)$ we get that it will be strictly increasing on $\mathbb{R}$. 
	
	\noindent\textbf{(2)} Let any real numbers $s_1<s_2$. We have the following cases:
	
	\begin{itemize}
		\item $s_1<s_2\leq\varepsilon\ \Rightarrow |\overline{b}(x,s_1)-\overline{b}(x,s_2)|=|s_1-s_2|\leq \overline{L}_0|s_1-s_2|$.
		
		\item $\delta\leq s_1<s_2\ \Rightarrow |\overline{b}(x,s_1)-\overline{b}(x,s_2)|=|s_1-s_2|\leq \overline{L}_0|s_1-s_2|$.
		
		\item $\varepsilon\leq  s_1<s_2\leq \delta\ \Rightarrow |\overline{b}(x,s_1)-\overline{b}(x,s_2)|=|b(x,s_1)-b(x,s_2)|\leq L_0|s_1-s_2|\leq \overline{L}_0|s_1-s_2|$.
		
		\item $s_1\leq \varepsilon<s_2\leq \delta\ \Rightarrow |\overline{b}(x,s_1)-\overline{b}(x,s_2)|=|b(x,\varepsilon)+s_1-\varepsilon -b(x,s_2)|\leq$ $|b(x,\varepsilon)-b(x,s_2)|+|s_1-\varepsilon|\leq L_0|\varepsilon-s_2|+|s_1-\varepsilon|\leq \overline{L}_0 |\varepsilon-s_2|+\overline{L}_0|s_1-\varepsilon|=\overline{L}_0 (s_2-s_1)= \overline{L}_0|s_1-s_2|$.
		
		\item $s_1\leq\varepsilon<\delta\leq s_2\ \Rightarrow\ |\overline{b}(x,s_1)-\overline{b}(x,s_2)|=\overline{b}(x,s_2)-\overline{b}(x,s_1)= -b(x,\varepsilon)-s_1+\varepsilon+b(x,\delta)+s_2-\delta=$ $b(x,\delta)-b(x,\varepsilon)-(\delta-\varepsilon)+s_2-s_1\leq  (L_0-1)(\delta-\varepsilon)+(s_2-s_1)\leq (\overline{L}_0-1)(\delta-\varepsilon)+(s_2-s_1)\leq (\overline{L}_0-1)(s_2-s_1)+(s_2-s_1)=  \overline{L}_0(s_2-s_1)=\overline{L}_0|s_1-s_2|$.
		
		\item $\varepsilon\leq s_1\leq \delta\leq s_2\ \Rightarrow\ |\overline{b}(x,s_1)-\overline{b}(x,s_2)|=\overline{b}(x,s_2)-\overline{b}(x,s_1)=b(x,\delta)+s_2-\delta-b(x,s_1)=b(x,\delta)-b(x,s_1)+s_2-\delta_0\leq L_0(\delta-s_1)+s_2-\delta\leq \overline{L}_0(\delta-s_1)+\overline{L}_0(s_2-\delta)=\overline{L}_0(s_2-s_1)=\overline{L}_0|s_1-s_2|$.
	\end{itemize}
	
	\noindent The proof in now complete.

	\noindent\textbf{(3)} Note that for any $r\in [1,\infty)$ and $V\in L^{r}(\Omega)$, we have that:
	
	\begin{align*}
		|\mathcal{N}_{\overline{b}}(V)(x)|&=\bigl |\overline{b}(x,V(x)) \bigr |=\begin{cases} |b(x,\varepsilon)+V(x)-\varepsilon|, & V(x)<\varepsilon\\ |b(x,V(x))|, & V(x)\in [\varepsilon,\delta] \\ |b(x,\delta)+V(x)-\delta|, & V(x)>\delta \end{cases}\\
		&\leq |V(x)|+\max\{\delta,\Vert b(\cdot,\varepsilon)\Vert_{L^{\infty}(\Omega)}, \Vert b(\cdot,\delta)\Vert_{L^{\infty}(\Omega)}\}\in L^{r}(\Omega).
	\end{align*}
	
	\noindent Thus $\mathcal{N}_{\overline{b}}(V)\in L^r(\Omega)$. Using \textbf{(2)} we get that:
	
	\begin{align*}
		\Vert \mathcal{N}_{\overline{b}}(V_1)-\mathcal{N}_{\overline{b}}(V_2)\Vert_{L^r(\Omega)}&=\left (\int_{\Omega} [\overline{b}(x,V_1(x))-\overline{b}(x,V_2(x))] ^r\ dx \right )^{\frac{1}{r}}\\
		& \leq \overline{L}_0\left (\int_{\Omega} [V_1(x)-V_2(x)]^r\ dx \right )^{\frac{1}{r}}=\overline{L}_0 \Vert V_1-V_2\Vert_{L^r(\Omega)}.
	\end{align*}
	
	\noindent\textbf{(4)} For any $r\in [1,\infty)$ and $v\in L^{r}\big ((0,T)\times\Omega\big )$, we have that:
	
	\begin{align*}
		|\mathcal{N}_{\overline{b}}(v)(t,x)|&=\bigl |\overline{b}(x,v(t,x)) \bigr |=\begin{cases} |b(x,\varepsilon)+v(t,x)-\varepsilon|, & v(t,x)<\varepsilon\\ |b(x,v(t,x))|, & v(t,x)\in [\varepsilon,\delta] \\ |b(x,\delta)+v(t,x)-\delta|, & v(t,x)>\delta \end{cases}\\
		&\leq |v(t,x)|+\max\{\delta, \Vert b(\cdot,\varepsilon)\Vert_{L^{\infty}(\Omega)}, \Vert b(\cdot,\delta)\Vert_{L^{\infty}(\Omega)}\}\in L^{r}(\Omega).
	\end{align*}
	
	\noindent Thus $\mathcal{N}_{\overline{b}}(v)\in L^r\big ((0,T)\times\Omega\big )$. Using \textbf{(2)} we get for any $v_1,v_2\in L^r\big ((0,T)\times\Omega\big )$ that:
	
	\begin{align*}
		\Vert \mathcal{N}_{\overline{b}}(v_1)-\mathcal{N}_{\overline{b}}(v_2)\Vert_{L^r((0,T)\times\Omega)}&=\left (\int_{(0,T)\times\Omega} [\overline{b}(x,v_1(t,x))-\overline{b}(x,v_2(t,x))] ^r\ dx \right )^{\frac{1}{r}}\\
		& \leq \overline{L}_0\left (\int_{(0,T)\times\Omega} [v_1(t,x)-v_2(t,x)]^r\ dx\ dt \right )^{\frac{1}{r}}\\
		&=\overline{L}_0 \Vert v_1-v_2\Vert_{L^r((0,T)\times\Omega)}.
	\end{align*}
	
\end{proof}

\subsection{Properties of $\overline{f}$}

\begin{proposition}\label{propoverf} The following properties of $\overline{f}:\overline{\Omega}\times\mathbb{R}\to \mathbb{R}$ hold:
	
	\begin{enumerate}
		\item[\textbf{(1)}] $\overline{f}$ is a Carath\'{e}odory function.
		
		\item[\textbf{(2)}] For any $\lambda\geq \lambda_0$ the function $\mathbb{R}\ni s\longmapsto\overline{f}(x,s)+\lambda\cdot\overline{b}(x,s)$ is strictly increasing for a.e. $x\in\Omega$. In particular  $\mathbb{R}\ni s\longmapsto g_0(x,s)$ is strictly increasing for a.e. $x\in\Omega$.
		
		\item[\textbf{(3)}] The following inequality holds:
		\begin{equation}
			|\overline{f}(x,s)|\leq \lambda_0 |s|+\dfrac{\lambda_0\varepsilon}{2}+\lambda_0\Vert b(\cdot,\delta)-b(\cdot,\varepsilon) \Vert_{L^{\infty}(\Omega)}
		\end{equation}
		
		\item[\textbf{(4)}] For any $V\in L^2(\Omega)$ we have that $\overline{f}(\cdot, V(\cdot))\in L^2(\Omega)$. This means that $\mathcal{N}_{\overline{f}}:L^2(\Omega)\to L^2(\Omega)$. In consequence, $\mathcal{N}_{\overline{f}}$ is continuous, meaning that if $V_n\to V$ in $L^2(\Omega)$ then $\overline{f}(\cdot, V_n)\to \overline{f}(\cdot, V)$ in $L^2(\Omega)$. This is also true if we replace $L^2(\Omega)$ by $L^r(\Omega)$, for any $r\in [1,\infty)$.
		
		\item[\textbf{(5)}] Similarly, for any $r\in [1,\infty)$ the Nemytskii operator $\mathcal{N}_{\overline{f}}:L^r\big ((0,T)\times\Omega\big )\to L^r\big ((0,T)\times\Omega\big )$ is continuous with respect to the norm of $L^r\big ((0,T)\times\Omega\big )$. This means that if $v_n\to v$ in $L^r\big ((0,T)\times\Omega\big )$ then:
		
		\begin{equation}
			\lim\limits_{n\to\infty} \int_{0}^T\int_{\Omega} |\overline{f}(x,v_n(t,x))-f(x,v(t,x))|^r\ dx\ dt=0.
		\end{equation}
		
	\end{enumerate}
	
\end{proposition}

\begin{proof} \noindent\textbf{(1)} If $s<0$ then $\Omega\ni x\mapsto\overline{f}(x,s)=f(x,\varepsilon)-\dfrac{\lambda}{2}(s-\varepsilon)$ which is measurable, knowing \textbf{(H8)}. If $s\in [\varepsilon,\delta]$, $\Omega\ni x\mapsto\overline{f}(x,s)=f(x,s)$ is also measurable from \textbf{(H8)}. For $s>\delta$ then $\Omega\ni x\mapsto\overline{f}(x,s)=f(x,\delta_0)-\tilde{\lambda}_0(s-\delta_0)$ is clearly measurable from the same hypothesis \textbf{(H8)}. For a.e. $x\in\Omega$ it is easy to see that $\mathbb{R}\ni s\mapsto \overline{f}(x,s)$ is continuous, using again \textbf{(H8)}. Thus $\overline{f}$ is a Carath\'{e}odory function.
	
	\noindent \textbf{(2)} For almost all $x\in\Omega$ the function:
	
	\begin{equation}
		\mathbb{R}\ni s\mapsto\overline{f}(x,s)+\lambda\cdot\overline{b}(x,s)=\begin{cases} f(x,\varepsilon)+\lambda b(x,\varepsilon)+\left (\lambda-\dfrac{\lambda_0}{2}\right )(s-\varepsilon), & s<\varepsilon\\[3mm] f(x,s)+\lambda_0 b(x,s)+(\lambda-\lambda_0) b(x,s), & s\in [\varepsilon,\delta]\\[3mm] f(x,\delta)+\lambda b(x,\delta)+(\lambda-\tilde{\lambda}_0)(s-\delta), & s>\delta \end{cases}
	\end{equation}
	
	\noindent is continuous and strictly increasing on each branch, from \textbf{(H8), (H9)} and \textbf{(H14)}. Here we have used that $\tilde{\lambda}_0<\lambda_0$. This proves the assertion.
	
	\noindent \textbf{(3)} First note that from \textbf{(H14)} we can write for any $(x,s)\in\overline{\Omega}\times [\varepsilon,\delta]$ that: $-\lambda_0\big (b(x,\delta)-b(x,\varepsilon)\big )\leq f(x,\varepsilon)+\lambda_0\big (b(x,\varepsilon)-b(x,s)\big )\leq f(x,s)\leq f(x,\delta)+\lambda_0\big ( b(x,\delta)-b(x,s)\big )\leq \lambda_0\big (b(x,\delta)-b(x,\varepsilon)\big )$. So:
	
	\[
	|\overline{f}(x,s)|\leq \lambda_0 \big (b(x,\delta)-b(x,\varepsilon)\big ),\ \forall\ (x,s)\in\overline{\Omega}\times[\varepsilon,\delta].
	\]
	
	\noindent Now, for any $(x,s)\in\overline{\Omega}\times\mathbb{R}$ we have that:
	
	\begin{align*}
		|\overline{f}(x,s)|&=\begin{cases} \left|f(x,\varepsilon)-\dfrac{\lambda_0}{2}(s-\varepsilon)\right |, & s<0\\[3mm] |f(x,s)|, & s\in [\varepsilon,\delta] \\[3mm] |f(x,\delta)-\tilde{\lambda}_0 (s-\delta)|, & s>\delta\end{cases} \ \ \leq \begin{cases} f(x,\varepsilon)+\dfrac{\lambda_0}{2}|s|+\dfrac{\lambda_0\varepsilon}{2}, & s<\varepsilon\\[3mm] \lambda_0\big (b(x,\delta)-b(x,\varepsilon) \big ), & s\in [\varepsilon,\delta] \nonumber\\[3mm] -f(x,\delta)+\tilde{\lambda}_0|s|-\delta\tilde{\lambda}_0, & s>\delta\end{cases}\\
		&\leq \lambda_0 |s|+\dfrac{\lambda_0\varepsilon}{2}+\lambda_0\Vert b(\cdot,\delta)-b(\cdot,\varepsilon) \Vert_{L^{\infty}(\Omega)}.
	\end{align*}

	\noindent \textbf{(4)} Just note that from \textbf{(3)} we have that for any $V\in L^r(\Omega),\ r>1$:
	
	\[
	\big |\mathcal{N}_{\overline{f}}(V) \big |=\big |\overline{f}(x,V(x))\big |\leq \lambda_0 |V|+\dfrac{\lambda_0\varepsilon}{2}+\lambda_0\Vert b(\cdot,\delta)-b(\cdot,\varepsilon) \Vert_{L^{\infty}(\Omega)}\in L^r(\Omega).
	\]
	
	\noindent The continuity of the Nemystkii operator $\mathcal{N}_{\overline{f}}:L^r(\Omega)\to L^{r}(\Omega)$ with respect to the norm of $L^r(\Omega)$ follows from \cite[Theorem 19.1, page 155]{vain}.

	\noindent\textbf{(5)} Let any $v\in L^r\big ((0,T)\times\Omega\big )$. Thefore, from \textbf{(3)} we have:
	
	\[
	|\mathcal{N}_{\overline{f}}(v)(t,x)|=|\overline{f}(x,v(t,x))|\leq \lambda_0 |v(t,x)|+\dfrac{\lambda_0\varepsilon}{2}+\lambda_0\Vert b(\cdot,\delta)-b(\cdot,\varepsilon) \Vert_{L^{\infty}(\Omega)}\in L^r\big ((0,T)\times\Omega\big ).
	\]
	
	\noindent The continuity of the Nemystkii operator $\mathcal{N}_{\overline{f}}:L^r\big ((0,T)\times\Omega\big )\to L^{r}\big ((0,T)\times\Omega\big )$ with respect to the norm of $L^r\big ((0,T)\times\Omega\big )$ follows again from \cite[Theorem 19.1, page 155]{vain}, but applied for the function $\widehat{h}:\big [(0,T)\times\Omega \big ]\times\mathbb{R}\to\mathbb{R},\ \widehat{h}(t,x,s)=\overline{f}(x,s)$. Therefore $\mathcal{N}_{\widehat{h}}(v)=\widehat{h}(t,x,v(t,x))=\overline{f}(x,v(t,x))$.
	
\end{proof}

\subsection{Properties of $\frak{B}$}

\begin{proposition}\label{propofrakB}
	The following assertions about $\frak{B}$ are true:
	
	\begin{enumerate}
		\item[\textnormal{\textbf{(1)}}] For any $(x,s)\in\overline{\Omega}\times [\varepsilon,\delta]$:
		
		\begin{equation}\label{ineqfrakB}
			0\leq \dfrac{\ell_0}{2}(s^2-\varepsilon^2)\leq B(x,s)\leq \dfrac{L_0}{2}(s^2-\varepsilon^2).
		\end{equation}

		\item[\textnormal{\textbf{(2)}}] $\frak{B}\in C\bigl(\overline{\Omega}\times [\varepsilon,\delta]\bigr)$ and $\dfrac{\partial \frak{B}}{\partial s}(x,s)=s\dfrac{\partial b}{\partial s}(x,s)$ for every $(x,s)\in\overline{\Omega}\times [\varepsilon,\delta]$. Furthermore $\dfrac{\partial \frak{B}}{\partial s}\in C\bigl(\overline{\Omega}\times [\varepsilon,\delta]\bigr)$ and $\left |\dfrac{\partial \frak{B}}{\partial s}(x,s) \right |\leq \delta L_0$ for any $(x,s)\in\overline{\Omega}\times [\varepsilon,\delta]$. In particular:
		
		\begin{equation}\label{fcbotosani1}
			|\frak{B}(x,s_1)-\frak{B}(x,s_2)|\leq \delta L_0 |s_1-s_2|,\ \forall\ (x,s_1),(x,s_2)\in \overline{\Omega}\times [\varepsilon,\delta].
		\end{equation}
		
		\item[\textnormal{\textbf{(3)}}] If $u:(0,T)\times\Omega\to [\varepsilon,\delta]$ and $u\in H^1\bigl((0,T);L^2(\Omega)\bigr)$, then $w:(0,T)\times\Omega\to\mathbb{R}$ given by $w(t,x)=b(x,u(t,x))$ has the same property, i.e. $w\in H^1\bigl((0,T);L^2(\Omega)\bigr)$ and moreover the following formula holds:
		
		\begin{equation}\label{bleah3}
			\dfrac{\partial w}{\partial t}(t,\cdot)=\dfrac{\partial b}{\partial s}(\cdot,u(t,\cdot))\cdot\dfrac{\partial u}{\partial t}(t,\cdot),\ \text{for a.e.}\ t\in (0,T).
		\end{equation}
		
		\item[\textnormal{\textbf{(4)}}] If $u:(0,T)\times\Omega\to [\varepsilon,\delta]$ and $u\in H^1\bigl((0,T);L^2(\Omega)\bigr)$, then $v:(0,T)\times\Omega\to\mathbb{R}$ given by $v(t,x)=\frak{B}(x,u(t,x))$ has the same property, i.e. $v\in H^1\bigl((0,T);L^2(\Omega)\bigr)$ and moreover the following formula holds:
		
		\begin{align}\label{bleah1}
			\dfrac{\partial v}{\partial t}(t,\cdot)&=\dfrac{\partial \frak{B}}{\partial s}(\cdot,u(t,\cdot))\cdot\dfrac{\partial u}{\partial t}(t,\cdot)=\dfrac{\partial b}{\partial s}(\cdot,u(t,\cdot))\cdot u(t,\cdot)\cdot\dfrac{\partial u}{\partial t}(t,\cdot)\nonumber\\
\eqref{bleah3}\ \ \			&=\dfrac{\partial b(\cdot,u(t,\cdot))}{\partial t}\cdot u(t,\cdot),\ \text{for a.e.}\ t\in (0,T).
		\end{align}
		
		\item[\textnormal{\textbf{(5)}}] If $u:(0,T)\times\Omega\to [\varepsilon,\delta]$ with $u\in H^1\bigl((0,T);L^2(\Omega)\bigr)\cap C\bigl([0,T];L^2(\Omega)\bigr)$\footnote{This is not an extra assumption, because $H^1\bigl((0,T);L^2(\Omega)\bigr)\hookrightarrow C\bigl([0,T];L^2(\Omega)\bigr)$ -- see for example \cite[Section 5.9.2, Theorem 2, page 302]{evans2022partial}.}, then for any subinterval $(t_1,t_2)\subseteq (0,T)$ the following formula holds:
		
		\begin{equation}\label{bleah2}
			\int_{t_1}^{t_2}\int_{\Omega} \dfrac{\partial b(x,u(t,x))}{\partial t}\cdot u(t,x)\ dx\ dt=\int_{\Omega}\frak{B}(x,u(t_2,x))\ dx-\int_{\Omega}\frak{B}(x,u(t_1,x))\ dx\ dt.
		\end{equation}
	\end{enumerate}
\end{proposition}

\begin{proof} \noindent\textbf{(1)} From Proposition \ref{propoext} \textbf{(3)} it is easy to remark that 
	
\begin{equation}
	B(x,s)=\displaystyle\int_{\varepsilon}^s r\dfrac{\partial b}{\partial s}(x,r)\ dr\leq L_0\int_{\varepsilon}^s r\ dr=\dfrac{L_0}{2}(s^2-\varepsilon^2),\ \forall\ (x,s)\in\overline{\Omega}\times [\varepsilon,\delta].
\end{equation}

\noindent Similarly, from \textbf{(H13)}

\begin{equation}
	B(x,s)=\displaystyle\int_{\varepsilon}^s r\dfrac{\partial b}{\partial s}(x,r)\ dr\geq \ell_0\int_{\varepsilon}^s r\ dr=\dfrac{\ell_0}{2}(s^2-\varepsilon^2),\ \forall\ (x,s)\in\overline{\Omega}\times [\varepsilon,\delta].
\end{equation}	
	
\bigskip
	
\noindent\textbf{(2)} From Remark \ref{remb} we know that $b\in C^1(\mathbb{R}^{N+1})$, and therefore $\dfrac{\partial b}{\partial s}\in C\bigl(\overline{\Omega}\times [\varepsilon,\delta]\bigr)$. Since the product of two continuous functions in also a continuous function we deduce that $s\cdot\dfrac{\partial b}{\partial s}\in C\bigl(\overline{\Omega}\times [\varepsilon,\delta]\bigr)$. Using now Proposition \ref{pcontB} from the Appendix the conclusion follows immediately. The required inequality follows from Proposition \ref{propoext} \textbf{(3)}, since:
	
	\begin{equation}
		\left |\dfrac{\partial \frak{B}}{\partial s}(x,s) \right |=s\left |\dfrac{\partial b}{\partial s}(x,s) \right |\leq \delta L_0,\ \forall\ (x,s)\in\overline{\Omega}\times [\varepsilon,\delta].
	\end{equation}
	
\bigskip

\noindent\textbf{(3)} From Remark \ref{remb} we have that $b\in C^1(\mathbb{R}^{N+1})$. Therefore we can apply the \textit{Chain Rule for Bochner-Sobolev spaces} (see Remark \ref{rembochner}) for $b:\overline{\Omega}\times [\varepsilon,\delta]\to\mathbb{R}$ and $u:(0,T)\times\Omega\to [\varepsilon,\delta]$. We will get that $w\in H^1\bigl((0,T);L^2(\Omega)\bigr)$ and \eqref{bleah3} holds.

\bigskip

\noindent\textbf{(4)} From \textbf{(2)} we have all that is needed to apply \textit{Chain Rule for Bochner-Sobolev spaces} -- i.e. Theorem \ref{thmbochner} \textbf{(1)} -- for $p=2$, $\frak{B}:\overline{\Omega}\times [\varepsilon,\delta]\to \mathbb{R}$ and $u:(0,T)\times\Omega\to [\varepsilon,\delta]$ with $u\in H^1\bigl((0,T);L^2(\Omega)\bigr)$. Hence we obtain that $v\in H^1\bigl((0,T);L^2(\Omega)\bigr)$ and \eqref{bleah1} holds.

\bigskip

\noindent\textbf{(5)} From \textbf{(4)} we know that $v\in H^1\bigl((0,T);L^2(\Omega)\bigr)$. Since $u\in C\bigl([0,T];L^2(\Omega)\bigr)$, we also have that $v\in C\bigl([0,T];L^2(\Omega)\bigr)$. Indeed, for any $t_0,t\in [0,T]$ we have that:

\begin{align*}
	\Vert v(t,\cdot)-v(t_0,\cdot)\Vert_{L^2(\Omega)}&=\Vert \frak{B}(\cdot,u(t,\cdot))-\frak{B}(\cdot,u(t_0,\cdot))\Vert_{L^2(\Omega)}\\
	&=\left(\int_{\Omega} |\frak{B}(x,u(t,x))-\frak{B}(x,u(t_0,x))|^2\ dx\right)^{\frac{1}{2}}\\
\eqref{fcbotosani1}\ \ \	&\leq \delta L_0 \left(\int_{\Omega} |u(t,x)-u(t_0,x)|^2\ dx\right)^{\frac{1}{2}}\\
&=\delta L_0\Vert u(t,\cdot)-u(t_0,\cdot)\Vert_{L^2(\Omega)}\stackrel{t\to t_0}{\longrightarrow}0.
\end{align*}

\noindent Therefore we have proved that $v\in H^1\bigl((0,T);L^2(\Omega)\bigr)\cap C\bigl([0,T];L^2(\Omega)\bigr)$, and in fact $v\in \textnormal{AC}\bigl([0,T];L^2(\Omega)\bigr)$, from Lemma \ref{lemmaACX}. Moreover, from Lemma \ref{lemmaACX}, we can write for any $0\leq t_1\leq t_2\leq T$ the following equality in $L^2(\Omega)$:

\begin{equation}\label{bleah4}
	\frak{B}(\cdot,u(t_2,\cdot))-\frak{B}(\cdot,u(t_1,\cdot))=v(t_2,\cdot)-v(t_1,\cdot)=\int_{t_1}^{t_2} \dfrac{\partial v}{\partial t}(t,\cdot)\ dt\stackrel{\eqref{bleah1}}{=} \int_{t_1}^{t_2}\dfrac{\partial b(\cdot,u(t,\cdot))}{\partial t}\cdot u(t,\cdot)\ dt.
\end{equation}

\noindent Consider now the following operator $\mathcal{H}:L^2(\Omega)\to\mathbb{R},\ \mathcal{H}(z):=\displaystyle\int_{\Omega} z(x)\ dx$. It is easy to see that $\mathcal{H}$ is a linear operator. Moreover $\mathcal{H}$ is also bounded since from \textit{Cauchy inequality} we get that:

\begin{equation}
	\mathcal{H}(z)=\int_{\Omega} z(x)\ dx\leq \left (\int_{\Omega} 1^2\ dx\right )^{\frac{1}{2}}\cdot \left (\int_{\Omega} z^2(x)\ dx\right )^{\frac{1}{2}}=|\Omega|^{\frac{1}{2}}\Vert z\Vert_{L^2(\Omega)},\ \forall\ z\in L^2(\Omega).
\end{equation}
	
\noindent This shows that $\mathcal{H}\in L^2(\Omega)^*$. Using now Lemma \ref{lemcomute} we obtain that:

\begin{align*}
	&\mathcal{H}\bigl(v(t_2,\cdot)-v(t_1,\cdot) \bigr)\stackrel{\eqref{bleah4}}{=}\mathcal{H}\left (\int_{t_1}^{t_2} \dfrac{\partial v}{\partial t}(t,\cdot)\ dt\right )=\int_{t_1}^{t_2}\mathcal{H}\left (\dfrac{\partial v}{\partial t}(t,\cdot)\right )\ dt\\
	\Longleftrightarrow\ &\mathcal{H}(v(t_2,\cdot))-\mathcal{H}(v(t_1,\cdot))=\int_{\Omega} v(t_2,x)\ dx-\int_{\Omega} v(t_1,x)\ dx=\int_{t_1}^{t_2}\int_{\Omega} \dfrac{\partial v}{\partial t}(t,x)\ dx\ dt\\
	\Longleftrightarrow\ & \int_{\Omega}\frak{B}(x,u(t_2,x))\ dx-\int_{\Omega}\frak{B}(x,u(t_1,x))\ dx\stackrel{\eqref{bleah1}}{=}\int_{t_1}^{t_2}\int_{\Omega}\dfrac{\partial b(x,u(t,x))}{\partial t}\cdot u(t,x) \ dx\ dt.
\end{align*}
\end{proof}

\section{Defining what we mean by a solution}\label{s4}

	\begin{definition}
	We say that $u\in C\big ( [0,T]; L^2(\Omega)\big )\cap H^1_{\textnormal{loc}}((0,T);L^2(\Omega)\big )\cap L^{1}\big (0,T;W^{1,p(x)}(\Omega)\big )$ is a \textbf{weak solution} of \eqref{eqdpg} if $u(0,x)=u_0(x)$ for a.e. $x\in\Omega$ and for a.e. $t\in (0,T)$ we have that:
	
	\begin{equation}
	\int_{\Omega} \dfrac{\partial \overline{b}(x,u(t,x))}{\partial t}\phi\ dx+\int_{\Omega} \mathbf{a}(x,\nabla u(t,x))\cdot\nabla\phi\ dx=\int_{\Omega} \overline{f}(x,u(t,x))\phi\ dx,\ \forall\ \phi\in W^{1,p(x)}(\Omega).
	\end{equation}
\end{definition}

\noindent This definition makes sense because for a.e. $t\in (0,T)$:
\begin{itemize}
	\item $u(t,\cdot)\in W^{1,p(x)}(\Omega)\ \Rightarrow\ \nabla u(t,\cdot)\in L^{p(x)}(\Omega)^N\ \Rightarrow\ \mathbf{a}(\cdot,\nabla u(t,\cdot))\in L^{p'(x)}(\Omega)^N$. But $\nabla\phi\in L^{p(x)}(\Omega)^N$, so the second term makes sense, from \textit{H\"{o}lder inequality}.
	
	\item $u(t,\cdot)\in L^2(\Omega)$ Therefor from Proposition \ref{propoverf} \textbf{(4)} we get that $\overline{f}(\cdot, u(t,\cdot))\in L^2(\Omega)$. But $\phi\in W^{1,p(x)}(\Omega)\hookrightarrow L^2(\Omega)$. So the third integral makes sense, from \textit{Cauchy inequality}.\footnote{Proposition \ref{propoverf} \textbf{(2)}.}
	
	\item From Proposition \ref{propoverb} \textbf{(3)}, the Nemytskii operator $\mathcal{N}_{\overline{b}}:L^2(\Omega)\to L^2(\Omega)$ is a Lipschitz operator. Therefore from \cite[Corollary 3.23, page 36]{kreuter2015sobolev}, since $u\in H^1\big ((a,b); L^2(\Omega)\big )$, we get that $\mathcal{N}_{\overline{b}}\circ u\in H^1\big ((a,b); L^2(\Omega)\big )$ for any $0<a<b<T$. So $\dfrac{\partial \overline{b}(\cdot,u(t,\cdot))}{\partial t}\in L^2(\Omega)$ for a.e. $t\in (0,T)$. Thus the first integral makes sense.

\end{itemize}

\section{The weak parabolic comparison principle}\label{s5}

	\begin{theorem}[\textbf{Weak parabolic comparison principle}]
	For any $\lambda\in\mathbb{R}$, let $v_1,v_2:(0,T)\to W^{1,p(x)}(\Omega)$ with $v_1,v_2\in C\big ([0,T];L^2(\Omega)\big )\cap  H^1_{\textnormal{loc}}((0,T);L^2(\Omega)\big )$, $v_2(0,x)\geq v_1(0,x)$ for a.e. $x\in\Omega$ and $g_1,g_2\in L^2\big (0,T;L^2(\Omega)\big )$ with $g_2\geq g_1$ a.e. on $(0,T)\times\Omega$, satisfying in a \textbf{weak sense} the following inequalities:
	
	\begin{equation}
		\begin{cases}\dfrac{\partial \overline{b}(x,v_2)}{\partial t}-\operatorname{div}\mathbf{a}(x,\nabla v_2)+\lambda \overline{b}(x,v_2)=g_2\geq \\[2mm] 
			\geq g_1=\dfrac{\partial \overline{b}(x,v_1)}{\partial t}-\operatorname{div}\mathbf{a}(x,\nabla v_1)+\lambda \overline{b}(x,v_1), \ (t,x)\in (0,T)\times\Omega,  \\[2mm] \mathbf{a}(x,\nabla v_2)\cdot\nu\geq \mathbf{a}(x,\nabla v_1)\cdot\nu,\ (t,x)\in (0,T)\times\partial\Omega\end{cases}
	\end{equation}
	
	\noindent i.e. for a.e. $t\in (0,T)$ and for any $\phi\in W^{1,p(x)}(\Omega)^+$:
	
	\begin{align*}
		&\int_\Omega\dfrac{\partial \overline{b}\big (x,v_2(t,x)\big )}{\partial t}\phi\ dx+\int_{\Omega}\mathbf{a}\big (x,\nabla v_2(t,x)\big )\cdot\nabla\phi\ dx+\lambda\int_{\Omega} \overline{b}\big (x,v_2(t,x)\big )\phi\ dx=\int_{\Omega} g_2(t,x)\phi(x)\ dx\\
		\geq &\int_{\Omega} g_1(t,x)\phi(x)\ dx= \int_\Omega\dfrac{\partial \overline{b}\big (x,v_1(t,x)\big )}{\partial t}\phi\ dx+\int_{\Omega}\mathbf{a}\big (x,\nabla v_1(t,x)\big )\cdot\nabla\phi\ dx+\lambda\int_{\Omega} \overline{b}\big (x,v_1(t,x)\big )\phi\ dx.
	\end{align*}
	
	\noindent Then for any $t\in [0,T]$ we have that $v_2(t,x)\geq v_1(t,x)$ for a.e. $x\in\Omega$.
\end{theorem}

\begin{proof} For a.e. $t\in (0,T)$ and for any $\phi\in W^{1,p(x)}(\Omega)^+$ it follows that:

\begin{align}\label{inegalitate1}
	&\int_{\Omega} \dfrac{\partial \big [\overline{b}(x,v_2(t,x))-\overline{b}(x,v_1(t,x))\big ]}{\partial t}\phi\ dx+\int_{\Omega}\big [\mathbf{a}(x,\nabla v_2(t,x))-\mathbf{a}(x,\nabla v_1(t,x))\big ]\cdot\nabla \phi\ dx\\ 
	&+\lambda\int_\Omega\big [\overline{b}(x,v_2(t,x))-\overline{b}(x,v_1(t,x))]\phi\ dx\geq 0.\nonumber
\end{align}

\noindent For a fixed $t\in (0,T)$ and any $\tau>0$, define: $\phi_{\tau}=\begin{cases} 1, &  v_1(t,\cdot)-v_2(t,\cdot)\geq \tau\\[3mm] \dfrac{v_1(t,\cdot)-v_2(t,\cdot)}{\tau}, & |v_1(t,\cdot)-v_2(t,\cdot)|<\tau \\[3mm] -1, & v_1(t,\cdot)-v_2(t,\cdot)\leq -\tau\end{cases}$. It is easy to remark that $\phi_{\tau}\in W^{1,p(x)}(\Omega)$. Indeed, since $\mathbb{R}\ni s\mapsto\begin{cases} 1, & s\geq \tau\\[2mm] \dfrac{s}{\tau}, & s\in (-\tau,\tau )\\[2mm] -1, & s\leq -\tau\end{cases}$ is a $\dfrac{1}{\tau}$ -- Lipschitz function and $v_1(t,\cdot)-v_2(t,\cdot)\in W^{1,p(x)}(\Omega)$, the conclusion follows from the \textit{Chain Rule}\footnote{See Theorem \ref{chainrulevar} from the Appendix.}. Moreover $\nabla\phi_{\tau}(x)=\begin{cases} \mathbf{0}, & \text{if}\ |v_1(t,x)-v_2(t,x)|\geq\tau\\ \dfrac{1}{\tau}\bigl [\nabla v_1(t,x)-\nabla v_2(t,x) \bigr ], & \text{if}\ |v_1(t,x)-v_2(t,x)|<\tau \end{cases}$.

\noindent From Theorem \ref{apthplus} given in the Appendix, we have that $\phi_{\tau}^+\in W^{1,p(x)}(\Omega)^+$ and

\begin{align*}
	\nabla\phi_{\tau}^+(x)&=\begin{cases}\mathbf{0}, & \text{if}\ \phi_{\tau}(x)\leq 0\\ \nabla \phi_{\tau}(x), & \text{if}\ \phi_{\tau}(x)>0 \end{cases}=\nabla\phi_{\tau}(x)\chi_{\phi_{\tau}>0}(x)\\
	&=\dfrac{1}{\tau}\bigl [\nabla v_1(t,x)-\nabla v_2(t,x) \bigr ]\chi_{|v_1(t,\cdot)-v_2(t,\cdot)|<\tau}(x)\chi_{v_{1}(t,\cdot)>v_2(t,\cdot)}(x)\\
	&=\dfrac{1}{\tau}\bigl [\nabla v_1(t,x)-\nabla v_2(t,x) \bigr ]\chi_{v_1(t,\cdot)-v_2(t,\cdot)\in (0,\tau)}(x).
\end{align*}

\noindent Choose $\phi=\phi_{\tau}^+$ as test function in \eqref{inegalitate1}. Using the \textbf{strict monotony} of $\mathbf{a}$\footnote{Take a look at Proposition 3.1 \textbf{(1)} from \cite{max4}.} we get that:

\begin{align}\label{inegalitatea2}
	&\int_{\Omega} \dfrac{\partial \big [\overline{b}(x,v_2)-\overline{b}(x,v_1)\big ]}{\partial t}\phi_{\tau}^+\ dx 
	+\lambda\int_\Omega\big [\overline{b}(x,v_2)-\overline{b}(x,v_1)]\phi_{\tau}^+\ dx\\
	&\geq \dfrac{1}{\tau}\int_{\Omega}\big [\mathbf{a}(x,\nabla v_2)-\mathbf{a}(x,\nabla v_1)\big ]\cdot\big (\nabla v_2-\nabla v_1\big )\chi_{v_1(t,\cdot)-v_2(t,\cdot)\in (0,\tau)}\ dx  \geq 0.\nonumber
\end{align}

\noindent Note that $\phi_{\tau}^+(x)=\begin{cases}\phi_{\tau}(x), & \phi_{\tau}(x)> 0 \\ 0, & \phi_{\tau}(x)\leq 0 \end{cases}=\begin{cases}1, & v_1(t,x)-v_2(t,x)\geq \tau\\[2mm]  \dfrac{v_1(t,x)-v_2(t,x)}{\tau}, & v_1(t,x)-v_2(t,x)\in (0,\tau) \\[2mm] 0, & v_1(t,x)-v_2(t,x)\leq 0\end{cases}$. We have $\lim\limits_{\tau\to 0^+}\phi^{+}_{\tau}(x)\to \begin{cases}1, & v_1(t,x)-v_2(t,x)>0\\[2mm] 0, & v_1(t,x)-v_2(t,x)\leq 0\end{cases}= \chi_{v_1(t,\cdot)>v_2(t,\cdot)}(x)$ pointwise a.e. on $\Omega$. Since $0\leq \phi_{\tau}^+\leq 1\in L^{2}(\Omega)$ for any $\tau>0$ we get from the \textit{Lebesgue Dominated Convergence Theorem} that:

\begin{equation}
\lim\limits_{\tau\to 0^+} \phi^{+}_{\tau}=\chi_{v_1(t,\cdot)>v_2(t,\cdot)}\ \text{in}\ L^2(\Omega).
\end{equation}

\noindent Taking into account that strong convergence in $L^2(\Omega)$ implies weak convergence in $L^2(\Omega)$, we can make $\tau\to 0^+$ in \eqref{inegalitatea2} to obtain, for a.e. $t\in (0,T)$, that: 

\begin{equation}\label{ecuatieinterm}
	\int_{\Omega} \dfrac{\partial \big [\overline{b}(x,v_2)-\overline{b}(x,v_1)\big ]}{\partial t}\chi_{v_1(t,\cdot)>v_2(t,\cdot)}(x)\ dx 
	+\lambda\int_\Omega\big [\overline{b}(x,v_2)-\overline{b}(x,v_1)\big ]\chi_{v_1(t,\cdot)>v_2(t,\cdot)}(x)\ dx\geq 0.
\end{equation}

\noindent Let us define $w:[0,T]\times\Omega\to\mathbb{R}$ by $w(t,x)=\overline{b}(x,v_2(t,x))-\overline{b}(x,v_1(t,x))$. So $w=\mathcal{N}_{\overline{b}}\circ v_2-\mathcal{N}_{\overline{b}}\circ v_1$. Since $\mathcal{N}_{\overline{b}}:L^2(\Omega)\to L^2(\Omega)$ is a Lipschitz operator -- from Proposition \ref{propoverb}
\textbf{(3)} -- and $v_1,v_2\in H^1_{\text{loc}}\big ((0,T);L^2(\Omega)\big )$, we deduce, using \cite[Corollary 3.23, page 36]{kreuter2015sobolev}, that $\mathcal{N}_{\overline{b}}\circ v_1,\mathcal{N}_{\overline{b}}\circ v_2\in H^1_{\text{loc}}\big ((0,T);L^2(\Omega)\big )$. Therefore $w\in H^1_{\text{loc}}\big ((0,T);L^2(\Omega)\big )$. Using the \textbf{strict monotony} of $\overline{b}$\footnote{See Proposition \ref{propoverb} \textbf{(1)}.} it is sufficient to show that $w\geq 0$ a.e. on $(0,T)\times\Omega$. 

\noindent Now \eqref{ecuatieinterm} becomes:

\begin{equation}\label{ecuatieinterme}
	\int_{\Omega} \dfrac{\partial w}{\partial t}(t,x)\chi_{v_1(t,\cdot)>v_2(t,\cdot)}(x)\ dx 
	+\lambda\int_\Omega w(t,x) \chi_{v_1(t,\cdot)>v_2(t,\cdot)}(x)\ dx\geq 0.
\end{equation}

\noindent Since $w=w^+-w^-$, where:

\begin{equation}
	w^{+}(t,\cdot)=\begin{cases} w(t,\cdot ), & w(t,\cdot)\geq 0\\[2mm] 0, & w(t,\cdot)<0\end{cases}=\begin{cases} w(t,\cdot), & v_2(t,\cdot)\geq v_1(t,\cdot)\\[2mm] 0, & v_1(t,\cdot)>v_2(t,\cdot)\end{cases}=w(t,\cdot)\big (1-\chi_{v_1(t,\cdot)>v_2(t,\cdot)}\big ),
\end{equation}

\noindent and:

\begin{equation}
	w^{-}(t,\cdot)=\begin{cases} -w(t,\cdot ), & w(t,\cdot)< 0\\[2mm] 0, & w(t,\cdot)\geq 0\end{cases}=\begin{cases} -w(t,\cdot), & v_1(t,\cdot)>v_2(t,\cdot)\\[2mm] 0, & v_1(t,\cdot)\leq v_2(t,\cdot)\end{cases}=-w(t,\cdot)\chi_{v_1(t,\cdot)>v_2(t,\cdot)},
\end{equation}

\noindent  we may write, using Theorem \ref{athminmax}, that $\dfrac{\partial w}{\partial t}(t,\cdot)=\dfrac{\partial w^+}{\partial t}(t,\cdot)-\dfrac{\partial w^-}{\partial t}(t,\cdot)$, where $w^+,w^-\in H^1_{\text{loc}}\big ((0,T);L^2(\Omega)\big )$ and:

\begin{align*}
	&\dfrac{\partial w^+}{\partial t}(t,\cdot)=\begin{cases}\dfrac{\partial w}{\partial t}(t,\cdot),\ \text{a.e. on}\ \{x\in\Omega\ |\ w(t,x)\geq 0\}=\{x\in\Omega\ |\ v_2(t,x)\geq v_1(t,x)\}\\[3mm] 0,\ \text{a.e. on}\ \{x\in\Omega\ |\ w(t,x)<0\}=\{x\in\Omega\ |\ v_2(t,x)< v_1(t,x)\}\end{cases} \\
	&\dfrac{\partial w^-}{\partial t}(t,\cdot)=\begin{cases}-\dfrac{\partial w}{\partial t}(t,\cdot),\ \text{a.e. on}\ \{x\in\Omega\ |\ w(t,x)<0\}=\{x\in\Omega\ |\ v_2(t,x)<v_1(t,x)\}\\[3mm] 0,\ \text{a.e. on}\ \{x\in\Omega\ |\ w(t,x)\geq 0\}=\{x\in\Omega\ |\ v_2(t,x)\geq v_1(t,x)\}\end{cases}.
\end{align*}

\noindent Therefore: $\dfrac{\partial w^+}{\partial t}(t,\cdot)=\dfrac{\partial w}{\partial t}(t,\cdot)\big (1-\chi_{v_1(t,\cdot)>v_2(t,\cdot)}\big ),\ \dfrac{\partial w^-}{\partial t}(t,\cdot)=-\dfrac{\partial w}{\partial t}(t,\cdot)\chi_{v_1(t,\cdot)>v_2(t,\cdot)}$ and \eqref{ecuatieinterme} gives us that:

\begin{align}
	&\int_{\Omega} \underbrace{\dfrac{\partial w^+(t,x)}{\partial t}\chi_{v_1>v_2}}_{=0}-\dfrac{\partial w^-(t,x)}{\partial t}\chi_{v_1>v_2} dx 
	+\lambda\int_\Omega \underbrace{w^+(t,x)\chi_{v_1>v_2}}_{=0}-w^-(t,x)\chi_{v_1>v_2} dx\geq 0 \nonumber\\
\Longrightarrow\	&\int_{\Omega} \dfrac{\partial w^-(t,x)}{\partial t} dx 
	+\lambda\int_\Omega w^-(t,x) dx\leq 0,\ \text{for a.e.}\ t\in (0,T).\label{eqintegrate}
\end{align}

\noindent We will show now that $w^{-}\in C\big ([0,T];L^2(\Omega)\big )$. Since $v_1,v_2\in C\big ([0,T];L^2(\Omega)\big )$ and $\mathcal{N}_{\overline{b}}:L^2(\Omega)\to L^2(\Omega)$ is a Lipschitz operator (hence continuous) we get that $w=\mathcal{N}_{\overline{b}}\circ v_2-\mathcal{N}_{\overline{b}}\circ v_1\in C\big ([0,T];L^2(\Omega)\big )$. Now remark that $w^{-}=-\min\{w,0\}$ is the composition between the functions $m:L^2(\Omega)\to L^2(\Omega),\ m(V)=-\min\{V,0\}$ and $w$. From the elementary inequality $|\min(s_1,0)-\min(s_2,0)|\leq |s_1-s_2|,\ \forall\ s_1,s_2\in\mathbb{R}$ we get that:

\begin{align*}
	\Vert m(V_1)-m(V_2)\Vert_{L^2(\Omega)}&=\left (\int_{\Omega} |m(V_1)-m(V_2)|^2\ dx\right )^{\frac{1}{2}}\\
	&=\left (\int_{\Omega} |\min\{V_1(x),0\}-\min\{V_2(x),0\}|^2\ dx\right )^{\frac{1}{2}}\\
	&\leq \left (\int_{\Omega} |V_1(x)-V_2(x)|^2\ dx\right )^{\frac{1}{2}}=\Vert V_1-V_2\Vert_{L^2(\Omega)},\ \forall\ V_1,V_2\in L^2(\Omega).
\end{align*}

\noindent Thus $m$ is a $1$ -- Lipschitz operator and hence continuous. As a consequence $w^{-}=m\circ w\in C\big ([0,T];L^2(\Omega)\big )$.

\noindent We consider the function $h:[0,T]\to [0,\infty),\ h(t)=\displaystyle\int_{\Omega} w^{-}(t,x)\ dx$.

\noindent It is easy to see that $h(0)=0$. Indeed, $v_2(0,x)\geq v_1(0,x)$ a.e. on $\Omega$ which means that $w(0,x)=\overline{b}(x,v_2(0,x))-\overline{b}(x,v_1(0,x))\geq 0$, from Proposition \ref{propoverb} \textbf{(1)}. Then $w^{-}(0,x)=0$ a.e. on $\Omega\ \Rightarrow\ h(0)=\displaystyle\int_{\Omega} w^{-}(0,x)\ dx=0$.

\noindent From the fact that $w^{-}\in H^1_{\text{loc}}\big ((0,T);L^2(\Omega)\big )$ we get that $w^-\in H^1\big ((a,b);L^2(\Omega)\big )$, for any $0<a<b<T$. Define $K:L^2(\Omega)\to \mathbb{R},\ K(V)=\displaystyle\int_{\Omega} V(x)\ dx$. Since $\Omega$ is bounded we have that $K$ is a bounded linear functional\footnote{From \textit{Cauchy inequality} we get that $|K(V_1)-K(V_2)|\leq \sqrt{|\Omega|}\cdot \Vert V_1-V_2\Vert_{L^2(\Omega)}$ for any $V_1,V_2\in L^2(\Omega)$.} and $h(t)=Kw^{-}(t,\cdot)$ for any $t\in [a,b]$. Notice that, $K:L^2(\Omega)\to\mathbb{R}$ and $w^{-}:[0,T]\to L^2(\Omega)$ are both continuous. Therefore $h=K\circ w^{-}:[0,T]\to\mathbb{R}$ is also continuous, i.e. $h\in C\big ([0,T]\big )$. Using now Theorem \ref{thmmaxpri} from the Appendix we get that $h\in H^1((a,b))$ and $h'(t)=\displaystyle\int_{\Omega}\dfrac{\partial w^{-}}{\partial t}(t,x)\ dx$ for a.e. $t\in (a,b)$. From $h\in H^1((a,b))$ and $h\in C\big ([a,b]\big )$ we deduce from \cite[Theorem 7.16, page 189]{leonibook} that $h\in\text{AC}\big ([a,b]\big )$.

\noindent Consider now $\tilde{h}:[0,T]\to [0,\infty),\ h(t)=e^{\lambda t}h(t)$. We immediately get that $\tilde{h}\in\text{AC}\big ([a,b]\big )\cap C\big ([0,T] \big )$, and from the \textit{product rule}\footnote{Absolutely continuous real functions are a.e. differentiable. See \cite[Theorem 12, page 107]{royden2022real}.} we get that for a.e. $t\in (0,T)$:

\begin{equation}
	\tilde{h}'(t)=e^{\lambda t}\big (h'(t)+\lambda h(t)\big )=e^{\lambda t}\left (\displaystyle\int_{\Omega}\dfrac{\partial w^{-}}{\partial t}(t,x)\ dx+\lambda\displaystyle\int_{\Omega} w^{-}(t,x)\ dx \right )\stackrel{\eqref{eqintegrate}}\leq 0.
\end{equation}

\noindent Furtermore, from \cite[Theorem 12, page 107]{royden2022real} we get that:

\begin{equation}\label{altaecint}
	\tilde{h}(b)-\tilde{h}(a)=\int_{a}^{b} \tilde{h}'(t)\ dt\leq 0.
\end{equation}

\noindent Therefore $\tilde{h}$ is decreasing on $(0,T)$ and continuous on $[0,T]$. We deduce
 that for any $t_0\in (0,T)$ we have that: $0\leq \tilde{h}(t_0)\leq \lim\limits_{t\to 0^+} \tilde{h}(t)=\tilde{h}(0)=h(0)=0$. So $0\leq \tilde{h}(t)\leq 0,\ \forall\ t\in (0,T)$, from where, by continuity: $\tilde{h}(t)=0$ for any $t\in [0,T]$, i.e. $h(t)=0$ for any $t\in [0,T]$. 

\noindent In conclusion for any $t\in [0,T]$ we have that $w^-(t,x)=0$ for a.e. $x\in\Omega$. This proves that $w(t,x)=w^+(t,x)\geq 0$ i.e. $v_2(t,x)\geq v_1(t,x)$, as needed.
\end{proof}

\begin{remark}\label{remremrem} We point out the following fact from the proof: If $v_1,v_2:(0,T)\to W^{1,p(x)}(\Omega)$ with $v_1,v_2\in C\big ([0,T];L^2(\Omega)\big )\cap  H^1_{\textnormal{loc}}((0,T);L^2(\Omega)\big )$, $v_2(0,x)\geq v_1(0,x)$ for a.e. $x\in\Omega$ and we manage somehow to prove that there is some $\lambda\in\mathbb{R}$ such that for a.e. $t\in (0,T)$ the following inequality holds 
	
	\begin{equation}
		\int_{\Omega} \dfrac{\partial \big [\overline{b}(x,v_2)-\overline{b}(x,v_1)\big ]}{\partial t}\chi_{v_1(t,\cdot)>v_2(t,\cdot)}(x)\ dx 
		+\lambda\int_\Omega\big [\overline{b}(x,v_2)-\overline{b}(x,v_1)\big ]\chi_{v_1(t,\cdot)>v_2(t,\cdot)}(x)\ dx\geq 0,
	\end{equation}
	
\noindent then we may conclude that for any $t\in [0,T]$ we have $v_2(t,x)\geq v_1(t,x)$ for a.e. $x\in\Omega$.
\end{remark}

\section{The auxiliary problem}\label{s6}

\noindent For every $\lambda>0$ and any $g\in\mathfrak{M}_{\lambda}^{[\varepsilon,\delta]}$, consider the following \textbf{auxiliary problem}:

\begin{equation}\label{eqdpgaux}\tag{A}
	\begin{cases}\dfrac{\partial \overline{b}(x,v(t,x))}{\partial t}-\operatorname{div}\mathbf{a}(x,\nabla v)+\lambda \overline{b}(x,v)=g(t,x), & (t,x)\in (0,T)\times\Omega\\[3mm] \mathbf{a}(x,\nabla v)\cdot\nu=0, & (t,x)\in (0,T)\times\partial\Omega\\[3mm] v(0,x)=u_0(x)\in [\varepsilon,\delta], & x\in\Omega\end{cases}
\end{equation}

\begin{definition}
		We say that $v\in C\big ([0,T]; L^2(\Omega)\big )\cap H^1_{\operatorname{loc}}\big ((0,T);L^2(\Omega)\big )\cap L^{1}\big (0,T;W^{1,p(x)}(\Omega)\big )$ is a \textbf{weak solution} of \eqref{eqdpgaux} if $v(0,x)=u_0(x)$ for a.e. $x\in\Omega$ and for a.e. $t\in (0,T)$ we have that:
	
	\begin{equation}
		\int_{\Omega} \dfrac{\partial \overline{b}(x,v(t,x))}{\partial t}\phi\ dx+\int_{\Omega} \mathbf{a}(x,\nabla v)\cdot\nabla\phi\ dx+\lambda\int_{\Omega} \overline{b}(x,v)\phi\ dx=\int_{\Omega} g(t,x)\phi\ dx,\ \forall\ \phi\in W^{1,p(x)}(\Omega).
	\end{equation}
\end{definition}

\subsection{Basic properties regarding the auxiliary problem}

\begin{proposition}\label{propauxepsdelta}
	If the auxiliary problem \eqref{eqdpgaux} admits a weak solution $v$, then $v\in\mathcal{V}_{[\varepsilon,\delta]}$.
\end{proposition}

\begin{proof} Just remark that for a.e. $t\in (0,T)$ and any $\phi\in W^{1,p(x)}(\Omega)^+$ we may write:
	
\begin{align*}
	&\int_{\Omega} \lambda b(x,\varepsilon)\phi\ dx=\int_{\Omega} \lambda \overline{b}(x,\varepsilon)\phi\ dx=\int_\Omega\underbrace{\dfrac{\partial \overline{b}(x,\varepsilon)}{\partial t}}_{=0}\phi\ dx+\int_{\Omega}\underbrace{\mathbf{a}(x,\nabla \varepsilon)}_{=\mathbf{0}}\cdot\nabla\phi\ dx+\lambda\int_{\Omega} \overline{b}(x,\varepsilon)\phi\ dx\\
	\leq &\int_{\Omega} g(t,x)\phi\ dx= \int_\Omega\dfrac{\partial \overline{b}(x,v)}{\partial t}\phi\ dx+\int_{\Omega}\mathbf{a}(x,\nabla v)\cdot\nabla\phi\ dx+\lambda\int_{\Omega} \overline{b}(x,v)\phi\ dx\\
	\leq &\int_{\Omega} \lambda b(x,\delta)\phi\ dx=\int_{\Omega} \lambda \overline{b}(x,\delta)\phi\ dx=\int_\Omega\underbrace{\dfrac{\partial \overline{b}(x,\delta)}{\partial t}}_{=0}\phi\ dx+\int_{\Omega}\underbrace{\mathbf{a}(x,\nabla \delta)}_{=\mathbf{0}}\cdot\nabla\phi\ dx+\lambda\int_{\Omega} \overline{b}(x,\delta)\phi\ dx.
\end{align*}
	
\noindent Applying now the \textit{Weak parabolic comparison principle} we deduce that $0\leq\varepsilon\leq v\leq \delta$ a.e. on $(0,T)\times\Omega$.
\end{proof}

\begin{theorem}\label{thmveryunique}
	The auxiliary problem \eqref{eqdpgaux} has \textbf{at most} one weak solution.
\end{theorem}

\begin{proof} Let $v_1,v_2$ be two weak solutions of \eqref{eqdpgaux}. Therefore we have for a.e. $t\in (0,T)$ and any $\phi\in W^{1,p(x)}(\Omega)^+$ that:
	
	\begin{align*}
		&\int_\Omega\dfrac{\partial \overline{b}(x,v_2)}{\partial t}\phi\ dx+\int_{\Omega}\mathbf{a}(x,\nabla v_2)\cdot\nabla\phi\ dx+\lambda\int_{\Omega} \overline{b}(x,v_2)\phi\ dx=\int_{\Omega} g(t,x)\phi\ dx\\
	\begin{cases}\geq \\ \leq\end{cases} &\int_{\Omega} g(t,x)\phi\ dx= \int_\Omega\dfrac{\partial \overline{b}(x,v_1)}{\partial t}\phi\ dx+\int_{\Omega}\mathbf{a}(x,\nabla v_1)\cdot\nabla\phi\ dx+\lambda\int_{\Omega} \overline{b}(x,v_1)\phi\ dx.
	\end{align*}
	
\noindent Thus, from the \textit{Weak parabolic comparison principle} we can conclude that $v_1\begin{cases} \geq\\ \leq\end{cases} v_2$ a.e. on $(0,T)\times\Omega$. In conclusion $v_1=v_2$ a.e. on $(0,T)\times\Omega$, i.e. $v_1\equiv v_2$. The uniqueness follows.
	
\end{proof}

\subsection{Constructing a weak solution to the auxiliary problem with regular initial data}

\noindent For every natural number $m\in\mathbb{N}$ we make the following construction of a \textit{time-discretization scheme}\footnote{This is called Rothe's method.}:
	
	\begin{itemize}
		\item \textbf{Step I:} Set $\Delta_m=\dfrac{T}{2^m}$ and $t_{m,n}=n\Delta_m$, for each $n\in\{0,1,\hdots,2^m\}$.
		
		\item \textbf{Step II:} Fix $v_{m,0}:=u_0\in\mathcal{U}_{[\varepsilon,\delta]}$.
		
		\item \textbf{Step III:} For each $n\in\{1,\hdots,2^m\}$ consider the following \textbf{doubly-nonlinear elliptic problem}, by which we define $v_{m,n}$ in terms of $v_{m,n-1}$:
	\end{itemize}
		
		\begin{equation}\label{eq1112}
			\begin{cases} \dfrac{\overline{b}(x,v_{m,n})-\overline{b}(x,v_{m,n-1})}{\Delta_m}-\operatorname{div} \mathbf{a}(x,\nabla v_{m,n})+\lambda \overline{b}(x,v_{m,n})=\underbrace{\dfrac{1}{\Delta_m}\displaystyle\int_{t_{m,n-1}}^{t_{m,n}} g(t,x)\ dt}_{:=g_{m,n}(x)}, & x\in \Omega \\[3mm] \mathbf{a}(x,\nabla v_{m,n})\cdot\nu=0, & x\in\partial\Omega\end{cases}
		\end{equation}
		
	\noindent\underline{\textbf{Claim 1:}} \textbf{For any $m\in\mathbb{N}$ and each $n\in\{1,2,\hdots,2^m\}$ we have that $g_{m,n}\in\mathcal{M}^{[\varepsilon,\delta]}_{\lambda}$.}
	
	\begin{proof} \noindent Indeed, since $g\in\mathfrak{M}^{[\varepsilon,\delta]}_{\lambda}$, we deduce that:
	
	\begin{equation*}
	\begin{cases} g_{m,n}\geq \dfrac{1}{\Delta_m}\displaystyle\int_{t_{m,n-1}}^{t_{m,n}}\lambda b(x,\varepsilon)\ dt= \lambda b(x,\varepsilon) \\[3mm] g_{m,n}\leq \dfrac{1}{\Delta_m}\displaystyle\int_{t_{m,n-1}}^{t_{m,n}}\lambda b(x,\delta)\ dt=\lambda b(x,\delta) \end{cases}\ \Longrightarrow\ g_{m,n}\in\mathcal{M}^{[\varepsilon,\delta]}_{\lambda}.
\end{equation*} 

\end{proof}
	
	\noindent We may rewrite \eqref{eq1112} as:
	
		\begin{equation}\label{eq1113}
		\begin{cases} -\operatorname{div} \mathbf{a}(x,\nabla v_{m,n})+\left (\lambda+\dfrac{1}{\Delta_m}\right ) \overline{b}(x,v_{m,n})=\underbrace{\dfrac{1}{\Delta_m}\displaystyle\int_{t_{m,n-1}}^{t_{m,n}} g(t,x)\ dt+\dfrac{1}{\Delta_m}\overline{b}(x,v_{m,n-1})}_{:=\tilde{g}_{m,n}(x)}, & x\in \Omega \\[3mm] \mathbf{a}(x,\nabla v_{m,n})\cdot\nu=0, & x\in\partial\Omega\end{cases}
	\end{equation}
	
	\noindent In order to study \eqref{eq1113} we need to recall the following result:
	
	\begin{proposition}\label{prop5551}
		For any $\lambda>0$ and any $G\in\mathcal{M}_{\lambda}^{[\varepsilon,\delta]}$ the problem:
		
		\begin{equation}
	\begin{cases}-\operatorname{div} \mathbf{a}\big (x,\nabla V(x)\big )+\lambda \overline{b}\big (x,V(x)\big )=G(x), & x\in\Omega \\[3mm] \mathbf{a}\big (x,\nabla V(x)\big )\cdot \nu=0, & x\in\partial\Omega  \end{cases},
		\end{equation}
		
		\noindent admits a unique weak solution $V\in W^{1,p(x)}(\Omega)$, i.e. for all $\phi\in W^{1,p(x)}(\Omega)$:
		
		\begin{equation}
			\int_{\Omega} \mathbf{a}\big (x,\nabla V(x)\big )\cdot\nabla\phi(x)\ dx+\lambda\int_{\Omega} \overline{b}\big (x,V(x)\big )\phi(x)\ dx=\int_{\Omega} G(x)\phi(x)\ dx.
		\end{equation}
		
		\noindent In addition, $V\in\mathcal{U}_{[\varepsilon,\delta]}$.\footnote{This is precisely Theorem 6.1 from \cite{max4} applied for $\hat{b}:\overline{\Omega}\times [0,\delta]\to\mathbb{R},\ \hat{b}(x,s)=\begin{cases} b(x,s), & s\in [\varepsilon,\delta]\\[2mm] b(x,\varepsilon)+s-\varepsilon, & s\in [0,\varepsilon) \end{cases}$, $\hat{f}:\overline{\Omega}\times [0,\delta]\to\mathbb{R},\ \hat{f}(x,s)=\begin{cases} f(x,s), & s\in [\varepsilon,\delta]\\[2mm] f(x,\varepsilon), & s\in [0,\varepsilon) \end{cases}$ and $G\in\mathcal{M}^{[\varepsilon,\delta]}_{\lambda}\subset \mathcal{M}^{[0,\delta]}_{\lambda}.$}
		
	\end{proposition}

	\bigskip

	\noindent \textbf{\underline{Claim 2:}} \textbf{For any $m\in\mathbb{N}$ and any $n\in\{1,2,\hdots, 2^m\}$, we have that $v_{m,n}\in W^{1,p(x)}(\Omega)\cap\mathcal{U}_{[\varepsilon,\delta]}$ is uniquely-defined and $\tilde{g}_{m,n}\in\mathcal{M}^{[\varepsilon,\delta]}_{\lambda+\frac{1}{\Delta_m}}$.}
	
	\bigskip
	
	\begin{proof}\noindent Indeed, since $v_{m,0}=u_0\in\mathcal{U}_{[\varepsilon,\delta]}$ it follows that:

	\begin{equation*}
		 \begin{cases} \tilde{g}_{m,1}\geq \dfrac{1}{\Delta_m}\displaystyle\int_{t_{m,0}}^{t_{m,1}}\lambda b(x,\varepsilon)\ dt+\dfrac{1}{\Delta_m}b(x,u_0)\geq \left (\lambda +\dfrac{1}{\Delta_m}\right )b(x,\varepsilon) \\[3mm] \tilde{g}_{m,1}\leq \dfrac{1}{\Delta_m}\displaystyle\int_{t_{m,0}}^{t_{m,1}}\lambda b(x,\delta)\ dt+\dfrac{1}{\Delta_m}b(x,u_0)\leq\left (\lambda +\dfrac{1}{\Delta_m}\right )b(x,\delta) \end{cases}. 
	\end{equation*}

	\noindent This shows that $\tilde{g}_{m,1}\in\mathcal{M}^{[\varepsilon,\delta]}_{\lambda+\frac{1}{\Delta_m}}$. Using now Proposition \ref{prop5551} we obtain that $v_{m,1}\in W^{1,p(x)}(\Omega)\cap\mathcal{U}_{[\varepsilon,\delta]}$ is the unique weak solution of \eqref{eq1113} for $n=1$. Now we have that $v_{m,1}\in\mathcal{U}_{[\varepsilon,\delta]}$. Thus:
	
	\begin{equation}
		\begin{cases} \tilde{g}_{m,2}\geq \dfrac{1}{\Delta_m}\displaystyle\int_{t_{m,1}}^{t_{m,2}}\lambda b(x,\varepsilon)\ dt+\dfrac{1}{\Delta_m}b(x,v_{m,1})\geq \left (\lambda +\dfrac{1}{\Delta_m}\right )b(x,\varepsilon) \\[3mm] \tilde{g}_{m,2}\leq \dfrac{1}{\Delta_m}\displaystyle\int_{t_{m,1}}^{t_{m,2}}\lambda b(x,\delta)\ dt+\dfrac{1}{\Delta_m}b(x,v_{m,1})\leq\left (\lambda +\dfrac{1}{\Delta_m}\right )b(x,\delta) \end{cases}.
	\end{equation}
	
	\noindent Using again Proposition \ref{prop5551} we obtain that $v_{m,2}\in W^{1,p(x)}(\Omega)\cap\mathcal{U}_{[\varepsilon,\delta]}$ is the unique weak solution of \eqref{eq1113} for $n=2$. Repeating the same process for $2^m-2$ more times we reach the conclusion of our claim.
\end{proof}

\begin{itemize}
	\item \textbf{Step IV:} We define the functions $v_m$ and $\tilde{v}_m$ for any $t\in [0,T]$ and a.e. $x\in\Omega$ via the relations:
	
	\begin{equation}
	\begin{cases}	v_m(t,x)=\displaystyle\sum_{n=1}^{2^m} v_{m,n}(x)\chi_{(t_{m,n-1},t_{m,n}]}(t),\quad v_m(0,x)=\tilde{v}_m(0,x):=v_{m,0}(x)=u_0(x) \\[5mm] b\big (x,\tilde{v}_{m}(t,x)\big )=\displaystyle\sum_{n=1}^{2^m} \left [b\big (x,v_{m,n-1}(x)\big )+\dfrac{t-t_{m,n-1}}{\Delta_m}\left (b\big (x,v_{m,n}(x)\big )-b\big (x,v_{m,n-1}(x)\big )\right )\right ]\chi_{(t_{m,n-1},t_{m,n}]}(t)\end{cases}.
	\end{equation}
\end{itemize}

\bigskip

\bigskip

\begin{proposition}\label{propovm} For each $m\in\mathbb{N}$ the following assertions are true:
	
	\begin{enumerate}
		\item[\textbf{(1)}] $v_m\in L^{\infty}\big (0,T;W^{1,p(x)}(\Omega)\big )$ and 
		
		\begin{equation}
			\Vert v_m\Vert_{L^{\infty}(0,T;W^{1,p(x)}(\Omega))}\leq \max_{n\in\overline{1,2^m}} \Vert v_{m,n}\Vert_{W^{1,p(x)}(\Omega)}.
		\end{equation}

		\item[\textbf{(2)}] $v_m\in \mathcal{V}_{[\varepsilon,\delta]}\subset L^{\infty}\big (0,T;L^\infty(\Omega)\big )\subset L^{\infty}\big (0,T;L^2(\Omega)\big )$ and 
		
		\begin{equation}
			\Vert v_m\Vert_{L^{\infty}(0,T;L^\infty(\Omega))}\leq \delta.
		\end{equation}
		
		\noindent \textbf{In particular the sequence $(v_m)_{m\geq 0}$ is bounded in $L^{\infty}\big (0,T;L^\infty(\Omega)\big )$.}
		
		\item[\textbf{(3)}] For any $t\in (0,T]$ and a.e. $x\in\Omega$ we have that $b(x,v_m(t,x))=\displaystyle\sum_{n=1}^{2^m} b(x,v_{m,n}(x))\chi_{(t_{m,n-1},t_{m,n}]}(t)$ and for $t=0$ we have that $b(x,v_m(0,x))=b(x,u_0(x))$ for a.e. $x\in\Omega$.
		
		\item[\textbf{(4)}] $b(\cdot,v_m)\in L^{\infty}\big (0,T;W^{1,p(x)}(\Omega)\big )$ and 
		
		\begin{equation}
			\Vert b(\cdot,v_m)\Vert_{L^{\infty}(0,T;W^{1,p(x)}(\Omega))}\leq \max_{n\in\overline{1,2^m}} \Vert b(\cdot,v_{m,n})\Vert_{W^{1,p(x)}(\Omega)}\leq L_0 \max_{n\in\overline{1,2^m}} \Vert v_{m,n}\Vert_{W^{1,p(x)}(\Omega)}+L_0\Vert 1\Vert_{L^{p(x)}(\Omega)}.
		\end{equation}

		\item[\textbf{(5)}] $b(\cdot,v_m)\in L^{\infty}\big (0,T;L^\infty(\Omega)\big )\subset L^{\infty}\big (0,T;L^2(\Omega)\big )$ and 
		
		\begin{equation}
			\Vert b(\cdot,v_m)\Vert_{L^{\infty}(0,T;L^\infty(\Omega))}\leq \Vert b(x,\delta)\Vert_{L^{\infty}(\Omega)}\leq L_0\delta.
		\end{equation}
		
		\noindent \textbf{In particular the sequence $\big (b(\cdot,v_m)\big )_{m\geq 0}$ is bounded in $L^{\infty}\big (0,T;L^\infty(\Omega)\big )$.}
		
	\end{enumerate}

\end{proposition}

\begin{proof} \noindent\textbf{(1)} For any $t\in (0,T]$ we have that there is exactly one $k\in\{1,2,\hdots, 2^m\}$ such that $t\in (t_{m,k-1},t_{m,k}]$. Therefore $v_{m}(t,\cdot)=v_{m,k}\in W^{1,p(x)}(\Omega)$, from Claim 2 and $v_m:(0,T]\to W^{1,p(x)}(\Omega)$ is a Bochner measurable function, since it is a simple function. Moreover
	
	\[
	\Vert v_{m}(t,\cdot)\Vert_{W^{1,p(x)}(\Omega)}=\Vert v_{m,k}\Vert_{W^{1,p(x)}(\Omega)}\leq \max_{n\in\overline{1,2^m}} \Vert v_{m,n}\Vert_{W^{1,p(x)}(\Omega)}.
	\]

\noindent This shows that $v_m\in L^{\infty}\big (0,T;W^{1,p(x)}(\Omega)\big )$.	
	
\medskip
	
\noindent\textbf{(2)} For any $t\in (0,T]$ we have from Claim 2 that 

\[
\varepsilon\leq \min_{n\in\overline{1,2^m}} v_{m,n}(x)\leq v_m(t,x)\leq \max_{n\in\overline{1,2^m}} v_{m,n}(x)\leq \delta.
\] 

\noindent This shows that $v_m\in \mathcal{V}_{[\varepsilon,\delta]}$.

\medskip

\noindent\textbf{(3)} For $t=0$ the equality is obvious. Fix any $k\in\{1,2,\hdots,2^m\}$ and any $t\in (t_{m,k-1},t_{m,k}]$. Thence

\[
b(x,v_{m}(t,x))=b(x,v_{m,k}(x))=b(x,v_{m,k}(x))\chi_{(t_{m,k-1},t_{m,k}]}(t)=\displaystyle\sum_{n=1}^{2^m} b(x,v_{m,n}(x))\chi_{(t_{m,n-1},t_{m,n}]}(t).
\]

\medskip

\noindent\textbf{(4)} For any $t\in (0,T]$ we have that there is exactly one $k\in\{1,2,\hdots, 2^m\}$ such that $t\in (t_{m,k-1},t_{m,k}]$. Therefore $v_{m}(t,\cdot)=v_{m,k}\in W^{1,p(x)}(\Omega)\cap\mathcal{U}_{[\varepsilon,\delta]}$ from Claim 2. Applying now Proposition \ref{propob} \textbf{(1)} we deduce that $b(\cdot,v_m(t,\cdot))=b(\cdot,v_{m,k})\in W^{1,p(x)}(\Omega)$. Also note that $b(\cdot,v_m):(0,T]\to W^{1,p(x)}(\Omega)$ is a Bochner measurable function, since it is a simple function. Moreover

\begin{align*}
\Vert b(\cdot,v_{m}(t,\cdot))\Vert_{W^{1,p(x)}(\Omega)}&=\Vert b(\cdot,v_{m,k})\Vert_{W^{1,p(x)}(\Omega)}\leq \max_{n\in\overline{1,2^m}} \Vert b(\cdot,v_{m,n})\Vert_{W^{1,p(x)}(\Omega)}\\
\text{Prop. \ref{propob} \textbf{(1)}}\ \ \ &\leq L_0 \max_{n\in\overline{1,2^m}} \Vert v_{m,n}\Vert_{W^{1,p(x)}(\Omega)}+L_0\Vert 1\Vert_{L^{p(x)}(\Omega)}.
\end{align*}

\noindent This shows that $b(\cdot,v_m)\in L^{\infty}\big (0,T;W^{1,p(x)}(\Omega)\big )$.

\medskip

\noindent\textbf{(5)} Note that, since $v_m\in\mathcal{V}_{[\varepsilon,\delta]}$, we can say that

\begin{equation}
	0\leq b(x,\varepsilon)\leq b(x,v_m(t,x))\leq b(x,\delta)\stackrel{\textbf{(H10)}}=b(x,\delta)-b(x,0)\stackrel{\text{Prop. \ref{propoext} \textbf{(3)}}}{\leq} L_0\delta_0.
\end{equation}

\noindent Taking into account that $b(\cdot,v_m)$ is also a measurable function (being a composition between measurable functions) we conclude that $b(\cdot,v_m)\in L^{\infty}\big (0,T;L^\infty(\Omega)\big )$ and $\Vert b(\cdot,v_m)\Vert_{L^{\infty}(0,T;L^\infty(\Omega))}\leq \Vert b(x,\delta)\Vert_{L^{\infty}(\Omega)}\leq L_0\delta$.
\end{proof}

\begin{proposition}\label{propobeta} For each $m\in\mathbb{N}$, let us define $\beta_m:[0,T]\times\Omega\to [0,\infty)$ by
	
	\begin{equation}
		\beta_m(t,x):=\displaystyle\sum_{n=1}^{2^m} \left [b\big (x,v_{m,n-1}(x)\big )+\dfrac{t-t_{m,n-1}}{\Delta_m}\left (b\big (x,v_{m,n}(x)\big )-b\big (x,v_{m,n-1}(x)\big )\right )\right ]\chi_{(t_{m,n-1},t_{m,n}]}(t),
	\end{equation}
	
	\noindent if $t\in (0,T]$ and for $t=0$ we set $\beta_m(0,x)=b(x,u_0(x))$ for a.e. $x\in\Omega$. 
	
	\medskip
	
	\noindent The following properties of $\beta_m$ hold:
	
	\begin{enumerate}
		\item[\textbf{(1)}] For any $t\in [0,T]$ and for a.e. $x\in \Omega$: $0\leq b(x,\varepsilon)\leq \beta_m(t,x)\leq b(x,\delta)\leq L_0\delta$. In particular $\beta_m\in L^{\infty}\big ((0,T)\times\Omega)$ and $\Vert\beta_m\Vert_{L^{\infty}((0,T)\times\Omega)}\leq L_0\delta$. Moreover we can say that $(x,\beta_m(t,x))\in \widehat{D}$ for any $t\in [0,T]$ and a.e. $x\in\Omega$.
		
		\item[\textbf{(2)}] $\beta_m\in C\big ([\Delta_m,T]; W^{1,p(x)}(\Omega)\big )$ and if $u_0\in W^{1,p(x)}(\Omega)$ then $\beta_m\in C\big ([0,T];W^{1,p(x)}(\Omega)\big )$.
		
		\item[\textbf{(3)}] $\beta_m:[0,T]\to L^{\infty}(\Omega)$ and $\beta_m\in C\big ([0,T];L^{\infty}(\Omega) \big )\subset C\big ([0,T];L^2(\Omega) \big )$.
		
		\item[\textbf{(4)}] $\beta_m\in L^{\infty}(\Delta_m,T; W^{1,p(x)}(\Omega))$ and $\beta_m\in L^{\infty}(0,T; W^{1,p(x)}(\Omega))$ if $u_0\in W^{1,p(x)}(\Omega)$. Moreover, in that case:
		
		\begin{equation}
			\Vert \beta_m\Vert_{L^{\infty}(0,T;W^{1,p(x)}(\Omega))}\leq \max_{n\in\overline{0,2^m}} \Vert b(\cdot,v_{m,n})\Vert_{W^{1,p(x)}(\Omega)}\leq L_0 \max_{n\in\overline{0,2^m}} \Vert v_{m,n}\Vert_{W^{1,p(x)}(\Omega)}+L_0\Vert 1\Vert_{L^{p(x)}(\Omega)}.
		\end{equation}
		
		\item[\textbf{(5)}] $\beta_m\in W^{1,\infty}\big ((0,T);L^{\infty}(\Omega)\big )\subset H^1((0,T);L^{\infty}(\Omega))\subset H^1((0,T);L^{2}(\Omega)) $ and for a.e. $t\in (0,T)$, and a.e. $x\in\Omega$ we have
		
		\begin{equation}
			\dfrac{\partial\beta_m}{\partial t}(t,x)=\displaystyle\sum_{n=1}^{2^m} \dfrac{b\big (x,v_{m,n}(x)\big )-b\big (x,v_{m,n-1}(x)\big )}{\Delta_m}\chi_{(t_{m,n-1},t_{m,n}]}(t).
		\end{equation}
		
		\noindent Moreover $\dfrac{\partial\beta_m}{\partial t}\in L^{\infty}((0,T)\times\Omega)$ and: $\left\Vert\dfrac{\partial\beta_m}{\partial t}\right\Vert_{L^{\infty}((0,T)\times\Omega)}\leq\dfrac{L_0|\delta-\varepsilon|}{\Delta_m}$.
		
	\end{enumerate}
	
\end{proposition}

\begin{proof} \noindent\textbf{(1)} For any $t\in [0,T]$ and a.e. $x\in\Omega$ we have, from the monotony of $b$ (i.e. \textbf{(H13)}) that:
	
	\begin{equation*}
		\beta_m(t,x)\leq \sum_{n=1}^{2^m} \max\{b(x,v_{m,n-1}(x)),b(x,v_{m,n}(x))\}\chi_{(t_{m,n-1},t_{m,n}]}(t)\leq \sum_{n=1}^{2^m} b(x,\delta)\chi_{(t_{m,n-1},t_{m,n}]}(t)=b(x,\delta),
	\end{equation*}
	
	\noindent and similarly:\footnote{$\beta_m$ is a convex combination between $b(\cdot,v_{m,n-1})$ and $b(\cdot,v_{m,n})$ on the interval $(t_{m,n-1},t_{m,n}]$.}
	
	\begin{equation*}
		\beta_m(t,x)\geq \sum_{n=1}^{2^m} \min\{b(x,v_{m,n-1}(x)),b(x,v_{m,n}(x))\}\chi_{(t_{m,n-1},t_{m,n}]}(t)\geq \sum_{n=1}^{2^m} b(x,\varepsilon)\chi_{(t_{m,n-1},t_{m,n}]}(t)=b(x,\varepsilon).
	\end{equation*}
	
	\noindent Note that, since $\varepsilon\geq 0$ we have that $b(x,\varepsilon)\geq 0$ and $b(x,\delta)\stackrel{\textbf{(H10)}}{=}b(x,\delta)-b(x,0)\stackrel{\text{Prop.\ref{propoext} \textbf{(3)}}}{\leq} L_0\delta$.
	
	\bigskip
	
	\noindent\textbf{(2)} For $t\geq\Delta_m=t_{m,1}$ we have that $\beta_m(t,\cdot)$ is a convex combination between $b(\cdot,v_{m,1}),\ b(\cdot,v_{m,2}),\ \hdots,$ $b(\cdot,v_{m,2^m})$. Since $v_{m,1},v_{m,2},\hdots,v_{m,2^m}\in W^{1,p(x)}(\Omega)\cap\mathcal{U}_{[\varepsilon,\delta]}$ we obtain from Proposition \ref{propob} \textbf{(1)} that $b(\cdot,v_{m,1}),\ b(\cdot,v_{m,2})$, $\ \hdots, b(\cdot,v_{m,2^m})\in W^{1,p(x)}(\Omega)$ and therefore $\beta_m(t,\cdot)\in W^{1,p(x)}(\Omega)$.
	
	\medskip
	
	\noindent Next, we will show that for any $t_0\in [\Delta_m,T]$: $\lim\limits_{t\to t_0} \beta_m(t,\cdot)=\beta_m(t_0,\cdot)$ in $W^{1,p(x)}(\Omega)$. If $t_0\notin\{\Delta_m=t_{m,1},\hdots,t_{m,2^{m}-1}\}$ then we can assume that there is exactly one $n\in \{2,\hdots, 2^m\}$ such that $t,t_0\in (t_{m,n-1},t_{m,n}]$. Therefore:
	
	\[
	\Vert \beta_m(t,\cdot)-\beta_m(t_0,\cdot)\Vert_{W^{1,p(x)}(\Omega)}=\dfrac{|t-t_0|}{\Delta_m}\Vert b(\cdot,v_{m,n})-b(\cdot,v_{m,n-1})\Vert_{W^{1,p(x)}(\Omega)}\stackrel{t\to t_0}{\longrightarrow} 0.
	\]
	
	\noindent If $t_0=t_{m,1}$ then $\beta_m(t_0,\cdot)=b(\cdot,v_{m,1})$, $t,t_0\in [t_{m,1},t_{m,2})$ and
	
	\[
	\Vert \beta_m(t,\cdot)-\beta_m(t_0,\cdot)\Vert_{W^{1,p(x)}(\Omega)}=\dfrac{|t-t_0|}{\Delta_m}\Vert b(\cdot,v_{m,2})-b(\cdot,v_{m,1})\Vert_{W^{1,p(x)}(\Omega)}\stackrel{t\to t_0}{\longrightarrow} 0.
	\]
	
	\noindent If $t_0=t_{m,n}$ for some $n\in\{2,\hdots,2^m-1\}$ we need to compute lateral limits:
	
	\begin{itemize}
		\item For $t\in (t_{m,n-1},t_{m,n}]$ we get that:
		
		\[
		\Vert \beta_m(t,\cdot)-\beta_m(t_0,\cdot)\Vert_{W^{1,p(x)}(\Omega)}=\dfrac{t_0-t}{\Delta_m}\Vert b(\cdot,v_{m,n})-b(\cdot,v_{m,n-1})\Vert_{W^{1,p(x)}(\Omega)}\stackrel{t\to t_0}{\longrightarrow} 0.
		\]
		
		\item For $t\in (t_{m,n},t_{m,n+1}]$, $\beta_m(t_0,\cdot)=b(\cdot,v_{m,n})$ and we get that:
		
		\[
		\Vert \beta_m(t,\cdot)-\beta_m(t_0,\cdot)\Vert_{W^{1,p(x)}(\Omega)}=\dfrac{t-t_0}{\Delta_m}\Vert b(\cdot,v_{m,n+1})-b(\cdot,v_{m,n})\Vert_{W^{1,p(x)}(\Omega)}\stackrel{t\to t_0}{\longrightarrow} 0.
		\]
		
	\end{itemize}
	
	\noindent Therefore $\beta_m\in C\big ([\Delta_m,T];W^{1,p(x)}(\Omega)\big )$.
	\medskip
	
	\noindent If $u_0=v_{m,0}\in W^{1,p(x)}(\Omega)$ then for any $t\in [0,T]$, $\beta_m(t,\cdot)$ is a convex combination between $b(\cdot,v_{m,0}),\ b(\cdot,v_{m,1}),\ b(\cdot,v_{m,2}),\ \hdots,$ $b(\cdot,v_{m,2^m})$. Since $v_{m,0},\ v_{m,1},v_{m,2},\hdots,v_{m,2^m}\in W^{1,p(x)}(\Omega)\cap\mathcal{U}_{[\varepsilon,\delta]}$ we obtain again from Proposition \ref{propob} \textbf{(1)} that $b(\cdot,v_{m,0}),\ b(\cdot,v_{m,1}),\ b(\cdot,v_{m,2})$, $\ \hdots, b(\cdot,v_{m,2^m})\in W^{1,p(x)}(\Omega)$ and therefore $\beta_m(t,\cdot)\in W^{1,p(x)}(\Omega)$. In the same way as above one can show that $\beta_m\in C\big ([0,T];W^{1,p(x)}(\Omega)\big )$.

	\noindent\textbf{(3)} The fact that for any $t\in [0,T]$: $\beta_m(t,\cdot)\in L^{\infty}(\Omega)$ follows from \textbf{(1)}. Note that for any $n\in\{1,2\hdots,2^m\}$ and any $t,t_0\in [t_{m,n-1},t_{m,n}]$
	
	\begin{equation}
		\Vert \beta_m(t,\cdot)-\beta_m(t_0,\cdot)\Vert_{L^{\infty}(\Omega)}=\dfrac{|t-t_0|}{\Delta_m}\Vert b(\cdot,v_{m,n})-b(\cdot,v_{m,n-1})\Vert_{L^{\infty}(\Omega)}\stackrel{t\to t_0}{\longrightarrow} 0.
	\end{equation} 
	
	\noindent So $\beta_m\in C\big ([0,T];L^{\infty}(\Omega)\big )$.
	
	\noindent\textbf{(4)} Note that for any $t\in [\Delta_m,T]$ we have that:
	
	\begin{align*}
		\Vert \beta_m(t,\cdot)\Vert_{W^{1,p(x)}(\Omega)}&\leq \sum_{n=2}^{2^m}\max\{\Vert b(\cdot,v_{m,n-1})\Vert_{W^{1,p(x)}(\Omega)}, \Vert b(\cdot,v_{m,n})\Vert_{W^{1,p(x)}(\Omega)} \}\chi_{(t_{m,n-1},t_{m,n}]}(t)\\
		&\leq \max_{n\in \overline{1,2^m}}\Vert b(\cdot,v_{m,n})\Vert_{W^{1,p(x)}(\Omega)}\\
		\text{Proposition \ref{propob} \textbf{(1)}}\ \ \ 		&\leq L_0\Vert 1\Vert_{L^{p(x)}(\Omega)}+L_0 \max_{n\in \overline{1,2^m}}\Vert v_{m,n}\Vert_{W^{1,p(x)}(\Omega)}.
	\end{align*}
	
	\noindent Therefore $\beta_m\in L^{\infty}\big (\Delta_m,T;W^{1,p(x)}(\Omega)\big )$. Similarly, if $u_0=v_{m,0}\in W^{1,p(x)}(\Omega)$ then for any $t\in [0,T]$
	
	\begin{align*}
		\Vert \beta_m(t,\cdot)\Vert_{W^{1,p(x)}(\Omega)}&\leq \sum_{n=1}^{2^m}\max\{\Vert b(\cdot,v_{m,n-1})\Vert_{W^{1,p(x)}(\Omega)}, \Vert b(\cdot,v_{m,n})\Vert_{W^{1,p(x)}(\Omega)} \}\chi_{(t_{m,n-1},t_{m,n}]}(t)\\
		&\leq \max_{n\in \overline{0,2^m}}\Vert b(\cdot,v_{m,n})\Vert_{W^{1,p(x)}(\Omega)}\\
\text{Proposition \ref{propob} \textbf{(1)}}\ \ \ 		&\leq L_0\Vert 1\Vert_{L^{p(x)}(\Omega)}+L_0 \max_{n\in \overline{0,2^m}}\Vert v_{m,n}\Vert_{W^{1,p(x)}(\Omega)}.
	\end{align*}
	
	\noindent So $\beta_m\in L^{\infty}\big (0,T;W^{1,p(x)}(\Omega)\big )$.

	\noindent\textbf{(5)} From \textbf{(1)} we already know that $\beta_m\in L^{\infty}\big ((0,T);L^{\infty}(\Omega)\big )\simeq L^{\infty}\big ((0,T)\times\Omega\big )$. 
	
	\noindent Now for any $\varphi\in C^{\infty}_{c}(\mathbb{R})$ with $\operatorname{supp}(\varphi)\subset (0,T)$ we have that
	
	\begin{align*}
		&\int_{0}^T\beta_m(t,\cdot)\varphi'(t)\ dt=\sum_{n=1}^{2^m} \int_{t_{m,n-1}}^{t_{m,n}} b(\cdot,v_{m,n-1})\varphi'(t)+\dfrac{b(\cdot,v_{m,n})-b(\cdot,v_{m,n-1})}{\Delta_m}(t-t_{m,n-1})\varphi'(t)\ dt \\
		=&\sum_{n=1}^{2^m}b(\cdot,v_{m,n-1})\int_{t_{m,n-1}}^{t_{m,n}}\varphi'(t)\ dt+\dfrac{b(\cdot,v_{m,n})-b(\cdot,v_{m,n-1})}{\Delta_m}\left (\int_{t_{m,n-1}}^{t_{m,n}}(t-t_{m,n-1})\varphi'(t)\ dt \right )\\[3mm]
		=&\sum_{n=1}^{2^m}b(\cdot,v_{m,n-1})\cdot \big (\varphi(t_{m,n})-\varphi(t_{m,n-1}) \big )+\dfrac{b(\cdot,v_{m,n})-b(\cdot,v_{m,n-1})}{\Delta_m}\bigg ((t-t_{m,n-1})\varphi(t)\bigg |_{t=t_{m,n-1}}^{t=t_{m,n}}\\
		&-\int_{t_{m,n-1}}^{t_{m,n}} \varphi(t)\ dt\bigg )\\
		=&-\int_{0}^T\left (\sum_{n=1}^{2^m}\dfrac{b(\cdot,v_{m,n})-b(\cdot,v_{m,n-1})}{\Delta_m}\chi_{(t_{m,n-1},t_{m,n}]}(t)\right )\varphi(t)\ dt\\
		&+\sum_{n=1}^{2^m}b(\cdot,v_{m,n-1})\cdot \big (\varphi(t_{m,n})-\varphi(t_{m,n-1}) \big )+\big (b(\cdot,v_{m,n})-b(\cdot,v_{m,n-1})\big )\varphi(t_{m,n})\\
		=&-\int_{0}^T\left (\sum_{n=1}^{2^m}\dfrac{b(\cdot,v_{m,n})-b(\cdot,v_{m,n-1})}{\Delta_m}\chi_{(t_{m,n-1},t_{m,n}]}(t)\right )\varphi(t)\ dt\\
		&+\sum_{n=1}^{2^m}b(\cdot,v_{m,n})\varphi(t_{m,n})-b(\cdot,v_{m,n-1})\varphi(t_{m,n-1})\\
		=&-\int_{0}^T\left (\sum_{n=1}^{2^m}\dfrac{b(\cdot,v_{m,n})-b(\cdot,v_{m,n-1})}{\Delta_m}\chi_{(t_{m,n-1},t_{m,n}]}(t)\right )\varphi(t)\ dt+b(\cdot,v_{m,2^m})\underbrace{\varphi(T)}_{=0}-b(\cdot,v_{m,0})\underbrace{\varphi(0)}_{=0}\\
		=&-\int_{0}^T\left (\sum_{n=1}^{2^m}\dfrac{b(\cdot,v_{m,n})-b(\cdot,v_{m,n-1})}{\Delta_m}\chi_{(t_{m,n-1},t_{m,n}]}(t)\right )\varphi(t)\ dt.
	\end{align*}
	
	\noindent Since for a.e. $t\in (0,T)$ and a.e. $x\in\Omega$
	
	\begin{align*}
	 &	\left|\sum_{n=1}^{2^m}\dfrac{b(x,v_{m,n}(x))-b(x,v_{m,n-1}(x))}{\Delta_m}\chi_{(t_{m,n-1},t_{m,n}]}(t)\right |\leq  \max_{n\in \overline{1,2^m}} \left|\dfrac{b(x,v_{m,n}(x))-b(x,v_{m,n-1}(x))}{\Delta_m}\right |\\
\stackrel{\text{Prop.\ref{propoext} \textbf{(3)}}}{\leq} & \max_{n\in \overline{1,2^m}} \dfrac{L_0}{\Delta_m}\left|v_{m,n}(x)-v_{m,n-1}(x)\right |\leq \dfrac{L_0(\delta-\varepsilon)}{\Delta_m}.
	\end{align*}
	
	 \noindent we get that $\displaystyle\sum_{n=1}^{2^m}\dfrac{b(\cdot,v_{m,n})-b(\cdot,v_{m,n-1})}{\Delta_m}\chi_{(t_{m,n-1},t_{m,n}]}(\cdot)\in L^{\infty}((0,T);L^{\infty}(\Omega))$ and the weak derivative of $\beta_m$ with respect to $t$ is
	 
	 \begin{equation}
	 	\dfrac{\partial\beta_m}{\partial t}(t,\cdot)=\sum_{n=1}^{2^m}\dfrac{b(\cdot,v_{m,n})-b(\cdot,v_{m,n-1})}{\Delta_m}\chi_{(t_{m,n-1},t_{m,n}]}(t),\ \text{for a.e.}\ t\in (0,T).
	 \end{equation}
	 
	 \noindent Thus $\beta_m\in W^{1,\infty}\big ((0,T);L^{\infty}(\Omega)\big )$. In particular, we have showed that $\left \Vert \dfrac{\partial \beta_m}{\partial t}\right \Vert_{L^{\infty}((0,T)\times\Omega)}\leq \dfrac{L_0(\delta-\varepsilon)}{\Delta_m}$.
\end{proof}

\begin{proposition}\label{propovtilde}The following assertions hold for any $m\in\mathbb{N}$:
	
	\begin{enumerate}
		
		\item[\textbf{(1)}] For any $t\in [0,T]$ and for a.e. $x\in\Omega$: $\tilde{v}_m(t,x)=\mathfrak{b}(x,\beta_m(t,x))$ and $\tilde{v}_m\in\mathcal{V}_{[\varepsilon,\delta]}$. In particular $\tilde{v}_m\in L^{\infty}\big ((0,T)\times\Omega\big )$ and $\Vert \tilde{v}_m\Vert_{L^{\infty}((0,T)\times\Omega)}\leq \delta$.
		
		\item[\textbf{(2)}] $\tilde{v}_m\in C\big ([\Delta_m,T];W^{1,p(x)}(\Omega) \big )$ and if $u_0\in W^{1,p(x)}(\Omega)$ then $\tilde{v}_m\in C\big ([0,T];W^{1,p(x)}(\Omega) \big )$. Moreover the following formula holds for any $t\in [\Delta_m,T]$ or any $t\in [0,T]$ -- if $u_0\in W^{1,p(x)}(\Omega)$ -- and for a.e. $x\in\Omega$:
		
		\begin{equation}
			\dfrac{\partial\tilde{v}_m}{\partial x_i}(t,x)=\dfrac{\partial\mathfrak{b}}{\partial x_i}(x,\beta_m(t,x))+\dfrac{\partial\mathfrak{b}}{\partial\sigma}(x,\beta_m(t,x))\dfrac{\partial\beta_m}{\partial x_i}(t,x).
		\end{equation}
		
		\item[\textbf{(3)}] $\tilde{v}_m:[0,T]\to L^{\infty}(\Omega)$ and $\tilde{v}_m\in C\big ([0,T];L^{\infty}(\Omega) \big )\subset C\big ([0,T];L^2(\Omega) \big )$.
		
		\item[\textbf{(4)}] $\tilde{v}_m\in L^{\infty}\big (\Delta_m,T;W^{1,p(x)}(\Omega)\big )$ and if $u_0\in W^{1,p(x)}(\Omega)$ then $\tilde{v}_m\in L^{\infty}(0,T;W^{1,p(x)}(\Omega))$. Moreover the following inequality holds
		
		\begin{equation}
			\Vert \tilde{v}_m\Vert_{L^{\infty}(0,T;W^{1,p(x)}(\Omega))}\leq \dfrac{1}{\ell_0}\Vert \beta_m\Vert_{L^{\infty}(0,T;W^{1,p(x)}(\Omega))}+\dfrac{L_0}{\ell_0}\Vert 1\Vert_{L^{p(x)}(\Omega)}+\dfrac{1}{\ell_0}\Vert b(\cdot,\varepsilon)-\ell_0\varepsilon\Vert_{L^{p(x)}(\Omega)}.
		\end{equation}

		\item[\textbf{(5)}] $\tilde{v}_m\in W^{1,\infty}\big ((0,T);L^{\infty}(\Omega)\big )\subset H^1\big ((0,T);L^{\infty}(\Omega)\big )\subset H^1\big ((0,T);L^2(\Omega) \big )$ and for a.e. $t\in (0,T)$ and a.e. $x\in\Omega$ we have that:
		
		\begin{equation}
			\dfrac{\partial\tilde{v}_m}{\partial t}(t,x)=\dfrac{\partial\mathfrak{b}}{\partial\sigma}(x,\beta_m(t,x))\cdot\dfrac{\partial\beta_m}{\partial t}(t,x).
		\end{equation}
		
		\noindent We have that $\dfrac{\partial\tilde{v}_m}{\partial t}\in L^{\infty}\big ((0,T)\times\Omega\big )$ and furthermore $\left\Vert\dfrac{\partial\tilde{v}_m}{\partial t} \right\Vert_{L^{\infty}((0,T)\times\Omega)}\leq\dfrac{L_0|\delta-\varepsilon|}{\ell_0\Delta_m}$.

	\end{enumerate}
	
\end{proposition}

\begin{proof} \textbf{(1)} From Proposition \ref{propobeta} \textbf{(1)} we know that for any $t\in [0,T]$ and a.e. $x\in \Omega$: $(x,\beta_m(t,x))\in \widehat{D}\subset\operatorname{cl}(D)$. But from Remark \ref{rembfraktilde} we have that $\mathfrak{b}:\widehat{D}\to[\varepsilon,\delta]$. Thus from the definition of $\tilde{v}_m$ we know that $b(x,\tilde{v}_m(t,x))=\beta_m(t,x)\in [b(x,\varepsilon),b(x,\delta)]$. Therefore, from the definition of $\mathfrak{b}$ we deduce that $\tilde{v}_m(t,x)=\mathfrak{b}(x,\beta_m(t,x))\in [\varepsilon,\delta]$. Being a composition of measurable functions we get that $\tilde{v}_m$ is measurable, and hence $\tilde{v}_m\in\mathcal{V}_{[\varepsilon,\delta]}$.
	
	\medskip
	
\noindent\textbf{(2)} We know from Proposition \ref{propobeta} \textbf{(2)} and \textbf{(1)} that for any $t\in [\Delta_m,T]$: $\beta_m(t,\cdot)\in W^{1,p(x)}(\Omega)$, and $b(x,\varepsilon)\leq\beta_m(t,\cdot)\leq b(x,\delta)$ a.e. on $\Omega$. Using now Proposition \ref{propobrussian} \textbf{(6)} we get that $\tilde{v}_m(t,\cdot)=\mathfrak{b}(\cdot,\beta_m(t,\cdot))\in W^{1,p(x)}(\Omega)\cap \mathcal{U}_{[\varepsilon,\delta]}$ and for each $i\in\overline{1,N}$ and a.e. $x\in\Omega$:

\begin{equation}
	\dfrac{\partial\tilde{v}_m}{\partial x_i}(t,x)=\dfrac{\partial\mathfrak{b}}{\partial x_i}(x,\beta_m(t,x))+\dfrac{\partial\mathfrak{b}}{\partial\sigma}(x,\beta_m(t,x))\dfrac{\partial\beta_m}{\partial x_i}(t,x).
\end{equation}

\noindent Consider now some arbitrary $t_0\in [\Delta_m,T]$ and any sequence $(t_n)_{n\geq 1}\subset [\Delta_m,T]$ with $t_n\to 0$. From Proposition \ref{propobeta} \textbf{(2)} we have that $\lim\limits_{n\to\infty} \beta_m(t_n,\cdot)=\beta_m(t_0,\cdot)$ in $W^{1,p(x)}(\Omega)$. Therefore, using Proposition \ref{propobrussian} \textbf{(7)}, we obtain that:

\begin{equation}
	\lim\limits_{n\to\infty} \tilde{v}_m(t_n,\cdot)=\lim\limits_{n\to\infty} \mathfrak{b}(\cdot,\beta_m(t_n,\cdot))=\mathfrak{b}(\cdot,\beta_m(t_0,\cdot))=\tilde{v}_m(t_0,\cdot)\ \text{in}\ W^{1,p(x)}(\Omega).
\end{equation}

\noindent This shows that $\tilde{v}_m\in C\big ([\Delta_m,T];W^{1,p(x)}(\Omega)\big )$. Similarly if $u_0\in W^{1,p(x)}(\Omega)$ then $\beta_m\in$ $C\big ([0,T];W^{1,p(x)}(\Omega)\big )$ from Proposition \ref{propobeta} \textbf{(2)} and therefore for any $t\in [0,T]$ we obtain via Proposition \ref{propobrussian} \textbf{(6)} that $\tilde{v}_m(t,\cdot)\in W^{1,p(x)}(\Omega)$ and from Proposition \ref{propobrussian} \textbf{(7)} we can conclude that $\tilde{v}_m\in C\big ([0,T];W^{1,p(x)}(\Omega)\big )$ as in the previous case.

\medskip

\noindent\textbf{(3)} From \textbf{(1)} we already know that for any $t\in [0,T]$: $\tilde{v}_m(t,\cdot)\in\mathcal{U}_{[\varepsilon,\delta]}\subset L^{\infty}(\Omega)$. Moreover, for any $t,t_0\in [0,T]$ and a.e. $x\in\Omega$, from Proposition \ref{propobrussian} \textbf{(5)}, we can write:

\begin{equation}
	|\tilde{v}_m(t,x)-\tilde{v}_m(t_0,x)|=|\mathfrak{b}(x,\beta_m(t,x))-\mathfrak{b}(x,\beta_m(t_0,x))|\leq \dfrac{1}{\ell_0}|\beta_m(t,x)-\beta_m(t_0,x)|.
\end{equation}

\noindent Therefore, for any $t,t_0\in [0,T]$

\begin{equation}
	\Vert \tilde{v}_m(t,\cdot)-\tilde{v}_m(t_0,\cdot)\Vert_{L^{\infty}(\Omega)}\leq\dfrac{1}{\ell_0}\Vert \beta_m(t,\cdot)-\beta_m(t_0,\cdot)\Vert_{L^{\infty}(\Omega)}\stackrel{t\to t_0}{\longrightarrow} 0,
\end{equation}

\noindent because from Proposition \ref{propobeta} \textbf{(3)} we know that $\beta_m\in C\big ([0,T];L^{\infty}(\Omega)\big )$. So we have proved that $\tilde{v}_m\in C\big ([0,T];L^{\infty}(\Omega)\big )$.

\medskip

\noindent\textbf{(4)} From \textbf{(2)} we know that for any $t\in [\Delta_m,T]$ we have that $\tilde{v}_m(t,\cdot)\in W^{1,p(x)}(\Omega)$. Using Proposition \ref{propobrussian} \textbf{(6)} we obtain that:

\begin{equation}
	\Vert \tilde{v}_m(t,\cdot)\Vert_{W^{1,p(x)}(\Omega)}\leq \dfrac{1}{\ell_0}\Vert \beta_m(t,\cdot)\Vert_{W^{1,p(x)}(\Omega)}+\dfrac{L_0}{\ell_0}\Vert 1\Vert_{L^{p(x)}(\Omega)}+\dfrac{1}{\ell_0}\Vert b(\cdot,\varepsilon)-\ell_0\varepsilon\Vert_{L^{p(x)}(\Omega)}.
\end{equation}

\noindent From Proposition \ref{propobeta} \textbf{(4)} we have that $\beta_m\in L^{\infty}\big (\Delta_m,T;W^{1,p(x)}(\Omega)\big )$, and from the above inequality we obtain that $\tilde{v}_m\in L^{\infty}\big (\Delta_m,T;W^{1,p(x)}(\Omega)\big )$. The proof is exactly the same when $u_0\in W^{1,p(x)}(\Omega)$, since in this case from Proposition \ref{propobeta} \textbf{(4)} we know that $\beta_m\in L^{\infty}\big (0,T;W^{1,p(x)}(\Omega)\big )$ and therefore: $\tilde{v}_m\in L^{\infty}\big (0,T;W^{1,p(x)}(\Omega)\big )$ from the above inequality which in this situation holds for any $t\in [0,T]$.

\medskip

\noindent\textbf{(5)} We want to apply here the \textit{Chain Rule for Bochner-Sobolev spaces} -- i.e. Theorem \ref{thmbochner} \textbf{(1)} from the Appendix. We have from Proposition \ref{propobeta} \textbf{(1)} that $\beta_m:(0,T)\times\Omega\to \displaystyle\bigcup_{x\in\Omega} [b(x,\varepsilon),b(x,\delta)]\subset [\underset{x\in\Omega}{\operatorname{ess\ inf}}\ b(x,\varepsilon),\underset{x\in\Omega}{\operatorname{ess\ sup}}\ b(x,\delta)]$ (a bounded interval, from \textbf{(H9)}). Note that from Proposition \ref{propobeta} \textbf{(5)} we have that $\beta_m\in W^{1,\infty}\big ((0,T);L^{\infty}(\Omega)\big )$. Also, from Remark \ref{rembrussian} we know that $\mathfrak{b}\in C^1(\mathbb{R}^{N+1})$. Since all the requirements of the \textit{Chain Rule for Bochner-Sobolev spaces} are fulfilled (see Remark \ref{rembochner}), we deduce that $\tilde{v}_m=\mathfrak{b}(\cdot,\beta_m)\in  W^{1,\infty}\big ((0,T);L^{\infty}(\Omega)\big )$ and for a.e. $(t,x)\in (0,T)\times\Omega$ we have that

\begin{equation}
	\dfrac{\partial\tilde{v}_m}{\partial t}(t,x)=\dfrac{\partial\mathfrak{b}}{\partial\sigma}(x,\beta_m(t,x))\cdot\dfrac{\partial\beta_m}{\partial t}(t,x).
\end{equation}

\noindent Moreover for a.e. $(t,x)\in (0,T)\times\Omega$ we have from Proposition \ref{propobrussian} \textbf{(5)} and Proposition \ref{propobeta} \textbf{(5)} that

\begin{align*}
	\left |	\dfrac{\partial\tilde{v}_m}{\partial t}(t,x) \right |=\left |\dfrac{\partial\mathfrak{b}}{\partial\sigma}\underbrace{(x,\beta_m(t,x))}_{\in\widehat{D}} \right |\cdot \left |\dfrac{\partial\beta_m}{\partial t}(t,x) \right |\leq \dfrac{1}{\ell_0}\cdot\dfrac{L_0|\delta-\varepsilon|}{\Delta_m}.
\end{align*}

\noindent From here we get that $\left\Vert\dfrac{\partial\tilde{v}_m}{\partial t} \right\Vert_{L^{\infty}((0,T)\times\Omega)}\leq\dfrac{L_0|\delta-\varepsilon|}{\ell_0\Delta_m}$ and the proof is now complete.

\end{proof}

\begin{proposition}\label{propovmn} If $u_0\in W^{1,p(x)}(\Omega)$ the following assertions are true:
	
	\begin{enumerate}
		\item[\textbf{(1)}]  For each $m\in\mathbb{N}$ we have that
		
		\begin{align}\label{est1}
			&\left (\lambda+\dfrac{1}{\Delta_m}\right )\sum_{n=1}^{2^m}\int_{\Omega} \big [b(x,v_{m,n})-b(x,v_{m,n-1}) \big ]\cdot (v_{m,n}-v_{m,n-1})\ dx \nonumber\\
			\leq&\mathcal{A}(u_0)+\lambda L_0(\delta-\varepsilon)\sum_{n=1}^{2^m} \int_{\Omega} |v_{m,n}-v_{m,n-1}|\ dx.
		\end{align}
		
		\item[\textbf{(2)}] For each $m\in\mathbb{N}$ we have that
		
		\begin{equation}\label{est2}
			\ell_0\left (\lambda+\dfrac{1}{\Delta_m}\right )\sum_{n=1}^{2^m}\int_{\Omega} (v_{m,n}-v_{m,n-1})^2\ dx \leq\mathcal{A}(u_0)+\lambda L_0(\delta-\varepsilon)\sum_{n=1}^{2^m} \int_{\Omega} |v_{m,n}-v_{m,n-1}|\ dx.
		\end{equation}
		
		\item[\textbf{(3)}] For each $m\in\mathbb{N}$ we have that
		
		\begin{equation}\label{est3}
			\dfrac{1}{L_0}\left (\lambda+\dfrac{1}{\Delta_m}\right )\sum_{n=1}^{2^m}\int_{\Omega} \big [b(x,v_{m,n})-b(x,v_{m,n-1}) \big ]^2\leq\mathcal{A}(u_0)+\lambda L_0(\delta-\varepsilon)\sum_{n=1}^{2^m} \int_{\Omega} |v_{m,n}-v_{m,n-1}|\ dx.
		\end{equation}
		
		\item[\textbf{(4)}] For each $m\in\mathbb{N}$ and for every $n\in\{1,2,\hdots,2^m\}$ the following inequality holds
		
		\begin{equation}\label{est4}
			\mathcal{A}(v_{m,n})\leq \mathcal{A}(u_0)+\dfrac{|\Omega|T\lambda^2L_0^2(\delta-\varepsilon)^2}{4\ell_0}.
		\end{equation}

		\noindent In particular there is a constant $C_0>0$ such that for each $m\in\mathbb{N}$ and for every $n\in\{1,2,\hdots,2^m\}$ we have that:
		
		\begin{equation}\label{est5}
			\Vert v_{m,n}\Vert_{W^{1,p(x)}(\Omega)}, \rho_{p(x)}(|\nabla v_{m,n}|)\leq C_0.
		\end{equation}
	\end{enumerate}
\end{proposition}

\begin{proof}\noindent\textbf{(1)} For each $n\in\{1,2,\hdots,2^m\}$ we have from \eqref{eq1112} that for any $\phi\in W^{1,p(x)}(\Omega)$ the following equality holds
	
\begin{equation}
	\dfrac{1}{\Delta_m}\int_{\Omega} \big [b(x,v_{m,n})-b(x,v_{m,n-1})]\cdot\phi\ dx+\int_{\Omega} \mathbf{a}(x,\nabla v_{m,n})\cdot\nabla\phi\ dx+\lambda\int_{\Omega} b(x,v_{m,n})\phi\ dx=\int_{\Omega} g_{m,n}(x)\phi\ dx.
\end{equation}

\noindent Choosing $\phi=v_{m,n}-v_{m,n-1}\in W^{1,p(x)}(\Omega)$\footnote{Here it is essential that $u_0\in W^{1,p(x)}(\Omega)$, because $v_{m,0}=u_0$.} the above relation will become

\begin{align*}
	&\left (\lambda+\dfrac{1}{\Delta_m}\right )\int_{\Omega} \big [b(x,v_{m,n})-b(x,v_{m,n-1})]\cdot (v_{m,n}-v_{m,n-1})\ dx+\int_{\Omega} \mathbf{a}(x,\nabla v_{m,n})\cdot\big (\nabla v_{m,n}-\nabla v_{m,n-1}\big )\ dx\\
	&=\int_{\Omega} \big [g_{m,n}-\lambda b(x,v_{m,n-1})\big ]\cdot (v_{m,n}-v_{m,n-1})\ dx.
\end{align*}

\noindent Using the inequality $\mathbf{a}(x,\xi_2)\cdot(\xi_2-\xi_1)\geq A(x,\xi_2)-A(x,\xi_2)$ which is true for a.e. $x\in\Omega$ and for any $\xi_1,\xi_2\in\mathbb{R}^N$\footnote{For a proof, see \cite[Proposition 3.1 \textbf{(3)}]{max4}.} we deduce that 

\[
\int_{\Omega} \mathbf{a}(x,\nabla v_{m,n})\cdot\big (\nabla v_{m,n}-\nabla v_{m,n-1}\big )\ dx\geq \int_{\Omega} A(x,\nabla v_{m,n})\ dx-\int_{\Omega} A(x,\nabla v_{m,n-1})\ dx=\mathcal{A}(v_{m,n})-\mathcal{A}(v_{m,n-1}).
\] 

\noindent Therefore we can write that

\begin{align*}
	&\left (\lambda+\dfrac{1}{\Delta_m}\right )\int_{\Omega} \big [b(x,v_{m,n})-b(x,v_{m,n-1})]\cdot (v_{m,n}-v_{m,n-1})\ dx+\mathcal{A}(v_{m,n})-\mathcal{A}(v_{m,n-1})\\
&\leq	\int_{\Omega} \big [g_{m,n}-\lambda b(x,v_{m,n-1})\big ]\cdot (v_{m,n}-v_{m,n-1})\ dx\\
&\leq  \int_{\Omega} |g_{m,n}-\lambda b(x,v_{m,n-1})\big |\cdot |v_{m,n}-v_{m,n-1}|\ dx\\
\text{(Claim 1)}\ \ \ &\leq  \lambda\int_{\Omega} |b(x,\delta)-b(x,\varepsilon)|\cdot |v_{m,n}-v_{m,n-1}|\ dx\\
\text{Prop. \ref{propoext} \textbf{(3)}}\ \ \ & \leq \lambda L_0(\delta-\varepsilon)\int_{\Omega} |v_{m,n}-v_{m,n-1}|\ dx.
\end{align*}

\noindent Adding up all these relations written for $n=1,2,\hdots, 2^m$ and using the simple fact that $\mathcal{A}(v_{m,2^m})\geq 0$ gives us the desired inequality.
	
\medskip

\noindent\textbf{(2)} This estimation comes from the fact that $b(x,\cdot)$ is strictly increasing on $[\varepsilon,\delta]$ and from \textbf{(H13)}: $0\leq \big [b(x,v_{m,n})-b(x,v_{m,n-1})\big ]\cdot (v_{m,n}-v_{m,n-1})=|b(x,v_{m,n})-b(x,v_{m,n-1})|\cdot |v_{m,n}-v_{m,n-1}|\geq \ell_0 (v_{m,n}-v_{m,n-1})^2$ for a.e. $x\in\Omega$.

\medskip

\noindent\textbf{(3)} This estimation comes from the fact that $b(x,\cdot)$ is strictly increasing on $[\varepsilon,\delta]$ and from Proposition \ref{propoext} \textbf{(3)}: $0\leq \big [b(x,v_{m,n})-b(x,v_{m,n-1})\big ]\cdot (v_{m,n}-v_{m,n-1})=|b(x,v_{m,n})-b(x,v_{m,n-1})|\cdot |v_{m,n}-v_{m,n-1}|\geq \dfrac{1}{L_0} [b(x,v_{m,n})-b(x,v_{m,n-1}]^2$ for a.e. $x\in\Omega$.

\medskip

\noindent\textbf{(4)} We know from the proof of \textbf{(1)} that for any $k\in\overline{1,n}$ the following inequality holds:

\begin{align*}
	&\left (\lambda+\dfrac{1}{\Delta_m}\right )\int_{\Omega} \big [b(x,v_{m,k})-b(x,v_{m,k-1})]\cdot (v_{m,k}-v_{m,k-1})\ dx+\mathcal{A}(v_{m,k})-\mathcal{A}(v_{m,k-1})\\
	\leq & \lambda L_0(\delta-\varepsilon)\int_{\Omega} |v_{m,k}-v_{m,k-1}|\ dx.
\end{align*} 

\noindent Summing these relations for $k\in\overline{1,n}$ produces the following inequality

\begin{align*}
	&\left (\lambda+\dfrac{1}{\Delta_m}\right )\sum_{k=1}^n \int_{\Omega} \big [b(x,v_{m,k})-b(x,v_{m,k-1})]\cdot (v_{m,k}-v_{m,k-1})\ dx+\mathcal{A}(v_{m,n})-\underbrace{\mathcal{A}(v_{m,0})}_{=\mathcal{A}(u_0)}\\
	\leq &\lambda L_0(\delta-\varepsilon)\sum_{k=1}^n\int_{\Omega} |v_{m,k}-v_{m,k-1}|\ dx.
\end{align*}

\noindent Therefore we may write that

\begin{align*}
&\mathcal{A}(v_{m,n})\leq \mathcal{A}(u_0) \\
& +\sum_{k=1}^n \int_{\Omega} \lambda L_0(\delta-\varepsilon)|v_{m,k}-v_{m,k-1}|-\left (\lambda+\dfrac{1}{\Delta_m}\right )\big [b(x,v_{m,k})-b(x,v_{m,k-1})]\cdot (v_{m,k}-v_{m,k-1})\ dx\\
&\stackrel{\textbf{(H13)}}{\leq}  \mathcal{A}(u_0)+\sum_{k=1}^n \int_{\Omega} \lambda L_0(\delta-\varepsilon)|v_{m,k}-v_{m,k-1}|-\ell_0\left (\lambda+\dfrac{1}{\Delta_m}\right ) (v_{m,k}-v_{m,k-1})^2\ dx\\
&\leq  \mathcal{A}(u_0)+\sum_{k=1}^n \int_{\Omega}\dfrac{\Delta_m\lambda^2 L_0^2(\delta-\varepsilon)^2}{4\lambda\ell_0\Delta_m+4\ell_0}\ dx=\mathcal{A}(u_0)+\dfrac{n\Delta_m|\Omega|\lambda^2 L_0^2(\delta-\varepsilon)^2}{4\lambda\ell_0\Delta_m+4\ell_0}\\
&\stackrel{n\leqslant 2^m}{\leq} \mathcal{A}(u_0)+\dfrac{2^m\Delta_m|\Omega|\lambda^2 L_0^2(\delta-\varepsilon)^2}{4\lambda\ell_0\Delta_m+4\ell_0}=\mathcal{A}(u_0)+\dfrac{T|\Omega|\lambda^2 L_0^2(\delta-\varepsilon)^2}{4\lambda\ell_0\Delta_m+4\ell_0}\\
&\leq \mathcal{A}(u_0)+\dfrac{T|\Omega|\lambda^2 L_0^2(\delta-\varepsilon)^2}{4\ell_0}.
\end{align*}

\noindent We have used above the elementary quadratic inequality 

\[
\lambda L_0(\delta-\varepsilon)\cdot a-\ell_0\left (\lambda+\dfrac{1}{\Delta_m}\right )\cdot a^2\leq \dfrac{\Delta_m\lambda^2 L_0^2(\delta-\varepsilon)^2}{4\lambda\ell_0\Delta_m+4\ell_0},\ \forall\ a\in\mathbb{R}.
\]

\noindent Next we will show the existence of a constant $C_0>0$ with $\rho_{p(x)}(|\nabla v_{m,n})\leq C_0$ for any $m\in\mathbb{N}$ and any $n\in\{1,2,\hdots,2^m\}$. Suppose that such a constant do not exist. Therefore we can find two sequences of natural numbers $(m_k)_{k\geq 1}$ and $(n_{k})_{k\geq 1}$ with $1\leq n_k\leq 2^{m_{k}}$ such that $\lim\limits_{k\to\infty} \rho_{p(x)}\big (|\nabla v_{m_{k},n_{k}}|\big )=\infty$. Using now the key hypothesis \textbf{(H7)} we deduce that $\lim\limits_{k\to\infty} \mathcal{A}(v_{m_k, n_k})=+\infty$. But this is a constradiction, because we have already proved that $\mathcal{A}(v_{m_k, n_k})\leq \mathcal{A}(u_0)+\dfrac{T|\Omega|\lambda^2 L_0^2(\delta-\varepsilon)^2}{4\ell_0}$, for any $k\geq 1$. Therefore we can find such a constant $C_0$.

\noindent Note that:

\begin{align*}
	\Vert v_{m,n}\Vert_{W^{1,p(x)}(\Omega)}&=\Vert v_{m,n}\Vert_{L^{p(x)}(\Omega)}+\Vert |\nabla v_{m,n}|\Vert_{L^{p(x)}(\Omega)}\stackrel{\text{Prop. \ref{propu2u1}}}{\leq} \delta \Vert 1\Vert_{L^{p(x)}(\Omega)}+\Vert |\nabla v_{m,n}|\Vert_{L^{p(x)}(\Omega)}\\
	&\leq \delta \Vert 1\Vert_{L^{p(x)}(\Omega)}+\max\{\rho_{p(x)}(|\nabla v_{m,n}|)^{\frac{1}{p^+}},\rho_{p(x)}(|\nabla v_{m,n}|)^{\frac{1}{p^-}}\}.
\end{align*}

\noindent We have used here \cite[Theorem 14, \textbf{(22),(23)}]{max2}. This shows that we can choose a bigger value for $C_0>0$ such that both inequalities will hold for all $m\in\mathbb{N}$ and all $n\in\{1,2,\hdots,2^m\}$.
\end{proof}

\noindent\underline{\textbf{Claim 3:}} \textbf{If $u_0\in W^{1,p(x)}\cap\mathcal{U}_{[\varepsilon,\delta]}$ then the sequence $\left (\dfrac{\partial b(\cdot,\tilde{v}_m)}{\partial t} \right )_{m\geq 0}=\left (\dfrac{\partial \beta_m}{\partial t} \right )_{m\geq 0}$ is bounded in $L^2\big (0,T;L^2(\Omega)\big )$.}

\begin{proof} From Proposition \ref{propobeta} \textbf{(5)} we know that for each $m\in\mathbb{N}$: $\beta_m\in L^{\infty}\big (0,T;L^{\infty}(\Omega)\big )\subset L^2\big (0,T;L^2(\Omega)\big )$ and for a.e. $(t,x)\in (0,T)\times\Omega$
	
	\begin{equation}
		\dfrac{\partial\beta_m}{\partial t}(t,x)=\displaystyle\sum_{n=1}^{2^m} \dfrac{b\big (x,v_{m,n}(x)\big )-b\big (x,v_{m,n-1}(x)\big )}{\Delta_m}\chi_{(t_{m,n-1},t_{m,n}]}(t).
	\end{equation}
	
\noindent Therefore

\begin{align*}
\left \Vert \dfrac{\partial\beta_m}{\partial t}\right\Vert_{L^2(0,T;L^2(\Omega))}^2	&=\int_{0}^T\left \Vert \dfrac{\partial\beta_m}{\partial t}(t,\cdot)\right\Vert_{L^2(\Omega)}^2\ dt=\sum_{n=1}^{2^m}\int_{t_{m,n-1}}^{t_{m,n}}\underbrace{\left \Vert \dfrac{\partial\beta_m}{\partial t}(t,\cdot)\right\Vert_{L^2(\Omega)}^2}_{\text{independent of }t}\ dt =\\
&=\sum_{n=1}^{2^m}\Delta_m\left \Vert \dfrac{b\big (\cdot,v_{m,n}\big )-b\big (\cdot,v_{m,n-1}\big )}{\Delta_m}\right\Vert_{L^2(\Omega)}^2\\
&=\sum_{n=1}^{2^m}\Delta_m\int_{\Omega} \left [\dfrac{b\big (x,v_{m,n}(x)\big )-b\big (x,v_{m,n-1}(x)\big )}{\Delta_m} \right ]^2 dx\\
&=\sum_{n=1}^{2^m}\dfrac{1}{\Delta_m}\int_{\Omega}\left [b\big (x,v_{m,n}(x)\big )-b\big (x,v_{m,n-1}(x)\big )\right ]^2 dx\\
\eqref{est3}\ \ \ &\leq\dfrac{L_0}{\lambda\Delta_m+1}\left [\mathcal{A}(u_0)+\lambda L_0(\delta-\varepsilon)\sum_{n=1}^{2^m} \int_{\Omega} |v_{m,n}-v_{m,n-1}|\ dx\right ] \\
\text{\textbf{(H13)}} \ \ \ &\leq\dfrac{L_0}{\lambda\Delta_m+1}\left [\mathcal{A}(u_0)+\dfrac{\lambda L_0(\delta-\varepsilon)}{\ell_0}\sum_{n=1}^{2^m} \int_{\Omega} |b(x,v_{m,n})-b(x,v_{m,n-1})|\ dx\right ]  .
\end{align*}

\noindent Remark that

\begin{align*}
\left \Vert \dfrac{\partial\beta_m}{\partial t}\right\Vert_{L^1(0,T;L^1(\Omega))}	&=\int_{0}^T\left \Vert \dfrac{\partial\beta_m}{\partial t}(t,\cdot)\right\Vert_{L^1(\Omega)}\ dt=\sum_{n=1}^{2^m}\int_{t_{m,n-1}}^{t_{m,n}}\underbrace{\left \Vert \dfrac{\partial\beta_m}{\partial t}(t,\cdot)\right\Vert_{L^1(\Omega)}}_{\text{independent of }t}\ dt =\\
&=\sum_{n=1}^{2^m}\Delta_m\left \Vert \dfrac{b\big (\cdot,v_{m,n}\big )-b\big (\cdot,v_{m,n-1}\big )}{\Delta_m}\right\Vert_{L^1(\Omega)}\\
&=\sum_{n=1}^{2^m}\Delta_m\int_{\Omega} \left |\dfrac{b\big (x,v_{m,n}(x)\big )-b\big (x,v_{m,n-1}(x)\big )}{\Delta_m}\right |\ dx\\
&=\sum_{n=1}^{2^m}\int_{\Omega} |b\big (x,v_{m,n}(x)\big )-b\big (x,v_{m,n-1}(x)\big )|\ dx
\end{align*}

\noindent Also, from \textit{Cauchy inequality} and \textit{Tonelli's theorem}\footnote{See Proposition 5.2.1, page 147 from \cite{Cohn}.} we have that

\begin{align*}
	\left \Vert \dfrac{\partial\beta_m}{\partial t}\right\Vert_{L^1(0,T;L^1(\Omega))}&=\int_{(0,T)\times\Omega}\left |\dfrac{\partial\beta_m}{\partial t} \right |\ dx\ dt\leq \left (\int_{(0,T)\times\Omega} 1^2\ dx\ dt \right )^{\frac{1}{2}}\cdot\left [\int_{(0,T)\times\Omega}\left (\dfrac{\partial\beta_m}{\partial t} \right )^2\ dx\ dt \right ]^{\frac{1}{2}}\\
	&=\sqrt{T|\Omega|}\cdot \left \Vert \dfrac{\partial\beta_m}{\partial t}\right\Vert_{L^2(0,T;L^2(\Omega))}.
\end{align*}

\noindent Combining the last three relations leads us to

\begin{align*}
	\left \Vert \dfrac{\partial\beta_m}{\partial t}\right\Vert_{L^2(0,T;L^2(\Omega))}^2&\leq \dfrac{L_0}{\lambda\Delta_m+1}\left [\mathcal{A}(u_0)+\dfrac{\lambda L_0(\delta-\varepsilon)}{\ell_0}\left \Vert \dfrac{\partial\beta_m}{\partial t}\right\Vert_{L^1(0,T;L^1(\Omega))}\right ]\\
	&\leq \dfrac{L_0}{\lambda\Delta_m+1}\left [\mathcal{A}(u_0)+\dfrac{\lambda L_0(\delta-\varepsilon)\sqrt{T|\Omega|}}{\ell_0}\cdot \left \Vert \dfrac{\partial\beta_m}{\partial t}\right\Vert_{L^2(0,T;L^2(\Omega))}\right ]\\
	&\leq L_0\mathcal{A}(u_0)+\dfrac{\lambda L_0^2(\delta-\varepsilon)\sqrt{T|\Omega|}}{\ell_0}\cdot \left \Vert \dfrac{\partial\beta_m}{\partial t}\right\Vert_{L^2(0,T;L^2(\Omega))}.	
\end{align*}

\noindent Solving this quadratic inequality gives us $\forall\ m\in\mathbb{N}$ that

\begin{equation}
\left \Vert \dfrac{\partial\beta_m}{\partial t}\right\Vert_{L^2(0,T;L^2(\Omega))}\leq
\frac{\lambda L_0^2(\delta-\varepsilon)\sqrt{T|\Omega|}}{2\ell_0}+\sqrt{\frac{\lambda^2 L_0^4(\delta-\varepsilon)^2 T|\Omega|}{4\ell_0^2}+L_0\,\mathcal{A}(u_0)}:=C_1.
\end{equation}

\noindent Since the right-hand side is independent of $m\in\mathbb{N}$ we deduce that the sequence $\left (\dfrac{\partial \beta_m}{\partial t} \right )_{m\geq 0}$ is bounded in $L^2\big (0,T;L^2(\Omega)\big )$.

\end{proof}

\noindent\underline{\textbf{Claim 4:}} \textbf{If $u_0\in W^{1,p(x)}\cap\mathcal{U}_{[\varepsilon,\delta]}$ then the sequence $\left (\dfrac{\partial\tilde{v}_m}{\partial t} \right )_{m\geq 0}$ is bounded in $L^2\big (0,T;L^2(\Omega)\big )$.}

\begin{proof} From Proposition \ref{propovtilde} \textbf{(5)} we know that $\dfrac{\partial\tilde{v}_m}{\partial t}\in L^\infty\big (0,T;L^\infty(\Omega)\big )\subset L^2\big (0,T;L^2(\Omega)\big )$ and for a.e. $(t,x)\in (0,T)\times\Omega$ we have that
	
\begin{equation}
	\dfrac{\partial\tilde{v}_m}{\partial t}(t,x)=\dfrac{\partial\mathfrak{b}}{\partial\sigma}(x,\beta_m(t,x))\cdot\dfrac{\partial\beta_m}{\partial t}(t,x).
\end{equation}

\noindent Therefore, from Proposition \ref{propobrussian} \textbf{(5)} we get for a.e. $(t,x)\in\Omega$ that

 \[
 \left |\dfrac{\partial\tilde{v}_m}{\partial t}(t,x) \right |=\dfrac{\partial\mathfrak{b}}{\partial\sigma}\underbrace{(x,\beta_m(t,x))}_{\in\widehat{D}} \cdot \left |\dfrac{\partial\beta_m}{\partial t}(t,x) \right |\leq \dfrac{1}{\ell_0}\cdot \left |\dfrac{\partial\beta_m}{\partial t}(t,x) \right |.
 \] 

\noindent In particular this shows that $\left \Vert\dfrac{\partial\tilde{v}_m}{\partial t} \right\Vert_{L^2(0,T;L^2(\Omega))}\leq \dfrac{1}{\ell_0}\cdot \left \Vert\dfrac{\partial\beta_m}{\partial t} \right \Vert_{L^2(0,T;L^2(\Omega))}\leq\dfrac{C_1}{\ell_0}$. Thus for every $m\in\mathbb{N}$, from Claim 3 we get that

\begin{equation}
\left \Vert\dfrac{\partial\tilde{v}_m}{\partial t} \right\Vert_{L^2(0,T;L^2(\Omega))}\leq	\frac{\lambda L_0^2(\delta-\varepsilon)\sqrt{T|\Omega|}}{2\ell_0^2}+\sqrt{\frac{\lambda^2 L_0^4(\delta-\varepsilon)^2 T|\Omega|}{4\ell_0^4}+\dfrac{L_0}{\ell_0^2}\,\mathcal{A}(u_0)}.
\end{equation}

\noindent Since the right-hand side is independent of $m\in\mathbb{N}$ we deduce that the sequence $\left (\dfrac{\partial \tilde{v}_m}{\partial t} \right )_{m\geq 0}$ is bounded in $L^2\big (0,T;L^2(\Omega)\big )$.	
\end{proof}

\noindent\underline{\textbf{Claim 5:}} \textbf{As $m\to\infty$ we have that $b(\cdot,\tilde{v}_m)-b(\cdot,v_m)\to 0$ in $L^{\infty}\big (0,T;L^2(\Omega)\big )$.}

\begin{proof} For any. $t\in (0,T]$ we may write that
	
	\begin{align*}
		\Vert b(\cdot,\tilde{v}_m(t,\cdot))-b(\cdot,v_m(t,\cdot))\Vert_{L^2(\Omega)}&=\left\Vert\sum_{n=1}^{2^m}\left (1-\dfrac{t-t_{m,n-1}}{\Delta_m}\right )\cdot\big (b(\cdot,v_{m,n})-b(\cdot,v_{m,n-1}) \big )\chi_{(t_{m,n-1},t_{m,n}]}(t) \right\Vert_{L^2(\Omega)}\\
		&\leq \sum_{n=1}^{2^m}\underbrace{\left (1-\dfrac{t-t_{m,n-1}}{\Delta_m}\right )}_{\leq 1}\Vert b(\cdot,v_{m,n})-b(\cdot,v_{m,n-1})\Vert_{L^2(\Omega)}\chi_{(t_{m,n-1},t_{m,n}]}(t)\\
		&\leq \sum_{n=1}^{2^m}\Vert b(\cdot,v_{m,n})-b(\cdot,v_{m,n-1})\Vert_{L^2(\Omega)}\chi_{(t_{m,n-1},t_{m,n}]}(t)\\
		&\leq \max_{n\in\overline{1,2^m}} \Vert b(\cdot,v_{m,n})-b(\cdot,v_{m,n-1})\Vert_{L^2(\Omega)}\\
		&\leq \sqrt{\sum_{n=1}^{2^m}\Vert b(\cdot,v_{m,n})-b(\cdot,v_{m,n-1})\Vert_{L^2(\Omega)}^2}\\
		&=\sqrt{\sum_{n=1}^{2^m} \int_{\Omega} \big [b(x,v_{m,n}(x))-b(x,v_{m,n-1}(x))\big ]^2\ dx}\\
		&=\sqrt{\Delta_m}\cdot\sqrt{\sum_{n=1}^{2^m} \dfrac{1}{\Delta_m}\int_{\Omega} \big [b(x,v_{m,n}(x))-b(x,v_{m,n-1}(x))\big ]^2\ dx}\\
		&\stackrel{\text{Claim 3}}{=}\sqrt{\Delta_m}\cdot \left\Vert\dfrac{\partial\beta_m}{\partial t}\right\Vert_{L^2(0,T;L^2(\Omega))}\leq C_1\sqrt{\Delta_m}.
	\end{align*}
	
	\noindent This means that $\Vert b(\cdot,\tilde{v}_m)-b(\cdot,v_m)\Vert_{L^\infty(0,T;L^2(\Omega))}\leq  C_1\sqrt{\Delta_m}$ for every $m\in\mathbb{N}$. But, $\lim\limits_{m\to\infty} C_1\sqrt{\Delta_m}=\lim\limits_{m\to\infty}\dfrac{C_1\sqrt{T}}{2^{m/2}}=0$. In conclusion $\lim\limits_{m\to\infty} \Vert b(\cdot,\tilde{v}_m)-b(\cdot,v_m)\Vert_{L^\infty(0,T;L^2(\Omega))}=0$.
	
\end{proof}

\noindent\underline{\textbf{Claim 6:}} \textbf{As $m\to\infty$ we have that $\tilde{v}_m-v_m\to 0$ in $L^{\infty}\big (0,T;L^2(\Omega)\big )$.}

\begin{proof}\noindent Just note that from \textbf{(H13)} we have for a.e. $(t,x)\in (0,T)\times\Omega$ that

\begin{equation}
	|\tilde{v}_m(t,x)-v_m(t,x)|\leq \dfrac{1}{\ell_0}|b(x,\tilde{v}_m(t,x))-b(x,v_m(t,x))|.
\end{equation}

\noindent Therefore for any $t\in (0,T]$

 $$\Vert \tilde{v}_m(t,\cdot)-v_m(t,\cdot) \Vert_{L^2(\Omega)}\leq \dfrac{1}{\ell_0}\Vert b(\cdot,\tilde{v}_m(t,\cdot))-b(\cdot,v_m(t,\cdot))\Vert_{L^2(\Omega)}\leq \dfrac{C_1}{\ell_0}\sqrt{\Delta}\stackrel{m\to\infty}{\longrightarrow}0,$$

\noindent which shows that $\lim\limits_{m\to\infty} \Vert \tilde{v}_m-v_m\Vert_{L^\infty(0,T;L^2(\Omega))}=0$.
\end{proof}

\noindent\underline{\textbf{Claim 7:}} \textbf{If $u_0\in W^{1,p(x)}(\Omega)\cap\mathcal{U}_{[\varepsilon,\delta]}$, the sequences $(v_m)_{m\geq 0},\ \big (b(\cdot,v_m) \big )_{m\geq 0},\ (\tilde{v}_m)_{m\geq 0}$ and $\big (b(\cdot,\tilde{v}_m) \big )_{m\geq 0}=\big (\beta_m \big )_{m\geq 0}$ are all bounded in $L^{\infty}\big (0,T;W^{1,p(x)}(\Omega)\big )$.}

\begin{proof} \noindent$\bullet$ We already know that $(v_m)_{m\geq 0}\subset L^{\infty}\big (0,T;W^{1,p(x)}(\Omega)\big )$ from Proposition \ref{propovm} \textbf{(1)} and we have that:
	
	\[
	\Vert v_m\Vert_{L^{\infty}(0,T;W^{1,p(x)}(\Omega))}\stackrel{\text{Prop. \ref{propovm} \textbf{(1)}}}{\leq} \max_{n\in \overline{1,2^m}} \Vert v_{m,n}\Vert_{W^{1,p(x)}(\Omega)} \stackrel{\text{Prop. \ref{propovmn} \textbf{(4)}}}{\leq} C_0.
	\]

\noindent$\bullet$ From Proposition \ref{propovm} \textbf{(4)} we know that $\big (b(\cdot,v_m)\big )_{m\geq 0}\subset L^{\infty}\big (0,T;W^{1,p(x)}(\Omega)\big )$ and

\begin{align*}
	\Vert b(\cdot,v_m)\Vert_{L^{\infty}(0,T;W^{1,p(x)}(\Omega))}&\leq \max_{n\in\overline{1,2^m}} \Vert b(\cdot,v_{m,n})\Vert_{W^{1,p(x)}(\Omega)}\leq L_0 \max_{n\in\overline{1,2^m}} \Vert v_{m,n}\Vert_{W^{1,p(x)}(\Omega)}+L_0\Vert 1\Vert_{L^{p(x)}(\Omega)}\\
	&\stackrel{\eqref{est5}}{\leq} L_0C_0+L_0\Vert 1\Vert_{L^{p(x)}(\Omega)},\ \forall\ m\in\mathbb{N}.
\end{align*}

\noindent$\bullet$ From Proposition \ref{propobeta} \textbf{(4)} we have that $(\beta_m)_{m\geq 0}\subset L^{\infty}\big (0,T;W^{1,p(x)}(\Omega)\big )$ and

	\begin{align*}
		\Vert \beta_m\Vert_{L^{\infty}(0,T;W^{1,p(x)}(\Omega))}&\leq \max_{n\in\overline{1,2^m}} \Vert b(\cdot,v_{m,n})\Vert_{W^{1,p(x)}(\Omega)}\leq L_0 \max_{n\in\overline{1,2^m}} \Vert v_{m,n}\Vert_{W^{1,p(x)}(\Omega)}+L_0\Vert 1\Vert_{L^{p(x)}(\Omega)}\\
		&\stackrel{\eqref{est5}}{\leq} L_0C_0+L_0\Vert 1\Vert_{L^{p(x)}(\Omega)},\ \forall\ m\in\mathbb{N}.
	\end{align*}
	
\noindent $\bullet$ From Proposition \ref{propovtilde} \textbf{(4)} we have that $(\tilde{v}_m)_{m\geq 0}\subset L^{\infty}\big (0,T;W^{1,p(x)}(\Omega)\big )$ and

	\begin{align*}
	\Vert \tilde{v}_m\Vert_{L^{\infty}(0,T;W^{1,p(x)}(\Omega))}&\leq \dfrac{1}{\ell_0}\Vert \beta_m\Vert_{L^{\infty}(0,T;W^{1,p(x)}(\Omega))}+\dfrac{L_0}{\ell_0}\Vert 1\Vert_{L^{p(x)}(\Omega)}+\dfrac{1}{\ell_0}\Vert b(\cdot,\varepsilon)-\ell_0\varepsilon\Vert_{L^{p(x)}(\Omega)}\\
	&\leq \dfrac{L_0C_0}{\ell_0}+2\dfrac{L_0}{\ell_0}\Vert 1\Vert_{L^{p(x)}(\Omega)}+\dfrac{1}{\ell_0}\Vert b(\cdot,\varepsilon)-\ell_0\varepsilon\Vert_{L^{p(x)}(\Omega)} .
\end{align*}
	
\end{proof}

\noindent\textbf{\underline{Claim 8:}} \textbf{If $u_0\in W^{1,p(x)}(\Omega)$ then there are some $v,\tilde{v}\in L^{\infty}\big (0,T;W^{1,p(x)}(\Omega)\big )$ such that for the same subsequence $v_{m_k}\rightharpoonup v$ and $\tilde{v}_{m_k}\rightharpoonup \tilde{v}$ in $L^r\big (0,T;W^{1,p(x)}(\Omega)\big )$ for each $r\in [1,\infty)$. Moreover $v_{m_k}\weakstar v$ and $\tilde{v}_{m_k}\weakstar \tilde{v}$ in $L^\infty\big (0,T;W^{1,p(x)}(\Omega)\big )$. Furthermore one gets that $\Vert v\Vert_{L^{\infty}(0,T;W^{1,p(x)}(\Omega))}\leq\displaystyle\liminf_{k\to\infty}\Vert v_{m_k}\Vert_{L^{\infty}(0,T;W^{1,p(x)}(\Omega))}$ and $\Vert \tilde{v}\Vert_{L^{\infty}(0,T;W^{1,p(x)}(\Omega))}\leq\displaystyle\liminf_{k\to\infty}\Vert \tilde{v}_{m_k}\Vert_{L^{\infty}(0,T;W^{1,p(x)}(\Omega))}$.}

\begin{proof} From \cite[Theorem 3.1]{kovacik} we have that $W^{1,p(x)}(\Omega)$ is a reflexive and separable Banach space. Therefore, using \textbf{Claim 7} and Proposition \ref{propocompactness1} we find some $v\in L^{\infty}\big (0,T;W^{1,p(x)}(\Omega)\big )$ and a subsequence with the following two properties $v_{m_k}\rightharpoonup v$ in $L^r\bigl(0,T;W^{1,p(x)}(\Omega)\bigr)$ and $v_{m_k}\weakstar v$ in $L^\infty\bigl(0,T;W^{1,p(x)}(\Omega)\bigr)$ . Applying one more time Proposition \ref{propocompactness1}, but now for $(\tilde{v}_{m_k})_{k\geq 1}$, we find some function $\tilde{v}\in L^{\infty}\big (0,T;W^{1,p(x)}(\Omega)\big )$ and a further subsequence (still denoted for simplicity with $m_k$ instead of $m_{k_\ell}$) such that $\tilde{v}_{m_k}\rightharpoonup \tilde{v}$ in $L^r\big (0,T;W^{1,p(x)}(\Omega)\big )$ and $\tilde{v}_{m_k}\weakstar \tilde{v}$ in $L^\infty\big (0,T;W^{1,p(x)}(\Omega)\big )$. The two inequalities follow from the standard properties of the weak* convergence.\footnote{See Proposition 3.13 (iii) in \cite[page 63]{brezis2011functional}.}
	
\end{proof}

\noindent\textbf{\underline{Claim 9:}} \textbf{$\tilde{v}=v$}.

\begin{proof} From hypothesis \textbf{(H2)} we know that $W^{1,p(x)}(\Omega)\hookrightarrow L^2(\Omega)$. Since from \textbf{Claim 8}, $v_{m_k}\weak v$ and $\tilde{v}_{m_k}\weak \tilde{v}$ in $L^2\bigl(0,T;W^{1,p(x)}(\Omega)\bigr)$ as $k\to\infty$. We are now in position to apply Lemma \ref{lemmaweakbochner} to obtain that $v_{m_k}\weak v$ and $\tilde{v}_{m_k}\weak \tilde{v}$ in $L^2\bigl(0,T;L^{2}(\Omega)\bigr)$. Therefore $\tilde{v}_{m_k}-v_{m_k}\weak \tilde{v}-v$ in $L^2\bigl(0,T;L^{2}(\Omega)\bigr)$. 
	
\noindent Note that, from \textbf{Claim 6}, we know that $\tilde{v}_{m_k}-v_{m_k}\to 0$ in $L^{\infty}\bigl(0,T;L^2(\Omega)\bigr)$. Also, Proposition 2.2.5 from \cite[page 128]{gasinski2005nonlinear} implies that $L^{\infty}\bigl(0,T;L^2(\Omega)\bigr)\hookrightarrow L^2\bigl(0,T;L^2(\Omega)\bigr)$. Combining these two facts we get that $\tilde{v}_{m_k}-v_{m_k}\to 0$ in $L^{2}\bigl(0,T;L^2(\Omega)\bigr)$. Now, since strong convergence implies weak convergence, we will also get that $\tilde{v}_{m_k}-v_{m_k}\weak 0$ in $L^{2}\bigl(0,T;L^2(\Omega)\bigr)$. But we have already shown that $\tilde{v}_{m_k}-v_{m_k}\weak \tilde{v}-v$ in $L^2\bigl(0,T;L^{2}(\Omega)\bigr)$. From the uniqueness of the weak limit we conclude that $\tilde{v}=v$.
	
\end{proof}

\noindent\textbf{\underline{Claim 10:}} \textbf{If $u_0\in W^{1,p(x)}(\Omega)$, then $v\in H^1\bigl((0,T);L^2(\Omega)\bigr)$ and on a further subsequence, denoted for simplicity in the same way, we have that $\dfrac{\partial\tilde{v}_{m_k}}{\partial t}\weak \dfrac{\partial v}{\partial t}$ in $L^2\bigl(0,T;L^2(\Omega)\bigr)$.}

\begin{proof} From Proposition 2.2.3 (c) given in \cite[page 127]{gasinski2005nonlinear} we get that $L^2\bigl(0,T;L^2(\Omega)\bigr)$ is a reflexive Banach space. Also, from \textbf{Claim 4} we know that the sequence $\left (\dfrac{\partial\tilde{v}_{m_k}}{\partial t}\right )_{k\geq 1}$ is bounded in $L^2\bigl(0,T;L^2(\Omega)\bigr)$. Here we have used the fact that $u_0\in W^{1,p(x)}(\Omega)$. Therefore, using \textit{Eberlein-\v{S}mulian Theorem}\footnote{For a proof, see \cite[Theorem 3.4.16, page 221]{papageorgiou2018applied}.}, we get that there is a function $w\in L^2\bigl(0,T;L^2(\Omega)\bigr)$ such that on a further subsequence, denoted in the same way for simplicity, we have $\dfrac{\partial\tilde{v}_{m_k}}{\partial t}\weak w$ in $L^2\bigl(0,T;L^2(\Omega)\bigr)$.
	
\noindent Now fix any $\varphi\in C^{\infty}_{c}(0,T)$. Using Proposition \ref{propovtilde} \textbf{(5)} and the definition of the weak derivative we have for each $k\geq 1$ that:

\begin{equation}\label{weakderivativeeq1}
\int_{0}^T\tilde{v}_{m_k}(t,\cdot)\varphi'(t)\ dt=-\int_{0}^T \dfrac{\partial \tilde{v}_{m_k}}{\partial t}(t,\cdot)\varphi(t)\ dt.
\end{equation}

\noindent Define the operators $\mathcal{T}_{\varphi},\mathcal{T}_{\varphi'}:L^2\bigl(0,T;L^2(\Omega)\bigr)\to L^2(\Omega)$ via the following Bochner integrals:

\begin{equation}
	\mathcal{T}_{\varphi}(u)=\int_{0}^T u(t,\cdot)\varphi(t)\ dt,\ \text{and}\ \mathcal{T}_{\varphi'}(u)=\int_{0}^T u(t,\cdot)\varphi'(t)\ dt,\ \forall\ u\in L^2\bigl(0,T;L^2(\Omega)\bigr).
\end{equation}

\noindent From Lemma \ref{bochnercontinuous} we have that $\mathcal{T}_{\varphi}$ and $\mathcal{T}_{\varphi'}$ are well-defined bounded linear operators. This is because $\varphi,\varphi'\in C^{\infty}_c(0,T)\subset L^2(0,T)$. Now \eqref{weakderivativeeq1} rewrites as:

\begin{equation}
	\mathcal{T}_{\varphi'}(\tilde{v}_{m_k})=-\mathcal{T}_{\varphi}\left (\dfrac{\partial \tilde{v}_{m_k}}{\partial t} \right), \ \forall\ k\geq 1.
\end{equation}

\noindent Using that $\tilde{v}_{m_k}\weak v$ and $\dfrac{\partial\tilde{v}_{m_k}}{\partial t}\weak w$ in $L^2\bigl(0,T;L^2(\Omega)\bigr)$ and taking into account that continuous linear operators preserve weak convergence (i.e. Lemma \ref{lemmaboundedweak} from the Appendix) we finally get that:

\begin{equation}
	\int_{0}^T v(t,\cdot)\varphi'(t)\ dt=\mathcal{T}_{\varphi'}(v)=\lim\limits_{k\to\infty} \mathcal{T}_{\varphi'}(\tilde{v}_{m_k})=-\lim\limits_{k\to\infty} \mathcal{T}_{\varphi}\left (\dfrac{\partial \tilde{v}_{m_k}}{\partial t} \right)=-\mathcal{T}_{\varphi}(w)=-\int_{0}^T w(t,\cdot)\varphi(t)\ dt.
\end{equation}

\noindent Since $\varphi\in C^\infty_{c}(0,T)$ was fixed arbitrarily we get that the $v\in L^2\bigl(0,T;L^2(\Omega)\bigr)$ has a weak derivative and $\dfrac{\partial v}{\partial t}=w\in L^2\bigl(0,T;L^2(\Omega)\bigr)$. This completely shows that $v\in H^1\bigl((0,T);L^2(\Omega)\bigr)$ and $\dfrac{\partial\tilde{v}_{m_k}}{\partial t}\weak \dfrac{\partial v}{\partial t}$ in $L^2\bigl(0,T;L^2(\Omega)\bigr)$.
\end{proof}	

\noindent\textbf{\underline{Claim 11:}} \textbf{If $u_0\in W^{1,p(x)}(\Omega)$ then $(\tilde{v}_{m_k})_{k\geq 1}\subset C([0,T];L^2(\Omega)),\ v\in C([0,T];L^2(\Omega))$ and on a further subsequence, still denoted by $(\tilde{v}_{m_k})_{k\geq 1}$, we have that $\tilde{v}_{m_k}\to v\ \text{in}\ C\bigl([0,T];L^2(\Omega)\bigr)$. In particular $\tilde{v}_{m_k}\to v$ in $L^2\bigl(0,T;L^2(\Omega)\bigr)$.}

\begin{proof} Recall that:
	
\begin{equation}
	\begin{cases} W^{1,p(x)}(\Omega)\stackrel{c}{\hookrightarrow} L^2(\Omega)\hookrightarrow L^2(\Omega),\ \text{from \textbf{(H2)}}\\ r:=2>1\\ (\tilde{v}_{m_k})_{k\geq 1}\subset L^{\infty}\bigl(0,T;W^{1,p(x)}(\Omega)\bigr),\ \text{from Proposition \ref{propovtilde} \textbf{(4)}}\\[3mm] (\tilde{v}_{m_k})_{k\geq 1}\ \text{is bounded in }L^{\infty}\bigl(0,T;W^{1,p(x)}(\Omega)\bigr),\ \text{from \textbf{Claim 7}}\\[3mm] \left(\dfrac{\partial \tilde{v}_{m_k}}{\partial t}\right)_{k\geq 1}\subset L^{2}\bigl(0,T;L^2(\Omega)\bigr),\ \text{from \textbf{Claim 4}}\\[5mm] \left(\dfrac{\partial \tilde{v}_{m_k}}{\partial t}\right)_{k\geq 1}\ \text{is bounded in }L^{2}\bigl(0,T;L^2(\Omega)\bigr),\ \text{from \textbf{Claim 4}}\end{cases}.
\end{equation}

\noindent Therefore using \textit{Aubin--Lions--Simon's compactness lemma} which can be found as Theorem \ref{lemmaaubin} in the Appendix, we obtain that $(\tilde{v}_{m_k})_{k\geq 1}\subset C([0,T];L^2(\Omega))$, and there is some function $\tilde{w}\in C([0,T];L^2(\Omega))$ such that on a further subsequence, still denoted by $(\tilde{v}_{m_k})_{k\geq 1}$, we have that $\tilde{v}_{m_k}\to \tilde{w}$ in $C([0,T];L^2(\Omega))$. We will show that $\tilde{w}=v$. Indeed, from the fact that $\tilde{v}_{m_k}\to \tilde{w}$ in $C([0,T];L^2(\Omega))$ we get that $\tilde{v}_{m_k}\to \tilde{w}$ in $L^2(0,T;L^2(\Omega))$. This is because:

\begin{align*}
	\Vert \tilde{v}_{m_k}-\tilde{w}\Vert_{L^2(0,T;L^2(\Omega))}&=\left (\int_{0}^T \Vert \tilde{v}_{m_k}(t,\cdot)-\tilde{w}(t,\cdot)\Vert_{L^2(\Omega)}^2\ dt\right )^{\frac{1}{2}}\\
	&\leq \sqrt{T}\cdot \sup_{t\in [0,T]}\Vert \tilde{v}_{m_k}(t,\cdot)-\tilde{w}(t,\cdot)\Vert_{L^2(\Omega)}\\
	&=\sqrt{T}\cdot\Vert \tilde{v}_{m_k}-\tilde{w}\Vert_{C([0,T];L^2(\Omega))}\stackrel{k\to\infty}{\longrightarrow}0.
\end{align*}

\noindent Because strong convergence in $L^2(0,T;L^2(\Omega))$ implies weak convergence in $L^2(0,T;L^2(\Omega))$, we get that $\tilde{v}_{m_k}\weak \tilde{w}$ in $L^2(0,T;L^2(\Omega))$. But, from \textbf{Claim 8} and \textbf{Claim 9} we know that $\tilde{v}_{m_k}\weak v$ in $L^2(0,T;L^2(\Omega))$. Now, from the uniqueness of the weak-limit, we conclude that $\tilde{w}=v$, and the proof of the claim is now complete.
	
\end{proof}

\noindent\textbf{\underline{Claim 12:}} \textbf{If $u_0\in W^{1,p(x)}(\Omega)$ then $v\in\mathcal{V}_{[\varepsilon,\delta]}$}.

\begin{proof} In the proof of \textbf{Claim 11} we showed that $\tilde{v}_{m_k}\to v$ in $L^2(0,T;L^2(\Omega))\simeq L^2((0,T)\times\Omega)$. Using now \cite[Corollary of Riesz-Fischer theorem, page 234]{Jones} we get on a further subsequence, still denoted by $(\tilde{v}_{m_k})_{k\geq 1}$, that:
	
	\begin{equation}\label{rieszfischereq1}
	\lim\limits_{k\to\infty} \tilde{v}_{m_k}(t,x)=v(t,x),\ \text{for a.e.}\ (t,x)\in (0,T)\times\Omega.
	\end{equation}
	
\noindent Now, from Proposition \ref{propovtilde} \textbf{(1)}, we know that $\tilde{v}_{m_k}\in \mathcal{V}_{[\varepsilon,\delta]}$ for each $k\geq 1$. Thence for a.e. $(t,x)\in (0,T)\times\Omega:\ \varepsilon\leq \tilde{v}_{m_k}(t,x)\leq \delta$. Making $k\to\infty$ in this inequality gives us (combined with \eqref{rieszfischereq1}) that $\varepsilon\leq v(t,x)\leq\delta$ for a.e. $(t,x)\in (0,T)\times\Omega$, as needed.
	
\end{proof}

\begin{remark}\label{remvsolfinal} At this point we can state that:
	
$$v\in C\bigl([0,T];L^2(\Omega)\bigr)\cap H^1\bigl((0,T);L^2(\Omega)\bigr)\cap L^{\infty}\bigl(0,T;W^{1,p(x)}(\Omega)\bigr)\cap \mathcal{V}_{[\varepsilon,\delta]}$$.
\end{remark}

\noindent\textbf{\underline{Claim 13:}} \textbf{If $u_0\in W^{1,p(x)}(\Omega)$, then $b(\cdot,\tilde{v}_{m_k})\to b(\cdot,v)$ in $C([0,T];L^2(\Omega))$}.

\begin{proof} First note that, since $v\in C([0,T];L^2(\Omega))\subset L^2(0,T;L^2(\Omega))$, from Proposition \ref{propob} \textbf{(4)} we have that $b(\cdot,v)\in L^2(0,T;L^2(\Omega))$. Moreover for any $t\in [0,T]$ we have that:
	
	\begin{align*}
		\Vert b(\cdot,v(t+h,\cdot))-b(\cdot,v(t,\cdot))\Vert_{L^2(\Omega)}&=\left (\int_{\Omega} |b(x,v(t+h,x))-b(x,v(t,x))|^2\ dx\right)^{\frac{1}{2}}\\
			\eqref{blipschitz}\ \ \ &\leq L_0\left (\int_{\Omega} |v(t+h,x)-v(t,x)|^2\ dx\right)^{\frac{1}{2}}\\
	&=L_0\Vert v(t+h,\cdot)-v(t,\cdot)\Vert_{L^2(\Omega)}\\
		&\stackrel{h\to 0}{\longrightarrow} 0,
	\end{align*}
	
	\noindent because $v\in C([0,T];L^2(\Omega))$. This proves that $b(\cdot,v)\in C([0,T];L^2(\Omega))$.
	
	\noindent Similarly, taking into account that from Proposition \ref{propovtilde} \textbf{(3)} we know that $\tilde{v}_{m_k}\in C([0,T];L^2(\Omega))$ for any $k\geq 1$, and then repeating the above argument we will get that $b(\cdot,\tilde{v}_{m_k})\in C([0,T];L^2(\Omega))$ for each $k\geq 1$. 
	
	\noindent Finally, remark that:
	
	\begin{align*}
		\Vert b(\cdot,\tilde{v}_{m_k})-b(\cdot,v)\Vert_{C([0,T];L^2(\Omega))}&=\sup_{t\in [0,T]} \Vert b(\cdot,\tilde{v}_{m_k}(t,\cdot))-b(\cdot,v(t,\cdot))\Vert_{L^2(\Omega)}\\
		&=\sup_{t\in [0,T]} \left (\int_{\Omega} |b(x,\tilde{v}_{m_k}(t,x))-b(x,v(t,x))|^2\ dx\right)^{\frac{1}{2}}\\
		\eqref{blipschitz}\ \ \ &\leq L_0\sup_{t\in [0,T]}\left (\int_{\Omega} |\tilde{v}_{m_k}(t,x)-v(t,x)|^2\ dx\right)^{\frac{1}{2}}\\
		&=L_0\sup_{t\in [0,T]} \Vert \tilde{v}_{m_k}(t,\cdot)-v(t,\cdot)\Vert_{L^2(\Omega)}\\
		&=L_0	\Vert \tilde{v}_{m_k}-v\Vert_{C([0,T];L^2(\Omega))}\stackrel{k\to\infty}{\longrightarrow}0,
	\end{align*}
	
	\noindent because, from \textbf{Claim 11} we already know that $\tilde{v}_{m_k}\to v$ in $C([0,T];L^2(\Omega))$. We conclude that $b(\cdot,\tilde{v}_{m_k})\to b(\cdot,v)$ in $C([0,T];L^2(\Omega))$.
	
\end{proof}

\noindent\textbf{\underline{Claim 14:}} \textbf{If $u_0\in W^{1,p(x)}(\Omega)$ then $b(\cdot,v_{m_k})\rightharpoonup b(\cdot,v)$ and $b(\cdot,\tilde{v}_{m_k})\rightharpoonup b(\cdot,v)$ in $L^r\big (0,T;W^{1,p(x)}(\Omega)\big )$ for each $r\in [1,\infty)$. Moreover $b(\cdot,v_{m_k})\weakstar b(\cdot,v)$ and $b(\cdot,\tilde{v}_{m_k})\weakstar b(\cdot,v)$ in $L^\infty\big (0,T;W^{1,p(x)}(\Omega)\big )$}.

\begin{proof} From \textbf{Claim 7} we have that $(b(\cdot,v_{m_k}))_{k\geq 1}, (b(\cdot,\tilde{v}_{m_k}))_{k\geq 1}\subset L^{\infty}(0,T;W^{1,p(x)}(\Omega))$. Applying the same argument as in the proof of \textbf{Claim 8} we will get that there are some $w,\tilde{w}\in L^{\infty}(0,T;W^{1,p(x)}(\Omega))$ such that: 
	
	\begin{equation}
	\begin{cases} b(\cdot,v_{m_k})\weak w\ \text{and}\ b(\cdot,\tilde{v}_{m_k})\weak \tilde{w}, \text{in}\  L^{r}(0,T;W^{1,p(x)}(\Omega))	\\ b(\cdot,v_{m_k})\weakstar w\ \text{and}\ b(\cdot,\tilde{v}_{m_k})\weakstar \tilde{w}, \text{in}\  L^{\infty}(0,T;W^{1,p(x)}(\Omega))\end{cases}.
	\end{equation}
	
	\noindent We will show next that $w=\tilde{w}=b(\cdot,v)$. Indeed, from \textbf{Claim 5}, we know that $b(\cdot,\tilde{v}_{m_k})-b(\cdot,v_{m_k})\to 0$ in $L^{\infty}(0,T;L^2(\Omega))$ But $L^{\infty}(0,T;L^2(\Omega))\hookrightarrow L^2(0,T;L^2(\Omega))$, from Proposition 2.2.5 in \cite[page 128]{gasinski2005nonlinear}. Thus $b(\cdot,\tilde{v}_{m_k})-b(\cdot,v_{m_k})\to 0$ in $L^{2}(0,T;L^2(\Omega))$. It follows that $b(\cdot,\tilde{v}_{m_k})-b(\cdot,v_{m_k})\weak 0$ in $L^{2}(0,T;L^2(\Omega))$. But $b(\cdot,\tilde{v}_{m_k})-b(\cdot,v_{m_k})\weak \tilde{w}-w$ in $L^{2}(0,T;L^2(\Omega))$. Thence $\tilde{w}=w$.
	
	\noindent Finally, from \textbf{Claim 13} we have that $b(\cdot,\tilde{v}_{m_k})\to b(\cdot,v)$ in $C([0,T];L^2(\Omega))$. Using now Lemma \ref{lemmacontlp} we get that $b(\cdot,\tilde{v}_{m_k})\to b(\cdot,v)$ in $L^2(0,T;L^2(\Omega))$, and from here $b(\cdot,\tilde{v}_{m_k})\weak b(\cdot,v)$ in $L^2(0,T;L^2(\Omega))$. But we showed above that $b(\cdot,\tilde{v}_{m_k})\weak \tilde{w}$ in $L^2(0,T;L^2(\Omega))$. From the uniqueness of the weak limit we derive that $\tilde{w}=b(\cdot,v)$. Now the proof is complete.
	
\end{proof}

\noindent\textbf{\underline{Claim 15:}} \textbf{If $u_0\in W^{1,p(x)}(\Omega)$, then $b(\cdot,v)\in H^1\bigl((0,T);L^2(\Omega)\bigr)$ and on a further subsequence, denoted for simplicity in the same way, we have that $\dfrac{\partial b(\cdot,\tilde{v}_{m_k})}{\partial t}\weak \dfrac{\partial b(\cdot,v)}{\partial t}$ in $L^2\bigl(0,T;L^2(\Omega)\bigr)$.}.

\begin{proof} From Proposition \ref{propobeta} \textbf{(5)} we know that $b(\cdot,\tilde{v}_{m_k})=\beta_{m_k}\in H^1\bigl((0,T);L^2(\Omega)\bigr)$ for each $k\geq 1$. Moreover, from \textbf{Claim 3}, we also know that the sequence $\left (\dfrac{\partial b(\cdot,\tilde{v}_{m_k})}{\partial t}\right)_{k\geq 1}$ is bounded in $L^2\bigl(0,T;L^2(\Omega)\bigr)$. Also, from \textbf{Claim 14} we have that $b(\cdot,\tilde{v}_{m_k})\weak b(\cdot,v)$ in $L^2\bigl(0,T;L^2(\Omega)\bigr)$. Now we have all the ingredients to make a verbatim repetition of the same argument given in the proof of \textbf{Claim 10} in order to get the desired conclusion.
	
\end{proof}

\begin{remark}\label{rembo} Now we can also state that:
	
	$$b(\cdot,v)\in C\bigl([0,T];L^2(\Omega)\bigr)\cap H^1\bigl((0,T);L^2(\Omega)\bigr)\cap L^{\infty}\bigl(0,T;W^{1,p(x)}(\Omega)\bigr)$$.
\end{remark}

\noindent\textbf{\underline{Claim 16:}} \textbf{If $u_0\in W^{1,p(x)}(\Omega)$ then $v_{m_k}\to v$ in $L^{\infty}\bigl(0,T;L^2(\Omega)\bigr)$.}

\begin{proof} We already know from Proposition \ref{propovm} \textbf{(2)} that $(v_{m_k})_{k\geq 1}\subset L^{\infty}\bigl(0,T;L^{2}(\Omega)\bigr)$. Also $v\in C([0,T];L^2(\Omega))\subset L^{\infty}(0,T;L^2(\Omega))$ from Lemma \ref{lemmacontlp}. 
	
\noindent Now, combining \textbf{Claim 6} and \textbf{Claim 11} we will get that:

\begin{align*}
	\Vert v_{m_k}-v\Vert_{L^{\infty}(0,T;L^2(\Omega))}&\leq \Vert v_{m_k}-\tilde{v}_{m_k}\Vert_{L^{\infty}(0,T;L^2(\Omega))}+\Vert \tilde{v}_{m_k}-v\Vert_{L^{\infty}(0,T;L^2(\Omega))}\\
	&=\Vert v_{m_k}-\tilde{v}_{m_k}\Vert_{L^{\infty}(0,T;L^2(\Omega))}+\Vert \tilde{v}_{m_k}-v\Vert_{C(0,T;L^2(\Omega))}\stackrel{k\to\infty}{\longrightarrow}0.
\end{align*}
	
\end{proof}

\bigskip

\noindent\textbf{\underline{Claim 17:}} \textbf{If $u_0\in W^{1,p(x)}(\Omega)$ then $\bigl(\mathbf{a}(\cdot,\nabla v_{m_k})\bigr)_{k\geq 1}\subset L^{\infty}\bigl(0,T;L^{p'(x)}(\Omega)^N\bigr)$ is bounded. Moreover there is a function $\boldsymbol{\eta}\in L^{\infty}\bigl(0,T;L^{p'(x)}(\Omega)^N\bigr)$ such that on a further subsequence, denoted in the same way for simplicity, the following convergence holds $\mathbf{a}(\cdot,\nabla v_{m_k})\weakstar \boldsymbol{\eta}\ \text{in}\  L^{\infty}\bigl(0,T;L^{p'(x)}(\Omega)^N\bigr)$. In particular $\mathbf{a}(\cdot,\nabla v_{m_k})\weak \boldsymbol{\eta}\ \text{in}\  L^{r}\bigl(0,T;L^{p'(x)}(\Omega)^N\bigr),$ for each $r\in [1,\infty)$.}

\begin{proof} From \textbf{Claim 7} we know that $(v_{m_k})_{k\geq 1}$ is a bounded sequence from $L^{\infty}\bigl(0,T;W^{1,p(x)}(\Omega)\bigr)$. Thus, from Proposition \ref{propoa1} \textbf{(2)} we deduce that the sequence $\bigl(\mathbf{a}(\cdot,\nabla v_{m_k})\bigr)_{k\geq 1}$ is bounded in $L^{\infty}\bigl(0,T;L^{p'(x)}(\Omega)^N\bigr)$.

\noindent Now, since $L^{p'(x)}(\Omega)^N$ is a reflexive and separable Banach space (see Proposition \ref{propospatii} \textbf{(2)}), we are in position to apply Proposition \ref{propocompactness1} for this space and for the bounded sequence $\bigl(\mathbf{a}(\cdot,\nabla v_{m_k})\bigr)_{k\geq 1}\subset L^{\infty}\bigl(0,T;L^{p'(x)}(\Omega)^N\bigr)$. Therefore there is some $\boldsymbol{\eta}\in L^{\infty}\bigl(0,T;L^{p'(x)}(\Omega)^N\bigr)$ such that on a further subsequence, denoted in the same way for simplicity, we have that: $\mathbf{a}(\cdot,\nabla v_{m_k})\weakstar \boldsymbol{\eta}\ \text{in}\  L^{\infty}\bigl(0,T;L^{p'(x)}(\Omega)^N\bigr)$ and $\mathbf{a}(\cdot,\nabla v_{m_k})\weak \boldsymbol{\eta}\ \text{in}\  L^{r}\bigl(0,T;L^{p'(x)}(\Omega)^N\bigr),$ for each $r\in [1,\infty)$.
	
\end{proof}

\begin{itemize}
	\item \textbf{Step V:} We define for each $m\geq 0$ the function $g_{m}:(0,T]\times\Omega\to\mathbb{R}$ by
\end{itemize}

\begin{equation}
	g_m(t,x)=\sum_{n=1}^{2^m} g_{m,n}(x)\chi_{(t_{m,n-1},t_{m,n}]}(t)=\sum_{n=1}^{2^m}\left [\dfrac{1}{\Delta_m}\int_{t_{m,n-1}}^{t_{m,n}}g(\tau,x)\ d\tau\right ]\chi_{(t_{m,n-1},t_{m,n}]}(t).
\end{equation}

\begin{proposition}\label{propogm} For each $m\in\mathbb{N}$, the following properties of $g_m$ hold:
	
	\begin{enumerate}
		\item[\textnormal{\textbf{(1)}}] $g_m\in\mathfrak{M}_{\lambda}^{[\varepsilon,\delta]}$.
		
		\item[\textnormal{\textbf{(2)}}] $g_m\in L^{\infty}\big (0,T;L^{\infty}(\Omega)\big )$ and $\Vert g_m\Vert_{L^{\infty}(0,T;L^{\infty}(\Omega))}\leq \Vert g\Vert_{L^{\infty}(0,T;L^{\infty}(\Omega))}$.
	\end{enumerate}

\end{proposition}

\begin{proof} Fix some arbitrary $t\in (0,T]$. There is a unique $n\in\{1,2,\hdots,2^m\}$ such that $t\in (t_{m,n-1},t_{m,n}]$. 
	
\noindent\textbf{(1)} Therefore $g_m(t,x)=g_{m,n}(x)=\dfrac{1}{\Delta_m}\displaystyle\int_{t_{m,n-1}}^{t_{m,n}} g(\tau,x)\ d\tau\in [\lambda b(x,\varepsilon),\lambda b(x,\delta)]$ for a.e. $x\in\Omega$, because $g\in \mathfrak{M}_{\lambda}^{[\varepsilon,\delta]}$. Thus $g_m\in \mathfrak{M}_{\lambda}^{[\varepsilon,\delta]}$.
	
\noindent\textbf{(2)} Similarly $g_m(t,x)=g_{m,n}(x)=\dfrac{1}{\Delta_m}\displaystyle\int_{t_{m,n-1}}^{t_{m,n}} g(\tau,x)\ d\tau\leq \dfrac{1}{\Delta_m}\displaystyle\int_{t_{m,n-1}}^{t_{m,n}} \Vert g\Vert_{L^{\infty}(0,T;L^{\infty}(\Omega))}\ d\tau=\Vert g\Vert_{L^{\infty}(0,T;L^{\infty}(\Omega))}$ for any $t\in (0,T]$ and for a.e. $x\in\Omega$. The conclusion follows.
	
\end{proof}

\begin{lemma}\label{lemmagm}
	For any $r\in [1,\infty)$ we have for a.e. $t\in (0,T)$ that
	
	\begin{equation}
		\lim\limits_{m\to\infty} g_m(t,\cdot)= g(t,\cdot)\ \text{in}\ L^r(\Omega).
	\end{equation}
	
	\noindent Moreover one has that $\lim\limits_{m\to\infty} g_m= g\ \text{in}\ L^r\big ((0,T)\times\Omega\big )$.
\end{lemma}

\begin{proof} Consider the extended function $\overline{g}:\mathbb{R}\times\Omega\to\mathbb{R},\ \overline{g}(t,x)=\begin{cases}\ g(t,x), & t\in (0,T) \\[2mm] 0, & t\notin (0,T)\end{cases}$. Fix some $t\in (0,T]$. First note that $g_m(t,\cdot), g(t,\cdot)\in L^{\infty}(\Omega)\subset L^r(\Omega)$. It is obvious to deduce that $\overline{g}(t,\cdot)\in L^{\infty}(\Omega)\subset L^r(\Omega)$. Moreover $\overline{g}\in L^\infty (\mathbb{R}\times\Omega)$ and $\overline{g}\in L^1(\mathbb{R}\times\Omega)$. There is a unique $n\in\{1,2,\hdots,2^m\}$ such that $t\in (t_{m,n-1},t_{m,n}]$ and therefore we may write for a.e. $x\in\Omega$ that
	\begin{align*}
		|g_m(t,x)-g(t,x)|&=\left |\dfrac{1}{\Delta_m}\int_{t_{m,n-1}}^{t_{m,n}}g(\tau,x)\ d\tau-g(t,x) \right |=\left |\dfrac{1}{\Delta_m}\int_{t_{m,n-1}}^{t_{m,n}}g(\tau,x)-g(t,x)\ d\tau \right |\\
		&\leq \dfrac{1}{\Delta_m}\int_{t_{m,n-1}}^{t_{m,n}} \big |g(\tau,x)-g(t,x)\big |\ d\tau\\
		&\leq  \dfrac{1}{\Delta_m}\int_{t-\Delta_m}^{t+\Delta_m} \big |\overline{g}(\tau,x)-\overline{g}(t,x)\big |\ d\tau=2\cdot \dfrac{1}{2\Delta_m}\int_{t-\Delta_m}^{t+\Delta_m} \big |\overline{g}(\tau,x)-\overline{g}(t,x)\big |\ d\tau
	\end{align*}

\noindent Since $\overline{g}(\cdot,x)\in L^1(\mathbb{R})$ we deduce from \textit{Lebesgue's theorem}\footnote{See \cite[page 456]{Jones} for statement and complete proof.} that for a.e. $t\in\mathbb{R}$, and hence for a.e. $t\in (0,T)$, that $\lim\limits_{\Delta_m\to 0}\dfrac{1}{2\Delta_m}\displaystyle\int_{t-\Delta_m}^{t+\Delta_m} \big |\overline{g}(\tau,x)-\overline{g}(t,x)\big |\ d\tau=0$, i.e. $\lim\limits_{m\to \infty} \dfrac{1}{2\Delta_m}\displaystyle\int_{t-\Delta_m}^{t+\Delta_m} \big |\overline{g}(\tau,x)-\overline{g}(t,x)\big |\ d\tau=0$, because $\Delta_m=\dfrac{T}{2^m}$. This shows that:

\begin{equation}\label{convergencelimit1}
	\text{for a.e.}\ x\in\Omega:\ \lim\limits_{m\to\infty} |g_m(t,x)-g(t,x)|=0,\ \text{for a.e.}\ t\in (0,T).
\end{equation}
	
\noindent Consider the following set: $E:=\bigl\{(t,x)\in(0,T)\times\Omega\ |\
g_m(t,x)\not\longrightarrow g(t,x)\bigr\}$. The set $E$ is measurable. Indeed,
\begin{equation}
	E=\bigcup_{k=1}^{\infty}\left (\bigcap_{N=1}^{\infty}\bigcup_{m\ge N} \left\{(t,x)\ |\ |g_m(t,x)-g(t,x)|>\frac1k\right\}\right ),
\end{equation}

\noindent which is obtained from measurable sets (from the definition of the fact that $g_m-g$ is a measurable function) by countable unions and intersections. Therefore the characteristic function $\chi_{E}:(0,T)\times\Omega\to \{0,1\}$ is measurable (and positive). Therefore we can apply \textit{Tonelli's Theorem} and obtain that:

\begin{equation}\label{fubini1}
	\int_{(0,T)\times\Omega} \chi_{E}(t,x)\ dt\ dx=\int_{(0,T)}\left (\int_{\Omega} \chi_{E}(t,x)\ dx \right )\ dt=\int_{\Omega}\left (\int_{(0,T)} \chi_{E}(t,x)\ dt \right )\ dx.
\end{equation}

\noindent For any $x\in\Omega$ we denote the time-section $E_x:=\bigl\{t\in (0,T)\ |\ (t,x)\in E\bigr\}$. This set is measurable for a.e. $x\in\Omega$\footnote{See the Theorem given at page 270 in \cite{Jones}.}. Similarly for any $t\in (0,T)$ we denote the space section $E^t:=\bigl\{x\in\Omega\ |\ (t,x)\in E\bigr \}$. This set is also measurable for a.e. $t\in (0,T)$.\footnote{The idea is that the product of the Lebesgue $\sigma$-algebras of $\mathbb{R}$ and $\mathbb{R}^N$ is not the $\sigma$-algebra of $\mathbb{R}^{N+1}$, because the product of two complete $\sigma$-algebras need not be complete.} 

\noindent From \eqref{convergencelimit1} we get that for a.e. $x\in\Omega$, $|E_x|=0$. Therefore \eqref{fubini1} becomes:

\begin{equation}
	|E|=\int_{0}^T |E^t|\ dt=\int_{\Omega} |E_x|\ dx=0.
\end{equation}

\noindent Since $|E|=0$ it follows that for a.e. $t\in (0,T)$: $|E^t|=0$, i.e.

\begin{equation}
	\text{for a.e.}\ t\in (0,T):\ \lim\limits_{m\to\infty} |g_m(t,x)-g(t,x)|=0,\ \text{for a.e.}\ x\in\Omega.
\end{equation}

\noindent Thus, for a.e. $t\in (0,T)$: $g_m(t,\cdot)\to g(t,\cdot)$ pointwise a.e. on $\Omega$. From Proposition \ref{propogm} \textbf{(2)} we get for a.e. $x\in\Omega$ that:

\begin{equation}
	|g_m(t,x)-g(t,x)|^r\leq \bigl(|g_m(t,x)|+|g(t,x)|\bigr)^r\leq 2^r\Vert g\Vert_{L^{\infty}((0,T)\times\Omega)}^r\in L^1(\Omega).
\end{equation}

\noindent Using now \textit{Lebesgue dominated convergence theorem} we conclude that

\begin{equation}
	\lim\limits_{m\to\infty} \Vert g_m(t,\cdot)-g(t,\cdot)\Vert_{L^r(\Omega)}^r=\lim\limits_{m\to\infty}\int_{\Omega} |g_m(t,x)-g(t,x)|^r\ dx=0,
\end{equation}

\noindent and thence $g_m(t,\cdot)\to g(t,\cdot)$ in $L^r(\Omega)$ for a.e. $t\in (0,T)$, as needed.

\noindent In an analog fashion, since from $|E|=0$ we get that $g_m\to g$ pointwise a.e. on $(0,T)\times\Omega$ and $|g_m-g|^r\leq 2^r\Vert g\Vert_{L^{\infty}((0,T)\times\Omega)}^r\in L^1\bigl((0,T)\times\Omega\bigr)$, applying  \textit{Lebesgue dominated convergence theorem} we conclude that $g_m\to g$ in $L^r\bigl ((0,T)\times\Omega\bigr)$.

\end{proof}

\begin{theorem}\label{thmauxiliar} Problem \eqref{eqdpgaux} has exactly one weak solution $u$ for any regular initial data $u_0\in W^{1,p(x)}(\Omega)\cap\mathcal{U}_{[\varepsilon,\delta]}$. Moreover $u\in C\bigl([0,T];L^2(\Omega)\bigr) \cap H^1\bigl((0,T); L^2(\Omega)\bigr)\cap L^{\infty}\bigl(0,T;W^{1,p(x)}(\Omega)\bigr)\cap \mathcal{V}_{[\varepsilon,\delta]}$.
\end{theorem}

\begin{proof} Notice that, for each $m\in\mathbb{N}$, the above constructed functions $v_m$ and $\tilde{v}_m$ satisfy in the weak sense:
	
	\begin{equation}\label{proof1}
		\begin{cases}\dfrac{\partial b\big (x,\tilde{v}_m\big )}{\partial t}-\operatorname{div}\mathbf{a}(x,\nabla v_m)+\lambda b(x,v_m)=g_m(t,x), & (t,x)\in (0,T)\times\Omega\\[3mm] \mathbf{a}(x,\nabla v_m)\cdot\nu=0, & (t,x)\in (0,T)\times\partial\Omega\\[3mm] \tilde{v}_m(0,x)=v_m(0,x)=u_0(x), & x\in\Omega \end{cases},
	\end{equation}
	
\noindent meaning that for a.e. $t\in (0,T)$ and for each test function $\phi\in W^{1,p(x)}(\Omega)$ we have that:

\begin{equation}\label{important1}
	\int_\Omega\dfrac{\partial b(x,\tilde{v}_m(t,x))}{\partial t}\phi\ dx+\int_{\Omega}\mathbf{a}(x,\nabla v_m(t,x))\cdot\nabla\phi\ dx+\lambda\int_{\Omega}b(x,v_m(t,x))\phi\ dx=\int_{\Omega} g_m(t,x)\phi\ dx.
\end{equation}

\noindent This equality makes sense because of the following facts: from Proposition \ref{propobeta} \textbf{(5)} we know that $\dfrac{\partial b(\cdot,\tilde{v}_m(t,\cdot))}{\partial t}\in L^{\infty}(\Omega)$, from hypothesis \textbf{(H6)} and Claim 2 we get that $\mathbf{a}(\cdot,\nabla v_m(t,\cdot))\in L^{p'(x)}(\Omega)^N$, from Proposition \ref{propovm} \textbf{(5)} it follows that $b(\cdot,v_m(t,\cdot))\in L^{\infty}(\Omega)$, from Proposition \ref{propogm} \textbf{(2)} we have that $g_m(t,\cdot)\in L^{\infty}(\Omega)$ and finally from $\phi\in W^{1,p(x)}(\Omega)\subset L^{p(x)}(\Omega)$ and $\nabla\phi\in L^{p(x)}(\Omega)^N$.

\noindent Now this equality is true because if we fix any $t\in (0,T]$ we get that there is a unique $n\in\{1,2,\hdots, 2^m\}$ such that $t\in (t_{m,n-1},t_{m,n}]$, and this equality follows directly from \eqref{eq1112}. 

\begin{remark}\label{remimp1} Note that the left-hand side and the right-hand side of \eqref{important1} represent measurable simple functions defined on the interval $(0,T)$. Therefore they are Lebesgue integrable on any interval $(t_1,t_2)\subset (0,T)$. Indeed:
	
	\begin{equation}
		\begin{cases} \displaystyle\int_\Omega\dfrac{\partial b(x,\tilde{v}_m(t,x))}{\partial t}\phi(x)\ dx=\sum_{n=1}^{2^m}\left (\int_{\Omega}\dfrac{b(x,v_{m,n}(x))-b(x,v_{m,n-1}(x))}{\Delta_m}\phi(x)\ dx \right )\chi_{(t_{m,n-1},t_{m,n}]}(t)\\[3mm]
		\displaystyle\int_{\Omega}\mathbf{a}(x,\nabla v_m(t,x))\cdot\nabla\phi(x)\ dx=\sum_{n=1}^{2^m} \left (\int_{\Omega}\mathbf{a}(x,\nabla v_{m,n}(x))\cdot\nabla\phi(x)\ dx \right )\chi_{(t_{m,n-1},t_{m,n}]}(t)\\[3mm] 
		\displaystyle\int_{\Omega}b(x,v_m(t,x))\phi(x)\ dx=\sum_{n=1}^{2^m}\left (\int_{\Omega}b(x, v_{m,n}(x))\phi(x)\ dx \right )\chi_{(t_{m,n-1},t_{m,n}]}(t)\\[3mm]
		\displaystyle\int_{\Omega}g_m(t,x)\phi(x)\ dx=\sum_{n=1}^{2^m}\left [\int_{\Omega} \left(\dfrac{1}{\Delta_m}\int_{t_{m,n-1}}^{t_{m,n}} g(t,x)\ dt \right)\phi(x)\ dx\right ]\chi_{(t_{m,n-1},t_{m,n}]}(t).
		\end{cases}
	\end{equation}
\end{remark}

\noindent Now we need the following lemma:

\begin{lemma}\label{lemmatotalintegral}
	For any $(t_1,t_2)\subseteq (0,T)$, for any $\Theta\in L^{2}\bigl(0,T;W^{1,p(x)}(\Omega)\bigr)$ and for each $m\in\mathbb{N}$ we have that:
	
	\begin{equation}\label{important2}
		\int_{t_1}^{t_2}\int_\Omega\dfrac{\partial b(x,\tilde{v}_{m})}{\partial t}\Theta\ dx\ dt+\int_{t_1}^{t_2}\int_{\Omega}\mathbf{a}(x,\nabla v_{m})\cdot\nabla\Theta\ dx\ dt+\lambda\int_{t_1}^{t_2}\int_{\Omega}b(x,v_{m})\Theta\ dx\ dt=\int_{t_1}^{t_2}\int_{\Omega} g_{m}\Theta\ dx\ dt.
	\end{equation}
	
\end{lemma}

\begin{proof} Fix any $\Theta\in L^{2}\bigl(0,T;W^{1,p(x)}(\Omega)\bigr)$. Fix also some $t\in (0,T)$ for which \eqref{important1} holds. By setting $\phi:=\Theta(t,\cdot)\in W^{1,p(x)}(\Omega)$ we obtain that:
	
\begin{align}\label{important3}
	&\int_\Omega\dfrac{\partial b(x,\tilde{v}_m(t,x))}{\partial t}\Theta(t,x)\ dx+\int_{\Omega}\mathbf{a}(x,\nabla v_m(t,x))\cdot\nabla\Theta(t,x)\ dx+\lambda\int_{\Omega}b(x,v_m(t,x))\Theta(t,x)\ dx \nonumber\\
	=&\int_{\Omega} g_m(t,x)\Theta(t,x)\ dx.
\end{align}
	
\noindent Note that \eqref{important3} holds for a.e. $t\in (0,T)$. The next step is to explain why all four terms of this equality are functions from $L^1(0,T)$, so that we may integrate \eqref{important3} on any subinterval $(t_1,t_2)\subseteq (0,T)$. 

\bigskip

\noindent $\bullet$ Indeed, from Proposition \ref{propobeta} \textbf{(5)}:

\begin{equation}
\begin{cases} \dfrac{\partial b(\cdot,\tilde{v}_m(\cdot,\cdot))}{\partial t}\in L^2\bigl(0,T;L^2(\Omega)\bigr)\simeq L^2\bigl((0,T)\times\Omega\bigr)\\ \Theta\in L^{2}\bigl(0,T;W^{1,p(x)}(\Omega)\bigr)\subset L^{2}\bigl(0,T;L^{2}(\Omega)\bigr)\simeq L^2\bigl((0,T)\times\Omega\bigr).  \end{cases}.
\end{equation}

\noindent Then, from \textit{Cauchy inequality} we get that $\dfrac{\partial b(\cdot,\tilde{v}_m(\cdot,\cdot))}{\partial t}\Theta\in L^1\bigl((0,T)\times\Omega\bigr)$. Now using \textit{Fubini's Theorem}\footnote{For statement and proof see \cite[Theorem 2.3.50, page 126]{papageorgiou2018applied}.} we conclude that the function $(0,T)\ni t\mapsto \displaystyle\int_\Omega\dfrac{\partial b(x,\tilde{v}_m(t,x))}{\partial t}\Theta(t,x)\ dx$ is from $L^1(0,T)$.

\bigskip

\noindent $\bullet$ Similarly, from Proposition \ref{propovm} \textbf{(5)}:

\begin{equation}
	\begin{cases} b(\cdot,v_m(\cdot,\cdot))\in L^2\bigl(0,T;L^2(\Omega)\bigr)\simeq L^2\bigl((0,T)\times\Omega\bigr)\\ \Theta\in L^{2}\bigl(0,T;W^{1,p(x)}(\Omega)\bigr)\subset L^{2}\bigl(0,T;L^{2}(\Omega)\bigr)\simeq L^2\bigl((0,T)\times\Omega\bigr).  \end{cases}.
\end{equation}

\noindent Then, from \textit{Cauchy inequality} we get that $b(\cdot,v_m(\cdot,\cdot))\Theta\in L^1\bigl((0,T)\times\Omega\bigr)$. Now using \textit{Fubini's Theorem} we conclude that the function $(0,T)\ni t\mapsto \displaystyle\int_\Omega b(x,v_m(t,x))\Theta(t,x)\ dx$ is from $L^1(0,T)$.

\bigskip

\noindent $\bullet$ Again, from Lemma \ref{lemmagm}:

\begin{equation}
	\begin{cases} g_m\in L^2\bigl((0,T)\times\Omega\bigr)\\ \Theta\in L^{2}\bigl(0,T;W^{1,p(x)}(\Omega)\bigr)\subset L^{2}\bigl(0,T;L^{2}(\Omega)\bigr)\simeq L^2\bigl((0,T)\times\Omega\bigr).  \end{cases}.
\end{equation}

\noindent Then, from \textit{Cauchy inequality} we get that $g_m\Theta\in L^1\bigl((0,T)\times\Omega\bigr)$. Now using \textit{Fubini's Theorem} we conclude that the function $(0,T)\ni t\mapsto \displaystyle\int_\Omega g_m(t,x)\Theta(t,x)\ dx$ is from $L^1(0,T)$.

\bigskip

\noindent $\bullet$ Since $v_m\in L^\infty\bigl(0,T;W^{1,p(x)}(\Omega)\bigr)$ (see Proposition \ref{propovm} \textbf{(1)}) we deduce from Proposition \ref{propoa1} \textbf{(1)} that $\mathbf{a}(\cdot,\nabla v_m(\cdot,\cdot))\in L^{\infty}\bigl(0,T;L^{p'(x)}(\Omega)^N\bigr)$.

\noindent Also, $\Theta\in L^2\bigl(0,T;W^{1,p(x)}(\Omega)\bigr)$ and then we have from Proposition \ref{propotheta} that $\nabla \Theta\in L^2\bigl(0,T;L^{p(x)}(\Omega)^N\bigr)$. Therefore, from Proposition \ref{propolplprim} we deduce that the function $(0,T)\times\Omega\ni (t,x)\mapsto h(t,x):=\mathbf{a}(x,\nabla v_m(t,x))\cdot\nabla \Theta(t,x)\in\mathbb{R}$ is in $L^2\bigl(0,T;L^1(\Omega)\bigr)\subset L^1\bigl(0,T;L^1(\Omega)\bigr)=L^1\bigl((0,T)\times\Omega\bigr)$. \textit{Fubini's theorem} applied for function $h$ gives us that the function $(0,T)\ni t\mapsto\displaystyle\int_{\Omega} h(t,x)\ dx$ is from $L^1(0,T)$.

\bigskip

\noindent Finally, integrating \eqref{important3} on $(t_1,t_2)\subseteq (0,T)$ gives us \eqref{important2}. The proof is complete.

\end{proof}

\noindent $\blacktriangleright$ From \textbf{Claim 15} we know that $\dfrac{\partial b(\cdot,\tilde{v}_{m_k})}{\partial t}\weak \dfrac{\partial b(\cdot,v)}{\partial t}$ in $L^2\bigl(0,T;L^2(\Omega)\bigr)$. This means that for every $\psi\in L^2(0,T;L^2(\Omega))$ one has that:

\begin{equation}
	\lim\limits_{k\to\infty}\int_{0}^T\int_{\Omega} \dfrac{\partial b(\cdot,\tilde{v}_{m_k})}{\partial t}\psi(t,x)\ dx\ dt=\int_{0}^T\int_{\Omega} \dfrac{\partial b(\cdot,v)}{\partial t}\psi(t,x)\ dx\ dt.
\end{equation}

\noindent Therefore, taking into account that $\Theta\in L^2\bigl(0,T;W^{1,p(x)}(\Omega)\bigr)\subset L^2\bigl(0,T;L^2(\Omega)\bigr)$\footnote{This follows from Proposition 2.2.5 given in \cite[page 128]{gasinski2005nonlinear}, using \textbf{(H2)}.} and choosing in the above equality $\psi(t,x)=\chi_{(t_1,t_2)}(t)\Theta(t,x)\in L^2(0,T;L^2(\Omega))$, we conclude that:

\begin{equation}\label{i3}
	\lim\limits_{k\to\infty}\int_{t_1}^{t_2}\int_{\Omega} \dfrac{\partial b(\cdot,\tilde{v}_{m_k})}{\partial t}\Theta(t,x)\ dx\ dt=\int_{t_1}^{t_2}\int_{\Omega} \dfrac{\partial b(\cdot,v)}{\partial t}\Theta(t,x)\ dx\ dt.
\end{equation}

\bigskip

\noindent $\blacktriangleright$ From \textbf{Claim 14} we have that $b(\cdot,v_{m_k})\weak b(\cdot,v)$ in $L^2(0,T;W^{1,p(x)}(\Omega))$ and from \textbf{(H2)} we know that $W^{1,p(x)}(\Omega)\hookrightarrow L^2(\Omega)$. Therefore, using Lemma \ref{lemmaweakbochner}, we get that $b(\cdot,v_{m_k})\weak b(\cdot,v)$ in $L^2(0,T;L^2(\Omega))$. This means that for every $\psi\in L^2(0,T;L^2(\Omega))$ one has that:

\begin{equation}
	\lim\limits_{k\to\infty}\int_{0}^T\int_{\Omega} b(x,v_{m_k}(t,x))\psi(t,x)\ dx\ dt=\int_{0}^T\int_{\Omega} b(x,v(t,x))\psi(t,x)\ dx\ dt.
\end{equation}

\noindent Therefore, taking into account that $\Theta\in L^2\bigl(0,T;W^{1,p(x)}(\Omega)\bigr)\subset L^2\bigl(0,T;L^2(\Omega)\bigr)$ and choosing in the above equality $\psi(t,x)=\chi_{(t_1,t_2)}(t)\Theta(t,x)\in L^2(0,T;L^2(\Omega))$, we conclude that:

\begin{equation}\label{i2}
	\lim\limits_{k\to\infty} \int_{t_1}^{t_2}\int_{\Omega}  b(x,v_{m_k}(t,x))\Theta(t,x)\ dx\ dt=\int_{t_1}^{t_2}\int_{\Omega}  b(x,v(t,x))\Theta(x)\ dx\ dt.
\end{equation}

\bigskip

\noindent $\blacktriangleright$ From Lemma \ref{lemmagm} we know that $g_{m_k}\to g$ in $L^2(0,T;L^2(\Omega))$. Since the strong convergence implies the weak convergence, we also get that $g_{m_k}\rightharpoonup g$ in $L^2(0,T;L^2(\Omega))$. This means that for every $\psi\in L^2(0,T;L^2(\Omega))$ one has that:

\begin{equation}
	\lim\limits_{k\to\infty}\int_{0}^T\int_{\Omega} g_{m_k}(t,x)\psi(t,x)\ dx\ dt=\int_{0}^T\int_{\Omega} g(t,x)\psi(t,x)\ dx\ dt.
\end{equation}

\noindent Therefore, taking into account that $\Theta\in L^2\bigl(0,T;W^{1,p(x)}(\Omega)\bigr)\subset L^2\bigl(0,T;L^2(\Omega)\bigr)$ and choosing in the above equality $\psi(t,x)=\chi_{(t_1,t_2)}(t)\Theta(t,x)\in L^2(0,T;L^2(\Omega))$, we conclude that

\begin{equation}\label{i4}
	\lim\limits_{k\to\infty} \int_{t_1}^{t_2}\int_{\Omega} g_{m_k}(t,x)\Theta(t,x)\ dx\ dt=\int_{t_1}^{t_2}\int_{\Omega} g(t,x)\Theta(t,x)\ dx\ dt.
\end{equation}

\bigskip

\noindent From Lemma \ref{lemmatotalintegral} we get that for any $\Theta\in L^2\bigl (0,T;W^{1,p(x)}(\Omega)\bigr)$, any subinterval $(t_1,t_2)\subseteq (0,T)$ and any $k\geq 1$ we have that:

\begin{equation}\label{i1}
	\int_{t_1}^{t_2}\int_\Omega\dfrac{\partial b(x,\tilde{v}_{m_k})}{\partial t}\Theta\ dx\ dt+\int_{t_1}^{t_2}\int_{\Omega}\mathbf{a}(x,\nabla v_{m_k})\cdot\nabla\Theta\ dx\ dt+\lambda\int_{t_1}^{t_2}\int_{\Omega}b(x,v_{m_k})\Theta\ dx\ dt=\int_{t_1}^{t_2}\int_{\Omega} g_{m_k}\Theta\ dx\ dt.
\end{equation}

\noindent At this point, making $k\to\infty$ in \eqref{i1}, and using \eqref{i3},\eqref{i2} and \eqref{i4}, gives us that:

\begin{align}
	&\lim\limits_{k\to\infty} \int_{t_1}^{t_2}\int_{\Omega}\mathbf{a}(x,\nabla v_{m_k})\cdot\nabla\Theta\ dx\ dt=\nonumber\\
	=&\int_{t_1}^{t_2}\int_{\Omega} g\Theta\ dx\ dt-\int_{t_1}^{t_2}\int_\Omega\dfrac{\partial b(x,v)}{\partial t}\Theta\ dx\ dt-\lambda\int_{t_1}^{t_2}\int_{\Omega}b(x,v)\Theta\ dx\ dt.\label{i5}
\end{align}

\noindent Now set $\Theta=v_{m_k}-v\in L^{\infty}\bigl(0,T;W^{1,p(x)}(\Omega)\bigr)\subset L^2\bigl(0,T;W^{1,p(x)}(\Omega)\bigr)$ in \eqref{i1} and obtain that:

\begin{align}\label{i6}
	&\int_{t_1}^{t_2}\int_\Omega\dfrac{\partial b(x,\tilde{v}_{m_k})}{\partial t}(v_{m_k}-v)\ dx\ dt+\int_{t_1}^{t_2}\int_{\Omega}\mathbf{a}(x,\nabla v_{m_k})\cdot(\nabla v_{m_k}-\nabla v)\ dx\ dt\nonumber\\
	+&\lambda\int_{t_1}^{t_2}\int_{\Omega}b(x,v_{m_k})(v_{m_k}-v)\ dx\ dt=\int_{t_1}^{t_2}\int_{\Omega} g_{m_k}(v_{m_k}-v)\ dx\ dt.
\end{align}

\noindent From \textbf{Claim 16} we know that $v_{m_k}\to v$ in $L^{\infty}\bigl(0,T;L^2(\Omega)\bigr)\hookrightarrow L^2\bigl(0,T;L^2(\Omega)\bigr)$. Thence:

\bigskip

\noindent$\blacktriangleright$ Using \textbf{Claim 3} and \textit{Cauchy inequality} we get:

\begin{align}\label{i7}
	&\left |\int_{t_1}^{t_2}\int_\Omega\dfrac{\partial b(x,\tilde{v}_{m_k})}{\partial t}(v_{m_k}-v)\ dx\ dt\right |\leq \left\Vert\dfrac{\partial b(\cdot,\tilde{v}_{m_k}(\cdot,\cdot))}{\partial t} \right\Vert_{L^2(0,T;L^2(\Omega))}\cdot \Vert v_{m_k}-v\Vert_{L^2(0,T;L^2(\Omega))}\nonumber \\
	&\leq \underbrace{\sup_{k\geq 1}\left\Vert\dfrac{\partial b(\cdot,\tilde{v}_{m_k}(\cdot,\cdot))}{\partial t} \right\Vert_{L^2(0,T;L^2(\Omega))}}_{<\infty}\cdot \Vert v_{m_k}-v\Vert_{L^2(0,T;L^2(\Omega))}\stackrel{k\to\infty}{\longrightarrow}0.
\end{align}

\noindent$\blacktriangleright$ From \textbf{Claim 7} we get that $\bigl(b(\cdot,v_{m_k})\bigr)_{k\geq 1}$ is bounded in $L^{\infty}\bigl(0,T;W^{1,p(x)}(\Omega)\bigr)\hookrightarrow L^2\bigl(0,T;W^{1,p(x)}(\Omega)\bigr)\hookrightarrow L^2\bigl(0,T;L^2(\Omega)\bigr)$. So $\bigl(b(\cdot,v_{m_k})\bigr)_{k\geq 1}$ is also bounded in $L^2\bigl(0,T;L^2(\Omega)\bigr)$.  Using \textit{Cauchy inequality} we get:

\begin{align}\label{i8}
	&\left |\int_{t_1}^{t_2}\int_\Omega b(x,v_{m_k})(v_{m_k}-v)\ dx\ dt\right |\leq \left\Vert b(\cdot,v_{m_k}(\cdot,\cdot)) \right\Vert_{L^2(0,T;L^2(\Omega))}\cdot \Vert v_{m_k}-v\Vert_{L^2(0,T;L^2(\Omega))}\nonumber\\
	&\leq \underbrace{\sup_{k\geq 1}\left\Vert b(\cdot,v_{m_k}(\cdot,\cdot)) \right\Vert_{L^2(0,T;L^2(\Omega))}}_{<\infty}\cdot \Vert v_{m_k}-v\Vert_{L^2(0,T;L^2(\Omega))}\stackrel{k\to\infty}{\longrightarrow}0.
\end{align}

\noindent$\blacktriangleright$ From Proposition \ref{propogm} \textbf{(2)} we get that:

\begin{align}\label{i9}
	&\left |\int_{t_1}^{t_2}\int_\Omega g_m(v_{m_k}-v)\ dx\ dt\right |\leq \int_{t_1}^{t_2}\int_\Omega \left |g_m\right |\cdot \left |v_{m_k}-v\right |\ dx\ dt\nonumber\\
	&\leq \int_{t_1}^{t_2}\int_\Omega \left \Vert g\right \Vert_{L^{\infty}(0,T;L^{\infty}(\Omega))}\cdot \left |v_{m_k}-v\right |\ dx\ dt\nonumber\\
	&\leq \left \Vert g\right \Vert_{L^{\infty}(0,T;L^{\infty}(\Omega))}\sqrt{(t_2-t_1)|\Omega|}\cdot \Vert v_{m_k}-v\Vert_{L^2(0,T;L^2(\Omega))}\stackrel{k\to\infty}{\longrightarrow}0.
\end{align}

\noindent Thus making $k\to\infty$ in \eqref{i6} and using \eqref{i7}, \eqref{i8} and \eqref{i9} we get that

\begin{equation}\label{i10}
	\lim\limits_{k\to\infty} \int_{t_1}^{t_2}\int_{\Omega}\mathbf{a}(x,\nabla v_{m_k})\cdot(\nabla v_{m_k}-\nabla v)\ dx\ dt=0.
\end{equation}

\noindent Applying now Proposition \ref{propoa2} for $v\in L^{\infty}\bigl(0,T;W^{1,p(x)}(\Omega)\bigr)$ we get that $\ell_v:L^2\bigl(0,T;W^{1,p(x)}(\Omega)\bigr)\to\mathbb{R},\ \ell_v(w)=\displaystyle\int_{t_1}^{t_2}\int_{\Omega}\mathbf{a}(x,\nabla v)\cdot\nabla w\ dx\ dt$ is a bounded linear operator, i.e. $\ell\in L^2\bigl(0,T;W^{1,p(x)}(\Omega)\bigr)^*$.

\noindent But, from \textbf{Claim 8} we know that $v_{m_k}\weak v$ in $L^2\bigl(0,T;W^{1,p(x)}(\Omega)\bigr)$. Using the fact that continuous linear operators preserve weak convergence (see Lemma \ref{lemmaboundedweak} from the Appendix) we get that $\ell(v_{m_k})\weak \ell(v)$ in $\mathbb{R}$ -- which is a finite dimensional normed space, in which weak and strong convergence coincide. So we have proved that $\ell(v_{m_k})\to \ell(v)$, i.e.

\begin{equation}\label{i11}
	\lim\limits_{k\to\infty} \int_{t_1}^{t_2}\int_{\Omega}\mathbf{a}(x,\nabla v)\cdot\nabla v_{m_k}\ dx\ dt=\int_{t_1}^{t_2}\int_{\Omega}\mathbf{a}(x,\nabla v)\cdot\nabla v\ dx\ dt.
\end{equation}

\noindent Combining \eqref{i10} and \eqref{i11} results in:

\begin{equation}\label{i12}
	\lim\limits_{k\to\infty} \int_{t_1}^{t_2}\int_{\Omega}\bigl(\mathbf{a}(x,\nabla v_{m_k})-\mathbf{a}(x,\nabla v_{m_k})\bigr)\cdot(\nabla v_{m_k}-\nabla v)\ dx\ dt=0.
\end{equation}

\noindent If we set for the moment $t_1=0$ and $t_2=T$ we obtain that:

\begin{equation}\label{i13}
	\lim\limits_{k\to\infty} \int_{0}^{T}\int_{\Omega}\bigl(\mathbf{a}(x,\nabla v_{m_k})-\mathbf{a}(x,\nabla v_{m_k})\bigr)\cdot(\nabla v_{m_k}-\nabla v)\ dx\ dt=0.
\end{equation}

\bigskip 

\noindent It is now the time to dig deeper into the meanings of \eqref{i5}. From \textbf{Claim 17} we know that $\mathbf{a}(\cdot,\nabla v_{m_k})\weak \boldsymbol{\eta}$ in $L^2\bigl(0,T;L^{p'(x)}(\Omega)^N\bigr)$. Note that, from \cite[Theorem 1.10.13]{megginson2012introduction} we have that $\bigl[L^{p'(x)}(\Omega)^N\bigr]^*=L^{p(x)}(\Omega)^N$ which is reflexive, being a finite cartesian product of reflexive spaces. In particular $L^{p(x)}(\Omega)^N$ has the Radon-Nikodym property. Therefore, using \textit{Riesz Representation Theorem for the Lebesgue-Bochner spaces} given as \cite[Theorem 2.2.9]{gasinski2005nonlinear} we have that $L^2\bigl(0,T;L^{p'(x)}(\Omega)^N\bigr)^*=L^2\bigl(0,T;L^{p(x)}(\Omega)^N\bigr)$.

\noindent So, for any $\boldsymbol{\psi}\in L^2\bigl(0,T;L^{p(x)}(\Omega)^N\bigr)$ we have that:

\begin{equation}\label{frumos1}
	\lim\limits_{k\to\infty}\int_{0}^T\int_{\Omega}\mathbf{a}(x,\nabla v_{m_k}(t,x))\cdot\boldsymbol{\psi}(t,x)\ dx\ dt=\int_{0}^T\int_{\Omega}\boldsymbol{\eta}(t,x)\cdot\boldsymbol{\psi}(t,x)\ dx\ dt.
\end{equation}

\noindent Since $\Theta\in L^2\bigl(0,T;W^{1,p(x)}(\Omega)\bigr)$ we have from Proposition \ref{propotheta} that $\nabla \Theta\in L^2\bigl(0,T;L^{p(x)}(\Omega)^N\bigr)$. Thence we may choose $\boldsymbol{\psi}:=\chi_{(t_1,t_2)}(\cdot)\nabla\Theta(\cdot,\cdot)\in  L^2\bigl(0,T;L^{p(x)}(\Omega)^N\bigr)$ in \eqref{frumos1} and obtain the following important relation:

\begin{equation}\label{frumos}
	\lim\limits_{k\to\infty}\int_{t_1}^{t_2}\int_{\Omega}\mathbf{a}(x,\nabla v_{m_k}(t,x))\cdot\nabla\Theta(t,x)\ dx\ dt=\int_{t_1}^{t_2}\int_{\Omega}\boldsymbol{\eta}(t,x)\cdot\nabla\Theta(t,x)\ dx\ dt.
\end{equation}

\noindent Thus \eqref{i5} becomes:

\begin{align}
	&\int_{t_1}^{t_2}\int_{\Omega}\boldsymbol{\eta}(t,x)\cdot\nabla\Theta(t,x)\ dx\ dt=\lim\limits_{k\to\infty} \int_{t_1}^{t_2}\int_{\Omega}\mathbf{a}(x,\nabla v_{m_k})\cdot\nabla\Theta\ dx\ dt=\nonumber\\
	=&\int_{t_1}^{t_2}\int_{\Omega} g\Theta\ dx\ dt-\int_{t_1}^{t_2}\int_\Omega\dfrac{\partial b(x,v)}{\partial t}\Theta\ dx\ dt-\lambda\int_{t_1}^{t_2}\int_{\Omega}b(x,v)\Theta\ dx\ dt.\label{i5frumos}
\end{align}

\noindent By setting now $\Theta=v\in L^{\infty}\bigl(0,T;W^{1,p(x)}(\Omega)\bigr)\subset L^{2}\bigl(0,T;W^{1,p(x)}(\Omega)\bigr)$ in \eqref{i5frumos}, we obtain that:

\begin{align}
	&\int_{t_1}^{t_2}\int_{\Omega}\boldsymbol{\eta}(t,x)\cdot\nabla v(t,x)\ dx\ dt=\lim\limits_{k\to\infty} \int_{t_1}^{t_2}\int_{\Omega}\mathbf{a}(x,\nabla v_{m_k})\cdot\nabla v\ dx\ dt=\nonumber\\
	=&\int_{t_1}^{t_2}\int_{\Omega} gv\ dx\ dt-\int_{t_1}^{t_2}\int_\Omega\dfrac{\partial b(x,v)}{\partial t}v\ dx\ dt-\lambda\int_{t_1}^{t_2}\int_{\Omega}b(x,v)v\ dx\ dt.\label{i5frumosv}
\end{align}

\noindent Replacing \eqref{i5frumos} in \eqref{i10} gives us that:

\begin{align}
	&\lim\limits_{k\to\infty} \int_{t_1}^{t_2}\int_{\Omega} \mathbf{a}(x,\nabla v_{m_k})\cdot\nabla v_{m_k}\ dx\ dt=\lim\limits_{k\to\infty} \int_{t_1}^{t_2}\int_{\Omega} \mathbf{a}(x,\nabla v_{m_k})\cdot\nabla v\ dx\ dt=\nonumber\\
	=&\int_{t_1}^{t_2}\int_{\Omega} gv\ dx\ dt-\int_{t_1}^{t_2}\int_\Omega\dfrac{\partial b(x,v)}{\partial t}v\ dx\ dt-\lambda\int_{t_1}^{t_2}\int_{\Omega}b(x,v)v\ dx\ dt=\nonumber\\
	=&\int_{t_1}^{t_2}\int_{\Omega}\boldsymbol{\eta}(t,x)\cdot\nabla v(t,x)\ dx\ dt,\ \forall\ (t_1,t_2)\subseteq (0,T)\label{eqminty}.
\end{align}

\noindent\textbf{Next, we will show that $\boldsymbol{\eta}=\mathbf{a}(\cdot,\nabla v(\cdot,\cdot))$ a.e. on $(0,T)\times\Omega$.}

\noindent Set for the moment $t_1=0$ and $t_2=T$ and consider any direction $\mathbf{w}\in L^{\infty}\bigl(0,T;L^{p(x)}(\Omega)^N\bigr)$. Because, from Proposition \ref{propoa0} \textbf{(4),(5)} the Nemytskii operator $\mathcal{N}_{\mathbf{a}}:L^{p(x)}(\Omega)^N\to L^{p'(x)}(\Omega)^N$ is continuous and bounded (in the nonlinear sense), we deduce using Proposition \ref{nonlinearlinfty}, that $\mathcal{N}_{\mathbf{a}}\circ \mathbf{w}=\mathbf{a}(\cdot,\mathbf{w})\in L^{\infty}\bigl(0,T;L^{p'(x)}(\Omega)^N\bigr)$.

\noindent From \textbf{Claim 17} we know that $\mathbf{a}(\cdot,\nabla v_{m_k})\in L^{\infty}\bigl(0,T;L^{p'(x)}(\Omega)^N\bigr)$, and hence $\mathbf{a}(\cdot,\nabla v_{m_k})-\mathbf{a}(\cdot,\mathbf{w})\in L^{\infty}\bigl(0,T;L^{p'(x)}(\Omega)^N\bigr)$, for each $k\geq 1$. 

\noindent Since, $\bigl(v_{m_k}\bigl)_{k\geq 1}\subset L^{\infty}\bigl(0,T;W^{1,p(x)}(\Omega)\bigr)$, we deduce from Proposition \ref{propotheta} that for each $k\geq 1$: $\nabla v_{m_k}\in L^{\infty}\bigl(0,T;L^{p(x)}(\Omega)^N\bigr)$. Thus $\nabla v_{m_k}-\mathbf{w}\in  L^{\infty}\bigl(0,T;L^{p(x)}(\Omega)^N\bigr)$ for each $k\geq 1$.

\noindent Applying now Proposition \ref{propolplprim} we get that the function:

\begin{equation}
	(0,T)\times\Omega\ni (t,x)\mapsto \bigl( \mathbf{a}(x,\nabla v_{m_k}(t,x)) -\mathbf{a}(x,\mathbf{w}(t,x))\bigr)\cdot (\nabla v_{m_{k}}(t,x)-\mathbf{w}(t,x))\in\mathbb{R},
\end{equation}

\noindent is from $L^{\infty}\bigl(0,T;L^1(\Omega)\bigr)\subset L^1\bigl(0,T;L^1(\Omega)\bigr)$. Now it makes sense to speak about the following integral, which from Remark \ref{remmax4} is positive:

\begin{equation}\label{benefic}
	\int_{0}^T\int_{\Omega} \underbrace{\bigl( \mathbf{a}(x,\nabla v_{m_k}) -\mathbf{a}(x,\mathbf{w})\bigr)\cdot (\nabla v_{m_{k}}-\mathbf{w})}_{\geq 0}\ dx\ dt\geq 0,\ \forall\ k\geq 1,\ \forall\ \mathbf{w}\in L^{\infty}\bigl(0,T;L^{p(x)}(\Omega)^N\bigr).
\end{equation}

\noindent Expanding \eqref{benefic} gives us that:

\begin{align}\label{benefic1}
	&\int_{0}^T\int_{\Omega}\mathbf{a}(x,\nabla v_{m_k})\cdot\nabla v_{m_k} \ dx\ dt+\int_{0}^T\int_{\Omega}\mathbf{a}(x,\mathbf{w})\cdot\mathbf{w} \ dx\ dt\nonumber\\ 
	-&	\int_{0}^T\int_{\Omega}\mathbf{a}(x,\nabla v_{m_k})\cdot\mathbf{w} \ dx\ dt -	\int_{0}^T\int_{\Omega} \mathbf{a}(x,\mathbf{w})\cdot\nabla v_{m_k}\ dx\ dt \geq 0.
\end{align}

\noindent$\blacktriangleright$ From \eqref{eqminty} we have that:

\begin{equation}\label{benefic2}
	\lim\limits_{k\to\infty} \int_{0}^T\int_{\Omega}\mathbf{a}(x,\nabla v_{m_k})\cdot\nabla v_{m_k} \ dx\ dt=\int_{0}^T\int_{\Omega}\boldsymbol{\eta}\cdot \nabla v\ dx\ dt.
\end{equation}

\bigskip
\noindent$\blacktriangleright$ From \textbf{Claim 17} we know that $\mathbf{a}(\cdot,\nabla v_{m_k})\weak\boldsymbol{\eta}$ in $L^2\bigl(0,T;L^{p'(x)}(\Omega)^N\bigr)$. This means that for any $\boldsymbol{\psi}\in L^2\bigl(0,T;L^{p'(x)}(\Omega)^N\bigr)^*=L^2\bigl(0,T;L^{p(x)}(\Omega)^N\bigr)$ we have that:

\begin{equation}
	\lim\limits_{k\to\infty} \int_{0}^T\int_{\Omega}\mathbf{a}(x,\nabla v_{m_k})\cdot\boldsymbol{\psi} \ dx\ dt=\int_{0}^T\int_{\Omega}\boldsymbol{\eta}\cdot\boldsymbol{\psi} \ dx\ dt.
\end{equation}

\noindent By choosing $\boldsymbol{\psi}:=\mathbf{w}\in L^\infty\bigl(0,T;L^{p(x)}(\Omega)^N\bigr)\subset L^2\bigl(0,T;L^{p(x)}(\Omega)^N\bigr)$ we obtain that:

\begin{equation}\label{benefic3}
	\lim\limits_{k\to\infty} \int_{0}^T\int_{\Omega}\mathbf{a}(x,\nabla v_{m_k})\cdot\mathbf{w} \ dx\ dt=\int_{0}^T\int_{\Omega}\boldsymbol{\eta}\cdot\mathbf{w} \ dx\ dt.
\end{equation}

\bigskip

\noindent$\blacktriangleright$ Remeber from the proof of Proposition \ref{propotheta} that $\overline{G}=\nabla:L^2\bigl(0,T;W^{1,p(x)}(\Omega)\bigr)\to L^{2}\bigl(0,T;L^{p(x)}(\Omega)^N\bigr)$ is a bounded linear operator. Using the fact that bounded linear operators preserve weak convergence\footnote{See Proposition \ref{lemmaboundedweak} from the Appendix.}, and $v_{m_k}\weak v$ in $L^2\bigl(0,T;W^{1,p(x)}(\Omega)\bigr)$ (see \textbf{Claim 8}), we get that $\nabla v_{m_k}\weak \nabla v$ in $L^{2}\bigl(0,T;L^{p(x)}(\Omega)^N\bigr)$. This means that for any $\boldsymbol{\psi}\in L^{2}\bigl(0,T;L^{p(x)}(\Omega)^N\bigr)^*=L^{2}\bigl(0,T;L^{p'(x)}(\Omega)^N\bigr)$ one may write:

\begin{equation}
		\lim\limits_{k\to\infty} \int_{0}^T\int_{\Omega}\boldsymbol{\psi}\cdot\nabla v_{m_k} \ dx\ dt=\int_{0}^T\int_{\Omega}\boldsymbol{\psi}\cdot\nabla v \ dx\ dt.
\end{equation}

\noindent Setting $\boldsymbol{\psi}:=\mathbf{a}(\cdot,\mathbf{w})\in L^{\infty}\bigl(0,T;L^{p'(x)}(\Omega)^N\bigr)\subset L^{2}\bigl(0,T;L^{p'(x)}(\Omega)^N\bigr)$, we get that:

\begin{equation}\label{benefic4}
	\lim\limits_{k\to\infty} \int_{0}^T\int_{\Omega}\mathbf{a}(x,\mathbf{w})\cdot\nabla v_{m_k} \ dx\ dt=\int_{0}^T\int_{\Omega}\mathbf{a}(\cdot,\mathbf{w})\cdot\nabla v \ dx\ dt.
\end{equation}

\noindent Now, everything is prepared so that we can make $k\to\infty$ in \eqref{benefic1}. Thus we obtain using \eqref{benefic2}, \eqref{benefic3} and \eqref{benefic4} that:

\begin{align}
	&\int_{0}^T\int_{\Omega}\boldsymbol{\eta}\cdot\nabla v \ dx\ dt+\int_{0}^T\int_{\Omega}\mathbf{a}(x,\mathbf{w})\cdot\mathbf{w} \ dx\ dt\nonumber\\ 
	-&	\int_{0}^T\int_{\Omega}\boldsymbol{\eta}\cdot\mathbf{w} \ dx\ dt -	\int_{0}^T\int_{\Omega} \mathbf{a}(x,\mathbf{w})\cdot\nabla v\ dx\ dt \geq 0,
\end{align}

\noindent which can be factored as:

\begin{equation}\label{beneficeta}
	\int_{0}^T\int_{\Omega} \underbrace{\bigl( \boldsymbol{\eta} -\mathbf{a}(x,\mathbf{w})\bigr)\cdot (\nabla v-\mathbf{w})}_{\geq 0}\ dx\ dt\geq 0,\ \forall\ k\geq 1,\ \forall\ \mathbf{w}\in L^{\infty}\bigl(0,T;L^{p(x)}(\Omega)^N\bigr).
\end{equation}

\bigskip

\noindent Take now any arbitrary direction $\mathbf{W}\in L^{\infty}\bigl(0,T;L^{p(x)}(\Omega)^N\bigr)$. Consider any sequence $(\epsilon_n)_{n\geq 1}\subset (0,\infty)$ with $\epsilon_n\to 0$, and for each $n\geq 1$, set $\mathbf{w}_n:=\nabla v-\epsilon_n \mathbf{W}\in L^{\infty}\bigl(0,T;L^{p(x)}(\Omega)^N\bigr)$. 

\noindent By substituting $\mathbf{w}=\mathbf{w}_n$ in \eqref{beneficeta} we get that:

\begin{equation}
	\epsilon_n\int_{0}^T\int_{\Omega} \bigl( \boldsymbol{\eta} -\mathbf{a}(x,\mathbf{w}_n)\bigr)\cdot \mathbf{W}\ dx\ dt\geq 0,\ \forall\ n\geq 1.
\end{equation}

\noindent Dividing by $\varepsilon_n>0$ allows us to write that:

\begin{equation}\label{grozavetan}
\int_{0}^T\int_{\Omega} \bigl( \boldsymbol{\eta} -\mathbf{a}(x,\mathbf{w}_n)\bigr)\cdot \mathbf{W}\ dx\ dt\geq 0,\ \forall\ n\geq 1.
\end{equation}

\noindent Now we have for each $n\geq 1$ that:

\begin{align*}
	&\left |\int_{0}^T\int_{\Omega} \bigl( \boldsymbol{\eta}-\mathbf{a}(x,\mathbf{w}_n)\bigr)\cdot\nabla \mathbf{W}\ dx\ dt -\int_{0}^T\int_{\Omega} \bigl(\boldsymbol{\eta}- \mathbf{a}(x,\nabla v)\bigr)\cdot\nabla\mathbf{W}\ dx\ dt \right |\\
	=&\left |\int_{0}^T\int_{\Omega} \bigl( \mathbf{a}(x,\nabla v) -\mathbf{a}(x,\mathbf{w}_n)\bigr)\cdot \mathbf{W}\ dx\ dt\right |\\
	\leq &\int_{0}^T\int_{\Omega} \left |\bigl( \mathbf{a}(x,\mathbf{w}_n)-\mathbf{a}(x,\nabla v)\bigr)\cdot \mathbf{W} \right | \ dx\ dt\\
	\stackrel{\text{Prop.} \ref{propolplprim}}{\leq}\ & 2\Vert \mathbf{a}(x,\mathbf{w}_n)-\mathbf{a}(x,\nabla v)\Vert_{L^{2}(0,T;L^{p'(x)}(\Omega)^N)}\cdot \Vert\mathbf{W}\Vert_{L^{\infty}(0,T;L^{p(x)}(\Omega)^N)}\stackrel{n\to\infty}{\longrightarrow} 0,
\end{align*}

\noindent because $\Vert\mathbf{w}_n-\nabla v\Vert_{L^{\infty}(0,T;L^{p'(x)}(\Omega))}=\epsilon_n\Vert W\Vert_{L^{\infty}(0,T;L^{p'(x)}(\Omega))}\stackrel{n\to\infty}{\longrightarrow} 0$ and from Proposition \ref{propoa3} we have that $\mathbf{a}(x,\mathbf{w}_n)\to\mathbf{a}(x,\nabla v)$ in $L^2\bigl(0,T;L^{p'(x)}(\Omega)^N\bigr)$.

\noindent Thus making $n\to\infty$ in \eqref{grozavetan} we will get that:

\begin{equation}\label{grozaveta}
	\int_{0}^T\int_{\Omega} \bigl( \boldsymbol{\eta} -\mathbf{a}(x,\nabla v)\bigr)\cdot \mathbf{W}\ dx\ dt\geq 0,\ \forall\ \mathbf{W}\in L^{\infty}\bigl(0,T;L^{p(x)}(\Omega)^N\bigr).
\end{equation}

\noindent Writing \eqref{grozaveta} for $-\mathbf{W}$ instead of $\mathbf{W}$, which will be also a member of $L^{\infty}\bigl(0,T;L^{p(x)}(\Omega)^N\bigr)$ we will obtain that:

\begin{equation}\label{grozaveta0}
	\int_{0}^T\int_{\Omega} \bigl( \boldsymbol{\eta} -\mathbf{a}(x,\nabla v)\bigr)\cdot \mathbf{W}\ dx\ dt\leq 0,\ \forall\ \mathbf{W}\in L^{\infty}\bigl(0,T;L^{p(x)}(\Omega)^N\bigr).
\end{equation}

\noindent Combining now \eqref{grozaveta} with \eqref{grozaveta0} gives us that:

\begin{equation}\label{grozavetasuper}
	\int_{0}^T\int_{\Omega} \bigl( \boldsymbol{\eta} -\mathbf{a}(x,\nabla v)\bigr)\cdot \mathbf{W}\ dx\ dt= 0,\ \forall\ \mathbf{W}\in L^{\infty}\bigl(0,T;L^{p(x)}(\Omega)^N\bigr).
\end{equation}

\noindent Now, if we choose in \eqref{grozavetasuper} $W:=\boldsymbol{\eta} -\mathbf{a}(\cdot,\nabla v)\in  L^{\infty}\bigl(0,T;L^{p(x)}(\Omega)^N\bigr)$ we get from \textit{Fubini's theorem} that:

\begin{equation}
 \int_{(0,T)\times\Omega} \underbrace{\bigl|\boldsymbol{\eta} -\mathbf{a}(x,\nabla v)\bigr |^2}_{\geq 0}\ dx\ dt=\int_{0}^T\int_{\Omega} \bigl|\boldsymbol{\eta} -\mathbf{a}(x,\nabla v)\bigr |^2\ dx\ dt=0.
\end{equation}

\noindent In conclusion $\boldsymbol{\eta}=\mathbf{a}(\cdot,\nabla v)$ a.e. on $(0,T)\times\Omega$.

\bigskip

\noindent\textbf{Now equation \eqref{i5frumos} says that for any subinterval $(t_1,t_2)\subseteq (0,T)$:}

\begin{align}
	&\int_{t_1}^{t_2}\int_{\Omega}\mathbf{a}(x,\nabla v(t,x))\cdot\nabla\Theta(t,x)\ dx\ dt=\lim\limits_{k\to\infty} \int_{t_1}^{t_2}\int_{\Omega}\mathbf{a}(x,\nabla v_{m_k})\cdot\nabla\Theta\ dx\ dt=\nonumber\\
	=&\int_{t_1}^{t_2}\int_{\Omega} g\Theta\ dx\ dt-\int_{t_1}^{t_2}\int_\Omega\dfrac{\partial b(x,v)}{\partial t}\Theta\ dx\ dt-\lambda\int_{t_1}^{t_2}\int_{\Omega}b(x,v)\Theta\ dx\ dt,\ \forall\ \Theta\in L^{2}\bigl(0,T;W^{1,p(x)}(\Omega)\bigr).\label{i5smecher}
\end{align}

\noindent\textbf{The last part of the proof:} We know from \cite[Theorem 3.1]{kovacik} that the Banach space $W^{1,p(x)}(\Omega)$ is separable. Therefore we can find a sequence $(\phi_n)_{n\geq 1}\subset W^{1,p(x)}(\Omega)$ which is dense in $W^{1,p(x)}(\Omega)$. For each $n\geq 1$ we denote $h_n:(0,T)\to\mathbb{R}$ the function given by:

\begin{equation}
h_n(t)=\int_{\Omega}\dfrac{\partial b(x,v)}{\partial t}\phi_n\ dx+\int_{\Omega}\mathbf{a}(x,\nabla v(t,x))\cdot\nabla\phi_n\ dx+\lambda\int_{\Omega}b(x,v)\phi_n\ dx-\int_{\Omega} g\phi_n\ dx.\
\end{equation}

\noindent Setting for each $n\geq 1$: $\Theta:=\phi_n\in W^{1,p(x)}(\Omega)\subset L^2\bigl(0,T;W^{1,p(x)}(\Omega)\bigr)$ we get that:

\begin{equation}
	\int_{t_1}^{t_2} h_n(t)\ dt=0,\ \forall\ (t_1,t_2)\subseteq (0,T).
\end{equation}

\noindent At this point we need the following lemma:

\begin{lemma}\label{intt1t2}
	If $h\in L^1(a,b)$ has the property that $\displaystyle\int_{t_1}^{t_2} h(\tau)\ d\tau=0$, for any subinterval $(t_1,t_2)\subseteq (a,b)$, then $h(t)=0$ for a.e. $t\in (a,b)$.
\end{lemma}

\begin{proof} We denote $H:(a,b)\to\mathbb{R},\ H(t)=\displaystyle\int_{a}^t h(\tau)\ d\tau$. From our hypothesis we have that $H(t)=0$ for any $t\in (a,b)$. Therefore $H'(t)=0$ for any $t\in (a,b)$. But from \textit{Lebesgue differentiation theorem} (see Theorem 6.3.6 in \cite[page 172]{Cohn}) we get that $H'(t)=h(t)$ for a.e. $t\in (a,b)$. Thus $h(t)=0$ for a.e. $t\in (a,b)$.
	
\end{proof}

\noindent Using Lemma \ref{intt1t2} for each $h_n\in L^1(0,T),\ n\geq 1$ we deduce that $h_n(t)=0$ for a.e. $t\in (0,T)$. Therefore there is a null-measure set $A_n\subset (0,T)$ such that $h_n(t)=0$ for any $t\in (0,T)\setminus A_n$. Let's denote $A=\displaystyle\bigcup_{n=1}^{\infty} A_n$. It follows from the countable subadditivity of the Lebesgue measure that $|A|\leq \displaystyle\sum_{n=1}^{\infty} |A_n|=0$. Thereof, for any $t\in (0,T)\setminus A$ we have that $h_n(t)=0$ for any $n\geq 1$. We can therefore state that $\text{for a.e.}\ t\in (0,T)$:

\begin{equation}\label{smecheriesep}
\int_{\Omega}\dfrac{\partial b(x,v)}{\partial t}\phi_n\ dx+\int_{\Omega}\mathbf{a}(x,\nabla v(t,x))\cdot\nabla\phi_n\ dx+\lambda\int_{\Omega}b(x,v)\phi_n\ dx=\int_{\Omega} g\phi_n\ dx,\forall\ n\geq 1. 
\end{equation}

\noindent Now fix any $\phi\in W^{1,p(x)}(\Omega)$. From the density of $(\phi_n)_{n\geq 1}$ in $W^{1,p(x)}(\Omega)$ we can find a subsequence $(\phi_{n_j})_{j\geq 1}$ such that:

\begin{equation}
	\lim\limits_{j\to\infty} \Vert \phi_{n_j}-\phi\Vert_{W^{1,p(x)}(\Omega)}=0\ \Longrightarrow\  \begin{cases}\nabla\phi_{n_j}\to\nabla \phi\ \text{in}\ L^{p(x)}(\Omega)^N\\ \phi_{n_j}\to \phi \ \text{in}\ L^2(\Omega) \end{cases}.
\end{equation}

\noindent Note that for a.e. $t\in (0,T)$ we have that for each $j\geq 1$:

\begin{align}\label{sep1}
	\blacktriangleright\left |\int_{\Omega}\dfrac{\partial b(x,v)}{\partial t}\phi_{n_j}\ dx-\int_{\Omega}\dfrac{\partial b(x,v)}{\partial t}\phi\ dx\right |&=\left |\int_{\Omega}\dfrac{\partial b(x,v)}{\partial t}\bigl(\phi_{n_j}-\phi\bigr)\ dx \right |\nonumber\\
	\text{(Cauchy ineq.)}\ \ \ &\leq\ \left \Vert \dfrac{\partial b(x,v(t,\cdot))}{\partial t}\right\Vert_{L^2(\Omega)}\cdot\Vert\phi_{n_j}-\phi\Vert_{L^{2}(\Omega)} \stackrel{j\to\infty}{\longrightarrow} 0.
\end{align}

\bigskip

\begin{align}\label{sep2}
	&\blacktriangleright\left |\int_{\Omega}\mathbf{a}(x,\nabla v(t,x))\cdot\nabla\phi_n\ dx-\int_{\Omega}\mathbf{a}(x,\nabla v(t,x))\cdot\nabla\phi\ dx\right |=\left |\int_{\Omega}\mathbf{a}(x,\nabla v(t,x))\cdot\ \bigl(\nabla\phi_{n_j}-\nabla\phi\bigr)\ dx \right |\nonumber\\
	&\leq \int_{\Omega}|\mathbf{a}(x,\nabla v(t,x))|\cdot\ \bigl|\nabla\phi_{n_j}-\nabla\phi\bigr|\ dx \nonumber\\
&\text{(H\"{o}lder ineq.)}\ \ \	\leq 2\left \Vert|\mathbf{a}(\cdot,\nabla v(t,\cdot))|\right\Vert_{L^{p'(x)}(\Omega)}\cdot\Vert |\nabla\phi_{n_j}-\nabla\phi|\Vert_{L^{p(x)}(\Omega)}\nonumber \\
	&=2\left \Vert \mathbf{a}(\cdot,\nabla v(t,\cdot))\right\Vert_{L^{p'(x)}(\Omega)^N}\cdot\Vert \nabla\phi_{n_j}-\nabla\phi\Vert_{L^{p(x)}(\Omega)^N} \stackrel{j\to\infty}{\longrightarrow} 0.
\end{align}

\bigskip

\begin{align}\label{sep3}
	\blacktriangleright\left |\int_{\Omega}b(x,v)\phi_{n_j}\ dx-\int_{\Omega}b(x,v)\phi\ dx\right |&=\left |\int_{\Omega}b(x,v)\bigl(\phi_{n_j}-\phi\bigr)\ dx \right |\nonumber\\
	\text{(Cauchy ineq.)}\ \ \ &\leq\ \left \Vert b(\cdot,v(t,\cdot))\right\Vert_{L^2(\Omega)}\cdot\Vert\phi_{n_j}-\phi\Vert_{L^{2}(\Omega)} \stackrel{j\to\infty}{\longrightarrow} 0.
\end{align}

\bigskip

\begin{align}\label{sep4}
	\blacktriangleright\left |\int_{\Omega}g\phi_{n_j}\ dx-\int_{\Omega}g\phi\ dx\right |&=\left |\int_{\Omega}g\bigl(\phi_{n_j}-\phi\bigr)\ dx \right |\nonumber\\
	\text{(Cauchy ineq.)}\ \ \ &\leq\ \left \Vert g(t,\cdot)\right\Vert_{L^2(\Omega)}\cdot\Vert\phi_{n_j}-\phi\Vert_{L^{2}(\Omega)} \stackrel{j\to\infty}{\longrightarrow} 0.
\end{align}

\noindent Making now $n\to\infty$ in \eqref{smecheriesep} and using \eqref{sep1}, \eqref{sep2}, \eqref{sep3} and \eqref{sep4} finally gives us that for a.e. $t\in (0,T)$:

\begin{equation}\label{smecheriesepfinal}
	\int_{\Omega}\dfrac{\partial b(x,v)}{\partial t}\phi\ dx+\int_{\Omega}\mathbf{a}(x,\nabla v(t,x))\cdot\nabla\phi\ dx+\lambda\int_{\Omega}b(x,v)\phi\ dx=\int_{\Omega} g\phi\ dx,\forall\ \phi\in W^{1,p(x)}(\Omega). 
\end{equation}

\noindent This shows that $v\in C\big ([0,T]; L^2(\Omega)\big )\cap H^1\big ((0,T);L^2(\Omega)\big )\cap L^{\infty}\big (0,T;W^{1,p(x)}(\Omega)\big )\cap\mathcal{V}_{[\varepsilon,\delta]}$ (see Remark \ref{remvsolfinal}) is indeed a solution for our auxiliary problem \eqref{eqdpgaux}. Uniqueness of the solution for \eqref{eqdpgaux} was proved before in Theorem \ref{thmveryunique}.

\end{proof}

\subsection{Existence of solution for the auxiliary problem with irregular initial data}

\begin{theorem}\label{theoremverygeneral} With the extra assumption that $p\in\mathcal{P}^{\text{log}}(\Omega)$ -- i.e. $p$ is a log-H\"{o}lder continuous exponent --, the auxiliary problem \eqref{eqdpgaux} has exactly one weak solution for any initial data $u_0\in \mathcal{U}_{[\varepsilon,\delta]}$.
\end{theorem}

\begin{proof} We know that $C^{\infty}_{c}(\Omega)\subset W^{1,p(x)}(\Omega)$, and that $C^{\infty}_c(\Omega)$ is dense in $L^r(\Omega)$ for any $r\in [1,\infty)$ with respect to its norm\footnote{See Corollary 4.23 from \cite{brezis2011functional}.}. Since $u_0\in \mathcal{U}_{[\varepsilon,\delta]}\subset L^{\infty}(\Omega)\subset L^2(\Omega)$, we get that there is a sequence $(\tilde{u}_{0,n})_{n\geq 1}\subset C^{\infty}_c(\Omega)\subset W^{1,p(x)}(\Omega)$ such that 
	
\begin{equation}
	\tilde{u}_{0,n}\stackrel{n\to\infty}{\longrightarrow} u_0\ \text{in}\ L^2(\Omega).
\end{equation}

\noindent Consider now the truncation function $T:\mathbb{R}\to [\varepsilon,\delta],\ T(s)=\begin{cases} \varepsilon, & s<\varepsilon\\ s, & s\in [\varepsilon,\delta]\\ \delta, & s>\delta\end{cases}$. It is easy to check that $T$ is a $1$ -- Lipschitz function. Therefore, using the \textit{Chain Rule for weak derivatives}\footnote{See Theorem 11.1.2 from \cite{Hasto}.} we deduce that $u_{0,n}:=T\circ\tilde{u}_{0,n}\in W^{1,p(x)}(\Omega)\cap\mathcal{U}_{[\varepsilon,\delta]}$ for every $n\geq 1$.

\noindent Note that $u_{0,n}\stackrel{n\to\infty}{\longrightarrow} u_0\ \text{in}\ L^2(\Omega)$, because:

\begin{align*}
	\Vert u_{0,n}-u_{0}\Vert_{L^2(\Omega)}&=\left (\int_{\Omega} |u_{0,n}(x)-u_0(x)|^2\ dx \right )^{\frac{1}{2}}=\left (\int_{\Omega} |T(\tilde{u}_{0,n}(x))-T(u_0(x))|^2\ dx \right )^{\frac{1}{2}}\\
	&\leq \left (\int_{\Omega} |\tilde{u}_{0,n}(x)-u_0(x)|^2\ dx \right )^{\frac{1}{2}}=\Vert \tilde{u}_{0,n}-u_{0}\Vert_{L^2(\Omega)}\stackrel{n\to\infty}{\longrightarrow} 0.
\end{align*}
	
\noindent Applying Theorem \ref{thmauxiliar}, we get for each $n\geq 1$ that there is a unique solution

\begin{equation}\label{relvninspatii}
	v_n\in C\big ([0,T]; L^2(\Omega)\big )\cap H^1\big ((0,T);L^2(\Omega)\big )\cap L^{\infty}\big (0,T;W^{1,p(x)}(\Omega)\big )\cap\mathcal{V}_{[\varepsilon,\delta]},
\end{equation}

\noindent for the following problem:

\begin{equation}\label{eqdpgauxn}\tag{$\textnormal{A}_n$}
	\begin{cases}\dfrac{\partial b(x,v_n(t,x))}{\partial t}-\operatorname{div}\mathbf{a}(x,\nabla v_n)+\lambda b(x,v_n)=g(t,x), & (t,x)\in (0,T)\times\Omega\\[3mm] \mathbf{a}(x,\nabla v_n)\cdot\nu=0, & (t,x)\in (0,T)\times\partial\Omega\\[3mm] v_n(0,x)=u_{0,n}(x)\in [\varepsilon,\delta], & x\in\Omega\end{cases}
\end{equation}

\noindent This means that for each $n\geq 1$, we have for a.e. $t\in (0,T)$ that:

\begin{equation}\label{ecuatialuivn}
	\int_{\Omega}\dfrac{\partial b(x,v_n)}{\partial t}\phi\ dx+\int_{\Omega}\mathbf{a}(x,\nabla v_n(t,x))\cdot\nabla\phi\ dx+\lambda\int_{\Omega}b(x,v_n)\phi\ dx=\int_{\Omega} g\phi\ dx,\forall\ \phi\in W^{1,p(x)}(\Omega). 
\end{equation}

\noindent If this relation does not hold in a null-measure set $I_n\subset (0,T)$, then if we denote $I=\displaystyle\bigcup_{n=1}^{\infty} I_n\subset (0,T)$ we have from the countable subadditivity of the Lebesgue measure that $|I|\leq \displaystyle\sum_{n=1}^{\infty} |I_n|=0$. Thereof, we have that for any $t\in (0,T)\setminus I$ (i.e. for a.e. $t\in (0,T)$) the equation \eqref{ecuatialuivn} holds for every $n\geq 1$ and any $\phi\in W^{1,p(x)}(\Omega)$.

\bigskip

\noindent\textbf{\underline{Fact I}: There is some function $v\in C\bigl([0,T];L^2(\Omega)\bigr)$ such that $v_n\to v$ in $C\bigl([0,T];L^2(\Omega)\bigr)$.} 

\begin{proof} For a.e. $t\in (0,T)$, for any $n,m\geq 1$ and any $\phi\in W^{1,p(x)}(\Omega)$ we have that
	
	\begin{align*}
		&\int_\Omega\dfrac{\partial b(x,v_n)}{\partial t}\phi\ dx+\int_{\Omega}\mathbf{a}(x,\nabla v_n)\cdot\nabla\phi\ dx+\lambda\int_{\Omega} b(x,v_n)\phi\ dx=\int_{\Omega} g\phi\ dx\\
		&\int_\Omega\dfrac{\partial b(x,v_m)}{\partial t}\phi\ dx+\int_{\Omega}\mathbf{a}(x,\nabla v_m)\cdot\nabla\phi\ dx+\lambda\int_{\Omega} b(x,v_m)\phi\ dx=\int_{\Omega} g\phi\ dx.
	\end{align*}

	\noindent Substracting these relations will give us
	
	\begin{align}
		&\int_\Omega\dfrac{\partial \big [b(x,v_n)-b(x,v_m)\big ]}{\partial t}\phi\ dx+\int_{\Omega}\big [\mathbf{a}(x,\nabla v_n)-\mathbf{a}(x,\nabla v_m)\big ]\cdot\nabla\phi\ dx+\lambda\int_{\Omega} \big [b(x,v_n)-b(x,v_m)\big ]\phi\ dx \nonumber\\
		&=0.
	\end{align}
	
	\noindent As in the proof of the \textit{weak comparison principle} we define for each $\tau>0$ the function $\phi_{\tau}=\begin{cases} 1, & v_n(t,\cdot)-v_m(t,\cdot)\geq\tau \\ \dfrac{v_n(t,\cdot)-v_m(t,\cdot)}{\tau}, & v_n(t,\cdot)-v_m(t,\cdot)<\tau \\ -1 & v_n(t,\cdot)-v_m(t,\cdot)\leq -\tau\end{cases}\in W^{1,p(x)}(\Omega)$ and select the test function $\phi=\phi_{\tau}^+\in W^{1,p(x)}(\Omega)$. Hence we obtain:
	
	\begin{align*}
		&\int_\Omega\dfrac{\partial \big [b(x,v_n)-b(x,v_m)\big ]}{\partial t}\phi_{\tau}^+\ dx+\lambda\int_{\Omega} \big [b(x,v_n)-b(x,v_m)\big ]\phi_{\tau}^+\ dx\\
		&=\dfrac{1}{\tau}\int_{\Omega}\big [\mathbf{a}(x,\nabla v_n)-\mathbf{a}(x,\nabla v_m)\big ]\cdot\big [\nabla v_n-\nabla v_m\big ]\chi_{v_m(t,\cdot)-v_n(t,\cdot)\in (0,\tau)}\ dx\\
		&\geq 0.
	\end{align*}
	
	\noindent Since $\lim\limits_{\tau\to 0^+} \phi_{\tau}^+\to\chi_{v_m(t,\cdot)-v_n(t,\cdot)}$ strongly in $L^2(\Omega)$ we deduce that
	
	\begin{align}\label{ecuatiaprincipala10}
		&\int_\Omega\dfrac{\partial \big [b(x,v_n)-b(x,v_m)\big ]}{\partial t}\chi_{v_m(t,\cdot)>v_n(t,\cdot)}\ dx+\lambda\int_{\Omega} \big [b(x,v_n)-b(x,v_m)\big ]\chi_{v_m(t,\cdot)>v_n(t,\cdot)}\ dx \geq 0.
	\end{align}

	\noindent We repeat the same process by switching the roles of $v_n$ and $v_m$ and obtain:
	
	\begin{align}\label{ecuatiaprincipala20}
		&\int_\Omega\dfrac{\partial \big [b(x,v_m)-b(x,v_n)\big ]}{\partial t}\chi_{v_n(t,\cdot)>v_m(t,\cdot)}\ dx+\lambda\int_{\Omega} \big [b(x,v_m)-b(x,v_n)\big ]\chi_{v_n(t,\cdot)>v_m(t,\cdot)}\ dx \geq 0.
	\end{align}
	
	\noindent Define the function $w_{n,m}(t,\cdot)=b(\cdot,v_n(t,\cdot))-b(\cdot,v_m(t,\cdot))\in C\bigl([0,T];L^2(\Omega)\bigr)$ -- see Remark \ref{rembo}. The inequality \eqref{ecuatiaprincipala10} may be writen as:
	
	\begin{equation}
		\int_{\Omega} \dfrac{\partial w_{n,m}}{\partial t}\chi_{v_m(t,\cdot)>v_n(t,\cdot)}\ dx 
		+\lambda\int_\Omega w_{n,m} \chi_{v_m(t,\cdot)>v_n(t,\cdot)}\ dx\geq 0.
	\end{equation}
	
	\noindent Using the strict monotony of $b$ we come at the following relation that holds for a.e. $t\in (0,T)$:
	
	\begin{equation}\label{ecuatieconti10}
		\int_{\Omega} \dfrac{\partial w_{n,m}^-(t,x)}{\partial t} dx 
		+\lambda\int_\Omega w_{n,m}^-(t,x) dx\leq 0.
	\end{equation}
	
	\noindent Similarly, the inequality \eqref{ecuatiaprincipala20} can be written as:
	
	\begin{equation}
		\int_{\Omega} -\dfrac{\partial w_{n,m}}{\partial t}\chi_{v_n(t,\cdot)>v_m(t,\cdot)}\ dx 
		-\lambda\int_\Omega w_{n,m} \chi_{v_n(t,\cdot)>v_m(t,\cdot)}\ dx\geq 0,
	\end{equation}
	
	\noindent i.e. for a.e. $t\in (0,T)$
	
	\begin{equation}\label{ecuatieconti20}
		\int_{\Omega} \dfrac{\partial w_{n,m}^+(t,x)}{\partial t} dx 
		+\lambda\int_\Omega w_{n,m}^+(t,x) dx\leq 0.
	\end{equation}
	
	\noindent Adding now \eqref{ecuatieconti10} and \eqref{ecuatieconti20} we get that
	
	\begin{equation}\label{ecuatieconti0}
		\int_{\Omega} \dfrac{\partial |w_{n,m}(t,x)|}{\partial t} dx 
		+\lambda\int_\Omega |w_{n,m}(t,x)| dx\leq 0.
	\end{equation}
	
	\noindent For each $n,m\geq 1$, consider now $h_{n,m}:[0,T]\to [0,\infty)$, $h_{n,m}(t)=\displaystyle\int_\Omega |w_n(t,x)|\ dx$. We have that, $h_{n,m}(0)=\displaystyle\int_\Omega |w_{n,m}(0,x)|\ dx=\displaystyle\int_\Omega |b(\cdot,u_{0,n})-b(\cdot,u_{0,m})|\ dx$ because $v_n(0,\cdot)=u_{0,n}$ and $v_{m}(0,\cdot)=u_{0,m}$ for any $n,m\geq 1$. The same arguments used in the proof of the \textit{weak comparison principle} allows us to say that $h_{n,m}\in C([0,T])$ is differentiable a.e. on $(0,T)$ and from Theorem \ref{thmmaxpri} $h_{n,m}'(t)=\displaystyle\int_{\Omega}\dfrac{\partial |w_{n,m}|}{\partial t}(t,x)\ dx$ for a.e. $t\in (0,T)$. Therefore \eqref{ecuatieconti0} rewrites as:
	
	\begin{equation}
		h_{n,m}'(t)+\lambda h_{n,m}(t)\leq 0,\ \text{for a.e.}\ t\in (0,T).
	\end{equation}
	
	\noindent From \textit{Gronwall's inequality -- differential form}\footnote{See Theorem \ref{gronwalldiff} from the Appendix.} we obtain for a.e. $t\in (0,T)$ that:
	
	\begin{align*}
		0\leq h_{n,m}(t)&\leq e^{-\lambda t}h_{n,m}(0)= e^{-\lambda t}\Vert b(\cdot,u_{0,n})-b(\cdot,u_{0,m})\Vert_{L^1(\Omega)}\\
		\text{(Cauchy ineq.)}\ \ \ 	&\leq e^{-\lambda t}\sqrt{|\Omega|}\cdot\Vert b(\cdot,u_{0,n})-b(\cdot,u_{0,m})\Vert_{L^2(\Omega)}\\
		\eqref{blipschitz}\ \ \	&\leq e^{-\lambda t} L_0\sqrt{|\Omega|}\cdot\Vert u_{0,n}-u_{0,m}\Vert_{L^2(\Omega)}.
	\end{align*}
	
	\noindent But, since $h_{n,m}$ is continuous on $[0,T]$ we find that the above inequality holds in fact for every $t\in [0,T]$, i.e. 
	
	\begin{align}\label{poc1}
		\displaystyle\int_\Omega |b(x,v_n(t,x))-b(x,v_m(t,x))|\ dx=h_{m,n}(t)&\leq e^{-\lambda t} L_0\sqrt{|\Omega|}\cdot\Vert u_{0,n}-u_{0,m}\Vert_{L^2(\Omega)}\nonumber \\
		&\leq L_0\sqrt{|\Omega|}\cdot\Vert u_{0,n}-u_{0,m}\Vert_{L^2(\Omega)} ,\ \forall\ t\in [0,T].
	\end{align}
	
	\noindent Also, from Remark \ref{rem23} and the fact that $v_n,v_m\in\mathcal{V}_{[\varepsilon,\delta]}$, it follows that:

	\begin{align}\label{poc2}
		\displaystyle\int_\Omega |b(x,v_n(t,x))-b(x,v_m(t,x))|\ dx&\geq \ell_0\int_{\Omega} |v_n(t,x)-v_m(t,x)|\ dx\nonumber \\
		&\geq\dfrac{\ell_0}{\delta-\varepsilon}\int_{\Omega} (\delta-\varepsilon)|v_n(t,x)-v_m(t,x)|\ dx\nonumber\\
		&\geq\dfrac{\ell_0}{\delta-\varepsilon}\int_{\Omega}|v_n(t,x)-v_m(t,x)|^2\ dx\\
		&=\dfrac{\ell_0}{\delta-\varepsilon}\Vert v_n(t,\cdot)-v_{m}(t,\cdot)\Vert_{L^2(\Omega)}^2,\ \forall\ t\in [0,T].
	\end{align}
	
	\noindent Henceforth, combining \eqref{poc1} and \eqref{poc2}, we get that for any $n,m\geq 1$:
	
	\begin{align*}
\Vert v_n(t,\cdot)-v_{m}(t,\cdot)\Vert_{L^2(\Omega)}\leq \sqrt{\dfrac{L_0\sqrt{|\Omega|}(\delta-\varepsilon)}{\ell_0}}\cdot \Vert u_{0,n}-u_{0,m}\Vert_{L^2(\Omega)}^{\frac{1}{2}},\ \forall\ t\in [0,T].
	\end{align*}
	
	\noindent Thereof:
	
	\begin{equation}
		\Vert v_{n}-v_{m}\Vert_{C([0,T];L^2(\Omega))}=\sup_{t\in [0,T]} \Vert v_n(t,\cdot)-v_{m}(t,\cdot)\Vert_{L^2(\Omega)}\leq \sqrt{\dfrac{L_0\sqrt{|\Omega|}(\delta-\varepsilon)}{\ell_0}}\cdot \Vert u_{0,n}-u_{0,m}\Vert_{L^2(\Omega)}^{\frac{1}{2}}.
	\end{equation}

	\noindent This relation show us that $\bigl(v_n\bigr)_{n\geq 1}$ is a Cauchy sequence from $C\bigl([0,T];L^2(\Omega)\bigr)$ -- which is a Banach space, i.e. a complete normed vector space. Thence it is convergent to some function $v\in C\bigl([0,T];L^2(\Omega)\bigr)$, i.e.
	
	\begin{equation}
		v_n\stackrel{n\to\infty}{\longrightarrow} v\ \text{in}\ C\bigl([0,T];L^2(\Omega)\bigr).
	\end{equation}
	
	\begin{remark}\label{vnvlr0tl2omega} In particular 
		\begin{equation}
			v_n\stackrel{n\to\infty}{\longrightarrow} v\ \text{in}\ L^r\bigl(0,T;L^2(\Omega)\bigr),\ \forall\ r\in [1,\infty].
		\end{equation}
	\end{remark}
	
	\begin{remark}\label{v0initial}
		From $v_n\to v$ in $C\bigl([0,T];L^2(\Omega)\bigr)$ we deduce that $\lim\limits_{n\to\infty}\displaystyle\sup_{t\in [0,T]}\Vert v_n(t,\cdot)-v(t,\cdot)\Vert_{L^2(\Omega)}=0$. In particular, for $t=0$, this means that $\lim\limits_{n\to\infty} \Vert v_n(0,\cdot)-v(0,\cdot)\Vert_{L^2(\Omega)}=0$, i.e. $u_{0,n}=v_n(0,\cdot)\to v(0,\cdot)$ in $L^2(\Omega)$. But from the construction we know that $u_{0,n}\to u_0$ in $L^2(\Omega)$, and in conclusion 
		
		\begin{equation}
			v(0,\cdot)=u_0.
		\end{equation}
	\end{remark}
\end{proof}

\bigskip 

\noindent\textbf{\underline{Fact II}: $v\in\mathcal{V}_{[\varepsilon,\delta]}$.}

\begin{proof} From Remark \ref{vnvlr0tl2omega} we have in particular that $v_n\to v$ in $L^2\bigl(0,T;L^2(\Omega)\bigr)\simeq L^2\bigl((0,T)\times\Omega\bigr)$. Using the Corollary from \cite[page 234]{Jones}, we get that there is a subsequence $(v_{n_k})_{k\geq 1}$ with $v_{n_k}\to v$ pointwise a.e. on $(0,T)\times\Omega$. So there is a null-measure set $A_0\subset (0,T)\times\Omega$, with $v_{n_k}(t,x)\to v(t,x)$ for any $(t,x)\in\bigl[(0,T)\times\Omega\bigr]\setminus A_0$. We know that $v_{n_k}\in\mathcal{V}_{[\varepsilon,\delta]}$ for any $k\geq 1$ and therefore there is a null-measure set $A_k\subset (0,T)\times\Omega$ such that:
	
\begin{equation}
	\varepsilon\leq v_{n_k}(t,x)\leq \delta,\ \text{for any}\ (t,x)\in \bigl [(0,T)\times\Omega\bigr]\setminus A_k.
\end{equation}

\noindent Consider $\tilde{A}:=\displaystyle\bigcup_{k=0}^{\infty} A_k$. Then from the countable subadditivity of the Lebesgue measure we get that $|\tilde{A}|\leq\displaystyle\sum_{k=0}^{\infty} |A_k|=0$. So $\tilde{A}\subset (0,T)\times\Omega$ is a null-measure set. Thence we can write that:

\begin{equation}
	\varepsilon\leq v_{n_k}(t,x)\leq \delta,\ \text{for any}\ (t,x)\in \bigl [(0,T)\times\Omega\bigr]\setminus \tilde{A}.
\end{equation}

\noindent Making $k\to\infty$ in this relation yields that $\varepsilon\leq v(t,x)\leq \delta$ for any $(t,x)\in \bigl [(0,T)\times\Omega\bigr]\setminus \tilde{A}$. In conclusion for a.e. $(t,x)\in (0,T)\times\Omega$ we have that $\varepsilon\leq v(t,x)\leq \delta$, i.e. $v\in\mathcal{V}_{[\varepsilon,\delta]}$.
	
\end{proof}

\noindent\textbf{\underline{Fact III}: The sequence of positive real numbers $\left(\displaystyle\int_{0}^T\mathcal{A}(v_n(t,\cdot))\ dt \right)_{n\geq 1}$ is bounded, and lies in some interval $[0,C_{\mathcal{A}}]$.}

\begin{proof} From the definition of $\mathcal{A}$ we have for every $n\geq 1$ and a.e. $t\in (0,T)$ that:
	
	\begin{align}\label{primaturamarginireaA}
	0\leq	\mathcal{A}(v_n(t,\cdot))&=\int_{\Omega} A(x,\nabla v_n(t,x))\ dx=\int_{\Omega}\int_{0}^{|\nabla v_n(t,x)|} \Phi(x,s)\ ds\ dx\nonumber\\
\textbf{(H5)}\ \ \ &\leq \int_{\Omega}\int_{0}^{|\nabla v_n(t,x)|} \Phi(x,|\nabla v_n(t,x)|)\ ds\ dx=\int_{\Omega} \Phi(x,|\nabla v_n(t,x)|)\cdot |\nabla v_n(t,x)|\ dx\nonumber\\
&=\int_{\Omega} \mathbf{a}(x,\nabla v_n(t,x))\cdot \nabla v_n(t,x)\ dx\nonumber\\
\eqref{ecuatialuivn}\ \ \ &=\int_{\Omega} \left (g(t,x)-\lambda b(x,v_n(t,x))-\dfrac{\partial b(x, v_n(t,x))}{\partial t}\right )v_n(t,x)\ dx\nonumber\\
&=\int_{\Omega} \bigl(g(t,x)-\lambda b(x,v_n(t,x))\bigr)v_n(t,x)\ dx-\int_{\Omega}\dfrac{\partial b(x, v_n(t,x))}{\partial t}v_n(t,x)\ dx\nonumber\\
&\leq \int_{\Omega} \bigl(|g(t,x)|+\lambda b(x,v_n(t,x))\bigr)v_n(t,x)\ dx-\int_{\Omega}\dfrac{\partial b(x, v_n(t,x))}{\partial t}v_n(t,x)\ dx\nonumber\\
\text{(Remark \ref{rem23})}\ \ \ &\leq \int_{\Omega}\delta\Vert g\Vert_{L^{\infty}((0,T)\times\Omega)}+\delta\lambda \Vert b(\cdot,\delta)\Vert_{L^{\infty}(\Omega)}\ dx-\int_{\Omega}\dfrac{\partial b(x, v_n(t,x))}{\partial t}v_n(t,x)\ dx\nonumber\\
&=\delta|\Omega|\bigl(\Vert g\Vert_{L^{\infty}((0,T)\times\Omega)}+ \lambda \Vert b(\cdot,\delta)\Vert_{L^{\infty}(\Omega)}\bigr)-\int_{\Omega}\dfrac{\partial b(x, v_n(t,x))}{\partial t}v_n(t,x)\ dx.
\end{align}

\noindent Now, since $\mathcal{A}:W^{1,p(x)}(\Omega)\to [0,\infty)$ is a continuous functional -- $\mathcal{A}\in C^1\bigl(W^{1,p(x)}(\Omega)\bigr)$ -- and $v_n:(0,T)\to W^{1,p(x)}(\Omega)$ is a strongly measurable function -- $v_n\in L^{\infty}\bigl(0,T;W^{1,p(x)}(\Omega)\bigr)$ -- we get that $\mathcal{A}\circ v_n:(0,T)\to [0,\infty)$ is a (strongly) measurable function.\footnote{Strong measurability and Lebesgue measurability are equivalent for a function $w:(0,T)\to\mathbb{R}$. See Proposition 9, page 60 and the \textit{Simple aproximation theorem} given at page 62, both in \cite{royden2022real}.} 

\medskip

\noindent Now, from Remark \ref{rembo} we have that $b(\cdot,v_n(\cdot,\cdot))\in H^1\bigl((0,T);L^2(\Omega)\bigr)$ and therefore we can write that:

\begin{equation}
	\begin{cases} \dfrac{\partial b(\cdot,v_n(\cdot,\cdot))}{\partial t}\in L^2\bigl(0,T;L^2(\Omega)\bigr)\simeq L^2\bigl((0,T)\times\Omega\bigr)\\ v_n\in L^{\infty}\bigl(0,T;W^{1,p(x)}(\Omega)\bigr)\subset L^{2}\bigl(0,T;W^{1,p(x)}(\Omega)\bigr)\subset L^{2}\bigl(0,T;L^{2}(\Omega)\bigr)\simeq L^2\bigl((0,T)\times\Omega\bigr).  \end{cases}.
\end{equation}

\noindent Then, from \textit{Cauchy inequality} we get that $\dfrac{\partial b(\cdot,v_n(\cdot,\cdot))}{\partial t}v_n\in L^1\bigl((0,T)\times\Omega\bigr)$. Now using \textit{Fubini's Theorem}\footnote{For statement and proof see \cite[Theorem 2.3.50, page 126]{papageorgiou2018applied}.} we conclude that the function $(0,T)\ni t\mapsto \displaystyle\int_\Omega\dfrac{\partial b(x,v_n(t,x))}{\partial t}v_n(t,x)\ dx$ is from $L^1(0,T)$. Henceforth, the right-hand side of \eqref{primaturamarginireaA} is a function from $L^1(0,T)$. Then from \eqref{primaturamarginireaA} we get that $\mathcal{A}\circ v_n\in L^1(0,T)$.

\medskip

\noindent Integrating \eqref{primaturamarginireaA} on $(0,T)$, and the using Proposition \ref{propofrakB} \textbf{(5)},\textbf{(1)} gives us that:

\begin{align*}
&0\leq \int_{0}^T\mathcal{A}(v_n(t,\cdot))\ dt\leq \delta T|\Omega|\bigl(\Vert g\Vert_{L^{\infty}((0,T)\times\Omega)}+ \lambda \Vert b(\cdot,\delta)\Vert_{L^{\infty}(\Omega)}\bigr)-\int_0^T\int_{\Omega}\dfrac{\partial b(x, v_n(t,x))}{\partial t}v_n(t,x)\ dx\ dt\\
&=\delta T|\Omega|\bigl(\Vert g\Vert_{L^{\infty}((0,T)\times\Omega)}+ \lambda \Vert b(\cdot,\delta)\Vert_{L^{\infty}(\Omega)}\bigr)-\int_{\Omega}\underbrace{\frak{B}(x,v_n(T,x))}_{\geq 0}\ dx+\int_{\Omega} \frak{B}(x,v_n(0,x))\ dx\\
&\leq \delta T|\Omega|\bigl(\Vert g\Vert_{L^{\infty}((0,T)\times\Omega)}+ \lambda \Vert b(\cdot,\delta)\Vert_{L^{\infty}(\Omega)}\bigr)+\int_{\Omega} \frak{B}(x,u_{0,n}(x))\ dx\\
&\leq \delta T|\Omega|\bigl(\Vert g\Vert_{L^{\infty}((0,T)\times\Omega)}+ \lambda \Vert b(\cdot,\delta)\Vert_{L^{\infty}(\Omega)}\bigr)+\int_{\Omega} \dfrac{L_0}{2}\bigl(u_{0,n}^2(x)-\varepsilon^2\bigr)\ dx\\
&\leq \delta T|\Omega|\bigl(\Vert g\Vert_{L^{\infty}((0,T)\times\Omega)}+ \lambda \Vert b(\cdot,\delta)\Vert_{L^{\infty}(\Omega)}\bigr)+\dfrac{L_0}{2}\Vert u_{0,n}\Vert^2_{L^2(\Omega)}\\
&\leq \delta T|\Omega|\bigl(\Vert g\Vert_{L^{\infty}((0,T)\times\Omega)}+ \lambda \Vert b(\cdot,\delta)\Vert_{L^{\infty}(\Omega)}\bigr)+\dfrac{L_0}{2}\sup_{n\geq 1} \Vert u_{0,n}\Vert^2_{L^2(\Omega)}:=C_{\mathcal{A}}<\infty,
\end{align*} 

\noindent because from the fact that $u_{0,n}\to u_0$ in $L^2(\Omega)$ we get that $(u_{0,n})_{n\geq 1}$ is a bounded sequence in $L^2(\Omega)$. The proof of this fact is now complete.
	
\end{proof}

\bigskip

\noindent\textbf{\underline{Fact IV}: For any $\tau\in (0,T)$, the sequence $\left (\dfrac{\partial v_n}{\partial t} \right )_{n\geq 1}$ is bounded in $L^2\bigl (\tau,T;L^2(\Omega)\bigr)$.}

\begin{proof} We start by using Proposition \ref{propomathcalA1} \textbf{(4)} to deduce that for every $w\in W^{1,p(x)}(\Omega)$, each $n\geq 1$ and any $t\in [0,T]$ one has that:
	
	\begin{equation}\label{mere1}
		\mathcal{A}(w)-\mathcal{A}(v_n(t,\cdot))\geq \int_{\Omega} \mathbf{a}(x,\nabla v_n(t,x))\cdot\bigl(\nabla w(x)-\nabla v_n(t,x)\bigr)\ dx.
	\end{equation}
	
\noindent Now if in equation \eqref{ecuatialuivn} we choose $\phi:=w-v_n(t,\cdot)\in W^{1,p(x)}(\Omega)$, we get that:

\begin{align}\label{mere2}
	&\int_{\Omega}\mathbf{a}(x,\nabla v_n(t,x))\cdot\bigl (\nabla w(x)-\nabla v_n(t,x)\bigr)\ dx=\int_{\Omega} g(t,x)\bigl (w(x)-v_n(t,x)\bigr)\ dx\nonumber\\
	&\ \ \ -\int_{\Omega}\dfrac{\partial b(x,v_n)}{\partial t}\bigl (w(x)-v_n(t,x)\bigr)\ dx-\lambda\int_{\Omega}b(x,v_n)\bigl (w(x)-v_n(t,x)\bigr)\ dx \nonumber\\
	=\ &\int_{\Omega} \left (g(t,x)-\lambda b(x,v_n(t,x))-\dfrac{\partial b(x,v_n(t,x))}{\partial t} \right )\cdot \bigl(w(x)-v_n(t,x)\bigr)\ dx.
\end{align}

\noindent Combining \eqref{mere1} and \eqref{mere2} leads us to:

\begin{equation}\label{mere3}
	\mathcal{A}(w)-\mathcal{A}(v_n(t,\cdot))\geq \int_{\Omega} \left (g-\lambda b(x,v_n(t,x))-\dfrac{\partial b(x,v_n)}{\partial t} \right )\cdot \bigl(w(x)-v_n(t,x)\bigr)\ dx.
\end{equation}

\noindent If we denonte $h_n=g-\lambda b(\cdot,v_n)-\dfrac{\partial b(\cdot,v_n)}{\partial t}\in L^2\bigl(0,T;L^2(\Omega)\bigr)$ (See Remark \ref{rembo}), we get that for a.e. $t\in (0,T)$, $h_n(t,\cdot)\in L^2(\Omega)$. Therefore, taking into account Proposition \ref{propoverA} \textbf{(3)}\footnote{Here we used the extra assumption of Theorem \ref{theoremverygeneral} that $p$ is a log-H\"{o}lder continuous exponent.}, relation \eqref{mere3} says literally that $h_n(t,\cdot)\in \partial\overline{\mathcal{A}}(v_n(t,\cdot))$ for a.e. $t\in (0,T)$.

\medskip
	
\noindent Recall that $\overline{\mathcal{A}}:L^2(\Omega)\to [0,\infty]$ is a convex, lower semicontinuous and proper functinal, as proved in Proposition \ref{propoverA}. Also, for each $n\geq 1$ we have that $v_n\in C\bigl([0,T];L^2(\Omega)\bigr)\cap H^1\bigl((0,T);L^2(\Omega)\bigr)\cap L^{\infty}\bigl(0,T;W^{1,p(x)}(\Omega)\bigr)$ -- see \eqref{relvninspatii} --. Now we have all we need in order to apply Proposition \ref{propoleoni} from the Appendix for $\overline{\mathcal{A}}$, $v_n$ and $h_n$. Thus we obtain that:

\begin{equation}\label{mere4}
	\begin{cases} \alpha_n:=\mathcal{A}\circ v_n=\overline{\mathcal{A}}\circ v_n\in \textnormal{AC}\bigl([0,T]\bigr) \ \text{and}\\[3mm]
	\alpha'_n(t)=\bigl(\mathcal{A}\circ v_n\bigr)'(t)=\bigl(\overline{\mathcal{A}}\circ v_n\bigr)'(t)=\displaystyle\int_{\Omega} h_n(t,x)\dfrac{\partial v_n}{\partial t}(t,x)\ dx\ \text{for a.e.}\ t\in (0,T).\end{cases}
\end{equation}

\noindent From \eqref{mere4} we can write for a.e. $t\in (0,T)$ that:

\begin{equation}\label{mere5}
\alpha_n'(t)+\int_{\Omega} \dfrac{\partial b(x,v_n(t,x))}{\partial t}\cdot\dfrac{\partial v_n}{\partial t}\ dx=\int_{\Omega} \bigl(g-\lambda b(x,v_n(t,x))\bigr)\cdot\dfrac{\partial v_n}{\partial t}\ dx.
\end{equation}

\noindent Multiplying by $t$ we get for a.e. $t\in (0,T)$ that:

\begin{equation}\label{mere6}
	t\alpha_n'(t)+\int_{\Omega} t\dfrac{\partial b(x,v_n(t,x))}{\partial t}\cdot\dfrac{\partial v_n}{\partial t}\ dx=\int_{\Omega} t\bigl(g-\lambda b(x,v_n(t,x))\bigr)\cdot\dfrac{\partial v_n}{\partial t}\ dx.
\end{equation}

\noindent Next, we will explain why each of the three terms in \eqref{mere6} is a function from $L^1(0,T)$, so that we may integrate on $(0,T)$ \eqref{mere6}.

\bigskip

\noindent$\blacktriangleright$ We know that  $\alpha_n\in \operatorname{AC}\bigl([0,T]\bigr)$. From Theorem 12 given in \cite[Section 6.5, page 107]{royden2022real} we get that $\alpha_n'\in L^1(0,T)$. Since the identity function $[0,T]\ni t\mapsto t\in [0,T]$ is from $L^{\infty}(\Omega)$, we deduce that $(0,T)\ni t\mapsto t\alpha_n'(t)$ is a measurable function (being a product of measurable functions) and in fact a member of $L^1(0,T)$, because $\displaystyle\int_{0}^T |t\alpha_n'(t)|\ dt\leq T\int_{0}^T|\alpha_n'(t)|\ dt=T\Vert\alpha_n'\Vert_{L^1(0,T)}<\infty$. So it makes sense to talk about $\displaystyle\int_{0}^T t\alpha'_n(t)\ dt$.

\bigskip

\noindent$\blacktriangleright$ We know that $\alpha_n\in \operatorname{AC}\bigl([0,T]\bigr)$, which in particular implies that $\alpha_n$ is continuous on the compact interval $[0,T]$. From \textit{Weierstrass Theorem} we get that $\alpha_n$ is a bounded function, i.e. $\alpha_n\in L^{\infty}(0,T)\subset L^1(0,T)$. So it makes sense to talk about $\displaystyle\int_{0}^T \alpha_n(t)\ dt$.

\begin{remark}\label{partsint1}
	Since $[0,T]\ni t\mapsto [0,T]$ and $\alpha_n$ are both from $\operatorname{AC}\bigl([0,T]\bigr)$, we deduce from the \textbf{Integration by parts formula for AC functions} -- see Theorem 19 in \cite[Section 6.5, page 112]{royden2022real} --  that
	
	\begin{equation}\label{intpartsformula1}
		\int_{0}^Tt\alpha_n'(t)\ dt=T\cdot\alpha_n(T)-0\cdot\alpha_n(0)-\int_{0}^Tt'\alpha_n(t)\ dt=T\cdot\alpha_n(T)-\int_{0}^T\alpha_n(t)\ dt.
	\end{equation}
\end{remark}
\bigskip

\noindent$\blacktriangleright$ We have that

\begin{equation}
	\begin{cases}(0,T)\times\Omega\ni (t,x)\mapsto t\in (0,T)\ \text{is from}\ L^{\infty}\bigl((0,T)\times\Omega\bigr)\\[3mm] \dfrac{\partial b(\cdot, v_n)}{\partial t}\in L^2\bigl(0,T;L^2(\Omega)\bigr)\simeq L^{2}\bigl((0,T)\times\Omega\bigr) \text{ -- see Remark \ref{rembo}}\\[3mm] \dfrac{\partial v_n}{\partial t}\in L^2\bigl(0,T;L^2(\Omega)\bigr)\simeq L^{2}\bigl((0,T)\times\Omega\bigr).\end{cases}
\end{equation}

\noindent Therefore, from \textit{Cauchy inequality} we deduce that $\dfrac{\partial b(\cdot, v_n)}{\partial t}\cdot\dfrac{\partial v_n}{\partial t}\in L^{1}\bigl((0,T)\times\Omega\bigr)$, and from \textit{H\"{o}lder inequality} we get that $t\dfrac{\partial b(\cdot, v_n)}{\partial t}\cdot\dfrac{\partial v_n}{\partial t}\in L^{1}\bigl((0,T)\times\Omega\bigr)$. So, from \textit{Fubini's theorem}, it makes sense to talk about $\displaystyle\int_{0}^T\int_{\Omega} t\dfrac{\partial b(x, v_n)}{\partial t}\cdot\dfrac{\partial v_n}{\partial t}\ dx\ dt$.

\bigskip

\noindent$\blacktriangleright$ Similarly

\begin{equation}
	\begin{cases}(0,T)\times\Omega\ni (t,x)\mapsto t\in (0,T)\ \text{is from}\ L^{\infty}\bigl((0,T)\times\Omega\bigr)\\[3mm] g-\lambda b(\cdot,v_n)\in L^2\bigl(0,T;L^2(\Omega)\bigr)\simeq L^{2}\bigl((0,T)\times\Omega\bigr) \text{ -- see Remark \ref{rembo}}\\[3mm] \dfrac{\partial v_n}{\partial t}\in L^2\bigl(0,T;L^2(\Omega)\bigr)\simeq L^{2}\bigl((0,T)\times\Omega\bigr).\end{cases}
\end{equation}

\noindent Thus, from \textit{Cauchy inequality} and \textit{H\"{o}lder inequality}, we deduce that  $t\bigl(g-\lambda b(\cdot,v_n)\bigr)\cdot\dfrac{\partial v_n}{\partial t}\in L^{1}\bigl((0,T)\times\Omega\bigr)$. So, from \textit{Fubini's theorem}, it makes sense to talk about $\displaystyle\int_{0}^T\int_{\Omega} t\bigl(g(t,x)-\lambda b(x,v_n)\bigr)\cdot\dfrac{\partial v_n}{\partial t}\ dx\ dt$.

\bigskip

\noindent Taking into account the above discussion we can integrate \eqref{mere6} on $(0,T)$ and obtain, using \eqref{intpartsformula1}, that

\begin{equation}\label{mere7}
	T\cdot\alpha_n(T)-\int_{0}^T\alpha_n(t)\ dt+\int_{0}^T\int_{\Omega} t\dfrac{\partial b(x,v_n(t,x))}{\partial t}\cdot\dfrac{\partial v_n}{\partial t}\ dx\ dt=\int_{0}^T\int_{\Omega} t\bigl(g-\lambda b(x,v_n(t,x))\bigr)\cdot\dfrac{\partial v_n}{\partial t}\ dx\ dt.
\end{equation}

\noindent From here we get that:

\begin{align}\label{mere8}
	\int_{0}^T\int_{\Omega} t\dfrac{\partial b(x,v_n)}{\partial t}\cdot\dfrac{\partial v_n}{\partial t}\ dx\ dt&=\int_{0}^T\int_{\Omega} t\bigl(g-\lambda b(x,v_n)\bigr)\cdot\dfrac{\partial v_n}{\partial t}\ dx\ dt+\underbrace{\int_{0}^T\alpha_n(t)\ dt}_{\leq C_{\mathcal{A}}}-T\cdot\underbrace{\alpha_n(T)}_{\geq 0}\nonumber\\
\text{(Fact III, Remark \ref{rem23})}\ \ \ 	&\leq C_{\mathcal{A}}+\int_{0}^T\int_{\Omega} \underbrace{\bigl(\Vert g\Vert_{L^{\infty}((0,T)\times\Omega)}+\lambda \Vert b(\cdot,\delta)\Vert_{L^{\infty}(\Omega)}\bigr)}_{:=\widetilde{C}}t\left |\dfrac{\partial v_n}{\partial t}\right |\ dx\ dt\nonumber\\
&=C_{\mathcal{A}}+\widetilde{C}\int_{0}^T\int_{\Omega} t\left |\dfrac{\partial v_n}{\partial t}\right |\ dx\ dt\\
\text{(AM-GM ineq.)}\ \ \ &\leq C_{\mathcal{A}}+\widetilde{C}\int_{0}^T\int_{\Omega} \epsilon t\left |\dfrac{\partial v_n}{\partial t}\right |^2+\dfrac{t}{4\epsilon}\ dx\ dt\nonumber\\
&=C_{\mathcal{A}}+\dfrac{\widetilde{C}T^2|\Omega|}{8\epsilon}+\widetilde{C}\epsilon\int_{0}^T\int_{\Omega} t\left |\dfrac{\partial v_n}{\partial t}\right |^2\ dx\ dt,\ \forall\ \epsilon>0.
\end{align}

\noindent On the other hand, from Proposition \ref{propofrakB} \textbf{(3)} we have that:

\begin{align}\label{mere9}
	\int_{0}^T\int_{\Omega} t\dfrac{\partial b(x,v_n)}{\partial t}\cdot\dfrac{\partial v_n}{\partial t}\ dx\ dt&=	\int_{0}^T\int_{\Omega} t\dfrac{\partial b}{\partial s}(x,v_n)\cdot\left |\dfrac{\partial v_n}{\partial t}\right|^2\ dx\ dt\nonumber\\
\textbf{(H13)}\ \ \ 	&\geq \ell_0\int_{0}^T\int_{\Omega} t\left |\dfrac{\partial v_n}{\partial t}\right|^2\ dx\ dt.
\end{align}

\noindent Putting together \eqref{mere8} and \eqref{mere9} leads us to:

\begin{equation}\label{mere10}
\ell_0\int_{0}^T\int_{\Omega} t\left |\dfrac{\partial v_n}{\partial t}\right|^2\ dx\ dt\leq C_{\mathcal{A}}+\dfrac{\widetilde{C}T^2|\Omega|}{8\epsilon}+\widetilde{C}\epsilon\int_{0}^T\int_{\Omega} t\left |\dfrac{\partial v_n}{\partial t}\right |^2\ dx\ dt,\ \forall\ \epsilon>0.
\end{equation}

\noindent We can rewrite \eqref{mere10} in the following form:

\begin{equation}\label{mere11}
\bigl(\ell_0-\widetilde{C}\epsilon\bigr)\int_{0}^T\int_{\Omega} t\left |\dfrac{\partial v_n}{\partial t}\right|^2\ dx\ dt\leq C_{\mathcal{A}}+\dfrac{\widetilde{C}T^2|\Omega|}{8\epsilon},\ \forall\ \epsilon>0.
\end{equation}

\noindent Choosing now $\epsilon:=\dfrac{\ell_0}{2\widetilde{C}}$ in \eqref{mere11} gives us that:

\begin{equation}\label{mere12}
	\dfrac{\ell_0}{2}\int_{0}^T\int_{\Omega} t\left |\dfrac{\partial v_n}{\partial t}\right|^2\ dx\ dt\leq C_{\mathcal{A}}+\dfrac{\widetilde{C}^2T^2|\Omega|}{4\ell_0}\ \Longrightarrow \ \int_{0}^T\int_{\Omega} t\left |\dfrac{\partial v_n}{\partial t}\right|^2\ dx\ dt\leq\dfrac{2C_{\mathcal{A}}}{\ell_0}+\dfrac{\widetilde{C}^2T^2|\Omega|}{2\ell_0^2}.
\end{equation}

\noindent For any $\tau\in (0,T)$, from \eqref{mere12} we deduce that:

\begin{equation}
	\tau\int_{\tau}^T\int_{\Omega}\left |\dfrac{\partial v_n}{\partial t}\right|^2\ dx\ dt\leq \int_{\tau}^T\int_{\Omega}t\left |\dfrac{\partial v_n}{\partial t}\right|^2\ dx\ dt\leq \int_{0}^T\int_{\Omega} t\left |\dfrac{\partial v_n}{\partial t}\right|^2\ dx\ dt\leq\dfrac{2C_{\mathcal{A}}}{\ell_0}+\dfrac{\widetilde{C}^2T^2|\Omega|}{2\ell_0^2}.
\end{equation}

\noindent Thus we proved that:

\begin{equation}\label{pere1}
	\left\Vert\dfrac{\partial v_n}{\partial t}\right\Vert_{L^2(\tau,T;L^2(\Omega))}^2\leq \dfrac{2C_{\mathcal{A}}}{\tau\ell_0}+\dfrac{\widetilde{C}^2T^2|\Omega|}{2\tau\ell_0^2},\ \forall\ n\geq 1.
\end{equation}

\noindent Relation \eqref{pere1} proves that $\left(\dfrac{\partial v_n}{\partial t}\right)_{n\geq 1}$ is a bounded sequence from $L^2\bigl(\tau,T;L^2(\Omega)\bigr)$ for any $\tau\in (0,T)$.

\end{proof}

\noindent\textbf{\underline{Fact V}: For any $\tau\in (0,T)$ we have that $\dfrac{\partial v}{\partial t}\in L^2\bigl(\tau,T;L^2(\Omega)\bigr)$ and moreover $\dfrac{\partial v_n}{\partial t}\weak \dfrac{\partial v}{\partial t}$ in $L^2\bigl(\tau,T;L^2(\Omega)\bigr)$. Thereby $v\in H^1_{\text{loc}}\bigl((0,T);L^2(\Omega)\bigr)$.}

\begin{proof}From Proposition 2.2.3 (c) given in \cite[page 127]{gasinski2005nonlinear} we get that $L^2\bigl(\tau,T;L^2(\Omega)\bigr)$ is a reflexive Banach space. Also, from \textbf{Fact IV} we know that the sequence $\left (\dfrac{\partial v_{n}}{\partial t}\right )_{n\geq 1}$ is bounded in $L^2\bigl(\tau,T;L^2(\Omega)\bigr)$. Take any weakly convergent subsequence to some function, say $w\in L^2\bigl(\tau,T;L^2(\Omega)\bigr)$. So we have $\dfrac{\partial v_{n_k}}{\partial t}\weak w$ in $L^2\bigl(\tau,T;L^2(\Omega)\bigr)$.
	
	\medskip
	
	\noindent Now fix any $\varphi\in C^{\infty}_{c}(\tau,T)$. Using the fact that $v_{n_k}\in H^1\bigl((0,T);L^2(\Omega)\bigr)\subset H^1\bigl((\tau,T);L^2(\Omega)\bigr)$ and the definition of the weak derivative we have for each $k\geq 1$ that:
	
	\begin{equation}\label{weakderivativeeq10}
		\int_{\tau}^T v_{n_k}(t,\cdot)\varphi'(t)\ dt=-\int_{\tau}^T \dfrac{\partial v_{n_k}}{\partial t}(t,\cdot)\varphi(t)\ dt.
	\end{equation}
	
	\noindent Define the operators $\mathcal{T}_{\varphi},\mathcal{T}_{\varphi'}:L^2\bigl(\tau,T;L^2(\Omega)\bigr)\to L^2(\Omega)$ via the following Bochner integrals:
	
	\begin{equation}
		\mathcal{T}_{\varphi}(u)=\int_{\tau}^T u(t,\cdot)\varphi(t)\ dt,\ \text{and}\ \mathcal{T}_{\varphi'}(u)=\int_{\tau}^T u(t,\cdot)\varphi'(t)\ dt,\ \forall\ u\in L^2\bigl(\tau,T;L^2(\Omega)\bigr).
	\end{equation}
	
	\noindent From Lemma \ref{bochnercontinuous} we have that $\mathcal{T}_{\varphi}$ and $\mathcal{T}_{\varphi'}$ are well-defined bounded linear operators. This is because $\varphi,\varphi'\in C^{\infty}_c(\tau,T)\subset L^2(\tau,T)$. Now \eqref{weakderivativeeq10} rewrites as:
	
	\begin{equation}
		\mathcal{T}_{\varphi'}(v_{n_k})=-\mathcal{T}_{\varphi}\left (\dfrac{\partial v_{n_k}}{\partial t} \right), \ \forall\ k\geq 1.
	\end{equation}
	
	\noindent From Remark \ref{vnvlr0tl2omega} we deduce that $v_{n_k}\to v$ in $L^2\bigl(0,T;L^2(\Omega)\bigr)$. In particular it follows that $v_{n_k}\to v$ in $L^2\bigl(\tau,T;L^2(\Omega)\bigr)$, and therefore $v_{n_k}\weak v$ in $L^2\bigl(\tau,T;L^2(\Omega)\bigr)$. We also know that $\dfrac{\partial\tilde{v}_{n_k}}{\partial t}\weak w$ in $L^2\bigl(\tau,T;L^2(\Omega)\bigr)$ and taking into account that continuous linear operators preserve weak convergence (i.e. Lemma \ref{lemmaboundedweak} from the Appendix) we finally get that:
	
	\begin{equation}
		\int_{\tau}^T v(t,\cdot)\varphi'(t)\ dt=\mathcal{T}_{\varphi'}(v)=\lim\limits_{k\to\infty} \mathcal{T}_{\varphi'}(v_{n_k})=-\lim\limits_{k\to\infty} \mathcal{T}_{\varphi}\left (\dfrac{\partial v_{n_k}}{\partial t} \right)=-\mathcal{T}_{\varphi}(w)=-\int_{\tau}^T w(t,\cdot)\varphi(t)\ dt.
	\end{equation}
	
	\noindent Since $\varphi\in C^\infty_{c}(\tau,T)$ was fixed arbitrarily we get that $v\in L^2\bigl(\tau,T;L^2(\Omega)\bigr)$ has a weak derivative and $\dfrac{\partial v}{\partial t}=w\in L^2\bigl(\tau,T;L^2(\Omega)\bigr)$. Since $v\in C\bigl([0,T];L^2(\Omega)\bigr)\subset L^2\bigl(0,T;L^2(\Omega)\bigr)\subset L^2\bigl(\tau,T;L^2(\Omega)\bigr)$, this completely shows that $v\in H^1\bigl((\tau,T);L^2(\Omega)\bigr)$ and $\dfrac{\partial v_{n_k}}{\partial t}\weak \dfrac{\partial v}{\partial t}$ in $L^2\bigl(\tau,T;L^2(\Omega)\bigr)$.
	
	\noindent We showed that any weakly convergent subsequence of $\left(\dfrac{\partial v_n}{\partial t}\right)_{n\geq 1}\subset L^2\bigl(\tau,T;L^2(\Omega)\bigr)$ weakly converges to $\dfrac{\partial v}{\partial t}\in  L^2\bigl(\tau,T;L^2(\Omega)\bigr)$. Using Lemma \ref{weakconvergencelemma} we conclude that $\dfrac{\partial v_{n}}{\partial t}\weak \dfrac{\partial v}{\partial t}$ in $L^2\bigl(\tau,T;L^2(\Omega)\bigr)$.

\end{proof}

\bigskip

\noindent\textbf{\underline{Fact VI}: For a.e. $t\in (0,T)$ we have that $v(t,\cdot)\in W^{1,p(x)}(\Omega)$ and moreover, for a.e. $t\in (0,T)$ there is a subsequence $(v_{n_k(t)}(t,\cdot))_{k\geq 1}$ such that $v_{n_k(t)}(t,\cdot)\weak v(t,\cdot)$ in $W^{1,p(x)}(\Omega)$ as $k\to\infty$.}

\begin{proof} Remeber that $\alpha_n:[0,T]\to [0,\infty),\ \alpha_n(t)=\mathcal{A}\bigl(v_n(t,\cdot)\bigr)$ is a function from $\operatorname{AC}\bigl([0,T]\bigr)$, for each $n\geq 1$. We define 
	
\begin{equation}
\alpha:[0,T]\to [0,\infty],\ \alpha(t)=\liminf\limits_{n\to\infty} \alpha_n(t),\ \forall\ t\in [0,T].
\end{equation}

\noindent First, we point out that $\alpha$ is a measurable function, being a limit inferior of a sequence of measurable functions. Then, \textit{Fatou inequality}\footnote{The statement and its proof can be found in \cite[Theorem 2.4.4, page 63]{Cohn}.} and \textbf{Fact III} gives us that:

\begin{equation}\label{lumanare3}
	\int_0^T\alpha(t)\ dt=\int_{0}^T\liminf\limits_{n\to\infty}\alpha_n(t)\ dt\leq \liminf\limits_{n\to\infty}\int_{0}^T \alpha_n(t)\ dt\leq C_{\mathcal{A}}<\infty.
\end{equation}

\noindent From this relation we find out that for a.e. $t\in (0,T)$ we have that $\alpha(t)\in [0,\infty)$. So we can find some representative $\alpha:(0,T)\to [0,\infty)$.

\noindent So, for a.e. $t\in (0,T)$, the definition of limit inferior tells us that there is a subsequence (which may depend on $t$) such that $\lim\limits_{k\to\infty}\alpha_{n_k(t)}(t)=\lim\limits_{k\to\infty}\mathcal{A}\bigl(v_{n_k(t)}(t,\cdot)\bigr)=\alpha(t)\in (0,\infty)$. Therefore the sequence $\left(\mathcal{A}\bigl(v_{n_k(t)}(t,\cdot)\bigr)\right)_{k\geq 1}$ is bounded. But, from Remark \ref{remboundednormmodular} this means that the sequence $\left(\rho_{p(x)}\bigl(|\nabla v_{n_k(t)}(t,\cdot)|\bigr)\right)_{k\geq 1}$ is also bounded i.e. $\left(\Vert\nabla v_{n_k(t)}(t,\cdot)\Vert_{L^{p(x)}(\Omega)^N}\right)_{k\geq 1}$ is bounded. Taking into account that $(v_n)_{n\geq 1}\subset C\bigl([0,T];L^2(\Omega)\bigr)\cap\mathcal{V}_{[\varepsilon,\delta]}$, from Lemma \ref{contl2allt} we obtain that for each $k\geq 1:\ \varepsilon\leq v_{n_k(t)}(t,x)\leq \delta$ for a.e. $x\in\Omega$. Hence we may write that:

\begin{align*}
	\Vert v_{n_k(t)}(t,\cdot)\Vert_{W^{1,p(x)}(\Omega)}&=	\Vert v_{n_k(t)}(t,\cdot)\Vert_{L^{p(x)}(\Omega)}+\Vert\nabla v_{n_k(t)}(t,\cdot)\Vert_{L^{p(x)}(\Omega)^N}\\
	&\leq \delta \Vert 1\Vert_{L^{p(x)}(\Omega)}+\Vert\nabla v_{n_k(t)}(t,\cdot)\Vert_{L^{p(x)}(\Omega)^N}\\
	&\leq \delta \Vert 1\Vert_{L^{p(x)}(\Omega)}+\sup_{k\geq 1}\Vert\nabla v_{n_k(t)}(t,\cdot)\Vert_{L^{p(x)}(\Omega)^N}<\infty.
\end{align*}

\noindent This relation shows us that the sequence $\bigl(v_{n_k(t)}(t,\cdot)\bigr)_{k\geq 1}$ is bounded in $W^{1,p(x)}(\Omega)$, which is a reflexive Banach space\footnote{See Proposition \ref{propospatii} \textbf{(1)}.}. In this context we want to apply Lemma \ref{weakconvergencelemma}. In that sense consider any subsequence $\left(v_{n_{k_\ell}(t)}(t,\cdot)\right)_{\ell\geq 1}$ that weakly converges to some element $z\in W^{1,p(x)}(\Omega)\hookrightarrow L^2(\Omega)$ (from \textbf{(H2)}) -- one such subsequence exists from \textit{Eberlein-\v{S}mulian theorem}. Applying Lemma \ref{lemmaboundedweak} we deduce that $v_{n_{k_\ell}(t)}(t,\cdot)\weak z$ in $L^2(\Omega)$.

\noindent Knowing that $v_n\to v$ in $C\bigl([0,T];L^2(\Omega)\bigr)$ we get that $\lim\limits_{n\to\infty} \displaystyle\sup_{\tau\in [0,T]} \Vert v_n(\tau,\cdot)-v(\tau,\cdot)\Vert_{L^2(\Omega)}=0$. So we also have that $\lim\limits_{n\to\infty} \Vert v_n(t,\cdot)-v(t,\cdot)\Vert_{L^2(\Omega)}=0$. This means that $v_n(t,\cdot)\to v(t,\cdot)$ in $L^2(\Omega)$, from where we deduce that $v_n(t,\cdot)\weak v(t,\cdot)$ in $L^2(\Omega)$. In conclusion $\begin{cases} v_{n_{k_\ell}(t)}(t,\cdot)\weak z,\ \text{in}\ L^2(\Omega)\\ v_{n_{k_\ell}(t)}(t,\cdot)\weak v(t,\cdot),\ \text{in}\ L^2(\Omega)\end{cases}$, and from the uniqueness of the weak limit, we get that $v(t,\cdot)=z\in W^{1,p(x)}(\Omega)$. This shows two things: one is that $v(t,\cdot)\in W^{1,p(x)}(\Omega)$ for a.e. $t\in (0,T)$, and the second one, using Lemma \ref{lemmaboundedweak}, is that $v_{n_k(t)}(t,\cdot)\weak v(t,\cdot)$ in $W^{1,p(x)}(\Omega)$. The proof is now complete.

\end{proof}

\bigskip

\noindent\textbf{\underline{Fact VII}: The function $v:(0,T)\to W^{1,p(x)}(\Omega)$ is strongly measurable.}

\begin{proof} We know from \textbf{Fact I} that $v:[0,T]\to L^2(\Omega)$ is a continuous function (hence it is Borel measurable), and from \textbf{(H2)} that $W^{1,p(x)}(\Omega)\hookrightarrow L^2(\Omega)$, where $W^{1,p(x)}(\Omega)$ and $L^2(\Omega)$ are both separable Banach spaces. Applying Lemma \ref{stronglymeasurableembedding} we deduce that $v:(0,T)\to W^{1,p(x)}(\Omega)$ is strongly measurable.

\end{proof}

\bigskip

\noindent\textbf{\underline{Fact VIII}: The sequence $(v_n)_{n\geq 1}$ is bounded in $L^1\bigl(0,T;W^{1,p(x)}(\Omega)\bigr)$.}

\begin{proof} We already know that $(v_n)_{n\geq 1}\subset L^{\infty}\bigl(0,T;W^{1,p(x)}(\Omega)\bigr)$. Therefore we may write for each $n\geq 1$ and a.e. $t\in (0,T)$ that:
	
	\begin{align*}
		\Vert v_n(t,\cdot)\Vert_{W^{1,p(x)}(\Omega)}&=\Vert v_n(t,\cdot)\Vert_{L^{p(x)}(\Omega)}+\Vert \nabla v_n\Vert_{L^{p(x)}(\Omega)^N}\\
\text{(Lemma \ref{contl2allt})}\ \ \		&\leq \delta\Vert 1\Vert_{L^{p(x)}(\Omega)}+\Vert \nabla v_n\Vert_{L^{p(x)}(\Omega)^N}\\
\text{(Prop. \ref{propomathcalA1} \textbf{(6)})}\ \ \ &\leq \delta\Vert 1\Vert_{L^{p(x)}(\Omega)}+C_*+C_*\mathcal{A}\bigl(v_n(t,\cdot)\bigr).
	\end{align*}
	
	\noindent Integrating this inequality on $(0,T)$ leads us to:
	
	\begin{align*}
		\Vert v_n\Vert_{L^1(0,T;W^{1,p(x)}(\Omega))}&=\int_{0}^T 	\Vert v_n(t,\cdot)\Vert_{W^{1,p(x)}(\Omega)}\ dt\leq \delta T\Vert 1\Vert_{L^{p(x)}(\Omega)}+C_*T+C_*\int_{0}^T\mathcal{A}\bigl(v_n(t,\cdot)\bigr)\  dt\\
(\textbf{Fact III})\ \ \		&\leq \delta T\Vert 1\Vert_{L^{p(x)}(\Omega)}+C_*T+C_*C_{\mathcal{A}},\ \forall\ n\geq 1.
	\end{align*}
	
	\noindent Therefore the sequence $(v_n)_{n\geq 1}$ is bounded in $L^1\bigl(0,T;W^{1,p(x)}(\Omega)\bigr)$.
	
\end{proof}

\bigskip

\noindent\textbf{\underline{Fact IX}: We have that $v\in L^1\bigl(0,T;W^{1,p(x)}(\Omega)\bigr)$.}

\begin{proof} From Proposition \ref{propoverA} \textbf(2) we know that $\overline{\mathcal{A}}:L^2(\Omega)\to [0,\infty]$ is lower semicontinuous. From $v_n\to v$ in $C\bigl([0,T];L^2(\Omega)\bigr)$ we get that for any $t\in [0,T]$: $v_n(t,\cdot)\to v(t,\cdot)$ in $L^2(\Omega)$. Taking into account that from \textbf{Fact VI} we know that $v(t,\cdot)\in W^{1,p(x)}(\Omega)$ for a.e. $t\in (0,T)$ we may write that:
	
	\begin{equation}\label{lumanare2}
	\mathcal{A}\bigl(v(t,\cdot)\bigr)=\overline{\mathcal{A}}\bigl(v(t,\cdot)\bigr)\leq \liminf\limits_{n\to\infty} \overline{\mathcal{A}}\bigl(\underbrace{v_n(t,\cdot)}_{\in W^{1,p(x)}(\Omega)}\bigr)=\liminf\limits_{n\to\infty} \mathcal{A}\bigl(v_n(t,\cdot)\bigr)=\liminf\limits_{n\to\infty} \alpha_n(t)=\alpha(t),
	\end{equation}
	
	\noindent for a.e. $t\in (0,T)$. Since, from \textbf{Fact VII} we know that $v:(0,T)\to W^{1,p(x)}(\Omega)$ is strongly measurable and $\mathcal{A}:W^{1,p(x)}(\Omega)\to [0,\infty)$ is a $C^1$ functional (in particular it is continuous), we deduce that the composition function $(0,T)\ni t\mapsto \mathcal{A}\bigl(v(t,\cdot)\bigr)\in [0,\infty)$ is measurable.
	
	\noindent Integrating \eqref{lumanare2} on $(0,T)$ gives us that:
	
	\begin{equation}\label{lumanare4}
		\int_{0}^T\mathcal{A}\bigl(v(t,\cdot)\bigr)\ dt\leq \int_{0}^T\alpha(t)\ dt\stackrel{\eqref{lumanare3}}{\leq} C_{\mathcal{A}}.
	\end{equation}
	
	\noindent Also, since from \textbf{Facts I, II} we have that $v\in C\bigl([0,T];L^2(\Omega)\bigr)\cap\mathcal{V}_{[\varepsilon,\delta]}$, we get from Lemma \ref{contl2allt} that for any $t\in [0,T]$: $\varepsilon\leq v(t,x)\leq\delta$ for a.e. $x\in\Omega$. Thence we may write:

	\begin{align*}
		\Vert v(t,\cdot)\Vert_{W^{1,p(x)}(\Omega)}&=\Vert v(t,\cdot)\Vert_{L^{p(x)}(\Omega)}+\Vert \nabla v(t,\cdot)\Vert_{L^{p(x)}(\Omega)^N}\\
		&\leq \delta\Vert 1\Vert_{L^{p(x)}(\Omega)}+\Vert \nabla v(t,\cdot)\Vert_{L^{p(x)}(\Omega)^N}\\
		\text{(Prop. \ref{propomathcalA1} \textbf{(6)})}\ \ \ &\leq \delta\Vert 1\Vert_{L^{p(x)}(\Omega)}+C_*+C_*\mathcal{A}\bigl(v(t,\cdot)\bigr).
	\end{align*}
	
	\noindent Integrating on $(0,T)$ finally gives us that:
	
	\begin{align*}
	\Vert v\Vert_{L^1(0,T;W^{1,p(x)}(\Omega))}&=\int_{0}^T 	\Vert v(t,\cdot)\Vert_{W^{1,p(x)}(\Omega)}\ dt\leq \delta T\Vert 1\Vert_{L^{p(x)}(\Omega)}+C_*T+C_*\int_{0}^T\mathcal{A}\bigl(v(t,\cdot)\bigr)\  dt\\
	\eqref{lumanare4}\ \ \		&\leq \delta T\Vert 1\Vert_{L^{p(x)}(\Omega)}+C_*T+C_*C_{\mathcal{A}}.
\end{align*}

\noindent This shows that $v\in L^1\bigl(0,T;W^{1,p(x)}(\Omega)\bigr)$, as required.

\end{proof}

\bigskip

\noindent\textbf{\underline{Fact X}: $b(\cdot,v)\in C\bigl([0,T];L^2(\Omega)\bigr)$ and $b(\cdot,v_n)\to b(\cdot,v)$ in $C\bigl([0,T];L^2(\Omega)\bigr)$ as $n\to\infty$. In particular $b(\cdot,v_n)\to b(\cdot,v)$ and $b(\cdot,v_n)\weak b(\cdot,v)$ in $L^2\bigl(0,T;L^2(\Omega)\bigr)$.}

\begin{proof} First note that, since $v\in C([0,T];L^2(\Omega))\subset L^2(0,T;L^2(\Omega))$, from Proposition \ref{propoverb} \textbf{(4)} we have that $b(\cdot,v)\in L^2(0,T;L^2(\Omega))$. Moreover for any $t\in [0,T]$ we have for sufficiently small $|h|$ that:
	
	\begin{align*}
		\Vert b(\cdot,v(t+h,\cdot))-b(\cdot,v(t,\cdot))\Vert_{L^2(\Omega)}&=\left (\int_{\Omega} |b(x,v(t+h,x))-b(x,v(t,x))|^2\ dx\right)^{\frac{1}{2}}\\
		\eqref{blipschitz}\ \ \ &\leq L_0\left (\int_{\Omega} |v(t+h,x)-v(t,x)|^2\ dx\right)^{\frac{1}{2}}\\
		&=L_0\Vert v(t+h,\cdot)-v(t,\cdot)\Vert_{L^2(\Omega)}\\
		&\stackrel{h\to 0}{\longrightarrow} 0,
	\end{align*}
	
	\noindent because $v\in C([0,T];L^2(\Omega))$. This proves that $b(\cdot,v)\in C([0,T];L^2(\Omega))$.
	
	\noindent Similarly, we know that $v_{n}\in C([0,T];L^2(\Omega))$ for any $n\geq 1$, and then repeating the above argument will give us that $b(\cdot,v_{n})\in C([0,T];L^2(\Omega))$ for each $n\geq 1$. 
	
	\noindent Finally, remark that:
	
	\begin{align*}
		\Vert b(\cdot,v_{n})-b(\cdot,v)\Vert_{C([0,T];L^2(\Omega))}&=\sup_{t\in [0,T]} \Vert b(\cdot,v_{n}(t,\cdot))-b(\cdot,v(t,\cdot))\Vert_{L^2(\Omega)}\\
		&=\sup_{t\in [0,T]} \left (\int_{\Omega} |b(x,v_{n}(t,x))-b(x,v(t,x))|^2\ dx\right)^{\frac{1}{2}}\\
		\eqref{blipschitz}\ \ \ &\leq L_0\sup_{t\in [0,T]}\left (\int_{\Omega} |v_{n}(t,x)-v(t,x)|^2\ dx\right)^{\frac{1}{2}}\\
		&=L_0\sup_{t\in [0,T]} \Vert v_{n}(t,\cdot)-v(t,\cdot)\Vert_{L^2(\Omega)}\\
		&=L_0	\Vert v_{n}-v\Vert_{C([0,T];L^2(\Omega))}\stackrel{n\to\infty}{\longrightarrow}0,
	\end{align*}
	
	\noindent because, from \textbf{Fact I} we already know that $v_{n}\to v$ in $C([0,T];L^2(\Omega))$. We conclude that $b(\cdot,v_{n})\to b(\cdot,v)$ in $C([0,T];L^2(\Omega))\hookrightarrow L^\infty\bigl(0,T;L^2(\Omega)\bigr)\hookrightarrow L^2\bigl(0,T;L^2(\Omega)\bigr)$. Thence $b(\cdot,v_{n})\to b(\cdot,v)$ in $L^2\bigl(0,T;L^2(\Omega)\bigr)$, from where $b(\cdot,v_{n})\weak b(\cdot,v)$ in $L^2\bigl(0,T;L^2(\Omega)\bigr)$.
	
\end{proof}

\bigskip

\noindent\textbf{\underline{Fact XI}: For any $\tau\in (0,T)$ we have that $\dfrac{\partial b(\cdot,v)}{\partial t}\in L^2\bigl(\tau,T;L^2(\Omega)\bigr)$ and moreover $\dfrac{\partial b(\cdot,v_n)}{\partial t}\weak \dfrac{\partial b(\cdot,v)}{\partial t}$ in $L^2\bigl(\tau,T;L^2(\Omega)\bigr)$. Thereby $b(\cdot,v)\in H^1_{\text{loc}}\bigl((0,T);L^2(\Omega)\bigr)$.}

\begin{proof} From Proposition \ref{propofrakB} \textbf{(3)} we get that for each $n\geq 1$, $b(\cdot,v_n)\in H^1\bigl((0,T);L^2(\Omega)\bigr)$ and 
	
	\begin{equation}
		\dfrac{\partial b(\cdot,v_n(t,\cdot))}{\partial t}=\dfrac{\partial b(\cdot,v_n(t,\cdot))}{\partial s}\cdot \dfrac{\partial v_n(t,\cdot)}{\partial t},\ \text{for a.e.}\ t\in (0,T).
	\end{equation}
	
	\noindent Moreover:
	
	\begin{align*}
		\left\Vert\dfrac{\partial b(\cdot,v_n(\cdot,\cdot))}{\partial t}\right\Vert_{L^2(\tau,T;L^2(\Omega))}&=\left(\int_{\tau}^T \left\Vert\dfrac{\partial b(\cdot,v_n(t,\cdot))}{\partial t}\right\Vert_{L^2(\Omega)}^2\ dt\right)^\frac{1}{2}\\
		&=\left(\int_{\tau}^T \left\Vert\dfrac{\partial b(\cdot,v_n(t,\cdot))}{\partial s}\cdot \dfrac{\partial v_n(t,\cdot)}{\partial t}\right\Vert_{L^2(\Omega)}^2\ dt\right)^\frac{1}{2}\\
		&=\left[\int_{\tau}^T \int_{\Omega}\left(\dfrac{\partial b(x,v_n(t,x))}{\partial s}\cdot \dfrac{\partial v_n(t,x)}{\partial t}\right)^2\ dx\ dt\right]^\frac{1}{2}\\
	\text{(Prop. \ref{propoext} \textbf{(3)})}\ \ \	&\leq L_0\left[\int_{\tau}^T \int_{\Omega}\left(\dfrac{\partial v_n(t,x)}{\partial t}\right)^2\ dx\ dt\right]^\frac{1}{2}\\
	&=L_0\left\Vert\dfrac{\partial v_n}{\partial t}\right\Vert_{L^2(\tau,T;L^2(\Omega))},\ \forall\ n\geq 1.
	\end{align*}
	
	\noindent Now, using \textbf{Fact IV}, we deduce that the sequence $\left(\dfrac{\partial b(\cdot,v_n(\cdot,\cdot))}{\partial t}\right)_{n\geq 1}$ is bounded in $L^2\bigl(\tau,T;L^2(\Omega)\bigr)$. We have now all that is needed to make a verbatim repetition of the argument given in the proof of \textbf{Fact V} for $b(\cdot,v)$ instead of $v$, so the conclusion follows.
\end{proof}

\bigskip

\begin{remark}
\noindent We know from \eqref{i5smecher} that for any $n\geq 1$, and any subinterval $(t_1,t_2)\subseteq (0,T)$:

\begin{align}
	&\int_{t_1}^{t_2}\int_\Omega\dfrac{\partial b(x,v_n)}{\partial t}\Theta\ dx\ dt+\int_{t_1}^{t_2}\int_{\Omega}\mathbf{a}(x,\nabla v_n(t,x))\cdot\nabla\Theta(t,x)\ dx\ dt+\lambda\int_{t_1}^{t_2}\int_{\Omega}b(x,v_n)\Theta\ dx\ dt=\nonumber\\
	=&\int_{t_1}^{t_2}\int_{\Omega} g\Theta\ dx\ dt,\ \forall\ \Theta\in L^{2}\bigl(0,T;W^{1,p(x)}(\Omega)\bigr).\label{i5smechern}
\end{align}
\end{remark}

\bigskip

\noindent\textbf{\underline{Fact XII}: Consider any $(t_1,t_2)\subseteq (0,T)$ with $t_1>0$. Then $\lim\limits_{n\to\infty} \displaystyle\int_{t_1}^{t_2}\int_{\Omega}\mathbf{a}(x,\nabla v_n(t,x))\cdot\nabla\Theta(t,x)\ dx\ dt=\int_{t_1}^{t_2}\int_{\Omega} \left( g-\lambda b(x,v)-\dfrac{\partial b(x,v)}{\partial t} \right )\Theta\ dx\ dt$, $\forall\ \Theta\in L^{2}\bigl(0,T;W^{1,p(x)}(\Omega)\bigr)$.}

\begin{proof} We have that:

\noindent $\blacktriangleright$ From \textbf{Fact XI}, applied for $\tau=t_1\in (0,T)$, we know that $\dfrac{\partial b(\cdot,v_{n})}{\partial t}\weak \dfrac{\partial b(\cdot,v)}{\partial t}$ in $L^2\bigl(t_1,T;L^2(\Omega)\bigr)$. This means that for every $\psi\in L^2(t_1,T;L^2(\Omega))$ one has that:

\begin{equation}
	\lim\limits_{k\to\infty}\int_{t_1}^T\int_{\Omega} \dfrac{\partial b(\cdot,v_{n})}{\partial t}\psi(t,x)\ dx\ dt=\int_{0}^T\int_{\Omega} \dfrac{\partial b(\cdot,v)}{\partial t}\psi(t,x)\ dx\ dt.
\end{equation}

\noindent Therefore, taking into account that $\Theta\in L^2\bigl(0,T;W^{1,p(x)}(\Omega)\bigr)\subset L^2\bigl(0,T;L^2(\Omega)\bigr)\subset L^2\bigl(t_1,T;L^2(\Omega)\bigr)$\footnote{This follows from Proposition 2.2.5 given in \cite[page 128]{gasinski2005nonlinear}, using \textbf{(H2)}.} and choosing in the above equality $\psi(t,x)=\chi_{(t_1,t_2)}(t)\Theta(t,x)\in L^2(t_1,T;L^2(\Omega))$, we conclude that:

\begin{equation}\label{i30}
	\lim\limits_{k\to\infty}\int_{t_1}^{t_2}\int_{\Omega} \dfrac{\partial b(\cdot,v_{n})}{\partial t}\Theta(t,x)\ dx\ dt=\int_{t_1}^{t_2}\int_{\Omega} \dfrac{\partial b(\cdot,v)}{\partial t}\Theta(t,x)\ dx\ dt.
\end{equation}

\bigskip

\noindent $\blacktriangleright$ From \textbf{Fact X} we have that $b(\cdot,v_{n})\weak b(\cdot,v)$ in $L^2(0,T;L^2(\Omega))$. This means that for every $\psi\in L^2(0,T;L^2(\Omega))$ one has that:

\begin{equation}
	\lim\limits_{k\to\infty}\int_{0}^T\int_{\Omega} b(x,v_{n}(t,x))\psi(t,x)\ dx\ dt=\int_{0}^T\int_{\Omega} b(x,v(t,x))\psi(t,x)\ dx\ dt.
\end{equation}

\noindent Therefore, taking into account that $\Theta\in L^2\bigl(0,T;W^{1,p(x)}(\Omega)\bigr)\subset L^2\bigl(0,T;L^2(\Omega)\bigr)$ and choosing in the above equality $\psi(t,x)=\chi_{(t_1,t_2)}(t)\Theta(t,x)\in L^2(0,T;L^2(\Omega))$, we conclude that:

\begin{equation}\label{i20}
	\lim\limits_{k\to\infty} \int_{t_1}^{t_2}\int_{\Omega}  b(x,v_{n}(t,x))\Theta(t,x)\ dx\ dt=\int_{t_1}^{t_2}\int_{\Omega}  b(x,v(t,x))\Theta(x)\ dx\ dt.
\end{equation}

\noindent Now we can make $n\to \infty$ in \eqref{i5smechern} and obtain that:

\begin{align}
	&\lim\limits_{n\to\infty} \int_{t_1}^{t_2}\int_{\Omega}\mathbf{a}(x,\nabla v_n(t,x))\cdot\nabla\Theta(t,x)\ dx\ dt=\nonumber\\
	=&\int_{t_1}^{t_2}\int_{\Omega} \left( g-\lambda b(x,v)-\dfrac{\partial b(x,v)}{\partial t} \right )\Theta\ dx\ dt,\ \forall\ \Theta\in L^{2}\bigl(0,T;W^{1,p(x)}(\Omega)\bigr).\label{i5smechernfinal}
\end{align}

\end{proof}

\begin{definition}\label{defhashh} Consider the function $h:(0,T)\times\Omega$,
	
	\begin{equation}
		h(t,x)=g(t,x)-\lambda b(x,v(t,x))-\dfrac{\partial b(x,v(t,x))}{\partial t},\ \text{for a.e.}\ (t,x)\in (0,T)\times\Omega.
	\end{equation}
	
\end{definition}

\begin{remark}\label{remhashh} Since $g\in L^{\infty}\bigl((0,T)\times\Omega\bigr)$, $\lambda b(\cdot,v)\in L^{\infty}\bigl((0,T)\times\Omega\bigr)$ and from \textnormal{\textbf{Fact XI}}: $\dfrac{\partial b(\cdot,v)}{\partial t}\in L^2\bigl(\tau,T;L^2(\Omega)\bigr),\ \forall\ \tau\in (0,T)$ we deduce that:
	
	\begin{equation}\label{has1}
		h\in L^2\bigl(\tau,T;L^2(\Omega)\bigr),\ \forall\ \tau\in (0,T).
	\end{equation}
	
\noindent Equation \eqref{i5smechernfinal} can be written for any $(t_1,t_2)\subseteq (0,T)$ with $t_1>0$ as:

\begin{equation}\label{i5smechernfinalh}
	\lim\limits_{n\to\infty} \int_{t_1}^{t_2}\int_{\Omega}\mathbf{a}(x,\nabla v_n(t,x))\cdot\nabla\Theta(t,x)\ dx\ dt=\int_{t_1}^{t_2}\int_{\Omega} h(t,x)\Theta(t,x)\ dx\ dt,\ \forall\ \Theta\in L^{2}\bigl(0,T;W^{1,p(x)}(\Omega)\bigr).
\end{equation}
	
\end{remark}

\bigskip
\noindent\textbf{\underline{Fact XIII}: The following relation holds for any $(t_1,t_2)\subseteq (0,T)$ with $t_1>0$:}

\begin{align}
	&\lim\limits_{n\to\infty} \int_{t_1}^{t_2}\int_{\Omega}\mathbf{a}(x,\nabla v_n(t,x))\cdot\nabla v_n(t,x)\ dx\ dt=\nonumber\\
	=&\int_{t_1}^{t_2}\int_{\Omega} \left( g-\lambda b(x,v)-\dfrac{\partial b(x,v)}{\partial t} \right )v\ dx\ dt=\int_{t_1}^{t_2}\int_{\Omega} h(t,x)v(t,x)\ dx\ dt.\label{i5smechernfinalvn}
\end{align}

\begin{proof} Setting $\Theta:=v_n\in L^{\infty}\bigl(0,T;W^{1,p(x)}(\Omega)\bigr)\subset L^{2}\bigl(0,T;W^{1,p(x)}(\Omega)\bigr)$ in \eqref{i5smechern} we get that:
	
\begin{equation}\label{mama1}
	\int_{t_1}^{t_2}\int_{\Omega}\mathbf{a}(x,\nabla v_n(t,x))\cdot\nabla v_n(t,x)\ dx\ dt=\int_{t_1}^{t_2}\int_{\Omega} \left( g-\lambda b(x,v_n)-\dfrac{\partial b(x,v_n)}{\partial t} \right )v_n\ dx\ dt,\ \forall\ n\geq 1.
\end{equation}

\noindent $\blacktriangleright$ Knowing that $g\in L^{\infty}\bigl((0,T)\times\Omega)\bigr)\subset L^2\bigl((0,T)\times\Omega\bigr)$ and $v_n\to v$ in $L^2\bigl(0,T;L^2(\Omega)\bigr)\simeq L^2\bigl((0,T)\times\Omega\bigr)$ (see Remark \ref{vnvlr0tl2omega}) we deduce that:

\begin{align}\label{mama4}
\left|\int_{t_1}^{t_2}\int_{\Omega} gv_n\ dx\ dt-\int_{t_1}^{t_2}\int_{\Omega}gv\ dx\ dt \right|&=\left|\int_{t_1}^{t_2}\int_{\Omega} g(v_n-v)\ dx\ dt \right|\nonumber\\
&\leq \int_{t_1}^{t_2}\int_{\Omega}|g|\cdot |v_n-v|\ dx\ dt\nonumber\\
\text{(Cauchy ineq.)}\ \ \ &\leq \Vert g\Vert_{L^2((0,T)\times\Omega)}\cdot \Vert v_n-v\Vert_{L^2((0,T)\times\Omega)}\stackrel{n\to\infty}{\longrightarrow} 0.
\end{align}

\bigskip

\noindent $\blacktriangleright$ Knowing from \textbf{Fact X} that $b(\cdot,v_n)\to b(\cdot,v)$ in $ L^2\bigl((0,T)\times\Omega\bigr)$ and $v_n\to v$ in $L^2\bigl((0,T)\times\Omega\bigr)$ we obtain from \textit{Cauchy inequality} that:

\begin{align}\label{mama3}
	&\left|\int_{t_1}^{t_2}\int_{\Omega} b(x,v_n)v_n\ dx\ dt-\int_{t_1}^{t_2}\int_{\Omega}b(x,v)v\ dx\ dt \right|\nonumber\\
	=\ & \left|\int_{t_1}^{t_2}\int_{\Omega} \bigl [b(x,v_n)-b(x,v)\bigr]\cdot v\ dx\ dt +\int_{t_1}^{t_2}\int_{\Omega} b(x,v_n)(v_n-v)\ dx\ dt \right|\nonumber\\
	=\ &\int_{t_1}^{t_2}\int_{\Omega} \bigl |b(x,v_n)-b(x,v)\bigr |\cdot |v|\ dx\ dt+\int_{t_1}^{t_2}\int_{\Omega} |b(x,v_n)|\cdot|v_n-v|\ dx\ dt\nonumber\\
\eqref{bequ1}	\leq\ &\int_{0}^{T}\int_{\Omega} \bigl |b(x,v_n)-b(x,v)\bigr |\cdot |v|\ dx\ dt+\int_{0}^{T}\int_{\Omega} \max\bigl\{\Vert b(\cdot,\varepsilon)\Vert_{L^{\infty}(\Omega)},\Vert b(\cdot,\delta)\Vert_{L^{\infty}(\Omega)}\bigr\}\cdot|v_n-v|\ dx\ dt\nonumber\\
\leq\ & \underbrace{\Vert b(\cdot,v_n)-b(\cdot,v)\Vert_{L^2((0,T)\times\Omega)}}_{\longrightarrow 0}\cdot \Vert v\Vert_{L^2((0,T)\times\Omega)}+\nonumber\\
 &+\max\bigl\{\Vert b(\cdot,\varepsilon)\Vert_{L^{\infty}(\Omega)},\Vert b(\cdot,\delta)\Vert_{L^{\infty}(\Omega)}\bigr\}\sqrt{T|\Omega|}\cdot\underbrace{\Vert v_n-v\Vert_{L^2((0,T)\times\Omega)}}_{\longrightarrow 0}\stackrel{n\to\infty}{\longrightarrow} 0.
\end{align}

\bigskip

\noindent $\blacktriangleright$ Knowing from \textbf{Fact XI} that $\dfrac{\partial b(\cdot,v_n)}{\partial t}\weak \dfrac{\partial b(\cdot,v)}{\partial t}$ in $ L^2\bigl((t_1,T)\times\Omega\bigr)$ and $v_n\to v$ in $L^2\bigl((0,T)\times\Omega\bigr)$ we obtain from \textit{Cauchy inequality} that:

\begin{align}\label{mama2}
	&\left|\int_{t_1}^{t_2}\int_{\Omega} \dfrac{\partial b(x,v_n)}{\partial t}v_n\ dx\ dt-\int_{t_1}^{t_2}\int_{\Omega}\dfrac{\partial b(x,v)}{\partial t}v\ dx\ dt \right|\nonumber\\
	=\ & \left|\int_{t_1}^{t_2}\int_{\Omega} \left [\dfrac{\partial b(x,v_n)}{\partial t}-\dfrac{\partial b(x,v)}{\partial t}\right ]\cdot v\ dx\ dt +\int_{t_1}^{t_2}\int_{\Omega} \dfrac{\partial b(x,v_n)}{\partial t}(v_n-v)\ dx\ dt \right|\nonumber\\
	\leq\ & \left|\int_{t_1}^{t_2}\int_{\Omega} \dfrac{\partial b(x,v_n)}{\partial t}v\ dx\ dt-\int_{t_1}^{t_2}\int_{\Omega} \dfrac{\partial b(x,v)}{\partial t} v\ dx\ dt\right |+\nonumber\\
	&+\left |\int_{t_1}^{t_2}\int_{\Omega} \dfrac{\partial b(x,v_n)}{\partial t}(v_n-v)\ dx\ dt \right|\nonumber\\
	\leq\ & \left|\int_{t_1}^{t_2}\int_{\Omega} \dfrac{\partial b(x,v_n)}{\partial t}v\ dx\ dt-\int_{t_1}^{t_2}\int_{\Omega} \dfrac{\partial b(x,v)}{\partial t} v\ dx\ dt\right |+\nonumber\\
	&+\int_{t_1}^{t_2}\int_{\Omega} \left |\dfrac{\partial b(x,v_n)}{\partial t}\right |\cdot |v_n-v|\ dx\ dt \nonumber\\
	\leq \ & \underbrace{\left|\int_{t_1}^{t_2}\int_{\Omega} \dfrac{\partial b(x,v_n)}{\partial t}v\ dx\ dt-\int_{t_1}^{t_2}\int_{\Omega} \dfrac{\partial b(x,v)}{\partial t} v\ dx\ dt\right |}_{\longrightarrow 0}+\nonumber\\
	&+\underbrace{\left\Vert\dfrac{\partial b(\cdot,v_n)}{\partial t}\right\Vert_{L^2((t_1,t_2)\times\Omega)}}_{\text{bounded}}\cdot \underbrace{\Vert v_n-v\Vert_{L^2((t_1,t_2)\times\Omega)}}_{\longrightarrow 0}\stackrel{n\to\infty}{\longrightarrow} 0.
\end{align}

\noindent Combining now relations \eqref{mama1}, \eqref{mama4}, \eqref{mama3} and \eqref{mama2} we get that:

\begin{align*}
	\lim\limits_{n\to\infty} \int_{t_1}^{t_2}\int_{\Omega}\mathbf{a}(x,\nabla v_n(t,x))\cdot\nabla v_n(t,x)\ dx\ dt&=\lim\limits_{n\to\infty}\int_{t_1}^{t_2}\int_{\Omega} \left( g-\lambda b(x,v_n)-\dfrac{\partial b(x,v_n)}{\partial t} \right )v_n\ dx\ dt\\
	&=\int_{t_1}^{t_2}\int_{\Omega} \left( g-\lambda b(x,v)-\dfrac{\partial b(x,v)}{\partial t} \right )v\ dx\ dt.
\end{align*}
	
\end{proof}

\bigskip
\noindent\textbf{\underline{Fact XIV}: For any $(t_1,t_2)\subseteq (0,T)$ with $t_1\in (0,T)$ we have that: $\mathcal{A}(v_n)\longrightarrow\mathcal{A}(v)$ in $L^1(t_1,t_2)$, i.e. $\lim\limits_{n\to\infty} \displaystyle\int_{t_1}^{t_2} \bigl |\mathcal{A}(v_n(t,\cdot))-\mathcal{A}(v(t,\cdot))\bigr|\ dt=0$. In particular $\lim\limits_{n\to\infty} \displaystyle\int_{t_1}^{t_2}\mathcal{A}(v_n(t,\cdot))\ dt=\displaystyle\int_{t_1}^{t_2}\mathcal{A}(v(t,\cdot))\ dt$.}
	
\begin{proof} In \textbf{Fact III} we proved that $(0,T)\ni t\mapsto\displaystyle\alpha_n(t)=\mathcal{A}(v_n(t,\cdot))$ is from $L^1(0,T)\subset L^1(t_1,t_2)$ for each $n\geq 1$. Also, in \textbf{Fact IX}, we showed that $(0,T)\ni t\mapsto \mathcal{A}(v(t,\cdot))$ is from $L^1(0,T)\subset L^1(t_1,t_2)$ -- see \eqref{lumanare4}. Moreover, \eqref{lumanare3} and \eqref{lumanare4} gives us that
	
	\begin{equation}
		\int_{0}^T\mathcal{A}\bigl(v(t,\cdot)\bigr)\ dt\leq \int_0^T\alpha(t)\ dt=\int_{0}^T\liminf\limits_{n\to\infty}\alpha_n(t)\ dt\leq \liminf\limits_{n\to\infty}\int_{0}^T \mathcal{A}\bigl(v_n(t,\cdot)\bigr)\ dt\leq C_{\mathcal{A}}.
	\end{equation}
	
	\noindent From \eqref{lumanare2} and \textit{Fatou's lemma} we obtain that:
	
	\begin{equation}\label{magazin1}
		\int_{t_1}^{t_2}\mathcal{A}\bigl(v(t,\cdot)\bigr)\ dt\leq \int_{t_1}^{t_2}\alpha(t)\ dt=\int_{t_1}^{t_2}\liminf\limits_{n\to\infty}\alpha_n(t)\ dt\leq \liminf\limits_{n\to\infty}\int_{t_1}^{t_2} \mathcal{A}\bigl(v_n(t,\cdot)\bigr)\ dt\leq C_{\mathcal{A}}.
	\end{equation}
	
	\noindent Using now Proposition \ref{propomathcalA1} \textbf{(4)} it follows for a.e. $t\in (0,T)$ that:
	
	\begin{equation}\label{magazin2}
		\mathcal{A}\bigl(v_n(t,\cdot)\bigr)\leq \mathcal{A}\bigl(v(t,\cdot)\bigr)+\int_{\Omega} \mathbf{a}(x,\nabla v_n(t,x))\cdot\bigl(\nabla v_n(t,x)-\nabla v(t,x)\bigr)\ dx.
	\end{equation}
	
	\noindent The problem with \eqref{magazin2} is that we cannot integrate it on $(t_1,t_2)$, because $(0,T)\ni t\mapsto \displaystyle\int_{\Omega} \mathbf{a}(x,\nabla v_n(t,x))\cdot\bigl(\nabla v_n(t,x)-\nabla v(t,x)\bigr)\ dx$, may fail to be in $L^1(0,T)$, because we only know that $v\in L^1\bigl(0,T;W^{1,p(x)}(\Omega)\bigr)$.
	
	\noindent However we can make a trick here. Proposition \ref{propomathcalA1} \textbf{(4)} allows us to write for any $\Theta\in L^{2}\bigl(0,T;W^{1,p(x)}(\Omega)\bigr)$ that:
	
	\begin{equation}\label{magazin3}
		\mathcal{A}\bigl(v_n(t,\cdot)\bigr)\leq\mathcal{A}\bigl(\Theta(t,\cdot)\bigr)+\int_{\Omega} \mathbf{a}(x,\nabla v_n(t,x))\cdot\bigl(\nabla v_n(t,x)-\nabla \Theta(t,x)\bigr)\ dx.
	\end{equation}
	
	\noindent Now $(0,T)\ni t\mapsto \displaystyle\int_{\Omega} \mathbf{a}(x,\nabla v_n(t,x))\cdot\bigl(\nabla v_n(t,x)-\nabla\Theta(t,x)\bigr)\ dx$ is from $L^1(0,T)\subset L^1(t_1,t_2)$, by Proposition \ref{propoa2}. This allows us to integrate \eqref{magazin3} on $(t_1,t_2)$ and get that:
	
	\begin{equation}\label{magazin4}
		\int_{t_1}^{t_2}\mathcal{A}\bigl(v_n(t,\cdot)\bigr)\ dt\leq\int_{t_1}^{t_2}\mathcal{A}\bigl(\Theta(t,\cdot)\bigr)\ dt +\int_{t_1}^{t_2}\int_{\Omega} \mathbf{a}(x,\nabla v_n(t,x))\cdot\bigl(\nabla v_n(t,x)-\nabla \Theta(t,x)\bigr)\ dx\ dt.
	\end{equation}
	
	\noindent The key step is here: using \textbf{Fact XIII} and Remark \ref{remhashh}\footnote{Here we use the fact that $t_1>0$.} -- especially \eqref{i5smechernfinalh} --  we obtain by passing \eqref{magazin4} to limit superior that:
	
	\begin{align}\label{magazin5}
		\limsup\limits_{n\to\infty} \int_{t_1}^{t_2}\mathcal{A}\bigl(v_n(t,\cdot)\bigr)\ dt\leq&\int_{t_1}^{t_2}\mathcal{A}\bigl(\Theta(t,\cdot)\bigr)\ dt+\nonumber\\
		&+\limsup\limits_{n\to\infty}\int_{t_1}^{t_2}\int_{\Omega} \mathbf{a}(x,\nabla v_n(t,x))\cdot\bigl(\nabla v_n(t,x)-\nabla \Theta(t,x)\bigr)\ dx\ dt\nonumber\\
		=&\int_{t_1}^{t_2}\mathcal{A}\bigl(\Theta(t,\cdot)\bigr)\ dt+\nonumber\\
		&+\lim\limits_{n\to\infty}\int_{t_1}^{t_2}\int_{\Omega} \mathbf{a}(x,\nabla v_n(t,x))\cdot\bigl(\nabla v_n(t,x)-\nabla \Theta(t,x)\bigr)\ dx\ dt\nonumber \\
		=&\int_{t_1}^{t_2}\mathcal{A}\bigl(\Theta(t,\cdot)\bigr)\ dt+\int_{t_1}^{t_2}\int_{\Omega} h(t,x)\cdot\bigl (v-\Theta)\ dx\ dt.
	\end{align}
	
	\noindent For each $m\geq 1$ consider the following set:
	
	\begin{equation}\label{magazin6}
		I_m:=\bigl\{t\in (t_1,t_2)\ |\ \Vert v(t,\cdot)\Vert_{W^{1,p(x)}(\Omega)}\leq m\bigr\}.
	\end{equation}
	
	\noindent Since $v:(0,T)\to W^{1,p(x)}(\Omega)$ is strongly measurable -- see \textbf{Fact VII}, it follows that $I_m$ is a measurable subset of $(t_1,t_2)$. This is because $(0,T)\ni t\mapsto\Vert v(t,\cdot)\Vert_{W^{1,p(x)}(\Omega)}$ is the composition between the norm of $W^{1,p(x)}(\Omega)$ -- which is a continuous function, even 1-Lipschitz -- and the strongly measurable function $v:(0,T)\to W^{1,p(x)}(\Omega)$.
	
	\noindent From the definition of the measurable set $I_m$ we deduce that $\Theta_m:=\chi_{I_m}(\cdot)v(\cdot,\cdot)\in L^{\infty}\bigl(0,T;W^{1,p(x)}(\Omega)\bigr)\subset L^2\bigl(0,T;W^{1,p(x)}(\Omega)\bigr)$. Also, from Remark \ref{vnvlr0tl2omega}, we get that $\Theta_m=\chi_{I_m}v\in L^2\bigl(0,T;L^2(\Omega)\bigr)$.
	
	\medskip
	
	\begin{align} 
	\blacktriangleright\lim\limits_{m\to\infty}\displaystyle\int_{t_1}^{t_2}\int_{\Omega} h(t,x)\cdot\bigl (v-\Theta_m)\ dx\ dt=&\int_{t_1}^{t_2} \left (\int_{\Omega} h(t,x)v(t,x)\ dx \right )\ dt\nonumber\\
	&-\lim\limits_{m\to\infty}\int_{t_1}^{t_2} \left (\int_{\Omega} h(t,x)\Theta_m(t,x)\ dx \right )\ dt\\
	=&\int_{t_1}^{t_2} \left (\int_{\Omega} h(t,x)v(t,x)\ dx \right )\ dt\nonumber\\
	&-\lim\limits_{m\to\infty}\int_{I_m} \left (\int_{\Omega} h(t,x)v(t,x)\ dx \right )\ dt\\
\text{(Lemma \ref{beppolevi1})}\ \ \	=&\int_{t_1}^{t_2} \left (\int_{\Omega} h(t,x)v(t,x)\ dx \right )\ dt\nonumber\\
	&-\int_{t_1}^{t_2} \left (\int_{\Omega} h(t,x)v(t,x)\ dx \right )\ dt=0.\label{magazin8}
	\end{align}
	
	\medskip
	
	\begin{align}\blacktriangleright\ \lim\limits_{m\to\infty}\displaystyle\int_{t_1}^{t_2} \mathcal{A}\bigl(\Theta_m(t,\cdot)\bigr)\ dt&=\lim\limits_{m\to\infty}\displaystyle\int_{I_m} \mathcal{A}\bigl(\Theta_m(t,\cdot)\bigr)+\displaystyle\int_{(t_1,t_2)\setminus I_m} \mathcal{A}\bigl(\Theta_m(t,\cdot)\bigr)\ dt\nonumber\\
	&=\lim\limits_{m\to\infty}\displaystyle\int_{I_m} \mathcal{A}\bigl(v(t,\cdot)\bigr)\ dt+\displaystyle\int_{(t_1,t_2)\setminus I_m} \underbrace{\mathcal{A}\bigl(0\bigr)}_{=0}\ dt\nonumber\\
	&=\lim\limits_{m\to\infty}\displaystyle\int_{I_m} \mathcal{A}\bigl(v(t,\cdot)\bigr)\ dt\nonumber\\
\text{(Lemma \ref{beppolevi1})}\ \ \	&=\int_{t_1}^{t_2} \mathcal{A}\bigl(v(t,\cdot)\bigr)\ dt.\label{magazin7}
	\end{align}
	
	\noindent Setting $\Theta=\Theta_m\in L^{\infty}\bigl(0,T;W^{1,p(x)}(\Omega)\bigr),\ m\geq 1$ in \eqref{magazin5} leads us to:
	
		\begin{equation}\label{magazin9}
		\limsup\limits_{n\to\infty} \int_{t_1}^{t_2}\mathcal{A}\bigl(v_n(t,\cdot)\bigr)\ dt\leq \int_{t_1}^{t_2}\mathcal{A}\bigl(\Theta_m(t,\cdot)\bigr)\ dt+\int_{t_1}^{t_2}\int_{\Omega} h(t,x)\cdot\bigl (v-\Theta_m)\ dx\ dt.
	\end{equation}
	
	\noindent Making $m\to\infty$ in \eqref{magazin9} and using \eqref{magazin8} and \eqref{magazin7} allows us to write:
	
	\begin{equation}\label{magazin10}
			\limsup\limits_{n\to\infty} \int_{t_1}^{t_2}\mathcal{A}\bigl(v_n(t,\cdot)\bigr)\ dt\leq \int_{t_1}^{t_2} \mathcal{A}\bigl(v(t,\cdot)\bigr)\ dt.
	\end{equation}
	
	\noindent We combine \eqref{magazin1} and \eqref{magazin10} and get that:
	
	\begin{equation}\label{magazin11}
		\int_{t_1}^{t_2}\mathcal{A}\bigl(v(t,\cdot)\bigr)\ dt\leq \liminf\limits_{n\to\infty}\int_{t_1}^{t_2} \mathcal{A}\bigl(v_n(t,\cdot)\bigr)\ dt\leq \limsup\limits_{n\to\infty} \int_{t_1}^{t_2}\mathcal{A}\bigl(v_n(t,\cdot)\bigr)\ dt\leq \int_{t_1}^{t_2} \mathcal{A}\bigl(v(t,\cdot)\bigr)\ dt.
	\end{equation} 
	
	\noindent This shows that:
	
	\begin{equation}\label{magazin12}
	\exists\  \lim\limits_{n\to\infty}\int_{t_1}^{t_2} \mathcal{A}\bigl(v_n(t,\cdot)\bigr)\ dt=\int_{t_1}^{t_2}\mathcal{A}\bigl(v(t,\cdot)\bigr)\ dt.
	\end{equation}
	
	\noindent Now, since $\liminf\limits_{n\to\infty}\mathcal{A}\bigl(v_n(t,\cdot)\bigr)\geq \mathcal{A}\bigl(v(t,\cdot)\bigr)$ for a.e. $t\in (0,T)$ -- see \eqref{lumanare2} --, we deduce from Lemma \ref{scheffeslemma}, taking into account \eqref{magazin12} that:
	
	\begin{equation}
		\lim\limits_{n\to\infty}\int_{t_1}^{t_2} \bigl|\mathcal{A}\bigl(v_n(t,\cdot)\bigr)-\mathcal{A}\bigl(v(t,\cdot)\bigr)\bigr|\ dt=0.
	\end{equation}

\end{proof}

\bigskip
\noindent\textbf{\underline{Fact XV}: There is a null-measure set $N_0\subseteq (0,T)$ with the property that:}

\begin{equation}
	\int_{\Omega} \mathbf{a}(x,\nabla v(t,x))\cdot\nabla\phi(x)\ dx=\int_{\Omega} h(t,x)\phi(x)\ dx,\ \forall\ \phi\in W^{1,p(x)}(\Omega),\ \forall\ t\in (0,T)\setminus N_0.
\end{equation}
	
\begin{proof} Proposition \ref{propomathcalA1} \textbf{(4)} allows us to write for any $\Theta\in L^{2}\bigl(0,T;W^{1,p(x)}(\Omega)\bigr)$ that:
	
	\begin{equation}\label{KAUF1}
	\mathcal{A}\bigl(\Theta(t,\cdot)\bigr)\geq 	\mathcal{A}\bigl(v_n(t,\cdot)\bigr)+\int_{\Omega} \mathbf{a}(x,\nabla v_n(t,x))\cdot\bigl(\nabla \Theta(t,x)-\nabla v_n(t,x)\bigr)\ dx.
	\end{equation}
	
	\noindent Now $(0,T)\ni t\mapsto \displaystyle\int_{\Omega} \mathbf{a}(x,\nabla v_n(t,x))\cdot\bigl(\nabla\Theta(t,x)-\nabla v_n(t,x)\bigr)\ dx$ is from $L^1(0,T)\subset L^1(t_1,t_2)$, by Proposition \ref{propoa2}. This allows us to integrate \eqref{KAUF1} on $(\tau,T)$, for any $\tau>0$, and get that:
	
	\begin{equation}\label{KAUF2}
		\int_{\tau}^{T}\mathcal{A}\bigl(\Theta(t,\cdot)\bigr)\ dt\geq\int_{\tau}^{T}\mathcal{A}\bigl(v_n(t,\cdot)\bigr)\ dt +\int_{\tau}^{T}\int_{\Omega} \mathbf{a}(x,\nabla v_n(t,x))\cdot\bigl(\nabla \Theta(t,x)-\nabla v_n(t,x)\bigr)\ dx\ dt.
	\end{equation}
	
	\noindent Making $n\to\infty$ in \eqref{KAUF2} and using \eqref{magazin12}, \eqref{i5smechernfinalh} and \textbf{Fact XIII} gives us for any $\Theta\in L^2\bigl(0,T;W^{1,p(x)}(\Omega)\bigr)$ that:
	
	\begin{equation}\label{KAUF3}
		\int_{\tau}^{T}\mathcal{A}\bigl(\Theta(t,\cdot)\bigr)\ dt\geq\int_{\tau}^{T}\mathcal{A}\bigl(v(t,\cdot)\bigr)\ dt +\int_{\tau}^{T}\int_{\Omega} h(t,x)\cdot\bigl(\Theta(t,x)-v(t,x)\bigr)\ dx\ dt.
	\end{equation}
	
	\noindent Fix any $z\in W^{1,p(x)}(\Omega)$ and any measurable set $I\subseteq (\tau,T)$. For each $m\geq 1$, consider as in \textbf{Fact XIV} the measurable sets (see \eqref{magazin6}):

	\begin{equation}
		I_m=\bigl\{t\in (\tau,T)\ |\ \Vert v(t,\cdot)\Vert_{W^{1,p(x)}(\Omega)}\leq m\}\subset (\tau,T),
	\end{equation}
	
	\noindent and define:
	
	\begin{equation}\label{KAUF4}
	\Theta_m:(0,T)\to W^{1,p(x)}(\Omega),\ \Theta_m(t,\cdot)=\begin{cases}z,& t\in I\\[3mm] v(t,\cdot), & \ t\in \bigl[(\tau,T)\setminus I\bigr]\cap I_m\\[3mm] 0, & t\in \bigl[(0,T)\setminus I\bigr]\setminus I_m  \end{cases}.
\end{equation}

\noindent Clearly $\Theta_m$ is a strongly measurable function, because it is the sum of two simple functions and a strongly measurable function. Also, it is easy to remark that:

\begin{align*}
	\Vert\Theta_m\Vert_{L^{\infty}(0,T;W^{1,p(x)}(\Omega))}&=\underset{t\in (0,T)}{\operatorname{ess\ sup}}\ \Vert\Theta_m(t,\cdot)\Vert_{W^{1,p(x)}(\Omega)}\leq \max\bigl\{\Vert z\Vert_{W^{1,p(x)}(\Omega)},	\Vert v\Vert_{L^{\infty}(I_m;W^{1,p(x)}(\Omega))}\}\\
	&\leq \max\bigl\{\Vert z\Vert_{W^{1,p(x)}(\Omega)},m\},\ \forall\ m\geq 1.
\end{align*}

\noindent Thus $\Theta_m\in L^{\infty}\bigl(0,T;W^{1,p(x)}(\Omega)\bigr)\subset L^2\bigl(0,T;W^{1,p(x)}(\Omega)\bigr)$. Choosing $\Theta=\Theta_m$ for each natural number $m\geq 1$, \eqref{KAUF3} becomes:

	\begin{equation}\label{KAUF5}
	\int_{\tau}^{T}\left [\mathcal{A}\bigl(\Theta_m(t,\cdot)\bigr)-\mathcal{A}\bigl(v(t,\cdot)\bigr)-\int_{\Omega} h(t,x)\cdot\bigl(\Theta_m(t,x)-v(t,x)\bigr)\ dx\right ]\ dt\geq 0.
\end{equation}

\noindent Splitting this integral into three parts gives us that:

\begin{align}\label{KAUF6}
	0&\leq \int_{\tau}^{T}\left [\mathcal{A}\bigl(\Theta_m(t,\cdot)\bigr)-\mathcal{A}\bigl(v(t,\cdot)\bigr)-\int_{\Omega} h(t,x)\cdot\bigl(\Theta_m(t,x)-v(t,x)\bigr)\ dx\right ]\ dt\nonumber\\
	&=\int_{I}\left [\mathcal{A}\bigl(z\bigr)-\mathcal{A}\bigl(v(t,\cdot)\bigr)-\int_{\Omega} h(t,x)\cdot\bigl(z(x)-v(t,x)\bigr)\ dx\right ] \ dt+\nonumber\\
	&\ \ \ +\int_{[(\tau,T)\setminus I]\cap I_m}\underbrace{\left [\mathcal{A}\bigl(v(t,\cdot)\bigr)-\mathcal{A}\bigl(v(t,\cdot)\bigr)-\int_{\Omega} h(t,x)\cdot\bigl(v(t,x)-v(t,x)\bigr)\ dx\right ]}_{=0} \ dt+\nonumber\\
	&\ \ \ +\int_{[(\tau,T)\setminus I]\setminus I_m}\left [\underbrace{\mathcal{A}\bigl(0\bigr)}_{=0}-\mathcal{A}\bigl(v(t,\cdot)\bigr)-\int_{\Omega} h(t,x)\cdot\bigl(0-v(t,x)\bigr)\ dx\right ] \ dt\nonumber\\
	&=\int_{I}\left [\mathcal{A}\bigl(z\bigr)-\mathcal{A}\bigl(v(t,\cdot)\bigr)-\int_{\Omega} h(t,x)\cdot\bigl(z(x)-v(t,x)\bigr)\ dx\right ] \ dt+\nonumber\\
	&\ \ \ +\int_{[(\tau,T)\setminus I]\setminus I_m}\left [-\mathcal{A}\bigl(v(t,\cdot)\bigr)+\int_{\Omega} h(t,x)\cdot v(t,x)\ dx\right ] \ dt.
\end{align}

\noindent From \textbf{Fact IX} we know that $v\in L^1\bigl(0,T;W^{1,p(x)}(\Omega)\bigr)$. Therefore:

\begin{equation}
	m\chi_{(\tau,T)\setminus I_m}(t)\leq \Vert v(t,\cdot)\Vert_{W^{1,p(x)}(\Omega)},\ \text{for a.e.}\ t\in (\tau,T),
\end{equation}

\noindent and then integrating on $(\tau,T)$ leads us to

\begin{equation}
 m|(\tau,T)\setminus I_m|\leq \int_{\tau}^T \Vert v(t,\cdot)\Vert_{W^{1,p(x)}(\Omega)}\ dt\leq \int_{0}^T \Vert v(t,\cdot)\Vert_{W^{1,p(x)}(\Omega)}\ dt=\Vert v\Vert_{L^{1}(0,T;W^{1,p(x)}(\Omega))}.
\end{equation}

\noindent So it follows that:

\begin{equation}\label{KAUF7}
	|[(\tau,T)\setminus I]\setminus I_m|\leq |(\tau,T)\setminus I_m|\leq\dfrac{\Vert v\Vert_{L^{1}(0,T;W^{1,p(x)}(\Omega))}}{m}\stackrel{m\to\infty}{\longrightarrow}0.
\end{equation}

\noindent Using Remark \ref{remhashh} we have that $h\in L^2\bigl(\tau,T;L^2(\Omega)\bigr)\simeq L^2\bigl((\tau,T)\times\Omega\bigr)$, and since $v\in L^2\bigl(0,T;L^2(\Omega)\bigr)\subset L^2\bigl(\tau,T;L^2(\Omega)\bigr)\simeq L^2\bigl((\tau,T)\times\Omega\bigr)$, we get from \textit{Cauchy inequality} that $hv\in L^1\bigl((\tau,T)\times\Omega\bigr)$. Now from \textit{Fubini's theorem} we obtain that the function $(\tau,T)\ni t\mapsto \displaystyle\int_{\Omega} h(t,x)\cdot v(t,x)\ dx$ is a function from $L^1(\tau,T)$. Taking into account that $(\tau,T)\ni t\mapsto -\mathcal{A}\bigl(v(t,\cdot)\bigr)$ is also from $L^1(\tau,T)$ we deduce that $(\tau,T)\ni t\mapsto -\mathcal{A}\bigl(v(t,\cdot)\bigr)+\displaystyle\int_{\Omega} h(t,x)\cdot v(t,x)\ dx$ is a function from $L^1(\tau,T)$.

\noindent Now, using this information and \eqref{KAUF7}, we deduce by applying Lemma \ref{beppolevi1} that:\footnote{This is because $[(\tau,T)\setminus I]\setminus I_m=(\tau,T)\setminus (I\cup I_m)$, and $I\cup I_m\subset (\tau,T)$ for each $m\geq 1$.}

\begin{equation}\label{KAUF8}
\lim\limits_{m\to\infty} \int_{[(\tau,T)\setminus I]\setminus I_m}\left [-\mathcal{A}\bigl(v(t,\cdot)\bigr)+\int_{\Omega} h(t,x)\cdot v(t,x)\ dx\right ] \ dt=0.
\end{equation}

\noindent By making $m\to\infty$ in \eqref{KAUF6} we get that for any measurable set $I\subseteq (\tau,T)$ and every $z\in W^{1,p(x)}(\Omega)$ the following inequality holds:

\begin{equation}\label{KAUF9}
	\int_{I}\left [\mathcal{A}\bigl(z\bigr)-\mathcal{A}\bigl(v(t,\cdot)\bigr)-\int_{\Omega} h(t,x)\cdot\bigl(z(x)-v(t,x)\bigr)\ dx\right ] \ dt\geq 0.
\end{equation}

\noindent Applying now Lemma \ref{intt1t2} we obtain that for every $z\in W^{1,p(x)}(\Omega)$ there exists a null-measure set $N_{z}\subset (\tau,T)$ such that:

\begin{equation}\label{KAUF10}
	\mathcal{A}\bigl(z\bigr)-\mathcal{A}\bigl(v(t,\cdot)\bigr)-\int_{\Omega} h(t,x)\cdot\bigl(z(x)-v(t,x)\bigr)\ dx \geq 0,\ \forall\ t\in (\tau,T)\setminus N_{z}.
\end{equation}

\noindent The problem is that this null-measure set $N_{z}$ depends on $z\in W^{1,p(x)}(\Omega)$, which is an uncountable set. But from Lemma \ref{propospatii} we have that $W^{1,p(x)}(\Omega)$ is a separable Banach space. Therefore we can find a sequence $(\phi_n)_{n\geq 1}\subset W^{1,p(x)}(\Omega)$ which is dense in $W^{1,p(x)}(\Omega)$ (with respect to its norm topology). We will define the following null-measure set:

\begin{equation}
	N_{\tau}:=\bigcup_{n=1}^{\infty} N_{\phi_n}\ \Longrightarrow\ |N_{\tau}|\leq \sum_{n=1}^{\infty} |N_{\phi_n}|=0.
\end{equation}

\noindent Thus, we know that:

\begin{equation}\label{KAUF11}
	\mathcal{A}\bigl(\phi_n\bigr)-\mathcal{A}\bigl(v(t,\cdot)\bigr)-\int_{\Omega} h(t,x)\cdot\bigl(\phi_n(x)-v(t,x)\bigr)\ dx \geq 0,\ \forall\ n\geq 1,\ \forall\ t\in (\tau,T)\setminus N_{\tau}.
\end{equation}

\noindent Thence we can find extract a subsequence 

\begin{equation}
	(\phi_{n_j})_{j\geq 1}\ \text{with}\ \lim\limits_{j\to\infty} \Vert \phi_{n_j}-z\Vert_{W^{1,p(x)}(\Omega)}=0.
\end{equation}

\noindent So, for any $t\in (\tau,T)\setminus N_{\tau}$ and every $j\geq 1$:

\begin{equation}\label{KAUF12}
	\mathcal{A}\bigl(\phi_{n_j}\bigr)-\mathcal{A}\bigl(v(t,\cdot)\bigr)-\int_{\Omega} h(t,x)\cdot\bigl(\phi_{n_j}(x)-v(t,x)\bigr)\ dx \geq 0.
\end{equation}

\noindent $\blacktriangleright$ From Proposition \ref{propomathcalA1} \textbf{(2)} we know that $\mathcal{A}\in C^1\bigl(W^{1,p(x)}(\Omega)\bigr)\subset C\bigl(W^{1,p(x)}(\Omega)\bigr)$ and since $\phi_{n_j}\to z$ in $W^{1,p(x)}(\Omega)$ we deduce that:

\begin{equation}\label{KAUF13}
	\lim\limits_{j\to\infty} \mathcal{A}\bigl(\phi_{n_j}\bigr)=\mathcal{A}\bigl(z\bigr).
\end{equation}

\noindent $\blacktriangleright$ From hypothesis \textbf{(H2)} we know that $W^{1,p(x)}(\Omega)\hookrightarrow L^2(\Omega)$. Since $\phi_{n_j}\to z$ in $W^{1,p(x)}(\Omega)$, it follows that $\phi_{n_j}\to z$ in $L^2(\Omega)$. In particular $\phi_{n_j}\weak z$ in $L^2(\Omega)$. Thence, knowing that $h(t,\cdot)\in L^2(\Omega)$, we can write that:

\begin{equation}\label{KAUF14}
	\lim\limits_{j\to\infty} \int_{\Omega} h(t,x)\cdot\phi_{n_j}(x)\ dx=\int_{\Omega} h(t,x)\cdot z(x)\ dx.
\end{equation}

\noindent Making now $j\to\infty$ in \eqref{KAUF12}, and using \eqref{KAUF13} and \eqref{KAUF14}, allows us to write that:

\begin{equation}\label{KAUF15}
	\mathcal{A}\bigl(z\bigr)-\mathcal{A}\bigl(v(t,\cdot)\bigr)-\int_{\Omega} h(t,x)\cdot\bigl(z(x)-v(t,x)\bigr)\ dx \geq 0, \ \forall\ z\in W^{1,p(x)}(\Omega),\ \forall\ t\in (\tau,T)\setminus N_{\tau}.
\end{equation}

\noindent Fix now any $\phi\in W^{1,p(x)}(\Omega)$ and choose in \eqref{KAUF15} $z:=v(t,\cdot)+s\phi\in W^{1,p(x)}(\Omega)$ for some fixed arbitrary $s>0$. Thus, \eqref{KAUF15} becomes:

\begin{equation}\label{KAUF16}
	\dfrac{\mathcal{A}\bigl(v(t,\cdot)+s\phi\bigr)-\mathcal{A}\bigl(v(t,\cdot)\bigr)}{s}\geq\int_{\Omega} h(t,x)\cdot\phi(x)\ dx,\ \forall\ t\in N_{\tau}.
\end{equation}

\noindent  Taking into account that $\mathcal{A}\in C^1\bigl(W^{1,p(x)}(\Omega)\bigr)$,\footnote{From $\mathcal{A}\in C^1\bigl(W^{1,p(x)}(\Omega)\bigr)$ it follows that $\mathcal{A}$ is G\^ateaux differentiable -- see \cite[page 59]{coleman}.} and using \eqref{laKAUF16}, produces the following relation when we make $s\to 0^+$ in \eqref{KAUF16}:

\begin{equation}\label{KAUF17}
\int_{\Omega}\mathbf{a}(x,\nabla v(t,x))\cdot\nabla\phi(x)\ dx=\langle\mathcal{A}'(v(t,\cdot)),\phi\rangle=\lim\limits_{s\to 0^+} \dfrac{\mathcal{A}\bigl(v(t,\cdot)+s\phi\bigr)-\mathcal{A}\bigl(v(t,\cdot)\bigr)}{s}\geq\int_{\Omega} h(t,x)\cdot\phi(x)\ dx, 
\end{equation}

\noindent for every $t\in N_{\tau}$.

\noindent Choose now in \eqref{KAUF15} $z:=v(t,\cdot)-s\phi\in W^{1,p(x)}(\Omega)$ for some fixed arbitrary $s>0$. Thus, \eqref{KAUF15} becomes:

\begin{equation}\label{KAUF18}
	\dfrac{\mathcal{A}\bigl(v(t,\cdot)-s\phi\bigr)-\mathcal{A}\bigl(v(t,\cdot)\bigr)}{-s}\leq\int_{\Omega} h(t,x)\cdot\phi(x)\ dx,\ \forall\ t\in N_{\tau}.
\end{equation}

\noindent Again, making $s\to 0^+$ in \eqref{KAUF18} yields:

\begin{equation}\label{KAUF19}
	\int_{\Omega}\mathbf{a}(x,\nabla v(t,x))\cdot\nabla\phi(x)\ dx=\langle\mathcal{A}'(v(t,\cdot)),\phi\rangle=\lim\limits_{s\to 0^+} \dfrac{\mathcal{A}\bigl(v(t,\cdot)-s\phi\bigr)-\mathcal{A}\bigl(v(t,\cdot)\bigr)}{-s}\leq\int_{\Omega} h(t,x)\cdot\phi(x)\ dx, 
\end{equation}

\noindent for every $t\in N_{\tau}$.

\medskip

\noindent Combining \eqref{KAUF17} and \eqref{KAUF19} gives us that:

\begin{equation}\label{KAUF20}
		\int_{\Omega}\mathbf{a}(x,\nabla v(t,x))\cdot\nabla\phi(x)\ dx=\int_{\Omega} h(t,x)\cdot\phi(x)\ dx, \forall\ \phi\in W^{1,p(x)}(\Omega),\ \forall\ t\in (\tau,T)\setminus N_{\tau}.
\end{equation}

\noindent Since $\tau>0$ was set arbitrarily, we obtain for each $n\geq1$ that for $\tau_{n}:=\dfrac{T}{n+1}>0$ there is a set of null-measure $N_{\tau_n}\subset (\tau_n,T)$ such that \eqref{KAUF20} holds for any $\phi\in W^{1,p(x)}(\Omega)$ and for every $t\in (\tau_{n},T)\setminus N_{\tau_n}$. We define:

\begin{equation}
	N_0=\bigcup_{n=1}^{\infty} N_{\tau_n}\ \Longrightarrow\ |N_0|\leq\sum_{n=1}^{\infty} |N_{\tau_n}|=0.
\end{equation}

\noindent Now for any fixed $t_0\in (0,T)\setminus N_0$ there is some $n_0\geq 1$ large enough such that $t\in (\tau_{n_0},T)$. Since $t\notin N_0$ it follows that $t\notin N_{\tau_{n_0}}$, and therefore $t\in (\tau_{n_0},T)\setminus N_{\tau_{n_0}}$, which proves that \eqref{KAUF20} holds for $t_0$ and every $\varphi\in W^{1,p(x)}(\Omega)$. In conclusion:

\begin{equation}\label{KAUF21}
	\int_{\Omega}\mathbf{a}(x,\nabla v(t,x))\cdot\nabla\phi(x)\ dx=\int_{\Omega} h(t,x)\cdot\phi(x)\ dx, \forall\ \phi\in W^{1,p(x)}(\Omega),\ \forall\ t\in (0,T)\setminus N_{0}.
\end{equation}

\end{proof}	
	
%
%
%
%
%
%

%
%
%
	
\noindent We conclude from \textbf{Fact XV} and from the definition of $h$ (see Definition \ref{defhashh}) that $v\in C\bigl([0,T];L^2(\Omega)\bigr)\cap H^1\bigl((\tau,T);L^2(\Omega)\bigr)\cap L^1\bigl(0,T;W^{1,p(x)}(\Omega)\bigr)\cap \mathcal{V}_{[\varepsilon,\delta]}$, for any $\tau\in (0,T)$ is a \textbf{weak solution} of the problem \eqref{eqdpgaux}, i.e. $v(0,\cdot)=u_0\in\mathcal{U}_{[\varepsilon,\delta]}$ and for any $t\in (0,T)\setminus N_0$ and every $\phi\in W^{1,p(x)}(\Omega)$ we have that:

	\begin{equation}
	\int_{\Omega} \dfrac{\partial b(x,v(t,x))}{\partial t}\phi\ dx+\int_{\Omega} \mathbf{a}(x,\nabla v)\cdot\nabla\phi\ dx+\lambda\int_{\Omega} b(x,v)\phi\ dx=\int_{\Omega} g\phi\ dx. 
\end{equation}	

\end{proof}

\section{The general doubly nonlinear parabolic problem}\label{s7}

\begin{proposition}\label{propepsdelta}
	If $u$ is a weak solution of \eqref{eqdpg}, then $u\in\mathcal{V}_{[\varepsilon,\delta]}$.
\end{proposition}

\begin{proof} We have for a.e. $t\in (0,T)$ and any $\phi\in W^{1,p(x)}(\Omega)$ that
	
	\begin{equation}\label{ecuatiecuu}
	\int_\Omega\dfrac{\partial \overline{b}(x,u(t,x))}{\partial t}\phi\ dx+\int_{\Omega}\mathbf{a}(x,\nabla u(t,x))\cdot\nabla\phi\ dx=\int_{\Omega} \overline{f}(t,u(t,x))\phi\ dx.
	\end{equation}
	
\noindent We cannot apply the \textit{Weak Comparison Principle} directly, but we will mimic its proof and see that, from a certain point onward, the arguments are exactly the same. 

\noindent\textbf{Part I: We show that $u\geq \varepsilon$ a.e. on $(0,T)\times\Omega$}

\noindent For a fixed $t\in (0,T)$, we define for any $\tau>0,\ \phi_{\tau}(x)=\begin{cases} 1, & \varepsilon-u(t,x)\geq \tau\\ \dfrac{\varepsilon-u(t,x)}{\tau}, & |\varepsilon-u(t,x)|<\tau \\ -1, & \varepsilon-u(t,x)\leq -\tau\end{cases}$. Since $\phi_{\tau}\in W^{1,p(x)}(\Omega)$ we take as test function in \eqref{ecuatiecuu} $\phi=\phi_{\tau}^+\in W^{1,p(x)}(\Omega)^+$ and obtain

\begin{align*}
	\int_\Omega\dfrac{\partial \big [\overline{b}(x,u(t,x))-\overline{b}(x,\varepsilon)\big]}{\partial t}\phi_{\tau}^+\ dx &=\int_{\Omega} \overline{f}(t,u(t,x))\phi_{\tau}^+\ dx+\dfrac{1}{\tau}\int_{\Omega}\underbrace{\mathbf{a}(x,\nabla u)\cdot\nabla u}_{\geq 0}\ \chi_{\varepsilon-u(t,\cdot)\in (0,\tau)} \ dx\\
	&\geq \int_{\Omega} \overline{f}(t,u(t,x))\phi_{\tau}^+\ dx.
\end{align*}

\noindent Because $\lim\limits_{\tau\to 0^+} \phi_{\tau}^+=\chi_{\varepsilon>u(t,\cdot)}$ in $L^2(\Omega)$ we get that:

\begin{align*}
	\int_\Omega\dfrac{\partial \big [\overline{b}(x,u(t,x))-\overline{b}(x,\varepsilon)\big]}{\partial t}\chi_{\varepsilon>u(t,\cdot)}\ dx & \geq \int_{\Omega} \overline{f}(t,u(t,x))\chi_{\varepsilon>u(t,\cdot)}\ dx\\
	&=\int_{\Omega} \left [\underbrace{f(x,\varepsilon)}_{\geq 0}+\dfrac{\lambda_0}{2}(\underbrace{\varepsilon-u(t,x)}_{\geq 0}) \right ]\chi_{\varepsilon>u(t,\cdot)}(x)\ dx\\
	&\geq 0.
\end{align*}

\noindent From Remark \ref{remremrem} (for $\lambda=0$) we deduce that for any $t\in [0,T]$ we have $u(t,x)\geq \varepsilon$ for a.e. $x\in\Omega$.

\bigskip

\noindent\textbf{Part II: We show that $u\leq \delta$ a.e. on $(0,T)\times\Omega$}

\noindent For a fixed $t\in (0,T)$, we define for any $\tau>0,\ \phi_{\tau}(x)=\begin{cases} 1, & u(t,x)-\delta\geq \tau\\ \dfrac{u(t,x)-\delta}{\tau}, & |u(t,x)-\delta|<\tau \\ -1, & u(t,x)-\delta\leq -\tau\end{cases}$. Since $\phi_{\tau}\in W^{1,p(x)}(\Omega)$ we take as test function in \eqref{ecuatiecuu} $\phi=\phi_{\tau}^+\in W^{1,p(x)}(\Omega)^+$ and obtain

\begin{align*}
	\int_\Omega\dfrac{\partial \big [\overline{b}(x,\delta)-\overline{b}(x,u(t,x))\big]}{\partial t}\phi_{\tau}^+\ dx &=-\int_{\Omega} \overline{f}(t,u(t,x))\phi_{\tau}^+\ dx+\dfrac{1}{\tau}\int_{\Omega}\underbrace{\mathbf{a}(x,\nabla u)\cdot\nabla u}_{\geq 0}\ \chi_{u(t,\cdot)\in (0,\tau)} \ dx\\
	&\geq -\int_{\Omega} \overline{f}(t,u(t,x))\phi_{\tau}^+\ dx.
\end{align*}

\noindent Because $\lim\limits_{\tau\to 0^+} \phi_{\tau}^+=\chi_{u(t,\cdot)>\delta}$ in $L^2(\Omega)$ we get that:

\begin{align*}
	\int_\Omega\dfrac{\partial \big [\overline{b}(x,\delta)-\overline{b}(x,u(t,x))\big]}{\partial t}\chi_{u(t,\cdot)>\delta}\ dx & \geq -\int_{\Omega} \overline{f}(t,u(t,x))\chi_{u(t,\cdot)>\delta}\ dx\\
	&=\int_{\Omega} \left [\underbrace{-f(x,\delta)}_{\geq 0}+\widetilde{\lambda}_0(\underbrace{u(t,x)-\delta}_{\geq 0}) \right ]\chi_{u(t,\cdot)>\delta}(x)\ dx\\
	&\geq 0.
\end{align*}

\noindent From Remark \ref{remremrem} (for $\lambda=0$) we deduce that for any $t\in [0,T]$ we have $u(t,x)\leq \delta$ for a.e. $x\in\Omega$.

\end{proof}

\begin{theorem}\label{teoremaunicitatii}
	If the problem \eqref{eqdpg} has a weak solution and the extra hypothesis \textnormal{\textbf{(EHU)}} holds, then it is unique.
\end{theorem}

\begin{proof} Let $u_1,u_2$ be two weak solutions of \eqref{eqdpg}. From Proposition \ref{propepsdelta} we get that $u_1,u_2\in\mathcal{V}_{[\varepsilon,\delta]}$. We may write for a.e. $t\in (0,T)$ and any $\phi\in W^{1,p(x)}(\Omega)$ that:
	
	\begin{equation}
		\begin{cases} \displaystyle\int_\Omega\dfrac{\partial b(x,u_1(t,x)}{\partial t}\phi\ dx+\int_{\Omega}\mathbf{a}(x,\nabla u_1(t,x))\cdot\nabla\phi\ dx=\int_{\Omega} f(x,u_1(t,x))\phi\ dx\\[3mm] \displaystyle\int_\Omega\dfrac{\partial b(x,u_2(t,x)}{\partial t}\phi\ dx+\int_{\Omega}\mathbf{a}(x,\nabla u_2(t,x))\cdot\nabla\phi\ dx=\int_{\Omega} f(x,u_2(t,x))\phi\ dx\end{cases}.
	\end{equation}
	
\noindent Substracting these relations leads us to:

\begin{align*}
	&\int_{\Omega} \dfrac{\partial \big [b(x,u_2(t,x))-b(x,u_1(t,x))\big ]}{\partial t}\phi\ dx+\int_{\Omega}\big [\mathbf{a}(x,\nabla u_2(t,x))-\mathbf{a}(x,\nabla u_1(t,x))\big ]\cdot\nabla \phi\ dx\\ 
	&=\int_{\Omega} \big [f(x,u_2(t,x))-f(x,u_1(t,x))\big ]\phi\ dx.
\end{align*}

\noindent Defining as in the proof of the \textit{weak comparison principle} for any $\tau>0$ the function  $\phi_{\tau}=\begin{cases} 1, &  u_1(t,\cdot)-u_2(t,\cdot)\geq \tau\\[3mm] \dfrac{u_1(t,\cdot)-u_2(t,\cdot)}{\tau}, & |u_1(t,\cdot)-u_2(t,\cdot)|<\tau \\[3mm] -1, & u_1(t,\cdot)-u_2(t,\cdot)\leq -\tau\end{cases}\in W^{1,p(x)}(\Omega)$ we can set as test function $\phi=\phi_{\tau}^+$ and get that:

\begin{align*}
	&\int_{\Omega} \dfrac{\partial \big [b(x,u_2)-b(x,u_1)\big ]}{\partial t}\phi_{\tau}^+\ dx\\ 
	&=\int_{\Omega} \big [f(x,u_1)-f(x,u_2)\big ]\phi_{\tau}^+\ dx+\int_{\Omega}\big [\mathbf{a}(x,\nabla u_2)-\mathbf{a}(x,\nabla u_1)\big ]\cdot\big [\nabla u_2-\nabla u_1\big ]\chi_{u_1(t,\cdot)-u_2(t,\cdot)\in (0,\tau)}\ dx\\
	&\geq \int_{\Omega} \big [f(x,u_2)-f(x,u_1)\big ]\phi_{\tau}^+\ dx.
\end{align*}

\noindent Since $\lim\limits_{\tau\to 0^+} \phi_{\tau}^+=\chi_{u_1(t,\cdot)>u_2(t,\cdot)}$ in $L^2(\Omega)$ we obtain that:

\begin{equation}
	\int_{\Omega} \dfrac{\partial \big [b(x,u_2)-b(x,u_1)\big ]}{\partial t}\chi_{u_1(t,\cdot)>u_2(t,\cdot)}\ dx\geq \int_{\Omega} \big [f(x,u_2)-f(x,u_1)\big ]\chi_{u_1(t,\cdot)>u_2(t,\cdot)}\ dx.
\end{equation}

\noindent Using now \textbf{(EHU)} and Remark \ref{rem23} we can write for a.e. $x\in\Omega$:

\begin{equation}
	|f(x,u_2)-f(x,u_1)|\leq L_f |u_2-u_1|\leq\dfrac{L_f}{\ell_0}|b(x,u_2)-b(x,u_1)|.
\end{equation}

\noindent Therefore $\big [f(x,u_2)-f(x,u_1)\big ]\chi_{u_1(t,\cdot)>u_2(t,\cdot)}\geq -\dfrac{L_f}{\ell_0}\big [b(x,u_1)-b(x,u_2)\big ]\chi_{u_1(t,\cdot)>u_2(t,\cdot)}$ a.e. on $\Omega$. So:

\begin{equation}
	\int_{\Omega} \dfrac{\partial \big [b(x,u_2)-b(x,u_1)\big ]}{\partial t}\chi_{u_1(t,\cdot)>u_2(t,\cdot)}\ dx-\dfrac{L_f}{\ell_0}\int_{\Omega}[b(x,u_2)-b(x,u_1)\big ]\chi_{u_1(t,\cdot)>u_2(t,\cdot)} \geq 0.
\end{equation}

\noindent Now Remark \ref{remremrem} (for $\lambda=-\dfrac{L_f}{\ell_0}$) gives us that for any $t\in [0,T]$, $u_2(t,x)\geq u_1(t,x)$ for a.e. $x\in\Omega$. Switching $u_1$ and $u_2$ and repeating the same process will result in the fact that for any $t\in [0,T]$, $u_1(t,x)\geq u_2(t,x)$ for a.e. $x\in\Omega$. In conclusion for any $t\in [0,T]$ we have that $u_1(t,x)=u_2(t,x)$ for a.e. $x\in\Omega$.
\end{proof}

\bigskip

\noindent For any $u\in\mathcal{V}_{[\varepsilon,\delta]}$ we know, from Proposition \ref{propauxepsdelta}, Theorem \ref{thmveryunique}, Theorem \ref{thmauxiliar} and Theorem \ref{theoremverygeneral}, that the following problem of type \eqref{eqdpgaux} has a unique weak solution $v\in C\big ( [0,T]; L^2(\Omega)\big )\cap H^1_{\textnormal{loc}}\big ((0,T);L^2(\Omega)\big )\cap L^{1}\big (0,T;W^{1,p(x)}(\Omega)\big )\cap\mathcal{V}_{[\varepsilon,\delta]}$ if $u_0\in \mathcal{U}_{[\varepsilon,\delta]}$ and a unique weak solution $v\in C\big ( [0,T]; L^2(\Omega)\big )\cap H^1\big ((0,T);L^2(\Omega)\big )\cap L^{\infty}\big (0,T;W^{1,p(x)}(\Omega)\big )\cap\mathcal{V}_{[\varepsilon,\delta]}$ if $u_0\in \mathcal{U}_{[\varepsilon,\delta]}\cap W^{1,p(x)}(\Omega)$:

\begin{equation}
	\begin{cases} \dfrac{\partial b(x,v)}{\partial t}-\operatorname{div}\mathbf{a}\big (x,\nabla v \big )+\lambda_0 b(x,v)=f(x,u)+\lambda_0 b(x,u), & (t,x)\in (0,T)\times\Omega\\[3mm] \mathbf{a}\big (x,\nabla v \big )\cdot \nu=0, & (t,x)\in (0,T)\times\partial\Omega\\[3mm] v(0,x)=u_0(x), & x\in\Omega\end{cases}.
\end{equation}

\noindent Thus we can define the operator $\mathcal{S}:\mathcal{V}_{[\varepsilon,\delta]}\to\mathcal{V}_{[\varepsilon,\delta]}$ by $\mathcal{S}(u)=v$.

\begin{proposition}\label{propoeS} The operator $\mathcal{S}:\mathcal{V}_{[\varepsilon,\delta]}\to\mathcal{V}_{[\varepsilon,\delta]}$ has the following properties:
	
	\begin{enumerate}
		\item[\textbf{(1)}] $\mathcal{S}$ is \textbf{strictly monotone}, i.e. for any $u_1,u_2\in\mathcal{V}_{[\varepsilon,\delta]},\ u_1\not\equiv u_2$ with $u_1\leq u_2$ a.e. on $(0,T)\times\Omega$ we have that $\mathcal{S}(u_1)\leq \mathcal{S}(u_2)$ a.e. on $(0,T)\times\Omega$ and $\mathcal{S}(u_1)\not\equiv\mathcal{S}(u_2)$.
		
		\item[\textbf{(2)}] $\mathcal{S}$ is \textbf{continuous with respect to the norm of $L^2\big ((0,T)\times\Omega\big )$}, i.e. for any sequence $(u_n)_{n\geq 1}\subset\mathcal{V}_{[\varepsilon,\delta]}$ and $u\in \mathcal{V}_{[\varepsilon,\delta]}$ with $u_n\to u$ in $L^2\big ((0,T)\times\Omega\big )$, we have that $\mathcal{S}(u_n)\to\mathcal{S}(u)$ in $L^2\big ((0,T)\times\Omega\big )$. The result still holds if we replace $L^2\big ((0,T)\times\Omega\big )$ by $L^r\big ((0,T)\times\Omega\big )$, for any $r\in [1,\infty)$.
		
	\end{enumerate}
	
\end{proposition}

\begin{proof}\noindent\textbf{(1)} Let us denote $v_1=\mathcal{S}(u_1)\in\mathcal{V}_{[\varepsilon,\delta]}$ and $v_2=\mathcal{S}(u_2)\in\mathcal{V}_{[\varepsilon,\delta]}$. Remark that for a.e. $t\in (0,T)$ and any $\phi\in W^{1,p(x)}(\Omega)^+$ we may write, using \textbf{(H14)}
	
	\begin{align*}
		&\int_\Omega\dfrac{\partial b(x,v_1)}{\partial t}\phi\ dx+\int_{\Omega}\mathbf{a}(x,\nabla v_1)\cdot\nabla\phi\ dx+\lambda_0\int_{\Omega} b(x,v_1)\phi\ dx=\int_{\Omega} \big [f(x,u_1)+\lambda_0 b(x,u_1)]\phi\ dx\\
	 & \leq \int_{\Omega} \big [f(x,u_2)+\lambda_0 b(x,u_2)\big]\phi\ dx=\int_\Omega\dfrac{\partial b(x,v_2)}{\partial t}\phi\ dx+\int_{\Omega}\mathbf{a}(x,\nabla v_2)\cdot\nabla\phi\ dx+\lambda_0\int_{\Omega} b(x,v_2)\phi\ dx.
	\end{align*}
	
	\noindent Applying now the \textit{Weak parabolic comparison principle} we deduce that $v_1\leq v_2$ a.e. on $(0,T)\times\Omega$. 
	
	\noindent If $v_1\equiv v_2$ we will get that for any $\phi\in W^{1,p(x)}(\Omega)^+$
	
	\[
	\int_{\Omega} \big [f(x,u_1)+\lambda_0 b(x,u_1)]\phi\ dx=\int_{\Omega} \big [f(x,u_2)+\lambda_0 b(x,u_2)\big ]\phi\ dx.
	\]
	
	\noindent Taking $\phi\equiv 1$ we get that $\displaystyle\int_{\Omega} \underbrace{g_0(x,u_2)-g_0(x,u_1)}_{\geq 0}\ dx=0$. Therefore $g_0(x,u_2)=g(x,u_1)$ for a.e. $t\in (0,T)$ and a.e. $x\in\Omega$. But from \textbf{(H14)} this means that $u_1\equiv u_2$, which is false. Therefore $v_1\not\equiv v_2$.
	
	\medskip
	
\noindent\textbf{(2)} Denote $v_n:=\mathcal{S}(u_n)\in\mathcal{V}_{[\varepsilon,\delta]},\ \forall\ n\geq 1$ and $v=\mathcal{S}(u)\in\mathcal{V}_{[\varepsilon,\delta]}$. For each $n\geq 1$, for a.e. $t\in (0,T)$ and any $\phi\in W^{1,p(x)}(\Omega)$ we have that

\begin{align*}
	&\int_\Omega\dfrac{\partial b(x,v_n)}{\partial t}\phi\ dx+\int_{\Omega}\mathbf{a}(x,\nabla v_n)\cdot\nabla\phi\ dx+\lambda_0\int_{\Omega} b(x,v_n)\phi\ dx=\int_{\Omega} \big [f(x,u_n)+\lambda_0 b(x,u_n)]\phi\ dx\\
	&\int_\Omega\dfrac{\partial b(x,v)}{\partial t}\phi\ dx+\int_{\Omega}\mathbf{a}(x,\nabla v)\cdot\nabla\phi\ dx+\lambda_0\int_{\Omega} b(x,v)\phi\ dx=\int_{\Omega} \big [f(x,u)+\lambda_0 b(x,u)]\phi\ dx.
\end{align*}

\noindent Substracting these relations will give us

\begin{align}
	&\int_\Omega\dfrac{\partial \big [b(x,v_n)-b(x,v)\big ]}{\partial t}\phi\ dx+\int_{\Omega}\big [\mathbf{a}(x,\nabla v_n)-\mathbf{a}(x,\nabla v)\big ]\cdot\nabla\phi\ dx+\lambda_0\int_{\Omega} \big [b(x,v_n)-b(x,v)\big ]\phi\ dx \nonumber\\
	&=\int_{\Omega} \big [g_0(x,u_n)-g_0(x,u)\big ]\phi\ dx.
\end{align}

 As in the proof of the \textit{weak comparison principle} we define for each $\tau>0$ the function $\phi_{\tau}=\begin{cases} 1, & v_n(t,\cdot)-v(t,\cdot)\geq\tau \\ \dfrac{v_n(t,\cdot)-v(t,\cdot)}{\tau}, & v_n(t,\cdot)-v(t,\cdot)<\tau \\ -1 & v_n(t,\cdot)-v(t,\cdot)\leq -\tau\end{cases}\in W^{1,p(x)}(\Omega)$ and select the test function $\phi=\phi_{\tau}^+\in W^{1,p(x)}(\Omega)$. Hence we obtain:

\begin{align*}
	&\int_\Omega\dfrac{\partial \big [b(x,v_n)-b(x,v)\big ]}{\partial t}\phi_{\tau}^+\ dx+\lambda_0\int_{\Omega} \big [b(x,v_n)-b(x,v)\big ]\phi_{\tau}^+\ dx\\
	&=\int_{\Omega} \big [g_0(x,u_n)-g_0(x,u)\big ]\phi_{\tau}^+\ dx+\dfrac{1}{\tau}\int_{\Omega}\big [\mathbf{a}(x,\nabla v_n)-\mathbf{a}(x,\nabla v)\big ]\cdot\big [\nabla v_n-\nabla v\big ]\chi_{v(t,\cdot)-v_n(t,\cdot)\in (0,\tau)}\ dx\\
	&\geq\int_{\Omega} \big [g_0(x,u_n)-g_0(x,u)\big ]\phi_{\tau}^+\ dx.
\end{align*}

\noindent Since $\lim\limits_{\tau\to 0^+} \phi_{\tau}^+\to\chi_{v(t,\cdot)-v_n(t,\cdot)}$ strongly in $L^2(\Omega)$, and $g_0(\cdot, u_n(t,\cdot))-g_0(\cdot, u(t,\cdot))\in L^2(\Omega)$ we deduce that

\begin{align}\label{ecuatiaprincipala1}
	&\int_\Omega\dfrac{\partial \big [b(x,v_n)-b(x,v)\big ]}{\partial t}\chi_{v(t,\cdot)>v_n(t,\cdot)}\ dx+\lambda_0\int_{\Omega} \big [b(x,v_n)-b(x,v)\big ]\chi_{v(t,\cdot)>v_n(t,\cdot)}\ dx \nonumber\\
	&\geq \int_{\Omega} \big [g_0(x,u_n)-g_0(x,u)\big ]\chi_{v(t,\cdot)>v_n(t,\cdot)}\ dx\nonumber \\
	&\geq -\int_{\Omega} \big |g_0(x,u_n)-g_0(x,u)\big |\ dx.
\end{align}

\noindent We repeat the same process by switching the roles of $v_n$ and $v$ and obtain:

\begin{align}\label{ecuatiaprincipala2}
	&\int_\Omega\dfrac{\partial \big [b(x,v)-b(x,v_n)\big ]}{\partial t}\chi_{v_n(t,\cdot)>v(t,\cdot)}\ dx+\lambda_0\int_{\Omega} \big [b(x,v)-b(x,v_n)\big ]\chi_{v_n(t,\cdot)>v(t,\cdot)}\ dx \nonumber\\
	&\geq -\int_{\Omega} \big |g_0(x,u)-g_0(x,u_n)\big |\ dx.
\end{align}

\noindent Define the function $w_n(t,\cdot)=b(\cdot,v_n(t,\cdot))-b(\cdot,v(t,\cdot))$. The inequality \eqref{ecuatiaprincipala1} may be writen as:

\begin{equation}
	\int_{\Omega} \dfrac{\partial w_n}{\partial t}\chi_{v(t,\cdot)>v_n(t,\cdot)}\ dx 
	+\lambda_0\int_\Omega w_n \chi_{v(t,\cdot)>v_n(t,\cdot)}\ dx\geq -\int_{\Omega} \big |g_0(x,u_n)-g_0(x,u)\big |\ dx.
\end{equation}

\noindent Using the strict monotony of $b$ we come at the following relation that holds for a.e. $t\in (0,T)$:

\begin{equation}\label{ecuatieconti1}
	\int_{\Omega} \dfrac{\partial w_n^-(t,x)}{\partial t} dx 
	+\lambda_0\int_\Omega w_n^-(t,x) dx\leq \int_{\Omega} \big |g_0(x,u(t,x))-g_0(x,u_n(t,x))\big |\ dx.
\end{equation}

\noindent Similarly, the inequality \eqref{ecuatiaprincipala2} can be written as:

\begin{equation}
	\int_{\Omega} -\dfrac{\partial w_n}{\partial t}\chi_{v_n(t,\cdot)>v(t,\cdot)}\ dx 
	-\lambda_0\int_\Omega w_n \chi_{v_n(t,\cdot)>v(t,\cdot)}\ dx\geq -\int_{\Omega} \big |g_0(x,u_n)-g_0(x,u)\big |\ dx,
\end{equation}

\noindent i.e. for a.e. $t\in (0,T)$

\begin{equation}\label{ecuatieconti2}
	\int_{\Omega} \dfrac{\partial w_n^+(t,x)}{\partial t} dx 
	+\lambda_0\int_\Omega w_n^+(t,x) dx\leq \int_{\Omega} \big |g_0(x,u(t,x))-g_0(x,u_n(t,x))\big |\ dx.
\end{equation}

\noindent Adding now \eqref{ecuatieconti1} and \eqref{ecuatieconti2} we get that

\begin{equation}\label{ecuatieconti}
	\int_{\Omega} \dfrac{\partial |w_n(t,x)|}{\partial t} dx 
	+\lambda_0\int_\Omega |w_n(t,x)| dx\leq 2\int_{\Omega} \big |g_0(x,u(t,x))-g_0(x,u_n(t,x))\big |\ dx.
\end{equation}

\noindent For each $n\geq 1$, consider now $h_n:[0,T]\to [0,\infty)$, $h_n(t)=\displaystyle\int_\Omega |w_n(t,x)| dx$. Of course, $h_n(0)=0$ because $v_n(0,x)=v(0,x)=u_0$ for any $n\geq 1$. The same arguments used in the proof of the \textit{weak comparison principle} allows us to say that $h_n$ is differentiable a.e. on $(0,T)$ and from Theorem \ref{thmmaxpri} $h_n'(t)=\displaystyle\int_{\Omega}\dfrac{\partial |w_n|}{\partial t}(t,x)\ dx$ for a.e. $t\in (0,T)$. Therefore \eqref{ecuatieconti} rewrites as:

\begin{equation}
h_n'(t)+\lambda_0 h_n(t)\leq 2\int_{\Omega} \big |g_0(x,u_n(t,x))-g_0(x,u(t,x))\big |\ dx,\ \text{for a.e.}\ t\in (0,T).
\end{equation}

\noindent From \textit{Gronwall's inequality -- differential form}\footnote{See Theorem \ref{gronwalldiff} from the Appendix.} we obtain for a.e. $t\in (0,T)$ that:

\begin{align*}
0\leq h_n(t)&\leq e^{-\lambda_0 t}\underbrace{h_n(0)}_{=0}+2\int_{0}^t e^{-\lambda_0(t-s)}\left (\int_{\Omega} \big |g_0(x,u_n(s,x))-g_0(x,u(s,x))\big |\ dx \right )\ ds\\
&=2\int_{0}^t \underbrace{e^{-\lambda_0(t-s)}}_{\leq 1}\left (\int_{\Omega} \big |g_0(x,u_n(s,x))-g_0(x,u(s,x))\big |\ dx \right )\ ds\\
&\leq 2\int_{0}^t \int_{\Omega} \big |g_0(x,u_n(s,x))-g_0(x,u(s,x))\big |\ dx \ ds\\
&\leq 2\int_{0}^T \int_{\Omega} \big |g_0(x,u_n(s,x))-g_0(x,u(s,x))\big |\ dx \ ds\\
&=2\Vert g_0(\cdot,u_n(\cdot,\cdot))-g_0(\cdot,u(\cdot,\cdot))\Vert_{L^1((0,T)\times\Omega)}\\
&=2\Vert \mathcal{N}_{g_0}(u_n)-\mathcal{N}_{g_0}(u)\Vert_{L^1((0,T)\times\Omega)}.
\end{align*}

\noindent Henceforth, for each $n\geq 1$

\begin{align*}
	\Vert \mathcal{N}_{b}(v_n)-\mathcal{N}_{b}(v)\Vert_{L^1((0,T)\times\Omega)}&=\int_{0}^T\int_{\Omega} |w_n(t,x)|\ dx\ dt=\int_{0}^T h_n(t)\ dt\\
	&\leq 2T\Vert \mathcal{N}_{g_0}(u_n)-\mathcal{N}_{g_0}(u)\Vert_{L^1((0,T)\times\Omega)}.
\end{align*}

\noindent From Proposition \ref{propoverf} \textbf{(5)} and Proposition \ref{propoverb} \textbf{(4)} we deduce that $\mathcal{N}_{g_0}(u_n)=\mathcal{N}_{f}(u_n)+\lambda_0\mathcal{N}_{b}(u_n)\to \mathcal{N}_{f}(u)+\lambda_0\mathcal{N}_{b}(u)=\mathcal{N}_{g_0}(u)$ in $L^1\big ((0,T)\times\Omega\big )$. Thus $\mathcal{N}_b(v_n)\to \mathcal{N}_b(v)$ in $L^1\big ((0,T)\times\Omega\big )$.

\noindent We will show that the sequence $\left (\Vert v_n-v\Vert_{L^r ((0,T)\times\Omega)}\right )_{n\geq 1}$ converges to $0$ by showing that for each subsequence of it $\left(\Vert v_{n_k}-v\Vert_{L^r ((0,T)\times\Omega)} \right )_{k\geq 1}$ there is a further subsequence with $\Vert v_{n_{k_\ell}}-v\Vert_{L^r ((0,T)\times\Omega)}\stackrel{\ell\to\infty}{\longrightarrow} 0$.

\noindent Indeed, since $\Vert b(\cdot,v_{n_k})-b(\cdot,v)\Vert_{L^1((0,T)\times\Omega)}\to 0$ we will get that there is a subsequence such that $b(x,v_{n_{k_\ell}}(t,x))\stackrel{\ell\to\infty}{\longrightarrow} b(x,v(t,x))$ pointwise a.e. on $(0,T)\times\Omega$. We will show next that $v_{n_{k_\ell}}(t,x)\to v(t,x)$ pointwise a.e. on $(0,T)\times\Omega$. Suppose the contrary. Then there is a set $\omega\subset (0,T)\times \Omega$ with $|\omega|>0$ such that $|v_{n_{k_\ell}}(t,x)-v(t,x)|\nrightarrow 0$ for each $(t,x)\in\omega$. Therefore there is for each $(t,x)\in\omega$ an $\varepsilon_{(t,x)}>0$ and a further subsequence -- still denoted by $\left (v_{n_{k_\ell}}\right )_{\ell\geq 1}$ -- such that $|v_{n_{k_\ell}}(t,x)-v(t,x)|>\varepsilon_{(t,x)}$ for any $\ell\geq 1$. There are two cases:

\begin{itemize}
	\item If $v_{n_{k_\ell}}(t,x)>v(t,x)+\varepsilon_{(t,x)}$ then from the strict monotonicity of $b(x,\cdot)$ we get $b(x,v_{n_{k_\ell}}(t,x))-b(x,v(t,x))>\underbrace{ b(x,v(t,x)+\varepsilon_{(t,x)})-b(x,v(t,x))}_{>0}\nrightarrow 0$, as $\ell\to\infty$.
	
	\item If $v_{n_{k_\ell}}(t,x)<v(t,x)-\varepsilon_{(t,x)}$ then from the strict monotonicity of $b(x,\cdot)$ we get $b(x,v_{n_{k_\ell}}(t,x))-b(x,v(t,x))>\underbrace{b(x,v(t,x))-b(x,v(t,x)-\varepsilon_{(t,x)})}_{>0}\nrightarrow 0$, as $\ell\to\infty$.
\end{itemize}

\noindent Thus we have obtained the desired contradiction: $b(x,v_{n_{k_\ell}}(t,x))-b(x,v(t,x))\nrightarrow 0$ for any $(t,x)\in\omega$.

\noindent Therefore $v_{n_{k_\ell}}(t,x)\to v(t,x)$ pointwise a.e. on $(0,T)\times\Omega$. Now since $|v_{n_{k_\ell}}|\leq \delta\in L^r\big ((0,T)\times\Omega\big )$ for each $\ell\geq 1$ we deduce from the \textit{Lebesgue Dominated Convergence Theorem} that $v_{n_{k_\ell}}\stackrel{\ell\to\infty}{\longrightarrow} v$ in $L^r\big ((0,T)\times\Omega\big )$, i.e. $\Vert v_{n_{k_\ell}}-v\Vert_{L^r ((0,T)\times\Omega)}\stackrel{\ell\to\infty}{\longrightarrow} 0$. In conclusion $\Vert v_{n}-v\Vert_{L^r ((0,T)\times\Omega)}\stackrel{n\to\infty}{\longrightarrow} 0$ and the proof is complete.\footnote{Note that we have proved this result without using \textbf{(H13)}. }

\end{proof}

\begin{theorem}\label{theexistenceresult}
	Problem \eqref{eqdpg} has a weak solution 
	
	\begin{equation}
	\begin{cases} u\in  C\big ( [0,T]; L^2(\Omega)\big )\cap H^1_{\textnormal{loc}}\big ((0,T);L^2(\Omega)\big )\cap L^{1}\big (0,T;W^{1,p(x)}(\Omega)\big )\cap\mathcal{V}_{[\varepsilon,\delta]}, & \text{if} \ u_0\in\mathcal{U}_{[\varepsilon,\delta]}\\ u\in  C\big ( [0,T]; L^2(\Omega)\big )\cap H^1\big ((0,T);L^2(\Omega)\big )\cap L^{\infty}\big (0,T;W^{1,p(x)}(\Omega)\big )\cap\mathcal{V}_{[\varepsilon,\delta]}, & \text{if} \ u_0\in\mathcal{U}_{[\varepsilon,\delta]}\cap W^{1,p(x)}(\Omega)\end{cases}.
	\end{equation}
\end{theorem}

\begin{proof} We will show that $\mathcal{S}:\mathcal{V}_{[\varepsilon,\delta]}\to\mathcal{V}_{[\varepsilon,\delta]}$ has a fixed point. Consider the following recurrence:
	
	\begin{equation}
		\begin{cases}u_{n+1}=\mathcal{S}(u_n),\ n\geq 0\\[3mm] u_0\in\mathcal{V}_{[\varepsilon,\delta]}\end{cases}.
	\end{equation}
	
	\noindent $\bullet$ If we take $u_0\equiv \delta\in\mathcal{V}_{[\varepsilon,\delta]}$ then $u_1=\mathcal{S}(u_0)\in\mathcal{V}_{[\varepsilon,\delta]}\ \Rightarrow\ u_1\leq \delta=u_0$ a.e. on $\Omega$. Using the monotony of $\mathcal{S}$ we get immediately by induction that $u_2=\mathcal{S}(u_1)\leq \mathcal{S}(u_0)=u_1$, ... , $u_{n+1}=\mathcal{S}(u_n)\leq \mathcal{S}(u_{n-1})=u_n$ for any $n\geq 1$. Also, from $\delta\equiv u_0\geq \varepsilon$ a.e. on $\Omega$ we obtain that $u_1=\mathcal{S}(u_0)\geq \mathcal{S}(\varepsilon)\geq \varepsilon$ because $\mathcal{S}(\varepsilon)\in\mathcal{V}_{[\varepsilon,\delta]}$. Repeating the same process $u_2=\mathcal{S}(u_1)\geq \mathcal{S}(\varepsilon)\geq \varepsilon$ and in general we get that $u_n\geq \varepsilon$ for any $n\geq 1$. Hence:
	
	\begin{equation}
		\delta \equiv u_0\geq u_1\geq u_2\geq\dots\geq u_n\geq u_{n+1}\dots\geq \varepsilon,\ \text{a.e. on } (0,T)\times\Omega.
	\end{equation}
	
	\noindent So, for a.e. $(t,x)\in (0,T)\times\Omega$, the sequence of real numbers $\big (u_n(t,x)\big )_{n\geq 1}\subset [\varepsilon,\delta]$ is decreasing and bounded. We conclude that it is convergent to some limit denoted by $\overline{u}(t,x)\in [\varepsilon,\delta]$. We can state that $u_n\to\overline{u}$ pointwise a.e. on $(0,T)\times\Omega$. Now, since $u_n=|u_n|\leq \delta\in L^2\big ((0,T)\times\Omega\big )$ for any $n\geq 1$\footnote{This is because $(0,T)\times\Omega$ is a bounded subset of $\mathbb{R}^{N+1}$.} we obtain from the \textit{Lebesgue dominated convergence theorem} that $u_n\to \overline{u}$ in $L^2\big ((0,T)\times \Omega\big )$.

	\noindent Using Proposition \ref{propoeS} \textbf{(2)} we get that $\mathcal{S}(u_n)\to\mathcal{S}(\overline{u})$ in $L^2\big ((0,T)\times\Omega\big )$. But $\mathcal{S}(u_n)=u_{n+1}\to \overline{u}$ in $L^2\big ((0,T)\times\Omega\big )$. From the uniqueness of the limit we conclude that $\mathcal{S}(\overline{u})=\overline{u}$. So, if $u_0\in\mathcal{U}_{[\varepsilon,\delta]}$, from Theorem \ref{theoremverygeneral} we get that $\overline{u}\in C\big ( [0,T]; L^2(\Omega)\big )\cap H^1_{\textnormal{loc}}\big ((0,T);L^2(\Omega)\big )\cap L^{1}\big (0,T;W^{1,p(x)}(\Omega)\big )\cap\mathcal{V}_{[\varepsilon,\delta]}$ and $\overline{u}$ is a weak solution of the problem \eqref{eqdpg}. Similarly, if $u_0\in\mathcal{U}_{[\varepsilon,\delta]}\cap W^{1,p(x)}(\Omega)$, from Theorem \ref{thmauxiliar} we get that $\overline{u}\in C\big ( [0,T]; L^2(\Omega)\big )\cap H^1\big ((0,T);L^2(\Omega)\big )\cap L^{\infty}\big (0,T;W^{1,p(x)}(\Omega)\big )\cap\mathcal{V}_{[\varepsilon,\delta]}$ and $\overline{u}$ is a weak solution of the problem \eqref{eqdpg}.
	
	\begin{remark} If $\delta\not\equiv\mathcal{S}(\delta)$ then $\delta\gneq u_1\gneq u_2\gneq\dots\gneq u_{n}\gneq u_{n+1}\gneq\dots\gneq \overline{u}$. This follows from the strict monotony of $\mathcal{S}$.
	\end{remark}
	
	\noindent $\bullet$ If we set $u_0\equiv \varepsilon$ then in the same manner as above we will get that:
	
	\begin{equation}
		\delta\geq\hdots u_{n+1}\geq u_n\geq\dots\geq u_2\geq u_1\geq u_0\equiv \varepsilon,\ \text{a.e. on }(0,T)\times\Omega.
	\end{equation}
	
	\noindent Thus for a.a. $(t,x)\in (0,T)\times \Omega$ the sequence $\big (u_n(t,x) \big )_{n\geq 1}\subset [\varepsilon,\delta]$ is increasing and bounded, hence convergent to some $\underline{u}(t,x)\in [\varepsilon,\delta]$. So $u_n\to\underline{u}$ pointwise a.e. on $(0,T)\times\Omega$. From $|u_n|=u_n\leq \delta\in L^2\big ((0,T)\times\Omega\big )$ we deduce from the \textit{Lebesgue dominated convergence theorem} that $u_n\to\underline{u}$ in $L^2\big ((0,T)\times\Omega\big )$. 
	
	\noindent Now from the continuity of $\mathcal{S}$ with respect to the $L^2\big ((0,T)\times\Omega\big )$-norm we conclude that $u_{n+1}=\mathcal{S}(u_n)\to\mathcal{S}(\underline{u})$ in $L^2\big ((0,T)\times\Omega\big )$. But since $u_{n+1}\to \underline{u}$ in $L^2\big ((0,T)\times\Omega\big )$ we deduce that $\underline{u}=\mathcal{S}(\underline{u})$. Thus $\underline{u}\in C\big ( [0,T]; L^2(\Omega)\big )\cap H^1_{\textnormal{loc}}((0,T);L^2(\Omega)\big )\cap L^{1}\big (0,T;W^{1,p(x)}(\Omega)\big )\cap\mathcal{V}_{[\varepsilon,\delta]}$ if $u_0\in\mathcal{U}_{[\varepsilon,\delta]}$, and $\underline{u}\in C\big ( [0,T]; L^2(\Omega)\big )\cap H^1\big ((0,T);L^2(\Omega)\big )\cap L^{\infty}\big (0,T;W^{1,p(x)}(\Omega)\big )\cap\mathcal{V}_{[\varepsilon,\delta]}$ if $u_0\in\mathcal{U}_{[\varepsilon,\delta]}\cap W^{1,p(x)}(\Omega)$ is also a weak solution of the problem \eqref{eqdpg}.
	
	\begin{remark} If $\varepsilon\not\equiv\mathcal{S}(\varepsilon)$ then $\varepsilon\lneq u_1\lneq u_2\lneq\dots\lneq u_{n}\lneq u_{n+1}\lneq\dots\lneq \underline{u}$. 
	\end{remark}
	
	\bigskip
	
	\noindent Now if $u$ is any weak solution of the problem \eqref{eqdpg} then: $\underline{u}(t,x)\leq u(t,x)\leq\overline{u}(t,x)$ for a.a. $(t,x)\in (0,T)\times \Omega$. Indeed, from $u\in\mathcal{V}_{[\varepsilon,\delta]}$,\footnote{See Proposition \ref{propepsdelta}.} we have that: $\varepsilon\leq u\leq \delta$ a.e. on $(0,T)\times\Omega$. Applying $\mathcal{S}$ gives us: $\mathcal{S}(\varepsilon)\leq\mathcal{S}(u)=u\leq \mathcal{S}(\delta)$. Repeating this process $n\geq 1$ times leads us to:
	
	\begin{equation}
		\mathcal{S}^n(\varepsilon)\leq u\leq\mathcal{S}^n(\delta),\ \text{pointwise a.e. on}\ (0,T)\times\Omega.
	\end{equation}
	
	\noindent Making $n\to\infty$ we get that $\underline{u}\leq u\leq\overline{u}$ a.e. on $(0,T)\times\Omega$. Of course, if $\underline{u}=\overline{u}$ then \eqref{eqdpg} has exactly one weak solution.	
	
\end{proof}

\section{Asymptotic behaviour of the solution}\label{s8}

\noindent In this section we assume the extra hypothesis \textbf{(EHU)}. Previously, we have shown that problem \eqref{eqdpg} has a unique weak solution for any $T>0$ and any $u_0\in\mathcal{U}_{[\varepsilon,\delta]}$. Therefore it has a \textbf{unique global weak solution} defined on $[0,\infty)\times\Omega$. 

\noindent For any $\psi\in \mathcal{U}_{[\varepsilon,\delta]}$ and any $t_0\in\mathbb{R}$ we will denote by $[t_0,\infty)\times\Omega\mapsto S(t,x;t_0,\psi)$ the unique solution of the following problem:

\begin{equation}\label{eqdpgpsit0}
	\begin{cases}\dfrac{\partial b(x,u(t,x))}{\partial t}-\operatorname{div}\mathbf{a}(x,\nabla u)=f\big (x,u(t,x)\big ), & (t,x)\in (t_0,\infty)\times\Omega\\[3mm] \mathbf{a}(x,\nabla u)\cdot\nu=0, & (t,x)\in (t_0,\infty)\times\partial\Omega\\[3mm] u(t_0,x)=\psi(x)\in [\varepsilon,\delta], & x\in\Omega\end{cases}
\end{equation}

\noindent This solution is well-defined because problem \eqref{eqdpgpsit0} becomes \eqref{eqdpg} via the time-translation $\tilde{u}(t,\cdot):=u(t+t_0,\cdot)$ for any $t\in [0,\infty)$ and a.e. on $\Omega$, where $\tilde{u}(0,x)=u(t_0,x)=\psi(x)\in [\varepsilon,\delta]$ a.e. on $\Omega$. Note that for any $t\geq 0$: $\varepsilon\leq S(t,\cdot;t_0,\psi)\leq \delta$ a.e. on $\Omega$, from Proposition \ref{propepsdelta}.

\medskip

\noindent Also, denote $u_{\varepsilon}(t,\cdot):=S(t,\cdot;0,\varepsilon)$ and $u_{\delta}(t,\cdot):=S(t,\cdot;0,\delta)$ for any $t\geq 0$.

\begin{proposition}\label{propos} The following properties of $S$ hold for any $\psi\in \mathcal{U}_{[\varepsilon,\delta]}$:
	
\begin{enumerate}
	\item[\textbf{(1)}] For any $t\geq t_0$ and any real number $h\in\mathbb{R}$ we have that $S(t,x;t_0,\psi)=S(t-h,x;t_0-h,\psi)$ a.e. on $\Omega$. In particular, for $h=t_0$ we have that: $S(t,x;t_0,\psi)=S(t-t_0,x;0,\psi)$ for any $t\geq t_0$ and for a.e. on $x\in\Omega$.
	
	\item[\textbf{(2)}] $S(t,x; t_1, S(t_1,x;t_0,\psi))=S(t,x; t_0,\psi)$ for any $t\geq t_1\geq t_0$ and for a.e. $x\in\Omega$.
	
	\item[\textbf{(3)}] If $\psi_1,\psi_2\in\mathcal{U}_{[\varepsilon,\delta]}$ and $\psi_1\leq \psi_2$ a.e. on $\Omega$ then for any $t\geq t_0$ we have that $S(t,\cdot;t_0,\psi_1)\leq S(t,\cdot;t_0,\psi_2)$ a.e. on $\Omega$.\footnote{We cannot hope to obtain a result like this: if $\psi_1\not\equiv\psi_2$ then $S(t,\cdot;t_0,\psi_1)\not\equiv S(t,\cdot;t_0,\psi_2)$ for any $t\geq t_0$, since we do not have backward uniqueness, and there are many types of parabolic problems which exhibit extinction in finite time, even though they have forward uniqueness. Many examples from ODE can be imported, by making them constant with respect to the spatial variable. See \cite{campbell2003exploring}.}. In particular, for any $t\geq t_0$ we have that $u_{\varepsilon}(t-t_0,\cdot)\leq S(t,\cdot;t_0,\psi)\leq u_{\delta}(t-t_0,\cdot)$ a.e. on $\Omega$.
	
	\item[\textbf{(4)}] For any $t\geq t_0$ and any $h\geq 0$ we have that $S(t,\cdot;t_0,\varepsilon)\leq S(t+h,\cdot;t_0,\varepsilon)$, a.e. on $\Omega$. In particular for any $t\geq 0$: $u_{\varepsilon}(t,\cdot)\leq u_{\varepsilon}(t+h,\cdot)$ a.e. on $\Omega$.
	
	\item[\textbf{(5)}] For any $t\geq t_0$ and any $h\geq 0$ we have that $S(t,\cdot;t_0,\delta)\geq S(t+h,\cdot;t_0,\delta)$, a.e. on $\Omega$. In particular for any $t\geq 0$: $u_{\delta}(t,\cdot)\geq u_{\delta}(t+h,\cdot)$ a.e. on $\Omega$.
	\item[\textbf{(6)}] If $(\psi_n)_{n\geq 1}\subset\mathcal{U}_{[\varepsilon,\delta]}$ and $\psi\in\mathcal{U}_{[\varepsilon,\delta]}$ such that $\psi_n\to\psi$ in $L^1(\Omega)$ then for any $t\geq t_0$ we have that $\lim\limits_{n\to\infty} S(t,\cdot;t_0,\psi_n)=S(t,\cdot;t_0,\psi)$ in $L^r(\Omega)$ for every $r\in [1,\infty)$.
\end{enumerate}
	
\end{proposition}

\begin{proof} \noindent\textbf{(1)} Let $u(t,\cdot)=S(t,\cdot;t_0,\psi)$, $v(t)=S(t-h,\cdot;t_0-h,\psi)$ and $w(t)=S(t,\cdot;t_0-h,\psi)$. We want to show that for any $t\geq t_0$: $u(t,x)=v(t,x)$ for a.e. $x\in\Omega$. Indeed: $v(t_0,x)=S(t_0-h,x;t_0-h,\psi)=w(t_0-h,x)=\psi(x)=u(t_0,x)$ for a.e. $x\in\Omega$, and for any $\phi\in W^{1,p(x)}(\Omega)$:
	
	\begin{align}\label{equ}
		&\int_\Omega\dfrac{\partial b(x,u(t,x))}{\partial t}\phi\ dx+\int_{\Omega}\mathbf{a}(x,\nabla u(t,x))\cdot\nabla\phi\ dx=\int_{\Omega} f(t,u(t,x))\phi\ dx,\ \text{for a.e.}\ t\in (t_0,\infty).\\
		&\int_\Omega\dfrac{\partial b(x,w(t,x))}{\partial t}\phi\ dx+\int_{\Omega}\mathbf{a}(x,\nabla w(t,x))\cdot\nabla\phi\ dx=\int_{\Omega} f(t,w(t,x))\phi\ dx,\ \text{for a.e.}\ t\in (t_0-h,\infty).\nonumber
	\end{align}

	\noindent Since $v(t,\cdot)=w(t-h,\cdot)$ we may write that: 
	
	\begin{equation}\label{eqv}
		\int_\Omega\dfrac{\partial b(x,v(t,x))}{\partial t}\phi\ dx+\int_{\Omega}\mathbf{a}(x,\nabla v(t,x))\cdot\nabla\phi\ dx=\int_{\Omega} f(t,v(t,x))\phi\ dx,\ \text{for a.e.}\ t\in (t_0,\infty).
	\end{equation}
	
	\noindent Comparing now \eqref{equ} and \eqref{eqv}, and taking into account that $u(t_0,\cdot)=v(t_0,\cdot)$ a.e. on $\Omega$, we observe that $u$ and $v$ are both weak solution to the same problem. Therefore, from Theorem \ref{teoremaunicitatii}, we get that for any $t\in [t_0,\infty)$: $u(t,x)=v(t,x)$ for a.e. $x\in\Omega$.

\noindent\textbf{(2)} Let us denote $u(t,\cdot)=S(t,\cdot; t_1, S(t_1,\cdot;t_0,\psi))$ and $v(t,\cdot)=S(t,\cdot; t_0,\psi)$. Remark that $u(t_1,\cdot)=S(t_1,\cdot;t_0,\psi)=v(t_1,\cdot)$ and for any $\phi\in W^{1,p(x)}(\Omega)$:

\begin{align*}
	&\int_\Omega\dfrac{\partial b(x,u(t,x))}{\partial t}\phi\ dx+\int_{\Omega}\mathbf{a}(x,\nabla u(t,x))\cdot\nabla\phi\ dx=\int_{\Omega} f(t,u(t,x))\phi\ dx,\ \text{for a.e.}\ t\in (t_1,\infty).\\
	&\int_\Omega\dfrac{\partial b(x,v(t,x))}{\partial t}\phi\ dx+\int_{\Omega}\mathbf{a}(x,\nabla v(t,x))\cdot\nabla\phi\ dx=\int_{\Omega} f(t,v(t,x))\phi\ dx,\ \text{for a.e.}\ t\in (t_0,\infty)\subset (t_1,\infty).\nonumber
\end{align*}

\noindent Thus $u$ and $v$ are weak solutions for the same problem. Theorem \ref{teoremaunicitatii} gives us that for any $t\in [t_1,\infty)$: $u(t,x)=v(t,x)$ for a.e. $x\in\Omega$, as needed.

\noindent\textbf{(3)} Let $u_1(t,\cdot)=S(t,\cdot; t_0,\psi_1)$ and $u_2(t,\cdot)=S(t,\cdot; t_0,\psi_2)$. Observe that $u_1(t_0,\cdot)=\psi_1\leq \psi_2=u_2(t_0,\cdot)$ a.e. on $\Omega$ and for any $\phi\in W^{1,p(x)}(\Omega)$ and a.e. $t\in (t_0,\infty)$:

\begin{align*}
	&\int_\Omega\dfrac{\partial b(x,u_1(t,x))}{\partial t}\phi\ dx+\int_{\Omega}\mathbf{a}(x,\nabla u_1(t,x))\cdot\nabla\phi\ dx=\int_{\Omega} f(t,u_1(t,x))\phi\ dx\\
	&\int_\Omega\dfrac{\partial b(x,u_2(t,x))}{\partial t}\phi\ dx+\int_{\Omega}\mathbf{a}(x,\nabla u_2(t,x))\cdot\nabla\phi\ dx=\int_{\Omega} f(t,u_2(t,x))\phi\ dx.\nonumber
\end{align*}

\noindent Using the same arguments as in Theorem \ref{teoremaunicitatii} we get that for a.e. $t\in (t_0,\infty)$

\begin{equation}
	\int_{\Omega} \dfrac{\partial \big [b(x,u_2)-b(x,u_1)\big ]}{\partial t}\chi_{u_1(t,\cdot)>u_2(t,\cdot)}\ dx-\dfrac{L_f}{\ell_0}\int_{\Omega}[b(x,u_2)-b(x,u_1)\big ]\chi_{u_1(t,\cdot)>u_2(t,\cdot)} \geq 0.
\end{equation}

\noindent Hence, from Remark \ref{remremrem}, we deduce that for any $t\in [t_0,\infty)$: $u_1(t,x)\leq u_2(t,x)$ for a.e. $x\in \Omega$, as claimed.

\medskip

\begin{remark} \noindent If there is some $t_1\geq t_0$ such that $S(t_1,\cdot;t_0,\psi_1)\equiv S(t_1,\cdot;t_0,\psi_2)$ we get from \textbf{(2)} that for any $t\geq t_1$: $S(t,\cdot;t_0,\psi_1)=S(t,\cdot;t_1,S(t_1,\cdot;t_0,\psi_1))=S(t,\cdot,t_1;S(t_1,\cdot;t_0,\psi_2))=S(t,\cdot;t_0,\psi_2)$.
\end{remark}

\medskip

\noindent Finally, since $\varepsilon\leq \psi\leq \delta$ we get that for any $t\in [t_0,\infty)$: $u_{\varepsilon}(t-t_0,\cdot)=S(t-t_0,\cdot;0,\varepsilon)=S(t,\cdot;t_0,\varepsilon)\leq S(t,\cdot;t_0,\psi)\leq S(t,\cdot;t_0,\delta)\stackrel{\textbf{(1)}}{=} S(t-t_0,\cdot;0,\delta)=u_{\delta}(t-t_0,\cdot)$ a.e. on $\Omega$.

\noindent\textbf{(4)} We have for any $t\geq t_0$ and $h\geq 0$ that:

\begin{align*}
	S(t+h,\cdot;t_0,\varepsilon)&\stackrel{\textbf{(2)}}{=} S(t+h,\cdot; t_0+h,S(t_0+h,\cdot;t_0,\varepsilon))\stackrel{\textbf{(1)}}{=} S(t,\cdot;t_0,S(t_0+h,\cdot;t_0,\varepsilon))\\
	&\stackrel{\textbf{(1)}}{=} S(t,\cdot;t_0,\underbrace{S(h,\cdot;0,\varepsilon)}_{\geq \varepsilon})\\
	&\stackrel{\textbf{(3)}}{\geq} S(t,\cdot;t_0,\varepsilon),\ \text{a.e. on }\Omega.
\end{align*}

\noindent For $t_0=0$ we get that $u_{\varepsilon}(t+h,\cdot)=S(t+h,\cdot;0,\varepsilon)\geq  S(t,\cdot;0,\varepsilon)=u_{\varepsilon}(t)$ a.e. on $\Omega$.

\noindent\textbf{(5)}  Similarly, for any $t\geq t_0$ and $h\geq 0$ we have that:

\begin{align*}
	S(t+h,\cdot;t_0,\delta)&\stackrel{\textbf{(2)}}{=} S(t+h,\cdot; t_0+h,S(t_0+h,\cdot;t_0,\delta))\stackrel{\textbf{(1)}}{=} S(t,\cdot;t_0,S(t_0+h,\cdot;t_0,\delta))\\
	&\stackrel{\textbf{(1)}}{=} S(t,\cdot;t_0,\underbrace{S(h,\cdot;0,\delta)}_{\leq \delta})\\
	&\stackrel{\textbf{(3)}}{\leq} S(t,\cdot;t_0,\delta),\ \text{a.e. on }\Omega.
\end{align*}

\noindent For $t_0=0$ we get that $u_{\delta}(t+h,\cdot)=S(t+h,\cdot;0,\delta)\leq  S(t,\cdot;0,\delta)=u_{\delta}(t)$ a.e. on $\Omega$.

\noindent\textbf{(6)} For each $n\geq 1$ let us denote for all $t\geq t_0$: $u_n(t,\cdot)=S(t,\cdot;t_0,\psi_n)\in \mathcal{U}_{[\varepsilon,\delta]}\subset L^{\infty}(\Omega)$ and let $u(t,\cdot)=S(t,\cdot;t_0,\psi)\in \mathcal{U}_{[\varepsilon,\delta]}\subset L^{\infty}(\Omega)$. Also define $w_n(t,\cdot)=b(\cdot,u_n(t,\cdot))-b(\cdot,u(t,\cdot))$. We can write for each $n\geq 1$ and any $\phi\in W^{1,p(x)}(\Omega)$:

\begin{align*}
	&\int_\Omega\dfrac{\partial b(x,u_n(t,x))}{\partial t}\phi\ dx+\int_{\Omega}\mathbf{a}(x,\nabla u_n(t,x))\cdot\nabla\phi\ dx=\int_{\Omega} f(t,u_n(t,x))\phi\ dx,\ \text{for a.e.}\ t\in (t_0,\infty)\\
	&\int_\Omega\dfrac{\partial b(x,u(t,x))}{\partial t}\phi\ dx+\int_{\Omega}\mathbf{a}(x,\nabla u(t,x))\cdot\nabla\phi\ dx=\int_{\Omega} f(t,u(t,x))\phi\ dx,\ \text{for a.e.}\ t\in (t_0,\infty)
\end{align*}

\noindent Using exactly the same arguments as in the proof of Proposition \ref{propoeS} \textbf{(2)} we obtain that

\begin{equation}\label{ecuatieconti3}
	\int_{\Omega} \dfrac{\partial |w_n(t,x)|}{\partial t} dx 
	+\lambda_0\int_\Omega |w_n(t,x)| dx\leq 2\int_{\Omega} \big |g_0(x,u_n(t,x))-g_0(x,u(t,x))\big |\ dx.
\end{equation}

\noindent For each $n\geq 1$, consider now $h_n:[t_0,\infty)\to [0,\infty)$, $h_n(t)=\displaystyle\int_\Omega |w_n(t,x)| dx$. Because $w_n(t_0,x)=b(x,\psi_n(x))-b(x,\psi(x))$, we get that:

\begin{equation}
	h_n(t_0)=\int_{\Omega} |b(x,\psi_n(x))-b(x,\psi(x))|\ dx, \text{for any}\ n\geq 1.
\end{equation}

\noindent The same arguments used in the proof of the \textit{weak comparison principle} allows us to say that $h_n$ is differentiable a.e. on $(t_0,\infty)$ and from Theorem \ref{thmmaxpri} $h_n'(t)=\displaystyle\int_{\Omega}\dfrac{\partial |w_n|}{\partial t}(t,x)\ dx$ for a.e. $t\in (t_0,\infty)$. Also, note that for a.e. $(t,x)\in (t_0,\infty)\times\Omega$ we have that:

\begin{align*}
|g_0(x,u_n(t,x))-g_0(x,u(t,x))\big |&=\left |f(x,u_n(t,x))-f(x,u(t,x))+\lambda_0 \big [b(x,u_n(t,x))-b(x,u(t,x))\big ]\right |\\
&\leq |f(x,u_n(t,x))-f(x,u(t,x))|+\lambda_0 |b(x,u_n(t,x))-b(x,u(t,x))|\\
\textbf{(EHU)} \ \ \ &\leq L_f |u_n(t,x)-u(t,x)|+\lambda_0 |b(x,u_n(t,x))-b(x,u(t,x))|\\
\text{Remark \ref{rem23}}\ \ \ &\leq \dfrac{L_f}{\ell_0} |b(x,u_n(t,x))-b(x,u(t,x))|+\lambda_0 |b(x,u_n(t,x))-b(x,u(t,x))|\\
&=\left (\dfrac{L_f}{\ell_0}+\lambda_0\right ) |w_n(t,x)|.
\end{align*}

\noindent Therefore \eqref{ecuatieconti3} rewrites as:

\begin{equation}
	h_n'(t)+\lambda_0 h_n(t)\leq 2\left (\dfrac{L_f}{\ell_0}+\lambda_0\right )\int_{\Omega} \big |w_n(t,x)\big |\ dx=2\left (\dfrac{L_f}{\ell_0}+\lambda_0\right )h_n(t),\ \text{for a.e.}\ t\in (t_0,\infty).
\end{equation}

\noindent So, we have that:

\begin{equation}
		h_n'(t)\leq \left (2\dfrac{L_f}{\ell_0}+\lambda_0\right )h_n(t),\ \text{for a.e.}\ t\in (t_0,\infty).
\end{equation} 

\noindent Applying \textit{Gronwall's inequality -- differential form}\footnote{This is Theorem \ref{gronwalldiff} from the Appendix.} we obtain for a.e. $t\in (t_0,\infty)$ that:

\begin{align*}
	0\leq h_n(t)&\leq e^{\left (2\frac{L_f}{\ell_0}+\lambda_0\right ) (t-t_0)}h_n(t_0).
\end{align*}

\noindent Since $h_n$ is a continuous function on $[t_0,\infty)$ we deduce that the above inequality holds in fact for every $t\in [t_0,\infty)$.

\noindent So far we have proved that for each $n\geq 1$ and for every $t\in [t_0,\infty)$:

\begin{align*}
	 \int_{\Omega} |b(x,u_n(t,x))-b(x,u(t,x))|\ dx &\leq e^{\left (2\frac{L_f}{\ell_0}+\lambda_0\right ) (t-t_0)} \int_{\Omega} |b(x,\psi_n(x))-b(x,\psi(x))|\ dx\\
\text{Prop. \ref{propoext} \textbf{(3)}}\ \ \ 	 &\leq L_0e^{\left (2\frac{L_f}{\ell_0}+\lambda_0\right ) (t-t_0)}\int_{\Omega}  |\psi_n(x)-\psi(x)|\ dx\\
&=L_0e^{\left (2\frac{L_f}{\ell_0}+\lambda_0\right ) (t-t_0)}\Vert \psi_n-\psi\Vert_{L^1(\Omega)}\stackrel{n\to\infty}{\longrightarrow} 0.
\end{align*}

\noindent Thus $\lim\limits_{n\to\infty} \Vert b(\cdot,u_n(t,\cdot))-b(\cdot,u(t,\cdot))\Vert_{L^1(\Omega)}=0$. We want to show that $\lim\limits_{n\to\infty} \Vert u_n(t,\cdot)-u(t,\cdot)\Vert_{L^r(\Omega)}=0$. To do this we will show that from any subsequence $\left ( \Vert u_{n_k}(t,\cdot)-u(t,\cdot)\Vert_{L^r(\Omega)} \right )_{k\geq 1}$ we can extract a further subsequence with $\lim\limits_{\ell\to\infty} \Vert u_{n_{k_{\ell}}}(t,\cdot)-u(t,\cdot)\Vert_{L^r(\Omega)}=0$.

\noindent Indeed, since $\lim\limits_{k\to\infty} \Vert b(\cdot,u_{n_k}(t,\cdot))-b(\cdot,u(t,\cdot))\Vert_{L^1(\Omega)}=0$ we get that on a subsequence that $\lim\limits_{\ell\to\infty} b(x,u_{n_{k_{\ell}}}(t,x))=b(x,u(t,x))$ pointwise for a.e. $x\in\Omega$.

\noindent  Next, we show that $u_{n_{k_\ell}}(t,\cdot)\to u(t,\cdot)$ pointwise a.e. on $\Omega$. Suppose the contrary. Then there is a set $\Omega_0\subset\Omega$ with $|\Omega_0|>0$ such that $|u_{n_{k_\ell}}(t,x)-u(t,x)|\nrightarrow 0$ for each $x\in\Omega_0$. Therefore there is for each $x\in\Omega_0$ an $\varepsilon_x>0$ and a further subsequence (still denoted by $u_{n_{k_\ell}}(t,x)$) such that $|u_{n_{k_\ell}}(t,x)-u(t,x)|>\varepsilon_x$ for any $\ell\geq 1$. There are two cases:

\begin{itemize}
	\item If $u_{n_{k_\ell}}(t,x)>u(t,x)+\varepsilon_x$ then from the strict monotonicity of $b(x,\cdot)$ we get that $b(x,u_{n_{k_\ell}}(t,x))-b(x,u(t,x))>\underbrace{\big (b(x,u(t,x)+\varepsilon_x)-b(x,u(t,x))\big )}_{>0}\nrightarrow 0$, as $\ell\to\infty$.
	
	\item If $u_{n_{k_\ell}}(t,x)<u(t,x)-\varepsilon_x$ then from the strict monotonicity of $b(x,\cdot)$ we get that $b(x,u_{n_{k_\ell}}(t,x))-b(x,u(t,x))>\underbrace{\big (b(x,u(t,x)-b(x,u(t,x)-\varepsilon_x)\big )}_{>0}\nrightarrow 0$, as $\ell\to\infty$.
\end{itemize}

\noindent Thus we have obtained the desired contradiction: $b(x,u_{n_{k_\ell}}(t,x))-b(x,u(t,x)) \nrightarrow 0$ for $x\in\Omega_0$. Up to this point we have that $u_{n_{k_\ell}}(t,\cdot)\to u(t,\cdot)$ pointwise a.e. on $\Omega$. But sice $|u_{n_{k_\ell}}(t,\cdot)|\leq\delta\in L^r(\Omega)$ we obtain from \textit{Lebesgue dominated convergence theorem} that $u_{n_{k_\ell}}(t,\cdot)\to u(t,\cdot)$ in $L^r(\Omega)$, i.e. $\Vert u_{n_{k_\ell}}(t,\cdot)-u(t,\cdot)\Vert_{L^r(\Omega)}\to 0$. So for the sequence $\big (\Vert u_n(t,\cdot)-u(t,\cdot)\Vert_{L^r(\Omega)}\big )_{n\geq 1}$ we have proved that any subsequence of it has a further subsequence that converges to $0$. In conclusion  $\lim\limits_{n\to\infty} \Vert u_n(t,\cdot)-u(t,\cdot)\Vert_{L^r(\Omega)}=0$.

\noindent

\end{proof}

\begin{theorem}\label{eqasym}
	If $u$ is the unique global weak solution of the problem \eqref{eqdpg}, and if the stationary problem \eqref{eqedg} has exactly one weak solution $U$ then $\lim\limits_{t\to\infty} \Vert u(t,\cdot)-U\Vert_{L^r(\Omega)}=0$ for any $r\in [1,\infty)$.\footnote{No matter what the initial value $u_0\in\mathcal{U}_{[\varepsilon,\delta]}$ is, the weak solution of \eqref{eqdpg} converges in $L^r(\Omega)$ to the unique steady-state of \eqref{eqedg}, as $t\to\infty$.}
\end{theorem}

\begin{proof} We have that $u=S(\cdot,\cdot;0,u_0)$. We divide the proof into 4 steps:
	
	\medskip
	
\noindent\textbf{Step I: Defining $u_{\infty}(x)$ and $u^{\infty}(x)$ for a.e. $x\in\Omega$.} 

\noindent From Proposition \ref{propos} \textbf{(4)} we get that there is a sequence of null-measure sets $(Z_n)_{n\geq 0}\subset\Omega$ such that:

\begin{align*}
	&\varepsilon=u_{\varepsilon}(0,x)\leq u_{\varepsilon}(1,x),\ \forall\ x\in\Omega\setminus Z_0
	&u_{\varepsilon}(1,x)\leq u_{\varepsilon}(2,x),\ \forall\ x\in \Omega\setminus (Z_0\cup Z_1)\\
	&u_{\varepsilon}(2,x)\leq u_{\varepsilon}(3,x),\ \forall\ x\in \Omega\setminus (Z_0\cup Z_1\cup Z_2)\\
	&\vdots\qquad\vdots\qquad\vdots\\
	&u_{\varepsilon}(n,x)\leq u_{\varepsilon}(n+1,x),\ \forall\ x\in \Omega\setminus\left (\bigcup_{k=0}^{n} Z_k  \right )\\
	&\vdots\qquad\vdots\qquad\vdots\\
\end{align*}

\noindent Therefore, for each $x\in \Omega\setminus Z$, where $Z:=\displaystyle\bigcup_{n=0}^{\infty} Z_n$ is also a null-measure set\footnote{$|Z|\leq \sum_{n=0}^{\infty}|Z_n|=0$.}, we have that:

\begin{equation}
	\varepsilon=u_{\varepsilon}(0,x)\leq u_{\varepsilon}(1,x)\leq u_{\varepsilon}(2,x)\leq\hdots\leq u_{\varepsilon}(n,x)\leq u_{\varepsilon}(n+1,x)\leq\hdots\leq \delta.
\end{equation}

\noindent So, for each $x\in\Omega\setminus Z$ we have that the sequence $(u_{\varepsilon}(n,x))_{n\geq 0}$ is increasing and bounded, therefore convergent to a limit that shall be denoted by $u_{\infty}(x)\in [\varepsilon,\delta]$. We have proved that for a.e. $x\in\Omega$ there exists $u_{\infty}(x)=\lim\limits_{n\to\infty} u(n,x)$. Being a limit of a sequence of measurable functions we get that $u_{\infty}$ is also a measurable function. Thus $u_{\infty}\in\mathcal{U}_{[\varepsilon,\delta]}$.

\noindent Similarly, now from Proposition \ref{propos} \textbf{(5)} we get that there is a sequence of null-measure sets, also denoted for simplicity by $(Z_n)_{n\geq 0}\subset\Omega$, such that:

\begin{align*}
	&\delta=u_{\delta}(0,x)\geq u_{\delta}(1,x),\ \forall\ x\in \Omega\setminus Z_0\\
	&u_{\delta}(1,x)\geq u_{\delta}(2,x),\ \forall\ x\in \Omega\setminus (Z_0\cup Z_1)\\
	&u_{\delta}(2,x)\geq u_{\delta}(3,x),\ \forall\ x\in \Omega\setminus (Z_0\cup Z_1\cup Z_2)\\
	&\vdots\qquad\vdots\qquad\vdots\\
	&u_{\delta}(n,x)\geq u_{\delta}(n+1,x),\ \forall\ x\in \Omega\setminus\left (\bigcup_{k=0}^{n} Z_k  \right )\\
	&\vdots\qquad\vdots\qquad\vdots\\
\end{align*}

\noindent Therefore, for each $x\in \Omega\setminus Z$, where $Z:=\displaystyle\bigcup_{n=0}^{\infty} Z_n$ is also a null-measure set, we have that:

\begin{equation}
	\delta=u_{\delta}(0,x)\geq u_{\delta}(1,x)\geq u_{\delta}(2,x)\geq\hdots\geq u_{\delta}(n,x)\geq u_{\geq}(n+1,x)\geq\hdots\geq \varepsilon.
\end{equation}

\noindent So, for each $x\in\Omega\setminus Z$ we have that the sequence $(u_{\delta}(n,x))_{n\geq 0}$ is decreasing and bounded, therefore convergent to a limit that shall be denoted by $u^{\infty}(x)\in [\varepsilon,\delta]$. We have showed that for a.e. $x\in\Omega$ there exists $u^{\infty}(x)=\lim\limits_{n\to\infty} u(n,x)$. Being a limit of a sequence of measurable functions we get that $u^{\infty}$ is also a measurable function. Thus $u^{\infty}\in\mathcal{U}_{[\varepsilon,\delta]}$.

\begin{remark}\label{remuinfty}
	Note that from Proposition \ref{propos} \textbf{(3)} we have that there is a null-measure set $Z\subset\Omega$ such that for any $n\in\mathbb{N}$ we have that $u_{\varepsilon}(n,x)=S(n,x;0,\varepsilon)\leq S(n,x;0,\delta)=u_{\delta}(n,x)$ for any $x\in\Omega\setminus Z$. Therefore we obtain by making $n\to\infty$ that $u_{\infty}(x)\leq u^{\infty}(x)$ for a.e. $x\in\Omega$. Therefore we may write for any $n\in\mathbb{N}$ that:
	
\begin{equation}
	u_{\varepsilon}(n,x)\leq u_{\infty}(x)\leq u^{\infty}(x)\leq u_{\delta}(n,x),\ \text{for a.e.}\ x\in\Omega.
\end{equation}
\end{remark}

\medskip

\noindent\textbf{Step II: Showing that $\lim\limits_{t\to\infty} u_{\varepsilon}(t,x)=u_{\infty}(x)$ and $\lim\limits_{t\to\infty} u_{\delta}(t,x)=u^{\infty}(x)$ pointwise for a.e. $x\in\Omega$.}

\noindent We will prove that for every sequence $(t_n)_{n\geq 1}\subset [0,\infty)$ we have that $\lim\limits_{n\to\infty} u_{\varepsilon}(t_n,x)=u_{\infty}(x)$ and $\lim\limits_{n\to\infty} u_{\delta}(t_n,x)=u^{\infty}(x)$ pointwise for a.e. $x\in\Omega$.

\noindent From Proposition \ref{propos} \textbf{(4)} we get that there is a sequence of null-measure sets $(Z_n)_{n\geq 1}$ such that:

\begin{align*}
	&u_{\varepsilon}(\lceil t_1\rceil,x)\geq u_{\varepsilon}(t_1,x)\geq u_{\varepsilon}(\lfloor t_1\rfloor,x),\ \forall\ x\in \Omega\setminus Z_1\\
	&u_{\varepsilon}(\lceil t_2\rceil,x)\geq u_{\varepsilon}(t_2,x)\geq u_{\varepsilon}(\lfloor t_2\rfloor,x),\ \forall\ x\in \Omega\setminus (Z_1\cup Z_2)\\
	&\vdots\qquad\vdots\qquad\vdots\\
	&u_{\varepsilon}(\lceil t_n\rceil,x)\geq u_{\varepsilon}(t_n,x)\geq u_{\varepsilon}(\lfloor t_n\rfloor,x),\ \forall\ x\in \Omega\setminus\left (\bigcup_{k=1}^{n} Z_k  \right )\\
	&\vdots\qquad\vdots\qquad\vdots\\
\end{align*}

\noindent Thus, for any $x\in\Omega\setminus Z$, where $Z=\displaystyle\bigcup_{n=1}^{\infty} Z_n$ we have that: $u_{\varepsilon}(\lceil t_n\rceil,x)\geq u_{\varepsilon}(t_n,x)\geq u_{\varepsilon}(\lfloor t_n\rfloor,x)$ for any $n\in\mathbb{N}^*$. Therefore, passing to limit superior we get that:

\begin{equation}
	u_{\infty}(x)=\lim\limits_{n\to\infty} u_{\varepsilon}(\lceil t_n\rceil,x)=\limsup_{n\to\infty} u_{\varepsilon}(\lceil t_n\rceil,x)\geq \limsup_{n\to\infty} u_{\varepsilon}(t_n,x),\ \text{for a.e. }x\in\Omega.
\end{equation}

\noindent Similarly, passing to limit inferior we get that:

\begin{equation}
	u_{\infty}(x)=\lim\limits_{n\to\infty} u_{\varepsilon}(\lfloor t_n\rfloor ,x)=\liminf_{n\to\infty} u_{\varepsilon}(\lfloor t_n\rfloor,x)\leq \liminf_{n\to\infty} u_{\varepsilon}(t_n,x),\ \text{for a.e. }x\in\Omega.
\end{equation}

\noindent Combining the last two relations we get that for a.e. $x\in\Omega$: $\exists\lim\limits_{n\to\infty} u_{\varepsilon}(t_n,x)=u_{\infty}(x)$.

\medskip

\noindent From Proposition \ref{propos} \textbf{(5)} we get that there is a sequence of null-measure sets $(Z_n)_{n\geq 1}$ such that:

\begin{align*}
	&u_{\delta}(\lceil t_1\rceil,x)\leq u_{\delta}(t_1,x)\leq u_{\delta}(\lfloor t_1\rfloor,x),\ \forall\ x\in \Omega\setminus Z_1\\
	&u_{\delta}(\lceil t_2\rceil,x)\leq u_{\delta}(t_2,x)\geq u_{\delta}(\lfloor t_2\rfloor,x),\ \forall\ x\in \Omega\setminus (Z_1\cup Z_2)\\
	&\vdots\qquad\vdots\qquad\vdots\\
	&u_{\delta}(\lceil t_n\rceil,x)\leq u_{\delta}(t_n,x)\leq u_{\delta}(\lfloor t_n\rfloor,x),\ \forall\ x\in \Omega\setminus\left (\bigcup_{k=1}^{n} Z_k  \right )\\
	&\vdots\qquad\vdots\qquad\vdots\\
\end{align*}

\noindent Thus, for any $x\in\Omega\setminus Z$, where $Z=\displaystyle\bigcup_{n=1}^{\infty} Z_n$ we have that: $u_{\delta}(\lceil t_n\rceil,x)\leq u_{\delta}(t_n,x)\leq u_{\delta}(\lfloor t_n\rfloor,x)$ for any $n\in\mathbb{N}^*$. Therefore, passing to limit inferior we get that:

\begin{equation}
	u^{\infty}(x)=\lim\limits_{n\to\infty} u_{\delta}(\lceil t_n\rceil,x)=\liminf_{n\to\infty} u_{\delta}(\lceil t_n\rceil,x)\leq \liminf_{n\to\infty} u_{\delta}(t_n,x),\ \text{for a.e. }x\in\Omega.
\end{equation}

\noindent Similarly, passing to limit superior we get that:

\begin{equation}
	u^{\infty}(x)=\lim\limits_{n\to\infty} u_{\delta}(\lfloor t_n\rfloor ,x)=\limsup_{n\to\infty} u_{\delta}(\lfloor t_n\rfloor,x)\geq \limsup_{n\to\infty} u_{\delta}(t_n,x),\ \text{for a.e. }x\in\Omega.
\end{equation}

\noindent Combining the last two relations we get that for a.e. $x\in\Omega$: $\exists\lim\limits_{k\to\infty} u_{\delta}(t_k,x)=u^{\infty}(x)$.

\noindent Moreover, from Remark \ref{remuinfty} we have for any $t\in [0,\infty)$ that:

\begin{equation}
	u_{\varepsilon}(t,x)\leq u_{\varepsilon}(\lceil t\rceil,x)\leq u_{\infty}(x)\leq u^{\infty}(x)\leq u_{\delta}(\lceil t\rceil,x)\leq  u_{\delta}(t,x),\ \text{a.e. on}\ \Omega.
\end{equation}

\noindent\textbf{Step III: Defining $u_{*}$ and $u^*$.}

\noindent Since $u_{\infty},u^{\infty}\in\mathcal{U}_{[\varepsilon,\delta]}$ we can define $u_{*}(t,\cdot):=S(t,\cdot;0,u_{\infty})$ and $u^{*}:=S(t,\cdot;0,u^{\infty})$ for any $t\geq 0$. So $u_{*}$ and $u^{*}$ are the global weak solutions of the following problems:

	\begin{equation}
		\begin{cases}\dfrac{\partial b(x,u_*(t,x))}{\partial t}-\operatorname{div}\mathbf{a}(x,\nabla u_*(t,x))=f\big (x,u_*(t,x)\big ), & (t,x)\in (0,\infty)\times\Omega\\[3mm] \mathbf{a}(x,\nabla u_*)\cdot\nu=0, & (t,x)\in (0,\infty)\times\partial\Omega\\[3mm] u_*(0,x)=u_{\infty}(x)\in [\varepsilon,\delta], & x\in\Omega\end{cases} 
	\end{equation}
	
	\noindent and
	
\begin{equation}
	\begin{cases}\dfrac{\partial b(x,u^*(t,x))}{\partial t}-\operatorname{div}\mathbf{a}(x,\nabla u^*(t,x))=f\big (x,u^*(t,x)\big ), & (t,x)\in (0,\infty)\times\Omega\\[3mm] \mathbf{a}(x,\nabla u^*)\cdot\nu=0, & (t,x)\in (0,\infty)\times\partial\Omega\\[3mm] u^*(0,x)=u^{\infty}(x)\in [\varepsilon,\delta], & x\in\Omega\end{cases}
\end{equation}

\medskip

\noindent Now it is time to use Proposition \ref{propos} \textbf{(6)}: first, observe that for any sequence $(t_n)_{n\geq 1}\subset [0,\infty)$, with $t_n\to\infty$, we have, from \textbf{Step II}, that $\psi_n:=u_{\varepsilon}(t_n,\cdot)\stackrel{n\to\infty}{\longrightarrow} u_{\infty}(\cdot):=\psi$ pointwise a.e. on $\Omega$. Since $\psi_n,\psi\in\mathcal{U}_{[\varepsilon,\delta]}$ we get from \textit{Lebesgue Dominated Convergence Theorem} that $\psi_n\to\psi$ in $L^1(\Omega)$. Therefore for any fixed $t\geq 0$ and each $n\in\mathbb{N}^*$ we have:

\begin{align}\label{equepss}
	u_{\varepsilon}(t+t_n,\cdot)&=S(t+t_n,\cdot;0,\varepsilon)\stackrel{\text{Prop.\ref{propos}\ \textbf{(2)}}}{=}S(t+t_n,\cdot;t_n,S(t_n,\cdot;0,\varepsilon))\nonumber\\
	&\stackrel{\text{Prop.\ref{propos}\ \textbf{(1)}}}{=}S(t,\cdot;0,S(t_n,\cdot;0,\varepsilon))=S(t,\cdot;0,u_{\varepsilon}(t_n,\cdot))=S(t,\cdot;0,\psi_n)
\end{align}

\noindent From \textbf{Step II} we get that $\lim\limits_{n\to\infty} u_{\varepsilon}(t+t_n,x)=u_{\infty}(x)$ pointwise for a.e. $x\in\Omega$, and since $|u_{\varepsilon}(t+t_n,\cdot)|\leq \delta\in L^{r}(\Omega)$ for any $n\geq 1$ we deduce that:

\begin{equation}\label{equepss1}
	\lim\limits_{n\to\infty} u_{\varepsilon}(t+t_n,\cdot)=u_{\infty}\ \text{in}\ L^r(\Omega).
\end{equation}

\noindent From equation \eqref{equepss}, using Proposition \ref{propos} \textbf{(6)}, and from equation \eqref{equepss1} we obtain that:

\begin{equation}\label{equepss2}
u_{\infty}=\lim\limits_{n\to\infty} u_{\varepsilon}(t+t_n,\cdot)=\lim\limits_{n\to\infty}S(t,\cdot;0,\psi_n)=S(t,\cdot;0,\psi)=S(t,\cdot;0,u_{\infty})=u_{*}(t,\cdot)\ \text{in}\ L^r(\Omega).
\end{equation}

\noindent We have proved that for any $t\geq 0$: $u_{*}(t,\cdot)=u_{\infty}$ a.e. on $\Omega$ (it is constant with respect to the time variable). This shows that $u_{\infty}\in W^{1,p(x)}(\Omega)$ is a weak solution of the stationary problem \eqref{eqedg}, i.e. a steady-state.

\medskip 

\noindent In a similar manner, from \textbf{Step II} we get that: $u_{\delta}(t_n,\cdot)\stackrel{n\to\infty}{\longrightarrow} u^{\infty}$ pointwise a.e. on $\Omega$, and since $|u_{\delta}(t_n,\cdot)|\leq \delta$ for each $n\geq 1$ we obtain from \textit{Lebesgue Dominated Convergence Theorem} that $\lim\limits_{n\to\infty} u_{\delta}(t_n,\cdot)=u^{\infty}$ in $L^1(\Omega)$. Therefore for any fixed $t\geq 0$ and each $n\in\mathbb{N}^*$ we have:

\begin{align}\label{equdels}
	u_{\delta}(t+t_n,\cdot)&=S(t+t_n,\cdot;0,\delta)\stackrel{\text{Prop.\ref{propos}\ \textbf{(2)}}}{=}S(t+t_n,\cdot;t_n,S(t_n,\cdot;0,\delta))\nonumber\\
	&\stackrel{\text{Prop.\ref{propos}\ \textbf{(1)}}}{=}S(t,\cdot;0,S(t_n,\cdot;0,\delta))=S(t,\cdot;0,u_{\delta}(t_n,\cdot)).
\end{align}

\noindent From \textbf{Step II} we get that $\lim\limits_{n\to\infty} u_{\delta}(t+t_n,x)=u^{\infty}(x)$ pointwise for a.e. $x\in\Omega$, and since $|u_{\delta}(t+t_n,\cdot)|\leq \delta\in L^{r}(\Omega)$ for any $n\geq 1$ we deduce that:

\begin{equation}\label{equdels1}
	\lim\limits_{n\to\infty} u_{\delta}(t+t_n,\cdot)=u^{\infty}\ \text{in}\ L^r(\Omega).
\end{equation}

\noindent From equation \eqref{equdels}, using Proposition \ref{propos} \textbf{(6)}, and from equation \eqref{equdels1} we obtain that:

\begin{equation}\label{equdels2}
	u^{\infty}=\lim\limits_{n\to\infty} u_{\delta}(t+t_n,\cdot)=\lim\limits_{n\to\infty}S(t,\cdot;0,u_{\delta}(t_n,\cdot))=S(t,\cdot;0,u^{\infty})=u^{*}(t,\cdot)\ \text{in}\ L^r(\Omega).
\end{equation}

\noindent We have proved that for any $t\geq 0$: $u^{*}(t,\cdot)=u^{\infty}$ a.e. on $\Omega$ (it is constant with respect to the time variable). This shows that $u^{\infty}\in W^{1,p(x)}(\Omega)$ is a weak solution of the stationary problem \eqref{eqedg}, i.e. a steady-state.

\medskip

\noindent But from the statement of the theorem we know that \eqref{eqedg} has a unique steady-state, which allows us to write that $u_{\infty}(x)=u^{\infty}(x)=U(x)$ for a.e. $x\in\Omega$.

\medskip

\noindent\textbf{Step IV: Getting the conclusion with the squeezing principle.} 

\noindent From Proposition \ref{propos} \textbf{(3)} we easily obtain that for each $n\geq 1$ and a.e. $x\in\Omega$:

\begin{align*}
	&u_{\varepsilon}(t_n,x)=S(t_n,x;0,\varepsilon)\leq S(t_n,x;0,u_0)=u(t,x)\leq S(t_n,x;0,\delta)=u_{\delta}(t_n,x)\\
\textbf{(Step II)}\ \ \ 	& \lim\limits_{n\to\infty} u_{\varepsilon}(t_n,x)=u_{\infty}(x)=U(x)=u^{\infty}(x)=\lim\limits_{n\to\infty} u_{\delta}(t_n,x)\ \text{pointwise}.
\end{align*}

\noindent Therefore, for any $n\in\mathbb{N}^*$ and for a.e. $x\in\Omega$: $\lim\limits_{n\to\infty} u(t,x)=U(x)$ pointwise. Since $|U|, |u(t_n,\cdot)|\leq\delta\in L^r(\Omega)$ for any $n\geq 1$, using one more time \textit{Lebesgue Dominated Convergence Theorem} we obtain that:

\begin{equation}
	\lim\limits_{n\to\infty} u(t_n,\cdot)=U\ \text{in}\ L^r(\Omega),\ \text{for any sequence}\ (t_n)_{n\geq 1}\subset [0,\infty),\ t_n\to\infty.
\end{equation}

\noindent In conclusion $\lim\limits_{t\to\infty} u(t,\cdot)=U$ in $L^r(\Omega)$ for any $r\in [1,\infty)$.

\end{proof}

\begin{remark} If \textnormal{\textbf{(EH$_\Phi$)}} or \textnormal{\textbf{(EH$_f$)}} hold and $\varepsilon>0$ then the stationary problem \eqref{eqedg} has exactly one weak solution. For more details, see Theorem 6.2 from \cite{max3} and Theorems 9.1 and Theorem 10.2 from \cite{max4}. 
\end{remark}

\begin{remark}
	In a future research we can study if $\lim\limits_{t\to\infty} u(t,\cdot)=U$ in $L^{\infty}(\Omega)$.
\end{remark}

\section{An example problem}

\noindent For $\varepsilon>0$ the choice $b:\overline{\Omega}\times [\varepsilon,\delta]\to\mathbb{R}$, $b(x,s)=s^{\theta(x)}$ where $\theta:\overline{\Omega}\to (0,\infty),\ \theta\in C^1(\Omega)\cap C(\overline{\Omega})$ has the property that $\displaystyle\lim_{\substack{x\to x_0\\ x\in\Omega}} \nabla\theta(x)
\quad \text{exists for every } x_0\in\partial\Omega$, verifies the hypotheses of Proposition \ref{propolip4}, and thus there is a $C^1\bigl(\mathbb{R}^{N+1}\bigr)$ extension of $b$, denoted by the same notation $b:\mathbb{R}^{N+1}\to\mathbb{R}$. Replacing eventually $b(x,s)$ by $b(x,s)-b(x,0)$ we get that $b$ satisfies \textbf{(H9), (H10), (H11), (H12)} and \textbf{(H13)}. So we have proved the existence of a strongly positive weak solution to the \textbf{porous medium equation with variable exponents and nonlinear heterogeneous diffusion}. Here \(\theta(x)\) plays the role of a spatially variable porous-medium exponent:

\begin{equation}\label{porousmeq}
	\begin{cases}
		\dfrac{\partial}{\partial t}\left( u^{\theta(x)}(t,x)\right)
		-\operatorname{div}\mathbf{a}(x,\nabla u)
		=
		f\bigl(x,u(t,x)\bigr),
		& (t,x)\in(0,T)\times\Omega,
		\\[3mm]
		\mathbf{a}(x,\nabla u)\cdot\nu=0,
		& (t,x)\in(0,T)\times\partial\Omega,
		\\[3mm]
		u(0,x)=u_0(x)\in[\varepsilon,\delta],
		& x\in\Omega.
	\end{cases}
\end{equation}

\newpage

\section*{Appendix}

\subsection*{Lipschitz domains}

\begin{proposition}\label{propolip1} If $\Omega_1\subset\mathbb{R}^{N_1}$ and $\Omega_2\subset\mathbb{R}^{N_2}$ are two open, bounded and connected Lipschitz domains, then $\Omega:=\Omega_1\times\Omega_2\subset\mathbb{R}^{N_1+N_2}$ is also an open, bounded and connected Lipschitz domain.\footnote{See \cite[Page 269]{diaconis2011geometric}.}
\end{proposition}

\begin{proposition}\label{propolip2} Let $\Omega\subset\mathbb{R}^N$ be an open, bounded and connected Lipschitz domain. Then there is some constant $c_{\Omega}\geq 1$ such that for any two points $x,y\in \Omega$ there is a curve $\gamma:[0,1]\to\Omega$ with $\gamma(0)=x$ and $\gamma(1)=y$ such that $\gamma\in C^1([0,1];\mathbb{R}^N)$ and
	
	\begin{equation}
		\int_{0}^{1}|\gamma'(s)|\ ds=\operatorname{length}(\gamma)\leq c_{\Omega}|x-y|.
	\end{equation}
	
\end{proposition}

\begin{proof} Since any open, bounded and connected Lipschitz domain is also an uniform $(\epsilon,\infty)$ domain for some $\epsilon\in (0,1]$.\footnote{For definition and proof see \cite[Definition 8.5.1, page 276]{Hasto}, \cite[Definition 3.1]{li2010unions} and \cite[Page 73]{jones1981quasiconformal}.} Setting $c_{\Omega}:=\dfrac{1}{\epsilon}\in [1,\infty)$ we get that for any fixed points $x,y\in\Omega$, there is a continuous rectifiable curve $\tilde{\gamma}:[0,1]\to \Omega$, with $\tilde{\gamma}(0)=x$, $\tilde{\gamma}(1)=y$ and $\ell(\tilde{\gamma})\leq c_{\Omega}|x-y|$.
	
\noindent Knowing that $\Omega$ is an open set we get that for any $s\in [0,1]$, there is some $r_s>0$ such that $B(\tilde{\gamma}(s),r_s)\subset\Omega$. Since $\Gamma:=\{\tilde{\gamma}(s)\ |\ s\in [0,1]\}\subset\Omega$ is a compact set and $\Gamma\subset\displaystyle\bigcup_{s\in [0,1]} B\left (\tilde{\gamma}(s),\frac{r_s}{2}\right )$, we get that there is some $n\in\mathbb{N}^*$ and some values $0\leq s_1<s_2<\hdots<s_n\leq 1$, such that $\Gamma\subset \displaystyle\bigcup_{k=1}^n B\left (\tilde{\gamma}(s_k),\frac{r_{s_k}}{2}\right )$. Setting $r=\displaystyle\min_{k\in\overline{1,n}} \dfrac{r_{s_{k}}}{2}$ we get that $\Gamma_{r}:=\{x\in\mathbb{R}^N\ |\ \operatorname{dist}(x,\Gamma)\leq r\}\subset \displaystyle\bigcup_{k=1}^n B\left (\tilde{\gamma}(s_k),r_{s_k}\right )\subset\Omega$. Indeed, if $z\in \Gamma_r$ we have that $\operatorname{dist}(z,\Gamma)\leq r$, and since $\Gamma$ is a closed set, we get that there is some $z_0\in\Gamma$ with $|z-z_0|=\operatorname{dist}(z,\Gamma)\leq r$.\footnote{See \cite[Theorem 2.1.5 (iii)]{durea2014introduction}.} As $z_0\in\Gamma$ we get that there is some $k\in\overline{1,n}$ such that $z_0\in B\left (\tilde{\gamma}(s_k),\frac{r_{s_k}}{2}\right )$. Therefore

\begin{equation}
	|z-\tilde{\gamma}(s_k)|\leq |z-z_0|+|z_0-\tilde{\gamma}(s_k)|\leq r+\frac{r_{s_k}}{2}\leq r_{s_k}\ \Longrightarrow\ z\in B(\tilde{\gamma}(s_k),r_{s_k})\subset\Omega.
\end{equation}

\noindent Hence we may write that $\displaystyle\bigcup_{s\in [0,1]} B(\tilde{\gamma}(s),r)\subset\Gamma_r\subset\Omega$.

\noindent As $[0,1]$ is a compact set and $\tilde{\gamma}$ is continuous on $[0,1]$, we get from \textit{Heine-Cantor theorem} that $\tilde{\gamma}$ is uniformly continuous on $[0,1]$. Thus there is some $\delta_r>0$ such that if $|s-\tilde{s}|\leq\delta_r$ then $|\tilde{\gamma}(s)-\tilde{\gamma}(\tilde{s})|\leq r$, i.e. $\tilde{\gamma}(\tilde{s})\in B(\tilde{\gamma}(s),r)$. But $\tilde{\gamma}(s)\in B(\tilde{\gamma}(s),r)$ (which is a convex set), and therefore the segment $[\tilde{\gamma}(s),\tilde{\gamma}(\tilde{s})]\subset B(\tilde{\gamma}(s),r)\subset\Gamma_r\subset \Omega$.

\noindent Now we consider an equidistant division of the interval $[0,1]$ with $\left\lceil\dfrac{1}{\delta_r}\right\rceil$ points, and define the curve $\gamma_r:[0,1]\to\Gamma_r\subset\Omega$ having as image the polygonal path starting with $\tilde{\gamma}(0)=x$ and ending with $\tilde{\gamma}(1)=y$ that has as vertices the images of the division points. Since $\tilde{\gamma}$ is rectifiable, we get from the definition of rectifiability as the supremum of the length taken for all polygonal paths with vertices on $\Gamma$ that

\begin{equation}
	\ell(\gamma_r)\leq \ell(\tilde{\gamma})\leq c_{\Omega}|x-y|.
\end{equation}

\noindent Rounding with circular arcs the inner vertices of the polygonal path $\gamma_r$ we obtain a $C^1$ curve $\gamma:[0,1]\to \Gamma_r\subset\Omega$ with $\ell(\gamma)\leq\ell(\gamma_r)\leq \ell(\tilde{\gamma})\leq c_{\Omega}|x-y|$ maintaining $\gamma(0)=x$ and $\gamma(1)=y$.\footnote{The length is smaller because of the elementary inequality $x<\tan(x)$ on $(0,\pi/2)$.}
\end{proof}

\subsection*{$C^1$ functions}

\begin{proposition}\label{propolip3} If $\Omega\subset\mathbb{R}^N$ is an open, bounded and connected Lipschitz domain and $f:\Omega\to\mathbb{R}$ is a $C^1$ function for which there exists some $L>0$ such that $|\nabla f(x)|\leq L$ for any $x\in\Omega$, then $f\in \operatorname{Lip}(\Omega)$. 	
\end{proposition}

\begin{proof} Let any $x,y\in\Omega$. From Proposition \ref{propolip2} we know that there is some $c_{\Omega}\geq 1$ and a curve $\gamma:[0,1]\to \Omega$ with $\gamma\in C^1([0,1];\mathbb{R}^N)$ such that $\gamma(0)=x$, $\gamma(1)=y$ and $\int_{0}^1 |\gamma'(s)|\ ds\leq c_{\Omega}|x-y|$. Note that $f\circ\gamma:[0,1]\to\mathbb{R}$ is a $C^1$ real function with $(f\circ\gamma)'(s)=\nabla f(\gamma(s))\cdot \gamma'(t)$ for any $s\in [0,1]$. Therefore we may write that:
	
\begin{align*}
	|f(y)-f(x)|&=|f(\gamma(1))-f(\gamma(0))|=\left |\int_{0}^1 \nabla f(\gamma(s))\cdot \gamma'(s)\ ds \right |\\
\text{(Cauchy ineq.)}\ \ \ 	&\leq \int_{0}^1 |\nabla f(\gamma(s))|\cdot |\gamma'(s)|\ ds \leq L \int_{0}^1 |\gamma'(s)|\ ds\\
&\leq Lc_{\Omega} |x-y|.
\end{align*}

\noindent Therefore $f\in\operatorname{Lip}(\Omega)$.
	
\end{proof}

\begin{lemma}\label{lemmaprelcont} Let $\Omega\subset\mathbb{R}^N$ be an open set and $f:\Omega\to\mathbb{R}$ a continuous function on $\Omega$ with the property that for each $x\in\partial\Omega$ the following limit exists and it is finite $\lim\limits_{y\to x,\ y\in\Omega} f(y)$. We introduce the function $g:\overline{\Omega}\to\mathbb{R},\ g(x)=\begin{cases} f(x), & x\in\Omega\\ \lim\limits_{y\to x,\ y\in\Omega} f(y), & x\in\partial\Omega\end{cases}$. Then $g\in C(\overline{\Omega})$.
\end{lemma}

\begin{proof} Since $g=f$ on $\Omega$ we already know that $g$ is continuous on $\Omega$. We only need to show that for each $x_0\in\partial\Omega$ and any sequence $(x_n)_{n\geq 1}\subset\Omega$ with $x_n\to x_0$, we have that $g(x_n)\to g(x_0)$.
	
\noindent For each $n\geq 1$ we have that $g(x_n)=\begin{cases} f(x_n), & x_n\in\Omega\\ \lim\limits_{y\to x_n,\ y\in\Omega} f(y), & x_n\in\partial\Omega\end{cases}$. Therefore, in each case, we can choose some $y_n\in\Omega$ such that $|y_n-x_n|<\dfrac{1}{n}$ and $|f(y_n)-g(x_n)|<\dfrac{1}{n}$.

\noindent Note that $0\leq |y_n-x_0|\leq |y_n-x_n|+|x_n-x_0|\leq \dfrac{1}{n}+|x_n-x_0|\longrightarrow 0$ as $n\to\infty$, so $y_n\to x_0$. By the definition of $g$ we get that $f(y_n)\to g(x_0)$. Now observe that:

\begin{equation}
	0\leq |g(x_n)-g(x_0)|\leq |g(x_n)-f(y_n)|+|f(y_n)-g(x_0)|\leq \dfrac{1}{n}+|f(y_n)-g(x_0)|\stackrel{n\to\infty}{\longrightarrow} 0.
\end{equation}

\noindent In conclusion $g(x_n)\to g(x_0)$ and we are done.
	
\end{proof}

\begin{theorem}[\textbf{Whitney's extension theorem -- $C^1$ version}]\label{thmwhit} Let $K\subset\mathbb{R}^N$ be a compact set and let $f:K\to\mathbb{R}$ and $L:K\to\mathbb{R}^N$ be continuous functions. Then there exists a function $\tilde{f}:\mathbb{R}^N\to\mathbb{R}$ with $\tilde{f}\in C^1(\mathbb{R}^N)$ such that:\footnote{For a clear proof, see \cite[Theorem 6.10, page 277]{evans2015measure}.}
	
	\begin{itemize}
		\item $\tilde{f}(x)=f(x),\ \forall\ x\in K$,
		
		\item $\nabla \tilde{f}(x)=L(x),\ \forall\ x\in K$,
	\end{itemize}
	
	\noindent if and only if
	
	\begin{equation}
		\lim\limits_{\underset{x\neq y; x,y\in K}{|x-y|\to 0}} \dfrac{|f(y)-f(x)-L(x)\cdot (y-x)|}{|y-x|}=0.
	\end{equation}
\end{theorem}

\begin{proposition}\label{propolip4} Let $\Omega\subset\mathbb{R}^N$ be an open, bounded and connected Lipschitz domain, and $f:\overline{\Omega}\to\mathbb{R}$ a function with $f\in C(\overline{\Omega})$ and $f\in C^1(\Omega)$. Moreover, we know that for any $x\in\partial\Omega$ the following limit exists $\lim\limits_{y\to x,\ y\in\Omega} \nabla f(y)\in\mathbb{R}^N$. We define the function
	
	\[
	L:\overline{\Omega}\to\mathbb{R}^N,\ L(x)=\begin{cases} \nabla f(x),& x\in\Omega\\[2mm] \lim\limits_{y\to x,\ y\in\Omega} \nabla f(y), & x\in\partial \Omega\end{cases}.
	\]
	
\noindent The following properties hold:

\begin{enumerate}
	\item[\textbf{(1)}] $L\in C(\overline{\Omega};\mathbb{R}^N)$.
	
	\item[\textbf{(2)}] $\lim\limits_{\underset{x,y\in\overline{\Omega},\ x\neq y}{|x-y|\to 0}} \dfrac{|f(y)-f(x)-L(x)\cdot (y-x)|}{|y-x|}=0$, i.e. for any $\epsilon>0$ there is some $\delta_{\epsilon}>0$ such that for any $x,y\in\overline{\Omega}$ with $|y-x|\leq \delta_{\epsilon}$ we have that:
	
	\[
	|f(y)-f(x)-L(x)\cdot (y-x)|\leq \epsilon |y-x|.
	\]
	
	\item[\textbf{(3)}] There is a function $f_{\textnormal{ext}}:\mathbb{R}^N\to\mathbb{R},\ f_{\textnormal{ext}}\in C^1(\mathbb{R}^N)$ such that 
	
	\begin{itemize}
		\item $f_{\textnormal{ext}}(x)=f(x)$ for any $x\in\overline{\Omega}$
		
		\item $\nabla f_{\textnormal{ext}}(x)=L(x)$ for any $x\in\overline{\Omega}$. In particular $\nabla f_{\textnormal{ext}}(x)=\nabla f(x)$ for any $x\in\Omega$.
	\end{itemize}  
\end{enumerate}
	
\end{proposition}

\begin{proof} \noindent \textbf{(1)} This follows directly by applying Lemma \ref{lemmaprelcont} to the functions $\dfrac{\partial f}{\partial x_i}\in C(\Omega),\ i\in\overline{1,N}$.
	
	\medskip
	
\noindent\textbf{(2)} We denote by $c_{\Omega}\geq 1$ the constant associated to $\Omega$ from Proposition \ref{propolip2}. Fix any $\epsilon>0$. From \textbf{(1)} we know that $L:\overline{\Omega}\to\mathbb{R}^N$ is continuous, and since $\overline{\Omega}$ is a compact set, we deduce via \textit{Heine-Cantor Theorem} that $L$ is uniformly continuous on $\overline{\Omega}$. Therefore there is some constant $\eta_{\epsilon}>0$ such that 

\begin{equation}\label{ecuatieLmare}
	|L(a)-L(b)|\leq \dfrac{\epsilon}{4c_{\Omega}},\ \text{for any}\ a,b\in\Omega\ \text{with}\ |a-b|\leq \eta_{\epsilon}.
\end{equation}

\noindent We define $\delta_{\epsilon}:=\dfrac{\eta_{\epsilon}}{2 c_{\Omega}}$ and we will show that for any $x_0,y_0\in\overline{\Omega}$ with $0<|y_0-x_0|\leq\delta_{\epsilon}$ the following inequality holds
\[
|f(y_0)-f(x_0)-L(x_0)\cdot (y_0-x_0)|\leq \epsilon |y_0-x_0|.
\]

\noindent In that sense let us denote $g:\overline{\Omega}\times\overline{\Omega}\to\mathbb{R},\ g(x,y)=
|f(y)-f(x)-L(x)\cdot (y-x)|$. From the continuity of $g$ we get that there are some points $x,y\in \Omega$ with the following properties

\begin{equation}\label{ecuatiecug}
|x-x_0|,|y-y_0|\leq \dfrac{|y_0-x_0|}{2} \ \text{and}\ g(x_0,y_0)<g(x,y)+\frac{\epsilon}{2}|y_0-x_0|.
\end{equation}

\noindent Next, from Proposition \ref{propolip2} we get that there is a curve $\gamma:[0,1]\to\Omega$ such that $\gamma(0)=x$ and $\gamma(1)=y$, with $\gamma\in C^1([0,1];\mathbb{R}^N)$ and 

\begin{equation}\label{ecuatiecucomega}
	\displaystyle\int_{0}^1|\gamma'(s)|\ ds\leq c_{\Omega}|x-y|\leq c_{\Omega}\big (|x-x_0|+|x_0-y_0|+|y_0-y| \big )\leq 2c_{\Omega}|x_0-y_0|\leq 2c_{\Omega}\delta_{\epsilon}=\eta_{\epsilon}.
\end{equation}

\noindent We may write that

\begin{align*}
	g(x,y)&=|f(y)-f(x)-L(x)\cdot (y-x)|=|f(y)-f(x)-L(x)\cdot (y-x)|\\
	&=\left |\int_{0}^1 \nabla f(\gamma(s))\cdot \gamma'(s)\ ds-\int_{0}^1 L(\gamma(0))\cdot \gamma'(s)\ ds \right |\\
	&=\left |\int_{0}^1 \big [L(\gamma(s))-L(\gamma(0))\big ]\cdot \gamma'(s)\ ds \right |\\
	&\leq \int_{0}^1 \left | L(\gamma(s))-L(\gamma(0))\right |\cdot |\gamma'(s)|\ ds
\end{align*}

\noindent Now observe for any $s\in [0,1]$ that $|\gamma(s)-\gamma(0)|=\left |\displaystyle\int_{0}^s\gamma'(\tau)\ d\tau \right |\leq \int_{0}^s |\gamma'(\tau)|\ d\tau\leq \int_{0}^1 |\gamma'(\tau)|\ d\tau\leq\eta_{\epsilon}$, from \eqref{ecuatiecucomega}. Therefore from \eqref{ecuatieLmare} we obtain that:

\begin{align*}
	g(x,y)&\leq \int_{0}^1 \left | L(\gamma(s))-L(\gamma(0))\right |\cdot |\gamma'(s)|\ ds\\
	&\leq \frac{\epsilon}{4c_{\Omega}}\int_{0}^1 |\gamma'(s)|\ ds\\
\eqref{ecuatiecucomega}\ \ \	&\leq \frac{\epsilon}{4c_{\Omega}}\cdot 2c_{\Omega}|y_0-x_0|=\frac{\epsilon}{2}|y_0-x_0|.
\end{align*}

\noindent Combining this relation with \eqref{ecuatiecug} given us that:

\[
|f(y_0)-f(x_0)-L(x_0)\cdot (y_0-x_0)|=g(x_0,y_0)\leq g(x,y)+\frac{\epsilon}{2}|y_0-x_0|\leq \epsilon |y_0-x_0|.
\]

\medskip

\noindent\textbf{(3)} This follows by applying Whitney's extension theorem for $f$, using \textbf{(2)}. See Theorem \ref{thmwhit}.
	
\end{proof}

\subsection*{Variable exponent spaces}

\begin{proposition}\label{propu2u1} Let $\Omega\subset\mathbb{R}^N$ be a measurable set and $p:\Omega\to [1,\infty)$ be a measurable and bounded exponent. If $u_1:\Omega\to\mathbb{R}$ is a measurable function, $u_2\in L^{p(x)}(\Omega)$ and $|u_1(x)|\leq |u_2(x)|$ for a.e. $x\in\Omega$, then $u_1\in L^{p(x)}(\Omega)$. Moreover $\rho_{p(x)}(u_1)\leq \rho_{p(x)}(u_2)$ and $\Vert u_1\Vert_{L^{p(x)}(\Omega)}\leq \Vert u_2\Vert_{L^{p(x)}(\Omega)}$.
\end{proposition}

\begin{proof}
	We have that $\rho_{p(x)}(u_1)=\displaystyle\int_{\Omega} |u_1(x)|^{p(x)}\ dx\leq \int_{\Omega} |u_2(x)|^{p(x)}\ dx=\rho_{p(x)}(u_2)<\infty$. This show that $u_1\in L^{p(x)}(\Omega)$. Now let us denote $M_1=\left \{\lambda>0\ \bigg |\ \rho_{p(x)}\left (\dfrac{u_1}{\lambda}\right )\leq 1\right \}$ and $M_2=\left \{\lambda>0\ \bigg |\ \rho_{p(x)}\left (\dfrac{u_2}{\lambda}\right )\leq 1\right \}$. If $\lambda\in M_2$ then $1\geq \rho_{p(x)}\left (\dfrac{u_2}{\lambda}\right )=\displaystyle\int_{\Omega}\left (\dfrac{u_2(x)}{\lambda}\right )^{p(x)}\ dx\geq \displaystyle\int_{\Omega}\left (\dfrac{u_1(x)}{\lambda}\right )^{p(x)}\ dx$. Thus $\lambda\in M_1$. This shows that $M_2\subset M_1$ and as a consequence $\Vert u_1\Vert_{L^{p(x)}(\Omega)}=\inf M_1\leq \inf M_2=\Vert u_2\Vert_{L^{p(x)}(\Omega)}$.
\end{proof}

\begin{lemma}\label{propoabel} Let any $a,b\geq 0$. The following elementary inequalities hold:
	
	\begin{enumerate}
		\item[\textbf{(1)}] For any $p\geq 1$ we have that $(a+b)^p\leq 2^{p-1}(a^p+b^p)$.
		
		\item[\textbf{(2)}] For any $p\in [0,1]$ we have that $(a+b)^p\leq a^p+b^p$ and $|a^p-b^p|\leq |a-b|^p$.
		
		\item[\textbf{(3)}] For any $p\in [0,\infty)$ we have that $(a+b)^p\leq \max\{2^{p-1},1\}(a^p+b^p)$.
		
		\item[\textbf{(4)}] For any $p\in [1,\infty)$ we have that $|a^p-b^p|\leq \max\{2^{p-2},1\}\big [|a-b|^p+pb^{p-1}|a-b|\big ]$.
	\end{enumerate}
	
\end{lemma}

\begin{proof} \noindent\textbf{(1)} The function $\theta:[0,\infty)\to\mathbb{R},\ \theta(x)=x^p$ is convex, since its second derivative is $\theta''(x)=p(p-1)x^{p-2}\geq 0$ for any $x\in (0,\infty)$. Therefore $\theta\left ( \dfrac{a+b}{2}\right )\leq \dfrac{\theta(a)+\theta(b)}{2}$, i.e. $(a+b)^p\leq 2^{p-1}(a^p+b^p)$.
	
\medskip
	
\noindent\textbf{(2)} Suppose that $b\geq a$ and define the function $h:[0,\infty)\to\mathbb{R},\ h(x)=(a+x)^p-a^p-x^p$. Observe that $h(0)=0$, $h$ is continuous on $[0,\infty)$ and differentiable on $(0,\infty)$ with $h'(x)=p\big [(a+x)^{p-1}-x^{p-1}\big ]\leq 0$, because $0\leq p\leq 1$. Thus $h$ is decreasing on $[0,\infty)$ and as a consequence $h(b-a)\leq h(0)$, i.e. $b^{p}-a^{p}\leq (b-a)^p$. Moreover $h(b)\leq 0$, i.e. $(a+b)^p\leq a^p+b^p$.

\medskip

\noindent\textbf{(3)} This follows immediately by combining \textbf{(1)} and \textbf{(2)}.

\medskip

\noindent\textbf{(4)} Suppose that $a\geq b$. Therefore:

\begin{align*}
	a^p-b^p&=p\int_{0}^1 \big (b+t(a-b) \big )^{p-1} (a-b)\ dt\\
	&\stackrel{\textbf{(3)}}{\leq}p\max\{2^{p-2},1\}\int_{0}^1 \big ( b^{p-1}+(a-b)^{p-1} t^{p-1}\big )\cdot (a-b)\ dt\\
	&=\max\{2^{p-2},1\} \left [ pb^{p-1}(a-b)+p(a-b)^p\int_{0}^1 t^{p-1}\ dt\right ]\\
	&=\max\{2^{p-2},1\}\big [ pb^{p-1}(a-b)+(a-b)^p\big ].
\end{align*}
	
\end{proof}

\begin{proposition}\label{lpxprimineq}
	Let $p:\Omega\to [1,\infty)$ be a measurable and bounded exponent and consider any function $u\in L^{p(x)}(\Omega)$. Therefore $|u|^{p(x)-1}\in L^{p'(x)}(\Omega)$ and moreover, if $p^{-}>1$ then the following inequality holds:
	
	\begin{equation}
		\left \Vert |u|^{p(x)-1}\right\Vert_{L^{p'(x)}(\Omega)}\leq \max\left\{\Vert u\Vert_{L^{p(x)}(\Omega)}^{p^--1},\Vert u\Vert_{L^{p(x)}(\Omega)}^{p^+-1} \right\}.
	\end{equation}
\end{proposition}

\begin{proof} Recall that $p'(x)=\begin{cases}\dfrac{p(x)}{p(x)-1}, & x\in\Omega\setminus\Omega_1\\ +\infty, & x\in \Omega_1:=\{x\in\Omega\ |\ p(x)=1\}\end{cases}$. Using the definition of the Luxemburg norm from $L^{p'(x)}(\Omega)$ given in \cite[Definition 2.6,page 17]{cruz2013variable}, we have that $|u|^{p(x)-1}\in L^{p'(x)}(\Omega)$ iff there is some $\lambda>0$ such that $\rho_{p'(x)}\left (\dfrac{|u|^{p(x)-1}}{\lambda}\right )<\infty$. But for $\lambda=1$ we get that:
	
	\begin{align*}
		\rho_{p'(x)}\left (|u|^{p(x)-1}\right )&=\int_{\Omega\setminus\Omega_1} \left [|u|^{p(x)-1}\right ]^{p'(x)}\ dx+\Vert |u|^{p(x)-1}\Vert_{L^{\infty}(\Omega_1)}\\
		&=\int_{\Omega\setminus\Omega_1} |u|^{p(x)}\ dx+1\leq \int_{\Omega} |u|^{p(x)}\ dx+1=\rho_{p(x)}(u)+1<\infty,
	\end{align*}
	
	\noindent because $u\in L^{p(x)}(\Omega)$. 
	
	\noindent Now we work in the case $p^->1$. First we denote $\mu=\Vert u\Vert_{L^{p(x)}(\Omega)}$. If $\mu=0$ then $u\equiv 0$ and there is nothing to prove. Consider $\mu>0$. Using Proposition 2.2.1 from \cite[page 24]{cruz2013variable} we find out that $\rho_{p(x)}\left (\dfrac{u}{\mu} \right )=1$. We distinsguish two separate cases:
	
	\begin{itemize}
		\item $0<\mu\leq 1$. In this case we may write $p(x)-1\geq p^--1>0$ a.e. on $\Omega$, from where $\mu^{p(x)-1}\leq\mu^{p^--1}:=\lambda$. Raising this inequality to the power $p'(x)$ gives us that

		\begin{equation}\label{neimportant3}
			\mu^{p(x)}\leq \lambda^{p'(x)},\ \text{a.e. on}\ \Omega. 
		\end{equation}
		
		\noindent Therefore
		
		\begin{align*}
			&\left (\dfrac{|u|^{p(x)-1}}{\lambda} \right )^{p'(x)}=\dfrac{|u|^{p(x)}}{\lambda^{p'(x)}}\stackrel{\eqref{neimportant3}}{\leq} \dfrac{|u|^{p(x)}}{\mu^{p(x)}}=\left (\dfrac{|u|}{\mu} \right)^{p(x)},\ \text{a.e. on}\ \Omega\\[3mm]
		\Longrightarrow\ \ \ 	& \rho_{p'(x)}\left(\dfrac{|u|^{p(x)-1}}{\lambda}\right)=\int_{\Omega} \left (\dfrac{|u|^{p(x)-1}}{\lambda} \right )^{p'(x)}\ dx\leq \int_{\Omega} \left (\dfrac{|u|}{\mu} \right)^{p(x)}\ dx=\rho_{p(x)}\left (\dfrac{u}{\mu} \right )=1.
		\end{align*}
		
		\noindent Since for $p^->1$ it follows that $p':\Omega\to (1,\infty)$, we may write: $\left\Vert |u|^{p(x)-1}\right\Vert_{L^{p'(x)}(\Omega)}=\inf\left\{\tilde{\lambda}>0\ \bigg |\ \rho_{p'(x)}\left(\dfrac{|u|^{p(x)-1}}{\tilde{\lambda}}\right)\leq 1\right\}$. Thence we proved that 
		
		\begin{equation}\label{neimportant4}
			\left\Vert |u|^{p(x)-1}\right\Vert_{L^{p'(x)}(\Omega)}\leq \lambda=\mu^{p--1}=\Vert u\Vert_{L^{p(x)}(\Omega)}^{p^--1},\ \text{for}\ \Vert u\Vert_{L^{p(x)}(\Omega)}\leq 1.
		\end{equation}
		
		\item $\mu>1$. In this case we may write $p(x)-1\leq p^+-1>0$ a.e. on $\Omega$, from where $\mu^{p(x)-1}\leq\mu^{p^+-1}:=\lambda$. Raising this inequality to the power $p'(x)$ gives us that
		
		\begin{equation}\label{neimportant5}
			\mu^{p(x)}\leq \lambda^{p'(x)},\ \text{a.e. on}\ \Omega. 
		\end{equation}
		
		\noindent Therefore
		
		\begin{align*}
			&\left (\dfrac{|u|^{p(x)-1}}{\lambda} \right )^{p'(x)}=\dfrac{|u|^{p(x)}}{\lambda^{p'(x)}}\stackrel{\eqref{neimportant5}}{\leq} \dfrac{|u|^{p(x)}}{\mu^{p(x)}}=\left (\dfrac{|u|}{\mu} \right)^{p(x)},\ \text{a.e. on}\ \Omega\\[3mm]
			\Longrightarrow\ \ \ 	& \rho_{p'(x)}\left(\dfrac{|u|^{p(x)-1}}{\lambda}\right)=\int_{\Omega} \left (\dfrac{|u|^{p(x)-1}}{\lambda} \right )^{p'(x)}\ dx\leq \int_{\Omega} \left (\dfrac{|u|}{\mu} \right)^{p(x)}\ dx=\rho_{p(x)}\left (\dfrac{u}{\mu} \right )=1.
		\end{align*}
		
		\noindent Now we may write: $\left\Vert |u|^{p(x)-1}\right\Vert_{L^{p'(x)}(\Omega)}=\inf\left\{\tilde{\lambda}>0\ \bigg |\ \rho_{p'(x)}\left(\dfrac{|u|^{p(x)-1}}{\tilde{\lambda}}\right)\leq 1\right\}$. Thence we proved that 
		
		\begin{equation}\label{neimportant6}
			\left\Vert |u|^{p(x)-1}\right\Vert_{L^{p'(x)}(\Omega)}\leq \lambda=\mu^{p^+-1}=\Vert u\Vert_{L^{p(x)}(\Omega)}^{p^+-1},\ \text{for}\ \Vert u\Vert_{L^{p(x)}(\Omega)}>1.
		\end{equation}
	\end{itemize}
	
	\noindent Combining now \eqref{neimportant4} and \eqref{neimportant6} completes the proof.

\end{proof}

\begin{proposition}[\textbf{Continuity of the modular}]\label{propomoducont} If $\Omega\subset\mathbb{R}^N$ is a measurable set and $p:\Omega\to [1,\infty)$ with $p^+<\infty$ is a bounded measurable exponent, then if $u_n\to u$ in $L^{p(x)}(\Omega)$ we have that $\rho_{p(x)}(u_n)\to\rho_{p(x)}(u)$. Moreover one has that:
	
	\begin{equation}
		\lim\limits_{n\to\infty} \int_{\Omega}\big | |u_n(x)|^{p(x)}-|u(x)|^{p(x)}\big |\ dx=0.
	\end{equation}
	
\end{proposition}

\begin{proof} $\Longrightarrow\ \ \ $ $u_n\to u$ in $L^{p(x)}(\Omega)$ we get that $\rho_{p(x)}(u_n-u)\to 0$, from \cite[Proposition 2.58]{cruz2013variable}. Moreover:
	
	\begin{align*}
		\big |\rho_{p(x)}(u_n)-\rho_{p(x)}(u)\big |&=\left | \int_{\Omega} |u_n(x)|^{p(x)}-|u(x)|^{p(x)}\ dx \right |\leq \int_{\Omega} \left | |u_n(x)|^{p(x)}-|u(x)|^{p(x)} \right |\ dx\\
	\text{Lemma \ref{propoabel} \textbf{(4)}}\ \ \ 	&\leq\max\{2^{p^+-2},1\}\int_{\Omega}\big [p^+|u|^{p(x)-1}|u_n-u|+|u_n-u|^{p(x)} \big ]\ dx\\
	&=p^+\max\{2^{p^+-2},1\}\int_{\Omega} |u|^{p(x)-1}|u_n-u|\ dx+\max\{2^{p^+-2},1\}\rho_{p(x)}(u_n-u)\\
\text{H\"{o}lder ineq.}\ \ \	&\leq p^+\max\{2^{p^+-2},1\}\left (\dfrac{1}{p^-}-\dfrac{1}{p^+}+2\right )\Vert |u|^{p(x)-1}\Vert_{L^{p'(x)}(\Omega)}\cdot \Vert u_n-u\Vert_{L^{p(x)}(\Omega)}\\
&\ \ \ +\max\{2^{p^+-2},1\}\rho_{p(x)}(u_n-u)\stackrel{n\to\infty}{\longrightarrow}\ 0.
	\end{align*}
	
\noindent Here we have used the \textit{H\"{o}lder inequality for variable exponent spaces}\footnote{See Theorem \cite[Theorem 2.26]{cruz2013variable}.} and the Proposition \ref{lpxprimineq}.

\end{proof}

\begin{proposition} Let $\Omega\subset\mathbb{R}^N$ be a measurable set and $p:\Omega\to [1,\infty)$ with $p^+<\infty$ be a bounded measurable exponent. Consider a sequence of functions $u_n:\Omega\to\mathbb{R}, n\geq 1$ from $L^{p(x)}(\Omega)$ that converges pointwise a.e. to a function $u\in L^{p(x)}(\Omega)$. If $\rho_{p(x)}(u_n)\to\rho_{p(x)}(u)$, then $u_n\to u$ in $L^{p(x)}(\Omega)$.
	
\end{proposition}

\begin{proof} Consider the sequence of measurable functions 
	
	\begin{equation}
		h_n=2^{p^+-1}\big (|u_n|^{p(x)}+|u|^{p(x)}\big )-|u-u_n|^{p(x)},\ n\geq 1.
	\end{equation}

\noindent We have that $h_n\to 2^{p^+}|u|^{p(x)}$ pointwise a.e. on $\Omega$. Also, from Lemma \ref{propoabel} \textbf{(1)} we get that $h_n\geq 0$ a.e. on $\Omega$ for each $n\geq 1$. Applying \textit{Fatou's Lemma} leads us to
	
	\begin{align*}
	2^{p^+}\rho_{p(x)}(u)&=2^{p^+}\int_{\Omega} |u|^{p(x)}\ dx=\int_{\Omega}\liminf_{n\to\infty} h_n(x)\ dx\leq \liminf_{n\to\infty} \int_{\Omega} h_n(x)\ dx\\
	&=\liminf_{n\to\infty} 2^{p^+-1}\int_{\Omega} |u_n|^{p(x)}\ dx+2^{p^+-1}\int_{\Omega}|u|^{p(x)}\ dx-\int_{\Omega} |u_n-u|^{p(x)}\ dx\\
	&=\liminf_{n\to\infty} 2^{p^+-1} \rho_{p(x)}(u_n)+2^{p^+-1} \rho_{p(x)}(u)-\rho_{p(x)}(u_n-u)\\
	&=2^{p^+}\rho_{p(x)}(u)+\liminf_{n\to\infty}-\rho_{p(x)}(u_n-u)\\
	&=2^{p^+}\rho_{p(x)}(u)-\limsup_{n\to\infty}\rho_{p(x)}(u_n-u).
	\end{align*}  
	
\noindent From here we deduce that

\begin{equation}
	0\leq\liminf_{n\to\infty} \rho_{p(x)}(u_n-u)\leq \limsup_{n\to\infty} \rho_{p(x)}(u_n-u)\leq 0,
\end{equation}

\noindent i.e. $\displaystyle\lim_{n\to\infty} \rho_{p(x)}(u_n-u)=0$. This shows that $u_n\to u$ in $L^{p(x)}(\Omega)$.
\end{proof}

\begin{theorem}[\textbf{General Lebesgue DCT for variable exponents}]\label{generaldct} Consider a measurable set $\Omega\subset\mathbb{R}^N$ and $p:\Omega\to [1,\infty)$ with $p^+<\infty$ a measurable and bounded exponent. Let $f_n:\Omega\to\mathbb{R}$ be a sequence of measurable functions such that there is a function $f:\Omega\to \mathbb{R}$ with $f_n\to f$ pointwise a.e. on $\Omega$. We know that there is a sequence of functions $(g_n)_{n\geq 1}\subset L^{p(x)}(\Omega)$ with the following properties
	
	\begin{itemize}
		\item For each $n\in\mathbb{N}^*$: $|f_n(x)|\leq |g_n(x)|$ for a.e. $x\in\Omega$.
		
		\item There is a function $g\in L^{p(x)}(\Omega)$ such that $g_n\to g$ in $L^{p(x)}(\Omega)$.
	\end{itemize}
	
\noindent Therefore $(f_n)_{n\geq 1}\subset L^{p(x)}(\Omega),\ f\in L^{p(x)}(\Omega)$ and $f_n\to f$ in $L^{p(x)}(\Omega)$.	
\end{theorem}

\begin{proof} From Theorem \cite[Theorem 2.58]{cruz2013variable} we know that convergence in the norm of $L^{p(x)}(\Omega)$ is equivalent to the convergence in the modular $\rho_{p(x)}$, thanks to the fact that $p^+<\infty$. We will use this throughout the proof.
	
\noindent From $|f_n(x)|\leq |g_n(x)|$ for a.e. $x\in\Omega$ we deduce that

\begin{equation}
	 \rho_{p(x)}(f_n)=\int_{\Omega} |f_n(x)|^{p(x)}\ dx\leq \int_{\Omega} |g_n(x)|^{p(x)}\ dx=\rho_{p(x)}(g_n)<\infty,\ \text{for each }n\in\mathbb{N}^*.
\end{equation}

\noindent Taking into account that $f_n$ is a measurable function, we deduce that $(f_n)_{n\geq 1}\subset L^{p(x)}(\Omega)$.

\noindent Consider now, for each $n\geq 1$ the function 

\begin{align*}
	h_n&:=\big | |g_n|^{p(x)}-|g|^{p(x)}\big |+|g|^{p(x)}+|f|^{p(x)}-\big ||f_n|^{p(x)}-|f|^{p(x)} \big |\\
	&\geq |g_n|^{p(x)}+|f|^{p(x)}-\big ||f_n|^{p(x)}-|f|^{p(x)} \big |\\
	&\geq |f_n|^{p(x)}+|f|^{p(x)}-\big ||f_n|^{p(x)}-|f|^{p(x)} \big |\geq 0.
\end{align*}

\noindent Since $f_n\to f$ pointwise a.e. on $\Omega$, we get that 

\begin{equation}
	\liminf_{n\to\infty} h_n=|g|^{p(x)}+|f|^{p(x)}+\liminf_{n\to\infty} \big | |g_n|^{p(x)}-|g|^{p(x)}\big |.
\end{equation}

\noindent On the other hand, since $g_n\to g$ in $L^{p(x)}(\Omega)$, we deduce, using Proposition 2.67 from the book \cite{cruz2013variable} that there is a subsequence $(g_{n_k})_{k\geq 1}$ that converges to $g$ pointwise a.e. on $\Omega$. Therefore 

\begin{equation}
	0\leq \liminf_{n\to\infty} \big | |g_n|^{p(x)}-|g|^{p(x)}\big |\leq \liminf_{k\to\infty} \big | |g_{n_k}|^{p(x)}-|g|^{p(x)}\big |=0,
\end{equation}

\noindent i.e. $\displaystyle\liminf_{n\to\infty} \big | |g_n|^{p(x)}-|g|^{p(x)}\big |=0$. Combining the last two relations results in: $\displaystyle\liminf_{n\to\infty} h_n=|g|^{p(x)}+|f|^{p(x)}$. Using now \textit{Fatou's Lemma}\footnote{See Theorem 2.4.4, page 63 from \cite{Cohn}.} we obtain that

\begin{align*}
	\int_{\Omega}|g|^{p(x)}+|f|^{p(x)}\ dx&=\int_{\Omega} \liminf_{n\to\infty} h_n(x)\ dx\leq \liminf_{n\to\infty} \int_{\Omega}h_n(x)\ dx\\
	&=\int_{\Omega}|g|^{p(x)}+|f|^{p(x)}\ dx+\liminf_{n\to\infty}\int_{\Omega}  \big | |g_n|^{p(x)}-|g|^{p(x)}\big |-\big ||f_n|^{p(x)}-|f|^{p(x)} \big |\ dx\\
\text{Proposition \ref{propomoducont}}\ \ \ &=\int_{\Omega}|g|^{p(x)}+|f|^{p(x)}\ dx+\liminf_{n\to\infty}\int_{\Omega}-\big ||f_n|^{p(x)}-|f|^{p(x)} \big |\ dx\\
&=\int_{\Omega}|g|^{p(x)}+|f|^{p(x)}\ dx-\limsup_{n\to\infty}\int_{\Omega}\big ||f_n|^{p(x)}-|f|^{p(x)} \big |\ dx.
\end{align*}

\noindent Thus 

\begin{equation}
0\leq\displaystyle\liminf_{n\to\infty}\int_{\Omega}\big ||f_n|^{p(x)}-|f|^{p(x)} \big |\ dx\leq\displaystyle\limsup_{n\to\infty}\int_{\Omega}\big ||f_n|^{p(x)}-|f|^{p(x)} \big |\ dx\leq 0,
\end{equation}

\noindent i.e. $\displaystyle\lim_{n\to\infty}\int_{\Omega}\big ||f_n|^{p(x)}-|f|^{p(x)} \big |\ dx=0$. In particular from this relation it follows that $\lim_{n\to\infty} \rho_{p(x)}(f_n)=\rho_{p(x)}(f)$. As a consequence, since $g_n\to g$ in $L^{p(x)}(\Omega)$ we can use Proposition \ref{propomoducont} to get that

\begin{equation}
	\rho_{p(x)}(f)=\lim_{n\to\infty} \rho_{p(x)}(f_n)\leq \lim_{n\to\infty} \rho_{p(x)}(g_n)=\rho_{p(x)}(g)<\infty.
\end{equation}

\noindent Knowing also that $f$ is a measurable function, we deduce that $f\in L^{p(x)}(\Omega)$. Let us define now for each $n\geq 1$ the function 

\begin{align*}
	\tilde{h}_n&=2^{p^+-1}\left [\big ||f_n|^{p(x)}-|f|^{p(x)}\big | +2 |f|^{p(x)}\right ]-|f_n-f|^{p(x)}\\
	&\geq 2^{p^+-1}\big [ |f_n|^{p(x)}+|f|^{p(x)}\big ]-|f_n-f|^{p(x)}\\
\text{Lemma \ref{propoabel}}\ \ \ 	&\geq 0.
\end{align*}

\noindent It is easy to see that $\tilde{h}_n\to 2^{p^+}|f|^{p(x)}$ pointwise a.e. on $\Omega$. Using again \textit{Fatou's lemma} we get that

\begin{align*}
	2^{p^+}\int_{\Omega} |f|^{p(x)}&=\int_{\Omega} \liminf_{n\to\infty} h_n(x)\ dx\leq \liminf_{n\to\infty} \int_{\Omega}h_n(x)\ dx=\\
	&=\liminf_{n\to\infty}2^{p^+-1} \int_{\Omega} \big ||f_n|^{p(x)}-|f|^{p(x)}\big |\ dx+2^{p^+}\int_{\Omega} |f|^{p(x)}\ dx-\int_{\Omega} |f_n-f|^{p(x)}\ dx\\
	&=2^{p^+}\int_{\Omega} |f|^{p(x)}\ dx+\liminf_{n\to\infty}-\int_{\Omega} |f_n-f|^{p(x)}\ dx\\
	&=2^{p^+}\int_{\Omega} |f|^{p(x)}\ dx-\limsup_{n\to\infty}\int_{\Omega} |f_n-f|^{p(x)}\ dx.
\end{align*}

\noindent Thus 

\begin{equation}
	0\leq \liminf_{n\to\infty}\int_{\Omega} |f_n-f|^{p(x)}\ dx\leq\limsup_{n\to\infty}\int_{\Omega} |f_n-f|^{p(x)}\ dx\leq 0, 
\end{equation}

\noindent i.e. $\displaystyle\lim_{n\to\infty}\int_{\Omega} |f_n-f|^{p(x)}\ dx=0$, which means that $\rho_{p(x)}(f_n-f)\to 0$. But this is equivalent to $f_n\to f$ in $L^{p(x)}(\Omega)$. The proof is now complete.

\end{proof}

\begin{proposition}[\textbf{Poincar\'{e} inequality for variable exponents}]\label{poincare}
	Let $\emptyset\neq\Omega\subset\mathbb{R}^N$ be an open, bounded and connected Lipschitz domain and $p\in\mathcal{P}^{\text{log}}(\Omega)$ a log-H\"{o}lder continuous exponent. Therefore there is a constant $C=C(\Omega,N,p)$ such that 
	
	\begin{equation}
		\Vert u-\langle u\rangle_{\Omega}\Vert_{L^{p(x)}(\Omega)}\leq C\Vert\nabla u\Vert_{L^{p(x)}(\Omega)^N},\ \forall\ u\in W^{1,p(x)}(\Omega),
	\end{equation}
	
	\noindent where $\langle u\rangle_{\Omega}:=\dfrac{1}{|\Omega|}\displaystyle\int_{\Omega} u(x)\ dx$.
\end{proposition}

\begin{proof} This is precisely Theorem 8.2.4 from \cite[page 255]{Hasto}. Just note that any bounded Lipschitz domain is also a John domain, as it is said at \cite[page 237]{Hasto}.
\end{proof}

\begin{lemma}\label{bijuteria} Let $\Omega\subset\mathbb{R}^N$ be an open, bounded and connected Lipschitz domain and $p\in\mathcal{P}^{\text{log}}(\Omega)$ a log-H\"{o}lder continuous exponent with $1<p^-\leq p^+<\infty$.
	
\noindent If we have a sequence $(w_n)_{n\geq 1}\subset W^{1,p(x)}(\Omega)$, a function $w\in L^2(\Omega)$ and some $\mathbf{W}\in L^{p(x)}(\Omega)^N$ such that:

\begin{equation}
	w_n\to w\ \text{in}\ L^2(\Omega)\ \text{and}\ \nabla w_n\weak \mathbf{W}\ \text{in}\ L^{p(x)}(\Omega)^N,
\end{equation}

\noindent then $w\in W^{1,p(x)}(\Omega),\ \nabla w=\mathbf{W}$ a.e. on $\Omega$ and moreover $w_n\weak w$ in $W^{1,p(x)}(\Omega)$.
	
\end{lemma}

\begin{proof} Since $\bigl (\nabla w_n\bigr)_{n\geq 1}$ is weakly convergent, we have from \cite[Proposition 3.3.13 (c)]{papageorgiou2018applied}, that $\bigl (\nabla w_n\bigr)_{n\geq 1}$ is bounded in $L^{p(x)}(\Omega)^N$, which is a reflexive Banach space, from Proposition \ref{propospatii} \textbf{(2)}. Therefore there is a constant $M_1\geq 0$ such that
	
	\begin{equation}\label{woo1}
		\Vert \nabla w_n\Vert_{L^{p(x)}(\Omega)^N}\leq M_1,\ \forall\ n\geq 1.
	\end{equation}
	
\noindent Also, from the fact that $w_n\to w$ in $L^2(\Omega)$ we obtain that the sequence $(w_n)_{n\geq 1}$ is bounded in $L^2(\Omega)$, i.e. there is some constant $M_2\geq 0$ such that

\begin{equation}\label{woo2}
	\Vert w_n\Vert_{L^2(\Omega)}\leq M_2,\ \forall\ n\geq 1.
\end{equation}

\noindent Let $\langle w_n\rangle_{\Omega}:=\dfrac{1}{|\Omega|}\displaystyle\int_{\Omega} w_n(x)\ dx$. From \textit{Cauchy inequality} we get that:

\begin{equation}\label{woo3}
	|\langle w_n\rangle_{\Omega}|=\dfrac{1}{|\Omega|}\left |\displaystyle\int_{\Omega} w_n(x)\ dx \right |\leq \dfrac{1}{|\Omega|}\cdot\sqrt{|\Omega|}\cdot \Vert w_n\Vert_{L^2(\Omega)}=|\Omega|^{-\frac{1}{2}}\Vert w_n\Vert_{L^2(\Omega)},\ \forall\ n\geq 1.
\end{equation}

\noindent Combining all these facts and applying \textit{Poincar\'{e} inequality for variable exponents} -- i.e. Proposition \ref{poincare} -- will lead us to:

\begin{align}\label{woo4}
	\Vert w_n\Vert_{L^{p(x)}(\Omega)}&\leq \Vert w_n-\langle w_n\rangle_{\Omega}\Vert_{L^{p(x)}(\Omega)}+\Vert \langle w_n\rangle_{\Omega}\Vert_{L^{p(x)}(\Omega)}\nonumber\\
	&=\Vert w_n-\langle w_n\rangle_{\Omega}\Vert_{L^{p(x)}(\Omega)}+|\langle w_n\rangle_{\Omega}|\Vert 1\Vert_{L^{p(x)}(\Omega)}\nonumber\\
\eqref{woo3}\ \ \	&\leq \Vert w_n-\langle w_n\rangle_{\Omega}\Vert_{L^{p(x)}(\Omega)}+|\Omega|^{-\frac{1}{2}}\Vert w_n\Vert_{L^2(\Omega)}\Vert 1\Vert_{L^{p(x)}(\Omega)}\nonumber\\
\eqref{woo2}\ \ \ &\leq  \Vert w_n-\langle w_n\rangle_{\Omega}\Vert_{L^{p(x)}(\Omega)}+|\Omega|^{-\frac{1}{2}}M_2\Vert 1\Vert_{L^{p(x)}(\Omega)}\nonumber\\
\text{(Poincar\'{e} ineq.)}\ \ \ &\leq C(\Omega,N,p)\Vert \nabla w_n\Vert_{L^{p(x)}(\Omega)^N}+|\Omega|^{-\frac{1}{2}}M_2\Vert 1\Vert_{L^{p(x)}(\Omega)}\nonumber\\
\eqref{woo1}\ \ \ &\leq C(\Omega,N,p)M_1+|\Omega|^{-\frac{1}{2}}M_2\Vert 1\Vert_{L^{p(x)}(\Omega)}:=M_3<\infty.
\end{align}

\noindent From \eqref{woo1} and \eqref{woo4} we get that:

\begin{equation}\label{woo5}
	\Vert w_n\Vert_{W^{1,p(x)}(\Omega)}=\Vert w_n\Vert_{L^{p(x)}(\Omega)}+\Vert\nabla w_n\Vert_{L^{p(x)}(\Omega)^N}\leq M_3+M_1,\ \forall\ n\geq 1.
\end{equation}

\noindent Relation \eqref{woo5} shows that $(w_n)_{n\geq 1}$ is a bounded sequence in $W^{1,p(x)}(\Omega)$, which from Proposition \ref{propospatii} \textbf{(1)} is a reflexive Banach space. Using now \textit{Eberlein-\v{S}mulian theorem} we get that there is some $\widetilde{w}\in W^{1,p(x)}(\Omega)$ and a subsequence $(w_{n_k})_{k\geq 1}$ such that:

\begin{equation}
	w_{n_k}\weak \widetilde{w}\ \text{in}\ W^{1,p(x)}(\Omega)\hookrightarrow L^{p(x)}(\Omega)\hookrightarrow L^{p^-}(\Omega)\hookrightarrow L^{\min\{p^-,2\}}(\Omega).
\end{equation}

\noindent All these embeddings are continuous and knowing that continuous linear operators preserve weak convergence (see Lemma \ref{lemmaboundedweak}) we deduce that 

\begin{equation}\label{woo6}
w_{n_k}\weak \widetilde{w}\ \text{in}\ L^{\min\{p^-,2\}}(\Omega).
\end{equation}

\noindent But from the statement of the lemma we know that $w_{n_k}\to w$ in $L^2(\Omega)\hookrightarrow L^{\min\{p^-,2\}}(\Omega)$. This means that $w_{n_k}\to w$ in $L^{\min\{p^-,2\}}(\Omega)$. Knowing that strong convergence implies weak convergence, we deduce that:

\begin{equation}\label{woo7}
	w_{n_k}\weak w\ \text{in}\ L^{\min\{p^-,2\}}(\Omega).
\end{equation}

\noindent From \eqref{woo6} and \eqref{woo7}, taking into account that the weak limit of a sequence is unique (if it exists), we conclude that $w=\widetilde{w}\in W^{1,p(x)}(\Omega)$.

\noindent Next, we'll show that $\mathbf{W}=\nabla w$ a.e. on $\Omega$. Indeed, since the gradient map $G:W^{1,p(x)}(\Omega)\to L^{p(x)}(\Omega)^N,\ G(u):=\nabla u$ is a continuous linear operator, and $w_{n_k}\weak w$ in $W^{1,p(x)}(\Omega)$ we infer, using Lemma \ref{lemmaboundedweak}, that 

\begin{equation}
	\nabla w_{n_k}=G(w_{n_k})\weak G(w)=\nabla w\ \text{in}\ L^{p(x)}(\Omega)^N.
\end{equation} 

\noindent But we already have that $\nabla w_{n_k}\weak \mathbf{W}$. Thus $\mathbf{W}=\nabla w$ a.e. on $\Omega$.

\noindent We have proved that $(w_n)_{n\geq 1}$ is a bounded sequence from $W^{1,p(x)}(\Omega)$ and that every weakly convergent subsequence $(w_{n_k})_{k\geq 1}$ must have the weak limit $w$. Thus from Lemma \ref{weakconvergencelemma} we conclude that $w_n\weak w$ in $W^{1,p(x)}(\Omega)$.

\end{proof}

\begin{theorem}[\textbf{Chain Rule for Sobolev spaces with variable exponents}]\label{chainrulevar}
	Let $p:\Omega\to [1,\infty)$ be a measurable exponent and $\Omega\subset\mathbb{R}^N$ an open set with finite measure. If $\eta:\mathbb{R}\to\mathbb{R}$ is a Lipschitz function and $u\in W^{1,p(x)}(\Omega)$, then $\eta\circ u\in W^{1,p(x)}(\Omega)$ and:
	
	\[
	\nabla \big (\eta\circ u\big )(x)=\eta'(u(x))\nabla u(x),\ \text{for a.e.}\ x\in\Omega.
	\]
\end{theorem}

\begin{proof} From \cite[page 604]{kovacik} we have that $u\in W^{1,p(x)}(\Omega)\hookrightarrow W^{1,1}(\Omega)$. Using now the classical chain rule for Sobolev spaces with constant exponent\footnote{See Exercise 11.51 (i), from \cite[page 340]{leonibook}.} we get that $\eta\circ u\in W^{1,1}(\Omega)$ and $\nabla \big (\eta\circ u\big )(x)=\eta'(u(x))\nabla u(x),\ \text{for a.e.}\ x\in\Omega$. Having that $\eta'\circ u\in L^{\infty}(\Omega)$ and $\nabla u\in L^{p(x)}(\Omega)^N$ we also get that $\nabla (\eta\circ u)\in L^{p(x)}(\Omega)^N$. 
	
\noindent Suppose now that $\eta$ is $L$--Lipschitz. Then for any fixed $x_0\in\Omega$: $|\eta(u(x))-\eta(u(x_0))|\leq L|u(x)-u(x_0)|$. From here:

\begin{align*}
	|\eta(u(x))|\leq L|u(x)-u(x_0)|+|\eta (u(x_0))|\leq L|u(x)|+L|u(x_0)|+|\eta (u(x_0))|\in L^{p(x)}(\Omega).
\end{align*}

\noindent Therefore $\eta\circ u\in L^{p(x)}(\Omega)$. In conclusion we have proved that $\eta\circ u\in W^{1,p(x)}(\Omega)$.
\end{proof}

\noindent We say that $u:\Omega\to\mathbb{R},\ \Omega\subseteq\mathbb{R}^N$ open, is an ACL function (absolutely continuous on lines), and write $u\in \text{ACL}(\Omega)$, if $u$ is absolutely continuous on almost every line segment in $\Omega$ parallel to a coordinate axis. For a measurable exponent $p:\Omega\to [1,\infty)$, we say that $u\in\text{ACL}^{p(x)}(\Omega)$ if $u\in \text{ACL}(\Omega)$ and $u,|\nabla u|\in L^{p(x)}(\Omega)$.

\begin{theorem}[\textbf{ACL Characterization of $W^{1,p(x)}(\Omega)$}] Let $\Omega\subseteq \mathbb{R}^N$ an open set and $p:\Omega\to [1,\infty)$ a measurable exponent. Then $\text{ACL}^{p(x)}(\Omega)=W^{1,p(x)}(\Omega)$. In particular any $u\in W^{1,p(x)}(\Omega)$ has classical partial derivative almost everywhere on $\Omega$ and these coincide with the weak partial derivatives of $u$.\footnote{This is precisely Theorem 11.1.12 from \cite[page 345]{Hasto}.}
\end{theorem}

\begin{theorem}[\textbf{General Chain Rule}]\label{thmchainrule} Let  $\emptyset\neq D\subset\mathbb{R}^m$ an open set. Consider $f:D\to\mathbb{R},\ f=f(y_1,y_2,\hdots,y_m)$ a function with $f\in C^1(D)\cap\operatorname{Lip}(D)$. Take also a mapping $\mathbf{u}:\Omega \to D, \ \mathbf{u}=(u_1,u_2,\hdots, u_m)\in W^{1,p(x)}(\Omega)^m$, where $\Omega\subset\mathbb{R}^N$ is also an open set. Then $f\circ \mathbf{u}\in W^{1,p(x)}(\Omega)$ and for each $i\in\overline{1,N}$:
	
	\begin{equation}
		\dfrac{\partial}{\partial x_i}\big (f\circ \mathbf{u}\big )(x)=\nabla f(\mathbf{u}(x))\cdot\dfrac{\partial \mathbf{u}}{\partial x_i}(x)=\displaystyle\sum_{j=1}^m \dfrac{\partial f}{\partial y_j}\big ( \mathbf{u}(x)\big )\dfrac{\partial u_j}{\partial x_i}(x),\ \text{for a.e.}\ x\in\Omega.
	\end{equation}
	
\noindent Moreover if $\mathbf{u}_n\to \mathbf{u}$ in $W^{1,p(x)}(\Omega)^m$, where $\mathbf{u}_n:\Omega\to D$ for each $n\geq 1$, then $f\circ\mathbf{u}_n\to f\circ\mathbf{u}$ in $W^{1,p(x)}(\Omega)$.\footnote{For more general versions of the chain rule for Lipschitz functions, instead of $C^1$, see \cite{ambrosio1990general} and \cite{murat2003chain}.}
	
\end{theorem}

\begin{proof} Since $W^{1,p(x)}(\Omega)=\text{ACL}^{p(x)}(\Omega)\subseteq\text{ACL}(\Omega)$, we deduce that $u_1,u_2,\hdots, u_m$ have classical partial derivatives almost everywhere on $\Omega$. So there is a set $N\subset\Omega$ of null measure such that the classical partial derivatives of $u_1, u_2,\hdots, u_m$ exist on $\Omega\setminus N$.
	
\noindent Fix any $x\in\Omega\setminus N$ and any $i\in\overline{1,N}$. Since $\Omega$ is open we can choose some $\tau>0$ such that $B(x,\tau)=\{\tilde{x}\in\mathbb{R}^N\ |\ |\tilde{x}-x|<\tau\}\subseteq\Omega$, and define: $\mathbf{r}:(-\tau,\tau)\to\mathbb{R}^m,\ \mathbf{r}(h)=\mathbf{u}(x+he_i)-\mathbf{u}(x)$. Knowing that the classical partial derivatives exist at $x$ we deduce that:

\begin{equation}\label{equation34}
	\lim\limits_{h\to 0} \dfrac{\mathbf{r}(h)}{h}=\dfrac{\partial \mathbf{u}}{\partial x_i}(x)\ \text{and}\ 	\lim\limits_{h\to 0} \mathbf{r}(h)=\mathbf{0}.
\end{equation}

\noindent From the definition of differentiability of $f$ in $\mathbf{u}(x)$ we can write:

\begin{equation}\label{equation35}
	f(\mathbf{u}(x+he_i))=f(\mathbf{u}(x)+\mathbf{r}(h))=f(\mathbf{u}(x))+\nabla f(\mathbf{u}(x))\cdot \mathbf{r}(h)+\alpha(\mathbf{r}(h)),\ \text{where}\ \lim\limits_{\mathbf{r}(h)\to \mathbf{0}} \dfrac{\alpha(\mathbf{r}(h))}{|\mathbf{r}(h)|}=0.
\end{equation}
	
\noindent Therefore we may write:

\begin{equation}\label{equation36}
	\dfrac{f(\mathbf{u}(x+he_i))-f(\mathbf{u}(x))}{h}=\nabla f(\mathbf{u}(x))\cdot \dfrac{\mathbf{r}(h)}{h}+\dfrac{\alpha(\mathbf{r}(h))}{|\mathbf{r}(h)|}\dfrac{|\mathbf{r}(h)|}{h}.
\end{equation}

\noindent Making $h\to 0$ in \eqref{equation36}, and using \eqref{equation34} and \eqref{equation35}, we obtain that:

\begin{equation}\label{equation37}
	\dfrac{\partial f\circ\mathbf{u}}{\partial x_i}(x)=\nabla f(\mathbf{u}(x))\cdot\dfrac{\partial \mathbf{u}}{\partial x_i}(x),\ \forall\ x\in\Omega\setminus{N}.
\end{equation}

\noindent Since $f$ is a Lipschitz function on $D$, we get that $f\circ\mathbf{u}$ is also absolutely continuous on lines. So $f\circ\mathbf{u}\in\text{ACL}(\Omega)=W^{1,1}(\Omega)$ and the classical partial derivatives of $f\circ\mathbf{u}$ that are defined almost everywhere on $\Omega$ by \eqref{equation37} are equal to the weak partial derivatives. So we only need to show that $f\circ\mathbf{u}\in L^{p(x)}(\Omega)$ and $\nabla f(\mathbf{u})\cdot\dfrac{\partial \mathbf{u}}{\partial x_i}\in L^{p(x)}(\Omega)$ in order to conclude that $f\circ\mathbf{u}\in\text{ACL}^{p(x)}(\Omega)=W^{1,p(x)}(\Omega)$.

\medskip

\noindent Knowing that $f\in C^1(D)\cap\operatorname{Lip}(D)$ we can find some $L>0$ such that $f$ is an $L$ -- Lipschitz function and $\left |\nabla f(y) \right |\leq L,\ \forall\ y\in D$. Let any $\mathbf{y}_0\in D$. Then for a.e. $x\in\Omega$:

\begin{equation*}
|f(\mathbf{u}(x))-f(\mathbf{y}_0)|\leq L |\mathbf{u}(x)-\mathbf{y}_0|\ \Longrightarrow\ |f(\mathbf{u}(x))|\leq L|\mathbf{u}(x)|+L|\mathbf{y}_0|\in L^{p(x)}(\Omega).
\end{equation*}

\noindent Thus $f\circ\mathbf{u}\in L^{p(x)}(\Omega)$. Furthermore, from \textit{Cauchy inequality for scalar product}, one gets that:

\begin{equation}
\left |\nabla f(\mathbf{u}(x))\cdot\dfrac{\partial \mathbf{u}}{\partial x_i}(x)\right |\leq |\nabla f(\mathbf{u}(x))|\cdot \left |\dfrac{\partial \mathbf{u}}{\partial x_i}(x) \right |\leq L\left |\dfrac{\partial \mathbf{u}}{\partial x_i}(x) \right |\in L^{p(x)}(\Omega).
\end{equation}

\noindent So $\nabla f(\mathbf{u})\cdot\dfrac{\partial \mathbf{u}}{\partial x_i}\in L^{p(x)}(\Omega)$ and the first part of the proof is complete.

\noindent Consider now a sequence $(\mathbf{u}_n)_{n\geq 1}\subset W^{1,p(x)}(\Omega)^m$ with $\mathbf{u}_n:\Omega\to D$, such that $\mathbf{u}_n\to\mathbf{u}$ in $W^{1,p(x)}(\Omega)^m$. We will show that $f\circ\mathbf{u_n}\to f\circ\mathbf{u}$ in $W^{1,p(x)}(\Omega)$. First, notice that:

\begin{equation*}
	\Vert f\circ\mathbf{u}_n-f\circ\mathbf{u}\Vert_{W^{1,p(x)}(\Omega)}=\Vert  f\circ\mathbf{u}_n-f\circ\mathbf{u}\Vert_{L^{p(x)}(\Omega)}+\Vert |\nabla (f\circ \mathbf{u}_n)-\nabla (f\circ\mathbf{u})|\Vert_{L^{p(x)}(\Omega)}.
\end{equation*}

\noindent Then, since $f$ is and $L$--Lipschitz function we can write for a.e. $x\in\Omega$:

\begin{equation}
	|f(\mathbf{u}_n(x))-f(\mathbf{u}(x))|\leq L|\mathbf{u}_n(x)-\mathbf{u}(x)|\ \Rightarrow\ |f(\mathbf{u}_n(x))-f(\mathbf{u}(x))|^{p(x)}\leq L^{p(x)}|\mathbf{u}_n(x)-\mathbf{u}(x)|^{p(x)}.
\end{equation}

\noindent Observe that $L^{p(x)}\leq \max\{L^{p^-},L^{p^+}\}$ for a.e. $x\in\Omega$, and therefore:

\begin{equation}
	\rho_{p(x)}(f\circ\mathbf{u}_n-f\circ\mathbf{u})\leq \max\{L^{p^-},L^{p^+}\}\rho_{L^{p(x)}(\Omega)}\big (|\mathbf{u}_n-\mathbf{u}|\big )\stackrel{n\to\infty}{\longrightarrow} 0.
\end{equation}

\noindent Thus $f\circ\mathbf{u}_n\to f\circ\mathbf{u}$ in $L^{p(x)}(\Omega)$. We have used that $|\mathbf{u}_n-\mathbf{u}|\to 0$ in $W^{1,p(x)}(\Omega)$.

\noindent Next, note that:

\begin{equation*}
	|\nabla(f\circ\mathbf{u}_n)(x)-\nabla (f\circ\mathbf{u})(x)|_{\mathbb{R}^N}\leq \sum_{i=1}^N \left |\dfrac{\partial f\circ\mathbf{u}_n}{\partial x_i}(x)-\dfrac{\partial f\circ\mathbf{u}}{\partial x_i}(x) \right |
\end{equation*}

\noindent Then from \textit{H\"{o}lder inequality} we obtain that:

\begin{equation}\label{equation41}
	\rho_{p(x)}\big (|\nabla(f\circ\mathbf{u}_n)(x)-\nabla (f\circ\mathbf{u})(x)| \big )\leq N^{p^+-1}\displaystyle\sum_{i=1}^N\rho_{p(x)}\left (\left |\dfrac{\partial f\circ\mathbf{u}_n}{\partial x_i}-\dfrac{\partial f\circ\mathbf{u}}{\partial x_i} \right |\right ).
\end{equation}

\noindent We will show that every term in the right hand side of \eqref{equation41} converges to $0$ as $n$ goes to infinity. We have for each $n\geq 1$ that:

\begin{align*}
	&\left |\dfrac{\partial f\circ\mathbf{u}_n}{\partial x_i}(x)-\dfrac{\partial f\circ\mathbf{u}}{\partial x_i}(x) \right | =\left |\nabla f(\mathbf{u}_n(x))\dfrac{\partial \mathbf{u}_n}{\partial x_i}(x)-\nabla f(\mathbf{u}(x))\dfrac{\partial \mathbf{u}}{\partial x_i}(x) \right |\\[3mm]
	=&\left |\nabla f(\mathbf{u}_n(x))\left (\dfrac{\partial \mathbf{u}_n}{\partial x_i}(x)-\dfrac{\partial \mathbf{u}}{\partial x_i}(x)\right ) +\big (\nabla f(\mathbf{u}_n(x))-\nabla f(\mathbf{u}(x))\big )\dfrac{\partial \mathbf{u}}{\partial x_i}(x) \right |\\[3mm]
\text{(Cauchy's ineq.)}\ \leq 	&\left |\nabla f(\mathbf{u}_n(x))\right |\cdot \left |\dfrac{\partial \mathbf{u}_n}{\partial x_i}(x)-\dfrac{\partial \mathbf{u}}{\partial x_i}(x)\right | +\big |\nabla f(\mathbf{u}_n(x))-\nabla f(\mathbf{u}(x))\big |\cdot\left |\dfrac{\partial \mathbf{u}}{\partial x_i}(x) \right |\\[3mm]
\leq & L\left |\dfrac{\partial \mathbf{u}_n}{\partial x_i}(x)-\dfrac{\partial \mathbf{u}}{\partial x_i}(x)\right | +\underbrace{\big |\nabla f(\mathbf{u}_n(x))-\nabla f(\mathbf{u}(x))\big |\cdot\left |\dfrac{\partial \mathbf{u}}{\partial x_i}(x) \right |}_{:=g_n(x)}.
\end{align*}

\noindent So:

\begin{equation}\label{equation42}
\rho_{p(x)}\left (\left |\dfrac{\partial f\circ\mathbf{u}_n}{\partial x_i}-\dfrac{\partial f\circ\mathbf{u}}{\partial x_i} \right | \right )\leq 2^{p^+-1}\max\{L^{p^+},L^{p^-}\}\rho_{p(x)}\left (\left |\dfrac{\partial \mathbf{u}_n}{\partial x_i}(x)-\dfrac{\partial \mathbf{u}}{\partial x_i}(x)\right | \right )+2^{p^+-1} \rho_{p(x)}(g_n).
\end{equation}

\noindent The first term in the right hand side of \eqref{equation42} goes to $0$ as $n\to \infty$ since $\dfrac{\partial\mathbf{u}_n}{\partial x_i}\to \dfrac{\partial\mathbf{u}}{\partial x_i}$ in $L^{p(x)}(\Omega)$.

\noindent For the second term, we will show that $g_n\to 0$ in $L^{p(x)}(\Omega)$ by showing that every subsequence of $g_n$, denoted by $g_{n_k}$ has a further subsequence that converges to $0$ in $L^{p(x)}(\Omega)$. Since $\mathbf{u}_{n_k}\to\mathbf{u}$ in $L^{p(x)}(\Omega)^m$ it follows that there is a further subsequence $\bigl (\mathbf{u}_{n_{k_\ell}}\bigl )_{\ell\geq 1}$ that converges pointwise to $\mathbf{u}$ a.e. on $\Omega$\footnote{See \cite{Jones}, page 234.}. Using the continuity of $\nabla f$ we find that $\big |\nabla f(\mathbf{u}_{n_{k_\ell}}(x))-\nabla f(\mathbf{u}(x))\big |\to 0$ a.e. on $\Omega$. Also, note that for each $n\geq 1$: 

\begin{equation}
	|g_n(x)|\leq 2L\left |\dfrac{\partial \mathbf{u}}{\partial x_i}(x) \right |\in L^{p(x)}(\Omega),
\end{equation}

\noindent Using now the \textit{Lebesgue dominated convergence theorem for variable exponents}\footnote{See Theorem 2.62 from \cite{cruz2013variable}.} we get that $g_{n_{k_\ell}}\to 0$ in $L^{p(x)}(\Omega)$. Thus $g_n\to 0$ in $L^{p(x)}(\Omega)$, which means that $\rho_{p(x)}(g_n)\to 0$ as $n\to\infty$.

\noindent Now, equations \eqref{equation42} and \eqref{equation41}, show that indeed $|\nabla (f\circ\mathbf{u}_n)-\nabla (f\circ\mathbf{u})|\to 0$ in $L^{p(x)}(\Omega)$, and the proof is now complete.

\end{proof}

\begin{theorem}\label{apthplus} Let $\Omega$ be an open and bounded domain, $p:\Omega\to [1,\infty)$ be a measurable exponent and $u\in W^{1,p(x)}(\Omega)$. Set $u^+=\max\{u,0\}$ and $u^{-}=-\min\{u,0\}$. Then $u^{+},u^{-}\in W^{1,p(x)}(\Omega)^+$ and moreover:
	
	\begin{equation}\label{appeqnab}
		\nabla u^{+}=\begin{cases}0,\ \text{a.e. on}\ \{x\in\Omega\ |\ u(x)\leq 0\}\\[3mm]\nabla u,\ \text{a.e. on}\ \{x\in\Omega\ |\ u(x)>0\} \end{cases},\ \nabla u^{-}=\begin{cases}-\nabla u, \ \text{a.e. on}\ \{x\in\Omega\ |\ u(x)<0\}\\[3mm] 0, \ \text{a.e. on}\ \{x\in\Omega\ |\ u(x)\geq 0\} \end{cases}
	\end{equation}
	
	\begin{equation}
		\text{and}\ \nabla |u|=\begin{cases}-\nabla u, \ \text{a.e. on}\ \{x\in\Omega\ |\ u(x)<0\}\\[3mm] 0,\ \text{a.e. on}\ \{x\in\Omega\ |\ u(x)= 0\}\\[3mm]\nabla u,\ \text{a.e. on}\ \{x\in\Omega\ |\ u(x)>0\} \end{cases}.
	\end{equation}
	
	\noindent As a consequence: $u=u^+-u^-,\ |u|=u^++u^-,\ \nabla u=\nabla u^+-\nabla u^{-}$, $u^+u^-=0=\nabla u^+\cdot\nabla u^-=0$ and $\big | \nabla |u|\big |=|\nabla u|$ a.e. on $\Omega$.
	
	\noindent Another important fact is that:
	
	\begin{equation}\label{aineqplusminus}
		\Vert u^+\Vert_{W^{1,p(x)}(\Omega)},\Vert u^-\Vert_{W^{1,p(x)}(\Omega)}\leq \Vert |u|\Vert_{W^{1,p(x)}(\Omega)}=\Vert u\Vert_{W^{1,p(x)}(\Omega)}.
	\end{equation}
\end{theorem}

\begin{proof} Since $p(x)\geq p^{-}\geq 1$ a.e. on $\Omega$ we have that $W^{1,p(x)}(\Omega)\hookrightarrow W^{1,p^-}(\Omega)$. Then it follows from \cite[Theorem 2.2, page 25]{kinnu} that $u^+,u^-\in W^{1,p^{-}}(\Omega)^+$ and \eqref{appeqnab} holds.
	
	\noindent Since $|u|=u^++u^{-}$ we have that $|u^+|,|u^-|\leq |u|$ a.e. on $\Omega$. Since both $u^+$ and $u^-$ are measurable functions we conclude that $u^+,u^-\in L^{p(x)}(\Omega)$. In the same manner $|\nabla u^+|,|\nabla u^{-}|\in L^{p^{-}}(\Omega)$, so they are measurable and from \eqref{appeqnab} we easily see that $|\nabla u^+|,|\nabla u^-|\leq |\nabla u|\in L^{p(x)}(\Omega)$. We get that $|\nabla u^+|,|\nabla u^-|\in L^{p(x)}(\Omega)$. In conclusion $u^+,u^-\in W^{1,p(x)}(\Omega)$ as claimed.
\end{proof}

\subsection*{Bochner-Lebesgue and Bochner-Sobolev spaces}

\begin{theorem}[\textbf{Chain Rule for Bochner-Sobolev Spaces}]\label{thmbochner} Let $p\in (1,\infty]$, $a,b\in\mathbb{R}$ with $a<b$ and $u:(a,b)\times\Omega\to I:=[c,d]\cap\mathbb{R}$, where $-\infty\leq c<d\leq+\infty$ with $u\in W^{1,p}\big ((a,b);L^p(\Omega)\big )$. Here $\Omega\subset\mathbb{R}^N,\ N\geq 1$ is an open, bounded and measurable set. Consider $f:\overline{\Omega}\times I\to\mathbb{R}$ a continuous function for which the partial derivative with respect to the second argument exists on $\Omega\times I$, i.e. $\dfrac{\partial f}{\partial s}:\Omega\times I\to\mathbb{R}$, and has the property that there is some $L>0$ with
	
	\begin{equation}\label{cocos1}
		\left |\dfrac{\partial f}{\partial s}(x,s) \right |\leq L,\ \forall\ (x,s)\in\Omega\times I.
	\end{equation}
	
\noindent We define $v:(a,b)\times\Omega\to\mathbb{R},\ v(t,x):=f(x,u(t,x))$ for a.e. $(t,x)\in (a,b)\times\Omega$. Therefore:

\begin{enumerate}
	\item[\textbf{(1)}] $v\in W^{1,p}\big ((a,b);L^p(\Omega)\big )$ and for a.e. $t\in (a,b)$
	
	\begin{equation}\label{chainbochner}
		\dfrac{\partial v}{\partial t}(t,\cdot)=\dfrac{\partial f}{\partial s}(\cdot,u(t,\cdot))\cdot\dfrac{\partial u}{\partial t}(t,\cdot).
	\end{equation}
	
	\item[\textbf{(2)}] For $p\in (1,\infty)$, with the extra assumption that $\dfrac{\partial f}{\partial s}\in C(\Omega\times I)$, we have that if $u_n\to u$ in $W^{1,p}\big ((a,b);L^p(\Omega)\big )$, where $u_n,u:(a,b)\times\Omega\to I$, and $v_n:(a,b)\times\Omega\to\mathbb{R},\ v_n(t,x):=f(x,u_n(t,x))$ for a.e. $(t,x)\in (a,b)\times\Omega$ then $v_n\to v$ in $W^{1,p}\big ((a,b);L^p(\Omega)\big )$.
	
	\item[\textbf{(3)}] For $p=\infty$, with the extra assumption that $f\in C^1(\mathbb{R}^{N+1})$, we have that  if $u_n\to u$ in $W^{1,\infty}\big ((a,b);L^{\infty}(\Omega)\big )$, where $u_n,u:(a,b)\times\Omega\to I$, and $v_n:(a,b)\times\Omega\to\mathbb{R},\ v_n(t,x):=f(x,u_n(t,x))$ for a.e. $(t,x)\in (a,b)\times\Omega$ then $v_n\to v$ in $W^{1,\infty}\big ((a,b);L^\infty(\Omega)\big )$.
\end{enumerate}

\end{theorem}

\begin{proof} First, using the scalar \textit{Mean value theorem} we get from \eqref{cocos1} that:
	
	\begin{equation}\label{cocos3}
		|f(x,s_1)-f(x,s_2)|\leq L|s_1-s_2|,\ \forall\ x\in\Omega,\ \forall\ s_1,s_2\in [c,d]\cap\mathbb{R}.
	\end{equation}

\noindent Consider the following extension of $f$ given by $\overline{f}:\overline{\Omega}\times\mathbb{R}\to\mathbb{R},\ \overline{f}(x,s)=\begin{cases}f(x,c), & s\leq c\\ f(x,s),& s\in (c,d)\\ f(x,d), & s\geq d \end{cases}$.

\noindent It is straightforward to check that $\overline{f}\in C\bigl(\overline{\Omega}\times\mathbb{R}\bigr)$ and we also have

\begin{equation}\label{cocos2}
	|\overline{f}(x,s_1)-\overline{f}(x,s_2)|\leq L|s_1-s_2|,\ \forall\ x\in\Omega,\ \forall\ s_1,s_2\in\mathbb{R}.
\end{equation}

\noindent \textbf{(1)} $\bullet$ First we consider the case $p\in (1,\infty)$. We divide the proof into several steps:

\medskip

\noindent\textbf{Step I: The Nemytskii operator $\mathcal{N}_{\overline{f}}:L^p(\Omega)\to L^p(\Omega)$ is $L$ -- Lipschitz.}

\noindent Let any $w\in L^{p}(\Omega)$. We need to show that $\mathcal{N}_{\overline{f}}(w)\in L^{p}(\Omega)$. Indeed, fix any $s_0\in (c,d)$, and write for a.e. $x\in\Omega$

\begin{align*}
	|\mathcal{N}_{\overline{f}}(w)(x)|&=|\overline{f}(x,w(x))|\leq |\overline{f}(x,w(x))-\overline{f}(x,s_0)|+|\underbrace{\overline{f}(x,s_0)}_{=f(x,s_0)}|\leq L|w(x)-s_0|+\Vert f(\cdot,s_0)\Vert_{L^{\infty}(\Omega)}\\
	&\leq L|w(x)|+L|s_0|+\Vert f(\cdot,s_0)\Vert_{L^{\infty}(\Omega)} \in L^p(\Omega).
\end{align*}

\noindent Here we used the fact that $f(\cdot,s_0):\overline{\Omega}\to\mathbb{R}$ is a continuous function on a compact subset of $\mathbb{R}^N$, and hence it is bounded.

\noindent Since $\mathcal{N}_{\overline{f}}(w)$ is a measurable function (being a composition between the continuous function $\overline{f}$ and the measurable mapping $\Omega\ni x\mapsto (x,w(x))\in\mathbb{R}^2$) we deduce that $\mathcal{N}_{\overline{f}}(w)\in L^p(\Omega)$.

\noindent Let now any $w_1,w_2\in L^p(\Omega)$. We know from \eqref{cocos2} that for any $x\in\Omega$ 

\begin{equation*}
	|\mathcal{N}_{\overline{f}}(w_1)(x)-\mathcal{N}_{\overline{f}}(w_2)(x)|=|{\overline{f}}(x,w_1(x))-{\overline{f}}(x,w_2(x))|\leq L|w_1(x)-w_2(x)|.
\end{equation*}

\noindent Therefore

\begin{equation}
	\Vert\mathcal{N}_{\overline{f}}(w_1)-\mathcal{N}_{\overline{f}}(w_2)\Vert_{L^{p}(\Omega)}\leq L\Vert w_1-w_2\Vert_{L^{p}(\Omega)}.
\end{equation}

\medskip

\noindent\textbf{Step II: $v\in W^{1,p}\big ((a,b);L^p(\Omega)\big )$.}

\noindent Note that $u:(a,b)\to L^p(\Omega),\ \mathcal{N}_{\overline{f}}:L^p(\Omega)\to L^{p}(\Omega)$, and $v=\mathcal{N}_f\circ u=\mathcal{N}_{\overline{f}}\circ u:(a,b)\to L^{p}(\Omega)$ is the composition between a Lipschitz mapping and a Sobolev function. Since $L^{p}(\Omega)$ is a reflexive Banach space for $p\in (1,\infty)$ (in particular it has the Radon-Nikodym property) we deduce from \cite[Corollary 3.14]{kreuter2015sobolev} that $v\in W^{1,p}\big ((a,b);L^p(\Omega)\big )$.

\medskip

\noindent\textbf{Step III: Establishing formula \eqref{chainbochner}.}

\noindent From the AC characterization of the Bochner-Sobolev space $W^{1,p}\big ((a,b);L^p(\Omega)\big )$ -- see Theorem \ref{thmACX} -- we know that there is a representative of $v$, still denoted by $v\in \text{AC}_\text{loc}\big ((a,b);L^p(\Omega)\big )$ that has a strong derivative a.e. on $(a,b)$ which equals a.e. on $(a,b)$ the weak derivative of $v$.

\noindent Since $u\in W^{1,p}\big ((a,b);L^p(\Omega)\big )$ we can assume, without losing the generality -- from Theorem \ref{thmACX} --, that $u\in \text{AC}_\text{loc}\big ((a,b);L^p(\Omega)\big )$ and $u$ has a strong derivative a.e. on $(a,b)$.

\noindent Take any sequence of real numbers $(h_n)_{n\geq 1}$ with $h_n\to 0$. Therefore, for a.e. $t\in (a,b)$ the weak derivative of $v$ satisfies

\begin{equation}
	\dfrac{\partial v}{\partial t}(t,\cdot)=\lim\limits_{n\to \infty} \dfrac{v(t+h_n,\cdot)-v(t,\cdot)}{h_n}.
\end{equation}

\noindent On the other hand we will show that for a.e. $t\in (a,b)$:

\begin{equation}
	\lim\limits_{n\to \infty} \dfrac{v(t+h_n,\cdot)-v(t,\cdot)}{h_n}=\lim\limits_{n\to \infty} \underbrace{\dfrac{f(\cdot,u(t+h_n,\cdot))-f(\cdot,u(t,\cdot))}{h_n}}_{:=g_n(\cdot)}=\dfrac{\partial f}{\partial s}(\cdot,u(t,\cdot))\dfrac{\partial u}{\partial t}(t,\cdot)\ \text{in}\ L^p(\Omega).
\end{equation}

\noindent Our approach will be to show that any subsequence $(g_{n_k})_{k\geq 1}$ has a further subsequence $(g_{n_{k_\ell}})_{\ell\geq 1}$ that converges to $\dfrac{\partial f}{\partial s}(\cdot,u(t,\cdot))\dfrac{\partial u}{\partial t}(t,\cdot)\ \text{in}\ L^p(\Omega)$.

\noindent Indeed, since we already know from the definition of the strong derivative of $u:(a,b)\to L^p(\Omega)$, that for any subsequence:

\begin{equation}
	\lim\limits_{k\to \infty} \left\Vert \dfrac{u(t+h_{n_k},\cdot)-u(t,\cdot)}{h_{n_k}}-\dfrac{\partial u}{\partial t}(t,\cdot)\right \Vert_{L^p(\Omega)}=0,
\end{equation}

\noindent i.e. $ \dfrac{u(t+h_{n_k},\cdot)-u(t,\cdot)}{h_{n_k}}\stackrel{k\to\infty}{\longrightarrow} \dfrac{\partial u}{\partial t}(t,\cdot)$ in $L^p(\Omega)$, we deduce from \cite[page 234]{Jones} that on a further subsequence $\dfrac{u(t+h_{n_{k_\ell}},\cdot)-u(t,\cdot)}{h_{n_{k_\ell}}}\stackrel{\ell\to\infty}{\longrightarrow} \dfrac{\partial u}{\partial t}(t,\cdot)$ pointwise a.e. on $\Omega$.

\noindent Since $f$ has partial derivative with respect to the second argument, we deduce that:

\begin{align*}
\lim\limits_{\ell\to \infty} \dfrac{f(\cdot,u(t+h_{n_{k_\ell}},\cdot))-f(\cdot,u(t,\cdot))}{h_{n_{k_\ell}}}&=\lim\limits_{\ell\to \infty} \dfrac{f(\cdot,u(t+h_{n_{k_\ell}},\cdot))-f(\cdot,u(t,\cdot))}{u(t+h_{n_{k_\ell}},\cdot)-u(t,\cdot)}\cdot\dfrac{u(t+h_{n_{k_\ell}},\cdot)-u(t,\cdot)}{h_{n_{k_\ell}}}\\
&=\dfrac{\partial f}{\partial s}(\cdot,u(t,\cdot))\dfrac{\partial u}{\partial t}(t,\cdot),\ \text{pointwise a.e. on}\ \Omega.
\end{align*}

\noindent Also from \eqref{cocos3}, note that $|g_{n_{k_\ell}}(x)|=\left |\dfrac{f(\cdot,u(t+h_{n_{k_\ell}},\cdot))-f(\cdot,u(t,\cdot))}{h_{n_{k_\ell}}} \right | \leq L\left |\dfrac{u(t+h_{n_{k_\ell}},\cdot)-u(t,\cdot)}{h_{n_{k_\ell}}} \right |\in L^p(\Omega)$ for each $\ell\geq 1$. Using now the \textit{General Lebesgue Dominated Convergence Theorem} -- i.e. Theorem \ref{generaldct} -- we get that $g_{n_{k_\ell}}\in L^{p}(\Omega)$ for each $\ell\geq 1$, $\dfrac{\partial f}{\partial s}(\cdot,u(t,\cdot))\dfrac{\partial u}{\partial t}(t,\cdot)\in L^p(\Omega)$ and $g_{n_{k_\ell}}\to \dfrac{\partial f}{\partial s}(\cdot,u(t,\cdot))\dfrac{\partial u}{\partial t}(t,\cdot)\ \text{in}\ L^p(\Omega)$. Therefore $\dfrac{\partial v}{\partial t}(t,\cdot )=\dfrac{\partial f}{\partial s}(\cdot,u(t,\cdot))\dfrac{\partial u}{\partial t}(t,\cdot)$ for a.e. $t\in (a,b)$.

\medskip

\noindent $\bullet$ For $p=\infty$, note that $u\in W^{1,\infty}\big ((a,b);L^{\infty}(\Omega)\big )$ from where

\begin{equation} \begin{cases} u\in L^{\infty}\big ((a,b);L^{\infty}(\Omega)\big )\simeq L^{\infty}((a,b)\times\Omega)\subset L^2((a,b)\times\Omega)\simeq L^2\big ((a,b);L^2(\Omega)\big )\\ \dfrac{\partial u}{\partial t}\in L^{\infty}\big ((a,b);L^{\infty}(\Omega)\big )\simeq L^{\infty}((a,b)\times\Omega)\subset L^2((a,b)\times\Omega)\simeq L^2\big ((a,b);L^2(\Omega)\big )\end{cases}.
\end{equation}

\noindent This shows that $u\in  W^{1,2}\big ((a,b);L^{2}(\Omega)\big )$. Therefore, using the previous case we deduce that $v\in  W^{1,2}\big ((a,b);L^{\infty}(\Omega)\big )$ and formula \eqref{chainbochner} holds. We only need to check that $\begin{cases} v\in L^{\infty}\big ((a,b);L^{\infty}(\Omega)\big )\\ \dfrac{\partial v}{\partial t}\in L^{\infty}\big ((a,b);L^{\infty}(\Omega)\big )\end{cases}$. Indeed, if we fix some $s_0\in I$, then for a.e. $(t,x)\in (a,b)\times\Omega$ we can write that

\begin{align*}
 |v(t,x)|&=|f(x,u(t,x))|\leq |f(x,u(t,x))-f(x,s_0)|+|f(x,s_0)|\leq L|u(t,x)-s_0|+\Vert f(\cdot,s_0)\Vert_{L^{\infty}(\Omega)}\\
 &\leq L|u(t,x)|+L|s_0|+\Vert f(\cdot,s_0)\Vert_{L^{\infty}(\Omega)}\\
 &\leq L\Vert u\Vert_{L^{\infty}((a,b)\times\Omega)}+L|s_0|+\Vert f(\cdot,s_0)\Vert_{L^{\infty}(\Omega)}\ \Rightarrow\ v\in L^{\infty}\big ((a,b);L^{\infty}(\Omega)\big ).
\end{align*}

\noindent Note that for a.e. $(t,x)\in (a,b)\times\Omega$

\begin{equation}
	\left |\dfrac{\partial v}{\partial t}(t,x) \right |=\left | \dfrac{\partial f}{\partial s}(x,\underbrace{u(t,x)}_{\in I})\right |\cdot\left |\dfrac{\partial u}{\partial t}(t,x) \right |\leq L\cdot \left |\dfrac{\partial u}{\partial t}(t,x) \right |\leq L\left\Vert\dfrac{\partial u}{\partial t}\right\Vert_{L^{\infty}((a,b)\times\Omega)},
\end{equation}

\noindent i.e. $\dfrac{\partial v}{\partial t}\in L^{\infty}\big ((a,b);L^{\infty}(\Omega)\big )$. Thus $v\in W^{1,\infty}\big ((a,b);L^{\infty}(\Omega)\big )$. The proof of \textbf{(1)} is done.

\medskip

\noindent\textbf{(2)} We need to show that $v_n\to v$ in $W^{1,p}\big ((a,b);L^{p}(\Omega)\big )$, i.e. $v_n\to v$ in $L^{p}\big ((a,b),L^p(\Omega)\big )$ and $\dfrac{\partial v_n}{\partial t}\to \dfrac{\partial v}{\partial t}$ in $L^{p}\big ((a,b),L^p(\Omega)\big )$. From \textbf{Step I} we get that:

\begin{align*}
	\Vert v_n-v\Vert_{L^{p}((a,b);L^p(\Omega))}&=\left (\int_{a}^b \Vert v_n(t,\cdot)-v(t,\cdot)\Vert^{p}_{L^{p}(\Omega)} \ dt\right )^{\frac{1}{p}}\\
	&=\left (\int_{a}^b \Vert \mathcal{N}_f(u_n(t,\cdot))-\mathcal{N}_f(u(t,\cdot))\Vert^{p}_{L^{p}(\Omega)} \ dt\right )^{\frac{1}{p}}\\
\eqref{cocos3}\ \ \	&\leq L\left (\int_{a}^b \Vert u_n(t,\cdot)-u(t,\cdot)\Vert^{p}_{L^{p}(\Omega)} \ dt\right )^{\frac{1}{p}}\\
	&=L\Vert u_n-u\Vert_{L^{p}((a,b);L^p(\Omega))}\stackrel{n\to\infty}{\longrightarrow} 0.
\end{align*}

\noindent Using formula \eqref{chainbochner} we also obtain for a.e. $t\in (a,b)$ that

\begin{align}\label{eqdvndv}
	&\left |\dfrac{\partial v_n}{\partial t}(t,\cdot)-\dfrac{\partial v}{\partial t}(t,\cdot)\right |=\left |\dfrac{\partial f}{\partial s}(\cdot, u_n(t,\cdot))\cdot\dfrac{\partial u_n}{\partial t}(t,\cdot)-\dfrac{\partial f}{\partial s}(\cdot, u(t,\cdot))\cdot\dfrac{\partial u}{\partial t}(t,\cdot) \right | \nonumber\\
	=&\left |\dfrac{\partial f}{\partial s}(\cdot, u_n(t,\cdot))\cdot\left (\dfrac{\partial u_n}{\partial t}(t,\cdot)-\dfrac{\partial u}{\partial t}(t,\cdot) \right ) +\left (\dfrac{\partial f}{\partial s}(\cdot, u_n(t,\cdot))-\dfrac{\partial f}{\partial s}(\cdot, u(t,\cdot)) \right )\cdot\dfrac{\partial u}{\partial t}(t,\cdot) \right |\nonumber \\
	\leq & \left |\dfrac{\partial f}{\partial s}(\cdot, u_n(t,\cdot))\right |\cdot\left |\dfrac{\partial u_n}{\partial t}(t,\cdot)-\dfrac{\partial u}{\partial t}(t,\cdot) \right | +\left |\dfrac{\partial f}{\partial s}(\cdot, u_n(t,\cdot))-\dfrac{\partial f}{\partial s}(\cdot, u(t,\cdot)) \right |\cdot\left |\dfrac{\partial u}{\partial t}(t,\cdot) \right |\nonumber \\
	\leq &L\left |\dfrac{\partial u_n}{\partial t}(t,\cdot)-\dfrac{\partial u}{\partial t}(t,\cdot) \right | +\underbrace{\left |\dfrac{\partial f}{\partial s}(\cdot, u_n(t,\cdot))-\dfrac{\partial f}{\partial s}(\cdot, u(t,\cdot)) \right |\cdot\left |\dfrac{\partial u}{\partial t}(t,\cdot) \right |}_{:=g_n(t,\cdot)}.
\end{align}

\noindent Therefore for a.e. $t\in (a,b)$

\begin{equation}
	\left \Vert \dfrac{\partial v_n}{\partial t}(t,\cdot)-\dfrac{\partial v}{\partial t}(t,\cdot)\right \Vert_{L^{p}(\Omega)}\leq L\left \Vert \dfrac{\partial u_n}{\partial t}(t,\cdot)-\dfrac{\partial u}{\partial t}(t,\cdot)\right \Vert_{L^{p}(\Omega)}+\Vert g_n(t,\cdot)\Vert_{L^{p}(\Omega)}.
\end{equation}

\noindent Raising to the power $p$ this relation, and using Lemma \ref{propoabel} \textbf{(1)} gives us that

\begin{equation}
	\left \Vert \dfrac{\partial v_n}{\partial t}(t,\cdot)-\dfrac{\partial v}{\partial t}(t,\cdot)\right \Vert_{L^{p}(\Omega)}^p\leq 2^{p-1}L^p\left \Vert \dfrac{\partial u_n}{\partial t}(t,\cdot)-\dfrac{\partial u}{\partial t}(t,\cdot)\right \Vert_{L^{p}(\Omega)}^p+2^{p-1}\Vert g_n(t,\cdot)\Vert_{L^{p}(\Omega)}^p.
\end{equation}

\noindent Integrating now on $(a,b)$ allows us to write that

\begin{equation}
	\left \Vert \dfrac{\partial v_n}{\partial t}-\dfrac{\partial v}{\partial t}\right \Vert_{L^{p}((a,b);L^{p}(\Omega))}^p\leq  2^{p-1}L^p\left\Vert \dfrac{\partial u_n}{\partial t}-\dfrac{\partial u}{\partial t}\right \Vert_{L^{p}((a,b);L^{p}(\Omega))}^p+2^{p-1}\Vert g_n\Vert_{L^p((a,b)\times\Omega)}^p
\end{equation}

\noindent Raising now to the power $\dfrac{1}{p}\leq 1$ and using Lemma \ref{propoabel} \textbf{(2)} yields

\begin{equation}
	\left \Vert \dfrac{\partial v_n}{\partial t}-\dfrac{\partial v}{\partial t}\right \Vert_{L^{p}((a,b);L^{p}(\Omega))}\leq  2^{1-\frac{1}{p}}L\left\Vert \dfrac{\partial u_n}{\partial t}-\dfrac{\partial u}{\partial t}\right \Vert_{L^{p}((a,b);L^{p}(\Omega))}+2^{1-\frac{1}{p}}\Vert g_n\Vert_{L^p((a,b)\times\Omega)}.
\end{equation}

\noindent Note that from $u_n\to u$ in $W^{1,p}\big ((a,b);L^p(\Omega)\big )$ we get that $\dfrac{\partial u_n}{\partial t}\to \dfrac{\partial u}{\partial t}$ in $L^{p}\big ((a,b);L^p(\Omega)\big )$, i.e. $\left\Vert \dfrac{\partial u_n}{\partial t}-\dfrac{\partial u}{\partial t}\right \Vert_{L^{p}((a,b);L^{p}(\Omega))}\to 0$ as $n\to\infty$.

\medskip

\noindent Also observe that $g_n:(a,b)\times\Omega\to [0,\infty)$ is a sequence of measurable functions for which we have $|g_n(t,x)|\leq 2L\left |\dfrac{\partial u}{\partial t}\right |\in L^p\big ((a,b)\times\Omega)$. We will show that $g_n\to 0$ in $L^{p}\big ((a,b)\times\Omega\big )$ by showing that for any subsequence $(g_{n_k})_{k\geq 1}$ we can find a further subsequence $(g_{n_{k_{\ell}}})_{\ell\geq 1}$ such that $g_{n_{k_{\ell}}}\to 0$ in $L^{p}\big ((a,b)\times\Omega\big )$ as $\ell\to\infty$.

\noindent Indeed, since $u_{n_{k}}\to u$ in $L^{p}\big ((a,b);L^p(\Omega)\big )\simeq L^{p}\big ((a,b)\times\Omega\big )$ we get that there is a further subsequence such that $u_{n_{k_{\ell}}}\to u$ pointwise a.e. on $(a,b)\times\Omega$. Using now the fact that $\dfrac{\partial f}{\partial s}\in C(\Omega\times\mathbb{R})$ we deduce that $g_{n_{k_\ell}}\to 0$ pointwise a.e. on $(a,b)\times\Omega$. From the \textit{Lebesgue dominated convergence theorem} we conclude that $g_{n_{k_\ell}}\to 0$ in $L^{p}\big ((a,b)\times\Omega\big )$. Therefore $g_n\to 0$ in $L^{p}\big ((a,b)\times \Omega\big )$. So $\dfrac{\partial v_n}{\partial t}\to \dfrac{\partial v}{\partial t}$ in $L^{p}\big ((a,b);L^p(\Omega)\big )$ and $v_n\to v$ in $L^{p}\big ((a,b);L^p(\Omega)\big )$, i.e. $v_n\to v$ in $W^{1,p}\big ((a,b);L^p(\Omega)\big )$.

\medskip

\noindent\textbf{(3)} For $p=+\infty$ the proof is pretty much similar. We need to show that $v_n\to v$ in $W^{1,\infty}\big ((a,b);L^{\infty}(\Omega)\big )$, i.e. $v_n\to v$ in $L^{\infty}\big ((a,b),L^\infty(\Omega)\big )$ and $\dfrac{\partial v_n}{\partial t}\to \dfrac{\partial v}{\partial t}$ in $L^{\infty}\big ((a,b),L^\infty(\Omega)\big )$. From \textbf{Step I} we get that:

\begin{align*}
	\Vert v_n-v\Vert_{L^{\infty}((a,b);L^\infty(\Omega))}&=\underset{t\in (a,b)}{\operatorname{ess\ sup}}\ \Vert v_n(t,\cdot)-v(t,\cdot)\Vert_{L^{\infty}(\Omega)} \\
	&= \underset{t\in (a,b)}{\operatorname{ess\ sup}}\ \Vert \mathcal{N}_f(u_n(t,\cdot))-\mathcal{N}_f(u(t,\cdot))\Vert_{L^{\infty}(\Omega)} \\
\eqref{cocos3}\ \ \	&\leq L\cdot \underset{t\in (a,b)}{\operatorname{ess\ sup}}\ \Vert u_n(t,\cdot)-u(t,\cdot)\Vert_{L^{\infty}(\Omega)}\\
	&=L\Vert u_n-u\Vert_{L^{\infty}((a,b);L^\infty(\Omega))}\stackrel{n\to\infty}{\longrightarrow} 0.
\end{align*}

\noindent Recall that in \eqref{eqdvndv} we have proved for a.e. $(t,x)\in (a,b)\times\Omega$ that

\begin{equation}
	\left |\dfrac{\partial v_n}{\partial t}(t,x)-\dfrac{\partial v}{\partial t}(t,x)\right |\leq L\left |\dfrac{\partial u_n}{\partial t}(t,x)-\dfrac{\partial u}{\partial t}(t,x) \right | +\underbrace{\left |\dfrac{\partial f}{\partial s}(x, u_n(t,x))-\dfrac{\partial f}{\partial s}(x, u(t,x)) \right |\cdot\left |\dfrac{\partial u}{\partial t}(t,x) \right |}_{:=g_n(t,x)}.
\end{equation}

\noindent Therefore

\begin{equation}
	\left \Vert \dfrac{\partial v_n}{\partial t}-\dfrac{\partial v}{\partial t}\right \Vert_{L^{\infty}((a,b)\times\Omega)}\leq L\left \Vert \dfrac{\partial u_n}{\partial t}-\dfrac{\partial u}{\partial t}\right \Vert_{L^{\infty}((a,b)\times\Omega)}+\Vert g_n\Vert_{L^{\infty}((a,b)\times\Omega)}.
\end{equation}

\noindent We have that $\dfrac{\partial u}{\partial t}\in L^{\infty}\big ((a,b);L^\infty(\Omega)\big )$ and $\dfrac{\partial u_n}{\partial t}\to\dfrac{\partial u}{\partial t}$ in $L^{\infty}\big ((a,b);L^\infty(\Omega)\big )$. Observe that $g_n:(a,b)\times\Omega\to [0,\infty)$ is a sequence of measurable functions for which we have $|g_n(t,x)|\leq 2L\left |\dfrac{\partial u}{\partial t}\right |\in L^\infty\big ((a,b)\times\Omega)$. We only need to show that $g_n\to 0$ in $L^{\infty}\big ((a,b)\times\Omega\big )$.

\noindent Since $u_n\to u$ in $L^{\infty}\big ((a,b);L^p(\Omega)\big )\simeq L^{\infty}\big ((a,b)\times\Omega\big )$ we get that for any $1>\tau>0$ there is an $n_{\tau}\geq 1$ such that for any $n\geq n_{\tau}$ we have that $\Vert u_n-u\Vert_{L^{\infty}((a,b)\times\Omega)}<\tau$, i.e. for a.e. $(t,x)\in (a,b)\times\Omega$ the inequality $|u_n(t,x)-u(t,x)|\leq \tau$ holds for any $n\geq n_{\tau}$.

\noindent Now using the extra assumption that $f\in C^1(\mathbb{R}^{N+1})$ we get that $\dfrac{\partial f}{\partial s}\in C\big (\overline{\Omega}\times [-\Vert u\Vert_{L^{\infty}((a,b)\times\Omega))}-1,\Vert u\Vert_{L^{\infty}((a,b)\times\Omega))}+1]\big )$. Because $\overline{\Omega}\times [-\Vert u\Vert_{L^{\infty}((a,b)\times\Omega))}-1,\Vert u\Vert_{L^{\infty}((a,b)\times\Omega))}+1]$ is a compact set, we get from the \textit{Heine-Cantor theorem}\footnote{See Theorem 9.1.5 from \cite{searcoid2007metric}.} that $\dfrac{\partial f}{\partial s}$ is uniformly continuous on $\overline{\Omega}\times [-\Vert u\Vert_{L^{\infty}((a,b)\times\Omega))}-1,\Vert u\Vert_{L^{\infty}((a,b)\times\Omega))}+1]$. Therefore, for any $\epsilon>0$ there is some $\delta_{\epsilon}>0$ such that for any $(x_1,s_1),(x_2,s_2)\in \overline{\Omega}\times [-\Vert u\Vert_{L^{\infty}((a,b)\times\Omega))}-1,\Vert u\Vert_{L^{\infty}((a,b)\times\Omega))}+1]$ with $|x_1-x_2|+|s_1-s_2|\leq \delta_{\varepsilon}$ we have that $\left |\dfrac{\partial f}{\partial s}(x_1,s_1)-\dfrac{\partial f}{\partial s}(x_2,s_2) \right |<\epsilon$.

\noindent Choosing $\tau\leq \delta_{\epsilon}$, then observing that $u_n(t,x),u(t,x)\in [-\Vert u\Vert_{L^{\infty}((a,b)\times\Omega))}-1,\Vert u\Vert_{L^{\infty}((a,b)\times\Omega))}+1]$ for a.e. $(t,x)\in (a,b)\times\Omega$ and $|u_n(t,x)-u(t,x)|\leq\tau\leq\delta_{\epsilon}$ lead us to $\left |\dfrac{\partial f}{\partial s}(x,u_n(t,x))-\dfrac{\partial f}{\partial s}(x,u(t,x)) \right |\leq \epsilon$, for each $n\geq n_{\tau}$, and for a.e. $(t,x)\in (a,b)\times\Omega$. Thus for any $\epsilon>0$, we have for each $n\geq n_{\tau}$ (which depends only on $\epsilon$) that $\Vert g_n\Vert_{L^{\infty}((a,b)\times\Omega)}\leq \epsilon  \left \Vert\dfrac{\partial u}{\partial t}\right \Vert_{L^{\infty}((a,b)\times\Omega)}$. 

\noindent Therefore $g_n\to 0$ in $L^{\infty}\big ((a,b)\times \Omega\big )$. So $\dfrac{\partial v_n}{\partial t}\to \dfrac{\partial v}{\partial t}$ in $L^{\infty}\big ((a,b);L^\infty(\Omega)\big )$ and $v_n\to v$ in $L^{\infty}\big ((a,b);L^\infty(\Omega)\big )$, i.e. $v_n\to v$ in $W^{1,\infty}\big ((a,b);L^\infty(\Omega)\big )$.
\end{proof}

\begin{remark}\label{rembochner} If $I=[c,d]$ is a bounded interval and $f:\overline{\Omega}\times I\to\mathbb{R}$ has a $C^1$ extension on $\mathbb{R}^{N+1}$ then all the hypotheses of Theorem \ref{thmbochner} are fullfilled. This is true because $\dfrac{\partial f}{\partial s}\in C\bigl(\overline{\Omega}\times [c,d]\bigr )$, and since $\overline{\Omega}\times [c,d]$ is a compact set, we get that $\dfrac{\partial f}{\partial s}$ is bounded, as needed.
\end{remark}

\begin{theorem}\label{athminmax} If $u\in H^1\big ((a,b);L^{p(x)}(\Omega)\big )$, where $p:\Omega\to (1,\infty)$ is a measurable exponent with $1<p^{-}\leq p^+<\infty$, then $u^+,u^{-}\in H^1\big ((a,b);L^{p(x)}(\Omega)\big )$ and the following formulas hold for a.e. $t\in (a,b)$:\footnote{This is a straightforward adaptation of \cite[Corollary 5.18]{kreuter2015sobolev} since $L^{p(x)}(\Omega)$ is a reflexive Banach space -- \cite[Corollary 2.7, page 600]{kovacik} -- which is also $\sigma$-Dedekind complete. For a good understanding it is recommended to read the entire chapter 5 of \cite{kreuter2015sobolev}.}
	
	\begin{equation}
		\dfrac{\partial u^+}{\partial t}(t,\cdot)=\begin{cases}\dfrac{\partial u}{\partial t},\ \text{a.e. on}\ \{x\in\Omega\ |\ u(t,x)\geq 0\}\\[3mm] 0,\ \text{a.e. on}\ \{x\in\Omega\ |\ u(t,x)< 0\}\end{cases}
	\end{equation}
	
	\noindent and:
	
	\begin{equation}
		\dfrac{\partial u^-}{\partial t}(t,\cdot)=\begin{cases}-\dfrac{\partial u}{\partial t},\ \text{a.e. on}\ \{x\in\Omega\ |\ u(t,x)<0\}\\[3mm] 0,\ \text{a.e. on}\ \{x\in\Omega\ |\ u(t,x)\geq 0\}\end{cases}.
	\end{equation}
	
\end{theorem}

	\begin{lemma}[\textbf{Commutation of a bounded linear operator with the Bochner integral}]\label{lemcomute}
	Let \(X\) and \(Y\) be real Banach spaces, \(A\in \mathcal L(X,Y)\) a bounded linear operator, and
	\(f\in L^{p}(I;X)\) for some $p\in [1,\infty)$, where $I\subset\mathbb{R}$ any interval. Then \(Af\in L^{p}(I;Y)\), where
	\[
	(Af)(t):=A(f(t)) \qquad \text{for a.e. } t\in I,
	\]
	and for every \(\alpha,\beta\in I\) with \(\alpha<\beta\) one has $\Vert Af\Vert_{L^p(\alpha,\beta;Y)}\leq \Vert A\Vert_{\mathcal{L}(X,Y)}\Vert f\Vert_{L^{p}(\alpha,\beta;X)}$ and
	\[
	A\!\left(\int_{\alpha}^{\beta} f(t)\,dt\right)
	=
	\int_{\alpha}^{\beta} A(f(t))\,dt .
	\]
\end{lemma}

\begin{proof}
	We divide the proof into several steps.
	
	\medskip
	
	\noindent
	\textbf{Step 1: \(Af\) is strongly measurable and belongs to \(L^{1}(\alpha,\beta;Y)\).}
	
	\noindent Since \(f\in L^{1}(I;X)\), the function \(f\) is strongly measurable.
	Because \(A:X\to Y\) is continuous, the composition \(Af=A\circ f\) is again
	strongly measurable. Moreover, for almost every \(t\in (\alpha,\beta)\),
	\[
	\|A(f(t))\|_{Y}\le \|A\|_{\mathcal L(X,Y)}\,\|f(t)\|_{X}.
	\]
	Hence
	\[
	\int_{\alpha}^{\beta} \|A(f(t))\|^p_{Y}\,dt
	\le
	\|A\|^p_{\mathcal L(X,Y)}
	\int_{\alpha}^{\beta} \|f(t)\|^p_{X}\,dt
	<\infty\ \Longrightarrow\ \begin{cases} Af\in L^{p}(\alpha,\beta;Y)\\ \Vert Af\Vert_{L^p(\alpha,\beta;Y)}\leq \Vert A\Vert_{\mathcal{L}(X,Y)}\Vert f\Vert_{L^{p}(\alpha,\beta;X)}\end{cases}.
	\] 
	
	\medskip
	
	\noindent
	\textbf{Step 2: The identity is true for simple functions.}
	
	\noindent Assume first that \(f\) is an \(X\)-valued simple function on \((\alpha,\beta)\), say $f(t)=\sum_{k=1}^{N} x_{k}\,\mathbf 1_{E_{k}}(t)$, where \(x_{1},x_2,\hdots, x_N\in X\) and \(E_{1},E_2,\hdots,E_N\subset (\alpha,\beta)\) are measurable sets. Then, by the definition of the Bochner integral for simple functions, $\int_{\alpha}^{\beta} f(t)\,dt= \sum_{k=1}^{N} |E_{k}|x_{k}$, where \(|E_{k}|\) denotes the Lebesgue measure of \(E_{k}\). Applying the linear operator \(A\), we obtain
	
	\[
	A\!\left(\int_{\alpha}^{\beta} f(t)\,dt\right)
	=
	A\!\left(\sum_{k=1}^{N} |E_{k}|x_{k}\right)
	=
	\sum_{k=1}^{N} |E_{k}|A(x_{k}).
	\]
	\noindent  On the other hand, $A(f(t))= \sum_{k=1}^{N} A(x_{k})\,\mathbf 1_{E_{k}}(t)$, and hence $\int_{\alpha}^{\beta} A(f(t))\,dt=\sum_{k=1}^{N} |E_{k}|A(x_{k})$. Comparing the last two expressions yields
	
	\[
	A\!\left(\int_{\alpha}^{\beta} f(t)\,dt\right)
	=
	\int_{\alpha}^{\beta} A(f(t))\,dt .
	\]
	
	\medskip
	
	\noindent
	\textbf{Step 3: Approximation of a general Bochner integrable function by simple functions.}
	
	\noindent Since \(f\in L^{1}(\alpha,\beta;X)\), by the definition of Bochner integrability
	there exists a sequence of \(X\)-valued simple functions \((f_{n})_{n\geq 1}\)
	such that $\|f_{n}-f\|_{L^{1}(\alpha,\beta;X)} \longrightarrow 0 \ \text{as } n\to\infty$. Because \(A\) is bounded,
	\[
	\|Af_{n}-Af\|_{L^{1}(\alpha,\beta;Y)}
	=
	\int_{\alpha}^{\beta} \|A(f_{n}(t)-f(t))\|_{Y}\,dt
	\le
	\|A\|_{\mathcal L(X,Y)}\,\|f_{n}-f\|_{L^{1}(\alpha,\beta;X)}\stackrel{n\to\infty}{\longrightarrow} 0.
	\]
	Therefore $Af_{n}\to Af \ \text{in } L^{1}(\alpha,\beta;Y)$.
	\medskip
	
	\noindent
	\textbf{Step 4: Passage to the limit in the integral identity.}
	
	\noindent For each \(n\geq 1\), Step 2 gives $A\!\left(\int_{\alpha}^{\beta} f_{n}(t)\,dt\right) =\int_{\alpha}^{\beta} A(f_{n}(t))\,dt$. We now want to make \(n\to\infty\).
	
	\noindent First,
	\[
	\left\|
	\int_{\alpha}^{\beta} f_{n}(t)\,dt
	-
	\int_{\alpha}^{\beta} f(t)\,dt
	\right\|_{X}
	=
	\left\|
	\int_{\alpha}^{\beta} (f_{n}(t)-f(t))\,dt
	\right\|_{X}
	\le
	\int_{\alpha}^{\beta} \|f_{n}(t)-f(t)\|_{X}\,dt=\|f_{n}-f\|_{L^{1}(\alpha,\beta;X)},
	\]
	hence
	\[
	\int_{\alpha}^{\beta} f_{n}(t)\,dt \longrightarrow \int_{\alpha}^{\beta} f(t)\,dt
	\qquad\text{in } X.
	\]
	By continuity of \(A\),
	\[
	A\!\left(\int_{\alpha}^{\beta} f_{n}(t)\,dt\right)
	\longrightarrow
	A\!\left(\int_{\alpha}^{\beta} f(t)\,dt\right)
	\qquad\text{in } Y.
	\]
	
	\noindent Next,
	\[
	\left\|
	\int_{\alpha}^{\beta} A(f_{n}(t))\,dt
	-
	\int_{\alpha}^{\beta} A(f(t))\,dt
	\right\|_{Y}
	\le
	\int_{\alpha}^{\beta} \|A(f_{n}(t))-A(f(t))\|_{Y}\,dt
	=
	\|Af_{n}-Af\|_{L^{1}(\alpha,\beta;Y)},
	\]
\noindent	and the right-hand side tends to \(0\) from Step 3. Therefore
	\[
	\int_{\alpha}^{\beta} A(f_{n}(t))\,dt
	\longrightarrow
	\int_{\alpha}^{\beta} A(f(t))\,dt
	\qquad\text{in } Y.
	\]
	
	\noindent Passing to the limit in
	\[
	A\!\left(\int_{\alpha}^{\beta} f_{n}(t)\,dt\right)
	=
	\int_{\alpha}^{\beta} A(f_{n}(t))\,dt,
	\]
	\noindent  we obtain
	\[
	A\!\left(\int_{\alpha}^{\beta} f(t)\,dt\right)
	=
	\int_{\alpha}^{\beta} A(f(t))\,dt.
	\]
	
	\noindent This proves the claim.
	\end{proof}

	\begin{theorem}[\textbf{AC Characterization of $W^{1,p}(I;X)$}]\label{thmACX} Let $X$ be a real Banach space, $p\in [1,\infty)$ and $I\subset \mathbb{R}$ any open interval.\footnote{This is a combination between the following result: \cite[Corollary 3.3, Theorem 3.7 and Proposition 3.8]{kreuter2015sobolev}. See also \cite[Proposition 1.1, page 104]{showalter2013monotone}. In \cite[Lemma from the appendix]{komura1967nonlinear} we are attentioned that for general Banach spaces $X$ it is not true that every $AC$ function is strongly differentiable almost everywhere.}

	\begin{enumerate}
		\item[\textbf{(I)}] \noindent If $u\in W^{1,p}(I;X)$, then there is some $\tilde{u}\in \operatorname{AC}_{\operatorname{loc}}\big (I;X\big )$\footnote{By this we mean the functions $v:I\to X$ that are absolutely continuous on any compact interval $[a,b]\subset I$.} for which the strong derivative exists a.e. on $I$ and $\dfrac{d\tilde{u}}{dt}\in L^p(I;X)$ such that:
		
		\begin{itemize}
			\item $u(t)=\tilde{u}(t)$ for a.e. $t\in I$.
			
			\item The weak derivative of $u$ coincides with the strong derivative of $\tilde{u}$ a.e. on $I$, i.e. $u'(t)=\dfrac{d\tilde{u}}{dt}(t)$ for a.e. $t\in I$.
			
			\item There is a set of null measure $N\subset I$ such that:
			
			\[
			u(t_2)=u(t_1)+\int_{t_1}^{t_2} u'(t)\ dt,\ \forall\ t_1,t_2\in I\setminus N.
			\]
			
			\item For any $A\in X^*$ we have that $A\circ u\in\operatorname{AC}_{\operatorname{loc}}\big (I\big )$.\footnote{This is the same thing as saying that $u$ is weakly absolutely continuous.}	
	\end{itemize} 
		\end{enumerate}
		
		\begin{enumerate}
			\item[\textbf{(II)}] If $u\in L^p(I;X)$ and there is a function $v\in L^{p}(I;X)$ such that for some $x_0\in X$:
			
			\[
			u(t)=x_0+\int_{t_0}^t v(\tau)\ d\tau,\ \text{for a.e.}\ t\in I,
			\]
			
			\noindent then $u\in W^{1,p}(I;X)$ and $u'(t)=v(t)$ for a.e. $t\in I$.\footnote{See \cite[Theorem 6.35, page 203]{hunter2014notes} or \cite[Theorem 8.55, page 230]{leonibook}.}

		\end{enumerate}

	\end{theorem}

\begin{theorem}\label{thmmaxpri} Let $X$ and $Y$ be two real Banach spaces and $A:X\to Y$ be a continuous (i.e. bounded) linear operator. If $u\in W^{1,p}((a,b);X)$ then $v:=Au\in W^{1,p}((a,b);Y)$ and $\dfrac{\partial v}{\partial t}(t,\cdot)=A\left (\dfrac{\partial u}{\partial t}(t,\cdot) \right )$ for a.e. $t\in (a,b)$. Here $p\in [1,\infty)$ and $a<b$ are any real numbers.
	
\end{theorem}

\begin{proof} From $u\in W^{1,p}(a,b;X)$ we get that $\begin{cases} u\in L^{p}(a,b;X)\\ u'\in L^p(a,b;X)\end{cases}$. Using Lemma \ref{lemcomute} it follows that $v=Au\in L^p(a,b;Y)$ and $Au'\in L^{p}(a,b;Y)$. 
	
\noindent Applying Theorem \ref{thmACX} \textbf{(I)} we get that there is some $t_0\in (a,b)$ such that $u(t)=u(t_0)+\displaystyle\int_{t_0}^{t} u'(\tau)\ d\tau$, for a.e. $t\in (a,b)$. Therefore, using again Lemma \ref{lemcomute} we obtain:

\begin{align*}
	v(t)=Au(t)&=A\left(u(t_0)+\displaystyle\int_{t_0}^{t} u'(\tau)\ d\tau \right)=Au(t_0)+A\left(\displaystyle\int_{t_0}^{t} u'(\tau)\ d\tau \right)\\
	&=Au(t_0)+\displaystyle\int_{t_0}^{t} Au'(\tau)\ d\tau,\ \text{for a.e.}\ t\in (a,b).
\end{align*}
	
\noindent Since $Au'\in L^{p}(a,b;Y)$, using Theorem \ref{thmACX} \textbf{(II)}, we conclude that $v\in W^{1,p}((a,b);Y)$ and $v'(t)=Au'(t)$ for a.e. $t\in (a,b)$. The proof is complete.

\end{proof}

\begin{proposition}\label{propocompactness1}
	Let $X$ be a separable and reflexive Banach space, and let the sequence of functions $(f_n)_{n\ge1}\subset L^\infty(0,T;X)$ be bounded, i.e. for some $M>0$,
	\[
	\sup_{n\ge1}\|f_n\|_{L^\infty(0,T;X)}\le M.
	\]
	Then there exist a subsequence $(f_{n_k})_{k\ge1}$ and a function $f\in L^\infty(0,T;X)$ such that
	\[
	\|f\|_{L^\infty(0,T;X)}\le M,
	\]
	and
	\begin{equation}
		f_{n_k}\weakstar f
		\qquad\text{in }L^\infty(0,T;X).
	\end{equation}
	Moreover, for every $r\in[1,\infty)$,
	\begin{equation}
		f_{n_k}\weak f
		\qquad\text{in }L^r(0,T;X).
	\end{equation}
	The same subsequence $(f_{n_k})_{k\ge1}$ works simultaneously for every finite $r\in [1,\infty)$.
\end{proposition}

\begin{proof} Because $X$ is a reflexive Banach space, we have that $X$ has the Radon-Nikodym property\footnote{Take a look at Corollary 2.21 from \cite[page 17]{kreuter2015sobolev}.}. Therefore from the \textit{Riesz Representation Theorem for the Lebesgue-Bochner Spaces} -- Theorem 2.2.9 from \cite[page 129]{gasinski2005nonlinear} -- we infer that: 
	
\begin{equation}\label{embedings1}
	\begin{cases} L^{\infty}\bigl(0,T;X\bigr)=L^{\infty}\bigl(0,T;X^{**}\bigr)=\left [L^{1}\bigl(0,T;X^*\bigr) \right]^*\\ L^{r}\bigl(0,T;X\bigr)=L^{r}\bigl(0,T;X^{**}\bigr)=\left [L^{r'}\bigl(0,T;X^*\bigr) \right]^*\end{cases}.  
\end{equation}
	
\noindent Also, since $X$ is a reflexive and separable Banach space we deduce that $X^*$ is also a separable Banach space (see Corollary 3.27 from \cite[page 73]{brezis2011functional}). Therefore, from Proposition 2.2.3 given in \cite[page 127]{gasinski2005nonlinear}, we get that $L^{1}\bigl(0,T;X^*\bigr)$ is also a separable Banach space.

\noindent Applying now Corollary 3.30 from \cite[page 76]{brezis2011functional} we get that there is a subsequence $(f_{n_k})_{k\geq 1}$ of the bounded sequence $(f_n)_{n\geq 1}\subset \left [L^{1}\bigl(0,T;X^*\bigr) \right]^*$ and some function $f\in \left [L^{1}\bigl(0,T;X^*\bigr) \right]^*=L^{\infty}\bigl(0,T;X\bigr)$ such that:

\begin{equation}
	f_{n_k}\weakstar f\ \text{in}\ L^{\infty}\bigl(0,T;X\bigr).
\end{equation}

\noindent Using Proposition 3.13 (iii) from \cite[page 63]{brezis2011functional} we deduce that

\begin{equation}
	\Vert f\Vert_{L^{\infty}(0,T;X)}\leq \liminf_{n\to\infty} \Vert f_n\Vert_{L^{\infty}(0,T;X)}\leq M.
\end{equation}

\noindent By the definition of the weak* convergence we get that for any $\varphi\in L^1\bigl(0,T;X^*\bigr)$:

\begin{equation}
	\lim\limits_{k\to\infty} \int_{0}^T\varphi(t)\bigl(f_{n_k}(t)\bigr)\ dt=\int_{0}^T\varphi(t)\bigl(f(t)\bigr)\ dt.
\end{equation}

\noindent From Proposition 2.2.5 given at page 128 in \cite{gasinski2005nonlinear} we get that $L^{r'}\bigl (0,T;X^*\bigr)\hookrightarrow L^{1}\bigl (0,T;X^*\bigr)$ and $L^{\infty}\bigl(0,T;X\bigr)\hookrightarrow L^{r}\bigl (0,T;X\bigr)$. So, for any, $\varphi\in L^{r'}\bigl (0,T;X^*\bigr)$, taking into account that $(f_n)_{n\geq 1},f\in L^{\infty}\bigl(0,T;X\bigr)\subset L^{r}\bigl (0,T;X\bigr)$, we have that:

\begin{equation}\label{varphil1}
	\lim\limits_{k\to\infty} \int_{0}^T\varphi(t)\bigl(f_{n_k}(t)\bigr)\ dt=\int_{0}^T\varphi(t)\bigl(f(t)\bigr)\ dt.
\end{equation}

\noindent But \eqref{varphil1} is precisely the definition of the weak convergence in $L^r\bigl(0,T;X\bigr)$. So $f_{n_k}\weak f$ in $L^r\bigl(0,T;X\bigr)$. The proof is now complete.

\end{proof}

\begin{lemma}\label{lemmaweakbochner}
	Let $X,Y$ be two real Banach spaces such that $X\hookrightarrow Y$ and $T>0$. Consider also any $p\in [1,\infty)$ and suppose that $f_n\weak f$ in $L^p(0,T;X)$. Therefore $f_n\weak f$ in $L^p(0,T;Y)$.
\end{lemma}

\begin{proof} From Proposition 2.2.5 from \cite[page 128]{gasinski2005nonlinear} we get that $L^p(0,T;X)\hookrightarrow L^p(0,T;Y)$. Therefore there is a continuous injective linear operator $\mathcal{I}:L^p(0,T;X)\to L^p(0,T;Y)$. From $f_n\weak f$ in $L^p(0,T;X)$, using Lemma \ref{lemmaboundedweak}, we get that $\mathcal{I}(f_n)\weak \mathcal{I}(f)$ in $L^p(0,T;Y)$, as needed.
\end{proof}

\begin{lemma}\label{lemmacontlp} Let $X$ be a real Banach space and $b>a$ two real numbers. If $f_n\to f$ in $C\bigl([a,b];X\bigr)$ then for any $r\in [1,\infty]$ we also have that $f_n\to f$ in $L^r\bigl(a,b;X\bigr)$. 
\end{lemma}

\begin{proof} If $r\in [1,\infty)$ just note that:
	
	\begin{align*}
		\Vert f_n-f\Vert_{L^r(a,b;X)}&=\left(\int_{a}^b \Vert f_n(t)-f(t)\Vert_X^r\ dt\right)^{\frac{1}{r}}\\
		&\leq (b-a)^{\frac{1}{r}}\sup_{t\in [a,b]} \Vert  f_n(t)-f(t)\Vert_X\\
		&= (b-a)^{\frac{1}{r}}\Vert f_n-f\Vert_{C([a,b];X)}\stackrel{n\to\infty}{\longrightarrow} 0,
	\end{align*}
	
	\noindent because $f_n\to f$ in $C\bigl([a,b];X\bigr)$. Thus $f_n\to f$ in $L^r\bigl(a,b;X\bigr)$.
	
\noindent If $r=\infty$ the proof is quite similar:

\begin{align*}
	\Vert f_n-f\Vert_{L^\infty(a,b;X)}&=\sup_{t\in [a,b]} \Vert  f_n(t)-f(t)\Vert_X\\
	&=\Vert f_n-f\Vert_{C([a,b];X)}\stackrel{n\to\infty}{\longrightarrow} 0.
\end{align*}
	
\end{proof}

\begin{proposition}\label{inducedbochner}
	Let $X,Y$ be two real Banach spaces, $a<b$ two real numbers and $p\in [1,\infty]$. Consider a bounded linear operator $A:X\to Y$ and then define for each $u\in L^p\bigl(a,b;X\bigr)$ the function $\overline{A}u:(a,b)\to Y$, given by $(\overline{A}u)(t):=A(u(t))$. Then $\overline{A}:L^p\bigl(a,b;X\bigr)\to L^p\bigl(a,b;Y\bigr)$ is a bounded linear operator and moreover the following inequality holds:
	
	\begin{equation}\label{neimportant2}
		\Vert \overline{A}u\Vert_{L^{p}(a,b;Y)}\leq \Vert A\Vert_{\mathcal{L}(X,Y)}\Vert u\Vert_{L^p(a,b;X)},\ \forall\ u\in L^p(a,b;X).
	\end{equation}
\end{proposition}

\begin{proof} Let any $u\in L^p(a,b;X)$. Thus $u$ is a strongly measurable function. Therefore the function $\overline{A}u=A\circ u$ is also a strongly measurable function, because $A$ is continuous.
	
\noindent We have that $\Vert Ax\Vert_{Y}\leq \Vert A\Vert_{\mathcal{L}(X,Y)}\cdot \Vert x\Vert_X,\ \forall\ x\in X$. Therefore for a.e. $t\in (a,b)$ we can write that:

\begin{equation}\label{neimportant1}
	\Vert A(u(t))\Vert_{Y}\leq \Vert A\Vert_{\mathcal{L}(X,Y)}\cdot \Vert u(t)\Vert_X.
\end{equation}

\noindent Now, if $p<\infty$, from \eqref{neimportant1} we get that:

\begin{equation}
	\Vert A(u(t))\Vert_{Y}^p\leq \Vert A\Vert_{\mathcal{L}(X,Y)}^p\cdot \Vert u(t)\Vert_X^p\ \Rightarrow\ \int_{a}^b\Vert A(u(t))\Vert_{Y}^p\ dt\leq \Vert A\Vert_{\mathcal{L}(X,Y)}^p\cdot\int_{a}^b \Vert u(t)\Vert_X^p\ dt.
\end{equation}

\noindent On the other side, if $p=\infty$ we may write from \eqref{neimportant1} that:

\begin{equation}
	\underset{t\in (a,b)}{\operatorname{ess\ sup}}\ \Vert A(u(t))\Vert_{Y}\leq \Vert A\Vert_{\mathcal{L}(X,Y)}\cdot \underset{t\in (a,b)}{\operatorname{ess\ sup}}\ \Vert u(t)\Vert_X.
\end{equation}

 \noindent The last two relations can be writen also as: $\Vert \overline{A}u\Vert_{L^{p}(a,b;Y)}\leq \Vert A\Vert_{\mathcal{L}(X,Y)}\Vert u\Vert_{L^p(a,b;X)}<\infty$. So in both cases $\overline{A}u\in L^p\bigl(a,b;X\bigr)$, which means that $\overline{A}:L^p\bigl(a,b;X\bigr)\to L^p\bigl(a,b;Y\bigr)$ is well-defined. Moreover \eqref{neimportant2} holds. Finally since $A$ is a linear operator it is straightforward to check that $\overline{A}$ is also a linear operator. From \eqref{neimportant2} (proved above) it follows that it is in fact a bounded linear operator.

\end{proof}

\begin{proposition}\label{nonlinearlinfty} Let $X,Y$ be two real Banach spaces and $a<b$ two real numbers. Consider a (nonlinear) operator $F:X\to Y$ that is norm-to-norm continuous and bounded in the sense that it maps bounded sets from $X$ into bounded sets from $Y$. Therefore:
	
	\begin{enumerate}
		\item[\textnormal{\textbf{(1)}}] For any  $u\in L^{\infty}\bigl(a,b;X\bigr)$ we have that $F\circ u\in L^{\infty}\bigl(a,b;Y\bigr)$.
			
		\item[\textnormal{\textbf{(2)}}] If moreover we know that there is a nondecreasing function $f:[0,\infty)\to [0,\infty)$ such that $\Vert F(x)\Vert_Y\leq f\left(\Vert x\Vert_X\right),\ \forall\ x\in X$ then we have that:
			
			\begin{equation}
				\Vert F\circ u\Vert_{L^{\infty}(a,b;Y)}\leq f\left(\Vert u\Vert_{L^{\infty}(a,b;X)}\right),\ \forall\ u\in L^{\infty}(a,b;X).
			\end{equation}

	\end{enumerate}
	
\end{proposition}

\begin{proof}\noindent\textbf{(1)} First note that since $F:X\to Y$ is continuous and $u:(a,b)\to X$ is strongly measurable (being from $L^{\infty}(a,b;X$), then their composition $F\circ u:(a,b)\to Y$ is also strongly measurable\footnote{For a more general result, see \cite[Corollary 1.1.11, page 7]{hytonen2016analysis}.}
	
\noindent From the fact that $u\in L^{\infty}\bigl(a,b;X\bigr)$ we have that $\Vert u(t)\Vert_X\leq \underset{\tau\in (a,b)}{\operatorname{ess\ sup}}\ \Vert u(\tau)\Vert_X:=M\in [0,\infty)$, for a.e. $t\in (a,b)$. On the other hand from the fact that $F$ is bounded, we get that there is a constant $C_M\in [0,\infty)$ such that $\Vert Fx\Vert_Y\leq C_M$ for any $x\in X$ with $\Vert x\Vert_X\leq M$. Therefore, combining these two facts we get that:

\begin{equation}
	\Vert (F\circ u)(t)\Vert_Y=\Vert F(u(t))\Vert_Y\leq C_M,\ \text{for a.e.} \ t\in (a,b).
\end{equation}

\noindent This shows that $F\circ u\in L^{\infty}\bigl(a,b;Y\bigr)$ and $\Vert F\circ u\Vert_{L^{\infty}(a,b;Y)}\leq C_M$.
	
\bigskip

\noindent\textbf{(2)} Just note that for any $u\in L^{\infty}\bigl(a,b;X\bigr)$ we may write that

\begin{equation}
\Vert (F\circ u)(t)\Vert_Y=\Vert F(u(t))\Vert_Y\leq f\left ( \Vert u(t)\Vert_X\right )\leq f\left(\Vert u\Vert_{L^{\infty}(a,b;X)}\right),\ \text{for a.e.} \ t\in (a,b).	
\end{equation}

\noindent Therefore: 

\begin{equation}
\Vert F\circ u\Vert_{L^{\infty}(a,b;Y)}=\underset{t\in (a,b)}{\operatorname{ess\ sup}}\ \Vert (F\circ u)(t)\Vert_Y\leq  f\left(\Vert u\Vert_{L^{\infty}(a,b;X)}\right).
\end{equation}

\noindent The proof is now complete.

\end{proof}

\begin{proposition}\label{propocompositionprinc} Let $X,Y$ be two real Banach spaces and $a<b$ two real constants. Consider $F:X\to Y$ to be a continuous (nonlinear) operator that is norm-to-norm continuous and bounded in the sense that is maps bounded subsets of $X$ into bounded subsets of $Y$. If $u_n\to u$ in $L^{\infty}(0,T;X)$, then for every $r\in [1,\infty)$ we have that $F\circ u_n\to F\circ u$ in $L^r(0,T;Y)$.
	
\end{proposition}

\begin{proof} $F\circ u_n,\ n\geq 1$ and $F\circ u$ are strongly measurable, from Proposition \ref{nonlinearlinfty} and moreover $F\circ u_n,\ F\circ u\in L^{\infty}(a,b;Y)\hookrightarrow L^{r}(a,b;Y)$.

\noindent From $u_n\to u$ in $L^{\infty}(a,b;X)$ we get that $\lim\limits_{n\to\infty} \Vert u_n-u\Vert_{L^{\infty}(a,b;X)}=0$, and therefore the sequence $\bigl(\Vert u_n-u\Vert_{L^{\infty}(a,b;X)}\bigr)_{n\geq 1}$ is bounded, i.e. there is some $M\geq 0$ such that $\Vert u_n-u\Vert_{L^{\infty}(a,b;X)}\leq M,\ \forall\ n\geq 1$.

\noindent Now we have for each $n\geq 1$ that:

\begin{equation}
	\Vert u_n\Vert_{L^{\infty}(a,b;X)}=\Vert u_n-u+u\Vert_{L^\infty(a,b;X)}\leq \Vert u_n-u\Vert_{L^{\infty}(a,b;X)}+\Vert u\Vert_{L^{\infty}(a,b;X)}\leq M+\Vert u\Vert_{L^{\infty}(a,b;X)},
\end{equation}

\noindent i.e. the sequence $(u_n)_{n\geq 1}$ is bounded in $L^{\infty}(a,b;X)$. From the fact that $F$ is a bounded (nonlinear) operator we get that $(F\circ u_n)_{n\geq 1}$ is also a bounded sequence from $L^{\infty}(a,b;Y)$. Thus we can find some $M_F\geq 0$ such that 

\begin{equation}\label{ener1}
	\Vert F\circ u_n\Vert_{L^{\infty}(a,b;Y)},\Vert F\circ u\Vert_{L^{\infty}(a,b;Y)}\leq M_F\ \text{for any}\ n\geq 1.
\end{equation}

\noindent Let's denote by $\epsilon_n:=\Vert u_n-u\Vert_{L^{\infty}(a,b;X)}$. We know the fact that $\lim\limits_{n\to\infty}\epsilon_n =0$. For each $n\geq 1$, from the definition of the essential supremum, there is some null-measure set $I_n\subset (a,b)$ such that $\Vert u_n(t)-u(t)\Vert_X\leq \epsilon_n$ for every $t\in (a,b)\setminus I_n$. From here we infer that $\Vert u_n(t)-u(t)\Vert_X\leq\epsilon_n,\ \forall\ n\geq 1,\ \forall\ t\in (a,b)\setminus I$, where $I:=\displaystyle\bigcup_{n=1}^{\infty} I_n$ is also a set of null-measure, because $|I|\leq \displaystyle\sum_{n=1}^{\infty} |I_{n}|=0$. This show that:

\begin{equation}
	u_n(t)\stackrel{n\to\infty}{\longrightarrow} u(t)\ \text{in}\  X,\ \text{for a.e.}\ t\in (a,b).
\end{equation}

\noindent Since $F:X\to Y$ is continuous we get that

\begin{equation}
	(F\circ u_n)(t)\stackrel{n\to\infty}{\longrightarrow} (F\circ u)(t)\ \text{in}\  Y,\ \text{for a.e.}\ t\in (a,b).
\end{equation}

\noindent Thus for any $r\in [1,\infty)$ we can write that:

\begin{equation}\label{ener2}
	\Vert (F\circ u_n)(t)-(F\circ u)(t)\Vert^r_{Y}\stackrel{n\to\infty}{\longrightarrow} 0, \ \text{for a.e.}\ t\in (a,b).
\end{equation}

\noindent Also note that for each $n\geq 1$, using \eqref{ener1}, we get that:

\begin{equation}\label{ener3}
	\Vert (F\circ u_n)(t)-(F\circ u)(t)\Vert_{Y}^r\leq \left (\Vert (F\circ u_n)(t)\Vert_{Y} +\Vert(F\circ u)(t)\Vert_{Y}\right )^r\leq (2M_F)^r\in L^1(a,b).
\end{equation}

\noindent Combining now \eqref{ener2} and \eqref{ener3}, we deduce from \textit{Lebesgue dominated convergence theorem} that 

\begin{equation}
	\lim\limits_{n\to\infty} \Vert F\circ u_n-F\circ u\Vert_{L^r(a,b;Y)}^r=\lim\limits_{n\to\infty} \int_{a}^b \Vert (F\circ u_n)(t)-(F\circ u)(t)\Vert^r_{Y}\ dt=0,
\end{equation}

\noindent and from here $\lim\limits_{n\to\infty} \Vert F\circ u_n-F\circ u\Vert_{L^r(a,b;Y)}=0$, i.e. $F\circ u_n\to F\circ u$ in $L^r(a,b;Y)$, as needed.

\end{proof}

\begin{theorem}[\textbf{Aubin--Lions--Simon compactness lemma, $p=\infty$ case}]\label{lemmaaubin}
	Let $X$, $B$, and $Y$ be real Banach spaces such that $X \stackrel{c}{\hookrightarrow} B \hookrightarrow Y$,
	where the embedding $X\stackrel{c}{\hookrightarrow}B$ is compact and the embedding
	$B\hookrightarrow Y$ is continuous. Let $T>0$ and let $r\in(1,\infty]$. Suppose that we have a sequence
	\[
	(u_n)_{n\geq 1}\subset L^\infty(0,T;X)\ \textit{that is bounded in }L^\infty(0,T;X)
	\]
	and that the sequence formed by the weak time derivatives (which are assumed to exist for each $n\geq 1$) satisfy
	\[
	\left(\frac{\partial u_n}{\partial t}\right)_{n\geq 1}
	\subset L^r(0,T;Y) \ \text{and it is bounded in}\ L^r(0,T;Y).
	\]

	\noindent Then, after eventually choosing representatives of the functions
	$u_n,\ n\geq 1$, we get that $(u_n)_{n\geq 1}\subset C([0,T];B)$ and there exist a subsequence $(u_{n_k})_{k\geq 1}$ and a
	function $u\in C([0,T];B)$ such that:
	\[
	u_{n_k}\longrightarrow u
	\quad\text{in}\ C([0,T];B),\ \text{i.e.}\ \lim_{k\to\infty}\sup_{t\in[0,T]} \|u_{n_k}(t)-u(t)\|_B=0.
	\]
\end{theorem}

\begin{proof} For a complete proof see the classic paper \cite[Section 8, Corollary 4, page 85]{simon1986compact}.
\end{proof}

\begin{lemma}\label{bochnercontinuous} Let $X$ be a real Banach space, $t_2>t_1$ two real numbers and $p\in [1,\infty)$. For each scalar function $a\in L^{p'}(t_1,t_2)$ we define the following operator $\mathcal{T}_a:L^p\bigl(t_1,t_2;X\bigr )\to X$ by the Bochner integral:
	
	\begin{equation}
		\mathcal{T}_a(f)=\int_{t_1}^{t_2} a(t)f(t)\ dt,\ \forall\ f\in L^p\bigl(0,T;X\bigr ).
	\end{equation}
	
	\noindent Therefore $\mathcal{T}_a$ is a well-defined bounded linear operator and moreover 
	
	$$\Vert \mathcal{T}_a(f)\Vert_{X}\leq \Vert a\Vert_{L^{p'}(t_1,t_2)}\cdot\Vert f\Vert_{L^p(t_1,t_2;X)}\ \text{for every}\ f\in L^p\bigl(t_1,t_2;X\bigr ).$$
\end{lemma}

\begin{proof} Fix some $f\in L^p\bigl(t_1,t_2;X\bigr)$. First, since $f$ is a strongly measurable function and $a$ is a scalar measurable function we immediately get that $af:(t_1,t_2)\to X$ is also a strongly measurable function. Moreover $\Vert af\Vert_X$ is also measurable.\footnote{See \cite[page 6]{Arendt2011VectorValued}.} Note that:
	
	\begin{align*}
		\int_{t_1}^{t_2} \Vert a(t)f(t)\Vert_{X}\ dt&=\int_{t_1}^{t_2} |a(t)|\cdot\Vert f(t)\Vert_X\ dt\\
		\text{(H\"{o}lder ineq.)}\ \ \ 	&\leq \Vert a\Vert_{L^{p'}(t_1,t_2)}\cdot \Vert f\Vert_{L^p(t_1,t_2;X)}<\infty.
	\end{align*}
	\noindent This shows that $af$ is indeed Bochner integrable, so $\mathcal{T}_a$ is well-defined. It is obvious to check that $\mathcal{T}_a$ is a linear operator. Finally, from Theorem 1.1.4 given in \cite[page 9]{Arendt2011VectorValued} we obtain that:
	
	\begin{equation}
		\Vert \mathcal{T}_a(f)\Vert_{X}=\left\Vert\int_{t_1}^{t_2} a(t)f(t)\ dt\right\Vert_X\leq \int_{t_1}^{t_2} \Vert a(t)f(t)\Vert_{X}\ dt\leq \Vert a\Vert_{L^{p'}(t_1,t_2)}\cdot \Vert f\Vert_{L^p(t_1,t_2;X)}.
	\end{equation}	
	
	\noindent This shows that $\mathcal{T}_a$ is a bounded linear operator. The proof is now complete.
	
\end{proof}

\begin{lemma}\label{lemmaacdeschisacinchis} If $X$ is a real Banach space then $\textnormal{AC}\bigl((0,T);X\bigr)\simeq\textnormal{AC}\bigl([0,T];X\bigr)$. Precisely if $u:(0,T)\to X,\ u\in \textnormal{AC}\bigl((0,T);X\bigr)$ then the following two limits exist in $X$:
	
	\begin{equation}
		\lim\limits_{t\searrow 0} u(t):=u_0\in X\ \text{and}\ \lim\limits_{t\nearrow T} u(t):=u_T\in X.
	\end{equation}
	
	\noindent Furthermore if we define $\tilde{u}:[0,T]\to X$ by $\widetilde{u}(t)=\begin{cases} u_0,& t=0\\ u(t), & t\in (0,T)\\ u_T, & t=T\end{cases}$, then $\widetilde{u}\in \textnormal{AC}\bigl([0,T];X\bigr)$.
\end{lemma}

\begin{proof} Since $u\in \textnormal{AC}\bigl((0,T);X\bigr)$ we get that $u$ is uniformly continuous on $(0,T)$. Using Theorem 2.85 given in \cite{ClarkGeneralTopology}, there is a unique continuous extension of $u$ on $[0,T]$ that we will denote by $\widetilde{u}$. It remains to show that $\widetilde{u}\in \textnormal{AC}\bigl([0,T];X\bigr)$.
	
\noindent Fix now $\varepsilon>0$. Since $u\in \operatorname{AC}((0,T);X)$, there exists $\delta>0$ such that for every finite family of pairwise non-overlapping intervals $[\alpha_j,\beta_j]\subset(0,T),\ j\in\overline{1,m}$ satisfying $\displaystyle\sum_{j=1}^m (\beta_j-\alpha_j)<\delta$ we have $\displaystyle\sum_{j=1}^m \|u(\beta_j)-u(\alpha_j)\|_X<\varepsilon$.
	
\noindent Now take an arbitrary finite family of pairwise non-overlapping intervals $[a_j,b_j]\subset[0,T],j\in\overline{1,m}$ such that $\sum_{j=1}^m (b_j-a_j)<\delta$. For each $k\in\mathbb{N}$, define

	\[
	a_j^{(k)}:=\max\left\{a_j,\frac1k\right\},
	\qquad
	b_j^{(k)}:=\min\left\{b_j,T-\frac1k\right\}.
	\]

\noindent Since the family of intervals is finite and $a_j<b_j$ for every $j\in\overline{1,m}$, for all sufficiently large $k$ we have $a_j^{(k)}<b_j^{(k)}$ and $[a_j^{(k)},b_j^{(k)}]\subset(0,T)$ for every $j\in\overline{1,m}$. Moreover, shrinking the intervals preserves the pairwise non-overlapping property, and
	\[
	\sum_{j=1}^m\left(b_j^{(k)}-a_j^{(k)}\right)\leq\sum_{j=1}^m (b_j-a_j)<\delta.
	\]
	
\noindent Hence
	\begin{equation}\label{ciuciu1}
	\sum_{j=1}^m \left\| u\bigl(b_j^{(k)}\bigr)-u\bigl(a_j^{(k)}\bigr) \right\|_X <\varepsilon.
	\end{equation}
	
\noindent As $k\to\infty$, $a_j^{(k)}\to a_j$ and $b_j^{(k)}\to b_j$ for every $j\in\overline{1,m}$. Since $\widetilde u\in C([0,T];X)$ and $\widetilde u(t)=u(t)$ for any $t\in(0,T)$, it follows that
	\[
	u\bigl(a_j^{(k)}\bigr)\to \widetilde u(a_j),
	\qquad
	u\bigl(b_j^{(k)}\bigr)\to \widetilde u(b_j).
	\]
	\noindent Therefore,
	\[
	\left\|
	u\bigl(b_j^{(k)}\bigr)
	-
	u\bigl(a_j^{(k)}\bigr)
	\right\|_X
	\longrightarrow
	\|\widetilde u(b_j)-\widetilde u(a_j)\|_X.
	\]
	\noindent Because the sum contains only finitely many terms, we may pass to the limit
	in \eqref{ciuciu1} and obtain
	\[
	\sum_{j=1}^m
	\|\widetilde u(b_j)-\widetilde u(a_j)\|_X
	\leq \varepsilon.
	\]
	
	\noindent Since the finite family of pairwise non-overlapping intervals  was arbitrary, we conclude that $\widetilde u\in \operatorname{AC}([0,T];X)$.
	
\end{proof}

\begin{proposition}\label{propoleoni}
	Let $\bigl(H,(\cdot,\cdot)_H\bigr)$ be a Hilbert space and $f:H\to (-\infty,\infty]$ be a convex, lower semicontinuous and proper (i.e. $f\not\equiv \infty$) functional. Consider now a function $u\in H^1\bigl((0,T);H\bigr)\cap C\bigl([0,T];H\bigr)$ with the property that $u(t)\in \operatorname{dom}\ f:=\bigl\{v\in H\ |\ f(v)<\infty\bigr\}$ for any $t\in [0,T]$.
	
	\noindent If there is a function $h\in L^2\bigl(0,T;H\bigr)$ such that for a.e. $t\in (0,T)$ one has that:
	
	\begin{equation}
		h(t,\cdot)\in\partial f\bigl(u(t)\bigr):=\bigl\{w\in H\ |\ f(v)\geq f(u)+(w,v-u)_H,\ \forall\ v\in H\},
	\end{equation}
	
	\noindent then $f\circ u:[0,T]\to \mathbb{R}$ is absolutely continuous and furthermore for every $w\in\partial f\bigl(u(t)\bigr)$ the following equality takes place:
	
	\begin{equation}
		(f\circ u)'(t)=\bigl(w,u'(t)\bigr)_H,\ \text{for a.e.}\ t\in (0,T).
	\end{equation} 
\end{proposition}

\begin{proof} This is precisely Proposition 59 from the excellent PDE's course \cite{Leoni2014PDEII} that can be found at page 41. I restricted somehow the statement, setting $u\in C\bigl([0,T];L^2(\Omega)\bigr)$ in order to have the function $f\circ u$ defined on all $[0,T]$. Taking into account Lemma \ref{lemmaacdeschisacinchis} the conclusion follows.
\end{proof}

\begin{lemma}\label{lemmaACX}
	Let $X$ be a reflexive real Banach space. If $u\in H^1\bigl((0,T);X\bigr)\cap C\bigl([0,T];X\bigr)$ then $u\in \textnormal{AC}\bigl([0,T];X\bigr)$ and for every $0\leq t_1\leq t_2\leq T$ the following formula holds:
	
	\begin{equation}\label{leoniAC}
		u(t_2)-u(t_1)=\int_{t_1}^{t_2} u'(t)\ dt.
	\end{equation}
\end{lemma}

\begin{proof} From the $AC$ characterization of the Bochner-Sobolev space $H^1\bigl((0,T);X\bigr)$ -- i.e. Theorem \ref{thmACX} \textbf{(I)} --  we get that there is some $\tilde{u}\in AC_{\text{loc}}\bigl((0,T);X\bigr )$ with $u(t)=\tilde{u}(t)$ for a.e. $t\in (0,T)$. Using Theorem 8.38 (i) given in \cite[page 225]{leonibook}, we have that $\tilde{u}:(0,T)\to X$ is a continuous function. Since $u:[0,T]\to X$ is also a continuous function and $u=\tilde{u}$ for a.e. $t\in (0,T)$ we deduce that $\tilde{u}(t)=u(t)$ for any $t\in (0,T)$ (if we remove a null-measure set from $(0,T)$ we get a dense subset of $(0,T)$). This shows that $u\in AC_{\text{loc}}\bigl((0,T);X\bigr )$.
	
\noindent From the same Theorem 8.38 (iii) in \cite{leonibook} we deduce that for any fixed $t_0\in (0,T)$ one has that:

\begin{equation}
	u(t)=u(t_0)+\int_{t_0}^t u'(\tau)\ d\tau,\ \forall\ t\in (0,T).
\end{equation}

\noindent Define the function $v:[0,T]\to X,\ v(t)=u(t_0)+\displaystyle\int_{t_0}^t u'(\tau)\ d\tau$. We show that $v\in \textnormal{AC}\bigl([0,T];X\bigr)\subset C\bigl([0,T];X\bigr)$.

\noindent Indeed, for any finite family of pairwise disjoint intervals $(a_k,b_k)\subset [0,T],\ k\in\overline{1,n}$, we have that:

\begin{align*}
	\sum_{k=1}^n\Vert v(b_k)-v(a_k)\Vert_X&=\sum_{k=1}^n \left \Vert \int_{a_k}^{b_k} u'(\tau)\ d\tau\right\Vert_X\\
\text{\cite[Theorem 8.9]{leonibook}}\ \ \ 	&\leq \sum_{k=1}^n \int_{a_k}^{b_k} \Vert u'(\tau)\Vert_X\ d\tau\\
\text{(Cauchy ineq. for integrals)}\ \ \ &\leq \sum_{k=1}^n (b_k-a_k)^{\frac{1}{2}}\cdot \left (\int_{a_k}^{b_k} \Vert u'(\tau)\Vert_X^2\ d\tau\right )^{\frac{1}{2}}\\
\text{(Cauchy ineq. for sum of products)}\ \ \ &\leq \left (\sum_{k=1}^n b_k-a_k \right )^{\frac{1}{2}}\cdot \left (\sum_{k=1}^n \int_{a_k}^{b_k} \Vert u'(\tau)\Vert_X^2\ d\tau \right )^{\frac{1}{2}}\\
&\leq \left (\sum_{k=1}^n b_k-a_k \right )^{\frac{1}{2}}\cdot\left (\int_{0}^T \Vert u'(\tau)\Vert_X^2\right )^{\frac{1}{2}}\\
&=\left (\sum_{k=1}^n b_k-a_k \right )^{\frac{1}{2}}\cdot\Vert u'\Vert_{L^2(0,T;X)}.
\end{align*} 

\noindent This shows that $v\in \textnormal{AC}\bigl([0,T];X\bigr)$. Having that $u(t)=v(t)$ for any $t\in (0,T)$ and $u,v\in C\bigl([0,T];X\bigr)$, we deduce that $u(t)=v(t)$ for any $t\in [0,T]$. Hence $u\in \textnormal{AC}\bigl([0,T];X\bigr)$. Moreover, for any $0\leq t_1\leq t_2\leq T$ we have that:

\begin{equation}
u(t_2)-u(t_1)=v(t_2)-v(t_1)=\int_{t_1}^{t_2} u'(\tau)\ d\tau,	
\end{equation}

\noindent as required.
	
\end{proof}

\begin{lemma}\label{stronglymeasurableembedding} Let $X$ and $Y$ be two separable real Banach spaces with $X\hookrightarrow Y$, i.e. there is a continuous and injective linear map $I:X\to Y$. Consider $a<b$ two real numbers and $f:(a,b)\to Y$ a Borel measurable function with the property that $f(t)\in I(X)$ for a.e. $t\in (a,b)$. Then there exists a strongly measurable function $\tilde{f}:(a,b)\to X$ with $I(\tilde{f}(t))=f(t)$ for a.e. $t\in (a,b)$.
\end{lemma}

\begin{proof} First note that $I:X\to I(X):=Z$ is a bijective continuous linear operator. So there exists $I^{-1}:Z\to X$. Since $X$, $Y$ are separable Banach spaces, they are Polish spaces. Taking into account that $I:X\to Y$ is injective and continuous, we may apply \textit{Lusin-Souslin theorem} -- see the proof in \cite[Theorem 8.3.7, page 260]{Cohn} -- to deduce that $I(X)=Z$ is a Borel subset of $Y$ and the inverse $I^{-1}:Z\to X$ is Borel measurable.
	
\noindent Since $f:(a,b)\to Y$ is Borel measurable and $Z\subset Y$ is a Borel set we deduce that: $E:=f^{-1}(Z)\subset (a,b)$ is a Borel set. From the statement we know that $f(t)\in Z$ for a.e. $t\in (a,b)$, i.e. $t\in f^{-1}(Z)=E$ for a.e. $t\in (a,b)$. Thus $|(a,b)\setminus E|=0$.

\noindent Now define the function:

\begin{equation}
	\tilde{f}:(a,b)\to X,\ \tilde{f}(t)=\begin{cases} I^{-1}(f(t)), & t\in E\\[3mm] 0_X, & t\in (a,b)\setminus E\end{cases}.
\end{equation}

\noindent $\bullet$ Observe that for any $t\in E=f^{-1}(Z)$ we have that $f(t)\in Z=I(X)$, so it makes sense to define $\tilde{f}(t)=I^{-1}(f(t))\in X$. Moreover from here we deduce that $(I\circ \tilde{f})(t)=f(t)$ for any $t\in E$, so for a.e. $t\in (a,b)$.

\medskip

\noindent $\bullet$ We show now that $\tilde{f}:(a,b)\to X$ is a Borel measurable function. Let $U\subset X$ be an open set. We need to prove that $\tilde{f}^{-1}(U)$ is a Borel subset of $(a,b)$. We can encounter one of the following two situations:

\begin{enumerate}
	\item[$\blacksquare$] If $0_X\notin U$ then 
	
	\begin{align*}
		\tilde{f}^{-1}(U)&=\{t\in (a,b)\ |\ \tilde{f}(t)\in U\}=\{t\in E\ |\ \tilde{f}(t)\in U\}\\
		&=E\cap f^{-1}\bigl(I(U)\bigr),
	\end{align*} 
	
	\noindent which is a Borel set, being the intersection of two Borel sets. Here we used the fact that $I^{-1}$ is Borel measurable, and since $U\subset X$ is open (hence a Borel set) we get that $I(U)=(I^{-1})^{-1}(U)$ is a Borel set from $Y$. But, since $f:(a,b)\to Y$ is Borel measurable, we deduce that $f^{-1}\bigl(I(U)\bigr)$ is a Borel set.
	
	\medskip
	
	\item[$\blacksquare$] If $0_X\in U$ then 
	
	\begin{align*}
		\tilde{f}^{-1}(U)&=\{t\in (a,b)\ |\ \tilde{f}(t)\in U\}=\{t\in E\ |\ \tilde{f}(t)\in U\}\cup \{t\in (a,b)\setminus E\ |\ \tilde{f}(t)\in U\} \\
		&=\left [E\cap f^{-1}\bigl(I(U)\bigr)\right ]\cup \bigl[(a,b)\setminus E\bigr],
	\end{align*} 
	
	\noindent which is again a Borel set, being a reunion of Borel sets. Here we used the fact that $E$ is a Borel subset of $(a,b)$ and therefore $(a,b)\setminus E$ is also a Borel set.
\end{enumerate} 

\medskip

\noindent $\bullet$ The last step is to show that the Borel measurable function $\tilde{f}:(a,b)\to X$ is also strongly measurable. Taking into account that $X$ is a separable Banach space, from \textit{Pettis measurability theorem} -- see Corollary 1.1.2 b) in \cite{Arendt2011VectorValued} -- we get that $\tilde{f}$ is strongly measurable iff $\tilde{f}$ is weakly measurable, i.e. for any $\ell\in X^*$ we have that $\ell\circ\tilde{f}:(a,b)\to \mathbb{R}$ is Lebesgue measurable.

\medskip

\noindent Indeed, since $\ell\in X^*$, we get that $\ell:X\to\mathbb{R}$ is continuous, and hence Borel measurable\footnote{See Corollary 2.2 given at page 44 in \cite{folland1999real} for the proof that any continuous  function between two topological spaces is a Borel measurable function.}. Therefore $\ell\circ\tilde{f}:(a,b)\to\mathbb{R}$ is also Borel measurable, being a composition of two Borel measurable functions (for details, see Proposition 2.6.1 in \cite[page 73]{Cohn}). But any Borel measurable function is also a Lebesgue measurable function\footnote{Read the warning given at the end of the proof of Corollary 2.2. from \cite[page 44]{folland1999real}.}, so $\ell\circ\tilde{f}$ is Lebesgue measurable. This shows that $\tilde{f}$ is weakly measurable, and from the above discussion it is in fact strongly measurable, as needed.
	
\end{proof}

\begin{lemma}\label{contl2allt}
	Let $\Omega\subset\mathbb{R}^N$ a bounded measurable set, $a<b$ and $\alpha<\beta$ some real numbers. Assume that $v\in C\bigl([a,b];L^2(\Omega)\bigr)$ and:
	
	\begin{equation}\label{programarermn1}
		\alpha\leq v(t,x)\leq \beta,\ \text{for a.e.}\ (t,x)\in (a,b)\times\Omega.
	\end{equation}
	
	\noindent Then, for \textbf{every} $t\in [a,b]$, we have that $\alpha\leq v(t,x)\leq\beta$ for a.e. $x\in\Omega$.
\end{lemma}

\begin{proof} $\bullet$ Consider the following set $M:=\bigl\{(t,x)\in (a,b)\times\Omega\ |\ v(t,x)<\alpha\ \text{or}\ v(t,x)>\beta\bigr\}$ and for any $t\in (a,b)$ we denote $M_t:=\{x\in\Omega\ |\ (t,x)\in M\}$. Relation \eqref{programarermn1} says that $|M|=0$. Also consider the characteristic function $\chi_M:(a,b)\times\Omega\to\{0,1\}$ given by $\chi_M(t,x)=\begin{cases} 1, & (t,x)\in M\\ 0, & (t,x)\notin M\end{cases}$. Since $\chi_M$ is a measurable positive function, we get from \textit{Tonelli's theorem} that 
	
	\begin{equation}
		\text{the function}\ (a,b)\ni t\mapsto \int_{\Omega} \chi_{M}(t,x)\ dx=|M_t|\geq 0\ \text{is measurable},
	\end{equation}
	
	\noindent and
	
	\begin{equation}
	0=|M|=\int_{(a,b)\times\Omega}\chi_M(t,x)\ dx\ dt=\int_{a}^b\left(\int_{\Omega} \chi_M(t,x)\ dx\right)\ dt=\int_{a}^b \underbrace{|M_t|}_{\geq 0}\ dt\geq 0.
	\end{equation}

	\noindent From here we deduce that there is a measurable set $E\subset (a,b)$, with $|(a,b)\setminus E|=0$ such that $|M_t|=0$ for any $t\in E$. This means that:
	
	\begin{equation}\label{lumanare1}
		\text{for any}\ t\in E:\ \alpha\leq v(t,x)\leq \beta,\ \text{for a.e.}\ x\in \Omega.
	\end{equation}
	
	\noindent $\bullet$ Consider now the following set $K:=\bigl\{u\in L^2(\Omega)\ |\ \alpha\leq u(x)\leq \beta,\ \text{for a.e.}\ x\in\Omega\bigr\}$. Next, we'll show that $K$ is a closed subset of $L^2(\Omega)$ with respect to the topology induced by its norm. Let $(u_n)_{n\geq 1}\subset K$ and $u\in L^2(\Omega)$ with $u_n\to u$ in $L^2(\Omega)$. We need to prove that $u\in K$.
	
	\noindent Since $u_n\to u$ in $L^2(\Omega)$ we get, from \cite[page 234]{Jones}, that there is a subsequence such that $u_{n_k}\to u$ pointwise a.e. on $\Omega$. So there is a null-measure set $\omega_0\subset\Omega$ such that:
	
	\begin{equation}
		\lim\limits_{k\to\infty} u_{n_k}(x)=u(x),\ \forall\ x\in\Omega\setminus\omega_0.
	\end{equation}
	
	\noindent For each $k\geq 1$ we know that $u_{n_k}\in K$, so there is a null-measure set $\omega_k\subset\Omega$ with the property that:
	
	\begin{equation}
		\alpha\leq u_{n_k}(x)\leq \beta,\ \forall\ x\in\Omega\setminus\omega_k.
	\end{equation}
	
	\noindent Denote $\omega:=\bigcup_{k=0}^{\infty} \omega_k$. So $|\omega|\leq\sum_{k=0}^{\infty} |\omega_k|=0$, i.e. $\omega$ is also a null-measure subset of $\Omega$. But for each $x\in\Omega\setminus\omega$ we can write that:
	
	\begin{equation}
		\alpha \leq u_{n_k}(x)\leq \beta,\ \forall\ k\geq 1\ \text{and}\ \lim\limits_{k\to\infty} u_{n_k}(x)=u(x).
	\end{equation}
	
	\noindent Making $k\to \infty$ in the above double inequality gives us that $\alpha\leq u(x)\leq \beta$ for any $x\in\Omega\setminus\omega$, so for a.e. $x\in\Omega$. This shows that $u\in K$, i.e. $K$ is closed.

	\noindent $\bullet$ Coming back to our problem, \eqref{lumanare1} says that $v(t,\cdot)\in K$ for any $t\in E$. We will show that for any $t\in [a,b]$ we have that $v(t,\cdot)\in K$. 
	
	\noindent Indeed, since $[a,b]\setminus E$ is a null-measure set, we get that $E$ is dense in $[a,b]$. Therefore we can find a sequence $(t_n)_{n\geq 1}\subset E$ with $t_n\to t$. Now since $v\in C\bigl([a,b];L^2(\Omega)\bigr)$ we deduce that $\lim\limits_{n\to\infty} \underbrace{v(t_n,\cdot)}_{\in K}=v(t,\cdot)$ in $L^2(\Omega)$, and $K$ is a closed subset of $L^2(\Omega)$, we finally get that $v(t,\cdot)\in K$ for any $t\in [a,b]$. This is precisely our conclusion.
	
\end{proof}

\subsection*{Functional analysis}

\begin{lemma}[\textbf{Bounded linear operators preserve weak convergence}]\label{lemmaboundedweak} Let $X$ and $Y$ be two Banach spaces and $A:X\to Y$ be a bounded linear operator. If $x_n\weak x$ in $X$, then $Ax_n\weak Ax$ in $Y$.
\end{lemma}

\begin{proof} From $x_n\weak x$ in $X$ we get that for each $\phi\in X^*$: $\lim\limits_{n\to\infty} \phi(x_n)=\phi(x)$. We need to prove that for each $\psi\in Y^*$ we have that $\lim\limits_{n\to\infty} \psi(Ax_n)=\psi(Ax)$. But this is true by just setting $\phi=\psi\circ A\in X^*$ (a composition between two continuous linear operators is a continuous linear operator) in the previous relation. Therefore $Ax_n\weak Ax$ in $Y$. 
\end{proof}

\begin{lemma}\label{weakconvergencelemma} Let $X$ be a reflexive Banach space, $(x_n)_{n\geq 1}\subset X$ a bounded sequence and $x\in X$ with the following property: Every weakly convergent subsequence of $(x_n)_{n\geq 1}$ weakly converges to $x$. Then: $x_n\weak x$ in $X$.
\end{lemma}

\begin{proof} Suppose that $x_n\not\weak x$ in $X$. Then there is some $f\in X^*$ such that $f(x_n)\not\to f(x)$ in $\mathbb{R}$. So, there is some $\epsilon>0$ and a subsequence $(x_{n_k})_{k\geq 1}$ such that:
	
	\begin{equation}\label{relatiedecontrazis}
		|f(x_{n_k})-f(x)|\geq \epsilon,\ \forall\ k\geq 1.
	\end{equation}
	
	\noindent Since $(x_{n_k})_{k\geq 1}$ is still a bounded sequence from the reflexive Banach space $X$, we get from \textit{Eberlein-\v{S}mulian theorem} that it has a further subsequence $(x_{n_{k_\ell}})_{\ell\geq 1}$ that is weakly convergent to some element $\tilde{x}\in X$. But from the property given in the statement it follows that $x_{n_{k_\ell}}\weak \tilde{x}=x$. In particular, since $f\in X^*$, we have that
	
	\begin{equation}
		f(x_{n_{k_{\ell}}})\to f(x).
	\end{equation}
	
	\noindent But this relation contradicts \eqref{relatiedecontrazis}. The conclusion follows.
	
\end{proof}

\subsection*{Other results}

\begin{lemma}\label{lemmaliminfsupcontmon}
	Let $(a_n)_{n\geq 1}\subset\mathbb{R}$ be a sequence such that $\displaystyle\liminf_{n\to\infty}a_n\in\mathbb{R}$, and let $\varphi:\mathbb{R}\to\mathbb{R}$ be a continuous and nondecreasing function. Then:
	\begin{equation}
		\liminf_{n\to\infty}\varphi(a_n)=\varphi\left(\liminf_{n\to\infty}a_n\right).
	\end{equation}
\end{lemma}

\begin{proof}
	Set $\ell:=\displaystyle\liminf_{n\to\infty}a_n\in\mathbb{R}$. We first prove that $\displaystyle\liminf_{n\to\infty}\varphi(a_n)\geq\varphi(\ell)$. Let $\varepsilon>0$. By the definition of the limit inferior, there exists $n_\varepsilon\in\mathbb{N}$ such that $a_n\geq \ell-\varepsilon$ for every $n\geq n_\varepsilon$. Since $\varphi$ is nondecreasing, it follows that $\varphi(a_n)\geq\varphi(\ell-\varepsilon)$ for every $n\geq n_\varepsilon$. Consequently,
	\[
	\liminf_{n\to\infty}\varphi(a_n)\geq\varphi(\ell-\varepsilon).
	\]
	Letting $\varepsilon\to0^+$ and using the continuity of $\varphi$, we obtain $\displaystyle\liminf_{n\to\infty}\varphi(a_n)\geq\varphi(\ell)$.
	
	\medskip
	
	\noindent On the other hand, by the definition of the limit inferior, there exists a subsequence $(a_{n_k})_{k\geq 1}$ such that $a_{n_k}\to\ell$ as $k\to\infty$. Since $\varphi$ is continuous, we have $\varphi(a_{n_k})\to\varphi(\ell)$. Moreover, the lower limit of a sequence is not larger than the lower limit of any of its subsequences, and therefore
	\[
	\liminf_{n\to\infty}\varphi(a_n)\leq\liminf_{k\to\infty}\varphi(a_{n_k})=\varphi(\ell).
	\]
	Combining the two inequalities yields
	\[
	\liminf_{n\to\infty}\varphi(a_n)=\varphi(\ell)=\varphi\left(\liminf_{n\to\infty}a_n\right),
	\]
	which proves the assertion.
\end{proof}

\begin{lemma}\label{beppolevi1}
	Let $I\subset\mathbb{R}$ be a measurable set and $(I_n)_{n\geq 1}$ be a sequence of measurable subsets of $I$ with $\lim\limits_{n\to\infty} |I\setminus I_n|=0$. Then, for any $f\in L^1(I)$:
	
	\begin{equation}
		\lim\limits_{n\to\infty} \int_{I_n} f(t)\ dt=\int_I f(t)\ dt.
	\end{equation}
	
	\noindent In other words $\lim\limits_{n\to\infty} \displaystyle\int_{I\setminus I_n} f(t)\ dt=0$.
\end{lemma}

\begin{proof} Using the absolute continuity of the Lebesgue integral, proved in \cite[Proposition 1.12 (ii), pages 66-67]{stein2009real}, we get that:

\begin{equation}
	\left |\int_{I_n} f(t)\ dt-\int_I f(t)\ dt \right |= \left |\int_{I\setminus I_n} f(t)\ dt \right |\leq \int_{I\setminus I_n} |f(t)|\ dt\stackrel{n\to\infty}{\longrightarrow}0,
\end{equation}

\noindent because $|I\setminus I_n|\to 0$.
	
\end{proof}

\begin{lemma}\label{scheffeslemma}
	Let $\Omega\subset\mathbb{R}^N, N\geq 1$ be measurable and let $f_n,f\in L^1(\Omega)$ satisfy $f_n\geq 0$ for each $n\geq 1$ and $f\geq 0$ a.e. in $\Omega$. Assume that $\displaystyle\liminf_{n\to\infty} f_n(x)\geq f(x)$ for a.e. $x\in\Omega$ and
	\begin{equation}
		\lim_{n\to\infty}\int_\Omega f_n(x)\,dx=\int_\Omega f(x)\,dx.
	\end{equation}
	Then
	\begin{equation}
		\lim_{n\to\infty}\int_\Omega |f_n(x)-f(x)|\,dx=0.
	\end{equation}
	In other words, $f_n\to f$ in $L^1(\Omega)$.
\end{lemma}

\begin{proof}
	For every $n\in\mathbb{N}, n\geq 1$, define $g_n:=\min\{f_n,f\}$. We first show that $g_n(x)\to f(x)$ for a.e. $x\in\Omega$. Indeed, let $x\in\Omega$ be such that $\displaystyle\liminf_{n\to\infty}f_n(x)\geq f(x)$. Since $g_n(x)\leq f(x)$, we have $\displaystyle\limsup_{n\to\infty}g_n(x)\leq f(x)$. On the other hand, from Lemma \ref{lemmaliminfsupcontmon}, applied for the continuous and nondecreasing function $s\mapsto \min\{s,f(x)\}$, we get that:
	\[
	\liminf_{n\to\infty}g_n(x)
	=
	\liminf_{n\to\infty}\min\{f_n(x),f(x)\}
	=
	\min\left\{\liminf_{n\to\infty}f_n(x),f(x)\right\}
	=
	f(x).
	\]
	Therefore, $g_n(x)\to f(x)$ for a.e. $x\in\Omega$.
	
\noindent	Moreover, $0\leq g_n\leq f$ a.e. in $\Omega$. Since $f\in L^1(\Omega)$, the \textit{Lebesgue dominated convergence theorem} yields:
	\begin{equation}
		\lim_{n\to\infty}\int_\Omega g_n(x)\,dx=\int_\Omega f(x)\,dx.
	\end{equation}
	
\noindent Since $f_n,f\geq0$, we have $|f_n-f|=f_n+f-2\min\{f_n,f\}=f_n+f-2g_n$ a.e. in $\Omega$. Hence
	\[
	\int_\Omega |f_n(x)-f(x)|\,dx
	=
	\int_\Omega f_n(x)\,dx+\int_\Omega f(x)\,dx-2\int_\Omega g_n(x)\,dx.
	\]
	Passing to the limit as $n\to\infty$ and using the assumptions together with the previous convergence, we obtain
	\[
	\lim_{n\to\infty}\int_\Omega |f_n(x)-f(x)|\,dx
	=
	\int_\Omega f(x)\,dx+\int_\Omega f(x)\,dx-2\int_\Omega f(x)\,dx
	=
	0.
	\]
	Thus $f_n\to f$ in $L^1(\Omega)$.
\end{proof}

\begin{proposition}\label{pcontB} Let $\Omega\subset\mathbb{R}^N$ be an open and bounded set and $a<b$ two real numbers. Suppose that $f\in C\bigl(\overline{\Omega}\times [a,b]\bigr)$. Define the function $F:\overline{\Omega}\times [a,b]\to\mathbb{R},\ F(x,s)=\displaystyle\int_{a}^s f(x,r)\ dr$. Then $F\in  C\bigl(\overline{\Omega}\times [a,b]\bigr)$, $\dfrac{\partial F}{\partial s}(x,s)=f(x,s)$ for any $(x,s)\in\overline{\Omega}\times [a,b]$, which in particular means that $\dfrac{\partial F}{\partial s}\in C\bigl(\overline{\Omega}\times [a,b]\bigr)$.
	
\end{proposition}

\begin{proof} Fix any $(x_0,s_0)\in\overline{\Omega}\times [a,b]$. We need to show that for each $\epsilon>0$ there is some $\eta>0$ such that for any $(x,s)\in \overline{\Omega}\times [a,b]$ with $|x-x_0|\leq\eta$ and $|s-s_0|\leq\eta$ we have that $|F(x,s)-F(x_0,s_0)|\leq\epsilon$.
	
\noindent From \textit{Heine-Cantor theorem} we get that $f$ is uniformly continuous on the compact set $\overline{\Omega}\times [a,b]$. This means that for any $\epsilon'>0$ there is some $\eta'>0$ such that if $(x_1,s_1),(x_2,s_2)\in \overline{\Omega}\times [a,b]$ with $|x_1-x_2|\leq\eta'$ and $|s_1-s_2|\leq\eta'$, then $|f(x_1,s_1)-f(x_2,s_2)|\leq\epsilon'$.

\noindent Choosing $\epsilon'=\dfrac{\epsilon}{2(b-a)}$ we get that

\begin{equation}\label{barza1}
	|f(x,r)-f(x_0,r)|<\epsilon',\text{for any } |x-x_0|\leq\eta'\ \text{and any }r\in [a,b].
\end{equation}

\noindent Setting now $\eta=\min\left\{\eta',\dfrac{\epsilon}{2\displaystyle\sup_{r\in [a,b]} |f(x,r)|+1}\right\}$, we may write that for any $(x,s)\in\overline{\Omega}\times [a,b]$ with $|x-x_0|\leq\eta$ and $|s-s_0|\leq\eta$ the following inequality holds:
	
	\begin{align*}
		|F(x,s)-F(x_0,s_0)|&=\left |\int_{a}^s f(x,r)\ dr-\int_{a}^{s_0} f(x_0,r)\ dr \right | \\
		&=\left |\int_{a}^s f(x,r)-f(x_0,r)\ dr+\int_{s_0}^s f(x_0,r)\ dr \right |\\
		&\leq \int_{a}^s |f(x,r)-f(x_0,r)|\ dr+\int_{s_0}^s |f(x_0,r)\ dr\\
		&\leq \epsilon'(s-a)+(s-s_0)\displaystyle\sup_{r\in [a,b]} |f(x,r)|\\
		&\leq \epsilon'(b-a)+\dfrac{\epsilon}{2}\cdot \dfrac{\displaystyle\sup_{r\in [a,b]} |f(x,r)|}{\displaystyle\sup_{r\in [a,b]} |f(x,r)|+1}\\
		&\leq \dfrac{\epsilon}{2}+\dfrac{\epsilon}{2}=\epsilon.
	\end{align*}
	
\noindent This proves that $F\in C\bigl(\overline{\Omega}\times [a,b]\bigr)$. Now, for any fixed $x\in\overline{\Omega}$ we have that $f(x,\cdot):[a,b]\to\mathbb{R}$ is a continuous function. Therefore, from the \textit{Fundamental Theorem of Calculus} we get that $F(x,\cdot)$ is differentiable on $[a,b]$ and $\dfrac{\partial F}{\partial s}(x,s)=f(x,s)$ for any $s\in [a,b]$. The conclusion follows with ease.

\end{proof}

\begin{theorem}[\textbf{Gronwall inequality -- integral form}] Let $a,b \in \mathbb{R}$ with $a<b$, let $\alpha \in L^\infty(a,b)$ and $\beta \in L^1(a,b)$ with $\beta(t)\ge 0$ for almost all $t\in[a,b]$, and let $u\in L^\infty(a,b)$. If,\footnote{For a proof, see \cite[Lemma A.53, page 692]{john2016finite}.}
	\begin{equation}
		u(t)\le \alpha(t)+\int_a^t \beta(s)u(s)\,ds,
		\quad \text{for a.e. } t\in[a,b]
	\end{equation}
	\noindent then for almost all $t\in[a,b]$
	\begin{equation*}
		u(t)\le \alpha(t)+\int_a^t
		e^{\int_s^t \beta(\tau)\,d\tau}\,\beta(s)\alpha(s)\,ds .
	\end{equation*}
	If $\alpha\in W^{1,1}(a,b)$, it follows
	\begin{equation*}
		u(t)\le
		e^{\int_a^t \beta(\tau)\,d\tau}
		\left(
		\alpha(a)+\int_a^t e^{-\int_a^s \beta(\tau)\,d\tau }\,\alpha'(s)\,ds
		\right).
	\end{equation*}
	Moreover, if $\alpha$ is a monotonically increasing, continuous function, it holds
	\begin{equation}
		u(t)\le \exp\left (\int_a^t \beta(\tau)\,d\tau\right)\,\alpha(t).
	\end{equation}
	
\end{theorem}

\begin{theorem}[\textbf{Gronwall inequality -- differential form}]\label{gronwalldiff} Let $a,b\in\mathbb{R}$, $a<b$, $u\in W^{1,1}(a,b)$ and $\alpha,\beta\in L^1(a,b)$. If:\footnote{See Lemma A.54, page 693 from \cite{john2016finite}.}
	
	\begin{equation}
		u'(t)\leq \alpha(t)+\beta(t)u(t),\ \text{for a.e.}\ t\in (a,b),
	\end{equation}

\noindent then

\begin{equation}
u(t)\leq u(a)\exp\left (\int_{a}^t \beta(\tau)\ d\tau \right)+\int_{a}^t\exp\left (\int_{s}^t \beta(\tau)\ d\tau \right )g(s)\ ds,\ \text{for a.e.}\ t\in (a,b).
\end{equation}
\end{theorem}

\begin{theorem}\label{topothm}
	Let $X,Y$ be two topological spaces, such that $X$ is compact, and consider a continuous function $f:X\times Y\to \mathbb{R}$. Then the function $h:Y\to \mathbb{R}, h(y)=\displaystyle\inf_{x\in X} f(x,y)$ is a continuous function.
\end{theorem}

\begin{proof}
	\noindent The proof proceeds in three parts.
	
	\noindent First, we show that $h(y) > -\infty$ for every $y\in Y$. Fix $y\in Y$. Then $f(\cdot,y):X\to\mathbb{R}$ is continuous, and since $X$ is compact, it is bounded and attains its infimum. Hence $h(y)$ is finite and well-defined.
	
	\noindent Next, since the sets $(-\infty,a)$ and $(b,\infty)$ form a subbasis for the topology of $\mathbb{R}$, it suffices to show that $h^{-1}((-\infty,a))$ and $h^{-1}((b,\infty))$ are open in $Y$.
	
	\noindent Let $\pi_Y:X\times Y\to Y$ denote the canonical projection, 
	which is continuous and open. Observe that
	
	\[
	h^{-1}((-\infty,a)) 
	= \pi_Y\big(f^{-1}((-\infty,a))\big).
	\]
	
	\noindent Since $f$ is continuous, the set $f^{-1}((-\infty,a))$ is open in 
	$X\times Y$, and because $\pi_Y$ is an open map, 
	$h^{-1}((-\infty,a))$ is open in $Y$.\footnote{Note that this part does not depend on compactness of $X$. Indeed, a minor modification of this argument shows that the pointwise infimum of any family of upper semicontinuous functions is upper semicontinuous.}

	\noindent Finally, we show that $h^{-1}((b,\infty))$ is open. Here compactness of $X$ is essential. Observe that $h(y) > b$ implies $f(x,y) > b$ for all $x\in X$. Equivalently,
	
	\[
	h(y) > b \implies (x,y) \in f^{-1}((b,\infty)) 
	\quad \text{for all } x\in X.
	\]
	
	\noindent Since $f^{-1}((b,\infty))$ is open in $X\times Y$, for each $x\in X$ there exist open sets 
	$U_{(x,y)}\subset X$ and $V_{(x,y)}\subset Y$ such that
	
	\[
	(x,y)\in U_{(x,y)}\times V_{(x,y)}
	\subset f^{-1}((b,\infty)).
	\]
	
	\noindent The family $\{U_{(x,y)}\}_{x\in X}$ is an open cover of $X$. By compactness of $X$, there exist finitely many points $x_1,\dots,x_k\in X$ such that
	\[
	X \subset \bigcup_{i=1}^k U_{(x_i,y)}.
	\]
	
	\noindent Let
	\[
	V := \bigcap_{i=1}^k V_{(x_i,y)},
	\]
	
	\noindent which is an open neighborhood of $y$ in $Y$. Then for every $y'\in V$ and every $x\in X$, we have $(x,y')\in f^{-1}((b,\infty))$, so $f(x,y')>b$ for all $x$, hence $h(y')>b$. Therefore $V\subset h^{-1}((b,\infty))$, proving that $h^{-1}((b,\infty))$ is open. Since preimages of subbasic open sets are open, 
	$h$ is continuous.
\end{proof}

	\newpage
	\bibliographystyle{apalike}
	\bibliography{doubly}
\end{document}